\documentclass[11pt]{article}
\usepackage{amssymb,amsmath}
\usepackage{makeidx}
\usepackage{hyperref}
\makeindex
\newenvironment
	{example}
	{\refstepcounter{exampNum} \NI {\bf Example} \arabic{exampNum}:
	\medskip \newline}
{\newline \NI \BX}

\newenvironment
	{paragr}
	{\refstepcounter{paragNum} \NI {\bf \S \arabic{paragNum}}
	}
{}
\newenvironment
	{pictr}[4]
	{
	\refstepcounter{picnum}
	\begin{picture}(#1,#2)
	\put(#3,#4){\MB{\left(\arabic{picnum}\right)}}
	}
	{
	\end{picture}
	}

\newtheorem{newthm}{Theorem}
\newtheorem{newlem}{Lemma}
\newtheorem{newdef}{Definition}
\newtheorem{newcor}{Corollary}

\DeclareMathOperator{\CK}{\text{Cosk}}
\DeclareMathOperator{\CPL}{\text{Cpl}}
\DeclareMathOperator{\DIM}{\text{dim}}
\DeclareMathOperator{\EHULL}{\text{e-hull}}
\DeclareMathOperator{\FE}{\text{Fund}}
\DeclareMathOperator{\FSIG}{\text{Fsig}}
\DeclareMathOperator{\MX}{\overline{\text{max}}}
\DeclareMathOperator{\OB}{\text{Ob}}

\DeclareMathOperator{\PRJ}{\text{proj}}
\DeclareMathOperator{\Sig}{\text{sig}}

\DeclareMathOperator{\STS}{\text{\em Sets}}
\DeclareMathOperator{\TR}{\text{Tr}}

\DeclareMathOperator{\VL}{\text{vlist}}
\def\ALL{\forall}
\def\ATRIO#1#2#3#4#5#6{\bigl\{ (#1,#2),(#3,#4), (#5,#6) \bigr\}}

\def\BIGP#1{\Bigl( #1 \Bigr)}
\def\BLT{\bullet}
\def\BOX#1#2#3{\Lambda^{#1}(#2)(#3)}
\def\BS{\bigskip}
\def\BX{\NI \rule{5pt}{5pt}}
\def\C{\mathcal{C}}
\def\CAN#1{$#1$-cancellative}
\def\COMP#1#2{\text{comp}_{#1,#2}}
\def\DDD#1{#1 \cdots #1}
\def\del{\partial}
\def\DELSUB#1{\text{dsub}(#1)}
\def\DET#1{$#1$-determinate}
\def\DETY#1{$#1$-determinacy}
\def\DFAS{\stackrel{\text{def}}{=}}
\def\DLI#1#2#3{#1_0,\cdots,\widehat{#1_#3},\cdots,#1_{#2}}
\def\DLIST#1#2#3{#1_0, \cdots ,#1_{#3},#1_{#3}, \cdots , #1_{#2}}
\def\DMAP#1#2#3{#1 : [#2] \to [#3]}
\def\DOT#1{\bullet_#1}
\def\DPQ#1#2#3{(#1,#2).{#3}}
\def\DPQW#1#2{\{ #1,#2 \}\text{-determinate}}
\def\DREF#1{diagram (\ref{#1}.\pageref{#1})}
\def\DS#1{$\displaystyle #1 $}
\def\DTSET#1{\text{DetSet}({#1})}
\def\ELEM#1{\text{Elem}({#1})}
\def\FILL#1#2{\text{Fill}_{{#1}}({#2})}
\def\FL{^\flat}
\def\FLT#1{{#1}^\flat}
\def\FNS#1{\footnotesize{#1}}
\def\GAM{\gamma}
\def\HD#1{\NI {\bf {#1}}}
\def\HM#1#2#3{H^{#1}_{#2}({#3})}
\def\HMU#1#2#3#4{H^{{#1},{#2}}_{#3}({#4})}
\def\HT#1{\widehat{{#1}}}
\def\IMP{\Rightarrow}
\def\IMPBY#1{\overset{#1}\Longrightarrow}
\def\ISIT{\stackrel{\text{?}}{=}}
\def\ISO{\cong}
\def\LAM{\lambda}
\def\LIST#1#2{(#1_0, \cdots , #1_{#2})}
\def\LVO#1{\underset{(#1)}{-}}
\def\MB#1{\makebox(0,0){${#1}$}}
\def\MBP#1#2#3{\put(#1,#2){\MB{{#3}}}}
\def\MBPS#1#2#3{\put(#1,#2){\MBS{{#3}}}}
\def\MBS#1{\footnotesize{\makebox(0,0){${#1}$}}}
\def\MC#1{\mathcal{#1}}
\def\MIN#1{{#1}^\text{min}}
\def\MS{\medskip}
\def\MT{\emptyset}
\def\NI{\noindent}
\def\NICOMP#1#2{$(#1,#2)$-composer}
\def\OM#1#2#3{#1_0,\cdots,\underset{{#3}}{-},\cdots,#1_{#2}}
\def\omit#1{\widehat{#1}}

\def\PAR#1{\mathring{#1}}
\def\PAREL#1#2{\text{Par}_{{#1}}({#2})}
\def\PARL#1{\text{Par}_{#1}}
\def\PART#1#2{$({#1},{#2})$-partition}
\def\PLIST#1#2{{#1}_0 \cdots {#1}_{#2}}
\def\PR#1{\text{pr}_{#1}}
\def\PROD#1#2#3{\prod_{#1}^{#2}{#3}}
\def\Proof{\NI {\bf Proof:}\MS}
\def\PSMP#1#2{$({#1},{#2})$-partial simplex}
\def\qed{\NI $\square$}
\def\refpage#1{{\bf \ref{#1}}-\pageref{#1}}
\def\SACP#1#2#3{{#1}\bullet--\overset{{#3}}{-}--\bullet #2}
\def\SBS{\subseteq}
\def\SBSI{subface-simplicial}
\def\SEQ#1#2{#1_0, \cdots , #1_{#2}}
\def\SET#1#2{\bigl\{ \; #1 \; : \; #2 \; \bigr\}}
\def\SETT#1{\bigl\{  #1   \bigr\}}
\def\SH{^\sharp}
\def\SHR#1{{#1}^\sharp}
\def\SI{strictly increasing}
\def\SIG{\sigma}
\def\SL{\makebox(0,0){ }}
\def\SMAP#1#2#3#4#5#6{
  \Bigl({#3} \SBS \prod_{{#2}=0}^{#1}d_{#2}({#3}) \Bigr)
  \xrightarrow{({#5},{#6})} 
  \Bigl({#4} \SBS \prod_{{#2}=0}^{#1}d_{#2}({#4})\Bigr)
  }
\def\SPON#1#2{\text{Spon}({#1})_{{#2}}}
\def\SPS{\supseteq}
\def\SRSET#1{\text{SurjSet}({#1})}
\def\SUR#1{$#1$-surjective}
\def\SURD#1#2{#1\text{-surj}.#2} 
\def\SURY#1{$#1$\text{-surjectivity}}
\def\T#1{\scriptstyle #1}
\def\TILDE#1{\overset{\sim}{#1}}
\def\VERT#1{\text{vert}_{#1}}

\def\XRA#1{\xrightarrow{{#1}}}
\def\ZZ{\mathbb{Z}}
\newcounter{exampNum}

\newcounter{paragNum}
\newcounter{picnum}

\begin{document}

\newlength{\width}
\newcommand{\setw}[1]{\settowidth {\width}{#1}}
\newcommand{\textbox}[1]{\framebox[1.5\width]{#1}}

\title{Generalized Composition via Nerves: Models and Algebra}
\author{Paul G. Glenn\footnote{Emails: glenn\symbol{64}cua.edu\; or paulglenn47\symbol{64}gmail.com}
}

\date{}

\maketitle

\begin{abstract}
   The well-known conditions which tell when a simplicial set is the nerve of a small category generalize with respect to two parameters: the dimension $n$ of the things which compose, and the position $i$ of the thing which is the result of the composition. 
   
   In the nerve of a small category, the dimension of the things which compose is $n=1$ (i.e. its 1-simplices, the maps), and its compositions are 2-simplices (i.e. its commutative triangles) in which the position of the composite is the 1-simplex opposite vertex $i=1$.
   
These conditions generalize to all $n>1$ and $i \in \{0 , ... , n+1\}$. In this paper, a simplicial set which satisfies the generalized conditions will be called an ``\NICOMP{n}{i}".
   
   This paper explores two aspects of composers: models and the algebra of such generalized composition.
    
   With regard to models: the two-parameter generalization of composition allows a two-parameter ``theory of a composer" which generalizes the ``theory of a category". We develop a family of set-based models for composers analogous to sets-and-functions models for ordinary categories.
    
   The key to obtaining set-based models of the theory of a composer is to observe that in the nerve of a category of sets-and-functions, a 1-simplex is a binary relation with certain properties on its two faces, that is, its codomain and domain. In the nerve of a model of an \NICOMP{n}{i}, each $k$-simplex is a $(k+1)$-ary relation on its faces having certain properties, with special attention to $k=n$ (generalization of a function) and $k=n+1$ (generalization of a composition of a functions).
      
   The development of these models is guided by two constraints: (i) that the usual sets-and-functions model, i.e. when $n=1$ and $i=1$, belongs to the family; (ii) that the things which compose have properties which resemble and generalize properties possessed by ordinary functions.
   
   For the algebra of composers, we explore some notions which are inspired by the algebra of ordinary composition. These include ``comma-composers" generalizing comma categories, a composer structure for function complexes, representables, and a generalization of universal map. This part of the paper, in its current form, is provisional in content and organization.
\end{abstract}

\tableofcontents

\section{Preliminaries}

\subsection{What is in this paper?}
\label{what-is-in-paper}
The {\bf axiomatic presentation of a {category}} speaks of {objects}, {arrows}, the operations of {domain-of}, {codomain-of}, {identity-map-of} and {composition}, and equations concerning these operations. 

As is well known, one may capture this axiomatic presentation in a standard diagram like this:

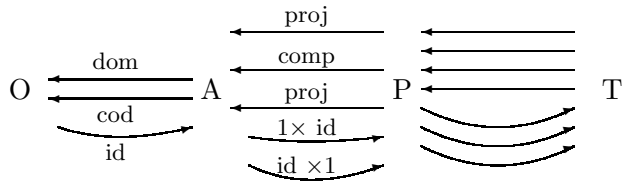
\setlength{\unitlength}{1in}
\begin{center}
\begin{picture}(3,.9)
\put(0,.5){\makebox(0,0){\cal{O}}}
\put(1,.5){\makebox(0,0){\cal{A}}}
\put(2,.5){\makebox(0,0){\cal{P}}}

\put(.9,.55){\vector(-1,0){.75}}
\put(.5,.65){\makebox(0,0){\footnotesize{dom}}}
\put(.9,.45){\vector(-1,0){.75}}
\put(.5,.37){\makebox(0,0){\footnotesize{cod}}}
\qbezier(.2,.3)(.5,.2)(.9,.3) \put(.9,.3){\vector(4,1){0}}
\put(.5,.17){\makebox(0,0){\footnotesize{id}}}

\put(1.9,.8){\vector(-1,0){.8}}
\put(1.5,.88){\makebox(0,0){\footnotesize{proj}}}
\put(1.9,.6){\vector(-1,0){.8}}
\put(1.5,.65){\makebox(0,0){\footnotesize{comp}}}
\put(1.9,.4){\vector(-1,0){.8}}
\put(1.5,.48){\makebox(0,0){\footnotesize{proj}}}
\qbezier(1.2,.25)(1.5,.2)(1.9,.25) \put(1.9,.25){\vector(4,1){0}}
\qbezier(1.2,.1)(1.5,-.05)(1.9,.1) \put(1.9,.1){\vector(4,1){0}}

\put(1.5,.3){\makebox(0,0){\footnotesize{$1 \times $ id}}}
\put(1.5,.1){\makebox(0,0){\footnotesize{id $\times 1$}}}

\put(3.1,.5){\makebox(0,0){\cal{T}}}
\put(2.9,.8){\vector(-1,0){.8}}
\put(2.9,.7){\vector(-1,0){.8}}
\put(2.9,.6){\vector(-1,0){.8}}
\put(2.9,.5){\vector(-1,0){.8}}
\qbezier(2.1,.4)(2.5,.2)(2.9,.4) \put(2.9,.4){\vector(3,1){0}}
\qbezier(2.1,.3)(2.5,.1)(2.9,.3) \put(2.9,.3){\vector(3,1){0}}
\qbezier(2.1,.2)(2.5,0)(2.9,.2) \put(2.9,.2){\vector(3,1){0}}

\end{picture}
\end{center}

\NI where \cal{O} represents {objects}, \cal{A} represents {arrows}, \cal{P} represents ``composable pairs'', \cal{T} represents ``composable trios'' 
and where {dom}, {cod}, {id} and {comp} are the operations. The ``arrows'' labelled {dom, cod, id} and {comp} in the diagram (distinct from axiomatized arrows) are a visual encoding of the operations. $P$, $T$ and the operations involving them are defined in terms of $O$, $A$, {dom}, {cod}, {id} and {comp}. 

The equations of the axiomatic presentation of a ``category'' can be stated in terms of ``following the `arrows' '' in the usual ways. These equations correspond to the {simplicial identities} in dimensions 0, 1, 2 and 3 regarding the diagram as a 3-truncated simplicial object.
\MS

This 3-truncated simplicial object can be extended to a full simplicial object/diagram amounting to a nerve-like presentation of the theory of a category. The essential information however is given by the 3-truncation and the corresponding simplicial identities.
The nerve of any specific small category has as its 3-truncation a simplicial set like the one in the diagram above.

For any small category of sets, the $0$-simplices are sets, the $1$-simplices are functions, the $2$-simplices are compositions (commutative triangles), the $3$-simplices are commutative tetrahedra (i.e. composable trios), etc. \MS
\MS

\NI {\bf Nerves of categories}
\MS

Simplicial sets which are nerves of small categories can be characterized very concisely, as follows. 

First, recall that given any simplicial set $C$, any $n \geq 1$ and any $i \in [n+1] \DFAS \SETT{0, \cdots ,n+1}$, one defines $\BOX{i}{n+1}{C}$ to consist of all 
\[
(x_0, \cdots, x_{i-1},-,x_{i+1}, \cdots , x_{n+1})
\]
where $x_j \in C_n$ for each $j \in [n+1]-\SETT{i}$ and for all $p<q$ in $[n+1]-\SETT{i}$, $d_p(x_q) = d_{q-1}(x_p)$. These are the so-called ``$i$-horns'' of dimension $n+1$.
Informally, the $x_j$'s fit together correctly so as to comprise all but the $i$'th face of a potential $(n+1)$-simplex of $C$. For each $m \geq 1$ and each $i \in [m]$ we denote by $\phi_{m,i}$ the map:
\[
\begin{array}{l}
\phi_{m,i} : C_m \to \BOX{i}{m}{C}\\
y \mapsto \Bigl(d_0(y), \cdots , d_{i-1}(y), \; \text{---} \; , d_{i+1}(y), \cdots, d_m(y)\Bigr)
\end{array}
\]
\MS

In the nerve of a small category $\C$, the set of composable pairs is $\BOX{1}{2}{\C}$ and the composition operation gives a map $\BOX{1}{2}{\C} \to \C_2$ which is an inverse of $\phi_{2,1} : \C_2 \to \BOX{1}{2}{\C}$. The composition operation is:
\[
\BOX{1}{2}{C} \xrightarrow{\phi_{2,1}^{-1}} C_2 \xrightarrow{d_1} C_1
\]
The usual picture for $\phi_{2,1}^{-1}$ is:

\begin{center}
\begin{picture}(2,.6)
\put(0,.5){\MB{V_0}}
\put(.5,.5){\MB{V_1}}
\put(.5,0){\MB{V_2}}
\put(.15,.5){\vector(1,0){.2}}
\put(.25,.6){\MBS{f_2}}
\put(.5,.35){\vector(0,-1){.2}}
\put(.6,.25){\MBS{f_0}}
\put(1,.25){\MB{\longmapsto}}
\put(1,.4){\MBS{\phi_{2,1}^{-1}}}
\put(1.5,.5){\MB{V_0}}
\put(2,.5){\MB{V_1}}
\put(2,0){\MB{V_2}}
\put(1.65,.5){\vector(1,0){.2}}
\put(1.75,.6){\MBS{f_2}}
\put(2,.35){\vector(0,-1){.2}}
\put(2.1,.25){\MBS{f_0}}
\put(1.6,.35){\vector(1,-1){.2}}
\put(1.5,.2){\MBS{f_0 \circ f_2}}
\end{picture}
\end{center}

The apparatus of simplicial algebra allows the following characterization of a (small) category, referring to the category and its nerve interchangeably:
\begin{quote}
{\bf
A {\em small category} is a simplicial set $C$ such that for all $m > 1$, the map
$\phi_{m,1} : C_m \to \BOX{1}{m}{C}$ is an isomorphism.
}
\end{quote}
\MS

\NI {\bf From categories to ``composers''}
\MS

The characterization above suggests the following generalization: regard ``1'' (the dimension of the things which ``compose'') and the other ``1'' (the position of the composite thing) as parameters and consider those simplicial sets $C$ which satisfy the following special family of {\em unique} Kan filler conditions parameterized by a given $n \geq 1$ and a given $i \in [n+1]$:
\begin{equation}
\label{intro-composer-condition}
\forall\; m>n \; \Bigl( \phi_{m,i}: C_m \to \BOX{i}{m}{C} \text{ is an isomorphism} \Bigr)
\end{equation}
In such a simplicial set $C$, $n$ is the dimension of the things that compose, and $i$ is the position of the composite thing in an $(n+1)$-simplex.
\MS

We have chosen to call any simplicial set which satisfies (\ref{intro-composer-condition}) an {\bf ``\NICOMP{n}{i}''}.
If $C$ is an \NICOMP{n}{i}, then its $(n+2)$-truncation gives a finite set of equations (namely, the simplicial identities) which generalize the equations which axiomatize a category. 
\label{axiomatic-equations}
This $(n+2)$-truncation contains all the essential information. The composite map
\[
\BOX{i}{n+1}{C} \xrightarrow{\phi_{n+1,i}^{-1}} C_{n+1} \xrightarrow{d_i} C_n
\]
generalizes composition of arrows. All this can be expressed axiomatically to present the ``theory of an \NICOMP{n}{i}''.
\MS

The family of conditions given in \eqref{intro-composer-condition} is a conceptual descendant of the family of Kan filler conditions (also unique) parameterized by the integer $n \geq 1$:
\begin{equation}
\label{hypergroupoid-condition}
\forall \; m>n, \; \forall \; i \in [m]\; \Bigl( \phi_{m,i} : C_m \to \BOX{i}{m}{C} \text{ is an isomorphism} \Bigr)
\end{equation}

Any simplicial object satisfying (\ref{hypergroupoid-condition}) is called an ``$n$-dimensional hypergroupoid'' (page \pageref{hypergroupoid}), a structure which arose as part of my dissertation research during 1973-76 under Jack Duskin's direction. (See \cite{duskin} and \cite{Glenn}).
\MS

\NI{\bf Important note:} 

There is a major distinction between the generalization of composition in the ``composer'' idea above and the comparatively much more important study of higher-dimensional categories ($\infty$-categories, quasi-categories). The family of Kan filler conditions for an $\infty$-category is quite different from the conditions in \eqref{intro-composer-condition} above: a simplicial set $C$ is an $\infty$-category iff for each $m \geq 2$ and each $i$ where $0<i<m$, the function $\phi_{m,i}: C_m \to \BOX{i}{m}{C}$ is surjective. That is, the fillers exist in all dimensions $\geq 2$ for all ``inner'' open $i$-horns, and the fillers are not in general unique. (See Lurie \cite{lurie}).
\BS

\HD{Note on internal references:}
\MS

\NI Most internal references in this paper take the form: 
\begin{quote}
{\tt [item][number]-[page number]}
\end{quote}
For example Theorem \ref{C^x-is-simplicial-set} on page \pageref{C^x-is-simplicial-set} is indicated as ``Theorem \refpage{C^x-is-simplicial-set}''.
\BS

\NI {\bf Models, sections \ref{defOfS} through \ref{model-construction}}
\MS

What about models for condition (\ref{intro-composer-condition})? Certainly, any $n$-dimensional hypergroupoid is also an \NICOMP{n}{i} (for {\em each} $i \in [n+1]$), and there are plenty of those.
\MS

The first main question addressed in this paper is: Given $n>1$ and $i \in [n+1]$, then is there a model of the theory of an \NICOMP{n}{i} which plays a role {analogous to that played by the category of sets} to the theory of a category? This paper develops a family of such models which generalize the category of sets as a model of the theory of $(1,1)$-composers.
\MS

\NI {\bf Composer algebra, sections \ref{ni-composers-algebra} through \ref{raw-material}}
\MS

The second main component of this paper is an exploration of the algebra associated with $(n,i)$-composition. This exploration is, for now, provisional in content and organization. It is guided by the search for analogs of well-known aspects of the algebra of ordinary categories. In the spirit of the somewhat informal nature of this part of the paper, we include a section of ``raw material'': ideas that look interesting with regard to $(n,i)$-composing but which need further development.

\subsection{Outline}

The key idea leading to an ``\NICOMP{n}{i} of sets'' is that an ordinary function $f : A \to B$, a $1$-simplex in the nerve of the category of sets, is a particular kind of binary relation $R_f \SBS B \times A$ on its $0$-dimensional faces $d_0(f)=B$ and $d_1(f) = A$. 

With care, one may generalize to higher dimensions the conditions characterizing those binary relations which are functions. The first step in doing so is to create a useful simplicial context for multiplace relations. 

We do this in section \ref{defOfS}, where we develop a simplicial set $\mathcal{S}$ of relations. Briefly, an $n$-simplex $y$ of $\mathcal{S}$ is an $(n+1)$-ary relation on the faces of $y$. We believe $\mathcal{S}$ may be of interest apart from providing a context for models of $(n,i)$-composition.

An ordinary, and arbitrary, $(n+1)$-place relation, $n \ge 0$, determines, and can be recovered from either of two special $n$-simplices of $\MC{S}$.
\MS

Next, (section \ref{comp}), given $n \geq 1$ and $i \in [n+1]$, we define a composition-like operation
\[
\Lambda^i(n+1)(\mathcal{S}) \to \MC{S}_{n+1} \xrightarrow{d_i} \MC{S}_n
\]
which generalizes the composition of binary relations on sets.
The key idea here, generalized to dimension $n+1$, is that a composition of ordinary functions $f_1 = f_0 \circ f_2$, corresponds to a $2$-simplex in the nerve of the category of sets:

\setlength{\unitlength}{.5in}
\begin{center}
\begin{picture}(1,1.3)
\put(0,1){\makebox(0,0){$V_0$}}
\put(1,1){\makebox(0,0){$V_1$}}
\put(.5,0){\makebox(0,0){$V_2$}}

\put(.2,1){\vector(1,0){.6}}
\put(.1,.8){\vector(1,-2){.3}}
\put(.9,.8){\vector(-1,-2){.3}}

\put(.5,1.2){\makebox(0,0){\footnotesize{$f_2$}}}
\put(.9,.5){\makebox(0,0){\footnotesize{$f_0$}}}
\put(.1,.5){\makebox(0,0){\footnotesize{$f_1$}}}
\end{picture}
\end{center}

\NI This 2-simplex is a {\em relation} on $f_0 \times f_1 \times f_2$ which expands as a relation on $(V_2 \times V_1) \; \times \; (V_2 \times V_0) \; \times \; (V_1 \times V_0)$ with certain special properties having evident generalizations to higher dimensions. That is, an ordinary composition is a certain $2$-simplex in $\mathcal{S}$. In the generalization to $n>1$, $i \in [n+1]$, an ``$(n,i)$-composition'' will be a certain kind of $(n+1)$-simplex of $\mathcal{S}$.
\MS

The technical core of this part of the paper, sections \ref{comp} -- \ref{determconds}, works out these properties and their implications, and establishes necessary conditions needed for obtaining an \NICOMP{n}{i} of sets as a subcomplex of $\mathcal{S}$. Section \ref{nisolution} then applies this material to deliver the promised models: \NICOMP{n}{i}s of sets.
\MS

Section \ref{furthermore} contains some consequences and corollaries of earlier results.
\MS

Section \ref{model-construction} examines more closely certain aspects of the models.
\MS

Sections \ref{ni-composers-algebra} through \ref{raw-material} develop some of the algebra of \NICOMP{n}{i}s.
\MS

This paper uses only basic simplicial algebra. An appendix, section \ref{simpAlg}, summarizes the ideas needed and notations.

\subsection{$n$-ary relations in this paper} 
\label{technote}

Let $n \ge 0$. In this paper, an {\bf $(n+1)$-ary relation of sets} will be understood to consist of the following information:
\begin{enumerate}
\item A function $f : [n] \to \MC{V}$, where $[n]$ is the ordered set $\{0,1,\cdots,n\}$ and $\MC{V}$ is some non-empty family of non-empty sets.
\item A subset $R \SBS \prod_{j=0}^n f(j)$.
\end{enumerate}

\MS

For brevity and readability we will speak of ``the relation $R$'', but
the phrase really includes the information above. When we mean to refer to $R$
itself as a set, we'll speak of the {\bf domain} of the relation. \index{domain of a relation}

While speaking of such relations, we will suppress the mention of
``$f$'' and use conventional notation. The function $f$ will be implicit in the notation; usually we will write something like $\prod_{j=0}^n V_j$. The relation will be written
\[
R \SBS \prod_{j=0}^n V_j \quad \text{or}  \quad R \SBS V_0 \cdots V_n
\]
Here $V_0 \cdots V_n$ is a space-saving abbreviation of $V_0 \times \cdots \times V_n$.\MS

\index{signature of a relation}
The {\bf signature} of the relation is the sequence $(f(0), \cdots ,
f(n))$. For notational compression, we will also refer to $\prod_{j=0}^n f(j)$ as the relation's signature and abbreviate this ordered product as $\Sig(R)$.
\MS

The ordering of the signature of $R \SBS \prod_{j=0}^n V_j$ is analogous to regarding a function $f : V_0 \to V_1$ as a relation on the {\em ordered} product $V_1 \times V_0$, where the ``ordering'' specifies domain $V_0$ and codomain $V_1$.
\MS

In this paper, we will be dealing with $(n+1)$-ary relations $R \SBS \prod_{j=0}^n V_j$ where
each $V_j$ is, itself, the domain of an $n$-ary relation. If $a \in R$
then $a=(a_0,\cdots,a_n) \in \prod_{j=0}^n V_j$ and, also, each $a_j$ is
$(a_{j0},\cdots,a_{j\,n-1})$. We'll examine the implications of this below in section \ref{arrays}.

\section{The simplicial set $\mathcal{S}$ of relations} \label{defOfS}

Let $\MC{V}$ be any non-empty family of non-empty sets. The simplicial set $\mathcal{S}$
defined below will refer to $\MC{V}$ but the definition will not be dependent
on any particular property of $\MC{V}$.

The models for \NICOMP{n}{i}s (see condition \eqref{intro-composer-condition} on page \pageref{intro-composer-condition}) that we seek will be subcomplexes of $\mathcal{S}$.

\subsection{$n$-simplices of $\mathcal{S}$ and face operators}

\begin{newdef}
\index{face operators of $\mathcal{S}$} 
\index{fundamental signature}
\index{$\FSIG(y)$}
\index{Fsig}
\label{Fsig-definition}
\label{simp-set-of-rels}
{\bf (Simplicial set $\MC{S}$ of relations)}
{\rm
\SL

Given any $n \ge 0$, an {\bf $n$-simplex} $y$ of $\mathcal{S}$ will be an $(n+1)$-ary relation which, in addition to its signature $\Sig(y)$ also has a so-called {\bf ``fundamental signature''} denoted $\FSIG(y)$.

\begin{enumerate} \index{signature} \index{fundamental signature}

\item {\bf Dimension 0:} The 0-simplices of $\mathcal{S}$ are unary relations $ y = (v \SBS V)$ where $V \in
\MC{V}$. By definition, $\Sig(y)=V$ and $\FSIG(y) = V$.

\item {\bf Dimension 1:} The 1-simplices of $\mathcal{S}$ are binary relations $y \SBS y_0 \times
y_1$ where $y_0$ and $y_1$ are the domains of 0-simplices $v_1 \SBS V_1$ and
$v_0 \SBS V_0$ of $\mathcal{S}$ where, for reasons which will be apparent below, $y_0=v_1$ and $y_1=v_0$. 
\[
\begin{aligned}
\Sig(y) &= y_0 \times y_1 = v_1 \times v_0\\
\FSIG(y) &= V_1 \times V_0
\end{aligned}
\]
We define
\[
\begin{aligned}
d_0(y \SBS y_0 \times y_1) &= (y_0 \SBS V_1) = (v_1 \SBS V_1)\\
d_1(y \SBS y_0 \times y_1) &=  (y_1 \SBS V_0) = (v_0 \SBS V_0)
\end{aligned}
\]

\item {\bf Dimension $n>1$:} The $n$-simplices of $\mathcal{S}$ are $(n+1)$-place
relations $y \SBS y_0 \times \cdots \times y_n$ where 
\begin{itemize}
\item $y_0, \cdots, y_n$ are the domains of $(n-1)$-simplices

\item for each $i<j$ in $[n]$,
$d_iy_j = d_{j-1}y_i$. 
\end{itemize}
We define 
\[
\begin{aligned}
\Sig(y) &= y_0 \times \cdots \times y_n\\
d_j(y \SBS \Sig(y)) &= y_j \SBS \Sig(y_j) \quad \text{each } j \in [n]\\
\FSIG(y) &=  \FSIG(y_0) \times \cdots \times \FSIG(y_n)
\end{aligned}
\]
\end{enumerate}
}
\BX
\end{newdef}
\MS

\NI {\bf Comments and notes:}
\begin{enumerate}

\item 
For notational brevity, we will usually write $y \in \mathcal{S}_n$ rather than
\[
\Bigl( y \SBS d_0(y) \times \cdots \times d_n(y) \Bigr) \in \mathcal{S}_n
\]
When referring to the set of elements of the domain of this relation, we will either use ``$y$'' when there is no ambiguity or else use ``$\ELEM{y}$'' to emphasize the distinction between the relation $y$ and the element set of its domain.
\index{element set of a simplex}
\index{Elem(y)}
\index{$\ELEM{y}$}

\item
The required face identities for $\MC{S}$, namely that for all $n>0$, all $i,j \in [n]$,
$i<j$ one has $d_i d_j (y) = d_{j-1} d_i (y)$, are built into the
definition. The degeneracy operators will be defined in the next
section.

The $i$'th vertex of $y \in \mathcal{S}_n$ is
\[
d_0 \cdots d_{i-1} d_{i+1} \cdots d_n(y) = 
d_0 \cdots \omit{d_i} \cdots d_n(y) = (v_i \SBS V_i)
\]

\item
\index{$\FSIG(V_0, \cdots , V_n)$}
\label{fsig-def-extension}
The fundamental signature is a kind of index-bookkeeping involving the vertices $(v_i \SBS V_i)$ of $y$.
Since fundamental signatures concern only the $V_i$'s, we will extend the ``Fsig'' notation and define
\[
\begin{aligned}
\FSIG(V_0) &= V_0\\
\FSIG(V_0,V_1) &= V_1 V_0\\
\FSIG(V_0, \cdots , V_n) &= \prod_{i=0}^n \FSIG(\OM{V}{n}{i})
\end{aligned}
\]
and keep in mind that the order of the factors is part of the definition.
\MS

Accordingly, if $y \in \mathcal{S}_n$ and the vertices of $y$ are $v_i \SBS V_i$, $i \in [n]$ then \DS{\FSIG(y) = \FSIG(\SEQ{V}{n})}

\item
One may, in fact, focus only on the sequence of subscripts in these ordered products and observe that subscripts in the ``Fsig'' construct are formed by an inductive concatenation:
\[
\begin{aligned}
p(0) &= 0\\
p(0,1) &= 1, 0 = p(1), p(0)\\
p(0,1,2) &= p(1,2)\; p(0,2)\; p(0,1) = 2, 1,\; 2, 0,\; 1, 0\\
\text{and}\\
p(0,1, \cdots , n) &= \bigoplus_{i=0}^n \; p(0, \cdots,\; \omit{i}\;, \cdots ,n)\\
\end{aligned}
\]
where ``$\bigoplus$'' is concatenation of sequences.

It follows directly that $\FSIG(\SEQ{V}{n})$ has $(n+1)!$ factors.

\item
Given $y \SBS \Sig(y) \in \mathcal{S}_n$ then any $a \in y$ expands to a list
\[
a = \LIST{a}{n} \SBS d_0(y) \DDD{\times} d_n(y)
\]
and each ``component'' $a_i$ also expands as
\[
a_i = (a_{i0}, \cdots , a_{i\, n-1}) \SBS
d_0 d_i(y) \DDD{\times} d_{n-1} d_i(y)
\]
The ultimate constituents of $a \in y$ are elements of the vertices $v_i \SBS V_i$ of $y$. If one were to list these ultimate entries in the correct product order, that array of ultimate entries would be an element of $\FSIG(y)$.

We call these ultimate constituents of $a$ the ``fundamental entries'' of $a$. The formal definition is below (page \pageref{fundamental-entries}).

\item
For each $n \geq 0$ and list $\SETT{\SEQ{V}{n}} \SBS \MC{V}$ there is an extremal $n$-simplex $y_\text{max}$ where $y_\text{max} = \Sig(y_\text{max})$ and for each subface $y'$ of $y_\text{max}$, $y' = \Sig(y')$.

It follows from the definitions that $y_\text{max} = \FSIG(y_\text{max}) = \FSIG \LIST{V}{n}$. If $y'$ is a subface with vertices $V_{i_0}, \cdots , V_{i_k}$ and $i_0 \DDD{<} i_k$ then $y' = \FSIG(V_{i_0}, \cdots , V_{i_k})$.

\end{enumerate}
\BX
\MS

\example{
\label{basic_2_simplex}
Consider $(y \SBS y_0 \times y_1\times  y_2) \in \mathcal{S}_2$. For each $i \in [2]$,  $y_i = d_i(y)$ and the vertices of $y$ are $(v_0 \SBS V_0) = d_1 d_2(y)$, $(v_1 \SBS V_1) = d_0d_2(y)$ and $(v_2 \SBS V_2) = d_0 d_1(y)$. 
\MS

By definition, for all $p<q$ in $[2]$, $d_p d_q(y) = d_{q-1} d_p(y)$. Therefore:
\[
\begin{aligned}
d_0(y) & \SBS  d_0 d_1(y) \times d_2 d_0 (y) \SBS V_2 \times V_1\\
d_1(y) & \SBS  d_0 d_1(y) \times d_1 d_2 (y) \SBS V_2 \times V_0\\
d_2(y) & \SBS  d_0 d_2(y) \times d_1 d_2 (y) \SBS V_1 \times V_0\\
\end{aligned}
\]
and
\[
\FSIG(y) = (V_2 \times V_1) \times (V_2 \times V_0) \times (V_1 \times V_0)
= \FSIG(V_0,V_1,V_2)
\]
A typical element of $y$ is $(a_0,a_1,a_2) \in d_0(y) \times d_1(y)\times d_2(y)$
which expands as
\[
\begin{aligned}
(a_0,a_1,a_2) &= ((a_{00},a_{01}),(a_{10},a_{11}),(a_{20},a_{21}))\\
& \in  \bigl(d_0 d_0(y) \times d_1 d_0 (y)\bigr) \times \bigl(d_0 d_1(y) \times d_1 d_1 (y)\bigr) 
\times \bigl(d_0 d_2(y) \times d_1 d_2 (y)\bigr)\\
& \SBS  (V_2 \times V_1) \times (V_2 \times V_0) \times (V_1 \times V_0)
\end{aligned}
\]
The ``fundamental entries'' of $a$ are $a_{00}$, $a_{01}$, $a_{10}$, $a_{11}$, $a_{20}$ and $a_{21}$.

Observe that $a_{00}, a_{10} \in V_2$, $a_{20}, a_{01} \in V_1$ and $a_{21}, a_{11} \in V_0$ but that {\bf it is generally not the case} that $a_{00}=a_{10}$ or $a_{20}= a_{01}$ or $a_{21}= a_{11}$. We will revisit this point below (page \pageref{arrays}).
} 

\BX
\BS

\begin{newdef}
\index{projection function $e_j$}
\index{$e_j$}
\index{e$_j$}
{\rm
\NI {\bf Projection functions $e_j^y:y \to d_j(y)$:}\MS
\SL

If $n>0$, $y \SBS \Sig(y)$ is an $n$-simplex and 
\[
a = (a_0, a_1, \cdots , a_n) \in d_0(y) \times \cdots \times d_n(y)
\]
then we'll write $e^y_j(a) \DFAS a_j$. That is, $e_j: y \to d_j(y)$ is the composite
\[
y \hookrightarrow y_0 \DDD{\times} y_n \xrightarrow{\text{proj}_j} y_j
\]
For notational brevity, we will usually just write $e_j$ instead of $e^y_j$.
\MS

\NI We will extend this notation as follows:\MS

Given any $n>0$ and sets $V_0, \cdots, V_n$ and recalling that
\[
\FSIG(V_0, \cdots V_n) \DFAS \prod_{i=0}^n \FSIG(\OM{V}{n}{i})
\]
we define
\[
e_j : \FSIG(V_0, \cdots V_n) \to \FSIG(\OM{V}{n}{j})
\]
to be projection to the $j$'th factor.
}

\BX
\end{newdef}

\BS

\begin{newdef}
{\rm
\SL

Suppose $y \in \mathcal{S}_n$ and $k \leq n$. Given 
\[
(j_n, j_{n-1},\cdots, j_k) \in [n] \times [n-1] \DDD{\times} [k]
\] 
then $d_{j_k} \cdots d_{j_n} (y)$ is an $(n-k+1)$-simplex which we will call a {\bf \em subface} of $y$.
}
\end{newdef}
\BX
\MS

\NI Some basic observations concerning $(y \SBS y_0 \times \cdots \times y_n) \in \mathcal{S}_n$:
\begin{enumerate}

\item For $n>1$, $a \in y$ and $i<j$ in $[n]$, it is, in general, the case that $e_i e_j(a) \neq e_{j-1}e_i(a)$ although both $e_i e_j(a)$ and $ e_{j-1}e_i(a)$ belong to the domain of $d_i d_j(y)$.

\item In general, $e_i : y \to d_i(y)$ is not surjective. That is, if $t$ is an element of some subface of $y$, then there need not be any element of $y$ in which $t$ appears.

\item Suppose $n \geq 1$ and $y \SBS y_0 \times \cdots \times y_n$ is an $n$-simplex. Although $a \in y$ implies $e_j(a) \in d_j(y)$, all $j \in [n]$, the elements of $y$ do not, in general, determine the elements of any of its subfaces. As an extreme case, $y$ could be the empty relation on $y_0 \times \cdots \times y_n$ while each $y_i =d_i y \SBS \Sig(d_i y)$ is non-empty.

On the other hand, if $y$ is non-empty then so are all of its subfaces.

\item Suppose $\LIST{y}{n} \in \Delta^\bullet(n)(\MC{S}) $. That is: $d_p(y_q) = d_{q-1}(y_p)$ for all $p<q$ in $[n]$.

Then {\em any} relation $y$ on $y_0 \times \cdots \times y_n$ is an $n$-simplex with $d_i(y) = y_i$, all $i \in [n]$. That is, $\mathcal{S}_n \to \Delta^\bullet(n)(\MC{S})$ is surjective.

\end{enumerate}

\subsection{Degeneracies}
\label{degen}

Let $y \in \mathcal{S}_n$ and $j \in [n]$. The simplicial identities for $d_i s_j$ inductively determine the signature of the (still to be defined) $s_j(y)$, namely
\[
\begin{aligned}
\label{sigsjy}
\Sig(s_j (y)) &= d_0 s_j (y) \times d_1 s_j (y) \times \cdots \times d_{n+1}s_j (y)\\
&= s_{j-1}d_0 (y) \times \cdots \times s_{j-1}d_{j-1} (y) \times y
\times y \times s_j d_{j+1} (y) \times \cdots \times s_j d_n (y)
\end{aligned}
\]

In order to define $s_j(y)$, we will use a family of functions defined as follows.

\begin{newdef}
\index{c$_j$}
\index{$c_j$}
{\rm 
\SL

Given $n \geq 0$ and any sets $V_0,V_1, \cdots, V_n$ from $\MC{V}$, then for each $j \in [n]$ the monic function
\[
c_j: \FSIG(\SEQ{V}{n}) \to \FSIG(\DLIST{V}{n}{j})
\]
is defined inductively on $n$ as follows.
\MS

\NI {\bf For $n=0$:}

\DS{c_0 : V_0  \to \FSIG(V_0,V_0), \quad t \mapsto (t,t)}
\MS

\NI {\bf For $n=1$:}

\DS{c_0: \FSIG(V_0,V_1) \to \FSIG(V_0,V_0,V_1), \quad 
(t_1,t_0) \mapsto (t_1,t_0,t_1,t_0,t_0,t_0)}

\DS{c_1:\FSIG(V_0,V_1) \to \FSIG(V_0,V_1,V_1), \quad 
(t_1,t_0) \mapsto (t_1,t_1,t_1,t_0,t_1,t_0)}
\MS

\NI {\bf For $n>1$:}

Let \DS{a = \LIST{a}{n} \in \FSIG(\SEQ{V}{n})} 
where for each $k \in [n]$:
\[a_k = e_k(a) \in \FSIG(\OM{V}{n}{k})
\]
Then
\DS{c_j :  \FSIG(\SEQ{V}{n}) \to \FSIG(\DLIST{V}{n}{j})} is defined by 
\[
\begin{aligned}
c_j(a) &= (c_{j-1}e_0(a), \cdots, c_{j-1}e_{j-1}(a),a,a,c_je_{j+1}(a),\cdots ,c_je_n(a))\\
&= (c_{j-1} (a_0), \cdots , c_{j-1} (a_{j-1}),a,a,c_j (a_{j+1}), \cdots , c_j (a_n))
\end{aligned}
\]
That is, $c_j$ is determined inductively by the requirement that 
\[
e_ic_j = 
\left\{
\begin{array}{ll}
c_{j-1}e_i & \text{when } i<j\\
1 & \text{when } i=j,j+1\\
c_j e_{i-1} &  \text{when } i>j+1
\end{array}
\right.
\]
}
\NI \BX
\end{newdef}

\begin{newdef}
\index{$s_j$}
\label{degen-def}
{\rm
\SL

Given $y \in \mathcal{S}_n$ and $\FSIG(y) = \FSIG(\SEQ{V}{n})$, then we define the $(n+1)$-simplex $s_j(y) \SBS \Sig(s_jy)$ by requiring the elements of $s_j(y)$ to be the image of the composite
\[
y \hookrightarrow \FSIG(\SEQ{V}{n}) \stackrel{c_j}{\longrightarrow}
\FSIG(\DLIST{V}{n}{j})
\]
and $c_j^y : y \to s_j(y)$ to be the corresponding surjection. 

\setlength{\unitlength}{1in}
\begin{center}
\begin{picture}(3,1)
\put(0,1){\makebox(0,0){$y$}}
\put(0,0){\makebox(0,0){$s_j(y)$}}
\put(-.5,0){\makebox(0,0){image $=$}}
\put(3,1){\makebox(0,0){$\FSIG(\SEQ{V}{n})$}}
\put(3,0){\makebox(0,0){$\FSIG(\DLIST{V}{n}{j})$}}
\put(.1,1){\vector(1,0){2.1}}
\put(0,.85){\vector(0,-1){.7}}
\put(3,.85){\vector(0,-1){.7}}
\put(.2,0){\vector(1,0){1.6}}
\put(3.1,.5){$c_j$}
\put(.1,.5){$c_j^y$}
\end{picture}
\end{center}

That is, $s_j(y)$ is the relation
\[
\SETT{c_j(a) : a \in y} \SBS \Sig(s_j(y))
\]
and we specify, as indicated above that
\[
d_ks_j(y) = 
\left\{
\begin{array}{ll}
s_{j-1}d_k(y) & \text{ if } k<j\\
y & \text{ if } k=j \text{ or } j+1\\
s_jd_{k-1}(y) & \text{ if } k>j
\end{array}
\right.
\]
}
\NI \BX
\end{newdef}
\MS

\begin{newlem}
{\rm
\SL

Suppose $n \geq 0$ and $\SEQ{V}{n}$ are sets from $\MC{V}$. Then
\begin{enumerate}

\item For each $j \in [n]$, $c_j : \FSIG [\SEQ{V}{n} ] \to \FSIG [V_0, \cdots, V_j,V_j, \cdots , V_n]$ is monic.

\item Given any $y \in \mathcal{S}_n$ and any $j \in [n]$, then $c^y_j : y \to s_j(y)$ is a bijection (of the domains of $y$ and $s_j(y)$).
\end{enumerate}
}
\end{newlem}

\Proof

These follow immediately from the definitions. Given any 
\[
a,a' \in \FSIG [ \SEQ{V}{n} ]
\]
then $c_j(a)=c_j(a')$ implies $e_j(c_j(a)) = e_j(c_j(a'))$ which implies $a=a'$. For the same reason, $c^y_j$ is monic, but it is also surjective, by definition.

\qed
\MS

\begin{newlem}
{\rm
\SL

Suppose $n>0$ and $i<j \in [n]$. Let $y \in \mathcal{S}_n$ and $a \in y$. Then
$c_j c_i (a) = c_i c_{j-1} (a)$.
}
\end{newlem}

\NI {\bf Proof:}\MS

The proof is by induction on $n$. The case $n=1$ is direct and trivial.
To show $c_j c_i (a) = c_i c_{j-1} (a)$ it
suffices to check that $e_k (c_j c_i (a)) = e_k (c_i c_{j-1} (a))$ for each $k \in [n+1]$. We break this down into cases and use the commutativity equations above.\MS

\NI {\bf Case} $0 \le k<i<j$:
\[
\begin{aligned}
e_k c_j c_i & =  c_{j-1}c_{i-1} e_k \text{ (identities above)}\\
&= c_{i-1}c_{j-2} e_k \text{ (induction)}\\
&= e_k c_i c_{j-1} \text{ (identities above)}
\end{aligned}
\]

\NI{\bf Case} $k=i<j$:\MS

$e_ic_jc_i = c_{j-1} e_i c_i = c_{j-1} = e_i c_i c_{j-1}$\MS

\NI{\bf Case} $i<i+1=k<j$:\MS

$e_{i+1}c_jc_i = c_{j-1} e_{i=1}c_i = c_{j-1} = e_{i+1}c_i c_{j-1}$\MS

\NI {\bf Case} $i<i+1=k=j$:\MS

$e_{i+1}c_ic_i = c_i = c_{j-1} = e_{i+1}c_ic_{j-1}$\MS

\NI {\bf Case} $i<i+2 \le k<j$:
\[
\begin{aligned}
e_k c_j c_i & =  c_{j-1}c_{i} e_{k-1} \text{ (identities above)}\\
&= c_{i}c_{j-2} e_{k-1} \text{ (induction)}\\
&= e_k c_i c_{j-1} \text{ (identities above)}
\end{aligned}
\]

\NI {\bf Case} $k=j$ or $k=j+1$:

$e_kc_j c_i=c_i = c_i e_{k-1} c_{j-1} = e_k c_i c_{j-1}$\MS

\NI {\bf Case} $k>j+1$:
\[
\begin{aligned}
e_k c_j c_i & =  c_{j}c_{i} e_{k-2} \text{ (identities above)}\\
&= c_{i}c_{j-1} e_{k-2} \text{ (induction)}\\
&= e_k c_i c_{j-1} \text{ (identities above)}
\end{aligned}
\]
\qed \BS

\begin{newthm}
{\rm
\SL

$\mathcal{S}$ is a simplicial set.
}
\end{newthm}

\NI {\bf Proof:} The verification of the simplicial identities follows
from the definitions and from the previous lemma.

\NI \qed

\subsection{Partial ordering of $n$-simplices}
\label{poset-of-n-simps}

\begin{newdef}
\label{po-def}
\index{partial order for simplices}
\index{$\leq$ for simplices}
{\rm
Suppose $y,y' \in \MC{S}_n$. We define $y \leq y'$ to mean that $\FSIG(y)=\FSIG(y')$, $\ELEM{y} \SBS \ELEM{y'}$ and that for all subface-permissible\footnote{See page \pageref{subface-permissible-for-n} for the definition of ``subface-permissible''.}
 \\$p_1 \DDD{<} p_k$,
\[
\ELEM{d_{p_1} \DDD{}d_{p_k}(y)} \SBS \ELEM{d_{p_1} \DDD{}d_{p_k}(y')}
\]
}
\BX
\end{newdef}

Given non-empty sets $\SEQ{V}{n}$ and any non-empty subset 
\[
T \SBS \FSIG \LIST{V}{n}
\]
the partially ordered set of all $n$-simplices $y$ with fundamental signature equal to $\FSIG \LIST{V}{n}$ and $\ELEM{y}=T$ has a least element $T^\text{min}$ and a greatest element $T^\text{max}$, defined below. The simplex $T^\text{min}$ will play a useful role below in section \refpage{model-construction}.
\MS

\NI{\bf Definition of $T^\text{min} \in \MC{S}_n$}

Set $\ELEM{T^\text{min}} \DFAS T$. For each subface-permissible sequence $p_1 \DDD{<} p_k$ we define
\[
\ELEM{d_{p_1} \DDD{} d_{p_k}(T^\text{min})} \DFAS
\SETT{e_{p_1} \DDD{} e_{p_k}(a):a \in y}
\]
which yields the surjective function
\[
e_{p_1} \DDD{} e_{p_k}: T^\text{min} \to d_{p_1} \DDD{} d_{p_k}(T^\text{min})
\]
\BX

\NI Note that 
\[
T^\text{min} \to \prod_{p=0}^n d_p(T^\text{min}), \qquad a \mapsto
(e_0(a) \DDD{,} e_n(a))
\]
is monic (with similar statements true for all subfaces of $T^\text{min}$), and therefore $T^\text{min} \in \MC{S}_n$. 

It is immediate that given any $y \in \MC{S}_n$ such that $\ELEM{y}=T$ and $\FSIG(y)=\FSIG \LIST{V}{n}$ then $T^\text{min} \leq y$.
\BS

\NI{\bf Definition of $T^\text{max} \in \MC{S}_n$}

Set $\ELEM{T^\text{max}} \DFAS T$. For each subface-permissible sequence $p_1 \DDD{<} p_k$ we define
\[
\ELEM{d_{p_1} \DDD{} d_{p_k}(T^\text{max})} \DFAS
V_0 \DDD{} \omit{V_{p_1}} \DDD{} \omit{V_{p_2}} \DDD{} \omit{V_{p_k}} \DDD{}  V_{n}
\]
That is, the product on the right omits factors $V_{p_1} \DDD{,} V_{p_k}$.
As in the definition of $T^\text{min}$
\[
T^\text{max} \to \prod_{p=0}^n d_p(y^\text{max}), \qquad a \mapsto
(e_0(a) \DDD{,} e_n(a))
\]
is monic (with similar statements true for all subfaces of $T^\text{max}$), and therefore $T^\text{max} \in \MC{S}_n$.

It follows directly from this definition that for all $y \in \MC{S}_n$ such that $\FSIG(y)=\FSIG \LIST{V}{n}$ and $\ELEM{y}=T$ then 
\[
\ELEM{d_{p_1} \DDD{} d_{p_k}(y)} \SBS \ELEM{d_{p_1} \DDD{} d_{p_k}(T^\text{max})}
\]
That is, $y \leq T^\text{max}$.

\BX

\subsection{Component-of-component arrays of elements of simplices}
\label{arrays}

Suppose $n \geq 1$, $y \in \mathcal{S}_n$ and $a \in y$. Then $a$ expands as a ``vector'' \DS{a = \LIST{a}{n}} with each $a_k = e_k(a) \in d_k (y)$. Similarly (if $n \geq 2$) each $a_k$ expands as \DS{a_k = (a_{k0}, \cdots , a_{k\, n-1})} with $a_{kj} = e_j e_k (a) \in d_j d_k(y)$.

Therefore, when $n \geq 2$ we may represent $a \in y$ as a rectangular array whose entries are indexed by $[n] \times [n-1]$:
\[
\left[
\begin{array}{cccc}
a_{0\,0} & a_{0\,1} & \cdots & a_{0\,n-1}\\
a_{1\,0} & a_{1\,1} & \cdots & a_{1\,n-1}\\
\vdots & \vdots & & \vdots\\
a_{n\,0} & a_{n\,1} & \cdots & a_{n\,n-1}
\end{array}
\right]
\]
For each pair of indices $k<j$, the entries $a_{jk} = e_ke_j(a)$
and $a_{j-1\,k}= e_{j-1} e_k (a)$ both belong to the same $(n-2)$-dimensional
subface $d_k d_j (y)$ of $y$ but, as noted in example \ref{basic_2_simplex}, they are {\em not necessarily}
equal. 
\MS

When $n>2$ then each of the $a_{jk}$ is itself an array so that, when fully expanded, $a$ is an $(n+1) \times n \times \cdots \times 2$
hyper-rectangular array with $(n+1)!$ entries, as observed earlier.\MS

\begin{newdef}
{\bf (Fundamental entries)}
\label{fundamental-entries}
\index{fundamental entries} 
\index{fundamental matrix}
{\rm 
\SL

The entries in the fully-expanded rectangular array of $a$ will be
called the {\bf fundamental entries} of $a$. That is, if $y \in \mathcal{S}_n$ and $a \in y$ then each fundamental entry of $a$ has the form
\[
e_{i_1} \cdots e_{i_j} \cdots e_{i_n}(a)\quad \text{where} \quad 
(i_1, \cdots, i_n) \in [1] \times \cdots \times [n]
\]
and belongs to one of the vertices of $y$.

This hyper-rectangular array
will be called the {\bf fundamental matrix} of $a$.

\NI \BX}
\end{newdef}\MS

Note that if $y \in \mathcal{S}_n$ and $a \in y$ then $e_0 \cdots \omit{e_i} \cdots e_n(a)$ belongs to the $i$'th vertex of $y$.

Also, if $d_{j_1} \cdots d_{j_n} = d_0 \cdots \omit{d_i} \cdots d_n$ is a face identity then $e_{j_1} \cdots e_{j_n}(a)$ belongs to the $i$'th vertex of $y$, though in general it is {\em not} equal to $e_0 \cdots\omit{e_i} \cdots e_n(a)$.\BS

\begin{newdef}
\index{subelement} 
\index{component}
\label{fundEntries}
{\rm
\SL

Given $y \in \mathcal{S}_n$ where $n \ge 1$ and $a \in y$, we'll call $e_{j_1} \cdots e_{j_k} (a) \in d_{j_1} \cdots d_{j_k} (y)$
a {\bf \em subelement} or {\bf \em component} of $a$.

Here, $(j_k, j_{k-1}, \cdots, j_1) \in [n] \times [n-1] \times \cdots \times [n-k+1]$. For notational brevity where convenient, we'll write $A = (j_k, j_{k-1}, \cdots, j_1)$ and the subelement as $e_A(a)$.
\NI \BX}
\end{newdef}\MS

\begin{newdef} 
\index{partial element}
\label{defPartialElement}
{\rm
\SL

A {\bf \em partial element}, or ``compatible family'', of $y$ is a
family of elements of some subfaces of $y$ such that there exists an
$\tilde{a} \in \FSIG(y)$ with the property that whenever $b \in d_{A}y$
belongs to the given family then $b = e_A \tilde{a}$

Note that, in general, even if it exists, that $\tilde{a}$ need not belong to $y$.

\NI \BX}
\end{newdef} \MS

\NI Notes:

\begin{enumerate}
\item There is {\em no essential difference} in speaking of $a \in y$ as a matrix  and as an element in the product $\FSIG(y)$. Each fundamental entry of $a$ is a list item in an element of $\FSIG(y)$.

\item One may visualize a partial element of $y \in \mathcal{S}_n$ as a
fundamental matrix with some missing entries.
\end{enumerate}

\section{Composition}
\label{comp}

In this section, guided by ordinary composition of functions (or binary relations), we define an ``$(n,i)$-composition'' as an $(n+1)$-simplex in $\MC{S}$ with certain properties, and develop some implications of those properties. We use the definition of $(n,i)$-composition to specify, in section \refpage{wantcompstruct}, exactly what is meant by a model of the ``theory of an \NICOMP{n}{i}'': a subcomplex of $\MC{S}$ with certain properties.

\subsection{Example from composition of binary relations}
\label{binaryrelcomp}

To anticipate the higher-dimensional notion of composition defined below, first consider ordinary composition of binary relations (and also
of functions) in terms of the ``composition relation''. That is, suppose
there are binary relations, 1-simplices of $\mathcal{S}$:
\[
\begin{aligned}
y_0  \SBS  d_0(y_0) \times d_1 (y_0) \SBS V_2 \times V_1\\
y_2  \SBS  d_0(y_2) \times d_1 (y_2) \SBS V_1 \times V_0
\end{aligned}
\]
where
$d_0(y_2) =d_1(y_0)$. The composite of these two relations, possibly empty, is:
\[
y_1 \SBS d_0(y_1) \times d_1(y_1) = d_0(y_0) \times d_1(y_2)
\]
where $y_1$ consists of all pairs $(a_2,a_0) \in V_2 \times V_0$ for which there
exists some $a_1 \in d_0(y_2)$ such that $(a_2,a_1)\in y_0$ and $(a_1,a_0) \in y_2$.\MS

In a conventional diagram:

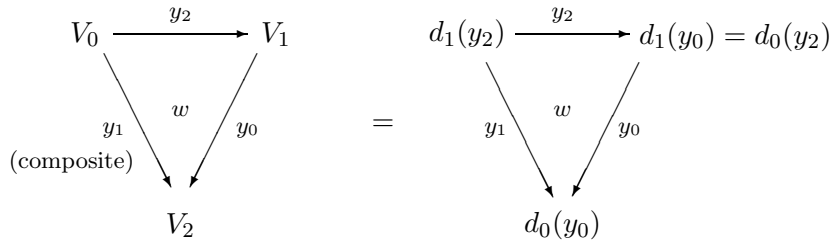
\setlength{\unitlength}{1in}
\begin{center}
\begin{picture}(3,1.2)
\put(0,1){\MB{V_0}}
\put(1,1){\MB{V_1}}
\put(.5,0){\MB{V_2}}
\put(.5,1.1){\MBS{y_2}}
\put(.85,.5){\MBS{y_0}}
\put(.15,1){\vector(1,0){.7}}
\put(.9,.9){\vector(-1,-2){.35}}
\put(.1,.9){\vector(1,-2){.35}}
\put(.5,.6){\MBS{w}}
\put(.15,.5){\MBS{y_1}}
\put(-.4,.3){\footnotesize{(composite)}}
\put(1.5,.5){$=$}

\put(2,1){\MB{d_1(y_2)}}
\put(2.9,.95){{$d_1(y_0) = d_0(y_2)$}}
\put(2.5,0){\MB{d_0(y_0)}}
\put(2.5,1.1){\MBS{y_2}}
\put(2.85,.5){\MBS{y_0}}
\put(2.25,1){\vector(1,0){.6}}
\put(2.9,.85){\vector(-1,-2){.35}}
\put(2.1,.85){\vector(1,-2){.35}}
\put(2.5,.6){\MBS{w}}
\put(2.15,.5){\MBS{y_1}}
\end{picture}
\end{center}

Expressed in terms of $\mathcal{S}$, the ``composition relation'' is a 2-simplex
\[
w \SBS d_0 (w) \times d_1 (w) \times d_2(w) = y_0 \times d_1 (w) \times y_2 \in \mathcal{S}_2
\]
with $d_1(w)$ equal to the composite of $y_0$ with $y_2$. A typical element $b=(b_0,b_1,b_2) \in w$ expands as
\[
b=(b_0,b_1,b_2)=
\left[
\begin{array}{cc}
b_{00} & b_{01}\\
b_{10} & b_{11}\\
b_{20} & b_{21}
\end{array}
\right]
= 
\left[
\begin{array}{cc}
a_2 & a_1 \\
a_2 & a_0 \\
a_1 & a_0
\end{array}
\right]
\begin{array}{c}
\text{ (element of } y_0 \text{)}\\
\text{ (element of } d_1(w) \text{)}\\
\text{ (element of } y_2 \text{)}
\end{array}
\]
and $w$ has the property that the 1-simplex $d_1(w) \SBS V_2 \times V_0$
contains $(a_2,a_0) \in V_2 \times V_0$ {\em if and only if} there exists
$b \in w$ such that $(a_2,a_0) = (b_{10},b_{11}) = e_1(b)$. 

This captures the conventional definition of {\em function} composition: the ordered pair $(a_0,a_2)$ belongs to the composite function $y_1 = y_0 \circ y_2$ if and only if there exists an element $a_1 \in V_1$ such that the ordered pair $(a_1,a_0)$ belongs to the function $y_2$ (that is, $y_2(a_0)=a_1$) and the ordered pair $(a_2,a_1)$ belongs to the function $y_0$ (that is, $y_0(a_1)=a_2$). 
\MS

As a 2-simplex of $\mathcal{S}$, $w$ is characterized by three special properties:
\begin{enumerate}
\item For all $b \in w$, it is the case that $b_{00}=b_{10}$,
$b_{01}=b_{20}$ and $b_{11}= b_{21}$. That is, for all $b \in w$ it is the case that for all $j<k$ in $[2]$, $e_j e_k(b) = e_{k-1} e_j(b)$.

\item $(a_2,a_0) \in d_1(w)$ if and only if there exists $b \in w$ such that
$b_{10} = a_2$ and $b_{11} = a_0$.
\item $w$ contains all $b=(b_0,b_1,b_2)$ such that $b_0 \in y_0, b_2 \in y_2$ and $b_{01}=b_{20}$.
\end{enumerate}
We also observe that even if some 2-simplex satisfies the first condition, then it might not satisfy the second or third one.\MS

In summary, ordinary composition of binary relations (dimension $1$ simplices of $\mathcal{S}$) gives rise to a simplex $w \in \mathcal{S}_2$ which has these three properties and which we might call the ``composition relation'' or ``composition simplex'' for
the composition.\MS

In the next section we will generalize to dimensions $n \geq 2$, the first of the properties of $w$ listed above.

\subsection{The ``component-simplicial'' and ``element-simplicial'' properties}
\label{simplicialProp}

When we defined the projection functions 
\[
e_k : y \to d_k(y), \quad a = \LIST{a}{n} \mapsto a_k
\]
we observed that, in general, the face-like identity ``$e_i e_j = e_{j-1} e_i$'' (when $i<j$) need {\em not} hold. As we have just seen however, these equations do hold for the ``composition'' simplex. In the next few sections, we will develop the consequences of imposing such equational conditions on $n$-simplices for $n \geq 2$.

\begin{newdef} 
\index{element-simplicial hull} 
\index{component-simplicial} 
\index{$A$-indexed partial element}
\index{A-indexed partial element}
\index{indexed partial element}
\index{element-simplicial simplex}
\index{e-simplicial}
\index{maximal e-simplicial}
\index{$M(\SEQ{x}{n})$}
\index{component-simplicial partial element}
\index{$\PARL{A}(y)$}
\label{defSimplicial}
{\rm
\SL

\begin{enumerate}

\item Suppose $n \geq 2$ and $y \in \mathcal{S}_n$. An element $a \in y$ is said to be {\bf component-simplicial} if for all $j<k \in [n]$, it is the case that $e_j e_k(a) = e_{k-1} e_j(a)$.

\item Suppose $n \geq 2$, $y \in \mathcal{S}_n$ and $A \subsetneq [n]$ contains at least two elements. Then a {\bf component-simplicial partial element of $y$} is a set
\[
\SETT{a_j \; : \; j \in A \text{ and } a_j \in d_j(y)}
\]
such that for all $p,q \in A$ with $p<q$ then $e_p(a_q) = e_{q-1}(a_p)$.

We will also refer to this as an {\bf $A$-indexed partial element}.
\MS

We denote the set of $A$-indexed partial elements of $y$ by $\PARL{A}(y)$.

\item Suppose $n \geq 2$ and $(x_0, \cdots,x_n) \in \Delta^\bullet(n)(\MC{S})$. That is, for all $p<q$ in $[n]$, $d_p(x_q) = d_{q-1}(x_p)$.
Then the {\bf element-simplicial hull} of $(x_0, \cdots,x_n)$ is
\[
\EHULL \LIST{x}{n} \DFAS 
\Bigl\{
a \in \prod_{j=0}^n x_j \; : \; 
\forall \; p<q \in [n], \; 
e_p e_q (a) = e_{q-1}e_p(a)
\Bigr\}
\]
In terms of arrays, $a \in \EHULL \LIST{x}{n}$ if and only if 
\[
a = 
\left[
\begin{array}{c}
a_0 \\ 
a_1 \\
\vdots \\
a_n
\end{array} 
\right]
=
\left[
\begin{array}{ccc}
a_{00} & \cdots & a_{0\, n-1}\\
a_{10} & \cdots & a_{1\, n-1}\\
\vdots & & \vdots \\
a_{n0} & \cdots & a_{n\, n-1}\\
\end{array}
\right]
\]
has $a_{qp} = a_{p \, q-1}$ whenever $0 \leq p<q \leq n$.

We can regard $\EHULL \LIST{x}{n}$ as an $n$-simplex of $\mathcal{S}$ where 
\[
d_j\BIGP{\EHULL \LIST{x}{n}} = x_j
\]

\item The non-empty simplex $y \in \mathcal{S}_n$ is {\bf \em element-simplicial} if every $a \in y$ is component-simplicial. Non-empty element-simplicial $y$ do exist; see below.
{\bf We will abbreviate the term ``element-simplicial'' by ``e-simplicial''.}
\index{e-simplicial}

\item
Given $y \in \mathcal{S}_n$, the {\bf e-simplicial hull of $y$} is defined by
\[
\EHULL (y) \DFAS \EHULL (d_0(y),\cdots ,d_n(y))
\]
$\EHULL(y)$ can be construed as an $n$-simplex with $\Sig(\EHULL(y)) = \Sig(y)$.

\item $y \in \mathcal{S}_n$ is {\bf maximal e-simplicial} if $y$ consists of all $a \in \Sig(y)$ which are component-simplicial. That is, $y = \EHULL(y)$ as $n$-simplices.
\end{enumerate}

\NI \BX}
\end{newdef}

\NI Note: The element-sets of $y$ and $\EHULL(y)$ are both subsets of $\Sig(y)$ but there is no other relationship, in general. Either could be empty. Even if both are non-empty, they could be disjoint (or not) and one could be a subset of the other (or not).

In dimensions 0 and 1, the condition for being ``e-simplicial''
is vacuous.
\MS

\begin{newlem}
\label{degenerates are maximal e-simplicial}
{\rm
Let $n \geq 2$, $y \in \MC{S}_n$ be e-simplicial and $k \in [n]$. Then $s_k(y)$ is maximal e-simplicial.
}
\end{newlem}

\Proof

It follows from the definition of $c_k:y \to s_k(y)$ that for all $a \in y$, $c_k(a)$ is component-simplicial i.e. $c_k(a) \in \EHULL(s_k(y))$. To complete the proof we will show that $\EHULL(s_k(y)) \SBS s_k(y)$.
\MS

Let $\LIST{b}{n+1} \in \EHULL(s_k(y))$.
Then $b_k, b_{k+1} \in y$ and, for $j \neq k, k+1$
\[
b_j = \left\{
\begin{array}{ll}
c_{k-1}(a^j) & \text{ for some }a^j \in y,\quad j<k\\
c_k(a^{j-1}) & \text{ for some }a^{j-1} \in y,\quad j>k+1
\end{array}
\right.
\]
From the conditions $e_i(b_j)=e_{j-1}(b_i)$ the components of $b_k$ are
\[
e_i(b_k) = \left\{
\begin{array}{ll}
e_{k-1}(b_i) = e_{k-1}c_{k-1}(a^i) = a^i & i<k\\
e_k(b_{i+1}) = e_k c_k(a^i) = a^i & i>k+1
\end{array}
\right.
\]
Therefore, 
\[
b_i = \left\{
\begin{array}{ll}
c_{k-1}(a^i) = c_{k-1}e_i(b_k) & i<k\\
c_k(a^{i-1}) = c_k e_{i-1}(b_k) & i>k+1
\end{array}
\right.
\]
Also
\[
e_i(b_{k+1}) = \left\{
\begin{array}{ll}
e_k(b_i) = e_k c_{k-1}(a^i) = a^i = e_i(b_k) & i<k\\
e_k(b_k) & i=k\\
e_k(b_{i+1}) = a^i = e_i(b_k) & i>k
\end{array}
\right.
\]
That is, $b_{k+1}=b_k$.
This shows:
\[
\LIST{b}{n+1} = 
\bigl(
c_{k-1}e_0(b_k), \cdots, c_{k-1}e_{k-1}(b_k),b_k,b_{k},c_ke_{k+1}(b_k), \cdots , c_ke_n(b)
\bigr) = c_k(b_k)
\]
\qed
\MS

\begin{example}
If $y \in \mathcal{S}_3$ and $y$ is e-simplicial, then the matrix of a typical
element $a \in y$ would look like this:
\[
\left[
\begin{array}{ccc}
a_{10} & a_{20} & a_{30}\\
a_{10} & a_{21} & a_{31}\\
a_{20} & a_{21} & a_{32}\\
a_{30} & a_{31} & a_{32}
\end{array}
\right]
\]
\end{example}

Note that when $n>3$, $i \in [n]$ and even if $y \in \mathcal{S}_n$
is e-simplicial, then $d_i (y)$ need not be e-simplicial.

The following fact, though immediate from the definition of component-simplicial, is a key to generalizing the notion of function composition.

\begin{newlem}
{\bf (Unique filler)}
\label{unique-determine-lemma}
\index{unique filler lemma}
\index{lemma!Unique filler}
{\rm
\SL

Suppose $y \in \mathcal{S}_n$, with $n \geq 2$, $i \in [n]$ and that $a = \LIST{a}{n} \in y$ is component-simplicial.
Then $a_i$ is determined uniquely by the $a_j$, $j \neq i$.
}
\end{newlem}

\Proof

For any $a = \LIST{a}{n} \in y$ and any $j \in [n]$,
\[
\begin{aligned}
a_i &= \Bigl( e_0(a_i), \cdots , e_{n-1}(a_i) \Bigr) \\
&= \Bigl( e_{i-1}(a_0), \cdots , e_{i-1}(a_{i-1}), e_i (a_{i+1}), \cdots, e_i(a_n) \Bigr)
\end{aligned}
\]

\qed

\subsection{Definition of $(n,i)$-composition simplices of $\mathcal{S}$}
\label{compSimps}

\begin{newdef}
{\bf ($(n,i)$-composition)}
\label{defCompSimps} 
\label{n-i-comp-definition}
\index{$(n,i)$-composition}
{\rm
\SL

Suppose $n\ge 1$ and $i \in [n+1]$. Then $w \in \mathcal{S}_{n+1}$ is said to be an  {\bf $(n,i)$-composition} if:
\begin{enumerate}

\item 
$w$ is e-simplicial.

\item \label{operation-condition}
Given any component-simplicial partial element $(\OM{b}{n+1}{i})$ of $w$, let $b_i \in \Sig(d_i(w))$ be the unique element defined by
\[
e_p(b_i) \DFAS 
\left\{
\begin{array}{cc}
e_{i-1}(b_p) & \text{ if } p<i\\
& \\
e_i(b_{p+1}) & \text{ if } p \geq i
\end{array}
\right.
\]
That is,
\[
b_i \DFAS \Bigl(e_{i-1}(b_0), \cdots , e_{i-1}(b_{i-1}), e_i(b_{i+1}), \cdots , e_i(b_{n+1})\Bigr)
\]
and $b = \LIST{b}{n+1} \in \Sig(w)$. Then $b \in w$.

\item \label{surjectivity-condition}
$a \in d_i(w)$ if and only if there exists $b \in w$ such that $e_i(b)=a$.

\end{enumerate}

(Make the obvious notational adjustments in condition \ref{operation-condition} for $i=0$ or $i=n+1$).

}
\BX
\end{newdef}
\NI {\bf Remarks unpacking this definition:}
\begin{enumerate} \label{unpack-comp-def}

\item Condition \ref{operation-condition} refers back to lemma \refpage{unique-determine-lemma}, which states that whenever 
\[
(\OM{b}{n+1}{i})
\]
is a component-simplicial partial element of $w$, there is a unique $b_i \in \Sig(d_i(w))$ such that 
\[
(b_0, \cdots ,  b_i , \cdots , b_{n+1}) \in \FSIG(w)
\]
is component-simplicial. 
In general, if $w$ is not an $(n,i)$-composition then $b_i$ need not be an element of $d_i(w)$. Moreover, even if $b_i \in d_i(w)$, then $b = \LIST{b}{{n+1}}$ need not be an element of $w$.

Condition \ref{operation-condition} is therefore an operation-like property of $w$ stating that $b \in w$ and $b_i \in d_i(w)$.

\item The ``only if'' part of condition \ref{surjectivity-condition} means that $e_i : w \to d_i(w)$ is surjective. This is an example of a ``surjectivity property'' which will play a role later.
\label{i-surj}

\item
A (1,1)-composition is a composition of binary relations as described in section \refpage{binaryrelcomp}.

\item
Given $n \geq 1$ and $i \in [n+1]$, there is a partial operation for $(n,i$)-composition analogous to (and generalizing) the partial operation of composing binary relations. See Definition \refpage{nicomp} below. 

\end{enumerate}
\BX
\MS

\begin{example}
To visualize $(n,i)$-composition in the case $(n,i)=(2,2)$, suppose $w \in
\mathcal{S}_3$ is a (2,2)-composition. Then $a = (a_0,a_1,a_2) \in d_2(w)$
iff there is an element $b \in w$ whose array
is
\[
\left[
\begin{array}{ccc}
b_{10} & b_{20} & b_{30}\\
b_{10} & b_{21} & b_{31}\\
b_{20} & b_{21} & b_{32}\\
b_{30} & b_{31} & b_{32} 
\end{array}
\right] = 
\left[
\begin{array}{ccc}
b_{10} & a_0 & b_{30}\\
b_{10} & a_1 & b_{31}\\
a_0 & a_1 & a_2\\
b_{30} & b_{31} & a_2
\end{array}
\right]
\begin{array}{c}
\text{ row 0} \in d_0(w)\\
\text{ row 1} \in d_1(w)\\
\text{ row 2} \in d_2(w)\\
\text{ row 3} \in d_3(w)
\end{array}
\]
From the partial-operation viewpoint, consider
\[
(y_0,y_1,-,y_3) \in \BOX{2}{3}{\MC{S}}
\]
Let $w \SBS y_0 \times y_1 \times d_2(w) \times y_3$
be the $(2,2)$-composition, where
\[
d_2(w) \SBS d_0 d_2(w) \times d_1 d_2 (w) \times d_2 d_2(w) = 
d_1(y_0) \times d_1 (y_1) \times d_2(y_3)
\]
The elements of $w$ are all $b=(b_0,b_1,b_2,b_3) \in \Sig(w)$ with matrix
\[
\left[
\begin{array}{ccc}
b_{10} & b_{20} & b_{30}\\
b_{10} & b_{21} & b_{31}\\
b_{20} & b_{21} & b_{32}\\
b_{30} & b_{31} & b_{32} 
\end{array}
\right]
\]
and
\[
d_2(w) = \SET{a}{\exists b \in w \text{ with } a = e_2(b)}
\]
That is, $a \in d_2(w)$ iff there exists $b \in w$ with $a=b_2 = (b_{20},
b_{21}, b_{32})$.
\end{example}

\begin{newthm}
{\bf (Maximal Property)}
\label{max-property} 
\index{maximal property theorem}
\index{Theorem!Maximal Property}
{\rm
\SL

If $w \in \mathcal{S}_{n+1}$ is an $(n,i)$-composition then $w$ is maximal e-simplicial.
}
\end{newthm}

\Proof

If $b \in \Sig(w)$ is component-simplicial i.e. $b \in \EHULL(w)$ then
\[
(\OM{b}{n+1}{i})
\]
is a component-simplicial partial element of $w$ and the component $b_i = e_i(b)$ is uniquely determined according to lemma \refpage{unique-determine-lemma}. Condition \ref{operation-condition} of the definition implies $b \in w$. Therefore $w$ is maximal e-simplicial.

\qed
\BS

\begin{example}
It is possible for $w \in S_n$ to be maximal e-simplicial but not an $(n,i)$-composition, even in dimension $n=1$. 

To illustrate this, suppose $(y_0,-,y_2) \in \Lambda^1(2)(\MC{S})$; that is, $d_0(y_2) = d_1(y_0)$. Suppose $b_0 \in y_0$ and $b_2 \in y_2$ with $e_0(b_2) = e_1(b_0)$. That is, $(b_0, - , b_2)$ is a component-simplicial partial element.

Holding $b_0,b_2$ {\em fixed}, let $y_1 \SBS d_0(y_0) \times d_1(y_2)$ be any 1-simplex with $(e_0(b_0),e_1(b_2)) \notin y_1$.
\MS

Now define $w \in \mathcal{S}_2$ by $w = \EHULL(y_0,y_1,y_2)$ with $d_j(w) \DFAS y_j$; so $w$ is maximal e-simplicial. A typical element of $w$ is $b' = (b'_0,b'_1,b'_2)$ where $b'$ is component-simplicial and $b'_j \in y_j$, all $j \in [2]$.

Note that $\Sig(d_1(w)) = d_0(y_0) \times d_1(y_2) = \Sig(y_1)$ and that $(b_0,-,b_2)$ is a component-simplicial partial element of $w$.

Then $(b_0,-,b_2)$ of $w$ has, as unique filler in $\Sig(w)$:
\[
\left[
\begin{array}{c}
b_0\\
b_1\\
b_2
\end{array}
\right] = 
\left[
\begin{array}{cc}
e_0(b_0) & e_1(b_0)\\
e_0(b_0) & e_1(b_2)\\
e_0(b_2) & e_1(b_2)
\end{array}
\right]
\]
Since $(e_0(b_0),e_1(b_2)) \notin y_1$ then $w$ fails to satisfy condition \ref{operation-condition} of the definition of $(1,1)$-composition.
\end{example}
\MS

The following definition generalizes function (or binary relation) composing to $(n,i)$-composition with $n>1$. The details are a direct consequence of the face identities, the requirement that an $(n,i)$-composition be e-simplicial and the requirement that $e_i: w \to  d_i(w)$ be surjective.

\begin{newdef}
\index{$(n,i)$-composing}
\index{$\COMP{n}{i}(\text{--})$}
\label{nicomp}
{\rm
\SL

Given a fixed $n \geq 1$ and fixed $i \in [n+1]$, suppose
\[
\Bigl(
\OM{y}{n+1}{i}
\Bigr) \in \Lambda^i(n+1)(\MC{S})
\]
Then there is an $(n,i)$-composition $w \in \MC{S}_{n+1}$ is defined as follows.
\begin{enumerate}

\item
For all $p \in [n+1]-\SETT{i}$, $d_p(w) = y_p$.

\item
\[
\begin{aligned}
\Sig(d_i(w)) & = d_0 d_i(w) \times \cdots \times d_n d_i(w) \\
& \DFAS  d_{i-1}(y_0) \times \cdots \times d_{i-1}(y_{i-1}) \times
d_i (y_{i+1}) \times \cdots \times d_i(y_{n+1})
\end{aligned}
\]

\item Whenever $(\OM{b}{n+1}{i}) \in \prod_{p \neq i} y_p$ such that for all $p<q$ in $[n+1]- \SETT{i}$, $e_p(b_q) = e_{q-1}(b_p)$ then
\[
a = \bigl(
e_{i-1}(b_0), \cdots , e_{i-1}(b_{i-1}), e_i(b_{i+1}, \cdots e_i(b_{n+1})
\bigr) \in d_i(w)
\]
and $(b_0, \cdots, b_{i-1},a,b_{i+1}, \cdots , b_{n+1}) \in w$.

The ``end cases'' are:

If $i=0$ then $a = (e_0(b_1), \cdots , e_0(b_{n+1}))$.

If $i=n+1$ then $a = (e_n(b_0), \cdots, e_n(b_n))$.

\item
$a \in d_i(w)$ iff there exists $b \in w$ such that $e_i(b) = a$.

\end{enumerate}

We denote $w$ by $\COMP{n}{i} ( \OM{y}{n+1}{i} )$ and call the operation
\begin{eqnarray*}
\COMP{n}{i} : \BOX{i}{n+1}{\MC{S}} & \to &  \MC{S}_{n+1}\\
(\OM{y}{{n+1}}{i}) & \mapsto & w
\end{eqnarray*}
just described {\bf \em $(n,i)$-composing}, where the factors are the $y_p$ ($p \neq i$), and $d_i(w)$ the resulting {\bf $(n,i)$-composite}.\\
}
\BX
\end{newdef}
\BS

Now suppose $w \in \mathcal{S}_{n+1}$ is e-simplicial but is {\em not} an $(n,i)$-composition. How then might $w$ relate to $\TILDE{w} = \COMP{n}{i}(d_0(w), \cdots , \underset{i}{\text{---}}, \cdots , d_{n+1}(w))$?  

Directly from the definitions it follows that, for all $p \neq i$, $d_p(\TILDE{w}) = d_p(w)$, and, for all $p \in [n]$, $d_p d_i(\TILDE{w}) = d_p d_i(w)$ and therefore, $d_i(\TILDE{w})$ and $d_i(w)$ have equal signatures. 

In general, $d_i(\TILDE{w})$ and $d_i(w)$ are unequal both as sets and as $n$-simplices. The inequality can arise in the following ways:
\begin{itemize}
\item 
Suppose $( \OM{b}{n+1}{i} )$ is any component-simplicial partial element of $w$. By lemma \refpage{unique-determine-lemma} above, it has a unique filler $b = \LIST{b}{n+1} \in \FSIG(w)$. By definition of $(n,i)$-composition, $b \in \TILDE{w}$. In general however, $b \notin w$ and $b_i = e_i(b) \notin d_i(w)$.

Thus, $\TILDE{w}$ may contain elements not in $w$ and $d_i(\TILDE{w})$ may contain elements not in $d_i(w)$.

\item
By definition of $(n,i)$-composition, each $t \in d_i(\TILDE{w})$ arises as $e_i(b)$ for some $b \in \TILDE{w}$. Since, in general, $w \xrightarrow{e_i} d_i(w)$ need {\em not}  be surjective then
$d_i(w)$ may contain an element $t'$ for which there is no $b' \in w$ with $e_i(b') = t'$. Thus, $d_i(w)$ may contain elements not in $d_i(\TILDE{w})$ and $w$ may contain elements not in $\TILDE{w}$.
\end{itemize}

If, however, we assume $w \in \mathcal{S}_{n+1}$ is e-simplicial {\em and} that $e_i : w \to d_i(w)$ is surjective then $w$ and $\TILDE{w}$ relate according to the following lemma.

\begin{newlem}
{\bf (Expansion)}
\label{sub-composition}
\index{Expansion lemma}
\index{lemma!Expansion lemma}
{\rm
\SL

Let $n \geq 2$ and $i \in [n+1]$. Suppose $w \in \mathcal{S}_{n+1}$ is e-simplicial and that $w \xrightarrow{e_i} d_i(w)$ is surjective. Let
\[
\TILDE{w} \DFAS \COMP{n}{i}\left(d_0(w), \cdots , \underset{i}{-} , \cdots , d_{n+1}(w)\right)
\]
Then $\Sig(d_i(w)) = \Sig(d_i(\TILDE{w}))$, $w \SBS \TILDE{w}$ and $d_i(w) \SBS d_i(\TILDE{w})$.
}
\end{newlem}
\MS

\Proof

As observed above, for each $p \in [n]$
\[
d_p d_i(w) = 
\left\{
\begin{array}{ll}
d_{i-1} d_p(w) =  d_{i-1} d_p(\TILDE{w}) = d_p d_i( \TILDE{w}) & \text{if } p<i\\
& \\
d_i d_{p+1}(w) = d_i d_{p+1}(\TILDE{w}) = d_p d_i( \TILDE{w})& \text{if } p \geq i
\end{array}
\right.
\]
That is, $\Sig(d_i(w)) = \Sig(d_i(\TILDE{w}))$.
\MS

Given any $t \in d_i(w)$, then the surjectivity assumption implies there exists $b = \LIST{b}{n+1} \in w$ such that $e_i(b) = b_i = t$. Since $w$ is e-simplicial then $(\DLI{b}{n+1}{i})$ is a component-simplicial partial element of $\TILDE{w}$ and, by definition of $\TILDE{w}$ and by the Unique filler lemma above, then $t = b_i \in d_i(\TILDE{w})$. Therefore $d_i(w) \SBS d_i(\TILDE{w})$.

By essentially the same observation, given any $b = \LIST{b}{n+1} \in w$ then $(\DLI{b}{n+1}{i})$ is a component-simplicial partial element of $\TILDE{w}$ filled uniquely by $b_i$. By definition of $(n,i)$-composition, $b \in \TILDE{w}$. Thus $w \SBS \TILDE{w}$, as claimed.

\qed

\subsection{Consequences of the e-simplicial property}

If a (for now hypothetical) subcomplex $\mathcal{S}'$ of $\mathcal{S}$ is an \NICOMP{n}{i} then all degenerate $(n+1)$-simplices of $\mathcal{S}'$ must be e-simplicial. This has consequences for lower-dimensional simplices of $\mathcal{S}'$. (Theorem \refpage{downwards-e-simplicial-thm} below).

In this section, we will continue from section \refpage{simplicialProp} to deduce the following:
\begin{itemize}
\item Whenever $n \ge 2$ and $y \in \mathcal{S}_n$ then $y$ is e-simplicial if and only if for all $j \in [n]$, $s_j(y)$ is e-simplicial.

\item If $w \in \mathcal{S}_{n+1}$ is an $(n,i)$-composition and if for all $j \neq i$, $d_j(w)$ is e-simplicial, then $d_i(w)$ is also e-simplicial.\MS
 
\end{itemize}
\BX

\begin{newlem}
{\bf (Component-simplicial lemma)}
\label{sptTechLemma}
\index{component-simplicial lemma}
\index{lemma!Component-simplicial lemma}

{\rm
\SL

Suppose $n \ge 2$, $y \in \mathcal{S}_n$ and $a \in y$.
\begin{enumerate}
\item For any $j \in [n]$, if $c_j(a)$ is component-simplicial then
for all $p<q \in [n] - \{j\}$, $e_p e_q (a) = e_{q-1} e_p (a)$.\MS

\item If $i,j,k \in [n]$ are distinct and $c_i(a), c_j(a),c_k(a)$ all
are component-simplicial, then $a$ is component-simplicial.
\end{enumerate}
}
\end{newlem}

\NI {\bf Proof:}\MS

There are three cases to consider for statement 1.\MS

\NI Case $p<q<j$:
\[
\begin{array}{l}
e_p e_q c_j (a) = e_p c_{j-1} e_q (a) = c_{j-2} e_p e_q (a)\\
e_{q-1} e_p c_j (a) = e_{q-1} c_{j-1} e_p (a) = c_{j-2} e_{q-1} e_p (a)
\end{array}
\]
Since $e_p e_q c_j(a) = e_{q-1} e_p c_j (a)$ by hypothesis and
$c_{j-2}$ is monic then $e_p e_q (a) = e_{q-1} e_p (a)$ in this case.
\MS

\NI Case $p<j<q$:
\[
\begin{array}{l}
e_p e_{q+1} c_j (a) = e_p c_j e_q (a) = c_{j-1} e_p e_{q} (a)\\
e_q e_p c_j (a) = e_q c_{j-1} e_p (a) = c_{j-1} e_{q-1}e_p (a)
\end{array}
\]
Since $e_p e_{q+1} c_j(a) = e_{q} e_p c_j (a)$ by hypothesis and
$c_{j-1}$ is monic then $e_p e_q (a) = e_{q-1} e_p (a)$ in this case.
\MS

\NI Case $j<p<q$:
\[
\begin{array}{l}
e_{p+1} e_{q+1} c_j(a) = e_{p+1} c_j e_q (a) = c_j e_p e_q (a)\\
e_q e_{p+1} c_j (a) = e_q c_j e_p (a) = c_j e_{q-1} e_p (a)
\end{array}
\]
Since $e_{p+1} e_{q+1} c_j(a) = e_{q} e_{p+1} c_j (a)$ by hypothesis and
$c_{j}$ is monic then $e_p e_q (a) = e_{q-1} e_p (a)$ in this case also.
\MS

For statement 2, we may suppose that $i<j<k$. Since $c_i (a)$ and $c_j(a)$
are component-simplicial, statement 1 implies that $e_p e_q (a) = e_{q-1} e_p (a)$
whenever $p<q \in [n]$ except possibly for $p=i$ and $q=j$. Statement 1 and
the assumption that $c_k(a)$ are component-simplicial imply that $e_i e_j(a) = e_{j-1} e_i (a)$.\\
\qed

\begin{newcor}
\label{arraymatchcor}
{\rm
\SL

Let $n \ge 2$, $y \in \mathcal{S}_n$ and $a \in y$. Then the following statements
are equivalent:
\begin{enumerate}
\item $a$ is component-simplicial
\item For all $j \in [n]$, $c_j(a)$ is component-simplicial.
\item There are three distinct indices $i,j,k \in [n]$ such that
$c_i(a), c_j(a), c_k(a)$ are component-simplicial.
\end{enumerate}
}
\end{newcor}

\NI {\bf Proof:}\MS

Clearly statement 2 implies statement 3, and the previous lemma
established that statement 3 implies statement 1. 

The proof that 
statement 1 implies statement 2 goes straight-forwardly by cases:\MS

Let $p<q \in [n+1]$.\MS

\NI {\bf Case} $p<q<j$:
\[
e_pe_q c_j(a) = e_p c_{j-1} e_q (a) = c_{j-2} e_p e_q(a) = c_{j-2} e_{q-1}e_p (a)
= 
\]
\[
=e_{q-1}c_{j-1} e_p (a) = e_{q-1} e_p c_j (a)
= e_{q-1} e_p c_j (a)
\]

\NI {\bf Cases} $p<q$ with $q=j$ or $q=j+1$:
\[
e_p e_q c_j (a) = e_p (a) = e_{q-1} c_{j-1} e_p (a) = 
e_{q-1} e_p c_j (a)
\]

\NI {\bf Case} $p<j<q-1$:
\[
e_p e_q c_j (a) = e_p c_j e_{q-1} (a) = c_{j-1} e_p e_{q-1} (a) =
c_{j-1} e_{q-2} e_p (a) = 
\]
\[
= e_{q-1} c_{j-1} e_p (a) =e_{q-1} e_p c_j (a)
\]

\NI {\bf Cases} $p<q$ and $p=j$ or $p=j+1$:
\[
e_p e_q c_j (a) = e_p c_j e_{q-1} (a) = e_{q-1} (a) = e_{q-1} e_p c_j (a)
\]

\NI {\bf Case} $j<p-1<p<q$:
\[
e_p e_q c_j (a) =  e_p c_j e_{q-1} (a) = c_j e_{p-1} e_{q-1} (a) = 
c_j e_{q-2} e_{p-1} (a) =
\]
\[
 = e_{q-1} c_j e_{p-1} (a) = e_{q-1} e_p c_j (a)
\]
\qed \MS

A key consequence:

\begin{newcor}
{\bf (Degeneracy lemma)}
\index{Degeneracy lemma}
\index{lemma!Degeneracy lemma}
\label{degenClosure}
{\rm 
\SL

Suppose $n \ge 2$. Then
$y \in \mathcal{S}_n$ is e-simplicial if and only if for each $i \in [n]$, $s_i(y)$ is e-simplicial.
}
\end{newcor}

\NI {\bf Proof:} $y$ is e-simplicial if and only if each $a \in y$ is component-simplicial
if and only if for each $i \in [n]$ $c_i(a)$ is component-simplicial if and only if for each
$i \in [n]$ $s_i(y)$ is e-simplicial.\\
\qed

\begin{newlem}
{\bf (e-simplicial composites)}
\label{simpclosure}
\index{e-simplicial composite lemma}
\index{lemma!e-simplicial composites}
{\rm 
\SL

Suppose $i \in [n+1]$, and $w \in \mathcal{S}_{n+1}$ is an $(n,i)$-composition with
$d_j(w)$ e-simplicial for all $j \in [n]-\{i\}$. Then $d_i(w)$ is also 
e-simplicial.
}
\end{newlem}

\MS

\NI {\bf Proof:} Suppose $a \in d_i(w)$. We must verify that $e_je_k(a)
= e_{k-1}e_j (a)$ for all $j<k \in [n]$. By definition of $(n,i)$-composition,
$a \in d_i(w)$ iff there exists $b \in w$ such that $e_i(b) = a$. Given such a $b \in w$, we then examine $e_je_k(a) = e_je_k(e_i(b))$ (by cases, somewhat
tediously).\MS

\NI {\bf Case} $j<k<i$:
\[
\begin{aligned}
e_je_k(a) &= e_je_k(e_i(b))\\
&= e_j e_{i-1} e_k (b)\qquad  \text{ (} w \text{ is e-simplicial)}\\
&= e_{i-2} e_j e_k(b)\qquad  \text{ (} d_kw \text{ is e-simplicial)}\\
&= e_{i-2} e_{k-1} e_j(b)\qquad  \text{ (} w \text{ is e-simplicial)}\\
&= e_{k-1} e_{i-1} e_j(b)\qquad  \text{ (} d_jw \text{ is e-simplicial)}\\
&= e_{k-1} e_j e_i(b)\qquad  \text{ (} w \text{ is e-simplicial)}\\
&= e_{k-1} e_j (a)
\end{aligned}
\]

\MS

\NI {\bf Case} $j<i \le k$:
\[
\begin{aligned}
e_je_k(a) &= e_je_k(e_i(b))\\
&= e_j e_{i} e_{k+1} (b)\qquad  \text{ (} w \text{ is e-simplicial)}\\
&= e_{i-1} e_j e_{k+1}(b)\qquad  \text{ (} d_{k+1}w \text{ is e-simplicial)}\\
&= e_{i-1} e_{k} e_j(b)\qquad  \text{ (} w \text{ is e-simplicial)}\\
&= e_{k-1} e_{i-1} e_j(b)\qquad  \text{ (} d_jw \text{ is e-simplicial)}\\
&= e_{k-1} e_j e_i(b)\qquad  \text{ (} w \text{ is e-simplicial)}\\
&= e_{k-1} e_j (a)
\end{aligned}
\]

\MS

\NI {\bf Case} $i \le j<k$:
\[
\begin{aligned}
e_je_k(a) &= e_je_k(e_i(b))\\
&= e_j e_{i} e_{k+1} (b)\qquad  \text{ (} w \text{ is e-simplicial)}\\
&= e_{i} e_{j+1} e_{k+1}(b)\qquad  \text{ (} d_{k+1}w \text{ is e-simplicial)}\\
&= e_{i} e_{k} e_{j+1}(b)\qquad  \text{ (} w \text{ is e-simplicial)}\\
&= e_{k-1} e_{i} e_{j+1}(b)\qquad  \text{ (} d_{j+1}w \text{ is e-simplicial)}\\
&= e_{k-1} e_j e_i(b)\qquad  \text{ (} w \text{ is e-simplicial)}\\
&= e_{k-1} e_j (a)
\end{aligned}
\]
\qed

\subsection{Wanted: an \NICOMP{n}{i} of sets}
\label{wantcompstruct} 

Starting with a fixed choice of $n>1$ and $i \in [n+1]$, we seek an {\em \NICOMP{n}{i} of sets}, by which we mean a subcomplex $\mathcal{S}'$ of $\mathcal{S}$ which satisfies condition \eqref{intro-composer-condition} (page \pageref{intro-composer-condition}). Namely, 
\[
\text{For all } m>n \quad \phi_{m,i}: \MC{S}'_m \to \Lambda^i(m)(\mathcal{S}') \text{ is an isomorphism}
\]
If $\MC{S}'$ satisfies condition \eqref{intro-composer-condition} then $\MC{S}'$ would also  be an \NICOMP{m}{i} for all $m \geq n+1$.

That is:
\begin{quote}
For each $m \geq n$, every $(m+1)$-simplex of $\MC{S}'$ is an $(m,i)$-composition.
\end{quote}
Each $(n+1)$-simplex of $\MC{S}'$ is an $(n,i)$-composition; each $(n+2)$-simplex is an $(n+1,i)$-composition, etc. \MS

If any such subcomplex $\MC{S}'$ exists then it has the properties:
\MS

\NI {\bf $\bullet$ Closure under degeneracies:}

\label{closuredegen}
For each $m \geq n$, every degenerate $(m+1)$-simplex of $\MC{S}'$ is an $(m,i)$-composition.
\MS

\NI {\bf $\bullet$ Closure under $(n,i)$-composition:}

\label{closurenicomp}
For all $m \geq n$, whenever $(\OM{y}{m+1}{i}) \in \Lambda^i(m+1)(\MC{S}')$ then the resulting $(m,i)$-composition
\[
w = \COMP{m}{i}(\OM{y}{m+1}{i}) 
\]
belongs to $\MC{S}'_{m+1}$ and $d_i(w) \in \MC{S}'_m$.\MS
\BS

$\mathcal{S}$ itself is {\em not} the structure we are looking for because not all degenerate $(n+1)$-simplices are $(n,i)$-compositions. Specifically, if $y \in \mathcal{S}_n$ and $y$ is not e-simplicial then (as a consequence of Corollary \refpage{degenClosure}) there exists at least one $j \in [n+1]$ such that $s_j(y)$ is not e-simplicial, hence not an $(n,i)$-composition.\MS

As a guide to finding $\MC{S}'$ (given $n$ and $i$) we ask that $\MC{S}'$ be {\em optimal} in the following senses.
\label{guidelines}
\begin{enumerate}
\item It is as large (inclusive) as possible. We should not impose conditions not strictly necessary.
\item It should be {\em reasonable} in the sense that the conditions we choose for the definition of $\MC{S}'$ should resemble what occurs in ordinary composition, the $(n,i)=(1,1)$ case. Also, the definition of $\MC{S}'$ should relate to $n$ and $i$ in as simple a way as possible. To state this negatively, we should not impose artificial conditions, conditions without motivation, or any conditions which involve $n$ and $i$ in a complicated or indirect way.
\end{enumerate}

Of course, ``reasonable'', ``resemble'', ``simple'', ``artificial'' etc. are judgments rather than precise criteria, and they are, therefore, guidelines rather than conditions.\MS

To illustrate the use of these guidelines, consider any small category of sets and functions, that is, $C$ is a \NICOMP{1}{1} where we require 1-simplices to be functions rather than general binary relations.
\label{function-characteristic}
Suppose $f:A \to B$ is any function of sets (1-simplex of $C$), and $\Gamma(f) \SBS B \times A$ denotes the graph of $f$. Two familiar properties characterize functions of sets and distinguish them from general binary relations.

\begin{enumerate}
\item[(i)] The composite $\Gamma(f) \hookrightarrow B \times A \stackrel{{\PRJ_A}}{\longrightarrow}A$ is surjective. That is: for every $a \in A$ there exists $b \in B$ such that $(b,a) \in \Gamma(f)$.

\item[(ii)] Given any $a \in A$ there exists a unique $(b,a) \in \Gamma(f)$. That is, the potential ordered pair $(-,a)$ uniquely determines $(b,a) \in \Gamma(f)$.
\end{enumerate}

The guideline of ``resemblance'' stated above will, at the least, support the reasonableness of any property we propose in the $n>1$ situation which bears some resemblance to properties (i) and (ii) for functions.

The same guideline would make us hesitate to propose a property completely unlike something which occurs in dimension 1.\MS

We will proceed as follows: first, we will deduce properties made necessary by closure under degeneracies and $(n,i)$-composition. Then we will complete the definition of $\MC{S}'$ {\em reasonably} in the sense just discussed.

\subsection{Surjectivity conditions and consequences}
\label{simplesurj}

By definition, if $w$ is an $(n,i)$-composition then $e_i^w : w \to d_i(w)$ is surjective. Other surjectivity conditions follow from this. This subsection focusses on a surjectivity property defined generally and develops some facts.

\begin{newdef}
{\bf ($i$-surjective)}
\index{$i$-surjective}
\index{surjectivity condition}
\label{surjCondDef}
{\rm
\SL

Given $n \geq 1$ and $y \in \mathcal{S}_n$ then $y$ will be said to be {\bf \em \SUR{i} } if $e_i^y : y \to d_i(y)$ is surjective.
}

\BX
\end{newdef}

In general, $n$-simplices of $\mathcal{S}$ are not \SUR{i}. However, $e_i : s_i (x) \to x$ and (if $i>0$) $e_i : s_{i-1} (x) \to x$ are in fact isomorphisms. Thus for any $x \in \mathcal{S}_n$ and $i>0$, $s_i(x)$ and $s_{i-1}(x)$ are \SUR{i}. $s_0(x)$ is \SUR{0} and \SUR{1}, and $s_n(x)$ is \SUR{n} and \SUR{(n+1)}.\MS

If a simplex is both \SUR{i} {and} e-simplicial, then some of its faces and simplicial images also have surjectivity conditions.
\MS

\begin{newlem}
{\bf (Basic Surjectivity)}
\label{basicSurjLemma}
\index{basic surjectivity lemma}
\index{lemma!basic surjectivity}
{\rm
\SL

\begin{enumerate}

\item Suppose $n \ge 2$, $p<q$ in $[n]$, $y \in \mathcal{S}_n$ and $y$ is e-simplicial. If $y$ is \SUR{q} and $d_q(y)$ is \SUR{p} then $d_p(y)$ is \SUR{{(q-1)}}.

\item Suppose $n \geq 1$, $y \in \mathcal{S}_n$. 
Then
\begin{enumerate}

\item If $i>1$, $s_0(y)$ is e-simplicial and \SUR{i} then $y$ is \SUR{(i-1)}.

\item If $i<n$, and $s_n(y)$ is e-simplicial and \SUR{i} then $y$ is \SUR{i}.
\end{enumerate}

\item Suppose $n \geq 2$, $i \in [n]$, $y \in \mathcal{S}_n$ and $y$ is \SUR{i}. Then:
\begin{enumerate}

\item If $j \leq i$ and $s_j(y)$ is e-simplicial then $s_j(y)$ is \SUR{(i+1)}.

\item If $j \geq i$ and $s_j(y)$ is e-simplicial, then $s_j(y)$ is \SUR{i}.

\end{enumerate}

\end{enumerate}

}
\end{newlem}

\NI {\bf Proof:}\MS

\NI (1) That $y$ is e-simplicial implies that the following square commutes:

\begin{center}
\setlength{\unitlength}{1in}
\begin{picture}(1,.6)
\put(0,.5){\MB{y}} 
\put(1,.5){\MB{d_qy}} 
\put(0,0){\MB{d_py}} 
\put(.8,0){{$d_pd_qy = d_{q-1}d_p y$}} 
\put(0.15,.5){\vector(1,0){0.7}} 
\put(0.5,.6){\MBS{e_q^y}} 
\put(0,0.35){\vector(0,-1){0.2}} 
\put(0.1,0.25){\MBS{e_p}} 
\put(1,0.35){\vector(0,-1){0.2}} 
\put(1.2,0.25){\MBS{e_p^{d_qy}}} 
\put(0.15,0){\vector(1,0){0.6}} 
\put(0.5,0.1){\MBS{e_{q-1}}} 
\end{picture}
\end{center}
\MS

\NI The surjectivity of $e_{q-1}$ follows immediately from the surjectivity of $e_q^y$
and $e_p^{d_qy}$.\MS

\NI (2) The conclusions follow from the commutativity of the following diagrams,
and that the vertical maps are isomorphisms.

\setlength{\unitlength}{1in}
\begin{picture}(4,.7)
\put(0,.5){\MB{s_0y}} 
\put(.85,.47){$d_is_0y = s_0 d_{i-1} y$} 
\put(0,0){\MB{y}} 
\put(1,0){\MB{d_{i-1}y}} 
\put(3,.5){\MB{s_ny}} 
\put(3.85,.47){$d_i s_n y = s_{n-1} d_i y$} 
\put(3,0){\MB{y}} 
\put(4,0){\MB{d_iy}} 
\put(0.15,.5){\vector(1,0){0.65}} 
\put(0.5,.6){\MBS{e_i}} 
\put(0,0.35){\vector(0,-1){0.2}} 
\put(0.1,0.25){\MBS{e_0}} 
\put(1,0.35){\vector(0,-1){0.2}} 
\put(1.1,0.25){\MBS{e_0}} 
\put(0.15,0){\vector(1,0){0.65}} 
\put(0.5,0.1){\MBS{e_{i-1}}} 
\put(3.15,.5){\vector(1,0){0.65}} 
\put(3.5,.6){\MBS{e_i}} 
\put(3,0.35){\vector(0,-1){0.2}} 
\put(3.1,0.25){\MBS{e_n}} 
\put(4,0.35){\vector(0,-1){0.2}} 
\put(4.2,0.25){\MBS{e_{n-1}}} 
\put(3.15,0){\vector(1,0){0.65}} 
\put(3.5,0.1){\MBS{e_i}} 
\put(1.6,.125){\MBS{(i>1)}}
\put(4.6,.125){\MBS{(i<n)}}
\end{picture}
\MS

\NI (3) If $j=i$ then then surjectivity claims follow immediately. In the cases when $j \neq i$, the hypotheses imply the following diagrams commute:

\setlength{\unitlength}{1in}
\begin{center}
\begin{picture}(4,.65)
\put(0,.5){\MB{s_jy}} 
\put(1,.5){\MB{d_{i+1}s_jy}} 
\put(0,0){\MB{y}} 
\put(1,0){\MB{d_iy}} 
\put(2,.5){\MB{s_jy}} 
\put(3,.5){\MB{d_is_jy}} 
\put(2,0){\MB{y}} 
\put(3,0){\MB{d_iy}} 
\put(0.15,.5){\vector(1,0){0.55}} 
\put(0.5,.6){\MBS{e_{i+1}}} 
\put(0,0.35){\vector(0,-1){0.2}} 
\put(.15,.25){\MBS{e_j^{s_jy}}}
\put(1,0.35){\vector(0,-1){0.2}} 
\put(1.25,0.25){\MBS{e_j^{d_{i+1}s_jy}}} 
\put(0.15,0){\vector(1,0){0.7}} 
\put(0.5,0.1){\MBS{e_i}} 
\put(2.15,.5){\vector(1,0){0.65}} 
\put(2.5,.6){\MBS{e_i}} 
\put(2,0.35){\vector(0,-1){0.2}} 
\put(2.15,0.25){\MBS{e_j^{s_jy}}} 
\put(3.1,0.35){\vector(0,-1){0.2}} 
\put(3.3,0.25){\MBS{e_{j-1}^{d_is_jy}}} 
\put(2.15,0){\vector(1,0){0.7}} 
\put(2.5,0.1){\MBS{e_i}} 
\put(-.25,.25){\footnotesize{iso}}
\put(.75,.25){\footnotesize{iso}}
\put(1.75,.25){\footnotesize{iso}}
\put(2.85,.25){\footnotesize{iso}}
\put(.5,-.2){\MBS{(j<i)}}
\put(2.5,-.2){\MBS{(j>i)}}
\end{picture}
\end{center}
\vskip.25in

If $j<i$ then $d_{i+1}s_j(y) = s_jd_i(y)$ and the conclusion follows from
$e_j^{s_jy}$ and $e_j^{d_{i+1}s_jy}$ being isomorphisms.

If $j>i$ then $d_is_j(y) = s_{j-1}d_i(y)$ and the conclusion follows from
$e_j^{s_jy}$ and $e_{j-1}^{d_is_jy}$ being isomorphisms.
\\

\qed

Surjectivity conditions propagate {\em downwards} in dimension under certain circumstances.

\begin{newcor}
{\bf (Downward Surjectivity)}
\label{downward-surj-corollary}
\index{Downward Surjectivity Corollary}
{\rm
\SL

Suppose $m \geq 2$, $i \in [m]$ and $\MC{S}'$ is a subcomplex of $\MC{S}$ such that every $y \in \MC{S}'_m$ is \SUR{i} and e-simplicial. Let $x \in \MC{S}'_{m-1}$.
\MS

1. If $i>1$ then $x$ is \SUR{(i-1)}.
\MS

2. If $i<m-1$ then $x$ is \SUR{i}.
}
\end{newcor}

\Proof

\NI 1. $s_0(x) \in \MC{S}'_m$ is \SUR{i} by hypothesis and the following diagram commutes, using that $d_i s_0(x) = s_0 d_{i-1}(x)$ when $i>1$:
\begin{center}
\begin{picture}(1,.6)
\put(0,.5){\MB{s_0(x)}}
\put(1,.5){\MB{s_0d_{i-1}(x)}}
\put(0,0){\MB{x}}
\put(1,0){\MB{d_{i-1}(x)}}
\put(0,.4){\vector(0,-1){.3}} 
\put(-.15,.25){\MBS{e_0}}
\put(1,.4){\vector(0,-1){.3}} 
\put(1.15,.25){\MBS{e_0}}
\put(.25,.5){\vector(1,0){.4}} 
\put(.4,.6){\MBS{e_i}}
\put(.25,0){\vector(1,0){.5}} 
\put(.4,-.1){\MBS{e_{i-1}}}
\end{picture}
\end{center}
Since $e_0 : s_0(x) \to x$ and $e_0 : s_0d_{i-1}(x) \to d_{i-1}(x)$ are bijections, and since $e_i: s_0(x) \to s_0 d_{i-1}(x)$ is surjective by hypothesis, then $e_{i-1} : x \to d_{i-1}(x)$ must also be surjective. 
\MS

\NI 2. $s_{m-1}(x) \in \MC{S'}_m$ is \SUR{i} and the following diagram commutes, using that $i<m-1$ implies $d_i s_{m-1}(x) = s_{m-2} d_{i}(x)$:
\begin{center}
\begin{picture}(1,.6)
\put(0,.5){\MB{s_{m-1}(x)}}
\put(1,.5){\MB{s_{m-2}d_{i}(x)}}
\put(0,0){\MB{x}}
\put(1,0){\MB{d_{i}(x)}}
\put(0,.4){\vector(0,-1){.3}} 
\put(-.15,.25){\MBS{e_{m-1}}}
\put(1,.4){\vector(0,-1){.3}} 
\put(1.2,.25){\MBS{e_{m-2}}}
\put(.25,.5){\vector(1,0){.4}} 
\put(.45,.6){\MBS{e_i}}
\put(.25,0){\vector(1,0){.5}} 
\put(.5,-.1){\MBS{e_{i}}}
\end{picture}
\end{center}
Since $e_{m-1} : s_{m-1}(x) \to x$ and $e_{m-2} : s_{m-2}d_{i}(x) \to d_i(x)$ are bijections, and since $e_i: s_{m-1}(x) \to s_{m-2} d_{i}(x)$ is surjective by hypothesis, then $e_{i} : x \to d_{i}(x)$ must also be surjective.

\qed
\MS

See section \refpage{furthermore} for a note on downward propagation of surjectivity conditions.

\section{The subface-simplicial property}
\label{subfacesimpprop}

We showed (corollary \refpage{degenClosure}) that $y \in \mathcal{S}_n$ is e-simplicial if and only if $s_j(y)$ is e-simplicial for all $j \in [n]$. By repeated use of this fact, we get:

\begin{newthm}
{\bf (e-simplicial downwards)}
\label{needSubfaceSimplicial} \label{downwards-e-simplicial-thm}
\index{e-simplicial downwards}
\index{Theorem!e-simplicial downwards}
{\rm 
\SL

Suppose $n \ge 1$ and $i \in [n+1]$. Let $\MC{S}'$ be any subcomplex of $\mathcal{S}$ such that all degenerate $(n+1)$-simplices of $\MC{S}'$ are e-simplicial. Then for all $k=2, \cdots,n$, the $k$-simplices of $\MC{S}'$ are e-simplicial.
}
\end{newthm}\MS

\NI {\bf Proof:} Since all degeneracies in $\mathcal{S}'_{n+1}$ are assumed e-simplicial, then all $n$-simplices of $\MC{S}'$ must be e-simplicial
by Corollary \refpage{degenClosure}. In particular, all the degenerate
$n$-simplices of $\MC{S}'$ must be e-simplicial, in which case, again by Corollary \refpage{degenClosure}, all $(n-1)$-simplices of $\MC{S}'$ must also be e-simplicial.
Thus, the hypothesis implies the ``e-simplicial'' property propagates
all the way down to dimension 2 (and to dimensions 1 and 0 vacuously).\\
\qed \MS

\begin{newdef}
{\bf (Subface-simplicial)}
\index{subface-simplicial}
\label{subfaceSimplicialDef}
{\rm
\SL

Let $k \ge 2$. The simplex $y \in \mathcal{S}_k$ will be called {\bf \em
subface-simplicial} if $y$ and every subface of $y$ is e-simplicial.\\
\NI \BX}
\end{newdef}

\begin{newlem}
\label{subSimpIsaComplex}
{\rm
\SL

The subface-simplicial simplices of $\mathcal{S}$ form a subcomplex of $\mathcal{S}$.
}
\end{newlem} \MS

\NI {\bf Proof:} It's immediate that they form a face complex. If $y \in
\mathcal{S}_1$ then it's vacuously e-simplicial and $s_0(y), s_1(y)$ are e-simplicial, hence also subface-simplicial. If $y \in \mathcal{S}_n$, is subface-simplicial
and $n>1$ then $s_j(y)$ is e-simplicial for all $j \in [n]$ and each face
$d_is_j(y)$ is either $y$ or the degenerate image of a face of $y$. Hence, by induction on dimension, each face of $s_j(y)$ is subface-simplicial. \\
\qed
\MS

In the terminology of this definition, the previous results imply:

\begin{newcor}
{\rm
\SL

Any subcomplex of $\mathcal{S}$ which satisfies closure-under-degeneracies and closure-under--$(n,i)$-composition (page \pageref{closuredegen}) has the property that all its simplices in dimensions $n+1$ and
below are subface-simplicial.
}

\qed
\end{newcor}

\subsection{Properties of subface-simplicial simplices}
\label{subface-simplicial-properties}

A \SBSI\ $n$-simplex has special properties which will be examined in this section. In particular, a \SBSI\ $n$-simplex is much simpler than a general $n$-simplex. 
\BS

Suppose $y \in \mathcal{S}_n$ is \SBSI. It follows that whenever
\[
(i_n, i_{n-1},\cdots, i_k), (j_n,j_{n-1}, \cdots, j_k) \in 
[n] \times [n-1] \times \cdots \times [k]
\]
and $d_{i_k} d_{i_{k+1}} \cdots d_{i_n} = d_{j_k} d_{j_{k+1}} \cdots
d_{j_n}$ is a simplicial face-identity, then for each $a \in y$
\[
e_{i_k} e_{i_{k+1}} \cdots  e_{i_n}(a) = 
e_{j_k} e_{j_{k+1}} \cdots e_{j_n}(a)
\]
Consider, in particular, the fundamental entries (see
definition \refpage{fundEntries}) of $a \in y$. Each has the form
\[
e_{i_n} e_{i_{n-1}} \cdots e_{i_1}(a) = 
e_0 \cdots \omit{e_j} \cdots e_n (a)
\]
for exactly one $j \in [n]$. Hence $a \in y$ has the special property
that its fundamental matrix, with its $(n+1)!$ entries, is formed from exactly $n+1$ fundamental
entries, namely $e_0 \cdots \omit{e_j} \cdots e_n (a)$ for $j=0, \cdots,
n$. We name this property (of an element of $y$) in the next definition.\MS

For the record:

\begin{newcor}
{\rm
\SL

Suppose $n \geq 1$, $y \in \mathcal{S}_n$ is subface-simplicial and that $a \in y$. Then each entry in the fundamental matrix of $a$ has, for some $j \in [n]$, the form  the form $e_0 \cdots \omit{e_j} \cdots e_n (a)$.
}
\qed
\end{newcor}

\NI The following definition sets up some convenient terminology:

\begin{newdef}
{\bf ($\FE(V_0 \cdots V_n)$)}
\index{subcomponent-simplicial}
\label{subcomp-simp}
\index{Fund$(V_0 \DDD{} V_n)$}
{\rm 
\SL 

Given any non-empty sets $\SEQ{V}{n}$ with $n \geq 2$

\begin{enumerate}

\item
We say that $a \in \FSIG(\SEQ{V}{n})$ is {\bf \em subcomponent-simplicial} if for any $k<n$ and whenever
\[
(i_n, i_{n-1},\cdots, i_k), (j_n,j_{n-1}, \cdots, j_k) \in 
[n] \times [n-1] \times \cdots \times [k]
\]
such that $d_{i_k} d_{i_{k+1}} \cdots d_{i_n} = d_{j_k} d_{j_{k+1}} \cdots
d_{j_n}$ is a simplicial face-identity, then
\[
e_{i_k} e_{i_{k+1}} \cdots  e_{i_n}(a) = 
e_{j_k} e_{j_{k+1}} \cdots e_{j_n}(a)
\]

\item 
Notations:
\index{$\FE(V_0 \cdots V_n)$}
\label{fund-elements}
\[
\begin{aligned}
\FE(\SEQ{V}{n}) & \DFAS
\SETT{a \in \FSIG(\SEQ{V}{n}):a \text{ is subcomponent-simplicial}}\\
\FE(y) & \DFAS  
\SET
{a \in y}
{
a \text{ is subcomponent-simplicial} 
}
\end{aligned}
\]

\end{enumerate}
\NI \BX}
\end{newdef}
\MS
\NI Observe that for all $n>0$, 
\[
\FE \LIST{V}{n} \SBS\prod_{i=0}^n \FE (\DLI{V}{n}{i})
\]
\MS

\begin{example}
Suppose $n=2$ and 
\[
a \in \FSIG(V_0,V_1,V_2) = V_2 V_1 V_2 V_0 V_1 V_0
\] 
If $a \in \FE(V_0,V_1,V_2)$ then its fundamental matrix is
\[
a = 
\left[
\begin{array}{cc}
a_{00} & a_{01}\\
a_{10} & a_{11}\\
a_{20} & a_{21}
\end{array}
\right]
=
\left[
\begin{array}{cc}
t_2 & t_1\\
t_2 & t_0\\
t_1 & t_0
\end{array}
\right]
\]
where the fundamental entries of $a$ are $a_{10}= e_0 e_1(a) = t_2 \in V_2, a_{20} = e_0 e_2(a) = t_1\in V_1$ and  $a_{21}  = e_1 e_2(a)=t_0 \in V_0$.
\end{example}
\MS

\NI {\bf A maximal \SBSI\ $n$-simplex:}
\label{maxsubsimp}
\index{maximal subface-simplicial $m$-simplex}
Given $\SETT{\SEQ{V}{m}} \SBS \MC{V}$, there is a corresponding maximal \SBSI\ $m$-simplex denoted $\MX(\SEQ{V}{m})$ and defined as follows.
\MS

\NI If $m=0$ then $\MX(V_0) \DFAS (V_0 \SBS V_0)$.

\NI If $m=1$ then $\MX(V_0,V_1) \DFAS (V_1 V_0 \SBS V_1 V_0)$.

\NI If $m>1$ then $\MX(\SEQ{V}{m})$ is the $m$-simplex
\[
\MX(\SEQ{V}{m}) \SBS \prod_{p=0}^m \MX(\OM{V}{m}{p})
\]
that is, $d_p(\MX(\SEQ{V}{m})) = \MX(\OM{V}{m}{p})$, 
where \\
$t \in \MX(\SEQ{V}{m})$ iff $t \in \FE(\SEQ{V}{m})$ i.e. $t$ is subcomponent-simplicial.

Note that $\FSIG(\MX(\SEQ{V}{m})) = \FSIG(\SEQ{V}{m})$. 
\MS

\NI It follows from the definitions that if $(y \SBS \Sig(y)) \in \mathcal{S}_m$ is \SBSI\ with fundamental signature $\FSIG(\SEQ{V}{m})$ then $y \SBS \MX(\SEQ{V}{m})$.

\BX

\subsection{Representation of subcomponent-simplicial elements}
\label{subcomp-simp-elem-rep}

Given an arbitrary $y \in \mathcal{S}_n$, the representation of $a \in y$ as an element of the product
$\FSIG(\SEQ{V}{n})$ is the higher-dimensional analog of listing the
entries of a rectangular array in row order.

Now if $y \in \mathcal{S}_n$ and $a \in y$ is subcomponent-simplicial then the row-order list of its fundamental entries $t_j = e_0 \cdots \omit{e_j} \cdots e_n(a)$ has a particularly simple description which we give inductively, as follows. \MS

\begin{newdef}
\index{$h_n$}
\index{h$_n$}
\label{hFunc}
{\rm
\SL

Let $n \ge 0$ and suppose $\SEQ{V}{n} \in \MC{V}$ are non-empty. For each $j \in [n]$, we will denote by $\text{pr}^{j}$ the projection function
\[
\begin{array}{l}
\prod_{i=0}^n V_i \to \prod_{i \neq j} V_i\\
\\
\LIST{t}{n} \mapsto
(\OM{t}{n}{j})
\end{array}
\]

Define the function
\[
h_n : \prod_{i=0}^n V_i \to \FSIG(\SEQ{V}{n})
\]
inductively (on $n$) as follows:

{\bf Case $n=0$:} $h_0 \DFAS 1_{V_0} : V_0 \to V_0$.\MS

{\bf Case $n>0$:} $h_n = ( h_{n-1} \circ \text{pr}^{0},\cdots , h_{n-1}\circ \text{pr}^{n} )$. That is, if $t = \LIST{t}{n} \in \prod_{i=0}^n V_i$, and $t^{j} \DFAS ( \OM{t}{n}{j} )=\text{pr}^{j}(t)$ then
\[
h_n(t) = \left( h_{n-1}(t^{0}), \cdots , h_{n-1}(t^{n}) \right)
\]

\NI \BX}
\end{newdef}

Note: This is related closely to the definition of $\FSIG$ back on page \pageref{Fsig-definition}.
\BS

\begin{example}
\[
\begin{aligned}
h_1(t_0, t_1) &= (h_0(t_1), h_0(t_0)) = (t_1, \hspace{.2cm} t_0)\\
h_2(t_0,t_1,t_2) &= (h_1(t_1,t_2),h_1(t_0,t_2),h_1(t_0,t_1))
= (t_2,t_1,\hspace{.2cm} t_2,t_0,\hspace{.2cm} t_1,t_0)
\end{aligned}
\]
\begin{multline*}
h_3(t_0,t_1,t_2,t_3) = \Bigl(h_2(t_1,t_2,t_3), h_2(t_0,t_2,t_3),h_2(t_0,t_1,t_3),h_2(t_0,t_1,t_2) \Bigr)\\
= (t_3,t_2,t_3,t_1,t_2,t_1,\hspace{.2cm} t_3,t_2,t_3,t_0,t_2,t_0,\hspace{.2cm} t_3,t_1,t_3,t_0,t_1,t_0,\hspace{.2cm} t_2,t_1,t_2,t_0,t_1,t_0)
\end{multline*}
\end{example}
\MS

It follows from the definition that if $a = h_n(\SEQ{t}{n})$ then
\[
\begin{aligned}
e_j(a) &= h_{n-1}(\OM{t}{n}{j})\\
c_j(a) &= h_{n+1}(t_0, \cdots , t_{j-1},t_j,t_j,t_{j+1}, \cdots ,t_n)
\end{aligned}
\]

\NI Following immediately from the definitions:\MS

\begin{newlem}
\label{hnIsBijec}
{\rm
\SL

For each $n \ge 0$, $h_n :\PROD{i=0}{n}{V_i} \to \FE(\SEQ{V}{n})$ is
bijective.
}
\end{newlem}
\qed

The main importance of these functions $h_n$, other than being descriptive, is simply that they are bijective.

\begin{newlem}
{\bf (Uniqueness)}
\label{uniquefundelem} 
\index{Uniqueness lemma}
\index{lemma!Uniqueness lemma}
{\rm
\SL 

Suppose $m \geq 2$, $p,q \in [m]$, $p<q$, $\SEQ{V}{m}$ are non-empty and
\[
\begin{aligned}
(t_0, \cdots , t_p, \cdots, \LVO{q}, \cdots , t_m) & \in  V_0 \cdots \widehat{V_q} \cdots V_m \\
(t_0, \cdots , \LVO{p}, \cdots, t_q, \cdots , t_m) & \in  V_0 \cdots \widehat{V_p} \cdots V_m
\end{aligned}
\] 

Let 
\[
\begin{aligned}
x &= h_{m-1}(t_0, \cdots , t_p, \cdots, \LVO{q}, \cdots , t_m)\\
y &= h_{m-1}(t_0, \cdots , \LVO{p}, \cdots, t_q, \cdots , t_m)\\
z &= h_m(t_0, \cdots , t_p, \cdots, t_q, \cdots , t_m)
\end{aligned}
\]

Then $x = e_q(z), y=e_p(z)$ and $z$ is unique with those properties.
}
\end{newlem}

\NI {\bf Proof:} 

That $x = e_q(z)$ and $y = e_p(z)$ follows directly from the definition of $h_m$.

If $e_q(z') = x = e_q(z)$ then
\[
h_{m-1}(t'_0, \cdots , t'_p, \cdots , \LVO{q}, \cdots , t'_m) =
h_{m-1}(t_0, \cdots , t_p, \cdots , \LVO{q}, \cdots , t_m)
\]
from which it follows that $t_i = t'_i$ for all $i \neq q$ because $h_{m-1}$ is monic. By the same reasoning, $e_q(z') = y = e_q(z)$ implies $t_i=t'_i$ for all $i \neq p$. The claim that $z=z'$ follows. 

\qed
\MS

Here is a low-dimensional example to illustrate.
\MS

\begin{example}
Suppose $y \in \mathcal{S}_3$ is subface-simplicial, $\FSIG(y) = \FSIG(V_0, V_1,V_2,V_3)$ and that $a \in y$ with fundamental entries $(t_0, t_1,t_2,t_3) \in V_0 V_1 V_2 V_3$. Then $a$ expands as follows:
\[
\begin{aligned}
a &= (a_0, a_1, a_2, a_3)\\
&= \Bigl( h_2(t_1,t_2,t_3),h_2(t_0,t_2,t_3),h_2(t_0,t_1,t_3),h_2(t_0,t_1,t_2) \Bigr)\\
\end{aligned}
\]

where these $e_j(a)$ form the rows of the following matrix for $a$:
\[
\left[
\begin{array}{l}
a_0\\
a_1\\
a_2\\
a_3\\
\end{array}
\right]
=
\left[
\begin{array}{r}
h_2(t_1,t_2,t_3)\\
h_2(t_0,t_2,t_3)\\
h_2(t_0,t_1,t_3)\\
h_2(t_0,t_1,t_2)\\
\end{array}
\right]
=
\left[
\begin{array}{ccc}
(t_3,t_2) & (t_3,t_1) & (t_2,t_1)\\
(t_3,t_2) & (t_3,t_0) & (t_2,t_0)\\
(t_3,t_1) & (t_3,t_0) & (t_1,t_0)\\
(t_2,t_1) & (t_2,t_0) & (t_1,t_0)\\
\end{array}
\right]
\]
Observe that {\em any two rows} of this matrix involve all the $t_j$. 
\end{example}
\MS

\begin{newcor} 
{\bf (2-lemma)}
\index{two-row lemma} 
\index{two-element uniqueness}
\index{2-lemma}
\index{lemma!2-lemma}
\label{two-element-uniqueness}
{\rm 
\SL

Assume $m \geq 2$, $y \in \mathcal{S}_m$ and $y$ is subface-simplicial. Then
\begin{enumerate}
\item {\bf (``Two-row lemma'')} Given $p<q$ in $[m]$, $a_p \in d_p(y)$, $a_q \in d_q(y)$ such that $e_p(a_q) = e_{q-1}(a_p)$  

-- that is, given a component-simplicial partial element $\SETT{a_p, a_q}$ of $y$ --

\NI then there exists a unique $b \in \FSIG(y)$ such that $e_p(b)=a_p$ and $e_q(b)=a_q$.

\item {\bf (``Two-element uniqueness'')} Suppose $a,a' \in y$ and for some $p \neq q$ one has $e_p(a) = e_p(a')$ and $e_q(a) = e_q(a')$. Then $a=a'$.
\end{enumerate}
}
\end{newcor}
\MS

Note that part (1) of the corollary can be interpreted as saying that given two rows, $a_p$ and $a_q$, of an array of a potential element of $y$, then there is a unique element of $\FSIG(y)$ which fills in that array.
\MS

\NI {\bf Proof:}

{\bf Part (1):} Suppose $\FSIG(y) = \FSIG(\SEQ{V}{m})$. We are given $p<q$. By the surjectivity of the ``$h$'' functions we know that
\[
\begin{aligned}
a_p &= h_{m-1}\left( t_0, \cdots , t_{p-1},t_{p+1}, \cdots ,t_q , \cdots ,t_m \right) \in V_0 \cdots V_{p-1} V_{p+1} \cdots V_m\\
a_q &= h_{m-1}\left( t'_0, \cdots , t'_p, \cdots, t'_{q-1},t'_{q+1}, \cdots, t'_m \right) \in V_0 \cdots V_{q-1}V_{q+1} \cdots V_m
\end{aligned}
\]
for some 
\[(t_0, \cdots , \LVO{p}, \cdots ,t_q , \cdots ,t_m) \in V_0 \cdots \omit{V_p} \cdots V_m
\]
and some
\[(t'_0, \cdots , t'_p, \cdots, \LVO{q}, \cdots, t'_m) \in V_0 \cdots \omit{V_q} \cdots V_m
\]
(Make the obvious adjustments if $p=0$ or $q=m$ or $q-1=p$). 

The condition $e_p(a_q) = e_{q-1}(a_p)$ translates to
\[
\begin{aligned}
e_p(a_q) &= h_{m-2} \left( t'_0 , \cdots , t'_{p-1},t'_{p+1} , \cdots , t'_{q-1},t'_{q+1}, \cdots , t'_m  \right) = \\
e_{q-1}(a_p) &= 
h_{m-2}
\left(
t_0, \cdots , t_{p-1},t_{p+1} , \cdots , t_{q-1},t_{q+1}, \cdots , t_m
\right)
\end{aligned}
\]

Then: the monotonicity of $h_{m-2}$ implies that $t_k = t'_k$ for all $k \neq p,q$.

Set $b = h_{m-1} \left( t_0, \cdots , t'_p, \cdots , t_q, \cdots , t_m \right) \in \FSIG(y)$. Then $b$ is unique  by the  uniqueness lemma (page \pageref{uniquefundelem}) with the property that $e_p(b) = a_p$ and $e_q(b) = a_q$.
\MS

{\bf Part (2):} Part (1) applies to $e_p(a)$ and $e_q(a)$ since $e_p e_q(a) = e_{q-1}e_p(a)$. That is, $a \in \FSIG(y)$ is uniquely determined by $e_p(a)$ and $e_q(a)$. Therefore, $e_p(a)=e_p(a')$ and $e_q(a) = e_q(a')$ can both occur only when $a=a'$. 

\qed

\subsection{The vertex-relation}
\label{vertex-relation}

A subface-simplicial $y \in \mathcal{S}_n$ determines an ordinary $(n+1)$-ary relation.

\begin{newdef}
{\bf (Vertex relation)}
\index{vertex relation of a simplex}
\label{vertRel}
{\rm
\SL

\NI Suppose $n \ge 1$, $y \in \mathcal{S}_n$ with fundamental signature $\FSIG(V_0,\cdots,V_n)$ and $y$ is subface-simplicial. Then the relation $R(y) \SBS V_0 \cdots V_n$ defined as the pullback below will be called the {\bf vertex-relation} of $y$.

\setlength{\unitlength}{1in}
\begin{center}
\begin{picture}(2,.6)
\put(0,.5){\MB{R(y)}} 
\put(2,.5){\MB{V_0 \cdots V_n}} 
\put(0,0){\MB{y}} 
\put(2,0){\MB{\FE(V_0, \cdots, V_n)}} 
\put(0.2,.5){\vector(1,0){1.45}} 
\put(0,0.35){\vector(0,-1){0.2}} 
\put(2,0.35){\vector(0,-1){0.2}} 
\put(0.2,0){\vector(1,0){1.1}} 
\put(2.15,.25){\MBS{h_n}}
\end{picture}
\end{center}

\NI \BX}
\end{newdef}
\NI That is, $\LIST{t}{n} \in R(y)$ iff $h_n \LIST{t}{n}\in y$.
\MS

A non-empty relation $R \SBS \PLIST{V}{n}$ corresponds to
\[
T(R) \DFAS \SETT{h_n \LIST{t}{n}: \LIST{t}{n} \in R}
\SBS \FSIG(\SEQ{V}{n})
\]
The definitions of $T^\text{min}$ and $T^\text{max}$ given in section \refpage{poset-of-n-simps} applied to $T(R)$ yield \SBSI\ $n$-simplices $T(R)^\text{min}$ and $T(R)^\text{max}$ whose vertex-relation is $R$. The simplex $T(R)^\text{min}$ will be denoted $y^R$ in section \refpage{model-construction} below. 

Summarizing:

\begin{newthm}
{\rm
\SL

To every $(n+1)$-ary relation $R \SBS V_0 \cdots V_n$ there
corresponds a subface-simplicial $y \in \mathcal{S}_n$ such that $R=R(y)$.
}
\end{newthm}
\qed

\subsection{Subface-simplicial $(n,i)$-compositions}

Our main goal (see section \refpage{wantcompstruct}) was: when given $n>1$ and $0 \leq i \leq n+1$, to identify subcomplexes $\MC{S}'$ of $\mathcal{S}$ in which all $(n+1)$-simplices are $(n,i)$-compositions. By Theorem \refpage{downwards-e-simplicial-thm} we know that every $(n,i)$-composition in such a subcomplex $\MC{S}'$ must be \SBSI. Here are two more deduced properties.

\begin{newthm}
{\bf (Necessity of surjectivity)}
\label{nicompSubfaceSimp} 
\label{necess-of-surj-thm}
\index{necessity of surjectivity theorem}
\index{Theorem!Necessity of surjectivity}
{\rm 
\SL

Fix $n \geq 2$, $i \in [n+1]$ and assume $\MC{S}'$ is a subcomplex of $\mathcal{S}$ with the property that each $(n+1)$-simplex of $\MC{S}'$ is an $(n,i)$-composition. Then
\begin{enumerate}

\item Every $k$-simplex of $\MC{S}'$, $k \geq 2$, is subface-simplicial.

\item If $i=0$ or $i=1$ then every $n$-simplex of $\MC{S}'$ is \SUR{i}.

\item If $1<i<n$, every $n$-simplex of $\MC{S}'$ is both \SUR{i} and \SUR{(i-1)}.

\item If $i=n$ or $i=n+1$ then every $n$-simplex of $\MC{S}'$ is \SUR{(i-1)}.

\end{enumerate}
}
\end{newthm}

\NI {\bf Proof:}

{Claim (1):} By hypothesis, every $(n+1)$-simplex of $\MC{S}'$ is e-simplicial. By Theorem \refpage{downwards-e-simplicial-thm}, it follows that every $k$-simplex of $\MC{S}'$, $k \geq 2$ is subface-simplicial.
\MS

The remaining claims follow from the Downward Surjectivity Corollary \refpage{downward-surj-corollary} and that every $(n,i)$-composition is \SUR{i}.

\qed
\MS

The next theorem is an essential step in characterizing \NICOMP{n}{i}s of sets.

\begin{newthm}
{\bf (Degenerate $(n,i)$-compositions) }
\label{condsNicomp} 
\label{degen-comps-thm} 
\index{degenerate compositions theorem}
\index{Theorem!degenerate compositions}
{\rm
\SL

Fix $n \geq 2$ and $i \in [n+1]$. Let $y \in \mathcal{S}_n$. Assume:
\begin{itemize}

\item 
$y$ is \SBSI.

\item
$y$ is \SUR{i} if $0 \leq i \leq n-1$.

\item $y$ is \SUR{(i-1)} if $2 \leq i \leq n+1$.
\end{itemize}
Then for all $k \in [n]$, $s_k(y)$ is an $(n,i)$-composition.
}
\end{newthm}

\Proof

The first step is to apply the surjectivity assumptions for $y$ to deduce that $s_k(y)$ is \SUR{i}. Then the Expansion lemma \refpage{sub-composition} applies.
\MS

If $i=0$: $y$ is assumed \SUR{0}. For $k=0$, $s_0(y) \xrightarrow{e_0} y$ is a bijection. For $k >0$ the Basic Surjectivity lemma \refpage{basicSurjLemma} implies that $s_k(y) \xrightarrow{e_0} y$ is surjective.
\MS

If $i=1$: $y$ is assumed \SUR{1}. For $k=0,1$, $s_k(y) \xrightarrow{e_1} y$ is a bijection. For $k>1$, the Basic Surjectivity lemma implies $s_k(y) \xrightarrow{e_0} y$ is surjective.
\MS

If $2 \leq i \leq n-1$: $y$ assumed \SUR{i} and \SUR{(i-1)}. For $k=i-1$ and $i$, $s_k(y) \xrightarrow{e_i} y$ is a bijection. The Basic Surjectivity lemma applies in the other cases. The diagrams are:

\setlength{\unitlength}{1in}
\begin{center}
\begin{picture}(3,1.2)
\put(.5,.15){\MB{0 \leq k \leq i-2}}
\put(0,.5){\MB{y}}
\put(0,1){\MB{s_k(y)}}
\put(1,.5){\MB{d_{i-1}(y)}}
\put(1,1){\MB{s_kd_{i-1}(y)}}
\put(.2,.5){\vector(1,0){.4}}
\put(.2,1){\vector(1,0){.4}}
\put(0,.85){\vector(0,-1){.25}}
\put(1,.85){\vector(0,-1){.25}}
\put(.4,1.1){\MBS{e_i}}
\put(.4,.4){\MBS{e_{i-1}}}
\put(.1,.75){\MBS{e_k}}
\put(1.1,.75){\MBS{e_k}}

\put(2.5,.15){\MB{i+2 \leq k \leq n}}
\put(2,.5){\MB{y}}
\put(2,1){\MB{s_k(y)}}
\put(3,.5){\MB{d_i(y)}}
\put(3,1){\MB{s_{k-1} d_i(y)}}
\put(2.2,.5){\vector(1,0){.4}}
\put(2.2,1){\vector(1,0){.4}}
\put(2,.85){\vector(0,-1){.25}}
\put(3,.85){\vector(0,-1){.25}}
\put(2.4,1.1){\MBS{e_i}}
\put(2.4,.4){\MBS{e_i}}
\put(2.1,.75){\MBS{e_k}}
\put(3.2,.75){\MBS{e_{k-1}}}
\end{picture}
\end{center}

\NI where the vertical maps are bijections and the lower horizontal maps are assumed surjective. It follows that the upper horizontal maps are also surjective, as claimed.
\MS

If $i=n$: $y$ is assumed \SUR{(n-1)}. Then for $k=n$ or $n-1$, $s_k(y) \xrightarrow{e_n} y$ is bijective. For $0 \leq k \leq n-2$ the relevant diagram is
\begin{center}
\begin{picture}(1,.7)
\put(0,0){\MB{y}}
\put(0,.5){\MB{s_k(y)}}
\put(1,0){\MB{d_{n-1}(y)}}
\put(1,.5){\MB{s_k d_{n-1}(y)}}
\put(.2,0){\vector(1,0){.4}}
\put(.2,.5){\vector(1,0){.4}}
\put(0,.35){\vector(0,-1){.25}}
\put(1,.35){\vector(0,-1){.25}}
\put(.4,.6){\MBS{e_n}}
\put(.4,-.1){\MBS{e_{n-1}}}
\put(.1,.25){\MBS{e_k}}
\put(1.1,.25){\MBS{e_k}}
\end{picture}
\end{center}
\vskip.1in

\NI where the vertical maps are bijections and the lower horizontal map is assumed surjective. It follows that the upper horizontal map is also surjective.
\MS

If $i=n+1$: $y$ is assumed \SUR{n}. When $k=n$ then $s_k(y) \xrightarrow{e_{n+1}} y$ is bijective. For $0 \leq k \leq n-1$ the relevant diagram is
\begin{center}
\begin{picture}(1,.7)
\put(0,0){\MB{y}}
\put(0,.5){\MB{s_k(y)}}
\put(1,0){\MB{d_{n}(y)}}
\put(1,.5){\MB{s_k d_{n}(y)}}
\put(.2,0){\vector(1,0){.4}}
\put(.2,.5){\vector(1,0){.4}}
\put(0,.35){\vector(0,-1){.25}}
\put(1,.35){\vector(0,-1){.25}}
\put(.4,.6){\MBS{e_{n+1}}}
\put(.4,-.1){\MBS{e_{n}}}
\put(.1,.25){\MBS{e_k}}
\put(1.1,.25){\MBS{e_k}}
\end{picture}
\end{center}

\NI where (again) the vertical maps are bijections and the lower horizontal map is assumed surjective. It follows that the upper horizontal map is also surjective.
\BS

At this stage of the proof we have verified that for all $k \in [n]$, $s_k(y)$ is \SUR{i}. Therefore the Expansion lemma is applicable to the situation where $w = \COMP{n}{i}(d_0 s_k(y), \cdots , \widehat{d_i s_k(y)}, \cdots , d_{n+1} s_k(y))$.

From the Expansion lemma it follows that $d_i s_k(y) \SBS d_i(w)$ and $s_k(y) \SBS w$. In order to prove that $s_k(y)$ is an $(n,i)$-composition, we will show that $d_i s_k(y) \supseteq d_i(w)$ and $s_k(y) \supseteq w$, as follows.
\MS

Start with any $b_i \in d_i(w)$ arising from $b \in w$ where $b = \LIST{b}{n+1}$. We aim to show that $b \in s_k(y)$ from which it will follow that $b_i \in d_i s_k(y)$.

For each $p \neq i$
\[
b_p  \in d_p(w) = d_p s_k(y) = \left\{
\begin{array}{ll}
s_{k-1} d_p(y) & \text{if } p<k\\
& \\
y & \text{if } p=k, k+1\\
&\\
s_k d_{p-1}(y) & \text{if } p>k+1
\end{array}
\right.
\]
Note that $e_k(b_{k+1}) = e_k(b_k)$ because $w$ is e-simplicial.
\MS

For the calculations below, choose any $p \in [n+1]-\SETT{i,k,k+1}$; this is possible because $n \geq 2$.
\MS

First we will show that $b_k = b_{k+1}$.

If $p<k$ then $b_p \in s_{k-1} d_p(y)$ implies $b_p = c_{k-1}(t)$ for some $t \in d_p(y)$. Then $e_p(b_k) = e_p(b_{k+1})$ because
\[
e_p(b_k) = e_{k-1}(b_p) = e_{k-1}(c_{k-1}(t)) = t = e_k(b_p) = e_p(b_{k+1})
\]
From $e_k(b_{k+1}) = e_k(b_k)$ and $e_p(b_k) = e_p(b_{k+1})$, the 2-lemma implies $b_k = b_{k+1}$.
\MS

If $p>k+1$ then $b_p \in s_k d_{p-1}(y)$ which implies that $b_p = c_k(t)$ for some $t \in d_{p-1}(y)$. Then $e_{p-1}(b_k) = e_{p-1}(b_{k+1})$ because
\[
e_{p-1}(b_k) = e_k(b_p) = e_k(c_k(t) = t = e_{k+1}(b_p) = e_{p-1}(b_{k+1})
\]
From $e_k(b_{k+1}) = e_k(b_k)$ and $e_{p-1}(b_k) = e_{p-1}(b_{k+1})$, the 2-lemma implies $b_k = b_{k+1}$.
\MS

If $k \neq i$ then $b_k \in d_k s_k(y) = y$. Otherwise $k+1 \neq i$ and $b_{k+1} \in d_{k+1} s_k(y) = y$. In either case, $b_k = b_{k+1} = a \in y$. Then the 2-lemma implies $c_k(a) = b$ because
$e_k(c_k(a)) = a = b_k = e_k(b)$ and $e_{k+1}(c_k(a)) = a = b_{k+1} = e_{k+1}(b)$. Therefore, $b \in s_k(y)$ and $b_i \in d_i s_k(y)$, as claimed.

\qed

\section{Determinacy conditions}
\label{determconds}

\subsection{Need for additional properties}
The property of closure under degeneracies (page \pageref{wantcompstruct}) for the sought-after \NICOMP{n}{i} of sets $\MC{S}'$ implies that all simplices of $\MC{S}'$ at and below dimension $n+1$ must be subface-simplicial (Theorem \refpage{downwards-e-simplicial-thm}) and that those simplices satisfy certain surjectivity conditions (Corollary \refpage{downward-surj-corollary}). The Degenerate Compositions theorem of the previous section gave conditions implying that {\em all degenerate $(n+1)$-simplices are $(n,i)$-compositions.}

Therefore, we now know that if there exists the sought-after subcomplex $\MC{S}'$ of $\mathcal{S}$ as described in section \refpage{wantcompstruct} then it will have to consist of subface-simplicial simplices satisfying certain surjectivity conditions in dimension $n$. 

However, the condition (given $n$ and $i$) that $\MC{S}'$ be closed under $(n,i)$-composition requires properties in addition to these. 
\MS

We begin with an observation made earlier on page \pageref{function-characteristic} about functions, the $(n,i)=(1,1)$ case. If $f:A \to B$ is a function of sets then any element $a \in A$ determines exactly one element $(a,b) \in \Gamma(f) \DFAS $ the graph of $f$.\MS

The analogous property in dimension $n>1$ is as follows. 
Suppose $y \in \mathcal{S}_n$ is subface-simplicial, $A \SBS [n]$ is any subset with at least two elements, and 
\[
\SETT{a_j : j \in A , a_j \in d_j(y)}
\]
is a component-simplicial partial element of $y$; that is, $e_p(a_q) = e_{q-1}(a_p)$ whenever $p<q$ and $p,q \in A$.

Then this set determines a unique $a \in \FE(V_0,\cdots,V_n)$ such that for all $j \in A$, $e_j(a) = a_j$. However, it does not follow in general that $a \in y$ because $y$ is, in general, a proper subset of $\FE(V_0,\cdots,V_n)$. In fact, even if for all $p \in [n]$ one has $e_p(a) \in d_p(y)$ it need not be the case that $a \in y$. As an extreme example, $y$ could be the empty relation with non-empty faces and non-empty e-simplicial hull (page \pageref{defSimplicial}).

If, however, $y$ happens to have the special property that each $a \in y$ is determined by $\SETT{e_j(a)\; : \; j \in A}$ we will say that $y$ ``is \DET{A}'' and that $y$ satisfies a ``determinacy'' property.
\MS

In this section we will define determinacy properties generally and establish some basic facts. Then in later sections, we will see how, given $n$ and $i$, the requirement of closure under $(n,i)$-composition points to sufficient determinacy conditions.

We will work out a specific example in section \refpage{specexample} with $(n,i)=(6,3)$. Theorems for the general $n$ and $i$ will be stated and proved in subsequent sections.

\subsection{General determinacy conditions}

\begin{newdef}
{\bf (Determinacy)}
\label{determinacy}
\index{$A$-indexed partial element}
\index{$A$-determinate}
\index{A-determinate}
\index{$A$-component}
\index{determinacy condition}
\index{$\DPQ{p}{q}{n}$}
\index{(p,q).n}
\index{$A\text{-det}[n]$}
{\rm
\SL

Suppose $m \geq 2$, $y \in S_m$ and $y$ is subface-simplicial.
 Suppose $A \subsetneq [m]$ has at least two elements. Recall (page \pageref{defSimplicial}) that an $A$-indexed partial element of $y$ is a component-simplicial partial element
\[
\SETT{a_j\, : \, j \in A \text{ and } a_j \in d_j(y)}
\]
and that the set of $A$-indexed partial elements of $y$ is denoted $\PARL{A}(y)$.

\begin{enumerate}

\item 
The simplex $y \in S_m$ will be said to be {\bf $A$-determinate} if for every $A$-indexed partial element $\SETT{a_j\, : \,j \in A \text{ and } a_j \in d_j(y)}$ there exists $a \in y$ such that $e_j(a) = a_j$ for each $j \in A$. Note that the element $a$ is {\em unique} according to the ``2-lemma''  (page \pageref{two-element-uniqueness}).

\item 
Given $m$ and $A$ as above, the condition of an $m$-simplex being $A$-determinate will be referred to as {\bf $A$-determinacy} and abbreviated $A$-det, with dimension $m$ understood from the context.

Later, we will specialize to $A = \SETT{p,q}$ and abbreviate the property of $A$-determinacy in dimension $n$ by the notation ``$\DPQ{p}{q}{n}$''.

\item 
Given $m,y$ and $A$ as above  and $a \in y$ then the set $\SET{e_j(a)}{j \in A}$ will be denoted $e_A(a)$ and be called the {\bf $A$-component} of $a \in y$.

\item 
If $\MC{S}'$ is a subcomplex of $\mathcal{S}$ with the property that, for some $n>2$ and some $A \SBS [n]$ with $2 \leq |A| < n+1$, every $y \in \MC{S}'_n$ is \DET{A}, then we will express this briefly as ``$\MC{S}'_n$ has $A$-determinacy''.
\end{enumerate}
}
\BX
\end{newdef}
\MS

It follows from the 2-lemma \refpage{two-element-uniqueness} that each $A$-indexed partial element $\SETT{a_j:j \in A}$ of $y$ determines a unique $b \in \Sig(y)$ such that for all $j \in A$, $e_j(b) = a_j$. That is, there is a monic from the set of $A$-indexed partial elements of $y$ to $\Sig(y)$.

The relationship between $\PARL{A}(y)$, the set $\ELEM{y}$ of elements of $y$ and $\Sig(y)$ is indicated in the following diagram:
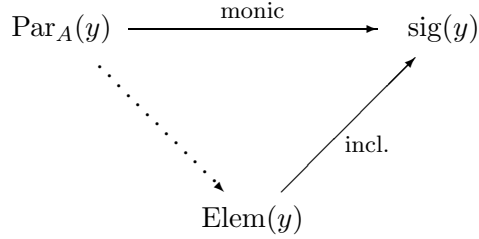
\begin{center}
\begin{picture}(2,1.1)
\put(0,1){\MB{\PARL{A}(y)}}
\put(2,1){\MB{\Sig(y)}}
\put(1,0){\MB{\ELEM{y}}}
\put(.35,1){\vector(1,0){1.3}}
\put(1.15,.15){\vector(1,1){.7}}
\multiput(.2,.8)(.05,-.05){13}{\MB{\cdot}}
\put(.85,.15){\vector(1,-1){0}}
\put(1.6,.4){\MBS{\text{incl.}}}
\put(1,1.1){\MBS{\text{monic}}}
\end{picture}
\end{center}
The dotted map exists iff $y$ is \DET{A}.

$A$-determinacy in dimension $m\geq 2$ may be construed as a kind of closure property or partial operation.

The definition of $y \in \mathcal{S}_m$ being \DET{A}\ applies to general $y$, without assuming $y$ is either \SBSI\ or component-simplicial. However, we will apply the idea only to \SBSI\ simplices.
\MS

\begin{example}
If $y \in \mathcal{S}_4$ is subface-simplicial, and $A=\SETT{1,3,4} \SBS [4]$, then a typical $A$-indexed partial-element of $y$, $\SETT{a_1,a_3,a_4}$, may be written as $(-,a_1,-,a_3,a_4)$ where $a_j \in d_j(y)$ for $j=1,3,4$. Since $y$ is subface-simplicial, there's a unique $a \in \FSIG(y)$ with $e_j(a) = a_j$, $j=1,3,4$. To say that $y$ is $\SETT{1,3,4}$-determinate is to say that for every such partial-element $(-,a_1,-,a_3,a_4)$ the unique $a=(a_0,a_1,a_2,a_3,a_4) \in \Sig(y)$ also belongs to $y$.
\end{example}\MS

\begin{newlem}
\label{subdet}
{\rm
\SL 

Suppose $m \geq 2$ and $A \SBS A' \subsetneq [m]$, where $A$ has at least two elements.
Then whenever $y \in \mathcal{S}_m$ is subface-simplicial and \DET{A}, then $y$ is also \DET{A'}.
}
\end{newlem}

\NI {\bf Proof:} If $\SET{a_j}{j \in A'}$ is an $A'$-indexed partial element of $y$ then the 2-row lemma, the assumption that $A \SBS A'$, and the $A$-determinacy of $y$ together imply that $a \in y$.\\
\qed
\MS

In the next section we will develop some facts about determinacy which we will use to specify models $\MC{S}'$ for \NICOMP{n}{i}s.

\subsection{Downward propagation of determinacy}
For subface-simplicial simplices, determinacy conditions imposed at dimension $m$ ``propagate downwards'' in dimension via degeneracies.
\MS

For $m \geq 1$, the goal in this section is to relate determinacy conditions for $x \in \mathcal{S}_m$ in terms of those for $s_k (x)$ and various $k \in [m]$.
\MS

First we consider a special case of \DETY{A} for $s_k(x)$ if it happens to be the case that $k$ or $k+1$ belongs to $A$.

\begin{newlem}
\label{restricted-degenerate-determinacy}
{\rm
\SL 

Suppose $x \in \MC{S}_m$ is \SBSI, $p \in [m]$, $A \subsetneq [m+1]$ where $|A| \geq 2$, and exactly one of $p$ and $p+1$ belongs to $A$. Then $s_p(x)$ is \DET{A}.
}
\end{newlem}

\Proof

We will consider the case when $p \in A$ and $p+1 \notin A$. The case when $p \notin A$ and $p+1 \in A$ is handled similarly.

Let $\SETT{b_t:t \in A}$ be any $A$-indexed partial element of $s_p(x)$. Since $p \in A$ by hypothesis then $b_p \in d_ps_p(x)=x$. We will show that $c_p(b_p)$ is the unique filler of $\SETT{b_t:t \in A}$ i.e. that $e_t c_p(b_p) = b_t$ for each $t \in A$. 

Given any $t \in A$, $t \neq p$ then either $t>p+1$ or $t<p$.

If $t>p+1$ then $b_t \in d_t s_p(x) = s_p d_{t-1}(x)$ implies $b_t = c_{p} (a)$ for some $a \in d_t(x)$. Then
\[
e_t c_p(b_p) = c_p e_{t-1}(b_p) = c_p e_p(b_t) = c_p e_p c_{p}(a) = b_t
\]
as required.

If $t<p$ then $b_t \in d_t s_p(x) = s_{p-1} d_t(x)$ implies $b_t = c_{p-1}(a)$ for some $a \in d_t(x)$. Then
\[
e_t c_p(b_p) = c_{p-1} e_t(b_p) = c_{p-1} e_{p-1}(b_t) = c_{p-1} e_{p-1}c_{p-1}(a) = b_t
\]
This shows that $c_p(b_p) \in s_p(x)$ is the unique filler of $\SETT{b_t:t \in A}$, and therefore that $s_p(x)$ is \DET{A}.

\qed
\BS

To illustrate this lemma: Given $p<q$ in $[m+1]$ then: $s_p(x)$ is \DET{(p,q)} if $p<q-1$ and $s_q(x)$ is \DET{(p,q)} if $q>p+1$.
\BS

Now if both $p$ {\em and} $p+1$ belong to $A$ then the conclusion of the previous lemma fails because given any $A$-indexed partial element $\SETT{b_t: t \in A}$ of $s_p(x)$ then both $b_p \in x$ and $b_{p+1} 
\in x$ with $e_p(b_{p+1}) = e_p(b_p)$. In general however, $b_p \neq b_{p+1}$.

Here is a low-dimensional illustration.
\MS

\begin{example}
{\rm
For notational convenience, we denote certain fundamental entries by 0, 1, 2, 3 and 4. Let $b_2=h_3(0,1,2,3)$ and $b_3 = h_3(0,1,4,3)$ and suppose the element set of $x \in \MC{S}_3$ is
\[
\ELEM{x} = \SETT{b_2,b_3} = \SETT{h_3(0,1,2,3), h_3(0,1,4,3)}
\]
Then
\[
\ELEM{s_2(x)} = \SETT{h_4(0,1,2,2,3), h_4(0,1,4,4,3)}
\]
$\SETT{b_2, b_3}$ is an $\{2,3\}$-indexed partial element of $s_2(x)$ since $b_2 \in d_2 s_2(x)=x$, $b_3 \in d_3 s_2(x)=x$ and $e_2(b_3) = h_2(0,1,3) = e_2(b_2)$. The unique filler of this partial element is $h_4(0,1,4,2,3)$ which is {\em not} in the image of $c_2 : x \to s_2(x)$. That is, $s_2(x)$ is not $\{2,3\}$-determinate.
}
\end{example}

We now continue with the case that $k, k+1 \notin A$.
\MS

\NI In preparation:
\MS

\NI {\bf Integer-interval notation:} 
\index{$[p,q]$}
\index{[p,q] integer interval}
We will use interval notation for subsets of integers. If $p,q \in \ZZ$ then 
\label{integer-interval-notation}
\[
[p,q] \DFAS
\left\{
\begin{array}{ll}
\SETT{p,p+1, \cdots , q} & \text{ if }p<q\\
\SETT{p} & \text{ if } p=q\\
\emptyset & \text{ if } p>q
\end{array}
\right.
\]
\MS

\NI {\bf Notation:} Suppose $m \geq 2$ and $A = \SETT{a_0 \DDD{<} a_n} \subsetneq [m+1]$ where $n \geq 1$. Then $A$ corresponds to the strictly increasing function
\[
\DMAP{\LAM_A}{n}{m+1}, \quad \LAM_A(t) \DFAS a_t
\]
Given any $k \in [m+1]$ such that $A \cap \{k,k+1\} = \MT$, then 
\[
[n] \XRA{\LAM_A} [m+1] \XRA{\sigma_k} [m]
\]
is also strictly increasing, and we will define
\[
\sigma_k(A) \DFAS \SETT{\sigma_k \LAM_A(t): t \in [n]} \SBS [m]
\]
Since
\[
\sigma_k\LAM_A(t) = \left\{
\begin{array}{ll}
\LAM_A(t) & \text{if } 0 \leq \LAM_A(t) \leq k-1\\
\LAM_A(t)-1 & \text{if } k+2 \leq\LAM_A(t) \leq m+1
\end{array}
\right.
\]
then
\[
\sigma_k(A) = \SETT{p: p \in A \cap [0,k-1]}\, \cup \, \SETT{p-1: p \in A \cap [k+2,m+1]}
\]
For example, $\SIG_6\big(\SETT{0,1,4,7,11,12}\big)=\SETT{0,1,4} \cup \SETT{6,10,11}$.

\BX
\MS

To set up the next lemma, we make two observations concerning partial elements. Suppose $m \geq 2$, $x \in \MC{S}_m$ where $x$ is \SBSI, $A \SBS [m+1]$ where $|A| \geq 2$, $k \in [m]$ and $A\cap \SETT{k,k+1} =\MT$.
\begin{itemize}
\item 
An $A$-indexed partial element of $s_k(x)$ is $\SETT{b_p \in d_p s_k(x): p \in A}$ such that for all $p<q$ in $A$, we have $e_p(b_q)=e_{q-1}(b_p)$. Since $|A|\geq 2$ and $s_k(x) $ is \SBSI\ then there is a unique $b \in \FSIG(s_k(x))$ such that for all $p \in A$, $b_p = e_p c_k(b)$. Using that
\[
e_p c_k(b) = \left\{
\begin{array}{ll}
 c_{k-1}e_p(b) & p \in A \cap [0,k-1]\\
 & \\
 c_k e_{p-1}(b) & p \in A \cap [k+2,m+1]
\end{array}
\right.
\]
it follows that the $A$-indexed partial element of $s_k(x)$ splits into
\[
\SETT{c_{k-1}e_p(b) : p \in A \cap [0,k-1]} \cup 
\SETT{c_k e_{p-1}(b) : p \in A \cap [k+2,m+1]}
\]

\item
An $\SIG_k(A)$-indexed partial element of $x$ is $\SETT{b_p \in d_p(x):p \in \SIG_k(A)}$. Since $|\SIG_k(A)| \geq 2$ also, there is a unique $b \in \FSIG(x)$ such that for all $p \in \SIG_k(A)$, $b_p = e_p(b)$. Therefore this $\SIG_k(A)$-indexed partial element of $x$ is
\[
\SETT{e_p(b): p \in A \cap [0,k-1]} \cup
\SETT{e_{p-1}(b):p \in A \cap [k+2,m+1]}
\]
\end{itemize}

\begin{newlem}
{\rm
\SL 

Suppose $m \geq 2$, $k \in [m]$, $A \SBS [m+1]- \{k,k+1\}$ and $|A| \geq 2$.
Let $x \in \mathcal{S}_m$ and suppose $x$ is subface-simplicial. Then there is a bijective correspondence between $A$-indexed elements of $s_k(x)$ and $\sigma_k(A)$-indexed elements of $x$:
\[
\SETT{e_p(b): p \in \sigma_k(A)} \leftrightarrow
\SETT{e_p c_k(b): p \in A}
\]
}
\end{newlem}

\Proof

First, suppose $\SETT{b_p \in d_p s_k(x): p \in A}$ is an $A$-indexed partial element of $s_k(x)$. As in the first observation above, there is a unique $b \in \FSIG(s_k(x))$ such that this partial element is 
\[
\SETT{c_{k-1}e_p(b) : p \in A \cap [0,k-1]} \cup 
\SETT{c_k e_{p-1}(b) : p \in A \cap [k+2,m+1]}
\]
which (since $c_{k-1}$ and $c_k$ are monic) determines the set
\[
\SETT{e_p(b) : p \in A \cap [0,k-1]} \cup 
\SETT{e_{p-1}(b) : p \in A \cap [k+2,m+1]} = \SETT{e_p(b): p \in \SIG_k(A)}
\]
that is, a $\SIG_k(A)$-indexed partial element of $x$.
\MS

Second, suppose $\SETT{b_p \in d_p(x):p \in \SIG_k(A)}$ is an $\SIG_k(A)$-indexed partial element of $x$. As noted in the second observation just above, there exists a unique $b \in \FSIG(x)$ such that for all $p \in \SIG_k(A)$, $e_p(b) = b_p$. This determines
\[
\SETT{c_{k-1}e_p(b) : p \in A \cap [0,k-1]} \cup 
\SETT{c_k e_{p-1}(b) : p \in A \cap [k+2,m+1]}
\]
which, as noted above, is an $A$-indexed partial element of $s_k(x)$.

It follows directly that the relationship between $A$-indexed partial elements of $s_k(x)$ and $\SIG_k(A)$-indexed partial elements of $x$ is bijective.

\qed

\begin{newcor}{\bf (Degenerate Determinacy)}
\label{degenDetlemma} 
\index{degenerate determinacy theorem}
\index{Theorem!Degenerate determinacy theorem}
{\rm  
\SL

\NI Assume:
\begin{itemize}
\item $m \geq 2$ and $k \in [m]$
\item $A \subsetneq [m+1]-\SETT{k,k+1}$ and $A$ has at least two elements.
\item $x \in S_m$ and $x$ is subface-simplicial.
\end{itemize}
Then $s_k(x)$ is $A$-determined if and only if $x$ is $\SIG_k (A)$-determined.
}
\end{newcor}

\Proof

Suppose $s_k(x)$ is $A$-determined. Given any $\SIG_k (A)$-indexed partial element of $x$, $\SETT{e_p(b): p \in \SIG_k (A)}$ for some $b \in \FSIG(x)$, we must show that $b \in x$. By the previous lemma, $\SETT{e_p c_k(b):p \in A}$ is an $A$-indexed partial element of $s_k(x)$. Since $s_k(x)$ is assumed $A$-determinate then $c_k(b) \in s_k(x)$, and therefore $b \in x$.

Conversely, suppose $x$ is $\SIG_k (A)$-determined. Given any $A$-indexed partial element of $s_k(x)$, $\SETT{e_p c_k(b) : p \in A}$ for some $b \in \FSIG(x)$, we must show $c_k(b) \in s_k(x)$. Again, by the previous lemma, $\SETT{e_p(b): p \in \SIG_k (A)}$ is a $\SIG_k (A)$-indexed partial element of $x$. Since $x$ is assumed $\SIG_k (A)$-determinate, $b \in x$ and therefore $c_k(b) \in s_k(x)$.

\qed
\BS

It follows as a consequence of this corollary that if $\MC{S}'$ has a hypothetical $(n,i)$-structure as described on page \pageref{wantcompstruct}, then requiring that {\em all}\, $n$-simplices be $A$-determined (in particular, all degenerate ones) implies that all $(n-1)$-simplices be $(\SIG_k A)$-determined whenever $k,k+1 \notin A$. The requirement for $(n-1)$-simplices propagates downwards in turn.
\BS

\subsection{Up and down translation of $\{p,q\}$-determinacy}

In this section we will specialize $A$-determinacy  to when $A=\SETT{p,q}$. 
\MS

\begin{newdef}
\label{pq-partial-element}
\index{(p,q)-indexed partial element}
\index{(p,q)-partial element}
{\rm
\SL

Suppose $n \geq 2$, $y \in \MC{S}_n$ is \SBSI\ and $A = \SETT{p,q}$ where $0 \leq p<q \leq n$. An $A$-indexed partial element of $y$ will be called a {\bf $(p,q)$-indexed partial element of $y$}. That is, it is a set $\SETT{b_p,b_q}$ where $b_p \in d_p(y)$, $b_q \in d_q(y)$ and $e_p(b_q)=e_{q-1}(b_p)$. For brevity: ``$(p,q)$-partial element of $y$''.
}

\BX
\end{newdef}

{\bf For brevity we will always assume $p<q$.} We will shorten notation and phrasing as follows.
\begin{itemize}

\item For ``$y \in \mathcal{S}_n$ is \DET{\SETT{p,q}}'' write ``$y$ is $\DPQ{p}{q}{n}$''. \index{$\DPQ{p}{q}{n}$} \index{(p,q).n}

\item Given a subcomplex $\MC{S}'$ of $S$, we will say ``$\MC{S}'$ has $\DPQ{p}{q}{n}$'' to mean ``all $n$-simplices of $\MC{S}'$ are \DET{\SETT{p,q}}''.

\end{itemize}
\MS
%
%

A $(p,q)$-partial element of $y$ determines partial elements in the faces $d_j(y)$ for $j \neq p,q$, as illustrated in the next example.
\MS

\begin{example}
Suppose $y \in \MC{S}_6$ and $y$ is \SBSI. A $(1,4)$-partial element of $y$ is $\SETT{b_1,b_4}$ where $b_1 \in d_1(y), b_4 \in d_4(y)$ and $e_1(b_4) = e_3(b_1)$. The array of such a partial element is:

\[
\begin{array}{c||c|c|c|c|c|c||c}
 & 0 & 1 & 2& 3& 4& 5 & \\
\hline
\hline
0 & b_{1,0} &  &  & b_{4,0} &  &  &  \\
\hline
1 & b_{1,0} & b_{1,1} & b_{1,2} & b_{1,3} & b_{1,4} & b_{1,5} & b_1\\
\hline
2 &  & b_{1,1} &  & b_{4,3} &  &  & \\
\hline
3 &  & b_{1,2} &  & b_{4,3} &  &  & \\
\hline
4 & b_{4,0} & b_{4,1} & b_{4,2} & b_{4,3} & b_{4,4} & b_{4,5} & b_4 \\
\hline
5 &  & b_{1,4} &   &   & b_{4,4} &  &  \\
\hline
6 &  & b_{1,5} &   &   & b_{4,5} &  &  \\
\hline
\end{array}
\]
where $b_{4,1} = e_1(b_4) = e_3(b_1) = b_{1,3}$. The partial element $\SETT{b_1, b_4}$ determines partial elements in the faces $d_j(y)$, $j=0,2,3,5,6$.
\end{example}
\MS

The next two theorems relate determinacy in dimension $n$ to: (i) determinacy in dimension $n+1$ (the ``Up Rule''), and (ii) determinacy in dimension $n-1$ (the ``Down Rule'').
\MS

\begin{newthm}
{\bf (Up Rule)}
\label{up rule}
\index{Up Rule}
{\rm

Suppose $n \geq 2$ and $0 \leq p<q \leq n$. If $z \in \MC{S}_{n+1}$ is maximal e-simplicial and for all $i \in [n+1]$, $d_i(z)$ is e-simplicial and \DET{(p,q)} then
\MS

\NI U1: $z$ is $\DPQ{p}{q}{(n+1)}$.
\MS

\NI U2: $z$ is $\DPQ{p}{q+1}{(n+1)}$.
\MS

\NI U3: $z$ is $\DPQ{p+1}{q+1}{(n+1)}$.
}
\end{newthm}

\Proof

For brevity, we'll denote $d_j(z)$ by $y_j$, $j \in [n+1]$.
\MS

\NI \fbox{That $z$ is \DET{(p,q)}}
\MS

Let $\{b_p, b_q\}$ be a $(p,q)$-partial element of $z$. That is, $b_p \in y_p =d_p(z)$, $b_q \in y_q =d_q(z)$ and $e_p(b_q)=e_{q-1}(b_p)$. The partial array of this partial element of $z$ is:

\[
\begin{array}{|c||c|c|c|c|c|c|c|c|c|c||c|}
\hline
\text{row} \downarrow & 0 & \cdots & p-1 & p & \cdots & q-1 & q & \cdots & n-1 & n & \\
\hline \hline
0 & & & 1 & & & 2 & & & & & \\
\hline
\vdots & & & \vdots & & & \vdots & & & & & \\
\hline
p & 1 & \cdots & 1 & 1 & \cdots & 1 & 1 & \cdots & 1 & 1 & b_p \\
\hline
\vdots & & & & \vdots & & \vdots & & & & & \\
\hline
q & 2 & \cdots & 2 & 1 & \cdots & 2 & 2 & \cdots & 2 & 2 & b_q \\
\hline
\vdots & & & & \vdots & & & \vdots & & & & \\
\hline
n+1 & & & & 1 & & & 2 & & & & \\
\hline
\end{array}
\]
The case when $p=0$ is similar.

Rows $q+1$ through $n+1$, of which there is at least one (namely row $q+1$), fill one at a time compatibly using the assumption that $y_{q+1}, \DDD{} , y_{n+1}$ are $\DPQ{p}{q}{n}$. The result is a partial array with filled entries looking like:

\[
\begin{array}{|c||c|c|c|c|c|c|c|c|c|c||c|}
\hline
\text{row} \downarrow & 0 & \cdots & p-1 & p & \cdots & q-1 & q & \cdots & n-1 & n & \\
\hline \hline
0 &  &  & 1 & & & 2 & / & / & / & / & \\
\hline
\vdots & & & \vdots & & & \vdots & / & / & / & / & \\
\hline
p & 1 & \cdots & 1 & 1 & \cdots & 1 & 1 & \cdots & 1 & 1 & b_p \\
\hline
\vdots & & & & \vdots & & \vdots & / & / & / & / & \\
\hline
q & 2 & \cdots & 2 & 1 & \cdots & 2 & 2 & \cdots & 2 & 2 &b_q\\
\hline
\vdots & / & \cdots & / & \vdots & / & / & \vdots & / & / & / & \\
\hline
n+1 & / & \cdots & / & 1 & / & / & 2 & / & / & / & \\
\hline
\end{array}
\]
where the new entries are compatible elements of various $y_j$, all indicated simply by ``/''. We note that all the entries in column 2 are present. 

Next, we use \DETY{(p,q)} to fills rows $p+1$ through $q$ one by one. These new entries fill in column $p$. In the case that $p+1=q$ then column $p$ has already been filled. Now, rows 0 through $p-1$ (if any) contain entries including those in columns $p$ and $q$. Again, these rows can be filled one by one compatibly using \DETY{(p,q)}. Now, for each $k$, the entries in row $k$ comprise an element $t_k$ of $y_k$ since $y_k$ is \DET{(p,q)}. Also the $t_k$ satisfy by construction the requirements $d_j(t_k)=d_{k-1}(t_j)$ for all $j<k$. That $z$ is maximal e-simplicial implies $\LIST{t}{n+1} \in z$, as claimed.
\BS

\NI \fbox{That $z$ is \DET{(p,q+1)}}
\MS

The reasoning is similar to the first case. Let $\{b_p, b_{q+1}\}$ be a $(p,q+1)$-partial element of $z$. The array for that partial element of $w$ (shown for when $0<p<q<n$):
\[
\begin{array}{|c||c|c|c|c|c|c|c|c|c|c|c||c|}
\hline
\text{row} \downarrow & 0 & \cdots & p-1 & p & \cdots & q-1 & q & q+1 & \cdots & n-1 & n & \\
\hline \hline
0 & & & 1 & & & & 2 & & & & & \\
\hline 
\vdots & & & \vdots & & & & \vdots & & & & & \\
\hline 
p & 1 & \cdots  & 1 & 1 & \cdots & 1 & 1 & 1 & \cdots & 1 & 1 & b_p\\
\hline 
\vdots & & & & \vdots & & & \vdots & & & & & \\
\hline 
q & & & & 1 & & & 2 & & & & & \\
\hline 
q+1 & 2 & \cdots & 2 & 1 & \cdots & 2 & 2 & 2 & \cdots & 2 & 2 & b_{q+1}  \\
\hline 
\vdots & & & & \vdots & & & &\vdots & & & & \\
\hline
n+1 & & & & 1 & & & & 2 & & & & \\
\hline
\end{array}
\]
Since $y_p$ and $y_{q+1}$ are \DET{(p,q)}, rows $p+1$ through $q$ fill, one by one, uniquely and compatibly. The resulting partial array becomes:
\[
\begin{array}{|c||c|c|c|c|c|c|c|c|c|c|c||c|}
\hline
\text{row} \downarrow & 0 & \cdots & p-1 & p & \cdots & q-1 & q & q+1 & \cdots & n-1 & n & \\
\hline  \hline
0 & & & 1 & / & & & 2 & & & & & \\
\hline 
\vdots & & & \vdots & / & & & \vdots & & & & & \\
\hline 
p & 1 & \cdots  & 1 & 1 & \cdots & 1 & 1 & 1 & \cdots & 1 & 1 & b_p\\
\hline 
\vdots & / & / & / & \vdots & \cdots & / & \vdots & / & \cdots & / & / & \\
\hline 
q & / & / & / & 1 & \cdots & / & 2 & / & \cdots & / & / & \\
\hline 
q+1 & 2 & \cdots & 2 & 1 & \cdots & 2 & 2 & 2 & \cdots & 2 & 2 & b_{q+1} \\
\hline 
\vdots & & & & \vdots & & & / &\vdots & & & & \\
\hline
n+1 & & & & 1 & & & / & 2 & & & & \\
\hline
\end{array}
\]
At this stage, all the other rows have filled entries in positions $p$ and $q$. From \DETY{(p,q)} it follows that these rows may also be filled, one by one, with elements $t_k$ of the corresponding $y_k$. And by the same reasoning as in the previous case, $\LIST{t}{n+1} \in z$, as claimed.
\BS

\NI \fbox{That $z$ is \DET{(p+1,q+1)}}
\MS

As in the previous cases, start with the $(p+1,q+1)$-partial element $\{b_{p+1},b_{q+1}\}$ and the corresponding partial array of $z$
\[
\begin{array}{|c||c|c|c|c|c|c|c|c|c|c|c|c||c|}
\hline
\text{row} \downarrow & 0 & \cdots & p-1 & p & p+1 & \cdots & q-1 & q & q+1 & \cdots & n-1 & n & \\
\hline \hline
0 & & & & 1 & & & & 2 & & & & & \\
\hline
\vdots & & & & \vdots & & & & \vdots & & & & & \\
\hline
p & & & & 1 & & & & 2 & & & & & \\
\hline
p+1 & 1 & \cdots & 1 & 1 & 1 & \cdots & 1 & 1 & 1 & \cdots & 1 & 1 & b_{p+1} \\
\hline
\vdots & & & & & \vdots & & & \vdots & & & & & \\
\hline
q & & & & & 1 & & & 2 & & & & & \\
\hline
q+1 & 2 & \cdots & 2 & 2 & 1 & \cdots & 2 & 2 & 2 & \cdots & 2 & 2 & b_{q+1} \\
\hline
\vdots & & & & & \vdots & & & & \vdots & & & & \\
\hline
n+1 & & & & & 1 & & & & 2 & & & & \\
\hline
\end{array}
\]
Rows 0 through $p$ (of which there is at least one) may be filled one by one uniquely and compatibly using \DETY{(p,q)} to obtain the partial array
\[
\begin{array}{|c||c|c|c|c|c|c|c|c|c|c|c|c||c|}
\hline
\text{row} \downarrow & 0 & \cdots & p-1 & p & p+1 & \cdots & q-1 & q & q+1 & \cdots & n-1 & n & \\
\hline \hline
0 & / & \cdots & / & 1 & / & \cdots & / & 2 & / & \cdots & / & / & \\
\hline
\vdots & / & \cdots & / & \vdots & / & \cdots & / & \vdots & / & \cdots & / & / & \\
\hline
p & / & \cdots & / & 1 & / & \cdots & / & 2 & / & \cdots & / & / & \\
\hline
p+1 & 1 & \cdots & 1 & 1 & 1 & \cdots & 1 & 1 & 1 & \cdots & 1 & 1 & b_{p+1} \\
\hline
\vdots & / & \cdots & / & / & \vdots & & & \vdots & & & & & \\
\hline
q & / & \cdots & / & / & 1 & & & 2 & & & & & \\
\hline
q+1 & 2 & \cdots & 2 & 2 & 1 & \cdots & 2 & 2 & 2 & \cdots & 2 & 2 & b_{q+1} \\
\hline
\vdots & / & \cdots & / & / & \vdots & & & & \vdots & & & & \\
\hline
n+1 & / & \cdots & / & / & 1 & & & & 2 & & & & \\
\hline
\end{array}
\]
This has the effect of filling in the entries in positions $p$ and $q$ in rows $p+1$ through $q$ (of which there is at least one) which yields
\[
\begin{array}{|c||c|c|c|c|c|c|c|c|c|c|c|c||c|}
\hline
\text{row} \downarrow & 0 & \cdots & p-1 & p & p+1 & \cdots & q-1 & q & q+1 & \cdots & n-1 & n & \\
\hline \hline
0 & / & \cdots & / & 1 & / & \cdots & / & 2 & / & \cdots & / & / & \\
\hline
\vdots & / & \cdots & / & \vdots & / & \cdots & / & \vdots & / & \cdots & / & / & \\
\hline
p & / & \cdots & / & 1 & / & \cdots & / & 2 & / & \cdots & / & / & \\
\hline
p+1 & 1 & \cdots & 1 & 1 & 1 & \cdots & 1 & 1 & 1 & \cdots & 1 & 1 & b_{p+1} \\
\hline
\vdots & / & \cdots & / & / & / & \cdots & /& / & / & \cdots &/ &/ & \\
\hline
q & / & \cdots & / & / & 1 & \cdots & / & 2 & / & \cdots & / & / & \\
\hline
q+1 & 2 & \cdots & 2 & 2 & 1 & \cdots & 2 & 2 & 2 & \cdots & 2 & 2 & b_{q+1} \\
\hline
\vdots & / & \cdots & / & / & / &  \cdots & / & / & \vdots & & & & \\
\hline
n+1 & / & \cdots & / & / & 1 & \cdots & / & / & 2 & & & & \\
\hline
\end{array}
\]
At this stage rows $q+2$ to $n+1$, if any, fill uniquely and compatibly by \DETY{(p,q)} yielding an element of $z$, as claimed.

\qed
\BS

\begin{newthm}
\label{down rule}
\index{Down Rule}
{\bf (Down Rule)}
{\rm

Suppose $n \geq 2$, $0 \leq p<q \leq n+1$ and $y \in \MC{S}_n$ is \SBSI. If for all $k \in [n]$, $s_k(y)$ is $\DPQ{p}{q}{n+1}$ then:
\MS

\NI D1: $p \geq 2 $ implies $y$ is $\DPQ{p-1}{q-1}{n}$.

\NI D2: $q-p \geq 3 $ implies $y$ is $\DPQ{p}{q-1}{n}$.

\NI D3: $q \leq n-1 $ implies $y$ is $\DPQ{p}{q}{n}$.
}
\end{newthm}

\Proof

This is a direct application of Corollary \refpage{degenDetlemma}. If $p \geq 2$ then $\{p,q\} \cap \{k,k+1\}=\MT$ when $k=0$. If $q-p \geq 3$ then $\{p,q\} \cap \{k,k+1\}=\MT$ when $k=p+1$. If $q \leq n-1$ then $\{p,q\} \cap \{k,k+1\}=\MT$ when $k=n$.
\MS

\NI For D1: If $p \geq 2$ then, by Corollary \refpage{degenDetlemma}, $s_0(y)$ is $\DPQ{p}{q}{n}$ iff $y$ is \DET{\SIG_0(\{p,q\})}, and $\SIG_0(\{p,q\})=\{p-1,q-1\}$.
\MS

\NI For D2: if $q-p \geq3$ then $s_{p+1}(y)$ is $\DPQ{p}{q}{n}$ iff $y$ is \DET{\SIG_{p+1}(\{p,q\})}, and $\SIG_{p+1}(\{p,q\})=\{p,q-1\}$.
\MS

\NI For D3: if $q \leq n$ then $s_{p+1}(y)$ is $\DPQ{p}{q}{n}$ iff $y$ is \DET{\SIG_{p+1}(\{p,q\})}, and $\SIG_{q+1}(\{p,q\})=\{p,q\}$.

\qed
\MS

Note that the Up Rule applies to any maximal e-simplicial $(n+1)$-simplex in $\MC{S}$; these include all degeneracies $s_k(y)$ for e-simplicial $n$-simplices $y$ (lemma \refpage{degenerates are maximal e-simplicial}), and to all $(n,i)$-compositions (Theorem \refpage{max-property}).
\MS

\NI {\bf Notes:} 

\begin{enumerate}

\item 
Eventually, we will be able to apply Up Rule and Down Rule to \NICOMP{n}{i}s, namely subcomplexes $\MC{S}'$ of $\mathcal{S}$ where all $n$ simplices of $\MC{S}'$ are subface-simplicial and where all $(n+1)$-simplices of $\MC{S}'$ are $(n,i)$-compositions. The somewhat technical quality of the Theorems \refpage{up rule} and \refpage{down rule} was necessary because we have not yet established that such $\MC{S}'$ exists. In fact these theorems will used below to verify that such $\MC{S}'$ exist.

\item
In section \refpage{furthermore}, we revisit the Up and Down Rules relaxing some hypotheses and examining determinacy relationships independent of $(n,i)$-composition.

\end{enumerate}

\subsection{Determinacy conditions and $(n,i)$-composition}

Suppose $n \geq 2$ and $w \in \mathcal{S}_{n+1}$ is e-simplicial. Fix any $i \in [n+1]$ and suppose $w$ is also \SUR{i}. 

Recall that if $(\OM{b}{n+1}{i})$ is a component-simplicial partial element of $w$ then there exists a unique $b_i \in \Sig(d_i(w))$ namely
\[
b_i = \bigl( e_{i-1}(b_0) \DDD{,} e_{i-1}(b_{i-1}),e_i(b_{i+1}) \DDD{,} e_i(b_{n+1}) \bigr)
\]
such that $b = (b_0, \cdots , b_i, \cdots , b_{n+1}) \in \Sig(w)$ is component-simplicial. Recall, again, that it need not be the case in general either that $b_i \in d_i(w)$ or $b \in w$.

Specifically, $w$ does not, in general, satisfy condition \ref{operation-condition} of definition \refpage{n-i-comp-definition} of $(n,i)$-composition.
\MS

Now suppose that in addition to the two given properties of $w$ above, we also assume that $w$ is \DET{A} for some $A \SBS [n+1]$ such that $|A| \geq 2$ and $i \notin A$.

\NI Then, given any component-simplicial partial element $(\OM{b}{n+1}{i})$ of $w$ then we have that 
\[
A \SBS [n+1]-\SETT{i} \text{ and }  w \text{ being \DET{A}\ }
\]
imply that $b = (b_0, \cdots, b_i, \cdots , b_{n+1}) \in w$. That is, with the extra condtion of $A$-determinacy, $w$ satisfies the ``if'' part of condition \ref{operation-condition} of the definition of $(n,i)$-composition. If we also assume that $w$ is \SUR{i} then $w$ then satisfies condition \ref{surjectivity-condition} as well.
\MS

Summarized as a lemma:
\begin{newlem} 
\label{comp-determinacy-lemma}
\index{Determinacy and Composition Lemma}
\index{lemma!Determinacy and Composition}
{\bf (Determinacy and Composition)}
{\rm
Given any $n \geq 2$ and any $i \in [n+1]$, suppose $w \in \mathcal{S}_{n+1}$ is:  e-simplicial, \SUR{i} and \DET{A} where $A \SBS [n+1]$, $|A| \geq 2$ and $i \notin A$. Then $w$ is an $(n,i)$-composition.
}
\end{newlem}
\qed

\section{Models of \NICOMP{n}{i}s}
\label{nisolution}

\subsection{Introduction}
Thus far in the search for a simplicial set $\MC{S}'$ which is the nerve of a structure with $(n,i)$-composition i.e. an \NICOMP{n}{i} of sets  (page \pageref{wantcompstruct}), we have used the closure-under-degeneracies requirement to deduce that all $m$-simplices of $\MC{S}'$ with $m \leq n+1$ must:
\begin{itemize}
\item be subface simplicial.
\item satisfy certain surjectivity conditions.
\end{itemize}

We will now examine the implications of requiring closure under $(n,i)$-composition. We have already seen by the ``e-simplicial composites'' lemma (page \pageref{simpclosure}) that if $w$ is an $(n,i)$-composition with each $d_jw$ ($j \neq i$) e-simplicial, then $d_i(w)$ is also e-simplicial. It follows immediately that if each $d_j(w)$ ($j \neq i$) subface-simplicial  then $d_i(w)$ is, in fact, also subface simplicial.

That $d_i(w)$ satisfies the required surjectivity conditions (when the $d_j(w)$ do for $j \neq i$) appears {\bf not} to follow, however, without the imposition of other conditions. The discussion to follow will show that imposing a determinacy condition in dimension $n$ suffices to yield the required surjectivity conditions. Of course, we will also have to show that the new condition is itself preserved by $(n,i)$-composition.

To illustrate this, we will first work out a specific example with small $n$.

\subsection{A specific example: $n=6$, $i=3$}
\label{specexample}

In this subsection, we will develop a \NICOMP{6}{3} $\MC{S}'$ and use this example to illustrate the general proof technique and notation to be applied later. \MS

In this example, $\mathcal{S}_{\leq 7}'$ will be a 7-truncated subcomplex of $\mathcal{S}$ where all simplices of $\MC{S}'_{\leq 7}$ are \SBSI\ and where $\MC{S}'_7$ consists of all (6,3)-compositions of 6-simplices. 
Previous results (the ``Necessity of surjectivity'' theorem, page \pageref{nicompSubfaceSimp} and the ``Degenerate $(n,i)$-compositions'' theorem, page \pageref{degen-comps-thm})
imply we must impose some surjectivity conditions on 6-simplices.
One determinacy condition suffices for closure under (6,3)-composition.
\MS

The conditions imposed on $\mathcal{S}_{6}'$ will be the following:
\begin{itemize}
\item All 6-simplices are \SUR3
\item All 6-simplices are \SUR2
\item All 6-simplices are \DET{\SETT{1,4}}. In the abbreviated notation: every $y \in \mathcal{S}_6'$ is $\DPQ{1}{4}{6}$.
\end{itemize}
As we have established earlier, these conditions imply related properties at dimensions below 6. We will not need to cite the lower-dimensional properties specifically in what follows. (See section \refpage{furthermore} for a discussion of this).
\MS

The Degenerate Compositions Theorem  \refpage{degen-comps-thm} implies that for all $y \in \MC{S}'_6$, and all $j \in [6]$ that $s_j(y)$ is a $(6,3)$-composition.\MS

To show that the conditions given above suffice, we will have to verify their preservation by $(6,3)$-composition. 
\MS

Specifically, if $w = \text{comp}_{6,3}(\OM{y}{7}{3})$ with all $y_j \in \MC{S}'_6$, $j \neq 3$, then we must verify that $d_3(w) \in \MC{S}'_6$ also.
\MS

In doing so, we will see how the requirement of $\DPQ{1}{4}{6}$ arises ``reasonably'' (in the sense of the introductory comments on page \pageref{guidelines}) from the need to establish closure of the two surjectivity conditions. 
Along the way, we will deduce that $\MC{S}'_6$ has $\DPQ{2}{4}{6}$ from $\DPQ{1}{4}{6}$. Then we will use $\DPQ{2}{4}{6}$ to establish closure of $\DPQ{1}{4}{6}$ under $(6,3)$-composition.\MS

Here is the sequence of verifications:
\begin{enumerate}
\item Closure of 3-surjectivity under (6,3)-composition.
\item Closure of 2-surjectivity under (6,3)-composition.
\item For any $(6,3)$-composition $w$ as above, $w$ is $\SETT{1,4}$-determinate and $\SETT{2,5}$-determinate.
\item All $y \in \MC{S}'_6$ are $\SETT{2,4}$-determinate. That is, $S'_6$ has $\DPQ{2}{4}{6}$.
\item Closure of $\DPQ{1}{4}{6}$ under $(6,3)$-composition
\end{enumerate}
\MS

\NI \fbox{\bf Closure of 3-surjectivity under $(6,3)$-composition:}\MS

Suppose $w = \COMP{6}{3}(y_0, y_1, y_2,-,y_4,y_5,y_6,y_7)$ 
is a (6,3)-composition where for each $j \neq 3$, $y_j \in \MC{S}'_6$. We must show that $y_3= d_3 (w)$ is also \SUR3. \MS

That is, given any arbitrary $t \in d_3 (y_3)$ we must show there exists an $a \in y_3$ with $e_3 (a) = t$. Since each element of $y_3$ comes from some element of $b \in w$ (by the 3-surjectivity of $w$ itself), it is both necessary and sufficient to show that there exists an element $b \in w$ such that $e_3 e_3 (b) = b_{33} = t$.\MS

Here is our method in outline, one which we will use repeatedly: start with $t$ as a component-of-component entry in the sought-for $b \in w$ (pictured below using the array notation developed earlier) and use some or all of the conditions above to build up an element $b \in w$ (i.e. fill in array entries) where each array entry belongs to an appropriate subface of the various $y_j$, and where the entries of row $j$ are those of an element of $y_j$ ($j \neq i$).\MS

For notational and visual convenience we will dispense with letters and write $\DOT1$ instead of $t$ or $b_{33}$, where the subscript of ``1'' simply indicates that this is the {\em first} thing we are given to work with.\MS

Here is the partial array for the sought-for element $b$, in effect a partial element of $b$:
\[
\begin{array}{|c||c|c|c|c|c|c|c||l|}
\hline
\text{row} \downarrow    & 0 & 1 & 2 & 3 & 4 & 5 & 6 & \text{comment:}\\
\hline
0  &   &   &   &   &   &   &   &   \\
\hline
1  &   &   &   &   &   &   &   &   \\
\hline
2  &   &   &   &   &   &   &   &   \\
\hline
*3*  &   &   &   & \DOT1  &   &   &   &  1: \text{ (given)} \\
\hline
4  &   &   &   & \DOT1  &   &   &   &   \\
\hline
5  &   &   &   &   &   &   &   &   \\
\hline
6  &   &   &   &   &   &   &   &   \\
\hline
7  &   &   &   &   &   &   &   &   \\
\hline
\end{array}
\]
where the second $\DOT1$, a component of $e_4 (b)$ (row 4), is present because $e_3 e_4 (b) = e_3 e_3 (b)$. Our goal is apply the various conditions to fill in the blanks of this array with appropriate entries. Here, ``filling in'' means that the appropriate elements exist either by virtue one of the given conditions or by $(6,3)$-composition.

In particular, as we fill in the entries of the array for $b$, we will use repeatedly that $w$ is subface-simplicial (i.e. that the $b$ we are trying to build is subcomponent-simplicial)
\MS

Since $y_4 = d_4 (w)$ is \SUR3 and $\DOT1 \in d_3 (y_4)$, we know that there is an element $u_4 \in y_4$ such that $e_3 (u_4) = \DOT1$. Therefore we may fill in more entries of $b$ to get:
\[
\begin{array}{|c||c|c|c|c|c|c|c||l|}
\hline
\text{row} \downarrow    & 0 & 1 & 2 & 3 & 4 & 5 & 6 & \text{comment:} \\
\hline
0  &   &   &   & \DOT2  &   &   &   &   \\
\hline
1  &   &   &   &\DOT2   &   &   &   &   \\
\hline
2  &   &   &   & \DOT2  &   &   &   &   \\
\hline
*3*  &   &   &   & \DOT1  &   &   &   &  1: \text{ (given)} \\
\hline
4  & \DOT2  & \DOT2  &\DOT2   & \DOT1  &\DOT2   & \DOT2  &\DOT2   & 2: \text{ \SUR3}  \\
\hline
5  &   &   &   &   & \DOT2  &   &   &   \\
\hline
6  &   &   &   &   & \DOT2  &   &   &   \\
\hline
7  &   &   &   &   & \DOT2  &   &   &   \\
\hline
\end{array}
\]
Row 4 of the array represents that element $u_4 \in y_4$, and in terms of the sought-after $b \in w$
\[
u_4 = (b_{40}, b_{41}, b_{42}, t, b_{44}, b_{45}, b_{46})
\]
but for brevity we will just use the symbol ``$\DOT2$'' in each slot of row 4 as a generic abbreviation for whatever sub-element of $u_4$ appears in that slot. (It is not really necessary to maintain the subscript-book-keeping here because the array display accomplishes that). The purpose of the subscript ``2'' is just to indicate these entries occurred as step 2 of the filling-in process. The positions of the ``$\DOT2$'s'' in the other rows is based on $w$ being subface-simplicial.\MS

Next, we may fill in row 1 with an element of $y_1$ using that $y_1 = d_1(w)$ is assumed \SUR{3}, and applying the same reasoning as above to obtain:
\[
\begin{array}{|c||c|c|c|c|c|c|c||l|}
\hline
\text{row} \downarrow & 0 & 1 & 2 & 3 & 4 & 5 & 6 & \text{comment:} \\
\hline
0  & \DOT3  &   &   & \DOT2  &   &   &   &   \\
\hline
1  & \DOT3  & \DOT3  & \DOT3  & \DOT2   & \DOT3  & \DOT3  & \DOT3  &  3: \text{ \SUR3} \\
\hline
2  &   & \DOT3  &   & \DOT2  &   &   &   &   \\
\hline
*3*  &   & \DOT3  &   & \DOT1  &   &   &   &  1: \text{ (given)} \\
\hline
4  & \DOT2  & \DOT2  &\DOT2   & \DOT1  &\DOT2   & \DOT2  &\DOT2   & 2: \text{ \SUR3}  \\
\hline
5  &   & \DOT3  &   &   & \DOT2  &   &   &   \\
\hline
6  &   & \DOT3  &   &   & \DOT2  &   &   &   \\
\hline
7  &   & \DOT3  &   &   & \DOT2  &   &   &   \\
\hline
\end{array}
\]
At this stage, we must apply another condition in order to continue to fill in the entries; 2-surjectivity and 3-surjectivity do not suffice. That condition is $\SETT{1,4}$-determinacy. This choice may seem arbitrary, although it arises naturally from the task of filling in the array above. We will address the motivation behind our choice of determinacy condition below (page \pageref{explain-choice}).
\MS

With $\DPQ{1}{4}{6}$ for $d_5(w)$ we may now fill in row 5 with an element of $y_5 = d_5 (w)$ to get:
\[
\begin{array}{|c||c|c|c|c|c|c|c||l|}
\hline
\text{row} \downarrow & 0 & 1 & 2 & 3 & 4 & 5 & 6 & \text{comment:} \\
\hline
0  & \DOT3  &   &   & \DOT2  & \DOT4  &   &   &   \\
\hline
1  & \DOT3  & \DOT3  & \DOT3  & \DOT2   & \DOT3  & \DOT3  & \DOT3  &  3: \text{ \SUR3} \\
\hline
2  &   & \DOT3  &   & \DOT2  & \DOT4  &   &   &   \\
\hline
*3*  &   & \DOT3  &   & \DOT1  & \DOT4  &   &   &  1: \text{ (given)} \\
\hline
4  & \DOT2  & \DOT2  &\DOT2   & \DOT1  &\DOT2   & \DOT2  &\DOT2   & 2: \text{ \SUR3}  \\
\hline
5  & \DOT4  & \DOT3  & \DOT4  & \DOT4  & \DOT2  & \DOT4  & \DOT4  & 4: \DPQ{1}{4}{6}  \\
\hline
6  &   & \DOT3  &   &   & \DOT2  & \DOT4  &   &   \\
\hline
7  &   & \DOT3  &   &   & \DOT2  & \DOT4  &   &   \\
\hline
\end{array}
\]
Recall (lemma \refpage{subdet}) that for any $A \subsetneq [6]$, if $\SETT{1,4} \SBS A$ then $\DPQ{1}{4}{6}$ implies $A$-determinacy. Therefore we can use the assumption of $\DPQ{1}{4}{6}$ to fill in row 2 next (here $A$ is $\SETT{1,3,4}$) with the entries of an element of $y_2 = d_2(w)$: 
\[
\begin{array}{|c||c|c|c|c|c|c|c||l|}
\hline
\text{row} \downarrow    & 0 & 1 & 2 & 3 & 4 & 5 & 6 & \text{comment:} \\
\hline
0  & \DOT3  & \DOT5  &   & \DOT2  & \DOT4  &   &   &   \\
\hline
1  & \DOT3  & \DOT3  & \DOT3  & \DOT2   & \DOT3  & \DOT3  & \DOT3  &  3: \text{ \SUR3} \\
\hline
2  & \DOT5  & \DOT3  &  \DOT5  & \DOT2  & \DOT4  &  \DOT5  &  \DOT5  &  5: \DPQ{1}{4}{6}  \\
\hline
*3*  &   & \DOT3  &  \DOT5  & \DOT1  & \DOT4  &   &   &  1: \text{ (given)} \\
\hline
4  & \DOT2  & \DOT2  &\DOT2   & \DOT1  &\DOT2   & \DOT2  &\DOT2   & 2: \text{ \SUR3}  \\
\hline
5  & \DOT4  & \DOT3  & \DOT4  & \DOT4  & \DOT2  & \DOT4  & \DOT4  & 4: \DPQ{1}{4}{6}  \\
\hline
6  &   & \DOT3  &  \DOT5  &   & \DOT2  & \DOT4  &   &   \\
\hline
7  &   & \DOT3  &  \DOT5  &   & \DOT2  & \DOT4  &   &   \\
\hline
\end{array}
\]
At this point, rows 0, 6 and 7 have empty slots but have entries whose indices include 1 and 4. For example, row 6 is filled-in at slots 1, 2, 4 and 5. Therefore, each of those sub-elements of $b$ is uniquely determined from the assumption of $\DPQ{1}{4}{6}$. The element which results in row 3 is an element $a \in y_3 = d_3 (w)$, by definition of $(6,3)$-composition, and $e_3(a) = \DOT1 = t$. 

This completes verifying that $y_3$ is \SUR3 using 3-surjectivity and $\DPQ{1}{4}{6}$.\MS

\NI {\bf Notational matters:} Before continuing checking closure of various conditions, we observe that the essential information in the filling-in process just completed can be recorded unambiguously and more concisely by specifying the order in which the rows are filled in and citing what property is used for each row-filling. In the case just completed verifying closure of 3-surjectivity, the record would be summarized by the following table:
\begin{center}
\begin{tabular}{l|l}
step: row & property used:\\
\hline
1: $t \in e_3 y$ & (given)\\
2: fill row 4 & 3-surj\\
3: fill row 1 & 3-surj\\
4: fill row 5 & $\DPQ{1}{4}{6}$\\
5: fill row 2 & $\DPQ{1}{4}{6}$\\
6: fill row 0 & $\DPQ{1}{4}{6}$\\
7: fill row 6 & $\DPQ{1}{4}{6}$\\
8: fill row 7 & $\DPQ{1}{4}{6}$
\end{tabular}
\end{center}

To further streamline the notation, we will replace ``$\DOT{j}$'' with simply ``$j$'' in denoting entries in these arrays.
\BS

\NI \fbox{\bf Closure of 2-surjectivity under $(6,3)$-composition:}\MS

This is a similar calculation, and we will retain the notation above: $w$ is an arbitrary $(6,3)$-composition with $y_3=d_3(w)$. We start with an arbitrary element $t \in d_2(y_3)$ and seek an element $a \in y_3$ such that $e_2(a)=t$. To do this, we, again, use that every element of $y_3$ comes from some $b \in w$.

We begin with the array:
\[
\begin{array}{|c||c|c|c|c|c|c|c||l|}
\hline
\text{row} \downarrow    & 0 & 1 & 2 & 3 & 4 & 5 & 6 & \\
\hline
0  &   &   &   &   &   &   &   &   \\
\hline
1  &   &   &   &   &   &   &   &   \\
\hline
2  &   &   &   1   &   &   &   &   &   \\
\hline
*3*  &   &   &   1   &   &   &   &   & (1)\; \text{given}  \\
\hline
4  &   &   &   &   &   &   &   &   \\
\hline
5  &   &   &   &   &   &   &   &   \\
\hline
6  &   &   &   &   &   &   &   &   \\
\hline
7  &   &   &   &   &   &   &   &   \\
\hline
\end{array}
\]
where ``$1$'' in row 3 column 2 is the place-holder for the given $t \in d_2 (y_3)$. The table for the filling-in process here is:
\begin{center}
\begin{tabular}{l|l}
step: row & property used:\\
\hline
1: $t \in d_2y_3, \text{ row 3}$ & (given)\\
2: fill row 2 & 2-surj \\
3: fill row 5 & 2-surj\\
4: fill row 1 & $\DPQ{1}{4}{6}$\\
5: fill row 4 & $\DPQ{1}{4}{6}$ \\
6: fill row 0 & $\DPQ{1}{4}{6}$ \\
7: fill row 6 & $\DPQ{1}{4}{6}$ \\
8: fill row 7 & $\DPQ{1}{4}{6}$ \\
\end{tabular}
\end{center}
For the record, the completed element of $b \in w$ has array:
\[
\begin{array}{|c||c|c|c|c|c|c|c||l|}
\hline
\text{row} \downarrow    & 0 & 1 & 2 & 3 & 4 & 5 & 6 &\\
\hline
0  & 4  & 2  & 6  & 5  & 3  & 6  & 6  &  6: \DPQ{1}{4}{6} \\
\hline
1  & 4  & 2  & 4  & 4  & 3  & 4  & 4  &  4: \DPQ{1}{4}{6} \\
\hline
2  & 2  & 2  & 1  & 2  & 2  & 2  & 2  &  2: \text{ 2-surj} \\
\hline
*3*  & 6  & 4  & 1  & 5  & 3  & 7  & 8  &  1: \text{ (given)} \\
\hline
4  & 5  & 4  & 2  & 5  & 3  & 5  & 5  &  5: \DPQ{1}{4}{6} \\
\hline
5  & 3  & 3  & 2  & 3  & 3  & 3  & 3  &  3: \text{ 2-surj} \\
\hline
6  & 6  & 4  & 2  & 7  & 5  & 3  & 7  &  7:  \DPQ{1}{4}{6}\\
\hline
7  & 6  & 4  & 2  & 8  & 5  & 3  & 7  &  8:  \DPQ{1}{4}{6}\\
\hline
\end{array}
\]
\BS

\NI \fbox{ \bf
$w$ is $\SETT{1,4}$-determinate and $\SETT{2,5}$-determinate:
}
\MS

These properties follow directly from the Up Rule (Theorem \refpage{up rule}).
\MS

\NI \fbox{\bf All $y \in \MC{S}'_6$ are $\SETT{2,4}$-determinate:}
\MS

Given any $y \in \MC{S}'_6$, the assumed surjectivity conditions imply that $s_3(y)$ is a $(6,3)$-composition. The previous calculation showed that $s_3(y)$ is $\SETT{2,5}$-determinate. The Down Rule (Theorem \refpage{down rule}) applies, with $p=2$ and $q=5$.
\BS

\NI \fbox{\bf
Closure of $\SETT{1,4}$-determinacy under $(6,3)$-composition:
}\MS

Given $w = \COMP{6}{3}( \OM{y}{7}{3} )$ where $y_k \in \MC{S}'_6$, for each $k \neq 3$ we must verify that $d_3(w)$ is $\SETT{1,4}$-determinate. We begin with an arbitrary $\SETT{1,4}$-indexed partial element $(-,t_1,-,-,t_4,-,-)$ of $d_3(w)$. The goal is to find an element of $b_3 \in d_3(w)$ such that $e_1(b_3)=t_1$ and $e_4(b_3)=t_4$. We will do this by finding an element $b \in w$ such that $e_3(b)$ works. The initial partial array for $b$ is:
\[
\begin{array}{|c||c|c|c|c|c|c|c||c|}
\hline
\text{row} \downarrow & 0 & 1 & 2 & 3 & 4 & 5 & 6 & \\
\hline
0 & & & & & & & & \\
\hline
1 & & &1 & & & & & \\
\hline
2 & & & & & & & & \\
\hline
*3* & & 1 & & & 1 & & & (1)\; \text{given} \\
\hline
4 & & & & & & & & \\
\hline
5 & & & & 1 & & & & \\
\hline
6 & & & & & & & & \\
\hline
7 & & & & & & & & \\
\hline
\end{array}
\]
By $3$-surjectivity, there is an element of $y_5$ which fills row $5$ (step ``2'') to give:
\[
\begin{array}{|c||c|c|c|c|c|c|c||l|}
\hline
\text{row} \downarrow & 0 & 1 & 2 & 3 & 4 & 5 & 6 & \\
\hline
0 & & & & & 2 & & & \\
\hline
1 & & & 1 & & 2 & & & \\
\hline
2 & & & & & 2 & & & \\
\hline
*3* & & 1 & & & 1 & & & (1)\; \text{given} \\
\hline
4 & & & & & 2 & & & \\
\hline
5 & 2 & 2 & 2 & 1 & 2 & 2 & 2 & (2)\; 3\text{-sur} \\
\hline
6 & & & & & & 2 & & \\
\hline
7 & & & & & & 2 & & \\
\hline
\end{array}
\]
Next, by $\SETT{2,4}$-determinacy, there is an element of $y_1$ which fills row 1 (step ``3''):
\[
\begin{array}{|c||c|c|c|c|c|c|c||l|}
\hline
\text{row} \downarrow & 0 & 1 & 2 & 3 & 4 & 5 & 6 & \\
\hline
0 & 3 & & & & 2 & & & \\
\hline
1 & 3 & 3 & 1 & 3 & 2 & 3 & 3 & (3)\; \DPQ{2}{4}{6} \\
\hline
2 & & 3 & & & 2 & & & \\
\hline
*3* & & 1 & & & 1 & & & (1)\; \text{given} \\
\hline
4 & & 3 & & & 2 & & & \\
\hline
5 & 2 & 2 & 2 & 1 & 2 & 2 & 2 & (2)\; 3\text{-sur} \\
\hline
6 & & 3 & & & & 2 & & \\
\hline
7 & & 3 & & & & 2 & & \\
\hline
\end{array}
\]
Use $\SETT{1,4}$-determinacy again to fill row 4 (step ``4''):
\[
\begin{array}{|c||c|c|c|c|c|c|c||l|}
\hline
\text{row} \downarrow & 0 & 1 & 2 & 3 & 4 & 5 & 6 & \\
\hline
0 & 3 & & & 4 & 2 & & & \\
\hline
1 & 3 & 3 & 1 & 3 & 2 & 3 & 3 & (3)\; \DPQ{2}{4}{6} \\
\hline
2 & & 3 & & 4 & 2 & & & \\
\hline
*3* & & 1 & & 4 & 1 & & & (1)\; \text{given} \\
\hline
4 & 4 & 3 & 4 & 4 & 2 & 4 & 4 & (4)\; \DPQ{1}{4}{6} \\
\hline
5 & 2 & 2 & 2 & 1 & 2 & 2 & 2 &(2)\; 3\text{-sur} \\
\hline
6 & & 3 & & & 4 & 2 & & \\
\hline
7 & & 3 & & & 4 & 2 & & \\
\hline
\end{array}
\]
At this stage, rows 2, 6 and 7 may be filled by $\SETT{1,4}$-determinacy, after which row 0 can be filled by the same condition, and row 3 is determined by $(6,3)$-composition. The completed array is:
\[
\begin{array}{|c||c|c|c|c|c|c|c||l|}
\hline
\text{row} \downarrow & 0 & 1 & 2 & 3 & 4 & 5 & 6 & \\
\hline
0 & 3 & 5 & 8 & 4 & 2 & 6 & 7 & (8)\;  \DPQ{1}{4}{6}\\
\hline
1 & 3 & 3 & 1 & 3 & 2 & 3 & 3 &(3)\;  \DPQ{2}{4}{6} \\
\hline
2 & 5 & 3 & 5 & 4 & 2 & 5 & 5 &(5)\; \DPQ{1}{4}{6} \\
\hline
*3* & 8 & 1 & 5 & 4 & 1 & 6 & 7 &(1)\; \text{given} \\
\hline
4 & 4 & 3 & 4 & 4 & 2 & 4 & 4 &(4)\;  \DPQ{1}{4}{6} \\
\hline
5 & 2 & 2 & 2 & 1 & 2 & 2 & 2 & (2)\; 3\text{-sur} \\
\hline
6 & 6 & 3 & 5 & 6 & 4 & 2 & 6 & (6) \DPQ{1}{4}{6} \\
\hline
7 & 7 & 3 & 5 & 7 & 4 & 2 & 6 & (7)\DPQ{1}{4}{6} \\
\hline
\end{array}
\]
which represents an element $b \in w$ such that $b_3 = e_3(b) \in d_3(w)$ has the claimed property that $e_1(b_3)$ and $e_4(b_3)$ are the originally given components $t_1$ and $t_4$.
\BS

Now we can specify a {\bf \NICOMP{6}{3}} as follows:\MS

We start with the 6-truncated complex $\MC{S}'$ as above, and define $S'_7$ to consist of all $(6,3)$-compositions involving $y \in S'_6$. The closure calculations we carried out show that $S'_7$ includes all the degeneracies $s_k(y)$, and that $d_3: S'_7 \to S'_6$ is well-defined.\MS

In fact, we can continue this process, with essentially the same verifications, to define $S'_8$ to consist of all $(7,3)$-compositions involving simplices of $\MC{S}'_7$, define $\MC{S}'_9$ to consist of all $(8,3)$-compositions involving simplices of $\MC{S}'_8$, etc.
\MS

Rather than work out the details for this example, we will develop the general construction below.
\BS

\label{explain-choice}
\NI {\bf Comment on the choice of $\{1,4\}$-determinacy:} As noted above it may have seemed arbitrary to have required $\SETT{1,4}$-determinacy in the definition of the \NICOMP{6}{3} $\MC{S}'$. It {\em was} arbitrary in the sense that other choices of a determinacy condition also work. On the other hand, this  particular one works for general combinatorial reasons, and can be applied in {\em all} higher-dimensional cases $(n,i)$ where $n>2$ and $1<i<n$. The cases $i=0,1,n,n+1$ are similar in style but have to be treated separately. In the next six sections, we will prove the corresponding theorems.
\MS

{\bf In each of those sections,} the work will be split into two theorems, the first of which establishes facts about properties of $(n,i)$-compositions under various assumptions, and the second of which specifies the model.

\subsection{Preparatory notes on the model theorems}

Each model theorem in the subsections to follow assumes $n \geq 2$ and $i \in [n+1]$. The starting point in each of the theorems is a given $n$-truncated subcomplex $\MC{S}'_{\leq n}$ of \SBSI\ simplices such that, dependent on $i$, certain surjectivity and determinacy conditions are assumed in dimension $n$. These conditions imply surjectivity and determinacy conditions in lower dimensions but those won't be referred to in the proofs.
\MS

The existence of such truncated complexes will be examined in Section \refpage{model-construction}.
\MS

Each proof will define sets $\MC{S}'_m \SBS \MC{S}_m$ inductively for each $m \geq n+1$ where $z \in \MC{S}'_m$ iff
\[
z = \COMP{m-1}{i}(\OM{u}{m}{i})
\]
with $u_j \in \MC{S}'_{m-1}$ for all $j \neq i$. To complete the reasoning, that is, to show that $\MC{S}'_m \ISO \BOX{i}{m}{\MC{S}'}$, we will have to show that for each $m \geq n+1$, that $\MC{S}'_{\leq m}$ is actually an $m$-truncated complex. This reduces to showing two things: (1) That for all $w \in \MC{S}'_m$, $d_i(w) \in \MC{S}'_{m-1}$. \; (2) That for all $y \in \MC{S}'_{m-1}$ and all $k \in [m-1]$, then $s_k(y) \in \MC{S}'_m$. Verifying item (1) involves some calculations with partial arrays, and the essential details will be treated is separate theorems which will be referred to as ``parameter'' theorems.

\subsection{\NICOMP{n}{i}s for $1<i<n$}
\label{solution-n-i}

\begin{newthm} 
{\bf ($(n,i)$-parameter theorem)}
\label{n-i-parameter-theorem}
\index{Theorem!$(n,i)$-parameter theorem}
{\rm
\SL

Fix $n$ and $i$ such that $1<i<n$. Suppose $\MC{S}'_{\leq n}$ is an $n$-truncated subcomplex of $\mathcal{S}$ such that all $y \in \MC{S}'_n$ are subface-simplicial, \SUR{i}, \SUR{(i-1)} and $\DPQW{i-2}{i+1}$.
\MS

\NI Let $w = \text{comp}_{n,i}(\OM{y}{n+1}{i})$ where $y_j \in \MC{S}'_n$ for all $j \neq i$.  
Then:
\begin{enumerate}
\item $w$ is $\DPQW{i-2}{i+1}$.
\item $w$ is $\DPQW{i-1}{i+2}$. 
\item All $y \in \MC{S}'_n$ are $\DPQW{i-1}{i+1}$. 
\item $d_i(w) \in \MC{S}'_n$ \label{d_i(w)}
\item $w$ is \SUR{i} and \SUR{(i-1)}.
\end{enumerate}
Note: The parameters in this theorem are $n$ and $i$.
}
\end{newthm}

\MS

\Proof

\NI {\bf Proofs of 1 and 2:}

That $w$ is \DET{\{i-2,i+1\}} and \DET{\{i-1,i+2\}} follow from the Up Rule, Theorem \refpage{up rule}.
\MS

\NI {\bf Proof of 3 (that all $y \in \MC{S}'_n$ are $\DPQW{i-1}{i+1}$):}

By the Degenerate $(n,i)$-compositions theorem (page \pageref{degen-comps-thm}), the hypotheses on $\MC{S}'_n$ imply that for all $y \in \MC{S}'_n$, it is the case that for each $k \in [n]$, $s_k(y) \in \mathcal{S}_{n+1}$ is an $(n,i)$-composition. In particular, for each $y \in \MC{S}'_n$, $s_i(y)$ is an $(n,i)$-composition. Therefore, by part 2, $s_i(y)$ is $\DPQ{i-1}{i+2}{{(n+1)}}$.

By The Down Rule, Theorem \refpage{down rule} applied to $p=i-1, q=i+2$ (so $q-p>2$), $s_i(y)$ being $\DPQ{i-1}{i+2}{(n+1)}$ implies $y$ is $\DPQ{i-1}{i+1}{n}$.
\MS

\NI {\bf Proof of 4 (that $d_i(w) \in \MC{S}'_n$):}

Suppose $y_i = d_i(w)$, given $w = \text{comp}_{n,i}( \OM{y}{n+1}{i})$. We must verify that $y_i$ is subface-simplicial, $i$-surjective, $(i-1)$-surjective and $\DPQ{i-2}{i+1}{n}$.
\MS

\NI {\bf That $y_i$ is subface-simplicial:} The e-simplicial composites Lemma \refpage{simpclosure} implies that $y_i$ is e-simplicial. Since all the faces of $y_i$ are faces of the other $y_j$ (which are all subface-simplicial) then $y_i$ is subface-simplicial also.
\MS

\NI {\bf That $y_i$ is \SUR{i}}: Starting with an arbitrary element $t  \in d_i(y_i) = d_i d_i(w)$ we must find an element $b_i \in y_i$ such that $e_i(b_i)=t$. Since every element of $y_i = d_i(w)$ arises as $e_i(b)$ for some $b \in w$, it suffices to find $b \in w$ such that $e_i(e_i(b)) = t$. The following array calculation shows how to do that. 

We start with the element $t \in d_i(y_i)$ in row $i$ column $i$ labelled ``$1$'' in the array below, then fill in row $i+1$ by $i$-surjectivity. 
\[
\begin{array}{|c||c|c|c|c|c|c|c|c|c||l|}
\hline
\text{row} \downarrow & 0 & \cdots & i-3 & i-2 & i-1 & i & i+1 & \cdots & n & \\
\hline
0 & & & & & & 2 & & & & \\
\hline
\vdots & & & & & & \vdots & & & & \\
\hline
i-2 & & & & & & 2 & & & & \\
\hline
i-1 & & & & & & 2 & & & & \\
\hline
*\;i\;* & & & & & & 1 & & & & (1)\; t\in d_i(y_i) \text{, given}\\
\hline
i+1 & 2 & \cdots & 2 & 2 & 2 & 1 & 2 & \cdots &2 & (2)\; $i$\text{-sur} \\
\hline
\vdots & & & & & & & \vdots & & & \\
\hline
n+1 & & & & & & & 2 & & & \\
\hline
\end{array}
\]
Next, we fill in row $i-2$ (step ``$3$'') by $i$-surjectivity. Now the array looks like this:
\[
\begin{array}{|c||c|c|c|c|c|c|c|c|c||l|}
\hline
\text{row} \downarrow & 0 & \cdots & i-3 & i-2 & i-1 & i & i+1 & \cdots & n & \\
\hline
0 & & & 3 & 4 & & 2 & & & & \\
\hline
\vdots & & & \vdots &\vdots & & \vdots & & & & \\
\hline
i-2 & 3 & \cdots & 3 & 3 & 3 & 2 & 3 & \cdots & 3 & (3)\; $i$\text{-sur} \\
\hline
i-1 &  &  &  & 3 &  & 2 & & & & \\
\hline
*\;i\;* & & & & 3 &  & 1 & & & & (1)\; t\in d_i(y_i) \text{, given}\\
\hline
i+1 & 2 & \cdots & 2 & 2 & 2 & 1 & 2 & \cdots &2 & (2)\; $i$\text{-sur} \\
\hline
\vdots & & & & \vdots & & & \vdots & & & \\
\hline
n+1 & & & & 3 &  & & 2 & & & \\
\hline
\end{array}
\]
Now, rows $i+2$ through $n+1$ can be filled using $\DPQ{i-2}{i+1}{n}$. The entries in row $i+2$ fill column $i+1$ in rows $0$ through $i$. Now the array is:
\[
\begin{array}{|c||c|c|c|c|c|c|c|c|c||l|}
\hline
\text{row} \downarrow & 0 & \cdots & i-3 & i-2 & i-1 & i & i+1 & \cdots & n & \\
\hline
0 & & & 3 &  & & 2 & / & \cdots & / & \\
\hline
\vdots & & & \vdots & & & \vdots & / & \cdots & / & \\
\hline
i-2 & 3 & \cdots & 3 & 3 & 3 & 2 & 3 & \cdots & 3 & (3)\; $i$\text{-sur} \\
\hline
i-1 &  &  &  & 3 &  & 2 & / & \cdots & / & \\
\hline
*\;i\;* & & & & 3 &  & 1 & / & \cdots & / & (1)\; t\in d_i(y_i) \text{, given}\\
\hline
i+1 & 2 & \cdots & 2 & 2 & 2 & 1 & 2 & \cdots &2 & (2)\; $i$\text{-sur} \\
\hline
\vdots & / & / & / & \vdots & / & / & \vdots & \cdots & / & \T\DPQ{i-2}{i+1}{n} \\
\hline
n+1 & / & / & / & 3 & / & / & 2 & / & / & \T\DPQ{i-2}{i+1}{n}\\
\hline
\end{array}
\]
Row $i-1$ can now be filled (step ``$4$'') using $\DPQ{i-2}{i+1}{n}$, and those entries fill column $i-2$ in rows $0$ through $i-3$. The result is:
\[
\begin{array}{|c||c|c|c|c|c|c|c|c|c||l|}
\hline
\text{row} \downarrow & 0 & \cdots & i-3 & i-2 & i-1 & i & i+1 & \cdots & n & \\
\hline
0 & & & 3 & 4  & & 2 & / & \cdots & / & \\
\hline
\vdots & & & \vdots & & & \vdots & / & \cdots & / & \\
\hline
i-2 & 3 & \cdots & 3 & 3 & 3 & 2 & 3 & \cdots & 3 & (3)\; $i$\text{-sur} \\
\hline
i-1 & 4  & \cdots  & 4  & 3 & 4  & 2 & / & \cdots & / & \\
\hline
*\;i\;* & & & & 3 & 4  & 1 & / & \cdots & / & (1)\; t\in d_i(y_i) \text{, given}\\
\hline
i+1 & 2 & \cdots & 2 & 2 & 2 & 1 & 2 & \cdots &2 & (2)\;  $i$\text{-sur} \\
\hline
\vdots & / & / & / & \vdots & / & / & \vdots & \cdots & / & \T\DPQ{i-2}{i+1}{n} \\
\hline
n+1 & / & / & / & 3 & / & / & 2 & / & / & \T\DPQ{i-2}{i+1}{n}\\
\hline
\end{array}
\]
Now the array can be completed using $\DPQ{i-2}{i+1}{n}$ for filling rows $0$ through $i-3$, and using $(n,i)$-composition to fill row $i$. We have therefore shown that an element $b \in w$ exists such that $e_i(b) = b_i$ and $e_i(b_i) = t$, as required.
\MS

\NI {\bf That $y_i$ is $(i-1)$-surjective:}

The approach is similar to $i$-surjectivity. Given an arbitrary $t \in d_{i-1}(y_i)$, we fill an array to obtain $b \in w$ such that $e_i(b) \in  y_i$ and $e_{i-1}(e_i(b)) = t$.

The initial steps are to apply $(i-1)$-surjectivity in rows $i-1$ and $i+2$. The partial array is:
\MS

\[
\begin{array}{|c||c|c|c|c|c|c|c|c|c|c||l|}
\hline
\text{row} \downarrow & 0 & \cdots & i-3 & i-2 & i-1 & i & i+1 & i+2 & \cdots & n & \\
\hline
0 & & & & 2 & & & 3 & & & & \\
\hline
\vdots & & & &\vdots & & &\vdots & & & & \\
\hline
i-2 & & & & 2 & & & 3 & & & & \\
\hline
i-1 & 2 & \cdots & 2 & 2 & 1 & 2 & 2 & 2 & \cdots & 2 & (2)\; (i-1)\text{-sur}\\
\hline
*\;i\;* & & & & & 1 & & 3 & & & & (1)\; t \in d_{i-1}(y_i) \text{, given}\\
\hline
i+1 & & & & & 2 & & 3 & & & & \\
\hline
i+2 & 3 & \cdots & 3 & 3 & 2 & 3 & 3 & 3 & \cdots & 3 & (3)\; (i-1)\text{-sur}\\
\hline
\vdots & & & & & \vdots & & &\vdots & & & \\
\hline
n+1 & & & & & 2 & & & 3 & & & \\
\hline
\end{array}
\]
Next, use $\DPQ{i-2}{i+1}{n}$ to fill rows $0$ through $i-2$; this also fills column $i-2$ in rows $j =i , \cdots, n+1$. Then one may complete row $i+1$ (as step ``$4$'') to give:
\MS

\[
\begin{array}{|c||c|c|c|c|c|c|c|c|c|c||l|}
\hline
\text{row}  \downarrow & 0 & \cdots & i-3 & i-2 & i-1 & i & i+1 & i+2 & \cdots & n & \\
\hline
0 & / & \cdots & / & 2 & / & / & 3 & / & \cdots & / & \\
\hline
\vdots & / & \cdots & / &\vdots & / & / &\vdots & / & \cdots & / & \\
\hline
i-2 & / & / & / & 2 & / & / & 3 & / & \cdots & / & \\
\hline
i-1 & 2 & \cdots & 2 & 2 & 1 & 2 & 2 & 2 & \cdots & 2 & (2)\; (i-1)\text{-sur}\\
\hline
*\;i\;* & / & / & / & / & 1 & 4 & 3 & & & & (1)\; t \in d_{i-1}(y_i) \text{, given}\\
\hline
i+1 & / & / & / & / & 2 & 4 & 3 & 4 & \cdots & 4 & (4)\; \T\DPQ{i-2}{i+1}{n} \\
\hline
i+2 & 3 & \cdots & 3 & 3 & 2 & 3 & 3 & 3 & \cdots & 3 & (3)\; (i-1)\text{-sur}\\
\hline
\vdots & / & / & / & / & \vdots & & &\vdots & & & \\
\hline
n+1 & / & / & / & / & 2 & & 4 & 3 & & & \\
\hline
\end{array}
\]
Finally, rows $i+3$ through $n+1$ can be filled using $\DPQ{i-2}{i+1}{n}$ and row $i$ is determined by $(n,i)$-composition. This array represents an element $b \in w$ such that $e_i(b) \in y_i$ and $e_{i-1}(e_i(b)) = t$, as required.
\MS

\NI {\bf That $y_i$ is $\DPQ{i-2}{i+1}{n}$:}

Suppose $( \cdots, t_{i-2}, - , - , t_{i+1}, \cdots )$ is an $\SETT{i-2,i+1}$-indexed partial element of $y_i$. We must show there exists $t \in y_i$ (necessarily unique by the 2-lemma) such that $e_{i-2}(t) = t_{i-2}$ and $e_{i+1}(t) = t_{i+1}$. The partial element entries of $y_i$ form entries in the array for a partial element of $w$. As with the other steps, we will use the hypotheses to prove that an element $b \in w$ exists such that $e_i(b)$ is the required $t$. The elements $t_{i-2}$ and $t_{i+1}$ are displayed as ``$1$'' in the array below, and we use $i$-surjectivity to fill row $i+2$ (step ``$2$'') to get:
\MS

\[
\begin{array}{|c||c|c|c|c|c|c|c|c|c|c||l|}
\hline
\text{row} \downarrow & 0 & \cdots & i-3 & i-2 & i-1 & i & i+1 & i+2 & \cdots & n & \\
\hline
0 & & & & & & & 2 & & & & \\
\hline
\vdots & & & & & & & \vdots & & & & \\
\hline
i-2 & & & & & 1 & & 2 & & & & \\
\hline
i-1 & & & & & & & 2 & & & & \\
\hline
*\;i\;* & & & & 1 & & & 1 & & & & (1) \text{ given} \\
\hline
i+1 & & & & & & & 2 & & & & \\
\hline
i+2 & 2 & \cdots & 2 & 2 & 2 & 1 & 2 & 2 & \cdots & 2 & (2)\;  i\text{-sur} \\
\hline
\vdots & & & & & & & & \vdots & & & \\
\hline
n+1 & & & & & & & & 2 & & & \\
\hline
\end{array}
\]
We may now fill in row $i-2$ (step ``$3$'') using $\DPQ{i-1}{i+1}{n}$, which we established earlier in the proof. Doing so determines entries in columns $i-2$ and $i+1$ of row $i+1$. Therefore, we may fill row $i+1$ (step ``$4$'') using $\DPQ{i-2}{i+1}{n}$. At this stage the array is:
\MS

\[
\begin{array}{|c||c|c|c|c|c|c|c|c|c|c||l|}
\hline
\text{row} \downarrow & 0 & \cdots & i-3 & i-2 & i-1 & i & i+1 & i+2 & \cdots & n & \\
\hline
0 & & & 3 & & & 4 & 2 & & & & \\
\hline
\vdots & & &\vdots & & &\vdots & \vdots & & & & \\
\hline
i-2 & 3 & \cdots & 3 & 3 & 1 & 3 & 2 & 3 & \cdots & 3 & (3)\; \T\DPQ{i-1}{i+1}{n} \\
\hline
i-1 & & & & 3 & & 4 & 2 & & & & \\
\hline
*\;i\;* & & & & 1 & & 4 & 1 & & & & (1) \text{ given}\\
\hline
i+1 & 4 & \cdots & 4 & 3 & 4 & 4 & 2 & 4 & \cdots & 4 & (4)\; \T\DPQ{i-2}{i+1}{n} \\
\hline
i+2 & 2 & \cdots & 2 & 2 & 2 & 1 & 2 & 2 & \cdots & 2 & (2)\; i\text{-sur} \\
\hline
\vdots & & & & \vdots & & & \vdots & \vdots & & & \\
\hline
n+1 & & & & 3 & & & 4 & 2 & & & \\
\hline
\end{array}
\]
Row $i-1$ now has entries in columns $i-2$ and $i+1$ and therefore may be filled (step ``$5$'') by $\DPQ{i-2}{i+1}{n}$ to give:
\MS

\[
\begin{array}{|c||c|c|c|c|c|c|c|c|c|c||l|}
\hline
\text{row} \downarrow & 0 & \cdots & i-3 & i-2 & i-1 & i & i+1 & i+2 & \cdots & n & \\
\hline
0 & & & 3 & 5 & & 4 & 2 & & & & \\
\hline
\vdots & & &\vdots &\vdots & &\vdots & \vdots & & & & \\
\hline
i-2 & 3 & \cdots & 3 & 3 & 1 & 3 & 2 & 3 & \cdots & 3 & (3)\; \T\DPQ{i-1}{i+1}{n} \\
\hline
i-1 & 5 & \cdots & 5 & 3 & 5 & 4 & 2 & 5 & \cdots & 5 & (5)\; \T\DPQ{i-2}{i+1}{n} \\
\hline
*\;i\;* & & & & 1 & 5 & 4 & 1 & & & &(1) \text{ given} \\
\hline
i+1 & 4 & \cdots & 4 & 3 & 4 & 4 & 2 & 4 & \cdots & 4 & (4)\; \T\DPQ{i-2}{i+1}{n} \\
\hline
i+2 & 2 & \cdots & 2 & 2 & 2 & 1 & 2 & 2 & \cdots & 2 & (2)\; i\text{-sur} \\
\hline
\vdots & & & & \vdots & \vdots & & \vdots & \vdots & & & \\
\hline
n+1 & & & & 3 & 5 & & 4 & 2 & & & \\
\hline
\end{array}
\]
We may now complete the array using $\DPQ{i-2}{i+1}{n}$ for rows $0$ through $i-3$ and rows $i+3$ to $n+1$; row $i$ is then determined by $(n,i)$-composition. This array represents an element $b \in w$ such that $e_i(b)= t$, and $e_{i-2}(t) = t_{i-2}$ and $e_{i+1}(t) = t_{i+1}$ as required.
\MS

\NI {\bf Proof of 5 (that $w$ is \SUR{i} and \SUR{(i-1)}):}

$w$ is \SUR{i} by definition. The proof that $w$ is \SUR{(i-1)} is by the following array argument. Let $w = \COMP{n}{i}( \OM{y}{{n+1}}{i})$ and suppose an arbitrary element of $b_{i-1} \in y_{i-1}$ is given. That is shown as row $i-1$ in the following array:
\[
\begin{array}{|c||c|c|c|c|c|c|c|c|c|c||l|}
\hline
\text{row} \downarrow & 0 & \cdots & i-3 & i-2 & i-1 & i & i+1 & i+2 & \cdots & n & \\
\hline 
0 & & & & 1 & & & & & & & \\
\hline 
\vdots & & & & \vdots & & & & & & & \\
\hline 
i-2 & & & & 1 & & & & & & & \\
\hline 
i-1 & 1 & \cdots & 1 & 1 & 1 & 1 & 1 & 1 & \cdots & 1 & (1)\; b_{i-1}\\
\hline 
*\;i\;* & & & & & 1 & & & & & & \\
\hline 
i+1 & & & & & 1 & & & & & & \\
\hline 
i+2 & & & & & 1 & & & & & & \\
\hline 
\vdots & & & & & \vdots & & & & & & \\
\hline 
n+1 & & & & & 1 & & & & & & \\
\hline 
\end{array}
\]
Next, use the $(i-1)$-surjectivity of $y_{i+2}$ to fill row $i+2$:
\[
\begin{array}{|c||c|c|c|c|c|c|c|c|c|c||l|}
\hline
\text{row} \downarrow & 0 & \cdots & i-3 & i-2 & i-1 & i & i+1 & i+2 & \cdots & n & \\
\hline 
0 & & & & 1 & & & 2 & & & & \\
\hline 
\vdots & & & & \vdots & & & \vdots & & & & \\
\hline 
i-2 & & & & 1 & & & 2 & & & & \\
\hline 
i-1 & 1 & \cdots & 1 & 1 & 1 & 1 & 1 & 1 & \cdots & 1 & (1)\;b_{i-1}\\
\hline 
*\;i\;* & & & & & 1 & & 2 & & & & \\
\hline 
i+1 & & & & & 1 & & 2 & & & & \\
\hline 
i+2 & 2 & 2 & 2 & 2 & 1 & 2 & 2 & 2 & 2 & 2 & (2)\; \T{(i-1)\text{-sur}} \\
\hline 
\vdots & & & & & \vdots & & & \vdots & & & \\
\hline 
n+1 & & & & & 1 & & & 2 & & & \\
\hline 
\end{array}
\]
Next, rows $0$ through $i-2$ may be filled using $\DPQ{i-2}{i+1}{n}$:
\[
\begin{array}{|c||c|c|c|c|c|c|c|c|c|c||l|}
\hline
\text{row} \downarrow & 0 & \cdots & i-3 & i-2 & i-1 & i & i+1 & i+2 & \cdots & n & \\
\hline 
0 & / & \cdots & / & 1 & / & / & 2 & / & \cdots & / & \T\DPQ{i-2}{i+1}{n} \\
\hline 
\vdots & / & \cdots & / & \vdots & / & / & \vdots & / & \cdots & / &  \T\DPQ{i-2}{i+1}{n}\\
\hline 
i-2 & / & \cdots & / & 1 & / & / & 2 & / & \cdots & / &  \T\DPQ{i-2}{i+1}{n}\\
\hline 
i-1 & 1 & \cdots & 1 & 1 & 1 & 1 & 1 & 1 & \cdots & 1 & (1)\; b_{i-1}\\
\hline 
*\;i\;* & & & & / & 1 & / & 2 & / & \cdots & / & \\
\hline 
i+1 & & & & / & 1 & / & 2 & / & \cdots & / & \\
\hline 
i+2 & 2 & 2 & 2 & 2 & 1 & 2 & 2 & 2 & 2 & 2 & (2)\; \T{(i-1)\text{-sur}} \\
\hline 
\vdots & & & & / & \vdots & / & / & \vdots & \cdots & / & \\
\hline 
n+1 & & & & / & 1 & / & / & 2 & \cdots & / & \\
\hline 
\end{array}
\]
Now the array can be completed in rows $i+1$ through $n+1$ by $\DPQ{i-2}{i+1}{n}$, and row $i$ by $(n,i)$-composition. This shows that $w$ is \SUR{(i-1)}.

\qed
\MS

\NI {\bf Note:} Similar reasoning to the proof that $w$ is \SUR{(i-1)} shows that $w$ is also \SUR{(i+1)}.

\begin{newthm}
\label{composer-model;1<i<n}
\index{Theorem!Existence of \NICOMP{n}{i}s}
{\rm
{\bf (Existence of \NICOMP{n}{i}s when $1<i<n$)}

Fix $n$ and $i$ where $1<i<n$. Suppose $\MC{S}'_{\leq n}$ is an $n$-truncated subcomplex of $\mathcal{S}$ such that all $y \in \MC{S}'_n$ are subface-simplicial, \SUR{i}, \SUR{(i-1)} and $\DPQW{i-2}{i+1}$. Then there is a subcomplex $\MC{S}' \SBS \MC{S}$ of $S$ such that $\TR^n(\MC{S}')$ is the given one and which is an \NICOMP{n}{i}.
}
\end{newthm}

\Proof

We will define $\MC{S}'_m$ for $m \geq n+1$ by induction on $m$. At each stage we will have to show that $\MC{S}'_{\leq m}$ is an $m$-truncated subcomplex of $\MC{S}$. This will yield a subcomplex $\MC{S}' \SBS \MC{S}$ such that $\TR^n(\MC{S})$ is the given $n$-truncated complex. By definition,  $\MC{S'}_m \ISO \BOX{m}{i}{\MC{S}'}$ for all $m \geq n+1$.
\MS

We begin by defining $\MC{S}'_{n+1}$ to consist of all $(n,i)$-compositions 
\[
w = \COMP{n}{i}(\OM{y}{n+1}{i})
\]
such that for each $j \neq i$, $y_j \in \MC{S}'_n$. In order to show that $\MC{S}'_{n+1}$ extends $\MC{S}'_{\leq n}$ to dimension $n+1$ we have to show that $d_i(w) = \MC{S}'_n$ and that for all $y \in \MC{S}'_n$ and all $k \in [n]$, $s_k(y) \in \MC{S}'_{n+1}$.
\MS

That $d_i(w) \in \MC{S}'_n$ was proved in item \ref{d_i(w)} of Theorem \refpage{n-i-parameter-theorem}. 
\MS

Let $y \in \MC{S}'_n$ and $k \in [n]$. Since $y$ is assumed \SUR{i} and \SUR{(i-1)} and since $1<i<n$ then the Degenerate Compositions Theorem \refpage{degen-comps-thm} implies that $s_k(y)$ is an $(n,i)$-composition. That the faces of $s_k(y)$ belong to $\MC{S}'_n$ implies that $s_k(y) \in \MC{S}'_{n+1}$.
\MS

Therefore $\MC{S}'_{\leq n+1}$ is an $(n+1)$-truncated extension of $\MC{S}'_{\leq n}$.
\MS

Now do induction on $m \geq n+2$. Define $\MC{S}'_m$ to consist of all $(m-1,i)$-compositions
\[
z = \COMP{m-1}{i}(\OM{u}{m}{i})
\]
where for each $j \neq i$, $u_j \in \MC{S}'_{m-1}$. The induction hypothesis is that each $u_j$ is subface-simplicial, \SUR{i} and \SUR{(i-1)}, and therefore Theorem \refpage{n-i-parameter-theorem} applies. By exactly the same reasoning as in the case $m=n+1$ above, it follows that $d_i(z) \in \MC{S}'_{m-1}$ and for all $k \in [m-1]$ and all $u \in \MC{S}'_{m-1}$, $s_k(u) \in \MC{S}'_m$.

We therefore obtain a subcomplex $\MC{S}' \SBS \MC{S}$ such that for all $m \geq n+1$, $\MC{S}'_m \ISO \BOX{m}{i}{\MC{S}'}$ and $\TR^n(\MC{S}
)$ is the given subcomplex.

\qed

\subsection{\NICOMP{n}{0}s}
\label{solution-n-0}

\begin{newthm}
{\bf ($(n,0)$-parameter theorem)}
\label{n-0-parameter-theorem}
\index{Theorem!$(n,0)$-parameter theorem}
{\rm
\SL

Fix $n \geq 2$. Suppose $w = \COMP{n}{0}(-,y_1, \cdots, y_{n+1})$ is an $(n,0)$-composition and that for each $ k \in [n+1]-\SETT{0}$, $y_k$ is \SBSI, \SUR{0}, \SUR{1} and $\DPQW{1}{2}$. Then:
\begin{enumerate}
\item $d_0(w)$ is \SUR{0}.
\item $d_0(w)$ is \SUR{1}.
\item $d_0(w)$ is $\DPQW{1}{2}$.
\item $w$ is \SUR{1}.
\item $w$ is $\DPQ{1}{2}{(n+1)}$.
\end{enumerate}
 Note: The ``parameter'' in this theorem is $n$.
}
\end{newthm}

\Proof

\begin{enumerate}
\item 
To prove $d_0(w)$ is \SUR{0}: Let $t \in d_0d_0(w)$ and consider the partial element of $w$ with array where the ``$1$'' is $t$:
\[
\begin{array}{|c||c|c|c|c|c|c||l|}
\hline
\text{row} \downarrow & 0 & 1 & 2 & 3 & \cdots & n & \\
\hline
*0* & 1 & & & & & & (1) \text{ given}\\
\hline
1 & 1 & & & & & & \\
\hline
2 & & & & & & & \\
\hline
3  & & & & & & & \\
\hline
\vdots & & & & & & & \\
\hline
n+1 & & & & & & & \\
\hline 
\end{array}
\]
By the 1-surjectivity of $y_1$, there is an element of $y_1$ which fills row 1 to give (as step 2):
\[
\begin{array}{|c||c|c|c|c|c|c||l|}
\hline
\text{row} \downarrow & 0 & 1 & 2 & 3 & \cdots & n & \\
\hline
*0* & 1 & & & & & & (1) \text{ given}\\
\hline
1 & 1 & 2 & 2 &  2 & \cdots  & 2 & (2)\; 1\text{-sur} \\
\hline
2 & & 2 & & & & & \\
\hline
3  & & 2 & & & & & \\
\hline
\vdots & & \vdots & & & & & \\
\hline
n+1 & & 2 & & & & & \\
\hline 
\end{array}
\]
Next, by 1-surjectivity of $y_2$, fill row 2 (step 3):
\[
\begin{array}{|c||c|c|c|c|c|c||l|}
\hline
\text{row} \downarrow & 0 & 1 & 2 & 3 & \cdots & n & \\
\hline
*0* & 1 & 3 & & & & & (1) \text{ given}\\
\hline
1 & 1 & 2 & 2 &  2 & \cdots  & 2 & (2)\; 1\text{-sur} \\
\hline
2 & 3 & 2 & 3 & 3 & \cdots & 3 & (3)\; 1\text{-sur} \\
\hline
3  & & 2 & 3 & & &  &  \\
\hline
\vdots & & \vdots & \vdots & & & & \\
\hline
n+1 & & 2 & 3 & & & & \\
\hline 
\end{array}
\]
Using that $y_3, \cdots, y_{n+1}$ are \DET{\SETT{1,2}}, there are elements of those $n$-simplices filling rows $3$ through $n+1$. $(n,0)$-composition then determines row 0, and therefore there is a element $b \in w$ such that $e_0(b) \in d_0(w)$ has $e_0(e_0(b)) = t$, as required.

\item 
To prove $d_0(w)$ is \SUR{1}: Let $t \in d_1(d_0(w))$ and consider the partial element array of $w$, where the ``$1$'' is $t$:
\[
\begin{array}{|c||c|c|c|c|c|c||l|}
\hline
\text{row} \downarrow & 0 & 1 & 2 & 3 & \cdots & n & \\
\hline \hline
*0* & & 1 & & & & & (1)\; \text{given}\\
\hline
1 & & & & & & & \\
\hline
2 & 1 & & & & & & \\
\hline
3 & & & & & & & \\
\hline
\vdots & & & & & & & \\
\hline
n+1 & & & & & & & \\
\hline 
\end{array}
\]
Then there is an element of $y_2$ (by 0-surjectivity) which fills row 2 (step 2) and then, as step 3, there is a element of $y_1$ (by 1-surjectivity) which fills row 1. The resulting array now is:
\[
\begin{array}{|c||c|c|c|c|c|c||l|}
\hline
\text{row} \downarrow & 0 & 1 & 2 & 3 & \cdots & n & \\
\hline \hline
*0* & 3 & 1 & & & & & (1)\; \text{given}\\
\hline
1 & 3 & 2 & 3 & 3 & \cdots & 3 &(3)\; 1\text{-sur} \\
\hline
2 & 1 & 2 & 2 & 2 & \cdots & 2 & (2)\; 0\text{-sur} \\
\hline
3 & & 3 & 2 & & & & \\
\hline
\vdots & & \vdots & \vdots & & & & \\
\hline
n+1 & & 3 & 2 & & & & \\
\hline 
\end{array}
\]
Using that $y_3, \cdots, y_{n+1}$ are \DET{\SETT{1,2}}, there are elements of those $n$-simplices which fill rows 3 through $n+1$. Row 0 is filled by $(n,0)$-composition and therefore there is an element $b \in w$ such that $e_1(e_0(b)) = t$ as required.

\item 
To prove $d_0(w)$ is \DET{\SETT{1,2}}: Start with any $\SETT{1,2}$-indexed partial element $\SETT{t_1,t_2}$ of $d_0(w)$ and consider the corresponding partial array of an element of $w$ showing $t_1 \in d_1 d_0(w)$ and $t_2 \in d_2 d_0(w)$:
\[
\begin{array}{|c||c|c|c|c|c|c||l|}
\hline
\text{row} \downarrow & 0 & 1 & 2 & 3 & \cdots & n & \\
\hline
*0* & & 1 & 1 & & & & (1) \text{ given} \\
\hline
1 & & & & & & &  \\
\hline
2 & 1 & & & & & &  \\
\hline
3 & 1 & & & & & &  \\
\hline
\vdots & & & & & & &  \\
\hline
n+1 & & & & & & &  \\
\hline 
\end{array}
\]
Use the 0-surjectivity of $y_2$ to fill row 2 (step 2) and the 1-surjectivity of $y_1$ to fill row 1 (step 3):
\[
\begin{array}{|c||c|c|c|c|c|c||l|}
\hline
\text{row} \downarrow & 0 & 1 & 2 & 3 & \cdots & n & \\
\hline
*0* & 3 & 1 & 1 & & & & (1) \text{ given} \\
\hline
1 & 3 & 2 & 3 & 3 & \cdots & 2 & (3)\; 1\text{-sur}  \\
\hline
2 & 1 & 2 & 2 & 2 & \cdots & 2 & (2)\; 0\text{-sur} \\
\hline
3 & 1 & 3 & 2 & & & &  \\
\hline
\vdots & & \vdots & \vdots & & & &  \\
\hline
n+1 & & 3 & 2 & & & &  \\
\hline 
\end{array}
\]
Then rows 3 through $n+1$ may be filled using that $y_3, \cdots , y_{n+1}$ are \DET{\SETT{1,2}}; row 0 is determined by $(n,0)$-composition. Therefore there exists $b \in w$ such that $e_1(e_0(b)) = t_1$ and $e_2(e_0(b)) = t_2$, as required.

\item 
To prove $w$ is \SUR{1}, consider an arbitrary $b_1 \in y_1$ and the array for the corresponding partial element of $w$:
\[
\begin{array}{|c||c|c|c|c|c|c||l|}
\hline
\text{row} \downarrow & 0 & 1 & 2 & 3 & \cdots & n & \\
\hline
*0* & 1 & & & & & & \\
\hline
1 & 1 & 1 & 1 & 1 & \cdots & 1 & (1) \text{ given} \\
\hline
2 & & 1 & & & & & \\
\hline
3 & & 1 & & & & & \\
\hline
\vdots & & \vdots & & & & & \\
\hline
n+1 & & 1 & & & & & \\
\hline 
\end{array}
\]
Since $y_2$ is \SUR{1}, then there is an element of $y_2$ which fills row 2 (step 2) to give:
\[
\begin{array}{|c||c|c|c|c|c|c||l|}
\hline
\text{row} \downarrow & 0 & 1 & 2 & 3 & \cdots & n & \\
\hline
*0* & 1 & 2 & & & & & \\
\hline
1 & 1 & 1 & 1 & 1 & \cdots & 1 & (1) \text{ given} \\
\hline
2 & 2 & 1 & 2 & 2 & \cdots & 2 & (2)\; 1\text{-sur}\\
\hline
3 & & 1 & 2 & & & & \\
\hline
\vdots & & \vdots & \vdots & & & & \\
\hline
n+1 & & 1 & 2 & & & & \\
\hline 
\end{array}
\]
Again, using that $y_3, \cdots , y_{n+1}$ are \DET{\SETT{1,2}}, rows $3$ through $n+1$ can be filled, and row 0 is then determined by $(n,0)$-composition, filling out the array. Therefore, there exists an element $b \in w$ such that $e_1(b) = b_1$, as required.

\item 
That $w$ is \DET{\{1,2\}} follows from the Up Rule, Theorem \refpage{up rule}.

\end{enumerate}
\qed

\begin{newthm}
{\bf (Existence of \NICOMP{n}{0}s)}
\label{composer-model;i=0}
\index{Theorem!Existence of $(n,0)$-composers}
{\rm
\SL

Fix $n \geq 2$. Suppose $\MC{S}'_{\leq n}$ is an $n$-truncated subcomplex of $\MC{S}$ such that every $n$-simplex of $\MC{S}'_{\leq n}$ is \SBSI, \SUR{0}, \SUR{1} and $\DPQW{1}{2}$. Then $\MC{S}'_{\leq n}$ is the $n$-truncation of an \NICOMP{n}{0} $\MC{S}'$.
}
\end{newthm}

\Proof

We will define $\MC{S}'_m$ for $m \geq n+1$ by induction on $m$. At each stage we will have to show that $\MC{S}'_{\leq m}$ is an $m$-truncated subcomplex of $\MC{S}$. This will yield a subcomplex $\MC{S}' \SBS \MC{S}$ such that $\TR^n(\MC{S})$ is the given $n$-truncated complex. Then we will show that $\MC{S'}_m \ISO \BOX{m}{0}{\MC{S}'}$ for all $m \geq n+1$.
\MS

We begin by defining $\MC{S}'_{n+1}$ to consist of all $(n,0)$-compositions 
\[
w = \COMP{n}{0}(-, y_1 \DDD{,} y_{n+1})
\]
such that for each $j \in [n+1]-\{0\}$, $y_j \in \MC{S}'_n$. In order to show that $\MC{S}'_{n+1}$ extends $\MC{S}'_{\leq n}$ to dimension $n+1$ we have to show that $d_0(w) = \MC{S}'_n$ and that for all $y \in \MC{S}'_n$ and all $k \in [n]$, $s_k(y) \in \MC{S}'_{n+1}$.
\MS

That $d_0(w) \in \MC{S}'_n$ is an immediate consequence of Theorem \refpage{n-0-parameter-theorem}.
\MS

Let $y \in \MC{S}'_n$ and $k \in [n]$. Since $y$ is assumed \SUR{0}, it follows from the Degenerate Compositions Theorem that $s_k(y)$ is an $(n,0)$-composition. Since all the faces of $s_k(y)$ belong to $\MC{S}'_n$ then $s_k(y) \in \MC{S}'_{n+1}$.
\MS

Now we define $\MC{S}'_m$ for $m>n+1$ inductively. We define $\MC{S}'_m$ to consist of all $(m-1,0)$-compositions
\[
z = \COMP{m-1}{0}(-,u_1 \DDD{,} u_m)
\]
where for each $j \neq 0$, $u_j \in \MC{S}'_{m-1}$. The induction hypothesis implies that each $u_j$ is \SBSI, \SUR{0}, \SUR{1} and $\DPQ{1}{2}{(m-1)}$. The same reasoning as in the $m=n+1$ case shows that $d_0(z) \in \MC{S}'_{m-1}$ and for all $k \in [m-1]$ and all $u \in \MC{S}'_{m-1}$, $s_k(u) \in \MC{S}'_m$.

We therefore obtain a subcomplex $\MC{S}' \SBS \MC{S}$ such that for all $m \geq n+1$, $\MC{S}'_m \ISO \BOX{m}{0}{\MC{S}'}$ and $\TR^n(\MC{S}
)$ is the given subcomplex.

\qed

\subsection{\NICOMP{n}{1}s} 

\begin{newthm}
{\bf ($(n,1)$-parameter theorem)}
\label{n-1-parameter-theorem}
\index{Theorem!$(n,1)$-parameter theorem}
{\rm
\SL

Fix $n \geq 2$. Suppose $w = \COMP{n}{1}(y_0,-,y_2, \cdots, y_{n+1})$ is an $(n,1)$-composition and that each $y_k$, $k \neq 1$, is subface-simplicial, \SUR{1}, and $\DPQW{0}{2}$. Then:
\begin{enumerate}
\item $d_1(w)$ is \SUR{1}.
\item $d_1(w)$ is $\DPQW{0}{2}$.
\item $w$ is \DET{\SETT{0,2}}.
\end{enumerate}
Note: The parameter is $n$.
}
\end{newthm}

\Proof

\begin{enumerate}
\item 
Given an arbitrary $t \in d_1d_1(w)$, we must show there exists $b_1 \in d_1(w)$ such that $e_1(b_1) = t_1$. Start with the corresponding partial element of $w$ as the array:
\[
\begin{array}{|c||c|c|c|c|c|c||l|}
\hline
\text{row} \downarrow & 0 & 1 & 2 & 3 & \cdots & n & \\
\hline
0 & & & & & & & \\
\hline
*\; 1 \; * & & 1 & & & & & (1)\; t \text{, given} \\
\hline
2 & & 1 & & & & & \\
\hline
3 & & & & & & & \\
\hline
\vdots & & & & & & & \\
\hline
n+1 & & & & & & & \\
\hline
\end{array}
\]
Since $y_2$ is \SUR{1}, there exists an element of $y_2$ which fills row 2 (step 2) to give:
\[
\begin{array}{|c||c|c|c|c|c|c||l|}
\hline
\text{row} \downarrow & 0 & 1 & 2 & 3 & \cdots & n & \\
\hline
0 & & 2 & & & & & \\
\hline
*\; 1 \; * & & 1 & & & & & (1)\; t \text{, given} \\
\hline
2 & 2 & 1 & 2 & 2 & \cdots & 2 & (2)\; 1\text{-sur} \\
\hline
3 & & & 2 & & & & \\
\hline
\vdots & & & \vdots & & & & \\
\hline
n+1 & & & 2 & & & & \\
\hline
\end{array}
\]
Similarly, the 1-surjectivity of $y_0$ implies there is an element of $y_0$ which fills row 0 (step 3):
\[
\begin{array}{|c||c|c|c|c|c|c||l|}
\hline
\text{row} \downarrow & 0 & 1 & 2 & 3 & \cdots & n & \\
\hline
0 & 3 & 2 & 3 & 3 & \cdots & 3 & (3)\; 1\text{-sur} \\
\hline
*\; 1 \; * & 3 & 1 & & & & & (1)\; t \text{, given} \\
\hline
2 & 2 & 1 & 2 & 2 & \cdots & 2 & (2)\; 1\text{-sur} \\
\hline
3 & 3 & & 2 & & & & \\
\hline
\vdots & \vdots & & \vdots & & & & \\
\hline
n+1 & 3 & & 2 & & & & \\
\hline
\end{array}
\]
Now, using that $y_3, \cdots , y_{n+1}$ are all $\DPQW{0}{2}$, there are elements of those $n$-simplices which fill rows 3 through $n+1$. Row 1 is then determined by $(n,1)$-composition and the completed array represents an element $b \in w$. Then $b_1 = e_1(w) \in d_1(w)$ has $e_1(b_1) = t$, which proves that $d_1(w)$ is \SUR{1}.

\item 
To verify that $d_1(w)$ is $\DPQW{0}{2}$ suppose $t_0 \in d_0(d_1(w))$ and $t_2 \in d_2(d_1(w))$ are components of a $\SETT{0,2}$-indexed partial element of $d_1(w)$. Consider the corresponding partial element of $w$ as the array:
\[
\begin{array}{|c||c|c|c|c|c|c||l|}
\hline
\text{row} \downarrow & 0 & 1 & 2 & 3 & \cdots & n & \\
\hline
0 & 1 & & & & & & \\
\hline
*\; 1 \; * & 1 & & 1 & & & & (1)\; t_0,t_2 \text{ given} \\
\hline
2 & & & & & & & \\
\hline
3 & & 1 & & & & & \\
\hline
\vdots & & & & & & & \\
\hline
n+1 & & & & & & & \\
\hline
\end{array}
\]
By the 1-surjectivity of $y_3$ there exists an element of $y_3$ which fills row 3 (step 2):
\[
\begin{array}{|c||c|c|c|c|c|c||l|}
\hline
\text{row} \downarrow & 0 & 1 & 2 & 3 & \cdots & n & \\
\hline
0 & 1 & & 2 & & & & \\
\hline
*\; 1 \; * & 1 & & 1 & & & & (1)\; t_0,t_2 \text{ given} \\
\hline
2 & & & 2 & & & & \\
\hline
3 & 2 & 1 & 2 & 2 & \cdots & 2 & (2)\; 1 \text{-sur} \\
\hline
\vdots & & & & \vdots & & & \\
\hline
n+1 & & & & 2 & & & \\
\hline
\end{array}
\]
Next, row 0 may be filled (step 3) in because $y_0$ is $\DPQW{0}{2}$:
\[
\begin{array}{|c||c|c|c|c|c|c||l|}
\hline
\text{row} \downarrow & 0 & 1 & 2 & 3 & \cdots & n & \\
\hline
0 & 1 & 3 & 2 & 3 & \cdots & 3 & (3)\;  \DPQ{0}{2}{n}\\
\hline
*\; 1 \; * & 1 & & 1 & & & & (1)\; t_0,t_2 \text{ given} \\
\hline
2 & 3 & & 2 & & & & \\
\hline
3 & 2 & 1 & 2 & 2 & \cdots & 2 & (2)\; 1 \text{-sur} \\
\hline
\vdots & \vdots & & \vdots & & & & \\
\hline
n+1 & 3 & & 2 & & & & \\
\hline
\end{array}
\]
Row 2 and rows 4 through $n+1$ can then be filled using that the corresponding $y_k$ are all $\DPQW{0}{2}$; row 1 is filled by $(n,1)$-composition. The completed array represents an element $b \in w$, and $b_1 = e_1(w)$ is an element of $d_1(w)$ such that $e_0(b_1)=t_0$ and $e_2(b_1) = t_2$ as required.

\item 
That $w$ is \DET{\SETT{0,2}} follows from the Up Rule.

\end{enumerate}
\qed

\begin{newthm}
{\bf (Existence of \NICOMP{n}{1}s)}
\label{n-1-composer-theorem}
\label{composer-model;i=1}
\index{Theorem!Existence of \NICOMP{n}{1}s}
{\rm
\SL

Fix $n \geq 2$. Suppose $\MC{S}'_{\leq n}$ is an $n$-truncated subcomplex of $\MC{S}$ such that each $n$-simplex is \SBSI, \SUR{1} and $\DPQW{0}{2}$. Then there is a subcomplex $\MC{S}'$ of $S$ whose $n$-truncation is the given one and which is an \NICOMP{n}{1}.
}
\end{newthm}

\Proof

We will define $\MC{S}'_m$ for $m \geq n+1$ by induction on $m$. At each stage we will have to show that $\MC{S}'_{\leq m}$ is an $m$-truncated subcomplex of $\MC{S}$. This will yield a subcomplex $\MC{S}' \SBS \MC{S}$ such that $\TR^n(\MC{S})$ is the given $n$-truncated complex. Then we will show that $\MC{S'}_m \ISO \BOX{m}{1}{\MC{S}'}$ for all $m \geq n+1$.
\MS

We begin by defining $\MC{S}'_{n+1}$ to consist of all $(n,1)$-compositions 
\[
w = \COMP{n}{1}(y_0, - \DDD{,} y_{n+1})
\]
such that for each $j \neq 1$, $y_j \in \MC{S}'_n$. In order to show that $\MC{S}'_{n+1}$ extends $\MC{S}'_{\leq n}$ to dimension $n+1$ we have to show that $d_1(w) = \MC{S}'_n$ and that for all $y \in \MC{S}'_n$ and all $k \in [n]$, $s_k(y) \in \MC{S}'_{n+1}$.
\MS

That $d_1(w) \in \MC{S}'_n$ is an immediate consequence of Theorem \refpage{n-1-parameter-theorem}. Also, note that $w$ is $\DPQ{0}{2}{(n+1)}$.
\MS

Let $y \in \MC{S}'_n$ and $k \in [n]$. Since $y$ is assumed \SUR{1}, it follows from the Degenerate Compositions Theorem that $s_k(y)$ is an $(n,1)$-composition. Since all the faces of $s_k(y)$ belong to $\MC{S}'_n$ then $s_k(y) \in \MC{S}'_{n+1}$.
\MS

Now we define $\MC{S}'_m$ for $m>n+1$ inductively. We define $\MC{S}'_m$ to consist of all $(m-1,1)$-compositions
\[
z = \COMP{m-1}{1}(u_0,- \DDD{,} u_m)
\]
where for each $j \neq 1$, $u_j \in \MC{S}'_{m-1}$. The induction hypothesis implies that each $u_j$ is \SBSI, \SUR{1} and $\DPQ{0}{2}{(m-1)}$. The same reasoning as in the $m=n+1$ case shows that $d_1(z) \in \MC{S}'_{m-1}$ and for all $k \in [m-1]$ and all $u \in \MC{S}'_{m-1}$, $s_k(u) \in \MC{S}'_m$.

We therefore obtain a subcomplex $\MC{S}' \SBS \MC{S}$ such that for all $m \geq n+1$, $\MC{S}'_m \ISO \BOX{m}{1}{\MC{S}'}$ and $\TR^n(\MC{S}
)$ is the given subcomplex.

\qed

\subsection{\NICOMP{n}{n}s} 

\begin{newthm}
{\bf ($(n,n)$-parameter)}
\label{n-n-parameter-theorem}
\index{Theorem!$(n,n)$-parameter theorem}
{\rm
\SL

Fix $n \geq 2$. Suppose $w = \COMP{n}{n}(y_0,y_1, \cdots, y_{n-1},-, y_{n+1})$ is an $(n,n)$-composition and that each $y_k$, $k \neq n$, is subface-simplicial, \SUR{(n-1)}, and $\DPQW{n-2}{n}$. Then:
\begin{enumerate}
\item $d_n(w)$ is \SUR{(n-1)}.
\item $d_n(w)$ is $\DPQW{n-2}{n}$
\item $w$ is \SUR{(n-1)}.
\item $w$ is $\DPQW{n-2}{n+1}$.
\item $w$ is $\DPQW{n-1}{n+1}$.
\end{enumerate}
Note: The parameter is $n$.
}
\end{newthm}

\Proof

\begin{enumerate}
\item
\NI {\bf That $d_n(w)$ is \SUR{(n-1)}:}

Given $t \in d_{n-1}(d_n(w))$ consider the partial array of the corresponding potential element of $w$:
\[
\begin{array}{|c||c|c|c|c|c|c||l|}
\hline
\text{row} \downarrow & 0 & \cdots & n-3 & n-2 & n-1 & n & \\
\hline
0 & & & & & & & \\
\hline

\vdots & & & & & & & \\
\hline

n-2 & & & & & & & \\
\hline

n-1 & & & & & 1 & & \\
\hline

*\; n\; * & & & & & 1 & & (1) \quad $t$\text{, given} \\
\hline

n+1 & & & & & & & \\
\hline
\end{array}
\]
For step 2, fill in row $n-1$ using \SURY{(n-1)}; then for step 3, fill row $n+1$ again using \SURY{(n-1)}. The resulting array is:
\[
\begin{array}{|c||c|c|c|c|c|c||l|}
\hline
\text{row} \downarrow & 0 & \cdots & n-3 & n-2 & n-1 & n & \\
\hline
0 & & & & 2 & & 3 & \\
\hline

\vdots & & & & \vdots& &\vdots & \\
\hline

n-2 & & & & 2 & & 3 & \\
\hline

n-1 & 2 & \cdots & 2 & 2 & 1 & 2 & (2) \quad \SURD{(n-1)}{n} \\
\hline

*\; n\; * & & & & & 1 & 3 & (1) \quad $t$\text{, given} \\
\hline

n+1 & 3 & \cdots & 3 & 3 & 2 & 3 & (3) \quad  \SURD{(n-1)}{n}\\
\hline
\end{array}
\]
Now, rows 0 through $n-2$ may be filled using \DETY{\SETT{n-2,n}} and row $n$ is determined by $(n,n)$-composition. The completed array represents an element of $w$, and therefore $d_n(w)$ is \SUR{(n-1)}, as claimed.

\item
\NI {\bf That $d_n(w)$ is $\DPQW{n-2}{n}$:}

Suppose $(-, \cdots, -,t_{n-2},-,t_n,-)$ is an $\SETT{n-2,n}$-indexed partial element of $d_n(w)$. The corresponding partial array is:
\[
\begin{array}{|c||c|c|c|c|c|c||l|}
\hline
\text{row} \downarrow & 0 & \cdots & n-3 & n-2 & n-1 & n & \\
\hline
0 & & & & & & & \\
\hline

\vdots & & & & & & & \\
\hline

n-2 & & & & & 1 & & \\
\hline

n-1 & & & & & & & \\
\hline

*\; n\; * & & & & 1 & & 1 & (1) \quad \text{given} \\
\hline

n+1 & & & & & & 1 & \\
\hline
\end{array}
\]
For step 2, use \SURY{(n-1)} to fill row $n-2$ and, for step 3 use \DETY{\SETT{n-2,n}} to fill row $n+1$. The resulting array is:
\[
\begin{array}{|c||c|c|c|c|c|c||l|}
\hline
\text{row} \downarrow & 0 & \cdots & n-3 & n-2 & n-1 & n & \\
\hline
0 & & & 2 & & & 3 & \\
\hline

\vdots & & & \vdots & & & \vdots & \\
\hline

n-2 & 2 & \cdots & 2 & 2 & 1 & 2 & (2) \quad \SURD{(n-1)}{n} \\
\hline

n-1 & & & & 2 & & 3 & \\
\hline

*\; n\; * & & & & 1 & & 1 & (1) \quad \text{given} \\
\hline

n+1 & 3 & \cdots & 3 & 2 & 3 & 1 & (3) \quad \DPQ{n-2}{n}{n}  \\
\hline
\end{array}
\]
For step 4, fill row $n-1$ by \DETY{\SETT{n-2,n}} to give
\[
\begin{array}{|c||c|c|c|c|c|c||l|}
\hline
\text{row} \downarrow & 0 & \cdots & n-3 & n-2 & n-1 & n & \\
\hline
0 & & & 2 & 4 & & 3 & \\
\hline

\vdots & & & \vdots & \vdots & & \vdots & \\
\hline

n-2 & 2 & \cdots & 2 & 2 & 1 & 2 & (2) \quad \SURD{(n-1)}{n} \\
\hline

n-1 & 4 & \cdots & 4 & 2 & 4 & 3 & (4) \quad \DPQ{n-2}{n}{n} \\
\hline

*\; n\; * & & & & 1 & 4 & 1 & (1) \quad \text{given} \\
\hline

n+1 & 3 & \cdots & 3 & 2 & 3 & 1 & (3) \quad \DPQ{n-2}{n}{n}  \\
\hline
\end{array}
\]
Now, rows 0 through $n-3$ may be filled using \DETY{\SETT{n-2,n}} and row $n$ filled by $(n,n)$-composition. The completed array represents an element of $w$, and therefore $d_n(w)$ is \DET{\SETT{n-2,n}}, as claimed.

\item
\NI {\bf That $w$ is \SUR{(n-1)}:}

Given any element of $d_{n-1}(w)$ consider the corresponding partial array:
\[
\begin{array}{|c||c|c|c|c|c|c||l|}
\hline
\text{row} \downarrow & 0 & \cdots & n-3 & n-2 & n-1 & n & \\
\hline
0 & & & & 1 & & & \\
\hline

\vdots & & & & \vdots & & & \\
\hline

n-2 & & & & 1 & & & \\
\hline

n-1 & 1 & \cdots & 1 & 1 & 1 & 1 & (1) \quad \text{element given} \\
\hline

n & & & & & 1 & & \\
\hline

n+1 & & & & & 1 & & \\
\hline
\end{array}
\]
For step 2, fill in row $n+1$ using \SURY{(n-1)} to get
\[
\begin{array}{|c||c|c|c|c|c|c||l|}
\hline
\text{row} \downarrow & 0 & \cdots & n-3 & n-2 & n-1 & n & \\
\hline
0 & & & & 1 & & 2 & \\
\hline

\vdots & & & & \vdots & & \vdots & \\
\hline

n-2 & & & & 1 & & 2 & \\
\hline

n-1 & 1 & \cdots & 1 & 1 & 1 & 1 & (1) \quad \text{element given} \\
\hline

*n* & & & & & 1 & 2 & \\
\hline

n+1 & 2 & \cdots & 2 & 2 & 1 & 2 & (2) \quad \SURD{(n-1)}{n} \\
\hline
\end{array}
\]
Now, rows $0$ through $n-2$ can be filled using \DETY{\SETT{n-2,n}} and row $n$ is filled by $(n,n)$-composition. Therefore, $w$ is \SUR{(n-1)} as claimed.

\item

That $w$ is \DET{\SETT{n-2,n+1}} follows from the Up Rule.

\item
That $w$ is $\DPQW{n-1}{n+1}$ follows from the Up Rule.

\end{enumerate}
\qed

\begin{newthm} 
{\bf (Existence of \NICOMP{n}{n}s)}
\label{n-n-composer-theorem}
\label{composer-model;i=n}
\index{Theorem!Existence of \NICOMP{n}{n}s}
{\rm
\SL

Fix $n \geq 2$. Suppose $\MC{S}'_{\leq n}$ is an $n$-truncated subcomplex of $\MC{S}$ such that each $n$-simplex is \SBSI, \SUR{(n-1)} and $\DPQW{n-2}{n}$. Then there is a subcomplex $\MC{S}'$ of $S$ whose $n$-truncation is the given one and which is an \NICOMP{n}{n}.
}
\end{newthm}

\Proof

We begin by defining $\MC{S}'_{n+1}$ to consist of all $(n,n)$-compositions 
\[
w = \COMP{n}{n}(y_0 \DDD{,} y_{n-1},-,y_{n+1})
\]
such that for each $j \neq n$, $y_j \in \MC{S}'_n$. In order to show that $\MC{S}'_{n+1}$ extends $\MC{S}'_{\leq n}$ to dimension $n+1$ we have to show that $d_n(w) = \MC{S}'_n$ and that for all $y \in \MC{S}'_n$ and all $k \in [n]$, $s_k(y) \in \MC{S}'_{n+1}$.
\MS

That $d_n(w) \in \MC{S}'_n$ is an immediate consequence of Theorem \refpage{n-n-parameter-theorem}. Also, note that $w$ is  \SUR{(n-1)} and $\{n-2,n+1\}$-determinate.
\MS

Let $y \in \MC{S}'_n$ and $k \in [n]$. Since $y$ is assumed \SUR{(n-1)}, it follows from the Degenerate Compositions Theorem that $s_k(y)$ is an $(n,n)$-composition. Since all the faces of $s_k(y)$ belong to $\MC{S}'_n$ then $s_k(y) \in \MC{S}'_{n+1}$.
\MS

Now consider the truncated complex $\MC{S}'_{\leq n+1}$. According to Theorem \refpage{n-n-parameter-theorem}, every $w \in \MC{S}'_{n+1}$ is \SUR{(n-1)} and $\DPQW{n-2}{n+1}$. $w$ is also \SUR{n} because it is an $(n,n)$-composition. We may therefore apply  Theorem \refpage{composer-model;1<i<n}, to $n'=n+1$, $i'=n$ because $1<i'<n'$, and all $n'$-simplices of $\MC{S}'_{n'}$ are \SUR{(i'-1)}, \SUR{i'} and $\DPQW{i'-2}{i'+1}$. Thus there exists a subcomplex $\MC{S}'$ such that $\MC{S}'_{\leq n+1} = \TR^{n+1}(\MC{S}')$ and for all $m \geq n+1$, $\MC{S}'_m \ISO \BOX{i'}{m}{\MC{S}'} = \BOX{n}{m}{\MC{S}'}$. That is, $\MC{S}'$ is an \NICOMP{n}{n}.

\qed

\subsection{\NICOMP{n}{n+1}s} 

\begin{newthm}
{\bf ($(n,n+1)$-parameter)}
\label{n-n+1-parameter-theorem}
\index{Theorem!$(n,n+1)$-parameter theorem}
{\rm
\SL

Fix $n \geq 2$. Suppose $w = \COMP{n}{n+1} \Bigl(y_0, \cdots , y_n,\; -\; \Bigr)$
is an $(n,n+1)$-composition such that each $y_j$ is \SBSI, \SUR{n}, \SUR{(n-1)}, \DET{\SETT{n-2,n-1}} and \DET{\SETT{n-1,n}}. Then:
\begin{enumerate}
\item $d_{n+1}(w)$ is \SUR{n}.
\item $d_{n+1}(w)$ is \SUR{(n-1)}.
\item $d_{n+1}(w)$ is \DET{\SETT{n-2,n-1}}.
\item $d_{n+1}(w)$ is \DET{\SETT{n-1,n}}.
\item $w$ is \SUR{n}.
\item $w$ is \DET{\SETT{n-1,n+1}}.
\end{enumerate}
}
\end{newthm}

\NI Note: In this theorem, the ``extra'' determinacy condition $\DPQ{n-1}{n}{n}$ is used only to prove that $w$ is \DET{\SETT{n-1,n+1}}.
\MS

\Proof

\begin{enumerate}

\item
{\bf That $d_{n+1}(w)$ is \SUR{n}:}

Given $t \in d_nd_{n+1}(w)$ then the array corresponding to that partial element of $w$ is
\[
\begin{array}{|c||c|c|c|c|c|c|c||l|}
\hline
\text{row} \downarrow & 0 & \cdots & n-4 & n-3 & n-2 & n-1 &n & \\
\hline 
0  & & & & & & & & \\
\hline 
\vdots & & & & & & & & \\
\hline 
n-3 & & & & & & & & \\
\hline 
n-2 & & & & & & & & \\
\hline 
n-1 & & & & & & & & \\
\hline 
n & & & & & & & 1 & \\
\hline 
*\; n+1 \; * & & & & & & & 1 & (1)\; \text{given} \\
\hline 
\end{array}
\]
For step 2, fill row $n$ by the \SURY{n} of $y_n$ and, as step 3, fill row $n-1$ by the \SURY{(n-1)} of $y_{n-1}$. The array is now:
\[
\begin{array}{|c||c|c|c|c|c|c|c||l|}
\hline
\text{row} \downarrow & 0 & \cdots & n-4 & n-3 & n-2 & n-1 &n & \\
\hline 
0  & & & & & 3 & 2 & & \\
\hline 
\vdots & & & & & \vdots & \vdots & & \\
\hline 
n-3 & & & & & 3 & 2 & & \\
\hline 
n-2 & & & & & 3 & 2 & & \\
\hline 
n-1 & 3 & \cdots & 3 & 3 & 3 & 2 & 3 & (3)\; \SURD{(n-1)}{n} \\
\hline 
n & 2 & \cdots & 2 & 2 & 2 & 2 & 1 & (2)\; \SURD{n}{n} \\
\hline 
*\; n+1 \; * & & & & & & & 1 & (1)\; t\text{, given} \\
\hline 
\end{array}
\]
Then rows $0$ through $n-2$ are filled by the \DETY{\SETT{n-2,n-1}} of $y_0, \cdots ,y_{n-2}$. Row $n+1$ is determined by $(n,n+1)$-composition. The completed array represents an element of $w$ and row $n+1$ represents an element $b_{n+1} \in d_{n+1}(w)$ such that $e_n(b_{n+1})=t$, as required.

\item

{\bf That $d_{n+1}(w)$ is \SUR{(n-1)}:}

Given $t \in d_{n-1}d_{n+1}(w)$ then the array corresponding to that partial element of $w$ is
\[
\begin{array}{|c||c|c|c|c|c|c|c||l|}
\hline
\text{row} \downarrow & 0 & \cdots & n-4 & n-3 & n-2 & n-1 &n & \\
\hline 
0  & & & & & & & & \\
\hline 
\vdots & & & & & & & & \\
\hline 
n-3 & & & & & & & & \\
\hline 
n-2 & & & & & & & & \\
\hline 
n-1 & & & & & & & 1 & \\
\hline 
n & & & & & & & & \\
\hline 
*\; n+1 \; * & & & & & & 1 & & (1)\; t\text{, given} \\
\hline 
\end{array}
\]
For step 2, fill row $n-1$ by \SURY{n} of $y_{n-1}$ and, for step 3, fill row $n$ by the \SURY{(n-1)} of $y_n$. The array is now:
\[
\begin{array}{|c||c|c|c|c|c|c|c||l|}
\hline
\text{row} \downarrow & 0 & \cdots & n-4 & n-3 & n-2 & n-1 &n & \\
\hline 
0  & & & & & 2 & 3 & & \\
\hline 
\vdots & & & & & \vdots & \vdots & & \\
\hline 
n-3 & & & & & 2 & 3 & & \\
\hline 
n-2 & &  & & & 2 & 3 & & \\
\hline 
n-1 & 2 & \cdots & 2 & 2 & 2 & 2 & 1 & (2)\; \SURD{n}{n} \\
\hline 
n & 3 & \cdots & 3 & 3 & 3 & 2 & 3 & (3)\; \SURD{(n-1)}{n} \\
\hline 
*\; n+1 \; * & & & & & & 1 & 3 & (1)\; t \text{, given} \\
\hline 
\end{array}
\]
Now fill rows $0$ through $n-2$ using the \DETY{\SETT{n-2,n-1}} of $y_0, \cdots , y_{n-2}$, and fill row $n+1$ by $(n,n+1)$-composition. The completed array represents an element of $w$ and row $n+1$ represents $b_{n+1} \in d_{n+1}(w)$ such that $e_{n-1}(b_{n+1})=t$, as required.

\item
{\bf That $d_{n+1}(w)$ is \DET{\SETT{n-2,n-1}}:}

Suppose $\SETT{t_{n-2},t_{n-1}}$ is an $\SETT{n-2,n-1}$-indexed partial element of $d_{n+1}(w)$. The corresponding array is:
\[
\begin{array}{|c||c|c|c|c|c|c|c||l|}
\hline
\text{row} \downarrow & 0 & \cdots & n-4 & n-3 & n-2 & n-1 &n & \\
\hline 
0  & & & & & & & & \\
\hline 
\vdots & & & & & & & & \\
\hline 
n-3 & & & & & & & & \\
\hline 
n-2 & & & & & & & 1 & \\
\hline 
n-1 & & & & & & & 1 & \\
\hline 
n & & & & & & & & \\
\hline 
*\; n+1 \; * & & & & & 1 & 1 & & (1)\; \text{given} \\
\hline 
\end{array}
\]
For step 2, fill row $n-1$ using the \SURY{n} of $y_{n-1}$ and, for step 3, fill row $n$ using the \SURY{(n-1)} of $y_n$. The resulting array is:
\[
\begin{array}{|c||c|c|c|c|c|c|c||l|}
\hline
\text{row} \downarrow & 0 & \cdots & n-4 & n-3 & n-2 & n-1 &n & \\
\hline 
0  & & & & & 2 & 3 & & \\
\hline 
\vdots & & & & & \vdots & \vdots & & \\
\hline 
n-3 & & & & & 2 & 3 & & \\
\hline 
n-2 & & & & & 2 & 3 & 1 & \\
\hline 
n-1 & 2 & \cdots & 2 & 2 & 2 & 2 & 1 & (2)\; \SURD{n}{n} \\
\hline 
n & 3 & \cdots & 3 & 3 & 3 & 2 & 3 & (2) \SURD{(n-1)}{n} \\
\hline 
*\; n+1 \; * & & & & & 1 & 1 & 3 & (1)\; \text{given} \\
\hline 
\end{array}
\]
Again, rows $0$ through $n-2$ can be filled by the \DETY{\SETT{n-2,n-1}} of  $y_0, \cdots , y_{n-2}$, and row $n+1$ is filled by $(n,n+1)$-composition. The completed array represents an element $b \in w$ and $b_{n+1}= e_{n+1}(b)$ fills row $n+1$ with the required components.

\item
{\bf That $d_{n+1}(w)$ is \DET{\SETT{n-1,n}}:}

Suppose $\SETT{t_{n-1},t_{n}}$ is an $\SETT{n-1,n}$-indexed partial element of $d_{n+1}(w)$. The corresponding array is: 
\[
\begin{array}{|c||c|c|c|c|c|c|c||l|}
\hline
\text{row} \downarrow & 0 & \cdots & n-4 & n-3 & n-2 & n-1 &n & \\
\hline 
0  & & & & & & & & \\
\hline 
\vdots & & & & & & & & \\
\hline 
n-3 & & & & & & & & \\
\hline 
n-2 & & & & & & & & \\
\hline 
n-1 & & & & & & & 1 & \\
\hline 
n & & & & & & & 1 & \\
\hline 
*\; n+1 \; * & & & & & & 1 & 1 & (1)\; \text{given} \\
\hline 
\end{array}
\]
For step 2, fill row $n$ by the \SURY{n} of $y_n$, and then for step 3, fill row $n-1$ by the \DETY{\SETT{n-1,n}} of $y_{n-1}$. The resulting array is
\[
\begin{array}{|c||c|c|c|c|c|c|c||l|}
\hline
\text{row} \downarrow & 0 & \cdots & n-4 & n-3 & n-2 & n-1 &n & \\
\hline 
0  & & & & & 3 & 2 & & \\
\hline 
\vdots & & & & & \vdots & \vdots & & \\
\hline 
n-3 & & & & & 3 & 2 & & \\
\hline 
n-2 & & & & & 3 & 2 & & \\
\hline 
n-1 & 3 & \cdots & 3 & 3 & 3 & 2 & 1 & (3)\; \DPQ{n-1}{n}{n} \\
\hline 
n & 2 & \cdots & 2 & 2 & 2 & 2 & 1 & (2)\; \SURD{n}{n} \\
\hline 
*\; n+1 \; * & & & & & & 1 & 1 & (1)\; \text{given} \\
\hline 
\end{array}
\]
At this point, the rest of the array may be completed using the \DETY{\SETT{n-2,n-1}} of $y_0, \cdots , y_{n-2}$ and $(n,n+1)$-composition.
It represents an element $b \in w$ whose components include $\SETT{t_{n-1},t_{n}}$, as required.

\item
{\bf That $w$ is \SUR{n}:}

Given $b_n \in d_n(w)$, the corresponding array is
\[
\begin{array}{|c||c|c|c|c|c|c|c||l|}
\hline
\text{row} \downarrow & 0 & \cdots & n-4 & n-3 & n-2 & n-1 &n & \\
\hline 
0  & & & & & & 1 & & \\
\hline 
\vdots & & & & & & \vdots & & \\
\hline 
n-3 & & & & & & 1 & & \\
\hline 
n-2 & & & & & & 1 & & \\
\hline 
n-1 & & & & & & 1 & & \\
\hline 
n & 1 & \cdots & 1 & 1 & 1 & 1 & 1 & (1)\; \text{given} \\
\hline 
*\; n+1 \; * & & & & & & & 1 & \\
\hline 
\end{array}
\]
For step 2, fill row $n-1$ using the \SURY{(n-1)} of $y_{n-1}$ to get:
\[
\begin{array}{|c||c|c|c|c|c|c|c||l|}
\hline
\text{row} \downarrow & 0 & \cdots & n-4 & n-3 & n-2 & n-1 &n & \\
\hline 
0  & & & & & 2 & 1 & & \\
\hline 
\vdots & & & & & \vdots & \vdots & & \\
\hline 
n-3 & & & & & 2 & 1 & & \\
\hline 
n-2 & & & & & 2 & 1 & & \\
\hline 
n-1 & 2 & \cdots & 2 & 2 & 2 & 1 & 2 & (2)\; \SURD{(n-1)}{n} \\
\hline 
n & 1 & \cdots & 1 & 1 & 1 & 1 & 1 & (1)\; \text{given} \\
\hline 
*\; n+1 \; * & & & & & & 2 & 1 & \\
\hline 
\end{array}
\]
Then fill rows $0$ through $n-2$ using the \DETY{\SETT{n-2,n-1}} of $y_0, \cdots , y_{n-2}$ and row $n+1$ by $(n,n+1)$-composition. The completed array represents $b \in w$ such that $e_n(b)$ is the given $b_n$, as required.

\item
That $w$ is \DET{\SETT{n-1,n+1}} follows from the Up Rule.

\end{enumerate}
\qed

\begin{newthm} 
{\bf (Existence of \NICOMP{n}{n+1}s)}
\label{n-n+1-composer-theorem}
\label{composer-model;i=n+1}
\index{Theorem!Existence of \NICOMP{n}{n+1}s}
{\rm
\SL

Fix $n \geq 2$. Suppose $\MC{S}'_{\leq n}$ is an $n$-truncated subcomplex of $\MC{S}$ such that each $n$-simplex is \SBSI, \SUR{n}, \SUR{(n-1)}, \DET{\SETT{n-2,n-1}} and \DET{\SETT{n-1,n}}. Then there is a subcomplex $\MC{S}'$ of $\MC{S}$ whose $n$-truncation is the given one and which is an \NICOMP{n}{n+1}.
}
\end{newthm}

\Proof

We begin by defining  $\mathcal{S}'_{n+1}$ to consist of all $(n,n+1)$-compositions
\[
w = \COMP{n}{n+1} \left( y_0, \cdots, y_n, \; -\; \right)
\]
such that $y_j \in \MC{S}'_n$, all $j=0, \cdots, n$. In order to show that $\MC{S}'_{n+1}$ extends $\MC{S}'_{\leq n}$ to dimension $n+1$ we have to show that $d_{n+1}(w) \in \MC{S}'_n$, and show that for all $k \in [n]$, all $y \in \MC{S}'_n$, $s_k(y) \in \MC{S}'_{n+1}$.
\MS

The first four assertions of the $(n,n+1)$-parameter theorem (Theorem \refpage{n-n+1-parameter-theorem}) imply that $d_{n+1}(w) \in \MC{S}'_n$. Theorem \refpage{degen-comps-thm} implies that $s_j(y) \in \MC{S}'_{n+1}$ for all $j \in [n]$. Therefore, $\MC{S}'_{\leq n+1}$ forms an $(n+1)$-truncated subcomplex.
\MS

Next, set $n'=n+1=i'$. By definition, $w \in \MC{S}'_{n+1} = \MC{S}'_{n'}$ iff $w$ is \SUR{n} ($\iff$ \SUR{(n'-1)}) and $\{n-1,n+1\}$-determinate ($\iff$ $\{n'-2,n'\}$-determinate). Therefore Theorem \refpage{composer-model;i=n} (``Existence of \NICOMP{n}{n}s'') applies with parameter $n'$ to show there is a subcomplex $\MC{S}'$ such that for all $m \geq n'+1$, $\MC{S}'_m \ISO \BOX{i'}{m}{\MC{S}'} = \BOX{n+1}{m}{\MC{S}'}$. Since $\MC{S}'_{n+1} \ISO \BOX{n+1}{n+1}{\MC{S}'}$ (by definition), $\MC{S}'$ is an \NICOMP{n}{n+1}.

\qed

\section{Further results}
\label{furthermore}

\subsection{Derived conditions}

If $\MC{S}'$ is an \NICOMP{n}{i} of sets then the surjectivity and determinacy conditions which hold in dimension $n$ imply other such conditions. These other conditions follow from Downward Surjectivity (Corollary \refpage{downward-surj-corollary}), the Up Rule (Theorem \refpage{up rule}) and the Down Rule (Theorem \refpage{down rule}). 
\MS

Let $n \geq 2$. The results below apply to any subcomplex $\MC{S}'$ of $\MC{S}$ such that all $(n+1)$-simplices of $\MC{S}'$ are \SBSI\ and maximal e-simplicial. These conditions hold for all the models of \NICOMP{n}{i}s developed in section \refpage{nisolution}. \BS

\HD{Determinacy}
\MS

Starting with  \DETY{\{p,q\}} the results (ignoring repetitions) of applying the Up Rule and then the Down Rule are  as follows:
\[
\DPQ{p}{q}{n} \IMPBY{\text{Up }1} \DPQ{p}{q}{n+1} \IMPBY{\text{Down}}
\left\{
\begin{array}{ll}
\DPQ{p-1}{q-1}{n} & p \geq 2\\
\DPQ{p}{q-1}{n} & q-p \geq 3
\end{array}
\right.
\]
\[
\DPQ{p}{q}{n} \IMPBY{\text{Up }2} \DPQ{p}{q+1}{n+1} \IMPBY{\text{Down}}
\left\{
\begin{array}{ll}
\DPQ{p-1}{q}{n} & p \geq 2\\
\DPQ{p}{q+1}{n} & q+1 \leq n-1
\end{array}
\right.
\]
\[
\DPQ{p}{q}{n} \IMPBY{\text{Up }3} \DPQ{p+1}{q+1}{n+1} \IMPBY{\text{Down}}
\left\{
\begin{array}{ll}
\DPQ{p+1}{q}{n} & q-p \geq 3\\
\DPQ{p+1}{q+1}{n} & q+1 \leq n-1
\end{array}
\right.
\]
Summarized in one table:

\begin{newdef}
\index{Down-Up Rule}
{\rm
{\bf ``Down-Up'' Rule}:
\[
\DPQ{p}{q}{n} \IMP 
\left\{
\begin{array}{llr}
\DPQ{p-1}{q-1}{n} & \text{when }p \geq 2 & \text{DU }1 \\ 
\DPQ{p}{q-1}{n} & \text{when }q-p \geq 3 & \text{DU }2 \\
\DPQ{p-1}{q}{n} & \text{when }p \geq 2 & \text{DU }3 \\
\DPQ{p}{q+1}{n} & \text{when }q+1 \leq n-1 & \text{DU }4  \\
\DPQ{p+1}{q}{n} & \text{when }q-p \geq 3 & \text{DU }5  \\
\DPQ{p+1}{q+1}{n} & \text{when }q+1 \leq n-1 & \text{DU }6 
\end{array}
\right.
\]
}
\end{newdef}
\BX
\MS

We can use the Down-Up rule to systematically derive other determinacy conditions from a given one. The derivations are based applying one of the ``DU'' rules repeatedly to $\DPQ{p}{q}{n}$, $n \geq 2$. In the following calculations, we denote the implication ``DU t'' by ``$\IMPBY{t}$''and, in each line, we assume $p$ and $q$ satisfy the requirement of the applied rule.
\[
\begin{aligned}
\DPQ{p}{q}{n} & \IMPBY{1} \DPQ{p-1}{q-1 }{n}\DDD{\IMPBY{1}} \DPQ{1}{q-p+1}{n} \\
\DPQ{p}{q}{n} & \IMPBY{6} \DPQ{p+1}{q+1}{n} \DDD{\IMPBY{6}} \DPQ{n-1-q+p}{n-1}{n} \\
\DPQ{p}{q}{n} & \IMPBY{3} \DPQ{p-1}{q}{n} \DDD{\IMPBY{3}} \DPQ{1}{q}{n} \\
\DPQ{p}{q}{n} & \IMPBY{4} \DPQ{p}{q+1}{n} \DDD{\IMPBY{4}} \DPQ{p}{n-1}{n}
\end{aligned}
\]
Also, if $q-p \geq 3$ then DU 2 and DU 5 apply:
\[
\begin{aligned}
\DPQ{p}{q}{n} & \IMPBY{2} \DPQ{p}{q-1}{n} \DDD{\IMPBY{2}} \DPQ{p}{p+2}{n} \\
\DPQ{p}{q}{n} & \IMPBY{5} \DPQ{p+1}{q}{n} \DDD{\IMPBY{5}} \DPQ{q-2}{q}{n}
\end{aligned}
\]
It will be convenient, fixing $n$, to associate the condition $\DPQ{p}{q}{n}$ with the point $(p,q)$ in the integer-point $pq$-plane. The following theorem states that the Down-Up Rule applied to a given condition $\DPQ{p}{q}{n}$ generates a maximal set of conditions (i.e. ``closed'' under applications of the Down-Up Rule) which form an integer-point {\em right triangle} in the $pq$-plane.
\MS

The following example illustrates this; it's an instance of case 5 of the next theorem.
\MS

\begin{example}
{\rm
Take $(n,i)=(8,5)$ and the determinacy condition $\DPQ{3}{6}{8}$. Consider the sequence of applications of the Down-Up Rule starting at the point $(p,q)=(3,6)$ denoted $\DOT0$.
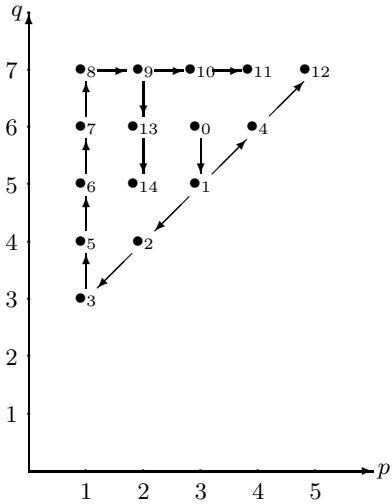
\setlength{\unitlength}{.3in}
\begin{center}
\begin{picture}(6,8)
\put(0,0){\vector(0,1){8}}
\put(-.2,8){\MBS{q}}
\put(0,0){\vector(1,0){6}}
\put(6.2,0){\MBS{p}}
\multiput(1,0)(1,0){5}{\MB{\cdot}}
\multiput(0,1)(0,1){8}{\MB{\cdot}}
\put(-.3,1){\MBS{1}}
\put(-.3,2){\MBS{2}}
\put(-.3,3){\MBS{3}}
\put(-.3,4){\MBS{4}}
\put(-.3,5){\MBS{5}}
\put(-.3,6){\MBS{6}}
\put(-.3,7){\MBS{7}}
\put(1,-.35){\MBS{1}}
\put(2,-.35){\MBS{2}}
\put(3,-.35){\MBS{3}}
\put(4,-.35){\MBS{4}}
\put(5,-.35){\MBS{5}}

\put(3,6){\MBS{\DOT0}}
\put(3,5.8){\vector(0,-1){.6}}

\put(3,5){\MBS{\DOT1}}
\put(2.8,4.8){\vector(-1,-1){.6}} 
\put(3.2,5.2){\vector(1,1){.6}} 

\put(2,4){\MBS{\DOT2}}
\put(1.8,3.8){\vector(-1,-1){.6}}

\put(1,3){\MBS{\DOT3}}
\put(1,3.2){\vector(0,1){.6}} 

\put(4,6){\MBS{\DOT4}}
\put(4.2,6.2){\vector(1,1){.6}}

\put(1,4){\MBS{\DOT5}}
\put(1,4.2){\vector(0,1){.6}} 

\put(1,5){\MBS{\DOT6}}
\put(1,5.2){\vector(0,1){.6}} 

\put(1,6){\MBS{\DOT7}}
\put(1,6.2){\vector(0,1){.6}} 

\put(1,7){\MBS{\DOT8}}
\put(1.2,7){\vector(1,0){.5}} 

\put(2,7){\MBS{\DOT9}}
\put(2.2,7){\vector(1,0){.5}} 
\put(2,6.8){\vector(0,-1){.6}} 

\put(3,7){\MBS{\DOT{{10}}}}
\put(3.2,7){\vector(1,0){.5}} 

\put(4,7){\MBS{\DOT{{11}}}}

\put(5,7){\MBS{\DOT{{12}}}}

\put(2,6){\MBS{\DOT{{13}}}}
\put(2,5.8){\vector(0,-1){.6}} 

\put(2,5){\MBS{\DOT{{14}}}}
\end{picture}
\end{center}
\setlength{\unitlength}{1in}
\MS

\NI The sequence of Down-Up rules is:
\[
\DOT0 \IMPBY{2} \DOT1  \IMPBY{1} \DOT2 \IMPBY{1} \DOT3 \IMPBY{4} \DOT5 \IMPBY{4} \DOT6
\IMPBY{4} \DOT7 \IMPBY{4} \DOT8 \IMPBY{5} \DOT9 \IMPBY{5} \DOT{{10}} \IMPBY{5} \DOT{{11}}
\]
\[
\DOT0 \IMPBY{2} \DOT{{1}} \IMPBY{6} \DOT4 \IMPBY{6} \DOT{{12}}
\]
and
\[
\DOT9 \IMPBY{2} \DOT{{13}} \IMPBY{2} \DOT{{14}}
\]
\SL
}
\end{example}
\MS

\HD{Notation: ``DetSet''} 
\index{DetSet}
\index{$\DTSET{\DPQ{p}{q}{n}}$}

Given a determinacy condition $\DPQ{p}{q}{n}$ then 
$\DTSET{\DPQ{p}{q}{n}} \DFAS$ the smallest set of determinacy conditions containing $\DPQ{p}{q}{n}$ closed under application of the Down-Up Rule.

\begin{newthm}
{\bf (Determinacy Set Theorem)}
\label{detset thm}
{\rm
\SL

\begin{enumerate}
\item 
If $n \geq 2$ then 
\[
\DTSET{\DPQ{1}{2}{n}} = \SETT{\DPQ{p}{q}{n}: 1 \leq p \leq n-2,\, p+1 \leq q \leq n-1 }
\]

\item
If $n \geq 2$ then 
\[
\DTSET{\DPQ{0}{2}{n}} = \SETT{\DPQ{p}{q}{n}: 1 \leq p \leq n-3,\, p+2 \leq q \leq n-1 }
\]

\item
If $n \geq 2$ then 
\[
\DTSET{\DPQ{n-2}{n}{n}} = \SETT{\DPQ{p}{q}{n}: 1 \leq p \leq n-2,\, p+2 \leq q \leq n }
\]

\item
If $n \geq 2$ then 
\[
\DTSET{\DPQ{n-1}{n}{n}} = \SETT{\DPQ{p}{q}{n}: 1 \leq p \leq n-1,\, p+1 \leq q \leq n }
\]

\item
If $n \geq 3$ and $1<i<n$ then 
\[
\DTSET{\DPQ{i-2}{i+1}{n}} = \SETT{\DPQ{p}{q}{n}: 1 \leq p \leq n-3,\, p+2 \leq q \leq n-1 }
\]
\end{enumerate}
}
\end{newthm}

\Proof

Each of the claimed generated sets of determinacy conditions forms a integer-point right triangle in the $pq$-plane with a vertical edge and a horizontal edge and a hypotenuse of slope 1. It will be convenient to refer to the conditions in the generated sets in those terms.
\MS

Note also that no Down-Up rule applied to any condition on the boundary of the claimed triangle yields a condition outside the triangle. Therefore it suffices to show that every condition in the triangle is deriveable from the generating condition and use the Down-Up rules.
\MS

The proof of claim 1 will serve as a model for the proofs of claims 2--5.
 
\begin{enumerate}
\item 

Rule DU 4 applied to $\DPQ{1}{2}{n}$:
\[
\DPQ{1}{2}{n} \DDD{\IMPBY{4}} \DPQ{1}{n-1}{n}
\]
generates the vertical edge of the claimed triangle.

Rule DU 6 applied to $\DPQ{1}{2}{n}$:
\[
\DPQ{1}{2}{n} \DDD{\IMPBY{6}} \DPQ{n-2}{n-1}{n}
\]
generates the hypotenuse of the claimed triangle.

Rule DU 5 applied to $\DPQ{1}{n-1}{n}$:
\[
\DPQ{1}{n-1}{n} \DDD{\IMPBY{5}} \DPQ{n-3}{n-1}{n}
\]
generates the horizontal edge of the claimed triangle.

For $p=2\DDD{,} n-3$, we have
\[
\DPQ{p}{n-1}{n} \IMPBY{2} \DPQ{p}{n-2}{n} \DDD{\IMPBY{2}} \DPQ{p}{p+2}{n}
\]
generating the conditions in the interior of the claimed triangle.

\item 
\[
\begin{aligned}
 \DPQ{0}{2}{n} & \DDD{\IMPBY{4}} \DPQ{0}{n-1}{n} & \text{ (vert edge)}\\
 \DPQ{0}{2}{n} &\DDD{\IMPBY{6}} \DPQ{n-3}{n-1}{n} & \text{ (hypot.)} \\
 \DPQ{0}{n-1}{n} &\DDD{\IMPBY{5}} \DPQ{n-4}{n-1}{n} & \text{ (horiz. edge)}
\end{aligned}
\]
The interior conditions are generated by DU 2 applied to the horizontal edge conditions $\DPQ{p}{n-1}{n}$, $1 \leq p \leq n-4$. 

\item 
\[
\begin{aligned}
\DPQ{n-2}{n}{n} & \DDD{\IMPBY{1}} \DPQ{1}{3}{n} & \text{ (hypot.)}\\
\DPQ{1}{3}{n} & \DDD{\IMPBY{4}} \DPQ{1}{n-1}{n} & \text{ (vert. edge)}\\
\DPQ{1}{n-1}{n} & \DDD{\IMPBY{5}} \DPQ{n-2}{n}{n} & \text{ (horiz. edge)}
\end{aligned}
\]
The interior conditions are generated by DU 2 applied to the horizontal edge conditions $\DPQ{p}{n-1}{n}$, $2 \leq p \leq n-3$. 

\item 
\[
\begin{aligned}
\DPQ{n-1}{n}{n} & \DDD{\IMPBY{1}} \DPQ{1}{2}{n} & \text{ (hypot.)}\\
\DPQ{n-1}{n}{n} & \DDD{\IMPBY{3}} \DPQ{1}{n}{n} & \text{ (horiz. edge)}\\
\DPQ{1}{2}{n} & \DDD{\IMPBY{4}} \DPQ{1}{n-1}{n} & \text{ (vert. edge)}
\end{aligned}
\]
The interior conditions are generated by DU 2 applied to the horizontal edge conditions $\DPQ{p}{n-1}{n}$, $2 \leq p \leq n-2$. 

\item 
\[
\begin{aligned}
\DPQ{i-2}{i+1}{n} & \IMPBY{2} \DPQ{i-2}{i}{n} & \\
\DPQ{i-2}{i}{n} & \DDD{\IMPBY{1}} \DPQ{1}{3}{n} & \text{ (lower hypot.)}\\
\DPQ{i-2}{i}{n} & \DDD{\IMPBY{6}} \DPQ{n-3}{n-1}{n} & \text{ (upper hypot.)}\\
\DPQ{1}{3}{n} & \DDD{\IMPBY{4}} \DPQ{1}{n-1}{n} & \text{ (vert. edge)}\\
\DPQ{1}{n-1}{n} & \DDD{\IMPBY{5}} \DPQ{n-3}{n-1}{n} & \text{ (horiz. edge)}
\end{aligned}
\]
The interior conditions are generated by DU 2 applied to the horizontal edge conditions $\DPQ{p}{n-1}{n}$, $2 \leq p \leq n-4$. 

\end{enumerate}

\qed
\MS

\BS

\HD{Surjectivity}
\MS

Repeated us of the Downward Surjectivity Corollary \refpage{downward-surj-corollary} applied to a subcomplex $\MC{S}'$ in which all $n$-simplices ($n \geq 2$) are \SBSI\ establishes lower-dimensional surjectivity conditions.
\MS

We will abbreviate ``\SUR{i} in dimension $n$'' by $\SURD{i}{n}$ and, where convenient, identify $\SURD{i}{n}$ with the integer point $(i,n)$in the $in$-plane. In what follows below, we will assume $n \geq 2$, $i \in [n]$ and $\MC{S}'$ is a complex with the property that all $n$-simplices are \SBSI\ and $\SURD{i}{n}$.
\MS

For $\MC{S}'$, the corollary may be restated as a ``rule'':

\begin{newdef}
\label{down surj rule}
\index{Downward surjectivity rule}
{\bf (Downward Surjectivity Rule)}
{\rm
\SL

Suppose $n \geq 2$ and $\MC{S}'$ is a subcomplex of $\MC{S}$ such that all $n$-simplices of $\MC{S}'$ are \SBSI\ and \SUR{i}. Then all $(n-1)$-simplices of $\MC{S}'$ are \SUR{(i-1)} if $i>1$ and \SUR{i} if $i<n-1$. That is:
\[
\SURD{i}{n} \IMP 
\left\{
\begin{array}{llr}
\SURD{(i-1)}{(n-1)} & \text{if } i>1 & \text{S1}\\
\SURD{i}{(n-1)} & \text{if } i<n-1 & \text{S2}
\end{array}
\right.
\]
}
\end{newdef}
\BX
\MS

Repeated use of this rule shows how $\SURD{i}{n}$ propagates down in dimension.
\MS

If $i>1$ then S1 gives:
\[
\SURD{i}{n} \IMPBY{\text{S}1} \SURD{(i-1)}{(n-1)} \DDD{\IMPBY{\text{S}1}} \SURD{1}{(n-i+1)}
\]
and if $i<n-1$ then S2 gives:
\[
\SURD{i}{n} \IMPBY{\text{S}2} \SURD{i}{(n-1)} \DDD{\IMPBY{\text{S}2}} \SURD{i}{(i+1)}
\]
Associating the condition $\SURD{i}{n}$ with the point $(i,n)$ define
\[
\begin{aligned}
\text{S1}(i,n) & \DFAS \SETT{(i-t,n-t): 0 \leq t \leq i-1}\\
\text{S2}(i,n) & \DFAS \SETT{(i,n-t): 0 \leq t \leq n-1-i}
\end{aligned}
\]
In the $in$-plane, $\text{S1}(i,n)$ consists of the set of (points representing) surjectivity conditions implied by $\SURD{i}{n}$ by rule S1. These lie on the line from $(1,n-i+1)$ to $(i,n)$.
Similarly, $\text{S2}(i,n)$ consists of the set of surjectivity conditions implied by $\SURD{i}{n}$ by rule S2 and consists of the points on the vertical line from $(i,n)$ to $(i,i+1)$.
\MS

Note that 
\begin{enumerate}
\item 
$\text{S1}(0,n) =\MT$

\item 
$\text{S1}(1,n) = \{ (1,n)\}$

\item 
$\text{S2}(n,n)=\MT$

\item
If $1<i<n-1$ then 
\[
\text{S1}(\text{S2}(i,n)) = \text{S2}(\text{S1}(i,n))
\]
because
\[
\SURD{i}{n} \IMPBY{S2} \SURD{i}{(n-1)} \IMPBY{S1} \SURD{(i-1)}{(n-2)}
\]
and
\[
\SURD{i}{n} \IMPBY{S1} \SURD{(i-1)}{(n-1)} \IMPBY{S2} \SURD{(i-1)}{(n-2)}
\]

\end{enumerate}
\MS

\HD{Notation: ``SurjSet''}
\index{SurjSet}
\index{$\SRSET{\SURD{i}{n}}$}

Given $(i,n)$ where $0 \leq i \leq n$ then $\SRSET{\SURD{i}{n}} \DFAS $ the smallest set of surjectivity conditions containing $\SURD{i}{n}$ closed under the Downward Surjectivity Rule. 

\begin{newthm}
\label{surjset thm}
{\rm
\SL

Let $n \geq 2$ and $i \in [n]$. Then
\[
\SRSET{\SURD{i}{n}}= \bigcup_{(j,m) \in \text{S1}(i,n)} \text{S2}(j,m) \quad \cup
\bigcup_{(j,m) \in \text{S2}(i,n)} \text{S1}(j,m)
\]
}
\end{newthm}

\Proof

Clearly, the union is a subset of $\SRSET{\SURD{i}{n}}$. On the other hand, given any $\SURD{j}{m} \in \SRSET{\SURD{i}{n}}$, then, by definition, it is derived by some combination of rules S1 and S2 and therefore belongs to at least one of $\bigcup_{(j,m) \in \text{S1}(i,n)} \text{S2}(j,m)$ and $\bigcup_{(j,m) \in \text{S2}(i,n)} \text{S1}(j,m)$.

\qed
\MS

\NI The points of $\SRSET{\SURD{i}{n}}$ consists of either a line or a parallelogram:

\begin{itemize}
\item 
$\SRSET{\SURD{0}{n}} = \text{S2}(0,n)$ is the vertical line from $(0,n)$ to $(0,1)$.

\item 
$\SRSET{\SURD{1}{n}} = \text{S2}(1,n)$ is the vertical line from $(1,n)$ to $(1,2)$.

\item 
$\SRSET{\SURD{n}{n}} = \text{S1}(n,n)$ is the line from $(n,n)$ to $(1,1)$.

\item 
$\SRSET{\SURD{n-1}{n}} = \text{S1}(n-1,n)$ is the line from $(n-1,n)$ to $(1,2)$.

\item 
For $1<i<n-1$, $\SRSET{\SURD{i}{n}}$ is the parallelogram with vertices $(i,n)$, $(i,i+1)$, $(1,n-i+1)$ and $(1,2)$. That is:
\[
\SRSET{\SURD{i}{n}} = \text{S2}(i,n) \cup \text{S2}(i-1,n-1) \DDD{\cup} \text{S2}(1,n-i+1)
\]

\end{itemize}

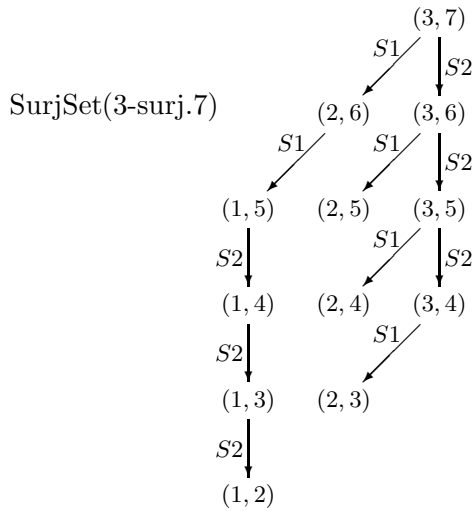
\begin{example}
{\rm
\setlength{\unitlength}{.5in}

\begin{center}
\begin{picture}(3,4.5)
\put(0,0){\MBS{(1,2)}}
\put(0,1){\MBS{(1,3)}}
\put(0,2){\MBS{(1,4)}}
\put(0,3){\MBS{(1,5)}}
\put(1,1){\MBS{(2,3)}}
\put(1,2){\MBS{(2,4)}}
\put(1,3){\MBS{(2,5)}}
\put(1,4){\MBS{(2,6)}}
\put(2,2){\MBS{(3,4)}}
\put(2,3){\MBS{(3,5)}}
\put(2,4){\MBS{(3,6)}}
\put(2,5){\MBS{(3,7)}}
\put(2,4.8){\vector(0,-1){.6}}
\put(2.2,4.5){\MBS{S2}}
\put(2,3.8){\vector(0,-1){.6}}
\put(2.2,3.5){\MBS{S2}}
\put(2,2.8){\vector(0,-1){.6}}
\put(2.2,2.5){\MBS{S2}}
\put(0,2.8){\vector(0,-1){.6}}
\put(-.2,2.5){\MBS{S2}}
\put(0,1.8){\vector(0,-1){.6}}
\put(-.2,1.5){\MBS{S2}}
\put(0,.8){\vector(0,-1){.6}}
\put(-.2,.5){\MBS{S2}}
\put(1.8,4.8){\vector(-1,-1){.6}}
\put(1.45,4.7){\MBS{S1}}
\put(0.45,3.7){\MBS{S1}}
\put(1.45,3.7){\MBS{S1}}
\put(1.45,2.7){\MBS{S1}}
\put(1.45,1.7){\MBS{S1}}
\put(0.8,3.8){\vector(-1,-1){.6}}
\put(1.8,3.8){\vector(-1,-1){.6}}
\put(1.8,2.8){\vector(-1,-1){.6}}
\put(1.8,1.8){\vector(-1,-1){.6}}
\put(-2.5,4){\SRSET{\SURD{3}{7}}}
\end{picture}
\end{center}
\MS

That is, 
\[
\SRSET{\SURD{3}{7}} = \text{S2}(3,7) \cup \text{S2}(2,6) \cup \text{S2}(1,5)
\]
}
\end{example}
\BS

\HD{Implications for the models}

The models developed in section \refpage{nisolution} cited sufficient conditions for a subcomplex to be an \NICOMP{n}{i}. From the discussion of determinacy and surjectivity just above, we can elaborate on the properties of those models using theorems \refpage{detset thm} and \refpage{surjset thm}.
\BS

\HD{\NICOMP{n}{0} model}
\MS

The requirements in Theorem \refpage{composer-model;i=0} were $\SURD{0}{n}$, $\SURD{1}{n}$ and $\DPQ{1}{2}{n}$. By Theorem \refpage{surjset thm} the set of surjectivity conditions is 
\[
\begin{aligned}
\SRSET{\SURD{0}{n}} \cup \SRSET{\SURD{1}{n}} = \\
\SETT{\SURD{0}{m} : 1 \leq m \leq n} \cup \SETT{\SURD{1}{m} : 2 \leq m \leq n}
\end{aligned}
\]
By Theorem \refpage{detset thm} the set of determinacy conditions is
\[
\DTSET{\DPQ{1}{2}{n}} = \SETT{\DPQ{p}{q}{n}: 1 \leq p \leq n-2, p+1 \leq 1 \leq n-1}
\]
\MS

\HD{\NICOMP{n}{1} model}
\MS

The requirements in Theorem \refpage{composer-model;i=1} were $\SURD{1}{n}$ and $\DPQ{0}{2}{n}$. Then:
\[
\SRSET{\SURD{1}{n}} = \SETT{\SURD{1}{m}: 2 \leq m \leq n}
\]
and
\[
\DTSET{\DPQ{0}{2}{n}} = \SETT{\DPQ{p}{q}{n}: 1 \leq p \leq n-3,\, p+2 \leq q \leq n-1 }
\]
\MS

\HD{\NICOMP{n}{i} model, $1<i<n-1$}
\MS

The requirements in Theorem \refpage{composer-model;1<i<n} were $\SURD{i}{n}$, $\SURD{(i-1)}{n}$ and $\DPQ{i-2}{i+1}{n}$.
\[
\SRSET{\SURD{i}{n}} = \bigcup_{t=0}^{i-1} \text{S2}(i-t,n-t)
\]
and
\[
\SRSET{\SURD{i-1}{n}} = \bigcup_{t=0}^{i-2} \text{S2}(i-1-t,n-t)
\]
Therefore
\[
\begin{aligned}
\SRSET{\SURD{i}{n}} \cup \SRSET{\SURD{i-1}{n}} = \\
\SRSET{\SURD{i}{n}} \cup \SETT{\SURD{(i-1-t)}{(n-t)}: 0 \leq t \leq i-2}
\end{aligned}
\]
And
\[
\DTSET{\DPQ{i-2}{i+1}{n}} = \SETT{\DPQ{p}{q}{n}: 1 \leq p \leq n-3,\, p+2 \leq q \leq n-1 }
\]
\MS

\HD{\NICOMP{n}{n} model}
\MS

The requirements in Theorem \refpage{composer-model;i=n} were $\SURD{(n-1)}{n}$ and $\DPQ{n-2}{n}{n}$. Then:
\[
\SRSET{\SURD{(n-1)}{n}} = \text{S1}(n-1,n) = \SETT{\SURD{(n-t-1)}{(n-t):0 \leq t \leq n-2}}
\]
and
\[
\DTSET{\DPQ{n-2}{n}{n}} = \SETT{\DPQ{p}{q}{n}: 1 \leq p \leq n-2,\, p+2 \leq q \leq n }
\]
\MS

\HD{\NICOMP{n}{n+1} model}
\MS

The requirements in Theorem \refpage{composer-model;i=n+1} were $\SURD{(n-1)}{n}$, $\SURD{n}{n}$, $\DPQ{n-1}{n}{n}$ and $\DPQ{n-2}{n-1}{n}$. Note that $\DPQ{n-1}{n}{n} \IMP \DPQ{n-2}{n-1}{n}$ by the Down-Up rule DU 1. Therefore the requirement of $\DPQ{n-2}{n-1}{n}$ was redundant.

The implied surjectivity conditions are
\[
\SRSET{\SURD{(n-1)}{n}} \cup \SRSET{\SURD{n}{n}}
\]
and the implied determinacy conditions are
\[
\DTSET{\DPQ{n-1}{n}{n}} = \SETT{\DPQ{p}{q}{n}: 1 \leq p \leq n-1,\, p+1 \leq q \leq n }
\]

\subsection{Succinct characterizations of \NICOMP{n}{i}s of sets}

The gist of each of the following theorems is that an \NICOMP{n}{i} of sets $\MC{S}'$ can be characterized by properties in dimension $n+1$. The theorems treat individual cases according to the value of $i \in [n+1]$.

In part \refpage{function-case}, specializing to $n=1$ and $i=1$, there is an analogous result characterizing the nerve of a category of sets and functions.

\subsubsection{Case $1<i<n-1$}

\begin{newthm}
{\rm
\SL

Suppose $1<i<n-1$. Suppose $\MC{S}'_{\leq n+1}$ is an $(n+1)$-truncated subcomplex of $\mathcal{S}$ such that: 
\begin{enumerate}
\item All simplices of $\MC{S}'_{\leq n+1}$ are \SBSI.
\item All $w \in \MC{S}'_{n+1}$ are \SUR{i}.
\item All $w \in \MC{S}'_{n+1}$ are \DET{\SETT{i-2,i+1}}.
\end{enumerate}
Then $\MC{S}'_{\leq n+1}$ is the $(n+1)$-truncation of an \NICOMP{n}{i} of sets.
}
\end{newthm}

\Proof

Down Rule (c), which applies to $\DPQ{i-2}{i+1}{(n+1)}$ because $i+1<n$, yields $\DPQ{i-2}{i+1}{n}$. Every $w \in \MC{S}'_{n+1}$ is \SUR{i} and therefore every $y \in \MC{S}'_n$ is both \SUR{i} and \SUR{(i-1)} by the Downward Surjectivity Corollary (page \pageref{downward-surj-corollary}). Therefore, by Theorem \refpage{composer-model;1<i<n}, $\MC{S}'_{\leq n+1}$ is the $(n+1)$-truncation of an \NICOMP{n}{i} of sets.

Note that the three assumed conditions imply that all $w \in \MC{S}'_{n+1}$ are $(n,i)$-compositions. (Lemma \refpage{comp-determinacy-lemma}).

\qed

\subsubsection{Case $n \geq 3$ and $i=n-1$}

\begin{newthm}
{\rm
\SL

Suppose $n \geq 3$. Suppose $\MC{S}'_{\leq n+1}$ is an $(n+1)$-truncated subcomplex of $\mathcal{S}$ such that
\begin{enumerate}
\item All simplices of $S'_{\leq n+1}$ are \SBSI.

\item All $w \in \MC{S}'_{n+1}$ are \SUR{(n-1)}.

\item All $w \in \MC{S}'_{n+1}$ are \DET{\SETT{n-3,n+1}}.
\end{enumerate}
Then $\MC{S}'_{\leq n+1}$ is the $(n+1)$-truncation of an \NICOMP{n}{n-1} of sets.
}
\end{newthm}

\Proof

Condition (2) implies all $y \in \MC{S}'_n$ are \SUR{(n-1)} and \SUR{(n-2)}. (Downward Surjectivity Corollary, page \pageref{downward-surj-corollary})

Condition (1) implies all $y \in \MC{S}'_n$ are \SBSI.

Down Rule (b) applies to $\DPQ{n-3}{n+1}{(n+1)}$ to imply $\DPQ{n-3}{n}{n}$.

It then follows from Theorem \refpage{composer-model;1<i<n} that $\MC{S}'_{\leq n+1}$ is the $(n+1)$-truncation of an \NICOMP{n}{n-1} of sets.

\qed

\subsubsection{Case $n \geq 2$ and $i=0$}

\begin{newthm}
{\rm
\SL

Suppose $n \geq 2$. Suppose $\MC{S}'_{\leq n+1}$ is an $(n+1)$-truncated subcomplex of $\mathcal{S}$ such that
\begin{enumerate}
\item All simplices of $\MC{S}'_{\leq n+1}$ are \SBSI.

\item All $w \in \MC{S}'_{n+1}$ are \SUR{0}.

\item All $w \in \MC{S}'_{n+1}$ are \SUR{1}.

\item All $w \in \MC{S}'_{n+1}$ are \DET{\SETT{2,3}}.
\end{enumerate}
Then $\MC{S}'_{\leq n+1}$ is the $(n+1)$-truncation of an \NICOMP{n}{0} of sets.
}
\end{newthm}

\Proof

Condition (1) implies all $y \in \MC{S}'_n$ are \SBSI.

The Downwards Surjectivity corollary and conditions (2) and (3) imply that all $y \in \MC{S}'_n$ are \SUR{0} and \SUR{1}.

Downrule (a) and condition (4) imply all $y \in \MC{S}'_n$ are \DET{\{1,2\}}.

Corollary \refpage{composer-model;i=0} implies that $\MC{S}'_{\leq n+1}$ is the $(n+1)$-truncation of an \NICOMP{n}{0} of sets, as claimed.

\qed

\subsubsection{Case $n \geq 2$ and $i=1$}

\begin{newthm}
{\rm
\SL

Suppose $n > 2$. Suppose $\MC{S}'_{\leq n+1}$ is an $(n+1)$-truncated subcomplex of $\mathcal{S}$ such that
\begin{enumerate}
\item All simplices of $S'_{\leq n+1}$ are \SBSI.

\item All $w \in \MC{S}'_{n+1}$ are \SUR{1}.

\item All $w \in \MC{S}'_{n+1}$ are \DET{\SETT{1,2}}.
\end{enumerate}
Then $\MC{S}'_{\leq n+1}$ is the $(n+1)$-truncation of an \NICOMP{n}{1} of sets.
}
\end{newthm}

\Proof

Condition (1) implies all $y \in \MC{S}'_n$ are \SBSI.

Condition (2) implies all $y \in \MC{S}'_n$ are \SUR{1}.

Down Rule (c), that $n>2$ and $\DPQ{0}{2}{(n+1)}$ imply $\DPQ{0}{2}{n}$. 

Therefore, by  Corollary \refpage{n-1-composer-theorem}, $\MC{S}'_{\leq n+1}$ is the $(n+1)$-truncation of an \NICOMP{n}{1} of sets, as claimed.

\qed

\subsubsection{Case $n \geq 2$ and $i=n$}

\begin{newthm}
{\rm
\SL

Suppose $n \geq 2$.  Suppose $\MC{S}'_{\leq n+1}$ is an $(n+1)$-truncated subcomplex of $\mathcal{S}$ such that
\begin{enumerate}
\item All simplices of $\MC{S}'_{\leq n+1}$ are \SBSI.

\item All $w \in \MC{S}'_{n+1}$ are \SUR{n}.

\item All $w \in \MC{S}'_{n+1}$ are \DET{\SETT{n-2,n+1}}.
\end{enumerate}
Then $\MC{S}'_{\leq n+1}$ is the $(n+1)$-truncation of an \NICOMP{n}{n} of sets.
}
\end{newthm}

\Proof

Condition (1) implies all $y \in \MC{S}'_n$ are \SBSI.

Condition (2) and that $n>1$ imply all $y \in \MC{S}'_n$ are \SUR{(n-1)}.

Condition (3) and Down Rule (b) imply $\DPQ{n-2}{n}{n}$. 

Therefore, by Corollary \refpage{n-n-composer-theorem}, $S'_{\leq n+1}$ is the $(n+1)$-truncation of an \NICOMP{n}{n} of sets, as claimed.

\qed

\subsubsection{Case $n \geq 3$ and $i=n+1$}

\begin{newthm}
{\rm
\SL

Suppose $n \geq 3$ and $\MC{S}'_{\leq n+1}$ is an $(n+1)$-truncated subcomplex of $\mathcal{S}$ such that
\begin{itemize}
\item 
All simplices of $\MC{S}'_{\leq n+1}$ are \SBSI.

\item
All $w \in S'_{n+1}$ are \SUR{n} and \SUR{(n+1)}.


\item
All $w \in S'_{n+1}$ are \DET{\SETT{n-1,n}}.
\end{itemize}
Then
\begin{enumerate}
\item 
All $w \in \MC{S}'_{n+1}$ are $(n,n+1)$-compositions.

\item 
All $y \in \MC{S}'_n$ are \SUR{n} and \SUR{(n-1)}.

\item
All $y \in \MC{S}'_n$ are \DET{\SETT{n-1,n}} and \DET{\SETT{n-2,n-1}}.

\item
$\MC{S}'_{\leq n+1}$ is the $(n+1)$-truncation of an \NICOMP{n}{n+1} of sets.
\end{enumerate}
}
\end{newthm}

\Proof

\begin{enumerate}
\item %
$w$ is \SUR{(n+1)}. Since $w$ is $\DPQ{n-1}{n}{(n+1)}$ and $n+1 \notin \{n-1,n\}$ then lemma \refpage{comp-determinacy-lemma} implies $w$ is an $(n,n+1)$-composition.

\item %
Given any $y \in \MC{S}'_n$ then $s_0(y)$ is, by hypothesis, both \SUR{n} and \SUR{(n+1)}. Since $n>2$ then Downwards Surjectivity (Corollary \refpage{downward-surj-corollary}) implies $y$ is both \SUR{n} and \SUR{(n-1)}.

\item %
Given any $y \in \MC{S}'_n$ then $s_0(y)$ is, by hypothesis, $\DPQ{n-1}{n}{(n+1)}$. Since $n>2$ then the Down Rule (page \pageref{down rule}) D1 implies $y$ is $\DPQ{n-1}{n}{n}$. It then follows from the Determinacy Set Theorem (page \pageref{detset thm}) that $y$ is also $\DPQ{n-2}{n-1}{n}$.

\item %
This follows from the previous items and Theorem \refpage{composer-model;i=n+1}.

\end{enumerate}

\qed

\subsubsection{Case $n=1$ and $i=1$: composition of functions}
\label{function-case}

In this section we will characterize the nerve $\MC{S}'$ of a category of sets and functions in terms of properties in dimension 2.
\MS

As observed earlier, a $(1,1)$-composition in $\mathcal{S}_2$ corresponds to ordinary composition of binary relations. If $w \in \mathcal{S}_2$ is a $(1,1)$-composition, then, by definition, $w$ is e-simplicial, $w$ is \DET{\SETT{0,2}}, and $(a_0,a_1) \in d_1(w)$ if and only if there exists $b \in w$ such that $e_1(b) = (a_0,a_1)$. A typical element of $w$ is
\[
b = 
\left[
\begin{array}{c}
b_0 \\
b_1 \\
b_2
\end{array}
\right]
=
\left[
\begin{array}{cc}
b_{00} & b_{01}\\
b_{00} & b_{21}\\
b_{01} & b_{21}
\end{array}
\right] =
\left[
\begin{array}{cc}
b_{00} & b_{01}\\
a_0 & a_1\\
b_{01} & b_{21}
\end{array}
\right]
\]
We note that in dimension 2, the properties e-simplicial and \SBSI\ coincide.

The two assumptions about $w$, namely that it is e-simplicial and \DET{\SETT{0,2}}, do not imply any properties for the faces of of $w$. In particular, the faces of $w$ (the factors and the resulting composite) need not be functions.
\MS

Now suppose $w $ is a $(1,1)$-composition whose factors {\em are} functions. Let $y_j = d_j(w)$, $j=0,1,2$. Then, in conventional function notation, $y_1 = y_0 \circ y_2$. For each
\[
\left[
\begin{array}{cc}
b_{00} & b_{01}\\
b_{00} & b_{21}\\
b_{01} & b_{21}
\end{array}
\right] \in w
\]
we have $y_2(b_{21}) = b_{01}$, $y_0(b_{01}) = b_{00}$ and $y_1(b_{21}) = y_0(y_2(b_{21})) = y_0(b_{01}) = b_{00}$. It follows that the entire array for $b \in w$ is determined by $e_1 e_2(w) = b_{21}$. This implies that $w$ is \DET{\SETT{2}}. By the same reasoning, it also implies $w$ is \SUR{2}.
\MS

Another note: in the definition of $A$-determinacy in dimension $m$ (page \pageref{determinacy}) the set $A \subsetneq [m]$ was required to have at least two elements. Here, we are stretching the terminology to allow for $A = \SETT{2}$.
\MS

Thus, if $w$ is a $(1,1)$-composition of {\em functions}, then $w$ is \DET{\SETT{2}} and \SUR{2}. In fact, these assumptions suffice to characterize \NICOMP{1}{1} of sets in which the 1-simplices are functions.

\begin{newthm}
{\rm
\SL

Suppose $\MC{S}'_{\leq 2}$ is a 2-truncated subcomplex of $\mathcal{S}$ such that 
\begin{enumerate}
\item All $w \in \MC{S}'_2$ are e-simplicial.
\item All $w \in \MC{S}'_2$ are \SUR{1}.
\item All $w \in \MC{S}'_2$ are \DET{\SETT{2}}.
\item All $w \in \MC{S}'_2$ are \SUR{2}.
\end{enumerate}
Then all $w \in \MC{S}'_2$ are $(1,1)$-compositions and all $y \in \MC{S}'_1$ are functions.
}
\end{newthm}

\Proof

That all $w \in \MC{S}'_2$ are \DET{\SETT{2}} implies that all $w \in \MC{S}'_2$ are \DET{\SETT{0,2}}. Therefore, properties (2) and (3) imply that $(a_0,a_1) \in d_1(w)$ if and only if there exists $b \in w$ such that $e_1(b) = (a_0,a_1)$ and where $b$ has the form
\[
b = 
\left[
\begin{array}{cc}
b_{00} & b_{01}\\
a_0 & a_1\\
b_{01} & b_{21}
\end{array}
\right]
\]
with $a_0=b_{00}$ and $a_1 = b_{21}$. That is, $w$ is a $(1,1)$-composition.
\MS

Next, suppose $y \in \MC{S}'_1$. Then $s_0(y) \in \MC{S}'_2$ and a typical element of $s_0(y)$ is
\[
c_0(a) = c_0(a_0,a_1) = 
\left[
\begin{array}{cc}
a_0 & a_1 \\
a_0 & a_1 \\
a_1 & a_1 
\end{array}
\right]
\]
That $s_0(y)$ is e-simplicial, \SUR{2} and \DET{\SETT{2}} implies that for every $a_1 \in d_1(y)$, there exists one (by 2-surjectivity) and only one (by $\SETT{2}$-determinacy) $a_0 \in d_0(y)$  such that $(a_0,a_1) \in y$. That is, $y$ is a function.
\qed

\section{Constructing truncations of models}
\label{model-construction}

\subsection{Truncated subcomplexes of $\MC{S}$}

The models of \NICOMP{n}{i}s described in section \refpage{nisolution} were based on certain $n$-truncated subcomplexes of $\MC{S}$. In this section we will show that such subcomplexes occur by developing some methods to produce examples of them.
\MS

Given a non-empty relation $R \SBS \PLIST{V}{n}$, we may associate to $R$ a specific \SBSI\ $n$-simplex denoted $y^R$ (defined below) from which $R$ may be recovered. The goal is to find enlargements $\bar{R}$ of $R$ such that the $n$-truncated subcomplex generated by $y^{\bar{R}}$ (i.e. the one consisting of all simplicial images of $y^{\bar{R}}$) satisfies the conditions for it to be the $n$-truncation of an \NICOMP{n}{i}. 

We denote the $n$-truncated subcomplex generated by $y^{\bar{R}}$ by $\MC{S}[y^{\bar{R}}]$.
\MS

Since the desired $\MC{S}[y^{\bar{R}}]$ must satisfy one or more determinacy conditions in dimension $n$, the Down Rule (Theorem \refpage{down rule}) implies certain determinacy conditions in dimensions below $n$. Therefore, to show the existence of $\bar{R}$, we must take into account a set of determinacy conditions.

\subsubsection{Minimal simplex of a relation}
\label{min-simplex-of-rel}

We begin by recalling a general observation concerning arbitrary $y \in \MC{S}_n$, $n \geq 1$. We use $\ELEM{y}$ to denote the set of elements of $y$, and assume here that $\ELEM{y} \neq \MT$.
\MS

\begin{newdef}
\label{general min simplex}
\index{minimal simplex of a simplex}
{\bf (Minimal simplex of a simplex)}
{\rm
\SL

Given $y \in \MC{S}_n$, $n \geq 1$ we define $\MIN{y} \in \MC{S}_n$ by: 
\[
\ELEM{\MIN{y}} = \ELEM{y}
\]
and, for all $k \leq n$ and all $p_1 \DDD{<} p_k$ subface permissible for $n$
\[
\ELEM{d_{p_1} \DDD{} d_{p_k}(\MIN{y})} \DFAS \SETT{e_{p_1} \DDD{} e_{p_k}(a): a \in \ELEM{y}}
\]
We say that $y \in \MC{S}_n$ is {\bf minimal} if $y = \MIN{y}$.
}
\BX
\end{newdef}

\NI Several consequences of this definition:
\begin{enumerate}

\item
A minimal $n$-simplex is \SUR{p} for all $p \in [n]$ and all its $m$-dimensional subfaces are \SUR{q} for all $q \in [m]$.

\item 
$\Sig(\MIN{y}) = \prod_{p=0}^n d_p(\MIN{y})$ where $d_p(\MIN{y}) \SBS d_p(y)$. If $y$ is not \SUR{p} then $d_p(\MIN{y}) \subsetneq d_p(y)$.

\item
For all $p \in [n]$, $d_p(\MIN{y}) \SBS \MIN{d_p(y)}$ and, if $y$ is not \SUR{p} then $d_p(\MIN{y}) \subsetneq \MIN{d_p(y)}$.

\item
All subfaces of a minimal simplex are minimal.

\end{enumerate}
\BX
\MS

Now we may specialize the minimal-simplex idea to the case when $y \in \MC{S}_n$ is \SBSI. To be specific, suppose $y$ has fundamental sets $V_0 \DDD{,} V_n$ and that its vertex-relation is $R=R(y) \SBS \PLIST{V}{n}$ where, by definition, $\LIST{t}{n} \in R \iff h_n \LIST{t}{n} \in \ELEM{y}$. Then the vertex-relation of $\MIN{y}$ is $R(y)$ and the vertex-relation of $d_p(\MIN{y})$ is $R(d_p(\MIN{y}))$ where
\[
(\OM{t}{n}{p}) \in R(d_p(\MIN{y})) \SBS V_0 \DDD{} \omit{V_p} \DDD{} V_n
\]
iff \;
$\exists t_p \in V_p$ such that $(t_0 \DDD{,} t_p \DDD{,} t_n) \in R(y)$.
\BS

The passage from \SBSI\ simplices to vertex-relations is not in general reversible since different \SBSI\ $n$-simplices may have the same vertex-relation. However, given a non-empty $R \SBS V_0 \DDD{}V_n$ there is a unique {\em minimal} $n$-simplex whose vertex-relation is $R$.

\begin{newdef}
\label{minimal-simplex}
\index{minimal simplex of a relation}
{\bf (Minimal simplex of a relation)}
{\rm
\SL

 Let $n \geq 0$, $R \SBS V_0 \DDD{} V_n$, and $R \neq \MT$. Then define the \SBSI\ $n$-simplex denoted $y^R$ with fundamental sets $\PLIST{V}{n}$ by
\[
\ELEM{y^R} \DFAS \SETT{h_n \LIST{t}{n} : \LIST{t}{n} \in R}
\]
and for each $p_1 \DDD{<} p_k$ subface-permissible for $n$,
\[
\ELEM{d_{p_1} \DDD{} d_{p_k}(y^R)} \DFAS 
\SETT{e_{p_1} \DDD{} e_{p_k} h_n \LIST{t}{n}: \LIST{t}{n} \in R}
\]
We call $y^R$ the {\bf minimal simplex of $R$}.
}

\BX
\end{newdef}

\NI That $y^R$ and all its subfaces are minimal follows from the definition.
\BS

The simplex structure of $y^R$ translates to one for $R$, as follows.

\begin{newdef}
{\bf (Simplex structure of a relation)}
\label{simp structure of a relation}
\index{simplex structure of a relation}
{\rm
\SL

Given $n \geq 1$ and $R \SBS \PLIST{V}{n}$ as above, the simplex structure of $y^R$ induces a simplex structure on $R$, as follows. 
\[
\begin{aligned}
d_p(R) & \DFAS \text{ the vertex-relation of }d_p(y^R) \\
s_p(R) & \DFAS \text{ the vertex-relation of } s_p(y^R)
\end{aligned}
\]
Note then that minimality implies $d_p(y^R) = y^{d_p(R)}$ and $s_p(y^R) = y^{s_p(R)}$
}

\BX
\end{newdef}

\NI According to this:
\[
\begin{aligned}
(\OM{t}{n}{p}) \in d_p(R) & \iff h_{n-1}(\OM{t}{n}{p}) \in d_p(y^R)\\
& \iff \exists t_p \; \Bigl( h_n(t_0 \DDD{,} t_p \DDD{,} t_n) \in y^R \Bigr) \\
&\iff \exists t_p\; \Bigl( (t_0 \DDD{,} t_p \DDD{,} t_n) \in R \Bigr)
\end{aligned}
\]
and
\[
s_p(R) = \SETT{(t_0 \DDD{,} t_p,t_p \DDD{,} t_n) : \LIST{t}{n} \in R}
\]
\BX
\MS

Suppose $y \in \MC{S}_n$ is \SBSI\ where $n \geq 2$ and $0\leq p<q \leq n$. Specializing the notion (page \pageref{defSimplicial}) of $A$-indexed partial elements to when $A=\{p,q\}$, recall that a {\bf $(p,q)$-indexed partial element of $y$} is a set $\SETT{b_p, b_q}$ where $b_p \in d_p(y), b_q \in d_q(y)$ and $e_p(b_q)=e_{q-1}(b_p)$. For brevity we will shorten this phrase to ``a $(p,q)$-partial element of $y$''.
\MS

Given a $(p,q)$-partial element $\SETT{b_p, b_q}$ of $y$ there is a unique ``filler'' $b \in \Sig(y)$ which may or may not be an element of $y$. 
\MS

\NI {\bf Notations:}
\begin{itemize}
\item 
$\PAREL{p,q}{y} \DFAS$ the set of $(p,q)$-partial elements of $y$.
\index{Par(p,q)(y)}
\index{$\PAREL{p,q}{y}$}

\item
$\FILL{p,q}{y} \DFAS$ the set of fillers of $(p,q)$-partial elements of $y$.
\index{Fill(p,q)(y)}
\index{$\FILL{p,q}{y}$}

\end{itemize}
In general, $\ELEM{y} \SBS \FILL{p,q}{y} \SBS \Sig(y)$. To say that $y$ is \DET{(p,q)} means that $\ELEM{y} = \FILL{p,q}{y}$.

\subsubsection{Partial elements and determinacy for relations}

\NI We will transfer various notions and notations concerning the $n$-simplex $y^R$ to $R$, as follows. 

\begin{itemize}
\item 
For each $p \in [n]$ and $\LIST{t}{n} \in \PLIST{V}{n}$
\[
\begin{aligned}
e_p : \PLIST{V}{n} \to V_0 \DDD{} \omit{V_p} \DDD{} V_n \\
e_p\LIST{t}{n} \DFAS (\OM{t}{n}{p})
\end{aligned}
\]
and
\[
\begin{aligned}
c_p: \PLIST{V}{n} \to V_0 \DDD{} V_p V_p  \DDD{} V_n\\
c_p \LIST{t}{n} \DFAS (t_0 \DDD{,} t_p, t_p \DDD{,} t_n)
\end{aligned}
\]

\item
A {\bf $(p,q)$-partial element of $R$} is $\SETT{b_p, b_q}$ where $p<q$ in $[n]$ such that $b_p \in d_p(R)$, $b_q \in d_q(R)$ and $e_p(b_q) = e_{q-1}(b_p)$. Denote by $\PAREL{p,q}{R}$ the set of $(p,q)$-partial elements of $R$.

Given $R$ and a \DETY{(p,q)} condition, we have $R \to \PAREL{p,q}{R}$ defined by $t \mapsto \SETT{e_p(t), e_q(t)}$.

Note also that a $(p,q)$-partial element of $R$ can have the form $\SETT{e_p(t), e_q(u)}$ where $t$ and $u$ need not belong $R$.

\item
Given $p<q$ in $[n]$ we define $\FILL{p,q}{R}$ to be the set of all $t \in \PLIST{V}{n}$ such that $e_p(t) \in d_p(R)$ and $ e_q(t) \in d_q(R)$. Clearly, $R \SBS \FILL{p,q}{R}$.

\item
{\bf $R \SBS V_0 \DDD{}V_n$ is \DET{(p,q)}} if $R = \FILL{p,q}{R}$.

\end{itemize}
\MS

\NI {\bf Determinacy in terms of coordinates}
\MS

Suppose $R \SBS \PLIST{V}{n}$ where $n \geq 2$. Suppose $u, t \in \PLIST{V}{n}$
\[
u = \LIST{u}{n}, \qquad t = \LIST{t}{n}
\]
and consider
\[
\begin{aligned}
e_p(u) &= (u_0 \DDD{,} \omit{u_p} \DDD{,} u_n)\\
e_q(t) &= (t_0\DDD{,} \omit{t_q} \DDD{,}, t_n)
\end{aligned}
\]
Then 
\[
\begin{aligned}
e_p e_q(t) &= (t_0 \DDD{,} \omit{t_p} \DDD{,} \omit{t_q} \DDD{,} t_n) \\
e_{q-1}e_p(u) = e_p e_q(u) &= (u_0 \DDD{,} \omit{u_p} \DDD{,} \omit{u_q} \DDD{,} u_n)
\end{aligned}
\]
If $e_p e_q(t) = e_{q-1} e_p(u)$ then for all $j \neq p \text{ or } q$, $u_j=t_j$ and the unique ``filler'' for $e_p(u)$ and $e_q(t)$ is
\[
v = \LIST{v}{n} = (u_0 \DDD{,} u_{p-1}, t_p, u_{p+1} \DDD{,} u_q \DDD{,} u_n) \in \PLIST{V}{n}
\]
since
\[
\begin{aligned}
e_p(v) &= (u_0 \DDD{,} u_{p-1}, \omit{t_p}, u_{p+1} \DDD{,} u_q \DDD{,} u_n) &= e_p(u) \\
e_q(v) &= (u_0 \DDD{,} u_{p-1}, t_p, u_{p+1} \DDD{,} \omit{u_q} \DDD{,} u_n) & \\
&= (t_0 \DDD{,} t_{p-1}, t_p, t_{p+1} \DDD{,} \omit{u_q} \DDD{,} t_n) &= e_q(t)
\end{aligned}
\]
That is, 
\[
\SETT{e_p(u), e_q(t)} = \SETT{e_p(v), e_q(v)}
\]

Therefore, $R$ is \DET{(p,q)} if and only if for all $u, t \in \PLIST{V}{n}$ such that $e_p(u) \in d_p(R)$, $e_q(t) \in d_q(R)$ and $e_p e_q(t) = e_p e_q(u)$ then
\[
v = (u_0 \DDD{,} u_{p-1}, t_p, u_{p+1} \DDD{,} u_{q-1},u_q, u_{q+1} \DDD{,} u_n) \in R
\]

\begin{example}
{\rm
Let $n=6$ and consider \DETY{(1,4)} in dimension 6 which is required for obtaining an example of a \NICOMP{6}{3}. Consider the relation 
\[
R = \SETT{r_1,r_2,r_3} \SBS \PLIST{V}{6}
\]
where
\[
\begin{aligned}
r_1 &= (0,1,2,3,4,5,6)\\
r_2 &= (0,1,2,7,4,5,6) \\
r_3 &= (0,8,2,3,4,9,6)
\end{aligned}
\]
and the $V_j$ are supersets of $\SETT{0},\SETT{1,8},\SETT{2},\SETT{3,7},\SETT{4},\SETT{5,9},\SETT{6}$.
\MS

A $(1,4)$-partial element of $R$ would have the form $\SETT{e_1(u), e_4(w)}$ for some
\[
u = (u_0,u_1,u_2,u_3,u_4,u_5,u_6), \qquad w = (w_0,w_1,w_2,w_3,w_4,w_5,w_6)
\]
with $e_1(u) \in d_1(R)$, $e_4(w) \in d_4(R)$ and $e_1 e_4(w)=e_3 e_1(u) = e_1 e_4(u)$.
We have:
\[
\begin{aligned}
d_4(R) &= \SETT{(0,1,2,3,5,6),(0,1,2,7,5,6),(0,8,2,3,9,6)}\\
d_1 d_4(R) &= \SETT{(0,2,3,5,6),(0,2,7,5,6),(0,2,3,9,6)}
\end{aligned}
\]
In this example the only occurence of $u, w \in R$ such that $e_1 e_4(w)= e_1 e_4(u)$ is when $w=u$. It follows that $R$ is \DET{(1,4)}. 
}
\end{example}

\begin{newlem}
{\rm
\SL 

Given $n \geq 1$ and relations $R, R' \SBS \PLIST{V}{n}$ then:

\begin{enumerate}
\item 
If $R \SBS R'$ then for all $j \in [n]$, $d_j(R) \SBS d_j(R')$ and $s_j(R) \SBS s_j(R')$.

\item
If $R \SBS R'$ then for all \DETY{(p,q)} conditions in dimension $n$, $\FILL{p,q}{R} \SBS \FILL{p,q}{R'}$.

\item
For all $p \in [n]$, $d_p(\FILL{p,q}{R}) = d_p(R)$. 

\item
For all $j \in [n]$, $d_j(R \cap R') \SBS d_j(R) \cap d_j(R')$ and $s_j(R \cap R') \SBS s_j(R) \cap s_j(R')$.

\item
$\FILL{p,q}{R \cap R'} \SBS \FILL{p,q}{R} \cap \FILL{p,q}{R'}$.

\item
$\FILL{p,q}{R}$ is \DET{(p,q)}.

\end{enumerate}
}
\end{newlem}

\Proof

Items (1) and (2) follow immediately from the definitions. 
\MS

Item (3): $R \SBS \FILL{p,q}{R}$ implies $d_p(R) \SBS d_p(\FILL{p,q}{R})$. And if $b_p \in d_p(\FILL{p,q}{R})$ then $b_p = e_p(b)$ for some $b \in \FILL{p,q}{R}$. That implies $e_p(b) = b_p \in d_p(R)$.
\MS

Items (4) and (5) follow from (1) and (2) and the inclusions $R \cap R' \SBS R$ and $R \cap R' \SBS R'$.
\MS

Item (6): Suppose $\SETT{b'_p, b'_q}$ is a $(p,q)$-partial element of $\FILL{p,q}{R}$ filled by $b' \in \FILL{p,q}{\FILL{p,q}{R}}$. Then $b'_p \in d_p(\FILL{p,q}{R}) = d_p(R)$ and $b'_q \in d_q(\FILL{p,q}{R}) = d_q(R)$ which implies that $\SETT{b'_p, b'_q}$ is a $(p,q)$-partial element of $R$. Therefore $b' \in \FILL{p,q}{R}$. That is, $\FILL{p,q}{R}$ is \DET{(p,q)}.

\qed

\subsubsection{Covers, anchors, determinacy condition sets and enlargements}

Given a composition of degeneracy and/or face operators in standard form
\[
f = s_{q_j} \DDD{} s_{q_1}\; d_{p_1} \DDD{} d_{p_k} \qquad \text{where }p_1 \DDD{<} p_k \text { and } q_1 \DDD{<} q_j
\]
we will denote by $\bar{f}$ the corresponding composition of functions
\[
\bar{f} = c_{q_j} \DDD{} c_{q_1}\; e_{p_1} \DDD{} e_{p_k}
\]
If $R \SBS \PLIST{V}{n}$ is any non-empty relation then $R$ and $f(R)$ relate as shown in the next diagram, where $f(R)$ is the vertex-relation of $f(y^R)$.

\setlength{\unitlength}{1in}
\begin{center}
\begin{picture}(3,1.1)
\put(.5,0){\MB{f(R)}}
\put(.5,1){\MB{R}}
\put(2,0){\MB{f(\PLIST{V}{n})}}
\put(2,1){\MB{\PLIST{V}{n}}}
\put(.5,.8){\vector(0,-1){.6}}
\put(.4,.5){\MB{\bar{f}}}
\put(2,.8){\vector(0,-1){.6}}
\put(2.1,.5){\MB{\bar{f}}}
\put(.8,0){\vector(1,0){.7}}
\put(1.1,.1){\makebox(0,0){\footnotesize{incl.}}}
\put(.8,1){\vector(1,0){.7}}
\put(1.1,1.1){\makebox(0,0){\footnotesize{incl.}}}
\end{picture}
\end{center}
\BS

\begin{newdef}
{\bf (Covers and anchors)}
\label{cover-and-anchor}
\index{cover of an element}
\index{anchor}
\index{anchored cover}
\index{anchor-complete}
{\rm
\SL

Suppose $R \SBS \PLIST{V}{n}$ is a non-empty relation and suppose 
\[
f = s_{q_k} \DDD{} s_{q_1}\, d_{p_1} \DDD{} d_{p_j}
\]
is in standard form where $j \geq 1$ and $k \leq j$. The simplicial image of $R$ is $f(R)$ which, since $k \leq j$, is the vertex-relation of $f(y^R)$ and belongs to the $n$-truncated complex $\MC{S}[y^R]$. 

\begin{enumerate}
\item 
Given $b \in f(\PLIST{V}{n})$ and $\bar{b} \in \PLIST{V}{n}$, then $\bar{b}$ is a {\bf cover of $b$} if $\bar{f}(\bar{b}) = b$. (Note: Even if $b \in f(R)$ there is no presumption that $\bar{b} \in R$).

\item
If $\SETT{b_p,b_q}$ is a $(p,q)$-indexed partial element of $f(R)$ then by a {\bf cover of the partial element} we mean a cover of its unique filler.

\item

Suppose $a = \LIST{a}{n} \in R$. Given $b \in f(\PLIST{V}{n})$, and a cover $\bar{b} = ( \bar{t}_0 \DDD{,} \bar{t}_n) \in \PLIST{V}{n}$ of $b$, then we say $\bar{b}$ is {\bf anchored to $a$} if for all $i \in \{p_1 \DDD{,} p_j\}$, $\bar{t}_i = a_i$.

\item $R$ is {\bf anchor-complete with respect to $a$} if for all $f$ (as above) and all $b \in f(R)$, $R$ contains the $a$-anchored cover $\bar{b}$ of $b$. That is, $R$ contains the $a$-anchored cover of every element in every simplicial image of $R$.

\end{enumerate}
}
\BX
\end{newdef}
\MS

\begin{example}
{\rm
Let $n=5$ and $V_0 \DDD{=} V_5 = \mathbb{N}$.
Consider the relation $R = \SETT{r_1,r_2,r_3,r_4} \SBS V_0 \DDD{} V_5$ where
\[
\begin{aligned}
r_1 &= (0,1,2,3,4,5)\\
r_2 &= (6,1,2,7,4,5)\\
r_3 &= (8,1,9,3,4,5)\\
r_4 &= (8,10,9,3,11,5)
\end{aligned}
\]
Suppose $f \DFAS s_3 d_1 d_4$. Then 
\[
\begin{aligned}
\bar{f}(\PLIST{V}{5}) &= V_0 V_2 V_3 V_5 V_5\\
\bar{f}(t_0,t_1,t_2,t_3,t_4,t_5) &= (t_0,t_2,t_3,t_5,t_5)\\
\bar{f}(R) &= \SETT{(0,2,3,5,5), (6,2,7,5,5), (8,9,3,5,5)}
\end{aligned}
\]
noting that $\bar{f}(r_3) = \bar{f}(r_4) = (8,9,3,5,5)$.
\MS

(1) By construction, each element of $\bar{f}(R)$ has at least one cover in $R$. In this example, $(8,9,3,5,5)$ happens to have two covers.
\MS

(2) Consider $(2,3,5,5) \in d_0 f(R)$ and $(6,2,5,5) \in d_2 f(R)$. These form a $(0,2)$-partial element of $\bar{f}(R)$
\[
\SETT{(2,3,5,5),\, (6,2,5,5)}
\]
since 
\[
e_0(6,2,5,5)=(2,5,5) = e_1(2,3,5,5)
\]
Its unique filler, $(6,2,3,5,5) \in \bar{f}(\PLIST{V}{5})$, does not belong to $\bar{f}(R)$ and it follows that $\bar{f}(R)$ is not \DET{(0,2)}.
\MS

(3) Continuing with $(6,2,3,5,5) \in \bar{f}(\PLIST{V}{5})$, a cover of this element is any element $(6,t_1,2,3,t_4,5) \in \PLIST{V}{5}$ and there are, in general, many such covers; in this example, infinitely many. If we choose and fix an element $a = \LIST{a}{5} \in \PLIST{V}{5}$ to serve as anchor then the unique cover of $(6,t_1,2,3,t_4,5)$ anchored to $a$ is 
$(6,a_1,2,3,a_4,5)$.
\MS

(4) Given a choice of anchor $a \in R$ it might not be the case that that the anchored cover of a given element of $f(R)$ belongs to $R$. For example, $(6,2,7,5,5) \in f(R)$ is covered by $r_2 = (6,1,2,7,4,5) \in R$. If we choose as anchor $a=(8,10,9,3,11,5)=r_4 \in R$ then the anchored cover of $(6,2,7,5,5)$ is $(6,\fbox{10},2,7,\fbox{11},5) \notin R$.
}
\end{example}
\BS

\NI {\bf Notes:}

\begin{enumerate}
\item 
Given $a = \LIST{a}{n} \in \PLIST{V}{n}$ as an anchor and any $b \in f(\PLIST{V}{n})$ then there is exactly one $a$-anchored cover $\bar{b}$ of $b$. If the entries in $b$ are $t_i \in V_i$ where $i \notin \SETT{p_1 \DDD{,}, p_j}$, and if $\bar{b} = ( \bar{t}_0 \DDD{,} \bar{t}_n)$ then
\[
\bar{t}_i = \left\{
\begin{array}{ll}
a_i & \text{if } i \in \SETT{p_1 \DDD{,} p_j} \\
& \\
\text{the entry of } b \text{ from } V_i & \text{if }i \notin \SETT{p_1 \DDD{,} p_j}
\end{array}
\right.
\]

\item 
A trivial but useful example of an anchor-complete relation on $\SEQ{V}{n}$ is $\PLIST{V}{n} \SBS \PLIST{V}{n}$. It is anchor-complete with respect to any choice of $a \in \PLIST{V}{n}$.

\item
Again, suppose $a = \LIST{a}{n} \in \PLIST{V}{n}$. Given $R \SBS \PLIST{V}{n}$, a non-empty relation, and $f = s_{q_k} \DDD{} s_{q_1} d_{p_1} \DDD{} d_{p_j}$ as above, suppose $b \in f(R)$ and that $\bar{b} \in \PLIST{V}{n}$ is the unique cover of $b$ anchored to $a$. Then $\bar{b}$ belongs to every enlargement $R' \SPS R$ which is anchor-complete with respect to $a$. Therefore, for all $m$, for all $f:\MC{S}_n \to \MC{S}_m$, and for all $b \in f(R)$, every $a$-anchored cover of $b$ belongs to every enlargement of $R$ which is anchor-complete with respect to $a$. Thus the relation $\bar{R}$ defined by
\[
\bar{R} \DFAS \bigcap R'
\]
where the intersection is over all enlargements of $R$ which are anchor-complete with respect to $a$, is the smallest enlargement of $R$ which is anchor-complete with respect to $a$. (The family of such $R'$ is non-empty, as noted just above).

\end{enumerate}

\BX

\subsubsection{Truncations of \NICOMP{n}{i}s; $0<i<n+1$}

Suppose $n \geq 3$. Consider the set of determinacy conditions
\[
\MC{D}(n) \DFAS \SETT{(p,q).n : 0 \leq p<q \leq n \text{ and } q-p=2 \text{ or } 3}
\]
\index{$\MC{D}(n)$}
The reason for this definition will be apparent below.
Order $\MC{D}(n)$ lexicographically with the first coordinate ``$p$'' as the leading symbol. For example, the ordering of $\MC{D}(5)$ is
\[
(0,2).5 \prec (0,3).5 \prec (1,3).5 \prec (1,4).5 \prec (2,4).5 \prec (2,5).5 \prec (3,5).5
\]
A ``subinterval'' of $\MC{D}(n)$ will refer to any non-empty subset of it which contains all the elements of $\MC{D}(n)$ between the endpoints. A one-element subset of $\MC{D}(n)$ will also count as a subinterval in this sense.

If $\MC{S}'$ is an $n$-truncated subcomplex of \SBSI\ simplices such that all $n$-simplices of $\MC{S}'$ are \DET{(p,p+3)} or \DET{(p,p+2)} then the Degenerate Determinacy Theorem (page \pageref{degenDetlemma}) and its corollary ``Down Rule'' (Theorem \refpage{down rule}) imply certain determinacy conditions in lower dimensions. To recall:
\begin{itemize}

\item 
{\bf D1:} $(p,q).n \IMP (p-1,q-1).(n-1)$ if $p \geq 2$.

\item
{\bf D2:} $(p,q).n \IMP (p,q-1).(n-1)$ if $q-p \geq 3$.

\item
{\bf D3:} $(p,q).n \IMP (p,q).(n-1)$ if $q \leq n-2$.
\end{itemize}

Observe that if $(p,q).n \in \MC{D}(n)$ and if at least one of the three rules applies to it then the result(s) form an subinterval of $\MC{D}(n-1)$. Also, given any two consecutive elements of $\MC{D}(n)$, then the Rules applied to them determine two possibly overlapping subintervals of $\MC{D}(n-1)$. And from this it follows:

\begin{newlem}
{\rm
\SL 

Given $n \geq 4$ then the Down Rules applied any subinterval of $\MC{D}(n)$ yield a subinterval of $\MC{D}(n-1)$.
}
\end{newlem}
\qed

In particular, given an interval 
\[
(p_1,q_1).n \DDD{\prec} (p_k,q_k).n
\] 
then the rules applied to these yield an interval in $\MC{D}(n-1)$ whose least element is either 
\[
\begin{aligned}
(p_1-1,q_1-1).(n-1) &\quad \text{ if } p_1 \geq 2 \\
(p_1,q_1-1).(n-1) &\quad \text{ if } q_1-p_1 \geq 3 \text{ and } p_1<2\\
(p_1,q_1).(n-1)&\quad \text{ if } q_1\leq n-2, p_1<2 \text{ and } q_1-p_1<3
\end{aligned}
\]
Similarly, the greatest of that subinterval in $\MC{D}(n-1)$ is
\[
\begin{aligned}
(p_k,q_k).(n-1) & \quad \text{ if } q_k\leq n-2 \\
(p_k,q_k-1).(n-1) &\quad \text{ if } q_k-p_k \geq 3 \text{ and } q_k>n-2\\
(p_k-1,q_k-1).(n-1) &\quad\text{ if } p_k \geq 2, q_k-p_k,3 \text{ and } q_k>n-2
\end{aligned}
\]
\MS

\begin{example}
{\rm
Consider $(2,5).7 \in \MC{D}(7)$. Rules D1, D2 and D3 apply to yield the following sequence of intervals:
\[
\begin{aligned}
(2,5).7 & \text{ in } \MC{D}(7) \\
(1,4).6 \prec (2,4).6 \prec (2,5).6 & \text{ in } \MC{D}(6) \\
(1,3).5 \prec (1,4).5 \prec (2,4).5 & \text{ in } \MC{D}(5) \\
(1,3).4 & \text{ in } \MC{D}(4)
\end{aligned}
\]
a sequence of intervals which terminates in dimension 4 because none of the Rules apply to $(1,3).4$. (See the diagram below).
}
\end{example}
\BS

\NI {\bf Some special cases of the Rules:}
\begin{enumerate} 

\item %
{\bf Repeated use of D3:}

Given $(p,q).n$ where $q \leq n-2$ then the conditions implied by applying Rule 3 repeatedly yields the sequence:
\[
\DPQ{p}{q}{n} \IMPBY{D3} \DPQ{p}{q}{(n-1)}\DDD{\IMPBY{D3}} \DPQ{p}{q}{(q+1)}
\]
Particular instances of this are:
\begin{enumerate} 
\item 
If $n>4$ then
\[
\DPQ{1}{3}{n} \IMPBY{D3} (1,3).(n-1) \DDD{\IMPBY{D3}} \DPQ{1}{3}{4}
\]
which terminates at $\DPQ{1}{3}{4}$ (to which no rules apply).

\item 
If $n>p+4$ then
\[
\DPQ{p}{p+3}{n} \DDD{\IMPBY{D3}} \DPQ{p}{p+3}{(p+4)}
\]

\item 
If $n>p+3$ then
\[
\DPQ{p}{p+2}{n} \DDD{\IMPBY{D3}} \DPQ{p}{p+2}{(p+3)}
\]

\end{enumerate}  

\item %
{\bf Repeated use of D1:}

Given $\DPQ{p}{q}{n}$ where $p \geq 2$ then the conditions implied by applying D1 repeatedly yields the sequence:
\[
\DPQ{p}{q}{n} \IMPBY{D1} \DPQ{p-1}{q-1}{(n-1)}\DDD{\IMPBY{D1}} \DPQ{1}{q-p+1}{(n-p+1)}
\]
Particular instances of this are
\begin{enumerate} 
\item 
If $p \geq 2$ and $n \geq p+2$ then
\[
\DPQ{p}{p+2}{n} \DDD{\IMPBY{D1}} \DPQ{1}{3}{(n-p+1)}(1,3).(n-p+1)
\]

\item 
If $p \geq 2$ and $n \geq p+3$ then
\[
\DPQ{p}{p+3}{n} \DDD{\IMPBY{D1}} \DPQ{1}{4}{(n-p+1)}
\]
and note that D2 gives $\DPQ{1}{4}{(n-p+1)} \IMPBY{D1} \DPQ{1}{3}{(n-p+1)}$.

\end{enumerate} 

\end{enumerate} 

\NI Here is a diagram showing how the rules apply, starting with $\MC{D}(6)$:
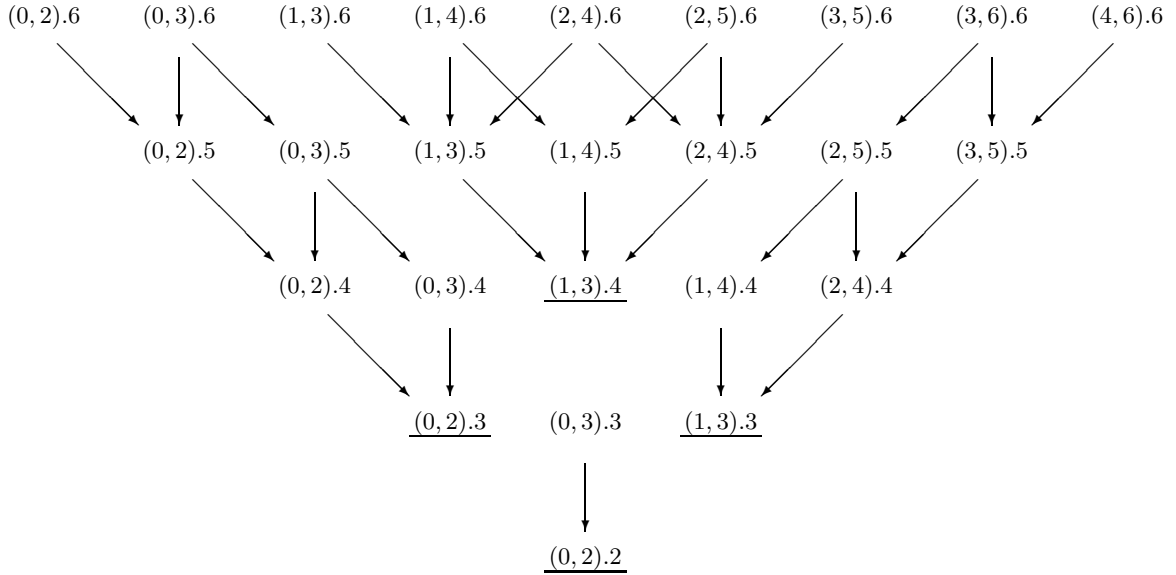
\index{Rule diagram for $\MC{D}(n)$}
\setlength{\unitlength}{1.8cm}
\begin{center}
\begin{picture}(8.5,4.1)
\put(0,4){\MBS{(0,2).6}}
\put(1,4){\MBS{(0,3).6}}
\put(2,4){\MBS{(1,3).6}}
\put(3,4){\MBS{(1,4).6}}
\put(4,4){\MBS{(2,4).6}}
\put(5,4){\MBS{(2,5).6}}
\put(6,4){\MBS{(3,5).6}}
\put(7,4){\MBS{(3,6).6}}
\put(8,4){\MBS{(4,6).6}}
\put(1,3){\MBS{(0,2).5}}
\put(2,3){\MBS{(0,3).5}}
\put(3,3){\MBS{(1,3).5}}
\put(4,3){\MBS{(1,4).5}}
\put(5,3){\MBS{(2,4).5}}
\put(6,3){\MBS{(2,5).5}}
\put(7,3){\MBS{(3,5).5}}
\put(2,2){\MBS{(0,2).4}}
\put(3,2){\MBS{(0,3).4}}
\put(4,2){\MBS{(1,3).4}}
\put(3.7,1.9){\line(1,0){.6}}
\put(5,2){\MBS{(1,4).4}}
\put(6,2){\MBS{(2,4).4}}
\put(3,1){\MBS{(0,2).3}}
\put(2.7,.9){\line(1,0){.6}}
\put(4,1){\MBS{(0,3).3}}
\put(5,1){\MBS{(1,3).3}}
\put(4.7,.9){\line(1,0){.6}}
\put(4,0){\MBS{(0,2).2}}
\put(3.7,-.1){\line(1,0){.6}}
\multiput(0.1,3.8)(1,-1){3}{\vector(1,-1){.6}}
\multiput(1.1,3.8)(1,-1){2}{\vector(1,-1){.6}}
\multiput(2.1,3.8)(1,-1){2}{\vector(1,-1){.6}}
\multiput(3.1,3.8)(1,-1){1}{\vector(1,-1){.6}}
\multiput(4.1,3.8)(1,-1){1}{\vector(1,-1){.6}}
\multiput(1,3.7)(1,-1){4}{\vector(0,-1){.5}}
\multiput(3,3.7)(1,-1){3}{\vector(0,-1){.5}}
\multiput(5,3.7)(1,-1){2}{\vector(0,-1){.5}}
\put(7,3.7){\vector(0,-1){.5}}
\put(3.9,3.8){\vector(-1,-1){.6}}
\put(4.9,3.8){\vector(-1,-1){.6}}
\multiput(5.9,3.8)(-1,-1){2}{\vector(-1,-1){.6}}
\multiput(6.9,3.8)(-1,-1){2}{\vector(-1,-1){.6}}
\multiput(7.9,3.8)(-1,-1){3}{\vector(-1,-1){.6}}
 
\end{picture}
\end{center}
In this diagram various determinacy conditions are arranged as follows: each horizontal row is $\MC{D}(n)$ in lexicographic order, $2 \leq n \leq 6$, starting with $(0,2).n$ and ending with $(n-2,n).n$. The arrows are determinacy implications, where northwest-southeast arrows are D3's, the north-south arrows are D2's and the northeast-southwest arrows are D1's. The underlined conditions are terminal in the sense that no Down Rule applies to them.

\BX
\MS

\begin{newlem}
{\rm
\SL

\begin{enumerate}
\item 
Every degenerate \SBSI\ 2-simplex of $\MC{S}$ is \DET{(0,2)}.

\item
Every degenerate \SBSI\ 3-simplex of $\MC{S}$ is \DET{(0,2)}.

\item
Every degenerate \SBSI\ 3-simplex of $\MC{S}$ is \DET{(1,3)}.

\item
Every degenerate \SBSI\ 4-simplex of $\MC{S}$ is \DET{(1,3)}.

\end{enumerate}

}
\end{newlem}

\Proof

These statements follow from the lemma \refpage{restricted-degenerate-determinacy}.

\qed

\begin{newthm}
\label{nondegen-fillers-implies-all-fillers}
{\rm
\SL

Suppose $n \geq 3$. Consider the following three cases.
\begin{enumerate}
\item 
$i=1$ and $\MC{D}$ is the set of determinacy conditions derived from $(0,2).n$ by the Down Rules.

\item
$1<i<n$ and $\MC{D}$ is the set of determinacy conditions derived from $(i-2,i+1).n$ by the Down Rules.

\item 
$i=n$ and $\MC{D}$ is the set of determinacy conditions derived from \\
$(n-2,n).n$ by the Down Rules.

\end{enumerate}

Suppose, for each case separately, $\MC{S}'$ is an $n$-truncated complex of \SBSI\ simplices such that for each $m \leq n$, each {\em non-degenerate} $y \in \MC{S}'_m$ satisfies the dimension $m$ conditions in $\MC{D}$. Then for each such $m$, {\em every} $m$-simplex of $\MC{S}'$ satisfies the dimension $m$ conditions in $\MC{D}$.
}
\end{newthm}

\Proof

Although the proof splits onto cases depending on $i$, each case applies the following step in common. Given $m<n$, and $x \in \MC{S}'_m$ where $x$ satisfies all the dimension $m$ conditions in $\MC{D}$, then for each $k \in [m]$ and each condition $(p,q).(m+1) \in \MC{D}$, $\{k,k+1\} \cap \{p,q\}$ is either
\begin{itemize}
\item 
empty, in which case the Degenerate Determinacy Theorem implies that $s_k(x)$ satisfies $(p,q).(m+1)$.

\item
or a singleton (since $q-p \geq 2$), in which case lemma \refpage{restricted-degenerate-determinacy} implies $s_k(x)$ satisfies $(p,q).(m+1)$.

\end{itemize}

\NI The cases are: (i) $i=1$ or $2$; (ii) $2<i < n-1$; (iii) $i=n-1$ or $n$ which determine the lowest dimensional determinacy condition in $\MC{D}$.
\MS

\NI {\bf Case: $i=1$ or $2$} 

In this case, the lowest dimensional condition in $\MC{D}$ is $(0,2).3$ (see the diagram above). As noted in the lemma above, all degenerate 3-simplices of $\MC{S}'$ are \DET{(0,2)}. Working up from dimension 3 to dimension $n$, suppose $x \in \MC{S}'_m$ where $3 \leq m<n$. Let $(p,q).(m+1)$ be any of the dimension $m+1$ determinacy conditions in $\MC{D}$. Let $k \in [m]$. The step-in-common given above proves that $s_k(x)$ satisfies $(p,q).(m+1)$.
\MS

\NI {\bf Case: $2<i <n-1$} 

In this case the lowest dimensional condition in $\MC{D}$ is $(1,3).4$. As noted in the lemma above, every degenerate 4-simplex in $\MC{S}'$ satisfies this determinacy condition. Working up from dimension 4 to dimension $n$, the same reasoning as in the previous case shows that for each $m$ from 4 to $n$, and every $(p,q).(m+1) \in \MC{D}$, every degenerate $(m+1)$-simplex of $\MC{S}'$ satisfies $(p,q).(m+1)$.
\MS

\NI {\bf Case: $i=n-1$ or $n$}

In this case, the lowest dimensional condition in $\MC{D}$ is $(1,3).3$. As noted in the lemma above, every degenerate 3-simplex in $\MC{S}'$ satisfies this determinacy condition. The reasoning of the previous two cases also applies to show that for each $3 \leq m < n$, every degenerate $(m+1)$-simplex of $\MC{S}'$ satisfies every dimension $m+1$ determinacy condition in $\MC{D}$.

\qed
\BS

In the previous theorem, the supposed $n$-truncated complex $\MC{S}'$ is assumed to satisfy a specific sets of determinacy conditions for certain of its non-degenerate simplices. The existence of such truncated complexes is implied by the next theorem and uses a standard least-upper-bound argument. 

\begin{newthm}
{\rm
\SL

Let $n \geq 4$ and let $\MC{D}$ be a set of determinacy conditions in dimensions $ \leq n$ such that for all $(p,q).m \in \MC{D}$, $q-p \geq 2$. Let $R \SBS \PLIST{V}{n}$ be a relation and $a = \LIST{a}{n} \in R$. Let 
\[
\SETT{R_\alpha \SBS \PLIST{V}{n} : \alpha \in \MC{A} }
\]
be a set of enlargements (supersets) of $R$ such that for all $\alpha \in \MC{A}$, $R_\alpha$ has the property:
\begin{quote}
(*) For all $(p,q).m \in \MC{D}$ and for all $(p,q)$-indexed partial elements of all $m$-dimensional non-degenerate subfaces $f(R_\alpha)$ of $R_\alpha$, there exists an $a$-anchored cover in $R_\alpha$.
\end{quote}
Then $\bar{R} \DFAS \bigcap_\alpha R_\alpha$ also satisfies (*).
}
\end{newthm}

\Proof

The family of $R_\alpha$'s is non-empty because it includes the maximal relation $\PLIST{V}{n} \SBS \PLIST{V}{n}$.

Given any $(p,q).m \in \MC{D}$ and any $m$-dimensional non-degenerate subface $f(\bar{R})$ of $\bar{R}$ where
\[
f(\bar{R}) = d_{r_1} \DDD{} d_{r_{n-m}}(\bar{R}) \qquad r_1 \DDD{<} r_{n-m}
\]
suppose $b$ is the unique filler of a $(p,q)$-partial element $\SETT{b_p,b_q}$ of $f(\bar{R})$. Since, $\bar{R} \SBS R_\alpha$ then $\SETT{b_p,b_q}$ is also a $(p,q)$-partial element of $f(R_\alpha)$. Property (*) of $R_\alpha$ says that there is an $a$-anchored cover $b' = \LIST{t}{n} \in R_\alpha$ which means that $\bar{f}(b')=b$, $t_j = a_j$ for each $j \in \SETT{r_1 \DDD{,} r_{n-m+1} }$, and $t_j = $ the corresponding fundamental entry of $b$ otherwise. Since this filler is uniquely determined by $b$, $\SETT{r_1 \DDD{,} r_{n-m} }$ and $a$, then $b' \in R_\alpha$ for all $\alpha \in \MC{A}$. That is, $b' \in \bar{R}$. This proves that $\bar{R}$ satisfies property (*).

\qed

\begin{newthm}
{\rm
\SL

Suppose $0<i<n+1$. Then there exists an $n$-truncated complex $\MC{S}'$ such that every $n$-simplex of $\MC{S}'$ is: 
\begin{enumerate}
\item 
$(0,2).n$ if $i=1$.

\item
$(i-2,i+1).n$ if $1<i<n$.

\item
$(n-2,n).n$ if $i=n$.
\end{enumerate}
}
\end{newthm}

\Proof

Let $\MC{D}$ be the set of determinacy conditions implied by any of the three cases. Let $R \SBS \PLIST{V}{n}$ be any non-empty relation and choose any $a \in R$. The previous theorem implies there exists an enlargement $\bar{R}$ of $R$ such that: for every $(p,q).m \in \MC{D}$, every composition of face operators $f = d_{p_1} \DDD{} d_{p_{n-m}}$ and for every $(p,q)$-partial element of $f(\bar{R})$ there exists an $a$-anchored cover $\bar{b} \in \bar{R}$ of that partial element. That is, $\bar{f}(\bar{b})$ is the unique filler of the partial element. Thus, setting $\MC{S}'$ to be $\MC{S}[y^{\bar{R}}]$, every {\em non-degenerate} simplex of $\MC{S}'$ satisfies every determinacy condition in $\MC{D}$ which applies to that simplex. By Theorem \refpage{nondegen-fillers-implies-all-fillers}, every simplex of $\MC{S}'$ satisfies every determinacy condition in $\MC{D}$ which applies to that simplex.

\qed
\MS

\begin{example}
{\rm
Let $n=6, i=3$. In order to form the 6-truncation for a \NICOMP{6}{3} the required determinacy condition is \DETY{(1,4)} in dimension 6. The Rules imply the set of determinacy conditions:
\[
\MC{D} = \SETT{(1,4).6,(1,3).5,(1,4).5,(1,3).4}
\]

We start with the following relation $R$:
\[
R = \SETT{
(0,1,2,3,4,5,6), (0,7,2,3,8,5,6), (0,9,2,10,4,5,6)
}
\]
and seek to enlarge $R$ appropriately to $\bar{R}$ so that every 6-simplex of the 6-truncated complex $\MC{S}[y^{\bar{R}}]$ is \DET{(1,4)}. $\MC{S}[y^{\bar{R}}]$ contains $3431$ simplices in all, of which 127 are non-degenerate.
\MS

\NI (Note: This example was calculated using custom software, not by hand.)
\MS

Initially:
\begin{itemize}
\item 
$y^R$ is not \DET{(1,4)} because the $(1,4)$-partial element
\[
\SETT{(0,2,3,4,5,6), (0,7,2,3,5,6)}
\]
has filler $(0,7,2,3,4,5,6) \notin R$. Similarly, the filler of the $(1,4)$-partial element $\SETT{(0,2,3,8,5,6),(0,1,2,3,5,6)}$ does not belong to $R$.

\item
$d_j d_k(y^R)$ is not \DET{(1,3)} for
\[
(j,k) = (2,5), (2,6), (3,5), (3,6), (4,5), (4,6) \text{ and } (5,6)
\]
and $d_k(y)$ is not \DET{(1,3)} for $k=2,3,4,5,6$.

For instance, 
\[
d_2 d_5(R) = \SETT{(0,1,3,4,6),(0,7,3,8,6),(0,9,10,4,6)}
\]
has a $(1,3)$-partial element
\[
\SETT{(0,3,4,6),(0,7,3,6)}
\]
with filler $(0,7,3,4,6) \notin R$.

\item 
$d_k(y)$ is not \DET{(1,3)} for $k=2,3,4,5,6$.

\end{itemize}

Now if we enlarge $R$ to $\bar{R}$ so that those seven 4-dimensional subfaces are \DET{(1,3)} using $(0,1,2,3,4,5,6)$ as the anchor element then $\bar{R}$ has the following nine elements:
\[
\begin{aligned}
(0,1,2,3,4,5,6),\quad & (0,7,2,3,8,5,6), & (0,9,2,10,4,5,6)\\
(0,7,2,3,4,5,6),\quad & (0,1,2,3,8,5,6), & (0,9,2,3,8,5,6) \\
(0,1,2,10,4,5,6),\quad & (0,7,2,10,4,5,6), &(0,9,2,3,4,5,6)
\end{aligned}
\]
This suffices: all 462 4-simplices of $\MC{S}[y^{\bar{R}}]$ are \DET{(1,3)}, all 924 5-simplices are \DET{(1,3)} and all 1716\; 6-simplices are \DET{(1,4)}.
Thus $\MC{S}[y^{\bar{R}}]$ is the 6-truncation of a \NICOMP{6}{3}.
\MS

\NI Variation \#1:

If we {\em change} the choice of anchor then the enlargement (in this example) is different. For instance, if we choose the anchor element to be $(0,7,2,3,8,5,6) \in R$ then the calculations go as follows: the anchored fillers for the missing $(1,4)$-partial elements of the same seven 4-dimensional subfaces enlarge $R$ to a relation $\tilde{R}$ with twelve elements:
\[
\begin{aligned}
(0,1,2,3,4,5,6),\quad & (0,7,2,3,8,5,6), & (0,9,2,10,4,5,6) \\          (0,7,2,3,4,5,6), \quad & (0,1,2,3,8,5,6), & (0,9,2,3,8,5,6) \\         (0,1,2,10,8,5,6), \quad & (0,7,2,10,8,5,6), & (0,9,2,3,4,5,6) \\          (0,1,2,10,4,5,6), \quad & (0,7,2,10,4,5,6), & (0,9,2,10,8,5,6)
\end{aligned}
\]
As with the first choice of anchor, $\MC{S}[y^{\tilde{R}}]$ also is a 6-truncation of a \NICOMP{6}{3}.
\MS

\NI Variation \#2:

Starting with the original three-element relation $R$, we may enlarge $R$ to $R'$ by adding all fillers of of $(1,4)$-partial elements, namely $(0,7,2,3,4,5,6)$ and $(0,1,2,3,8,5,6)$. However, the resulting 6-truncated complex $\MC{S}[y^{R'}]$ falls short of satisfying the required determinacy conditions. Specifically, writing $y=y^{R'}$:
\begin{itemize}
\item 
$s_2d_3(y)$ and eight other degenerate 6-dimensional images of $y$ are not \DET{(1,4)}.

\item
$d_3d_5(y)$ and four other 4-dimensional subfaces aren't \DET{(1,3)}.

\item
Nine 5-simplices including $d_j(y)$, $(j=3,4,5,6)$ are not \DET{(1,3)}, and $s_2d_3d_5(y)$ and four other degenerate 5-simplices are not \DET{(1,4)}.
\end{itemize}
\makebox{ }
}
\end{example}

\subsubsection{Truncations in the cases $i=0$ or $n+1$}

Fix $n \geq 3$. Given $y \in \MC{S}_n$, in order for the the $n$-truncated complex $\MC{S}[y]$ generated by $y$ to yield an \NICOMP{n}{i}, $\MC{S}[y]$ must satisfy a certain set $\MC{D}=\MC{D}(n,i)$ of determinacy conditions according to Down Rules D1, D2 and D3. The requirements are:
\begin{itemize}
\item 
{\bf Case $i=0$:}
\[
\MC{D} = \SETT{(1,2).n, (1,2).(n-1) \DDD{,} (1,2).3}
\]
where only D3 applies.

\item
{\bf Case $i=n+1$:}
\[
\begin{aligned}
\MC{D} = 
\SETT{
(n-2,n-1).n,(n-3,n-2).(n-1) \DDD{,} (1,2).3}\; \cup & \\
 \cup \;\SETT{
(n-1,n).n,(n-2,n-1).(n-1) \DDD{,} (1,2).2
}
\end{aligned}
\]
where only D1 applies.

\end{itemize}
\BS

In this subsection, we will establish necessary and sufficient conditions on a non-empty relation $R \SBS \PLIST{V}{n}$ such that $\MC{S}[y^R]$ is the $n$-truncation of an \NICOMP{n}{i} in the cases $i=0$ and $i=n+1$.
\MS

In the cases $i=0$ and $i=n+1$, the required determinacy conditions have the form $(p,p+1).m$ for various $p$ and  various $m \leq n$. Lemma \refpage{restricted-degenerate-determinacy} does {\em not} apply to these determinacy conditions for $s_p(x)$. Furthermore, certain choices of $y \in \MC{S}_n$ {\em cannot} be enlarged to yield a model of an \NICOMP{n}{i} when $i=0$ or $n+1$, as the next lemma shows.
\MS

\begin{newlem} 
\label{bad-pair-lemma}
{\rm
\SL 

Let $n \geq 3$, $y \in \MC{S}_n$ with $R=R(y)$ the vertex-relation of $y$. Suppose there is a pair of elements of $t=\LIST{t}{n},\, t'=\LIST{t'}{n} \in R$ and $j<k<m$ in $[n]$ such that $t_j=t'_j$ and $t_m=t'_m$ but $t_k \neq t'_k$. Then there are 3-simplices of $\MC{S}[y]$ which are not \DET{(1,2)}.

It follows from this that $\MC{S}[y]$ cannot be the $n$-truncation of either an \NICOMP{n}{0} or an \NICOMP{n}{n+1}.
}
\end{newlem}

\Proof

Given $t,t' \in R$ as stated let $f = d_0 \DDD{} \omit{d_j} \DDD{}\omit{d_k} \DDD{} \omit{d_m} \DDD{} d_n : \MC{S}_n \to \MC{S}_2$. Then
\[
\begin{aligned}
\bar{f}(t) = (t_j, t_k, t_m), & \qquad c_1 \bar{f}(t) = (t_j, t_k, t_k, t_m)\\
\bar{f}(t') = (t_j, t'_k, t_m), & \qquad c_1 \bar{f}(t) = (t_j, t'_k, t'_k, t_m)
\end{aligned}
\]
and $\SETT{\bar{f}(t),\bar{f}(t')}$ is a $(1,2)$-partial element of $s_1 f(R)$ with unique filler $(t_j,t'_k,t_k,t_m)$. Since $t_k \neq t'_k$ this filler cannot belong to $s_1 f(R)$; that is, $s_1 f(R)$ is a 3-simplex which is not \DET{(1,2)}.

The necessary determinacy conditions for truncations of \NICOMP{n}{0}s and \NICOMP{n}{n+1}s of the form $\MC{S}[y]$ include that all 3-simplices of $\MC{S}[y]$ be \DET{(1,2)}. The 3-simplex just described cannot be \DET{(1,2)} and therefore $S[y]$ cannot be the $n$-truncation of either an \NICOMP{n}{0} or an \NICOMP{n}{n+1}. 

\qed
\MS

\begin{example}
{\rm
Suppose $n=6$, $R \SBS \PLIST{V}{6}$ contains $(0,1,2,3,4,5,6)=t$ and $(a,1,2,b,c,5,d)=t'$. Set $y = y^R$. Let $f = d_0 d_1 d_4 d_6$. Then 
\[
\begin{aligned}
\bar{f}(t) &= (2,3,5) \in f(y) & \qquad c_1 \bar{f}(t) = (2,3,3,5) \in s_1 f(y)\\
\bar{f}(t') &= (2,b,5) \in f(y)& \qquad c_1 \bar{f}(t') = (2,b,b,5) \in s_1 f(y)
\end{aligned}
\]
and $\SETT{(2,3,5),(2,b,5)}$ is  $(1,2)$-partial element of $s_1 f(y)$ with filler $(2,b,3,5) \notin s_1 f(y)$.
}
\end{example}

\MS

The previous lemma shows that in order for $\MC{S}[y]$ to be the $n$-truncation for an \NICOMP{n}{0} or \NICOMP{n}{n+1}, the vertex-relation $R=R(y)$ must have the property:
\begin{quote}
For all $j<k<m$ in $[n]$ and all pairs $t,t' \in R$ it must be the case that whenever $t_j=t'_j$ and $t_m=t'_m$ then $t_k=t'_k$. 
\end{quote}

\NI This suggests the following definition.

\begin{newdef}
\label{overlap-definition}
\index{overlap on an interval}
{\rm
\SL

Suppose $n \geq 3$ and $t,t' \in \PLIST{V}{n}$. Given $j \leq k$ in $[n]$ we say that $t$ and $t'$ {\bf overlap on $[j,k]$} if for all $m \in [j,k]$, $t_m=t'_m$ {\em and} for all $m \in [0,j-1] \cup [k+1,n]$, $t_m \neq t'_m$.

Without explicit reference to $j$ and $k$, we say ``$t$ and $t'$ overlap on a single interval''. Also, for brevity, if $t$ and $t'$ do not overlap at all we will speak of them as ``overlapping'' in a single {\em empty} interval.
}

\BX
\end{newdef}

That is, $t,t'$ overlap on the single interval $[j,k]$ if
\[
\begin{aligned}
t & = (t_0 \DDD{,} t_{j-1}, u_j \DDD{,} u_k,t_{k+1} \DDD{,} t_n)\\
t' & = (t'_0 \DDD{,} t'_{j-1}, u_j \DDD{,} u_k,t'_{k+1} \DDD{,} t'_n)
\end{aligned}
\]
where $t_m \neq t'_m$, $m=0 \DDD{,} j-1,k+1 \DDD{,}n$. (Again, we understand the sequence $t_0 \DDD{,} t_{j-1}$ to be empty if $j=0$, and the sequence $t_{k+1} \DDD{,} t_n$ to be empty if $k=n$).
\MS

\begin{newlem} 
\label{image-overlap-lemma}
{\rm
\SL 

Suppose $n \geq 3$. If $t,t' \in \PLIST{V}{n}$ overlap in a single interval then for all $p \in [n]$, $e_p(t)$ and $e_p(t')$ overlap in a single interval and $c_p(t)$ and $c_p(t')$ overlap in a single interval.

Also, if $t$ and $t'$ overlap on the interval with right endpoint $n$ then for all $p \in [n]$, $e_p(t)$ and $e_p(t')$ overlap on an interval with right endpoint $n-1$ (or possibly the empty interval) and $c_p(t)$ and $c_p(t')$ overlap on an interval with right endpoint $n+1$.
}
\end{newlem}

\Proof

If $t,t'$ overlap in the interval $[j,k]$ of length $r=k-j+1$, then $e_p(t)$ and $e_p(t')$ overlap in a single interval of length $r$ if $p<j$ or $p>k$ and of length $r-1$ otherwise. $c_p(t)$ and $c_p(t')$ overlap in a single interval of length $r+1$ if $p \in [j,k]$ and of length $r$ otherwise.

The second claim follows directly.

\qed
\MS

It follows from this lemma that if $t,t' \in \PLIST{V}{n}$ overlap in a single interval then for all appropriate sequences $q_1 \DDD{<} q_m$ and $p_1 \DDD{<} p_m$, then
\[
\begin{aligned}
w & = c_{q_m} \DDD{} c_{q_1} e_{p_1} \DDD{} e_{p_m}(t) \\
w' & = c_{q_m} \DDD{} c_{q_1} e_{p_1} \DDD{} e_{p_m}(t')
\end{aligned}
\]
also overlap in a single interval.
\BS

\NI {\bf Models when $i=0$}
\MS

Suppose $R \SBS  \PLIST{V}{n}$ is non-empty. A $(1,2)$-partial element of $R$ has the form $\SETT{e_1(t), e_2(t')}$ for some pair of elements $t,t'$ of $R$ such that $e_1 e_2(t') = e_1 e_1(t)$. In general, if $R$ is not \DET{(1,2)} then there exists a pair $t,t' \in R$ such that no $t'' \in R$ exists such that $e_1(t'') = e_1(t)$ and $e_2(t'') = e_2(t')$. That is, $t,t'$ introduce a filler element for the partial element $\SETT{e_1(t),e_2(t')}$ which does not belong to $R$.

\begin{newthm} 
\label{overlap-theorem}
{\rm
\SL 

Suppose $n \geq 3$, $R \SBS \PLIST{V}{n}$ is non-empty and that $R$ satisfies property 
\begin{quote}
$(*)$ For all $t,t' \in R$ $t,t'$ overlap on a single (possibly empty) interval.
\end{quote}
Let $y \DFAS y^R$. Then every $n$-simplex in $\MC{S}[y]$ is \DET{(1,2)} and is the $n$-truncation of an \NICOMP{n}{0}.
}
\end{newthm}

\Proof

Suppose $z \in \MC{S}[y]_n$. That is, $z=y$ or $z$ is a (degenerate) simplicial image $f(y)$ where $f=s_{q_r} \DDD{} s_{q_1} d_{p_1} \DDD{} d_{p_r}$ for some $r \geq 1$.

Property $(*)$ and the previous lemma shows that every pair $t,t' \in z$ overlap in a single (possibly empty) interval.

We will show that for all $t,t' \in z$, $\SETT{e_1(t),e_2(t')}$ is a $(1,2)$-partial element of $z$ iff $t=t'$. That is, the only $(1,2)$-partial elements of $z$ are of the form $\SETT{e_1(t),e_2(t)}$ for some $t \in z$, which is equivalent to saying that $z$ is \DET{(1,2)}.

By hypothesis, given any $t,t' \in z$, either $t,t'$ overlap in a single interval, or else do not overlap at all (i.e. for all $m \in [n]$, $t_m \neq t'_m$).

If $t,t'$ don't overlap at all, then $e_1 e_2(t') \neq e_1 e_1(t)$ which means that $\SETT{e_1(t),e_2(t')}$ is not a $(1,2)$-partial element. Otherwise, $t,t'$ overlap in an interval (and are unequal) in two possible cases. 
\MS

\NI Case $j>0$:

Then
\[
\begin{aligned}
t & = (t_0 \DDD{,} t_{j-1}, u_j \DDD{,} u_k,t_{k+1} \DDD{,} t_n)\\
t' & = (t'_0 \DDD{,} t'_{j-1}, u_j \DDD{,} u_k,t'_{k+1} \DDD{,} t'_n)
\end{aligned}
\]
and
\[
\begin{aligned}
e_1 e_2(t') &= (t'_0 , \cdots )\\
e_1 e_1(t) &= (t_0 , \cdots )
\end{aligned}
\]
so that $e_1 e_2(t') \neq e_1 e_1(t)$ since $t_0 \neq t'_0$.
\MS

\NI Case $j=0$:

Then for some $k \geq j$
\[
\begin{aligned}
t & = (u_0 \DDD{,} u_k,t_{k+1} \DDD{,} t_n)\\
t' & = (u_0 \DDD{,} u_k,t'_{k+1} \DDD{,} t'_n)
\end{aligned}
\]
and
\[
\begin{aligned}
e_1 e_2(t') &= ( \cdots , t'_n )\\
e_1 e_1(t) &= (\cdots ,t_n )
\end{aligned}
\]
so that, again, $e_1 e_2(t') \neq e_1 e_1(t)$ because $t_n \neq t'_n$.
\MS

So, no matter how $t,t'$ overlap in a single interval, $\SETT{e_1(t),e_2(t')}$ cannot be a $(1,2)$-partial element unless $t=t'$. It follows that $z$ is \DET{(1,2)}. 

By Theorem \refpage{composer-model;i=0} it follows that $\MC{S}[y]$ is the $n$-truncation of an \NICOMP{n}{0}.

\qed
\MS

Combining Theorem \refpage{composer-model;i=0} and lemma \refpage{bad-pair-lemma} and Theorem \refpage{overlap-theorem} we get:

\begin{newthm}
{\rm
\SL

Suppose $n \geq 3$, $R \SBS \PLIST{V}{n}$ is non-empty and $y = y^R \in \MC{S}_n$. Then: $\MC{S}[y]$ is the $n$-truncation of an \NICOMP{n}{0} $\iff$ $R$ has the property that for all distinct $t,t' \in R$, $t$ and $t'$ overlap in a single (possibly empty) interval.
}
\end{newthm}

\Proof

\NI $\Longrightarrow$: \;If $\MC{S}[y^R]$ is the $n$-truncation of an \NICOMP{n}{0} then the Degenerate Determinacy Theorem implies that all degenerate 3-simplices of $\MC{S}[y]$ are \DET{(1,2)}. By lemma \refpage{bad-pair-lemma}, it follows that  all distinct $t,t' \in R$ either overlap not at all or do so in a single interval.
\MS

\NI $\Longleftarrow$:\; If $R$ has the stated overlap property then Theorem \refpage{overlap-theorem} states that all $n$-simplices of $\MC{S}[y]$ are \DET{(1,2)}. Therefore, by Theorem \refpage{composer-model;i=0}, $\MC{S}[y]$ is the $n$-truncation of an \NICOMP{n}{0}.

\qed
\BS

\NI {\bf Models when $i=n+1$}
\MS

Continuing with $n \geq 3$, non-empty $R \SBS \PLIST{V}{n}$ and $y \DFAS y^R$, $n$-truncated complexes $\MC{S}[y]$ which are models of \NICOMP{n}{n+1}s must satisfy (as noted above) the following set of determinacy conditions:
\[
\begin{aligned}
\MC{D} = 
\SETT{
(1,2).3, (2,3).4 \DDD{,} (n-2,n-1).n}\; \cup & \\
 \cup \;\SETT{
 (1,2).2,(2,3).3 \DDD{,} (n-1,n).n
}
\end{aligned}
\] 

\begin{newlem}
\label{bad-pair-lemma-2}
{\rm
\SL 

With $n$, $R$ and $y$ as above, suppose $t,t' \in R$ are distinct and overlap in a single interval $[j,k]$ where $k < n$. Then there exists a 2-simplex in $\MC{S}[y]$ which is not \DET{(1,2)}.
}
\end{newlem}

\Proof

By hypothesis,
\[
\begin{aligned}
t &= (t_0 \DDD{,} t_{j-1},u_j \DDD{,} u_k,t_{k+1} \DDD{,} t_n) \\
t' &= (t'_0 \DDD{,} t'_{j-1},u_j \DDD{,} u_k,t'_{k+1} \DDD{,} t'_n) 
\end{aligned}
\]
If $k<n-1$ then let $f = s_1 d_0 \DDD{}d_{k-1} d_{k+2} \DDD{} d_n$, and if $k=n-1$ then let $f = s_1 d_0 \DDD{}d_{k-1}$. In either case
\[
\begin{aligned}
\bar{f}(t) &= (u_k,t_{k+1},t_{k+1}) \\
\bar{f}(t') &= (u_k,t'_{k+1},t'_{k+1})
\end{aligned}
\]
and $\SETT{(u_k,t_{k+1}),(u_k,t'_{k+1})}$ is a $(1,2)$-partial element of $f(y)$ whose filler is $(u_k,t'_{k+1},t_{k+1}) \notin f(y)$. That is, $f(y)$ is a 2-simplex of $\MC{S}[y]$ which is not \DET{(1,2)}.

\qed
\MS

\begin{newthm}
\label{overlap-theorem-2}
{\rm
\SL 

Suppose $n \geq 3$, $R \SBS \PLIST{V}{n}$ is non-empty and $y \DFAS y^R$. Assume $R$ has the property:
\begin{quote}
$(*)$\; For all $t,t' \in R$ such that $t \neq t'$ then either $t$ and $t'$ do not overlap or else they overlap in a single interval with right endpoint $n$.
\end{quote}
Then all $n$-simplices in $\MC{S}[y]$ are \DET{(n-2,n-1)} and \DET{(n-1,n)} and $\MC{S}[y]$ is the $n$-truncation of an \NICOMP{n}{n+1}.
}
\end{newthm}

\Proof

Let $z \in \MC{S}[y]_n$. Then $z=y$ or $z$ is a degenerate simplicial image of $y$. By lemma \refpage{image-overlap-lemma}, for all $t$ and $t'$ in the vertex-relation for $z$, either they do not overlap at all or else overlap in a single interval with right endpoint $n$. In the first case $\SETT{e_{n-2}(t), e_{n-1}(t')}$ cannot form an $(n-2,n-1)$-partial element of $z$. In the second case,\[
\begin{aligned}
t &= (t_0 \DDD{,} t_{j-1},u_j \DDD{,} u_n) \\
t' &= (t'_0 \DDD{,} t'_{j-1},u_j \DDD{,} u_n)
\end{aligned}
\]
Then $\SETT{e_{n-2}(t), e_{n-1}(t')}$ cannot be an $(n-2,n-1)$-partial element of $z$ because $t_0 \neq t'_0$ implies
\[
\begin{aligned}
e_{n-2}e_{n-1}(t') &= (t'_0, \cdots ) \\
e_{n-2}e_{n-2}(t) &= (t_0, \cdots )
\end{aligned}
\]
are unequal. That is, the only $(n-2,n-1)$-partial elements of $z$ are $\SETT{e_{n-2}(t), e_{n-1}(t)}$ which, since the filler is $t$, implies that $z$ is \DET{(n-2,n-1)}.

By the same reasoning, $\SETT{e_{n-1}(t),e_n(t')}$ cannot be an $(n-1,n)$-partial element of $z$ which implies that $z$ is \DET{(n-1,n)}. Theorem \refpage{n-n+1-composer-theorem} implies that $\MC{S}[y]$ is the $n$-truncation of an \NICOMP{n}{n+1}.

\qed
\MS

Theorem \refpage{n-n+1-composer-theorem}, lemma \refpage{bad-pair-lemma-2} and Theorem \refpage{overlap-theorem-2} yield:

\begin{newthm}
{\rm
\SL

Suppose $n \geq 3$, $R \SBS \PLIST{V}{n}$ is non-empty and $y \DFAS y^R$. Then $\MC{S}[y]$ is the $n$-truncation of an \NICOMP{n}{n+1} $\iff$ $R$ has the property 
\begin{quote}
$(*)$\; For all $t,t' \in R$, $t \neq t'$ either $t$ and $t'$ do not overlap at all or else overlap in a single interval with right endpoint $n$.
\end{quote}
}
\end{newthm}

\Proof

\NI $\Longrightarrow$: If $\MC{S}[y]$ is the $n$-truncation of an \NICOMP{n}{n+1} then the Degenerate Determinacy Theorem implies that all 2-simplices of $\MC{S}[y]$ are \DET{(1,2)}. By lemmas \refpage{bad-pair-lemma} and \refpage{bad-pair-lemma-2}, if there existed $t,t' \in R$ which overlapped in more than one interval or overlapped in one interval $[j,k]$ where $k<n$ then $\MC{S}[y]$ would have a 2-simplex which is not \DET{(1,2)}. Therefore $R$ must satisfy the property $(*)$.
\MS

\NI $\Longleftarrow$: Theorem \refpage{overlap-theorem-2} shows that when $R$ satisfies property $(*)$ then $\MC{S}[y]$ has the determinacy properties in dimension $n$ which, by Theorem \refpage{n-n+1-composer-theorem}, implies $\MC{S}[y]$ is the $n$-truncation of an \NICOMP{n}{n+1}.

\qed

\subsection{Hypergroupoid models}

As mentioned briefly in section \refpage{what-is-in-paper} in equation \eqref{hypergroupoid-condition} (page \pageref{hypergroupoid-condition}), when given $n \geq 1$, an {\bf $n$-dimensional hypergroupoid} is a simplicial set $C$ such that
\[
\ALL m>n, \; \ALL i \in [m] \Bigl(
\phi_{m,i}:C_m \to \BOX{i}{m}{C} \text{ is an isomorphism}
\Bigr)
\]
This notion is developed in detail below in Definition \refpage{hypergroupoid}. Here we'll consider how the set/relation models developed in section \refpage{nisolution} can be specialized to set/relation models for $n$-dimensional hypergroupoids.
\BS

Let $\mathcal{V}$ be any non-empty family of non-empty sets and, relative to $\MC{V}$, let $\MC{S}$ be the simplicial set of relations as defined in section \refpage{defOfS}.

Given $n>1$ and any $i \in [n+1]$, let $\MC{S}^{n,i}$ denote the \NICOMP{n}{i} subcomplex of $\MC{S}$ as developed in Theorems \refpage{composer-model;1<i<n}, \refpage{n-1-composer-theorem}, \refpage{n-n-composer-theorem} and \refpage{n-n+1-composer-theorem}. Then the simplicial set $\MC{S}^{n}$
\[
\MC{S}^{n} \DFAS \bigcap_{i \in [n+1]} \MC{S}^{n,i} \SBS \mathcal{S}
\]
is an \NICOMP{n}{i} for each $i \in [n+1]$. Since for all $m>n$ and for all $i \in [m]$
\[
\phi_{m,i} : \MC{S}^n_m \to \BOX{i}{m}{\MC{S}^n}
\]
is an isomorphism, $\MC{S}^n$ is an $n$-dimensional hypergroupoid. 
\MS

Given any $y \in \MC{S}'_n$, then $y$ must satisfy those determinacy conditions specified in Theorems \refpage{composer-model;1<i<n}, \refpage{n-1-composer-theorem}, \refpage{n-n-composer-theorem} and \refpage{n-n+1-composer-theorem}, and simplices of $\MC{S}'$ of lower dimensions must satisfy the determinacy conditions spelled out by the rules above. 
\MS

These determinacy conditions, sorted by $i \in [n+1]$, are as follows.
\MS
\begin{itemize}
\item 
{\bf $i=0$:} Conditions
\[
(1,2).3 \DDD{,} (1,2).n
\]

\item
{\bf $i=1$:} Conditions
\[
(0,2).3 \DDD{,} (0,2).n
\]

\item
{\bf $i=n$:} Conditions
\[
(1,3).3,(2,4).4 \DDD{,} (n-2,n).n
\]

\item
{\bf $i=n+1$:} Conditions
\[
\begin{aligned}
(1,2).3,\; & (2,3).3, & (3,4).4 &\;\DDD{,}& (n-1,n).n\\
&\; (1,2).3, & (2,3).4 &\; \DDD{,}& (n-2,n-1).n
\end{aligned}
\]

\item
{\bf $1<i<n$:} Conditions $(i-2,i+1).n$ and all the lower dimensional conditions implied by the Rules.
\end{itemize}

Here is a construction of a finite relation $R \SBS \PLIST{V}{n}$ such that $\MC{S}[y^R]$ is the $n$-truncation of a $n$-dimensional hypergroupoid. We will assume that $n \geq 3$.

First, note that in order for $\MC{S}[y]$ to be an \NICOMP{n}{n+1}, $y=y^R$ must be \DET{(n-2,n-1)} and \DET{(n-1,n)}. By Theorem \refpage{overlap-theorem-2}, $R$ must have the following property:
\begin{quote}
$(*)$ For all pairs $t,t' \in R$, $t$ and $t'$ overlap not at all or overlap in a single interval with right endpoint $n$.
\end{quote}
Note that by Theorem \refpage{overlap-theorem}, if $R$ has property $(*)$ then $\MC{S}[y]$ is also be an \NICOMP{n}{0}. Assuming $t$ and $t'$ overlap on the interval $[j,n]$ where $j \geq 1$, then
\[
\begin{aligned}
t &= (t_0 \DDD{,} t_{j-1},u_j \DDD{,} u_n)\\
t' &= (t'_0 \DDD{,} t'_{j-1},u_j \DDD{,} u_n)
\end{aligned}
\]
and $t_m \neq t'_m$, $0 \leq m < j$. We will refer to the coordinates of $t$ in $[0,j-1]$ as the {\bf ``initial segment of $t$''}. It remains to verify that all the other dimension $n$ determinacy conditions listed above also hold for $y$ and for all dimension $n$ simplicial images of $y$.

Suppose $z$ is $y$ or is any dimension $n$ simplicial image of $y$. We will use that
\begin{itemize}
\item 
$z$ fails to be \DET{(p,q)} iff there exists $t \neq t'$ in $z$ such that $\SETT{e_p(t), e_q(t')}$ is a $(p,q)$-partial element.

\item 
$\SETT{e_p(t), e_q(t')}$ is a $(p,q)$-partial element iff $t$ and $t'$ are equal at all coordinates other than $p$ and $q$.
\end{itemize}

If $R$ satisfies condition $(*)$ then the proofs of Theorems \refpage{overlap-theorem} and \refpage{overlap-theorem-2} show that $\SETT{e_p(t), e_q(t')}$ cannot be a $(p,q)$-partial element for $(p,q)=(1,2), (n-2,n-1)$ or $(n-1,n)$ unless $t=t'$. Since the remaining dimension $n$ determinacy conditions $(p,q).n$ all have the property $q-p \geq 2$, the assumed coordinate pattern for $t$ and $t'$ implies that $t$ and $t'$ cannot be equal at all coordinates other $p$ and $q$.

To summarize:
\begin{newlem}
{\rm
\SL 

If $R \SBS \PLIST{V}{n}$ satisfies property $(*)$ then $\MC{S}[y^R]$ is the $n$-truncation of an $n$-dimensional hypergroupoid.
}
\end{newlem}
\qed

Next, here is a simple method to produce examples of $R$ satisfying condition $(*)$. For notational convenience we assume that the $V_j$ are all $\mathbb{N}$ and, as an example, we'll take $n=6$. The general case will be evident from this.
\MS

\begin{example}
We specify that $R \SBS \PLIST{V}{6}$ contains any/all of the following elements:
\[
\begin{aligned}
(0,1,2,3,4,5,6)\\
(7,1,2,3,4,5,6)\\
(8,9,2,3,4,5,6)\\
(10,11,12,3,4,5,6)\\
(13,14,15,16,4,5,6) \\
(17,18,19,20,21,5,6) \\
(22,23,24,25,26,27,6)
\end{aligned}
\]
It is clear that this $R$ satisfies property $(*)$ because the initial segments of any two of elements overlap not at all. We may extend $R$ by adding any and all of the following:
\[
\begin{aligned}
(28,29,30,31,32,33,34)\\
(35,29,30,31,32,33,34) \\
(36,37,30,31,32,33,34) \\
(38,39,40,31,32,33,34)\\
(41,42,43,44,32,33,34)\\
(45,46,47,48,49,33,34)\\
(50,51,52,53,54,55,34)
\end{aligned}
\]
The set of elements of this second grouping all satisfy $(*)$ and the union of the two groups also satisfies $(*)$ because no element in the second group overlaps at all with any member of the first. This process can be extended arbitrarily.
\MS

\NI {\bf Variation:} Each of the two ``blocks'' of 7 elements of $R$ in the example follows the same pattern and involves the integer ranges $[0,27]$ and $[28,55]$. This process also yields an infinite example of an $R$ such that $\MC{S}[y^R]$ is a 6-dimensional hypergroupoid. The $(k-1)$'st ``block'' of elements of the infinite version of $R$ would be
\[
\begin{aligned}
(28k, 28k+1, 28k+2, 28k+3, 28k+4, 28k+5, 28k+6)\\
(28k+7,28k+1, 28k+2, 28k+3, 28k+4, 28k+5, 28k+6)\\
(28k+8,28k+9, 28k+2, 28k+3, 28k+4, 28k+5, 28k+6)\\
(28k+10,28k+11, 28k+12, 28k+3, 28k+4, 28k+5, 28k+6)\\
(28k+13,28k+14, 28k+15, 28k+16, 28k+4, 28k+5, 28k+6)\\
(28k+17,28k+18, 28k+19, 28k+20, 28k+21, 28k+5, 28k+6)\\
(28k+22,28k+23, 28k+24, 28k+25, 28k+26, 28k+27, 28k+6)
\end{aligned}
\]
\end{example}

\subsection{Computation notes}

\subsubsection{Subcomplexes generated by an $n$-simplex}
Suppose $\MC{C}$ is a simplicial set $y \in \MC{C}_n$ and $\MC{C}[y]$ denotes the $n$-truncated complex consisting of $y$ and all the simplicial images of $y$ in dimensions $\leq n$. The number $|\MC{C}[y]|$ of simplices in $\MC{C}[y]$ depends on $n$, as follows.
\BS

Given $m \leq n$ then each non-degenerate $x\in \MC{C}[y]_m$ is determined by which vertices of $y$ are vertices of $x$. There are $\binom{n+1}{m+1}$ size $m+1$ subsets of the vertex set of $y$. That is: the number of non-degenerate $m$-simplices in $\MC{C}[y]$ is $\binom{n+1}{m+1}$ and therefore the total number of non-degenerate simplices of $\MC{C}[y]$ is
\[
\sum_{m=0}^n \binom{n+1}{m+1} = 2^{n+1}-1
\]

Now if $x\in \MC{C}[y]_m$ is non-degenerate then the set of degenerate images of $x$ of dimension $m+k$ is
\[
\SETT{
s_{q_k} \DDD{}s_{q_1} (x): q_1 \DDD{<} q_k \in [m+k-1]
}
\]
which has $\binom{m+k}{k}$ elements. Therefore, the total number of degenerate simplicial images of $x$ which are of dimension $\leq n$ is
\[
\sum_{k=1}^{n-m} \binom{m+k}{k}
\]
and the total number of degenerate simplices which come from dimension $m$ non-degenerate simplices is 
\[
\binom{n+1}{m+1} \sum_{k=1}^{n-m} \binom{m+k}{k}
\]
Using the standard identities
\[
\sum_{t=0}^q \binom{c+t}{t} = \binom{c+q+1}{q} \text{ and }
\binom{a}{b} = \binom{a}{a-b}
\]
we get
\[
\sum_{k=1}^{n-m} \binom{m+k}{k} = \binom{n+1}{n-m}-1
=\binom{n+1}{m+1}-1
\]
Therefore the total number of degenerate simplices of $\MC{C}[y]$ is
\[
\begin{aligned}
 \sum_{m=0}^{n-1} \Bigl[
\binom{n+1}{m+1} \sum_{k=1}^{n-m} \binom{m+k}{k}
\Bigr]
& = \sum_{m=0}^{n-1} \Bigl[
\binom{n+1}{m+1} \Bigl(\binom{n+1}{m+1}-1 \Bigr)
\Bigr]\\
&= \sum_{m=0}^{n-1} {\binom{n+1}{m+1}}^2 - \sum_{m=0}^{n-1}\binom{n+1}{m+1}\\
&=\sum_{m=0}^{n-1} {\binom{n+1}{m+1}}^2 - (2^{n+1}-2)
\end{aligned}
\]
Therefore, since $\MC{C}[y]$ has $2^{n+1}-1$ non-degenerate simplices then the total number of simplices of $\MC{C}[y]$ is
\[
\Bigl|\MC{C}[y]\Bigr| = 2^{n+1}-1 + \sum_{m=0}^{n-1} {\binom{n+1}{m+1}}^2 - (2^{n+1}-2)
= 1 + \sum_{m=0}^{n-1} {\binom{n+1}{m+1}}^2
\]
We get a slightly cleaner-looking formula by changing the sum variable to $p=m+1$:
\[
\Bigl|\MC{C}[y]\Bigr| =1 + \sum_{p}^{n-1} {\binom{n+1}{p}}^2
\]
Some sample values:
\[
\begin{array}{c|r}
n & 1+\sum_{p=1}^n \binom{n+1}{p}^2 \\
\hline
2 & 19\\
3 & 69\\
4 & 251 \\
5 & 922 \\
6 & 3431 \\
7 & 12869\\
8 & 48619
\end{array}
\]
\BX

\subsubsection{Checking subcomplexes efficiently}
\label{efficient model-checking}

Given $n$ and $i \in [n+1]$, the sufficient conditions spelled out in Theorems \refpage{composer-model;i=0}, \refpage{composer-model;i=1}, \refpage{composer-model;1<i<n}, \refpage{composer-model;i=n} and \refpage{composer-model;i=n+1} involve one determinacy condition and either one or two surjectivity conditions that must hold for the model $\MC{S}'$ obtained in any one of those theorems. The determinacy condition for a given $n$ and $i$ implies several lower-dimensional determinacy conditions simply using the Down Rule. The resulting set of determinacy conditions is structured so that the dimension $m<n$ conditions required of $m$-simplices imply that the {\em degeneracies} of those simplices satisfy all the dimension $m+1$ conditions. Now when $m<n$ each $x \in \MC{S}'_m$ is either non-degenerate or the degenerate image of a lower-dimensional non-degenerate simplex. Therefore, in order to verify that a $\MC{S}'$ satisfies the sufficent determinacy conditions, it suffices to verify that for each $m<n$  the {\em non-degenerate} $m$-simplices of that subcomplex satisfy the dimension $m$ conditions.
\MS

\begin{example}
{\rm
Let $(n,i)=(6,3)$. The sufficient determinacy condition is $\DPQ{1}{4}{6}$, and the Down Rule yields the following set of conditions: 
\[
\Bigl\{\DPQ{1}{3}{4}, \DPQ{1}{3}{5}, \DPQ{1}{4}{5}, \DPQ{1}{4}{6}\Bigr\}
\]
To verify that a 6-truncated complex $\MC{S}'_{\leq 6}$ satisfies the determinacy conditions to be the 6-truncation of a \NICOMP{6}{3}, it suffices to verify that all non-degenerate $x \in \MC{S}'_4$ are $\DPQ{1}{3}{4}$, all non-degenerate $y \in \MC{S}'_5$ are both $\DPQ{1}{3}{5}$ and $\DPQ{1}{4}{5}$ and that all non-degenerate $z \in \MC{S}'_6$ are $\DPQ{1}{4}{6}$.

This abbreviated verification saves a good deal of computation. Given a non-degenerate $z \in \MC{S}_6$, the 6-truncated subcomplex generated by $z$ has $3431$ simplices of which 
462 are of dimension 4,\, 924 are of dimension 5 and 1716 are of dimension 6. Only $2^7-1=127$ simplices are non-degenerate. Of these, there are $\binom{7}{5}=21$ non-degenerate 4-simplices, 7 non-degenerate 5-simplices and just one non-degenerate 6-simplex ($z$, itself).
}
\end{example}


\section{Algebra of \NICOMP{n}{i}s}
\label{ni-composers-algebra}

\subsection{Definitions and basic facts}

We begin by recalling the definition of \NICOMP{n}{i} and establishing some terminology. The rest of the section develops some basic facts about \NICOMP{n}{i}s.

In the definition and discussion below, $\BOX{i}{m}{C}$ denotes the set of open $i$-horns of the simplicial set $C$ (Appendix, page \pageref{open-i-horn}) and $\Delta^\bullet(n)(C)$ denotes the $n$'th simplicial kernel of the simplicial set $C$ (Appendix, page \pageref{simplicial-kernel}).

\begin{newdef}
{\bf (\NICOMP{n}{i})}
\label{nicatdef}
\index{\NICOMP{n}{i}}
\index{$(n,i)$-composition}
\index{$(n,i)$-factor}
\index{$(n,i)$-map}
{\rm  
\SL

\begin{enumerate}

\item Let $n \geq 1$ and $i \in [n+1]$. An {\bf \NICOMP{n}{i}} is a simplicial set $C$ such that for all 
$m>n$, 
\[
\begin{array}{c}
\phi_{m,i}: C_m \to \BOX{i}{m}{C} \\
y \mapsto (d_0y, \cdots, \widehat{d_iy}, \cdots , d_my)
\end{array}
\]
is an isomorphism. We will denote
\[
\phi_{m,i}^{-1}(\OM{y}{m}{i})
\]
by
\[
\COMP{m}{i}\left( \OM{y}{m}{i} \right)
\]
That is, for each $m \geq n+1$ and $w \in C_m$, then \[
w = \COMP{m-1}{i} \left(d_0(w), \cdots , \widehat{d_i(w)}, \cdots , d_m(w)\right)
\]

\item
Given an \NICOMP{n}{i} $C$, we will call $(n-1)$-simplices {\bf ``$(n,i)$-objects''}, call $n$-simplices {\bf ``$(n,i)$-factors''}, and call $(n+1)$-simplices {\bf ``$(n,i)$-compositions''}.

If $w \in C_{n+1}$
then $d_i(w)$ will be called the {\bf result} of the composition, or the {\bf composite}
of $(d_0(w), \cdots , \widehat{d_i(w)}, \cdots , d_{n+1}(w))$.

\item In anticipation of Theorem \refpage{function complex theorem} below: If $C$ and $D$ are \NICOMP{n}{i}s then an {\bf $(n,i)$-map from $C$ to $D$} is
a simplicial map $C \times \Delta[n-1] \to D$.

\end{enumerate}

}
\end{newdef}
\BX

\begin{newdef}
{\bf (Hypergroupoid)}
\label{hypergroupoid}
\index{hypergroupoid}
{\rm
\SL

Let $n \geq 1$. A simplicial set $C$ is an {\bf $n$-dimensional hypergroupoid} if for all $m>n$ and for all $i \in [m]$ the
map
\[
\phi_{m,i} : C_m \to \Lambda^i(m)(C)
\]
is an isomorphism. 
}

\BX
\end{newdef}
\BS

\NI Comment: The condition that $\phi_{m,i}: C_m \to \BOX{i}{m}{C}$ is
an isomorphism can be expressed as a unique extension condition:
\begin{quote}
For each $f : \Lambda^i[m] \to C$ there is a unique $g: \Delta[m] \to C$ such that
$f$ is $\Lambda^i[m] \stackrel{\text{incl}}{\to} \Delta[m] \stackrel{g}{\to} C$.
\end{quote}
(Appendix, page \pageref{standard-m-simplex} for definitions of $\Delta[m]$ and $\Lambda^i[m]$).
\BX
\MS

Here are some basic facts about \NICOMP{n}{i}s following directly from the definition:

\begin{newlem}
\label{immedlemma}
{\rm 
\SL 

\begin{enumerate}

\item An ordinary small category (identified with its nerve) is a \NICOMP{1}{1}.

\item If $C$ is an \NICOMP{n}{i} then for every $m > n$, $C$ is also an \NICOMP{m}{i}.

\item If $C$ is an $n$-dimensional hypergroupoid, then for every $m>n$, $C$ is also an $m$-dimensional hypergroupoid.

\item Suppose $C$ is an \NICOMP{n}{i}. Then: for each $m \geq n$, the map
\[
\BOX{i}{m+1}{C} \to \Delta^\bullet(m+1)(C)
\]
defined by $(\OM{y}{m+1}{i}) \mapsto \LIST{y}{m+1}$ is monic, where $y_i = d_i\bigl(\phi_{m,i}^{-1}(\OM{y}{m+1}{i}) \bigr)$

\end{enumerate}
}
\end{newlem}

\qed

\begin{newlem}
{\bf (Coskeleton lemma)}
\index{coskeleton lemma}
\index{lemma!Coskeleton lemma}
\label{trunclemma}
{\rm
\SL

\NI Let $C$ be a simplicial set. Then
\begin{enumerate}

\item  If $C$ is an \NICOMP{n}{i} then for all $m \ge n+2$, the map $C_m \to \Delta^\bullet(m)(C)$, $z \mapsto \Bigl( d_0(z), \cdots, d_m(z) \Bigr)$ is an isomorphism. Therefore,
any \NICOMP{n}{i} is determined by its truncation to dimension $n+1$ and $C \cong \text{Cosk}^{n+1}(C)$.

\item If $C$ is a simplicial set such that $C_{n+1} \ISO \BOX{i}{n+1}{C}$ and $C = \CK^{n+1}(C)$ then $C$ is a \NICOMP{n}{i}.

\item
If $C$ is an \NICOMP{n}{i} then for all $m \geq n+1$, any $z \in C_m$ is determined uniquely by its dimension $n$ subfaces.

\item Suppose $m \geq 1$ and $C = \CK^{m}(C)$. That is, for all $k>m$, $C_k = \Delta^\bullet(k)(C)$. Then $C$ is an $(m+1)$-dimensional hypergroupoid.

\item (``Truncation property'') 
\label{truncprop}
Given any simplicial set $X$, and \NICOMP{n}{i} $C$, then any truncated simplicial map $\TR_{n+1}(X) \to \TR_{n+1}(C)$ extends uniquely to a simplicial map $X \to C$.

\item (``Hypergroupoid property'')
If $C$ is an \NICOMP{n}{i} then $C$ is an $(n+2)$-dimensional hypergroupoid. 

\end{enumerate}
}
\end{newlem}
\MS

\Proof

\begin{enumerate}

\item

Since $C$ is an \NICOMP{m}{i} for each $m \geq n$, it suffices to consider just the case $m=n+2$. By the previous lemma, 
\[
\BOX{i}{n+2}{C} = C_{n+2} \XRA{\phi_{n+2}} \Delta^\bullet(C)
\]
is monic. 

Now, to show that $\phi_{n+2}$ is surjective,
suppose 
\[
\LIST{w}{n+2} \in \Delta^\bullet(n+2)(C)
\]
and 
\[
z = \COMP{n+1}{i}( w_0 \DDD{,} \widehat{w_i} \DDD{,} w_{n+2}) \in C_{n+2}
\]
with $d_i(z)$ denoted by $w'_i \in \C_{n+1}$ and 
\[
w'_i = \COMP{n}{i}\left(d_0(w'_i), \cdots , \widehat{d_i(w'_i)}, \cdots , d_{n+1}(w'_i) \right)
\]
It remains to show that $w'_i = w_i$.

Now when $p<i$ then
\[
d_p(w'_i) = d_{i-1}(w_p) = d_p(w_i)
\]
and when $p>i$ then
\[
d_p(w'_i) = d_i(w_{p+1}) = d_p(w_i)
\]
Therefore
\[
\begin{aligned}
w_i &=
\COMP{n}{i}\BIGP{
d_0(w_i), \cdots , \underset{i}{-}, \cdots , d_{n+1}(w_i)} \\
&=
\COMP{n}{i}\BIGP{
d_0(w'_i), \cdots , \underset{i}{-}, \cdots , d_{n+1}(w'_i)} \\
&= w'_i
\end{aligned}
\]
and thus $\phi_{n+2}$ is surjective.

\item 
We must show that for all $m \geq n+1$
\[
C_m \to \Lambda^i(m)(C), \qquad 
u \mapsto \bigl(d_0(u), \cdots, \widehat{d_i(u)}, \cdots , d_m(u) \bigr)
\]
is bijective. 

It is bijective for $m=n+1$, by hypothesis. When $m>n+1$ then, again by hypothesis, $C_m = \Delta^\bullet(m)(C)$. Therefore, $t \in C_m$ is
\[
t = \LIST{t}{m} = ( d_0(t) \DDD{,} d_m(t))
\]
and the map 
\[
\phi_{m,i}:C_m \to \Lambda^i(m)(C)
\]
is
\[
\bigl( t_0, \cdots , t_i, \cdots , t_m \bigr) \mapsto
\bigl( t_0, \cdots , \widehat{t_i}, \cdots , t_m \bigr)
\]
which is clearly monic.
\MS

Now given any $\bigl( \OM{t}{m}{i} \bigr) \in \Lambda^i(m)(C)$ there is a unique 
\[
t_i  = \Bigl(d_0(t_i) \DDD{,} d_{m-1}(t_i) \Bigr) \in C_{m-1}
\]
such that $\bigl(t_0, \cdots , t_i, \cdots , t_m \bigr) \in \Delta^\bullet(m)(C)$. It is determined uniquely by the face identities:
\[
d_p(t_i) = \left\{
\begin{array}{ll}
d_{i-1}(t_p) & \text{if } p<i\\
& \\
d_i(t_{p+1}) & \text{if } p \geq i
\end{array}
\right.
\]
When $m>n+1$ then $\LIST{t}{m} \in \Delta^\bullet(m)(C) = C_m$. Therefore, $C_m \to \Lambda^i(m)(C)$ is surjective.

\item 
Suppose $z \in C_m$ where $m \geq n+1$. If $m=n+1$ 
\[
\BOX{i}{n+1}{C} \to \Delta^\bullet(n+1)(C)
\]
is monic, as noted in the previous lemma. Therefore every $z \in C_{n+1}$ is determined uniquely by its dimension $n$ subfaces. 

Now reason inductively. $C$ is an \NICOMP{m}{i} for each $m>n$. Therefore, any $z \in C_{m+1}$ is determined uniquely by its dimension $m$ subfaces which, by induction, are each determined uniquely by their dimension $n$ subfaces.

\item
First we will show, for all $i \in [m+2]$, 
\[
\phi_{m+2,i}: C_{m+2} \to  \Lambda^i(m+2)(C)
\]
is an isomorphism. The same reasoning will apply $\phi_{n,i}$ for all $n>m+2$ and all $i \in [n]$.
\MS

Given $y \in C_{m+2}$ then, by hypothesis,
\[
y = \LIST{y}{m+2} \in \Delta^\bullet(m+2)(C) \qquad y_q = d_q(y)
\]
Fix any $i \in [m+2]$ and observe that 
\[
d_q(y_i) = \left\{
\begin{array}{ll}
d_{i-1}(y_q) & \text{if } 0 \leq q<i\\
&\\
d_i(y_{q+1}) & \text{if } i \leq q\leq m+1
\end{array}
\right.
\]
The map $\phi_{m+2,i} : C_{m+2} \to \Lambda^i(m+2)(C)$ defined by
\[
y = \LIST{y}{m+2} \mapsto (\OM{y}{m+2}{i})
\]
is clearly monic. The inverse is defined by
\[
\phi_{m+2,i}^{-1}(\OM{y}{m+2}{i}) = (y_0, \cdots , y_i, \cdots , y_{m+2})
\]
where
\[
y_i \DFAS \BIGP{
d_{i-1}(y_0), \cdots , d_{i-1}(y_{i-1}),\;  d_i(y_{i+1}), \cdots , y_{m+2}
}
\]
(with obvious adjustments in the cases $i=0$ and $i =m+2$).

Since this reasoning applies to each $i \in [m+2]$ we have proved that, for all $i \in [m+2]$, $\phi_{m+2,i}$ is an isomorphism. The same argument shows that $\phi_{n,i}$ is an isomorphism for all $n>m+2$ and all $i \in [n]$.

Therefore $C$ is an $(m+1)$-dimensional hypergroupoid.

\item This follows immediately from $C_{n+2}$ being the simplicial kernel of $d_0, \cdots , d_{n+1}: C_{n+1} \to C_n$.

\item
This follows from items 1 and 3.

\end{enumerate}
\qed
\MS

\NI {\bf Notes:}

(1) The Coskeleton lemma also applies to any $n$-dimensional hypergroupoid.
\MS

(2) It follows from the lemma on page \pageref{delta-n-is-hypergroupoid} that for every $m \geq 0$ and $n \geq 2$, then $\Delta[m]$ is an $n$-dimensional hypergroupoid.
\MS

\begin{newlem}
{\bf (Compositions lemma)}
\index{Compositions Lemma}
\index{lemma!Compositions lemma}
\label{compslemma}
{\rm
\SL

Let $C$ be an \NICOMP{n}{i}. Then
\begin{enumerate}

\item Given any sequence $w_0, \cdots , w_{i-1},\underset{i}{-},\; \underset{i+1}{-}, w_{i+1},\cdots , w_{n+2}$ of $(n+1)$-simplices of $C$ such that $d_p(w_q) = d_{q-1}(w_p)$ for all $p<q$ in $[n+1] - \SETT{i,i+1}$ then there is a unique $z \in C_{n+2}$ such that $d_p(z) = w_p$ for all $p \in [n+1]-\SETT{i,i+1}$.

\item $C_{n+2} \cong \BOX{i+1}{n+2}{C}$

\end{enumerate}
}
\end{newlem}

\NI {\bf Proof:} \MS

\begin{enumerate}

\item 
Suppose -- hypothetically -- that there exists $z \in C_{n+2}$ such that $d_p(z) = w_p$ for all $p \neq i,i+1$. Let $w_{i+1} \DFAS d_{i+1}(z)$. Then
\[
d_p(w_{i+1}) = \left\{
\begin{array}{ll}
d_i(w_p) & \text{if } p<i\\
& \\
d_{i+1}(w_{p+1}) & \text{if } p \geq i+1\\
\end{array}
\right.
\]
Therefore
\[
\bigl( \;d_0(w_{i+1}), \cdots, d_{i-1}(w_{i+1}),\; \underset{i}{\text{---}}\;, d_{i+1}(w_{i+1}), \cdots , d_{n+1}(w_{i+1})\; \bigr)
= 
\]
\[
\bigl( \; d_i(w_0), \cdots, d_i(w_{i-1}), \; \underset{i}{\text{---}} \;, d_{i+1}(w_{i+2}), \cdots, d_{i+1}(w_{n+2}) \; \bigr) \in \BOX{i}{n+1}{C}
\]
where the latter $i$-horn is determined by the given information without reference to the hypothetical $z$.

Then $(n,i)$-composition applied to this open $i$-horn determines a unique $w_{i+1} \in C_{n+1}$ such that
\[
(w_0, \cdots , w_{i-1}, \; \underset{i}{\text{---}}\; , w_{i+1}, w_{i+2}, \cdots , w_{n+2}) \in 
\BOX{i}{n+2}{C}
\]
From this, $(n+1,i)$-composition yields a unique $z\in C_{n+2}$ and $w_i \DFAS d_i(z) \in C_{n+1}$ such that
\[
z=(w_0, \cdots, w_{i-1}, w_{i},w_{i+1}, \cdots , w_{n+2}) \in C_{n+2}
\]

\item 
Given $\Bigl( w_0, \cdots , w_i, \; \underset{i+1}{\text{---}}\;, w_{i+2}, \cdots , w_{n+2} \Bigr) \in \BOX{i+1}{n+2}{C}$ then
\[
\Bigl(w_0, \cdots, w_{i-1},\underset{i}{-}, \underset{i+1}{-}, w_{i+2}, \cdots , w_{n+2} \Bigr) 
\]
satisfies the hypothesis of part (1). Therefore there exists a unique $z \in C_{n+2}$ such that $d_p(z) = w_p$ for all $p \in [n+2]-\SETT{i,i+1}$. It follows directly that the $(n,i)$-compositions $w_i$ and $d_i(z)$ have equal faces and therefore are equal. This verifies the claim.
\end{enumerate}
\qed

\NI See example \refpage{example-2-1-composer} below.
\MS

\begin{newcor}
\label{deduced-composer-structure}
{\rm
\SL

Suppose $C$ is an \NICOMP{n}{i}. Then 
\begin{enumerate}
\item 
$C$ is also an \NICOMP{n+1}{i+1}.

\item
For each $k \geq 0$ and each $j \in [k]$, $C$ is an \NICOMP{n+k}{i+j}.

That is, for all $m \geq n+k$ and $j \in [k]$, $C_m \cong \BOX{i+j}{m}{C}$.

\end{enumerate}
}
\end{newcor}

\Proof

1. If $C$ is an \NICOMP{n}{i} then $C$ is also an \NICOMP{m}{i} for all $m > n$. Part 2 of the Compositions Lemma states that $C_{m+2} \cong \BOX{i+1}{m+2}{C}$ holds for all $m>n$ i.e. $C_{m+1} \cong \BOX{i+1}{m+1}{C}$ for all $m > n+1$. That is, $C$ is an \NICOMP{n+1}{i+1}.
\MS

2. By repeated use of part 1 above: $C$ being an \NICOMP{n}{i} implies that for all $j \geq 0$, $C$ is an \NICOMP{n+j}{i+j}. It follows directly from the definition of \NICOMP{n}{i} that $C$ is an \NICOMP{n+k}{i+j} for all $k \geq j$.
%
%
%
%
%
%

\qed
\MS

To illustrate this corollary, suppose $C$ is a \NICOMP{6}{3}. Then: $C$ is a $(7,3)$- and \NICOMP{7}{4}. It is also an $(8,3)$-, $(8,4)$- and \NICOMP{8}{5}, etc.
%
%
%

\BX
\BS

Finally we note that if $C,D$ and $E$ are \NICOMP{n}{i}s, and $F:D \to C$, $G:E \to C$ are simplicial maps, then the pullback simplicial map $E \times_C D \to E$ is a simplicial map of \NICOMP{n}{i}s.
\begin{center}
\begin{picture}(1,.6)
\put(0,0){\MB{E}}
\put(1,0){\MB{C}}
\put(0,.5){\MB{E \times_C D}}
\put(1,.5){\MB{D}}
\put(.2,0){\vector(1,0){.6}}
\put(.35,.5){\vector(1,0){.45}}
\put(0,.375){\vector(0,-1){.25}}
\put(1,.375){\vector(0,-1){.25}}
\put(.5,.1){\MBS{G}}
\put(1.15,.25){\MBS{F}}
\end{picture}
\end{center}
\MS

In particular, any finite product of \NICOMP{n}{i}s is also an \NICOMP{n}{i}.

\subsection{Matrix notation for simplices}

Fix $n \ge 1$ and suppose $X$ is any simplicial set with $X = 
\text{Cosk}^n(X)$. Equivalently, for all $m>n$, $X_m = \Delta^\bullet(m)(X)$. 
It follows that if $y \in X_{m+1}$ with $m>n$ then 
\[
y = \bigl( d_0(y), \cdots , d_{m+1}(y)), \quad d_j(y) \in X_m
\]
and for each $j \in [m+1]$,
\[
d_j(y) = \bigl( d_0d_j (y), \cdots , d_md_j(y) \bigr) 
\]
Denote $d_p d_q(y)$ by $y_{q\,p}$.
We may then display $y$ with a faces-of-faces array:
\[
y = 
\left[
\begin{array}{cccc}
y_{00} & y_{01} & \cdots & y_{0\,m}\\
y_{10} & y_{11} & \cdots & y_{1\, m}\\
\vdots & \vdots & \ddots & \vdots \\
y_{m+1\, 0} & y_{m+1\, 1} & \cdots & y_{m+1\, m}
\end{array}
\right]
\]
where the entries occur in pairs of equal $(m-1)$-simplices according to the face identities
\[
\begin{array}{c}
d_p d_q (y) = d_{q-1} d_p (y)\\
y_{q\,p} = y_{p\,q-1}
\end{array}
\]
whenever $p<q$ in $[m+1]$.

For example, if $X = \text{Cosk}^1(X)$ and $y \in X_3$ then
\[
y = (y_0,y_1,y_2,y_3) = 
\left[
\begin{array}{ccc}
y_{10} & y_{20} & y_{30}\\
y_{10} & y_{21} & y_{31}\\
y_{20} & y_{21} & y_{32}\\
y_{30} & y_{31} & y_{32}
\end{array}
\right]
=
\left[
\begin{array}{ccc}
y_{00} & y_{01} & y_{02}\\
y_{00} & y_{11} & y_{12}\\
y_{01} & y_{11} & y_{22}\\
y_{02} & y_{12} & y_{22}
\end{array}
\right]
\]

\begin{example}
\label{example-2-1-composer}
To illustrate the use of this matrix notation, suppose $C$ is a \NICOMP{2}{1}. A typical
$w \in C_3$ is a $(2,1)$-composition
\[
w = \COMP{2}{1}(y_0, - , y_2, y_3) = (y_0,y_1,y_2,y_3) \in \Delta^\bullet(3)(C)
\]
where $y_1 = d_1(w) \in C_2$ defined by $\Lambda^1(3)(C) \xrightarrow{\phi_{2,1}^{-1}} C_3 \xrightarrow{d_1} C_2$

Now suppose $\bigl(\; w_0, \text{---}, \text{---}, w_3, w_4 \; \bigr)$ forms a compatible family of $(2,1)$-compositions i.e. for all $p<q$ in $[4]-\SETT{1,2}$, $d_p(w_q) = d_{q-1}(w_p)$. The
family may be displayed as a partially-filled matrix, where for each $p$, $w_p = (y_{p\;0}, y_{p\; 1}, y_{p\; 2},y_{p\; 3})$ and $d_p(w_q) = y_{q\, p}$.
\[
\left[
\begin{array}{c}
w_0 \\
- \\
- \\
w_3 \\
w_4\\
\end{array}
\right]
=
\left[
\begin{array}{cccc}
y_{00} & y_{01} & y_{02} & y_{03}\\
- & - & - & - \\
- & - & - & - \\
y_{30} & y_{31} & y_{32} & y_{33}\\
y_{40} & y_{41} & y_{42} & y_{43}\\
\end{array}
\right]
=
\left[
\begin{array}{cccc}
y_{00} & y_{01} & y_{02} & y_{03}\\
- & - & - & - \\
y_{01} & - & y_{32} & y_{42} \\
y_{30} & y_{31} & y_{32} & y_{33}\\
y_{40} & y_{41} & y_{42} & y_{43}\\
\end{array}
\right]
=
\left[
\begin{array}{cccc}
y_{00} & y_{01} & y_{02} & y_{03}\\
- & - & - & - \\
y_{01} & y_{21} & y_{32} & y_{42} \\
y_{30} & y_{31} & y_{32} & y_{33}\\
y_{40} & y_{41} & y_{42} & y_{43}\\
\end{array}
\right]
\]
By $(2,1)$-composition, the missing matrix entries are uniquely
determined as follows: row 2 $=w_2$ is determined by $(2,1)$-composition and the entire matrix is determined by $(3,1)$-composition. This is an instance of part (2) of the Compositions Lemma \refpage{compslemma}.
\end{example}

\subsection{Equational characteristics of \NICOMP{n}{i}s}
\index{$(n,i)$-associativity}

In this section we will examine aspects of $(n,i)$-composition related to associativity, units and cancellation.

\subsubsection{Associativity}
As noted in the introduction (page \pageref{axiomatic-equations}), the equational axioms for the theory of an \NICOMP{n}{i} $\C$ correspond to the simplicial identities of the $(n+2)$-truncation of $\C$.

For example, the axiom for associativity of composition of 1-simplices of a \NICOMP{1}{1} is the simplicial identity
\[
\C_3 \xrightarrow{d_2} \C_2 \xrightarrow{d_1} \C_1 =
\C_3 \xrightarrow{d_1} \C_2 \xrightarrow{d_1} \C_1
\]
In the usual notation, given a $(1,1)$-composition $w = (y_0,y_1,y_2) \in \C_2$ then $y_1 = y_0 \circ y_2$. Given any $(2,1)$-composition $z = (w_0, w_1, w_2, w_3) \in \C_3$, involving the $(1,1)$-compositions $w_j = d_j(z)$, then $d_p(w_q) = y_{q\,p}$ and $y_{q\, p} = y_{p\, q-1}$ whenever $p<q$. The array for $z$ is
\[
z =
\left[
\begin{array}{c}
w_0\\
w_1\\
w_2\\
w_3
\end{array}
\right]
=
\left[
\begin{array}{ccc}
y_{1\,0} & y_{2\,0} & y_{3\,0}\\
y_{1\,0} & y_{2\,1} & y_{3\,1}\\
y_{2\,0} & y_{2\,1} & y_{3\,2}\\
y_{3\,0} & y_{3\,1} & y_{3\,2}\\
\end{array}
\right]
\]
The identity $d_1 d_2(z) = y_{2\,1} = d_1 d_1(z)$ yields
\[
\begin{aligned}
y_{2\,1} &= y_{2\,0} \circ y_{3\,2} \qquad (\text{from } w_2)\\
&= (y_{1\,0} \circ y_{3\,0}) \circ y_{3\,2}\\
&= y_{1\,0} \circ y_{3\,1} \qquad (\text{from } w_1)\\
&= y_{1\,0} \circ ( y_{3\,0} \circ y_{3\,2})
\end{aligned}
\]

Now when $\C$ is an \NICOMP{n}{i} for arbitrary $n \geq 1$ and $i \in [n+1]$, the concept of associativity generalizes to the simplicial identity
\[
\C_{n+2} \xrightarrow{d_{i+1}} \C_{n+1} \xrightarrow{d_i} \C_n = 
\C_{n+2} \xrightarrow{d_{i}} \C_{n+1} \xrightarrow{d_i} \C_n
\]
Given any $(n+1,i)$-composition $z = \LIST{w}{n+2} \in \C_{n+2}$, involving $(n,i)$-compositions $w_j = (y_{j\,0}, \cdots , y_{j\,n+1})$, then we can draw a similar array for $z$. In the case when $0<i<n+1$ that array would look like the one below, emphasizing the entries in rows $i$ and $i+1$ with their corresponding column entries, and using ``$\bullet$'' for the other entries:
\[
z = \left[
\begin{array}{l}
w_0\\
\vdots \\
w_{i-1}\\
w_i\\
w_{i+1}\\
w_{i+2}\\
\vdots \\ 
w_{n+2}
\end{array}
\right]
=
\left[
\begin{array}{cccccccc}
\bullet & \cdots & \bullet & y_{i\,0} & \fbox{$y_{0\,i}$} & \bullet & \cdots & \bullet \\
\vdots & \cdots & \vdots & \vdots & \vdots & \vdots & \cdots & \vdots \\
\bullet & \cdots & \bullet & y_{i\, i-1} & \fbox{$y_{i-1\,i}$} & \bullet & \cdots & \bullet \\
y_{i\, 0} & \cdots & y_{i\, i-2}&  y_{i\, i-1}& \fbox{$y_{i\,i}$} & \fbox{$y_{i+2\,i}$} & \cdots &  \fbox{$y_{n+2\,i}$}\\
\fbox{$y_{0\, i}$} & \cdots & \fbox{$y_{i-2\, i}$}& \fbox{$y_{i-1\, i}$} & \fbox{$y_{i+1\, i}$} & y_{i+1\, i+1} & \cdots & y_{i+1\, n+1} \\
\bullet & \cdots &\bullet & \bullet & \fbox{$y_{i+2\, i}$} & \bullet & \cdots & \bullet \\
\vdots & \cdots & \vdots & \vdots & \vdots & \vdots & \cdots & \vdots \\
\bullet & \cdots &\bullet & \bullet & \fbox{$y_{n+2\, i}$} & \bullet & \cdots & \bullet \\
\end{array}
\right]
\]
Here, the boxed entries are $(n,i)$-composites. Written out, the simplicial identity $d_i d_{i+1} = d_i d_i$ says
\[
d_i( \fbox{$y_{0\,i}$}, \cdots, \fbox{$y_{i-1\,i}$}, y_{i+1\, i}, \cdots y_{i+1, n+1}) = 
d_i( 
y_{i\,0}, \cdots, y_{i\, n-2}, \fbox{$ y_{i\,i}$}, \fbox{$ y_{i+2\,i}$}, \cdots , \fbox{$ y_{n+2\,i}$}
)
\]

While the notation for ``$(n,i)$-associativity'' could be expanded and flattened to include all $\frac{(n+1)(n+2)}{2}$ $n$-simplices involved, analogous to the $n=1$, $i=1$ case, the result would be laborious to write and hard to read.

The array diagram above is perhaps the most concise and easiest visual display of the consequences of the simplicial identities.

\subsubsection{Cancellation}

In this section, we propose an extension to $(n,i)$-composition of the notions of monomorphism and epimorphism in ordinary categories. After the following definitions in an arbitrary simplicial set, we'll see how they relate to epics and monics in an arbitrary category.
\MS

In the following definitions $C$ is an arbitrary simplicial set, $n \geq 1$ and $A \subsetneq [n+1]$ is non-empty.

\begin{newdef}
\label{A-horn}
\index{A-horn}
\index{open A-horn}
\index{$\Lambda^A(n+1)(C)$}
{\bf ($A$-horn)}
{\rm
\SL

An {\bf  $A$-horn of dimension $n+1$} is any set 
\[
\SETT{y_t:t \in [n+1]-A,\; y_t \in C_n}
\]
which satisfies the {\bf face-compatibility conditions:} 
\[
\text{For all } p<q \text{ in } [n+1]-A, \quad d_p(y_q) = d_{q-1}(y_p)
\]
In effect, an $A$-horn consists of the $q$'th faces of a possible $(n+1)$-simplex, $q \notin A$.
\MS

\NI The set of $A$-horns of dimension $n+1$ will be denoted $\BOX{A}{n+1}{C}$.
\MS

Given an $A$-horn $\SETT{y_t:t \in [n+1]-A,\; y_t \in C_n}$, then we may also denote it by $(y_{t_0} \DDD{,} y_{t_k})$ where $\SETT{t_0 \DDD{<} t_k} = [n+1]-A$, or simply by $\PAR{w}$. If $\PAR{w} \in \BOX{A}{n+1}{C}$ then $A$ will be called the {\bf deletion subset} of $\PAR{w}$. 
\MS

\NI Notation: $A = \DELSUB{\PAR{w}}$.

\BX
}
\end{newdef}

\begin{newdef}
\label{deletion-subset}
\index{deletion subset}
\index{dsub}
{\bf ($A$-deletions and horn filling)}
{\rm
\SL

Define $\phi^A_{n+1} : C_{n+1} \to \BOX{A}{n+1}{C}$ by 
\[
\phi^A_{n+1}(w) \DFAS \SETT{d_p(w):p \in [n+1]-A}
\]
The $A$-horn $\phi^A_{n+1}(w)$ will be called the {\bf $A$-deletion of $w$.}
\MS

Suppose $B \SBS A \subseteq [n+1]$ and $B \neq \MT$. Then we define
\[
\phi^{B \SBS A}_{n+1} : \BOX{B}{n+1}{C} \to \BOX{A}{n+1}{C}
\]
by 
\[
\PAR{w} = \SETT{y_p:p \in [n+1]-B} \mapsto \SETT{y_p:p \in [n+1]-A}
\]
It follows that $\phi^A_{n+1} = \phi^{B \SBS A}_{n+1} \circ \phi^B_{n+1}$.
\MS

\label{fillable}
If $\PAR{w} \in \BOX{A}{n+1}{C}$, $w \in C_{n+1}$ and $\phi^A_{n+1}(w) = \PAR{w}$, then $w$ will be said to {\bf fill } $\PAR{w}$.
\index{fills a partial simplex}

\BX
}
\end{newdef}

\begin{newdef}
\label{cancellative}
{\bf ($A$-cancellative and $A$-cancelling)}
{\rm
\SL

Let $i \in A$.
The $A$-horn $\PAR{w} \in \BOX{A}{n+1}{C}$ is said to be {\bf \CAN{A} with respect to $i$} if whenever $w,w' \in C_{n+1}$ fill $\PAR{w}$ such that $d_i(w)=d_i(w')$ then $w=w'$. 

Also, if $w \in C_{n+1}$ and $\phi^A_{n+1}(w)$ is \CAN{A} with respect to $i$ then we'll say that {\bf $w$ is $A$-cancelling with respect to $i$}.

\BX
}
\end{newdef}
\MS

The motivation for these definitions comes from the following observations.
\MS

\begin{example}
{\rm
Let $C$ be a category i.e. a \NICOMP{1}{1} and $A = \SETT{0,1} \SBS [2]$. A typical $A$-horn in $\BOX{A}{2}{C}$ is $\PAR{w} = (-,-,y_2)$ where $y_2 \in C_1$ and the face compatibility requirements are vacuous. Suppose that $\PAR{w}$ has the following property: 
\begin{quote}
Whenever $w,w' \in C_2$ fill $\PAR{w}$ {\em and} $d_1(w) = d_1(w')$ then $w=w'$.
\end{quote}
That is, $\PAR{w}$ is $A$-cancelling with respect to $i=1$.
As diagrams, $w,w'$ are:

\begin{center}
\begin{picture}(3.5,.6)
\put(0,.25){\MB{w=}}
\put(.5,.5){\MB{v_0}}
\put(1,0){\MB{v_2}}
\put(1.5,.5){\MB{v_1}}
\put(.7,.5){\vector(1,0){.6}}
\put(1,.6){\MBS{y_2}}
\put(.6,.4){\vector(1,-1){.3}}
\put(.6,.25){\MBS{y_1}}
\put(1.4,.4){\vector(-1,-1){.3}}
\put(1.4,.25){\MBS{y_0}}

\put(2,.25){\MB{w'=}}
\put(2.5,.5){\MB{v_0}}
\put(3,0){\MB{v'_2}}
\put(3.5,.5){\MB{v_1}}
\put(2.7,.5){\vector(1,0){.6}}
\put(3,.6){\MBS{y_2}}
\put(2.6,.4){\vector(1,-1){.3}}
\put(2.55,.25){\MBS{y'_1}}
\put(3.4,.4){\vector(-1,-1){.3}}
\put(3.4,.25){\MBS{y'_0}}

\end{picture}
\end{center}
\MS

Here, $y_2=y'_2$ and the stated property of $(-,-,y_2)$ is that whenever $y_1=y'_1$ then $w=w'$.

When we translate this into familiar composition statements, it implies that whenever $y_0 y_2 = y_1=y'_1 = y'_0 y_2$ then $y_0=y'_0$. That is, $y_2$ is epic.

Conversely, if $y_2$ is epic then whenever it happens that you have $w = (y_0,y_0 \, y_2,y_2) \in C_2$ and $w' = (y'_0,y'_0 \, y_2,y_2) \in C_2$ with $y_0 y_2 = y'_0 y_2$ then $y_0=y'_0$, whence $w=w'$.
\MS

A similar exercise characterizes a monic $y_0 \in C_1$ in terms of 
$\BOX{\{1,2\}}{2}{C}$. 
}
\end{example}
\MS

This example suggests the following adaptation of the ``cancellation'' definition above.
\MS

\begin{newdef}
{\rm
\SL

Suppose $C$ is an \NICOMP{n}{i}, $n \geq 1$, $A \subsetneq [n+1]$ and $i \in A$.
\MS

We'll say $\PAR{w} \in \BOX{A}{n+1}{C}$ is {\bf $A$-cancellative with respect to $(n,i)$-composition} if it is \CAN{A} with respect to $i$.
\MS

If $w \in C_{n+1}$ and $\phi^A_{n+1}(w)$ is \CAN{A} with respect to $(n,i)$-composition then we'll say that $w$ is {\bf $A$-cancelling with respect to $(n,i)$-composition}.
}
\BX
\end{newdef}

\NI {\bf Comments:}
\MS

Let $C$ be an \NICOMP{n}{i}. 
\MS

1. If $A=\{i\} \subset [n+1]$ then $\BOX{A}{n+1}{C}= C_{n+1}$. Therefore every $\{i\}$-horn is $A$-cancellative with respect to $i$ and every $w \in C_{n+1}$ is $A$-cancelling with respect to $i$.
\MS

2. To rephrase ``$A$-cancelling'' for general $A \subsetneq [n+1]$ in the context of an \NICOMP{n}{i}: if $w,w' \in C_{n+1}$ such that $d_p(w)=d_p(w')$ for all $p \in \{i\} \cup \left([n+1]-A \right)$ then $w=w'$.
\BS

\begin{newlem}
{\rm
\SL

Let $C$ be a simplicial set, $n \geq 1$ and $B \SBS A \subsetneq [n+1]$ where $B \neq \MT$. Let $\PAR{w} \in \BOX{B}{n+1}{C}$. Then:
\MS

(i) $\PAR{w}$ is \CAN{B} implies $\phi^{B \SBS A}_{n+1}(\PAR{w})$ is \CAN{A}.
\MS

(ii) Let $w \in C_{n+1}$. If $w$ is $B$-cancelling then $w$ is $A$-cancelling.
}
\end{newlem}

\Proof

Suppose $w \in C_{n+1}$ fills $\PAR{w} \in \BOX{B}{n+1}{C}$ i.e. $\phi^B_{n+1}(w) = \PAR{w}$.
\setlength{\unitlength}{1in}
\begin{center}
\begin{picture}(2,.6)
\put(0,0){\MB{\BOX{B}{n+1}{C}}} 
\put(1,.5){\MB{C_{n+1}}}
\put(2,0){\MB{\BOX{A}{n+1}{C}}}
\put(.8,.4){\vector(-2,-1){.5}}
\put(.4,.4){\MBS{\phi^B_{n+1}}}
\put(1.2,.4){\vector(2,-1){.5}}
\put(1.6,.4){\MBS{\phi^A_{n+1}}}
\put(.5,0){\vector(1,0){1}}
\put(1,.125){\MBS{\phi^{B \SBS A}_{n+1}}}
\end{picture}
\end{center}
\BS

\NI Then $\phi^A_{n+1}(w) = \phi^{B \SBS A}_{n+1}\, \phi^B_{n+1}(w) = \phi^{B \SBS A}_{n+1}(\PAR{w})$. That is, $w$ fills $\phi^{B \SBS A}_{n+1}(\PAR{w})$. If $w$ is the only filler of $\PAR{w}$ then it is the only filler of $\phi^{B \SBS A}_{n+1}(\PAR{w})$.
\MS

\NI Therefore, as claimed: 
\MS

$w$ is $B$-cancelling implies $w$ is $A$-cancelling.
\MS

$\PAR{w}$ is \CAN{B} implies $\phi^{B \SBS A}_{n+1}(\PAR{w})$ is \CAN{A}.

\qed
\BS

Next we consider deletion subsets $A=\{i,j\}$ in the context of an \NICOMP{n}{i}.

\begin{newthm}
{\rm
\SL

Suppose $C$ is an \NICOMP{n}{i} and let $A = \SETT{i,j} \subsetneq [n+1]$ where $j \neq i$. Suppose $z \in C_{n+2}$. Then $d_i(z)$ is $A$-cancelling.
}
\end{newthm}

\Proof

Denote $d_p(z)$ by $w_p$.
\MS

To show that $w_i$ is $A$-cancelling suppose that $w'_i \in C_{n+1}$ is such that $d_p(w'_i)=d_p(w_i)$ for each $p \neq j$. We must show $w'_i=w_i$.
\MS

The rest of the proof splits into the cases $j>i$ and $j<i$.
\MS

\NI $\bullet$ Case: $j>i$:
\MS

Then
\[
\left(
w_0 \DDD{,} w'_i \DDD{,} w_j,-,w_{j+2} \DDD{,} w_{n+2}
\right) \in \BOX{j+1}{n+2}{C}
\]
because $w'_i$ satisfies all the face-compatibility conditions $w_i$ does except possibly for $d_j(w'_i) \overset{\text{??}}{=} d_i(w_{j+1}) = d_j(w_i)$. Since $C_{n+2} = \Delta^\bullet(n+2)(C)$, there is a unique $w'_{j+1} \in C_{n+1}$ such that
\[
z' = \left(
w_0 \DDD{,} w'_i \DDD{,} w_j,w'_{j+1},w_{j+1} \DDD{,} w_{n+2}
\right) \in C_{n+2}
\]
For all $p \neq i$,
\[
d_p(w'_{j+1}) = \left\{
\begin{array}{ll}
d_j(w_p) = d_p(w_{j+1}) & \text{if } p<j+1,\; p\neq i\\
& \\
d_{j+1}(w'_{p+1}) = d_{j+1}(w_{p+1}) & \text{if } p \geq j+1
\end{array}
\right.
\]
Therefore, by $(n,i)$-composition, $w'_{j+1}=w_{j+1}$, and by $(n+1,i)$-composition, $z'=z$ and $w'_i=w_i$.
\MS

\NI $\bullet$ Case $j<i$:
\MS

\NI The argument parallels the previous case.

First,
\[
\left(
w_0 \DDD{,} w_{j-1},-,w_{j+1} \DDD{,} w'_i,w_{i+1} \DDD{,} w_{n+2}
\right) \in \BOX{j}{n+2}{C}
\]
because $w'_i$ satisfies all the face-compatibility conditions $w_i$ does except possibly for $d_j(w'_i) \overset{\text{??}}{=} d_{i-1}(w_j) = d_j(w_i)$. Then there is a unique $w'_{j} \in C_{n+1}$ such that
\[
z' = \left(
w_0 \DDD{,} w_{j-1},w'_j,w_{j+1} \DDD{,} w'_i,w_{i+1} \DDD{,} w_{n+2}
\right) \in C_{n+2}
\]
For all $p \neq i$,
\[
d_p(w'_{j}) = \left\{
\begin{array}{ll}
d_{j-1}(w_p) = d_p(w_{j}) & \text{if } p<j\\
& \\
d_{j}(w_{p+1}) = d_{p}(w_j) & \text{if } p \geq j,\; p \neq i
\end{array}
\right.
\]
Therefore, by $(n,i)$-composition, $w'_{j+1}=w_{j+1}$, and by $(n+1,i)$-composition, $z'=z$ and $w'_i=w_i$.

\qed

\begin{newcor}
{\rm
\SL

Suppose $C$ is an \NICOMP{n}{i}, $k \geq 1$ and $A = \{i,j\} \SBS [n+k]$. Then:
Every $u \in C_{n+k}$ is $A$-cancelling. In particular, for all $q \in [n+k]$, $s_q(u)$ is $A$-cancelling.
}
\end{newcor}

\Proof

For each $k \geq 1$, $C$ is an \NICOMP{n+k}{i}. The previous theorem proved that for all $x \in C_{n+k+1}$ that $d_i(x)$ is $A$-cancelling. Therefore, for each $u \in C_{n+k}$, $u=d_is_i(u)$ is $A$-cancelling.

\qed

\subsection{Filling partial $m$-simplices; Extended Compositions Theorem}

The goal of this section is to prove that when $C$ is an \NICOMP{n}{i} or an $n$-dimensional hypergroupoid, then for certain dimensions $m>n$ and certain $B \subsetneq [m]$, partial $m$-simplices $\PAR z \in \BOX{B}{m}{C}$ can be filled uniquely (see page \pageref{fillable}) using the composition operation.

\begin{newdef}
\index{delsub}
\index{$\DELSUB{-}$}
\index{missing face of a partial simplex}
\index{known face of a partial simplex}
{\rm 
\SL

Let $C$ be a simplicial set and suppose $m \geq 2$. With reference to the  definition of horn (page \pageref{A-horn}), we will use the following terminology:
\begin{enumerate}
\item 
Given a proper non-empty subset $B \subsetneq [m]$, and a partial $m$-simplex $\PAR z \in \BOX{B}{m}{C}$ 
\[
\PAR z = \SETT{z_p : p \in [m]-B}
\]
then for each $p \in [m]-B$, we'll write $d_p(\PAR z) = z_p$ and refer to $z_p$ as a {\bf known face of $\PAR z$.}

\item
For each $p \in B$ we'll refer to $d_p(\PAR z)$ as a {\bf missing face of $\PAR z$} (without presuming that the missing face exists).

\end{enumerate}
\BX
}
\end{newdef}
\MS

To begin, we describe the relationship between the deletion subset $B$ of a $B$-horn and the deletion subsets of its missing faces. Consider the following situation when $m \geq 3$:
\[
\begin{aligned}
B & = \SETT{w_1 \DDD{<} w_k} \SBS [m] \text{ where } 2 \leq k \leq m \\
\PAR z & \in \BOX{B}{m}{C}
\end{aligned}
\]
For each $w_q \in B$, the missing face $d_{w_q}(\PAR z)$ is a partial $(m-1)$-simplex whose known faces and missing faces are determined by the face identities.

Using:
\[
d_t d_{w_q}(\PAR{z}) = \left\{
\begin{array}{ll}
 d_{w_q-1} d_t (\PAR{z}) & 0 \leq t<w_q\\
 & \\
 d_{w_q}d_{t+1} (\PAR{z}) & w_q \leq t \leq m-1
\end{array}
\right.
\]
it follows that $d_t d_{w_q}(\PAR{z})$ is missing either when $t<w_q$ and $t \in B$ or $t \geq w_q$ and $t+1 \in B$. Equivalently: $t \in B \cap [0,w_q-1]$ or $t+1 \in B \cap [w_q+1,m]$. 

That is:

\begin{newlem}
\label{delsub-lemma}
{\rm
Suppose $C$ is a simplicial set, $m \ge 3$, 
\[
B = \SETT{w_1 \DDD{<} w_k} \subsetneq [m] \quad \text{with } k \ge 2
\]
and $\PAR{z} \in \BOX{B}{m}{C}$. Then for each $w_q \in B$,
\begin{equation}
\DELSUB{d_{w_q}(\PAR{z})} = \Bigl( B \cap [0,w_q-1] \Bigr)
\cup \Bigl\{ r-1\; : r \in B \cap [w_q+1,m] \Bigr\}
\end{equation}
}
\end{newlem}
\qed
\MS

\begin{example}
{\rm
$m=6$, $B = \SETT{0,1,4,6} \in [6]$, $[6]-B = \SETT{2,3,5}$, $\PAR z \in \BOX{B}{6}{C}$ and the known faces of $\PAR z$ are $z_2, z_3, z_5 \in C_5$. The partial faces-of-faces array of $\PAR z$ is
\[
\begin{array}{|l|c|c|c|c|c|c||l|}
\hline
 & 0 & 1 & 2 & 3 & 4 & 5 & \text{deletion subset}\\
\hline
0 & & z_{20} & z_{30} &  & z_{50} & & \DELSUB{d_0(\PAR{z})} = \SETT{0,3,5}\\
\hline
1 &  & z_{21} & z_{31} & & z_{51} &  & \DELSUB{d_1(\PAR{z})} = \SETT{0,3,5}\\
\hline
2 & z_{20} & z_{21} & z_{22} & z_{23} & z_{24} & z_{25} & z_2\\
\hline
3 & z_{30} & z_{31} & z_{22} & z_{33} & z_{34} & z_{35} & z_3\\
\hline
4 &  & & z_{23} & z_{33} & z_{54} & &\DELSUB{d_4(\PAR{z})} = \SETT{0,1,5}\\
\hline
5 & z_{50} & z_{51} & z_{24} & z_{34} & z_{54} & z_{55} & z_5\\
\hline
6 & & & z_{25} & z_{35} & & z_{55} & \DELSUB{d_6(\PAR{z})} = \SETT{0,1,4}\\
\hline
\end{array}
\]
}
\end{example}
\MS

\NI{\bf Remarks:}
\MS

\NI 1. For each $w_q \in B$, $\DELSUB{d_{w_q}(\PAR{z})}$ has one less element than $\DELSUB{\PAR{z}}$.
\MS

\NI 2. Since $B \cap [0,w_1-1]=\MT$ and $B \cap [w_k+1,m] = \MT$ we have
\[
\DELSUB{d_{w_1}(\PAR{z})} = \SETT{w_2-1 \DDD{<} w_k-1}
\]
and 
\[
\DELSUB{d_{w_k}(\PAR{z})} = \SETT{w_1 \DDD{<} w_{k-1} }
\]

\NI 3. If $1<q<k$ then
\[
\begin{aligned}
B \cap [0,w_q-1] &= \SETT{w_1 \DDD{<} w_{q-1}}\\
B \cap [w_q+1,m] &= \SETT{w_{q+1} \DDD{<} w_k}
\end{aligned}
\]
which implies
\[
\DELSUB{d_{w_q}(\PAR{z})} = \SETT{w_1 \DDD{<} w_{q-1}} \cup
\SETT{w_{q+1}-1 \DDD{<} w_k-1}
\]

\NI 4. Items 2 and 3 can be combined by writing
\[
\DELSUB{d_{w_q}(\PAR{z})} = \SETT{w_1 \DDD{<} w_{q-1} < w_{q+1}-1 \DDD{<} w_k-1}
\]
where it is understood that $\SETT{w_1 \DDD{<} w_{q-1}} = \MT$ when $q=1$ and 
$\SETT{w_{q+1}-1 \DDD{<} w_k-1} = \MT$ when $q=k$.

\BX

\begin{newthm}
{\bf (Hypergroupoid filler)}
\label{hypgpd-filler-thm}
\index{Theorem!Hypergroupoid filler}
{\rm
\SL

Let $n \geq 1$, $k \geq 1$ and $C$ be an $n$-dimensional {\em hypergroupoid}. Suppose $B \SBS [n+k]$ and $1 \leq |B| \leq k$. Then every partial $(n+k)$-simplex $\PAR z \in \BOX{B}{n+k}{C}$ fills uniquely.

That is, $C_{n+k} \cong \BOX{B}{n+k}{C}$.
}
\end{newthm}

\Proof

We do induction on $k$.
\MS

\NI Case $k=1$:

In this case $B = \SETT{j}$ for some $j \in [n+1]$, and $\PAR z \in \BOX{j}{n+1}{C}$. Then $\PAR z$ fills uniquely by hypergroupoid composition.
\MS

\NI Case $k>1$:

In this case, we have $\DELSUB{\PAR z}= B = \SETT{q_1< \cdots < q_j}$ for some $j \leq k$. Let $q = q_j$. Then, by the previous lemma,
\[
\DELSUB{d_q(\PAR z)} = \SETT{q_1 < \cdots < q_{j-1}} = B' \SBS [n+(k-1)]
\]
That is, $B' \SBS [n+(k-1)]$, $|B'| = j-1 \leq k-1$ and $d_q(\PAR z) \in \BOX{B'}{n+(k-1)}{C}$. By induction on $k$, $d_q(\PAR z)$ fills uniquely and $\DELSUB{\PAR z}$ decreases to $\SETT{q_1 < \cdots < q_{j-1}}$. 
\MS

Repeat this reasoning for $q = q_{j-1}$ above:  thus $d_{q_{j-1}}(\PAR z)$ also fills uniquely. 
\MS

Repeat until the deletion subset of $\PAR z$ has been reduced to $\SETT{q_1}$. Then $\PAR z \in \BOX{q_1}{n+k}{C}$, and this fills uniquely by hypergroupoid composition.

\qed
\MS

\NI {\bf Alternative proof:}
\MS

Let $m>n$, $B \SBS [m]$ where $|B| \leq m-n$, and $\PAR z \in \BOX{B}{m}{C}$. Do induction on $|B|$.
\MS

If $|B|=1$ then $C_m \cong \BOX{B}{m}{C}$ by hypergroupoid composition.
\MS

If $|B|>1$ then for each $p \in B$, we have 
$\Bigl|\DELSUB{d_p(\PAR z}) \Bigr| = |B|-1$. By the induction hypothesis, $d_p(\PAR z)$ fills uniquely. This reduces the size of the deletion subset of $\PAR z$ by 1. After we repeat this reduction step for all but one of the remaining $p \in B$, the deletion subset of $\PAR z$ reduces to only one element. Then the hypergroupoid composition of $C$ fills $\PAR z$ uniquely.

\qed
\MS

\NI {\bf Comment:}
Note that in the last theorem, if $B \SBS [n+k] = [m]$ and $1 \leq |B| \leq k = m-n$ then $\Bigl|[n+k]-B \Bigr| \geq n+1$. Therefore we can paraphrase the last theorem for an $n$-dimensional hypergroupoid $C$ as follows: 
\begin{quote}
For any $m>n$, any partial $m$-simplex $\PAR z$  fills (uniquely) if at least $n+1$ of its faces are known.
\end{quote}

\BX
\MS

Suppose $C$ is an \NICOMP{n}{i}, $k \geq 2$, $m \geq n+k$ and $B \subsetneq [m]$ with $|B|=k$. Given $\PAR{z} \in \BOX{B}{m}{C}$, then as observed above:
\begin{itemize}
\item 
For each $p_1 \in \DELSUB{\PAR{z}}$, $d_{p_1}(\PAR z)$ is a partial $(m-1)$-simplex.

\item
For each $p_2 \in \DELSUB{d_{p_1}(\PAR z)}$, $d_{p_2} d_{p_1}(\PAR z)$ is a partial $(m-2)$-simplex.

\item
Etc.
\end{itemize}
If $\PAR y$ is any $(m-k+1)$-dimensional partial subface of $\PAR z$ then $\DELSUB{\PAR y}$ has just one element.
If that 1-element deletion subset happens to enable filling $\PAR y$ by $(n,i)$-composition, then the size of the deletion subsets of certain $(m-k+2)$-dimensional subfaces of $\PAR z$ decreases in size by 1. This sets up the possibility that $\PAR z$ itself fills by repeated use of $(n,i)$-composition applied to various partial subfaces.
\MS

\begin{example}
\label{extended-comp-example}
{\rm 
Suppose $C$ is a \NICOMP{3}{2}, $k=3$ and $B = \SETT{1,3,6} \SBS [3+3] \SBS [7]$ and $\PAR z \in \BOX{B}{7}{C}$, a partial $7$-simplex. Write $\PAR z = (z_0,-,z_2,-,z_4,z_5,-,z_7)$ to display the known faces of $\PAR z$. By Corollary \refpage{deduced-composer-structure}, $C$ is also a \NICOMP{n'}{2} for every $n' \geq 3$ and also a \NICOMP{n''}{3} for every $n'' \geq 4$.

We will calculate that $\PAR z$ fills by $(3,2)$-composition.
\MS

The partial array for $\PAR z$ is just below. Each entry is a known dimension 5 subface of $\PAR z$. For visual clarity, we indicate these known dimension 5 subfaces of $\PAR z$ simply by ``$\bullet$''.
\[
\PAR z = 
\begin{array}{|c|c|c|c|c|c|c|c||c|}
\hline
 & 0 & 1 & 2 & 3 & 4 & 5 & 6 & \\
\hline
0 & \BLT & \BLT & \BLT & \BLT & \BLT & \BLT & \BLT & \text{known}\\
\hline
(1) & \BLT & \BLT &      & \BLT & \BLT &      & \BLT & \DELSUB{d_1(\PAR{z})}=\SETT{2,5} \\
\hline
2 & \BLT & \BLT & \BLT & \BLT & \BLT & \BLT & \BLT & \text{known}\\
\hline
(3) & \BLT &      & \BLT & \BLT & \BLT &      & \BLT & \DELSUB{d_3(\PAR{z})}=\SETT{1,5}\\
\hline
4 & \BLT & \BLT & \BLT & \BLT & \BLT & \BLT & \BLT & \text{known}\\
\hline
5 & \BLT & \BLT & \BLT & \BLT & \BLT & \BLT & \BLT & \text{known}\\
\hline
(6) & \BLT &      & \BLT &      & \BLT & \BLT & \BLT & \DELSUB{d_6(\PAR{z})}=\SETT{1,3}\\
\hline
7 & \BLT & \BLT & \BLT & \BLT & \BLT & \BLT & \BLT & \text{known}\\
\hline
\end{array}
\]
Now consider the partial array for $d_1(\PAR z) \in \BOX{\{2,5\}}{6}{C}$:
\[
d_1(\PAR z) = 
\begin{array}{|c|c|c|c|c|c|c||c|}
\hline
  & 0 & 1 & 2& 3& 4& 5& \\
\hline
0 & \BLT & \BLT & \BLT & \BLT & \BLT & \BLT & \text{known}\\
\hline
1 & \BLT & \BLT & \BLT & \BLT & \BLT & \BLT & \text{known}\\
\hline
(2) & \BLT & \BLT & \BLT & \BLT &      & \BLT & \DELSUB{d_2 d_1(\PAR{z})} = \SETT{4}\\
\hline
3 & \BLT & \BLT & \BLT & \BLT & \BLT & \BLT & \text{known}\\
\hline
4 & \BLT & \BLT & \BLT & \BLT & \BLT & \BLT & \text{known}\\
\hline
(5) & \BLT & \BLT &      & \BLT & \BLT & \BLT & \DELSUB{d_5 d_1(\PAR{z})}= \SETT{2}\\
\hline
6 & \BLT & \BLT & \BLT & \BLT & \BLT & \BLT & \text{known}\\
\hline
\end{array}
\]
Observe that  $d_5 d_1(\PAR z)$ fills by $(4,2)$-composition and therefore 
\[
d_1(\PAR z) = \Bigl(d_0d_1(\PAR z), d_1d_1(\PAR z),-,d_3d_1(\PAR z),d_4d_1(\PAR z),d_5d_1(\PAR z),d_6d_1(\PAR z)
\Bigr) \in \BOX{2}{6}{C}
\]
fills by $(5,2)$-composition.

Returning to the partial array of $\PAR z$, the filling of $d_1(\PAR z)$ modifies the array for $\PAR z$ as follows:
\[
\PAR z = 
\begin{array}{|c|c|c|c|c|c|c|c||c|}
\hline
 & 0 & 1 & 2 & 3 & 4 & 5 & 6 & \\
\hline
0 & \BLT & \BLT & \BLT & \BLT & \BLT & \BLT & \BLT & \text{known}\\
\hline
1 & \BLT & \BLT & \BLT & \BLT & \BLT & \BLT & \BLT & \text{now known} = z_1\\
\hline
2 & \BLT & \BLT & \BLT & \BLT & \BLT & \BLT & \BLT & \text{known}\\
\hline
(3) & \BLT & \BLT & \BLT & \BLT & \BLT &      & \BLT & \DELSUB{d_3(\PAR{z})}=\SETT{5}\\
\hline
4 & \BLT & \BLT & \BLT & \BLT & \BLT & \BLT & \BLT & \text{known}\\
\hline
5 & \BLT & \BLT & \BLT & \BLT & \BLT & \BLT & \BLT & \text{known}\\
\hline
(6) & \BLT & \BLT & \BLT &      & \BLT & \BLT & \BLT & \DELSUB{d_6(\PAR{z})}=\SETT{3}\\
\hline
7 & \BLT & \BLT & \BLT & \BLT & \BLT & \BLT & \BLT & \text{known}\\
\hline
\end{array}
\]
The modified $d_6(\PAR z)$ fills (with $z_6 \in C_6$) by $(5,3)$-composition and so the original partial array $\PAR z$ now is
\[
\PAR z = ( z_0,z_1,z_2,-,z_4,z_5,z_6,z_7) \in \BOX{3}{7}{C}
\]
and this fills by $(6,3)$-composition.
}
\end{example}
\MS

The next definition and theorem address this kind of scenario precisely.
\MS

\begin{newdef}
{\bf (Type $p$ for $(n,i)$)}
\index{deletion subset of type $p$}
\index{type $p$ deletion subset}
\label{type-p-deletion-subset}
{\rm
\SL

Suppose $n \geq 1$, $i \in [n+1]$, $k \geq 2$ and 
\[
B = \SETT{ t_1 < \cdots < t_k} \SBS [n+k]
\] 
Let $p \in [1,k]$. We say that \fbox{$B$ {\bf is of type $p$ for $(n,i)$}} if $t_p = i+p-1$.
\MS

That is, $B$ is of type 1 if $t_1=i$, of type 2 if $t_2=i+1$, of type 3 if $t_3 = i+2$, etc.
}

\BX
\end{newdef}

\begin{newthm}
{\bf (Extended Compositions Theorem)}
\index{Extended Compositions Theorem}
\index{Theorem!Extended Compositions Theorem}
\index{Compositions Theorem}
\label{extended-comps-theorem}
{\rm 
\SL

Suppose $C$ is an \NICOMP{n}{i} where $n \geq 1$. Suppose $k \geq 2$, $m \geq n+k$ and $\PAR z \in \BOX{B}{m}{C}$ 
with 
\[
\DELSUB{\PAR z} = B = \SETT{t_1< \cdots < t_k} \SBS [n+k]\SBS [m]
\]
Then $(n,i)$-composition determines a unique filler for
$\PAR z$ if $B$ is of type $p$ for $(n,i)$ for some $p \in [1,k]$.
}
\end{newthm}

\Proof

First we will prove the claim for $p=1$. We do this by induction on $k \geq 2$.
\MS

\NI Case $k=2$:

In this case, $B = \SETT{i< t_2}$, and $\PAR z \in \BOX{B}{m}{C}$ for some $m \geq n+2$. Then $\DELSUB{d_{t_2}(\PAR z)} = \SETT{i}$; that is, $d_{t_2}(\PAR z) \in \BOX{i}{m}{C}$. So $d_{t_2}(\PAR z)$ fills uniquely by $(m-2,i)$-composition and the revised deletion subset of $\PAR z$ is $B' = \SETT{i}$. Therefore $\PAR z$ fills uniquely by $(m-1,i)$-composition. 
\MS

\NI Case $k>2$:

In this case, $B = \SETT{i < t_2< \cdots < t_k}$ and therefore
\[
\DELSUB{d_{t_k}(\PAR z)} = \SETT{i < t_2 < \cdots < t_{k-1}}
\]
By the induction hypothesis, $d_{t_k}(\PAR z)$ fills uniquely by $(m-2,i)$-composition and this revises the deletion subset of $\PAR z$ to
\[
B' = \SETT{i < t_2<t_3 < \cdots < t_{k-1}}
\]
Using that $B'$ is of type 1 for $(m-1,i)$ the induction hypothesis on $k$ implies that $\PAR z$ fills uniquely by $(m-1,i)$-composition. This proves the claim for $p=1$.
\BS

Next we prove the claim when $1<p \leq k$ by a terminating induction argument on $p$.
\MS

Suppose that $\PAR z \in \BOX{B}{m}{C}$ where $m \geq n+k$ and
\[
B = \SETT{t_1< \cdots < t_k} \text{ with } t_p = i+p-1
\]
Then
\[
\begin{aligned}
\DELSUB{d_{t_1}(\PAR z)} &= \SETT{t_2-1< \cdots < t_k-1}\\
&= \SETT{t'_1 < \cdots < t'_{k-1}}\\
&= B'
\end{aligned}
\]
where $t'_q = t_{q+1}-1$. In particular, $t'_{p-1} = t_p-1 = i + (p-1)-1$. That is, $B'$ is of type $p-1$ for $(n,i)$. By the induction hypothesis on $p$, $d_{t_1}(\PAR z)$ fills uniquely by $(m-2,i)$-composition and the revised deletion subset of $\PAR z$ becomes
\[
\begin{aligned}
B'' &= \SETT{t_2< \cdots < t_k}\\
&= \SETT{t''_1 < \cdots < t''_{k-1}}
\end{aligned}
\]
where $t''_q = t_{q+1}$.

In particular, $t''_{p-1} = t_p = i+p-1 = (i+1)+(p-1)+1$. Therefore $B''$ is of type $p-1$ for $(n+1,i+1)$. Again by the induction hypothesis on $p$, and using that $C$ is an \NICOMP{n+1}{i+1}, $\PAR{z}$ fills uniquely by $(n+1,i+1)$-composition.

\qed
\BS

In Example \refpage{extended-comp-example}, with $(n,i)=(3,2)$, the deletion subset $B=\SETT{1,3,6}= \SETT{t_1,t_2,t_3}$ is of type 2 for $(3,2)$ since $3=t_2 = i+p-1 = 2+2-1$.
\MS

One more consequence of the ``Hypergroupoid Filler'' theorem (Theorem \refpage{hypgpd-filler-thm}:

\begin{newthm}
{\rm
\SL

Suppose $C$ is an \NICOMP{n}{i} and $\PAR y \in \BOX{B}{m}{C}$ where $m>n+2$. Then $\PAR y$ fills uniquely by $(n,i)$-composition if $|B| \leq m-(n+2)$.
}
\end{newthm}

\Proof

$C$ is an $(n+2)$-dimensional hypergroupoid, by the ``Coskeleton lemma" lemma \refpage{trunclemma}. By Theorem \refpage{hypgpd-filler-thm}, $\PAR{y}$ fills uniquely if $[m]-B$ (the number of known faces of $y$) is at least $(n+2)+1$. That is $m+1-|B| \geq n+3$, equivalently $|B| \leq m-(n+2)$, as claimed.

\qed

\section{Function complexes}

The goal of this section is to prove:

\begin{newthm}
{\bf (Function Complex Theorem)}
\label{function complex theorem}
\index{function complex theorem}
\index{Theorem!Function Complex Theorem}
{\rm
\SL

\fbox{If $C$ and $D$ are \NICOMP{n}{i}s then so is $D^C$.}
}
\end{newthm}

\NI As the proof is long, we will present it in a sequence of numbered paragraphs.
\MS

\begin{tabular}{llc}
\S & Content & Page\\
\hline
\ref{step1} & Initial setup & \pageref{step1}\\
\ref{def-g(y,1)} & Definition of $g_m(-,1_m)$ & \pageref{def-g(y,1)}\\
\ref{def-g(-,del_i)} & Definition of $g_{m-1}(-,\del_i)$ & \pageref{def-g(-,del_i)}\\
\ref{def-g(-,sj di (1))} & Definition of $g_m(-,s_j d_i(1_m))$ & \pageref{def-g(-,sj di (1))}\\
\ref{tech-lemma} & Technical lemma & \pageref{tech-lemma}\\
\ref{g is simplicial} & Proof that $g$ is simplicial &\pageref{g is simplicial}\\
\end{tabular}
\vskip1cm

\begin{paragr}
\label{step1}
{\bf Initial set-up}
\MS

Suppose $m \geq n+1$ and $i \in [n+1]$. To prove that $D^C_{m} \ISO \BOX{i}{m}{D^C}$ is equivalent to proving that every simplicial map $f: C \times \Lambda^i[m] \to D$ extends along $1_C \times \text{incl} : C \times \Lambda^i[m] \to C \times \Delta[m]$ to a uniquely determined $g : C \times \Delta[m] \to D$.
\begin{center}
\setlength{\unitlength}{1in}
\begin{picture}(1,.8)
\put(0,0){\MB{C \times \Lambda^i[m]}}
\put(0,.6){\MB{C \times \Delta[m]}}
\put(1,.6){\MB{D}}
\put(0,.15){\vector(0,1){.3}}
\put(-.25,.3){\MBS{1 \times \text{incl}}}
\put(.36,.12){\vector(3,2){.55}}
\put(.7,.25){\MBS{f}}
\multiput(.36,.6)(.05,0){10}{\MB{.}}
\put(.86,.6){\vector(1,0){0}}
\put(.6,.7){\MBS{g}}
\end{picture}
\end{center}
\MS

\NI We will use the following background facts. (Also see the appendix).
\begin{enumerate}
\item 
$\Delta[m]$ is generated by $1_{[m]}$, which we will write simply as $1_m$.

\item
$\Lambda^i[m]$ is the subcomplex of $\Delta[m]$ generated by the $(m-1)$-simplices $\SETT{\del_p: p \in [m]-\{i\}}$. In some places we will write $d_p(1_m)$ for $\del_p$.

\item
For all $j \leq m-2$, $\Lambda^i[m]_j = \Delta[m]_j$. The reason is that for all $p \in [m-1]$, 
\[
d_p d_i(1_m) = 
\left\{
\begin{array}{ll}
d_{i-1} d_p(1_m) \in \Lambda^i[m]_{m-2} & \text{ if } p<i\\
d_i d_{p+1}(1_m) \in  \Lambda^i[m]_{m-2} & \text{ if } p \geq i
\end{array}
\right.
\]
Each $j$-simplex $h$ of $\Delta[m]$ has the form
\[
h = s_{p_1} \DDD{} s_{p_t} d_{q_1} \DDD{} d_{q_{u}}(1_m), \quad
u-t=m-j \geq 2
\]
with $d_{q_{u-1}} d_{q_u}(1_m) \in \Lambda^i[m]_{m-2}$.
Therefore $h \in \Lambda^i[m]_j$.

\item
The only $m$-simplices of $\Delta[m]$ which do not belong to $\Lambda^i[m]_m$ are $1_m$ and $s_j d_i(1_m)$ for $j \in [m-1]$. 

\end{enumerate}
\BX
\MS

\NI We will see that $g_{m-1}$ and $g_m$ are uniquely determined by the requirement that 
\[
g \circ (1_C \times \text{incl}) = f
\]
The long part of the proof is the verification that $g$ is a simplicial map.
\MS

\NI Clearly:\begin{itemize}
\item 
For all $j \leq m-2$, $g_j = f_j$ necessarily.

\item 
Whenever $ \gamma \in \Lambda^i[m]_m \SBS \Delta[m]_m$, then $g_m(y,\gamma)$ must be $f_m(y,\gamma)$.

\end{itemize}

Finally, since $D \ISO \CK^{n+2}(D)$ and $m \geq n+1$, the simplicial map $g$ is determined by $g_0 \DDD{,} g_m$.
\end{paragr}
\BS

\begin{paragr}
\label{def-g(y,1)}
{\bf Definition of $g_m(-,1_m)$}
\MS

Given $f : C \times \Lambda^i[m] \to D$ and $y \in C_m$, we note that 
\[
\SETT{f_{m-1}(d_p(y), d_p(1_m)):p \in [m]-\{i\}} \in \BOX{i}{m}{D}
\]
\NI{\bf Definition:}
\begin{equation}
\label{eqn-defining-g(y,1)}
g_m(y, 1_m) \DFAS \COMP{m-1}{i}
\Bigl( 
f_{m-1}(d_0(y), d_0(1_m)) \DDD{,} \underset{i}{-} \DDD{,} f_{m-1}(d_m(y), d_m(1_m))
\Bigr)
\end{equation}
\BX
\MS

\NI This definition forced by the requirement that 
\[
d_p \, g_m(y, 1_m) = g_{m-1}(d_p(y), d_p(1_m))
\]
It remains to show that $d_i \, g_m(y,1_m) \ISIT g_{m-1}(d_i(y), d_i(1_m))$ once we have defined $g_{m-1}$ below in \S \ref{def-g(-,del_i)}.

\begin{newlem}
\label{basic-lemma}
{\rm
\SL

Suppose $y \in C_m$. Then:
\MS

\NI (1)\; If $i>0$ then $d_i\, g_m(y,1_m) = d_i\, g_m(s_{i-1} d_i(y),1_m)$.
\MS

\NI (2)\; If $i<n+1$ then $d_i\, g_m(y,1_m) = d_i\, g_m(s_{i} d_i(y),1_m)$.
}
\end{newlem}

\Proof

Proof of (1):

For $p \neq i,i+1$ define $z_p \DFAS f_m(d_p s_{i-1}(y), d_p s_i(1_m))$. Then 
\[
\SETT{z_p:p \neq i,i+1} \in \BOX{i,i+1}{m+1}{D}
\]
because $f$ is simplicial. The Compositions Lemma \refpage{compslemma} implies there is a unique $z \in D_{m+1}$ such that for all $p \neq i,i+1$, $d_p(z) = z_p$.
\MS

Next we verify that $d_i(z)= g_m(y,1_m)$. Both $d_i(z)$ and $g_m(y,1_m)$ are $(m-1,i)$-compositions in $D$ and therefore to show they are equal it suffices to show that for all $p \neq i$, $d_p d_i(z) = d_p\, g_m(y,1_m)$.
\[
\begin{array}{l}
d_p d_i(z) =\left\{
\begin{array}{l}
\begin{aligned}
d_{i-1}(z_p) &= d_{i-1}\,f_m(d_p s_{i-1}(y), d_p s_i(1_m))\\
& = f_{m-1}(d_{i-1}d_p s_{i-1}(y), d_{i-1} d_p s_i(1_m))\\
&= f_{m-1}(d_p(y), d_p(1_m)) \quad \text{ if } p<i
\end{aligned}
\\
\\
\begin{aligned}
d_i(z_{p+1})&= d_i\, f_m(d_{p+1} s_{i-1}(y), d_{p+1} s_i(1_m)) \\
&= f_{m-1}(d_i d_{p+1} s_{i-1}(y), d_i d_{p+1} s_i(1_m)) \\
& = f_{m-1}(d_p(y), d_p(1_m)) \quad \text{ if } p>i
\end{aligned}
\end{array}
\right. 
\end{array}
\]
That is, for all $p \neq i$, 
\[
d_p d_i(z) = f_{m-1}(d_p(y), d_p(1_m)) = d_p\, g_m(y,1_m)
\]
Therefore $d_i(z) = g_m(y,1_m)$, as claimed.
\MS

By similar reasoning, we verify that $d_{i+1}(z)= g_m(s_{i-1} d_i(y),1_m)$. For all $p \neq i$:
\[
\begin{array}{l}
d_p d_{i+1}(z)
= \left\{
\begin{array}{l}
\begin{aligned}
d_i(z_p) &= d_i\, f_m(d_p s_{i-1}(y), d_p s_i(1_m)) \\
&= f_{m-1}(d_p d_{i+1} s_{i-1}(y), d_p(1_m)) \quad \text{ if } p<i\\
\\
\end{aligned}
\\
\begin{aligned}
d_{i+1}(z_{p+1}) &= d_{i+1}\, f_m(d_{p+1} s_{i-1}(y), d_{p+1} s_i(1_m)) \\
&= f_{m-1}(d_p d_{i+1} s_{i-1}(y), d_{p}(1_m)) \quad \text{ if } p>i
\end{aligned}
\end{array}
\right.
\end{array}
\]
That is, for all $p \neq i$,
\[
\begin{aligned}
d_p d_{i+1}(z) &= f_{m-1}(d_p d_{i+1} s_{i-1}(y), d_{p}(1_m)) 
= f_{m-1}(d_p s_{i-1}d_i(y), d_{p}(1_m))\\
& = d_p\, g_m(s_{i-1}d_i(y),1_m)
\end{aligned}
\]
Therefore, $d_{i+1}(z) = g_m(s_{i-1}d_i(y),1_m)$, as claimed.

Finally,
\[
d_i\, g_m(y,1_m) = d_i d_i(z) = d_i d_{i+1}(z) = d_i \, g_m(s_{i-1} d_i(y),1_m)
\]
which completes the proof of (1).
\BS

Proof of (2): We proceed as in the proof of (1). Here, we assume $i<n+1$.
\MS

For $p \neq i,i+1$, define $u_p \DFAS f_m(d_p s_{i+1}(y),d_p s_i(1_m))$. Then it follows that $\SETT{u_+p:p \neq i,i+1} \in \BOX{i,i+1}{m+1}{D}$ because $f$ is simplicial. Again, by the Compositions Lemma, there exists $u \in D_{m+1}$ such that for all $p \neq i,i+1$, $d_p(u) = u_p$.
\MS

\NI Claim: $d_i(u) = g_m(s_i d_i(y),1_m)$.
\MS

\NI Reason: For all $p \neq i$,
\[
\begin{array}{l}
d_p d_i(u) = \left\{
\begin{array}{l}
\begin{aligned}
d_{i-1}(u_p) &= d_{i-1}\, f_m(d_p s_{i+1}(y),d_p s_i(1_m))\\
&= f_{m-1}(d_p d_i s_{i+1}(y),d_p(1_m)) \quad \text{ if } p<i \\
\end{aligned}\\
\\
\begin{aligned}
d_i d_{p+1}(u) &= d_i\, f_m(d_{p+1} s_{i+1}(y),d_{p+1} s_i(1_m)) \\
&= f_{m-1}(d_{p} d_i s_{i+1}(y),d_{p}(1_m)) \quad \text{ if } p>i
\end{aligned}

\end{array}
\right.
\end{array}
\]
That is, for all $p \neq i$,
\[
\begin{aligned}
d_p d_i(u) &= f_{m-1}(d_{p} d_i s_{i+1}(y),d_{p}(1_m)) = 
d_p\, g_m(d_i s_{i+1}(y),1_m)\\
&= d_p\, g_m(s_i d_i(y),1_m)
\end{aligned}
\]
from which it follows that $d_i(u) = g_m(s_i d_i(y),1_m)$.
\MS

\NI Claim: $d_{i+1}(u) = g_m(s_i d_i(y),1_m)$.
\MS

\NI Reason: For all $p \neq i$,
\[
d_p d_{i+1}(u) = \left\{
\begin{array}{l}
\begin{aligned}
d_i(u_p) &= d_i\, f_m(d_p s_{i+1}(y),d_p s_i(1_m)) \\
&= f_{m-1}(d_i d_p s_{i+1}(y),d_i d_p s_i(1_m))\\
&= f_{m-1}(d_p(y),d_p(1_m)) \quad \text{ if } p<i
\end{aligned}
\\
\\
\begin{aligned}
d_{i+1}(u_{p+1}) &= d_{i+1}\, f_m(d_{p+1} s_{i+1}(y),d_{p+1} s_i(1_m)) \\
&= f_{m-1}(d_{i+1} d_{p+1} s_{i+1}(y),d_{i+1} d_{p+1} s_i(1_m))\\
&= f_{m-1}(d_{p}(y),d_{p}(1_m)) \quad \text{ if } p>i
\end{aligned}
 
\end{array}
\right.
\]
Therefore, for all $p \neq i$,
\[
\begin{aligned}
d_p d_{i+1}(u) &= f_{m-1}(d_{p}(y),d_{p}(1_m)) = d_p\, g_m(y,1_m)
\end{aligned}
\]
from which it follows that $d_{i+1}(u) = g_m(y,1_m)$.
\MS

Finally,
\[
d_i\, g_m(y,1_m) = d_i d_{i+1}(u) = d_i d_{i}(u) = d_i\, g_m(s_i d_i(y),1_m)
\]
as claimed.
\qed

\begin{newcor}
\label{basic-lemma-corollary}
{\rm
\SL

Suppose $y,y' \in C_m$ such that $d_i(y) = d_i(y')$. Then 
\[
d_i g_m(y,1_m) = d_i g_m(y',1_m)
\]
}
\end{newcor}

\Proof

If $i>0$ then 
\[
d_i g_m(y,1_m) = d_i g_m(s_{i-1}d_i(y),1_m) = d_i g_m(s_{i-1}d_i(y'),1_m) = d_i g_m(y',1_m)
\]

If $i<n+1$ then, using the simplicial identity $d_{i}s_{i+1} = s_i d_i$,
\[
d_i g_m(y,1_m) = d_i g_m(s_{i}d_i(y),1_m) = d_i g_m(s_{i}d_i(y'),1_m) = d_i g_m(y',1_m)
\]

\qed
\BS

\end{paragr}
\MS

\begin{paragr}
\label{def-g(-,del_i)}
{\bf Definition of $g_{m-1}(-,\del_i)$}
\MS

\begin{newdef}
{\rm
\SL

Given $x \in C_{m-1}$, and writing $d_i(1_m)$ as $\del_i$ we define
\begin{equation}
g_{m-1}(x,\del_i) \DFAS d_i \, g_m(s_i(x),1_m) 
\end{equation}
}
\BX
\end{newdef}

This definition is forced by the need for $g$ to be simplicial. That this is consistent follows from the corollary above: given any $y \in C_m$ such that $d_i(y) = x = d_i s_i(x)$, then $d_i \, g_m(s_i(x),1_m) = d_i \, g_m(y,1_m)$.
\MS

\begin{newcor}
\label{eqn3 for g(y,1)}
{\rm
\SL

For all $y \in C_m$, $d_i\, g_m(y,1_m) = g_{m-1}(d_i(y),d_i(1_m))$.
}
\end{newcor}

\Proof

By definition, $g_{m-1}(d_i(y),d_i(1_m)) = d_i \, g_m(s_i d_i(y),1_m)$. Since $d_i(y) = d_i s_i d_i(y)$, the lemma above implies
\[
g_{m-1}(d_i(y),d_i(1_m))= d_i \, g_m(s_i d_i(y),1_m) = d_i \, g_m(y,1_m)
\]
\qed

\end{paragr}
\BS

\begin{paragr}
\label{def-g(-,sj di (1))}
{\bf Definition of $g_m(-,s_j d_i(1_m))$}
\MS

As noted above, those $m$-simplices of $\Delta[m]$ which do not belong to $\Lambda^i[m]$ are $1_m$ and $s_j d_i(1_m)$ for each $j \in [m-1]$. 

\begin{newlem}
{\rm
\SL

Let $y \in C_m$ and $j \in [m-1]$. Then
\[
\SETT{g_{m-1}(d_p(y), d_p s_j d_i(1_m)) : p \neq i} 
\in \BOX{i}{m}{D}
\]
}
\end{newlem}

\Proof

First, observe that given any $t \in [m-1]$, and $x \in C_{m-1}$ then, by definition of $g_{m-1}(x,d_i(1_m))$,
\[
d_t\, g_{m-1}(x,d_i(1_m)) = d_t d_i\, g_m(s_i(x),1_m)
\]
Now if $t<i$ then
\[
\begin{aligned}
d_t d_i\, g_m(s_i(x),1_m) &= d_{i-1} d_{t}\,g_m(s_i(x),1_m) =
d_{i-1}\,f_{m-1}(d_{t} s_i(x),d_{t}(1_m)) \\
&= f_{m-2}(d_{i-1} d_{t} s_i(x),d_{i-1} d_{t}(1_m)) =
f_{m-2}(d_t(x),d_t d_i(1_m))
\end{aligned}
\]
Similarly, if $t>i$ then
\[
\begin{aligned}
d_t d_i\, g_m(s_i(x),1_m) &= d_{i} d_{t+1}\,g_m(s_i(x),1_m) =
d_{i}\,f_{m-1}(d_{t+1} s_i(x),d_{t+1}(1_m)) \\
&= f_{m-2}(d_{i} d_{t+1} s_i(x),d_{i} d_{t+1}(1_m)) =
f_{m-2}(d_t(x),d_t d_i(1_m))
\end{aligned}
\]
That is, 
\[
\ALL t \neq i, \quad d_t\, g_{m-1}(x,d_i(1_m)) = f_{m-2}(d_t(x), d_t d_i(1_m))
\]

\NI{\bf Notation:} For the remainder of the proof, abbreviate $g_{m-1}(d_p(y), d_p s_j d_i(1_m))$ by $u_p$ ($p \neq i$). That is
\[
u_p= \left\{
\begin{array}{ll}
f_{m-1}(d_p(y), d_p s_j d_i(1_m)) & \text{if } p \neq j,j+1\\
g_{m-1}(d_p(y), d_i(1_m)) & \text{if } p=j \text{ or } j+1
\end{array}
\right.
\]
The goal is to show that for all $p<q$ in $[m]-\{i\}$, $d_p(u_q)=d_{q-1}(u_p)$.
By cases:
\MS

\NI {\bf Case} $\{p,q\} \cap \{j,j+1\} = \MT$: Then $d_p(u_q) = d_{q-1}(u_p)$ because $f$ is simplicial.
\MS

\NI {\bf Case} $p<j$ and $q=j$ or $j+1$:
\[
d_p(u_q) = d_p\, g_{m-1}(d_q(y), d_i(1_m)) = f_{m-2}(d_p d_q(y), d_p d_i(1_m))
\]
equals
\[
\begin{aligned}
d_{q-1}(u_p) &= d_{q-1} \, f_{m-1}(d_p(y), d_p s_j d_i(1_m))=
f_{m-2}(d_{q-1} d_p(y), d_{q-1} d_p s_j d_i(1_m)) \\
&= f_{m-2}( d_p d_{q}(y), d_p d_i(1_m))
\end{aligned}
\]

\NI {\bf Case} $p=j$ and $q=j+1$
\[
\begin{aligned}
d_j(u_{j+1}) &= d_j \,g_{m-1}(d_{j+1}(y), d_i(1_m)) = 
f_{m-2}(d_j d_{j+1}(y), d_j d_i(1_m))\\
&= f_{m-2}(d_j d_{j+1}(y), d_j d_i(1_m)) = d_j(u_j)
\end{aligned}
\]

\NI{\bf Case} $p=j$ or $j+1$ and $q>j+1$:
\[
\begin{aligned}
d_p(u_q) &= d_p\; f_{m-1}(d_q(y), d_q s_j d_i(1_m)) =
f_{m-2}(d_p d_q(y), d_p d_q s_j d_i(1_m)) \\
&= f_{m-2}(d_p d_q(y), d_{q-1} d_i(1_m))
\end{aligned}
\]
(since $d_p s_j = 1$) and this equals
\[
d_{q-1}(u_p) = d_{q-1}\, g_{m-1}(d_p(y), d_i(1_m)) =
f_{m-2}(d_{q-1} d_p(y), d_{q-1} d_i(1_m))
\]
\qed

This lemma together with the requirement that $g$ be simplicial and that $g_m(y,s_j d_i(1_m))$ be an $(m-1,i)$-composition force the following definition.
\MS

\begin{newdef}
{\rm
\SL

Given $y \in C_m$,
\begin{equation}
\label{def-g(y,sj di(1))}
\begin{array}{l}
g_m(y, s_j d_i(1_m)) \DFAS \\ 
= \COMP{m-1}{i}\Bigl(
g_{m-1}(d_0(y),d_0 s_j d_i(1_m)) \DDD{,} \underset{i}{-} \DDD{,} g_{m-1}(d_m(y),d_m s_j d_i(1_m))
\Bigr)
\end{array}
\end{equation}
}
\BX
\end{newdef}

\MS

It remains to show that for each $j \in [m-1]$
\[
d_i \,  g_m(y, s_j d_i(1_m)) \ISIT g_{m-1}(d_i(y),d_i s_j d_i(1_m))
\]
\end{paragr}

\begin{paragr}
\label{tech-lemma}
{\bf Technical lemma}

\begin{newlem}
{\rm
\SL

Suppose $C$ and $D$ are \NICOMP{n}{i}s, $m \geq n+1$, $f : C \times \Lambda^i[m] \to D$ is a simplicial map, $y \in C_m$ and $j \in [m]$.

Then $d_i\, g_m(y,d_i s_j(1_m)) = d_i\, g_m(y,d_{i+1} s_j(1_m))$.
}
\end{newlem}

\Proof

Define $z_p \in D_m$ for $p \neq i,i+1$ as follows.
\[
z_p \DFAS \left\{
\begin{array}{ll}
f_m(d_p s_i(y), d_p s_j(1_m)) & \text{if } p \neq j,j+1\\
& \\
g_m(d_p s_i(y), 1_m) & \text{if } p=j \text{ or } p=j+1 
\end{array}
\right.
\]
{\bf Claim:} $\SETT{z_p:p \neq i,i+1} \in \BOX{i,i+1}{m+1}{D}$.
\MS

\NI{\bf  Verification: }

Suppose $p<q$ in $[m+1]-\{i,i+1\}$. We will verify $d_p(z_q) \ISIT d_{q-1}(z_p)$ by cases.
\MS

\NI{\bf Case} $\{p,q\} \cap \{j,i+1\} = \MT$. Then $d_p(z_q) = d_{q-1}(z_p)$ because $f$ is simplicial.
\MS

\NI{\bf Case} $p<q=j$
\MS

\NI 
$
\begin{array}{l}
d_p(z_j) = d_p\, g_m(d_j s_i(y), 1_m) = f_{m-1}(d_p d_j s_i(y), d_p(1_m)) 
\end{array}
$
\MS

\NI equals
\[
\begin{aligned}
d_{j-1}(z_p) &= d_{j-1}\, f_m(d_p s_i(y), d_p s_j(1_m)) =
f_{m-1}(d_{j-1} d_p s_i(y), d_{j-1} d_p s_j(1_m)) \\
&= f_{m-1}(d_{j-1} d_p s_i(y), d_p(1_m))
\end{aligned} 
\]

\NI{\bf Case} $p<j$ and $q=j+1$
\[
d_p(z_{j+1}) = d_p\,g_m(d_{j+1}s_i(y), 1_m) = f_{m-1}(d_p d_{j+1}s_i(y), d_p(1_m))
\]
equals
\[
\begin{aligned}
d_j(z_p) &= d_j\, f_m(d_p s_i(y), d_p s_j(1_m)) =
f_{m-1}(d_j d_p s_i(y), d_j d_p s_j(1_m)) \\
&= f_{m-1}(d_j d_p s_i(y), d_p(1_m))
\end{aligned}
\]

\NI{\bf Case} $p=j$ and $q=j+1$
\[
d_j(z_{j+1}) = d_j\, g_m(d_{j+1}s_i(y), 1_m) =
f_{m-1}(d_j d_{j+1}s_i(y), d_j(1_m)) = d_j(z_j)
\]

\NI{\bf Case} $p=j$ and $q>j+1$
\[
\begin{aligned}
d_j(z_q) &= d_j\, f_m(d_q s_i(y), d_q s_j(1_m)) =
f_{m-1}(d_j d_q s_i(y), d_j d_q s_j(1_m)) \\
&= f_{m-1}(d_j d_q s_i(y), d_{q-1}(1_m))
\end{aligned}
\]
equals
\[
d_{q-1}(z_j) = d_{q-1}\, g_m(d_j s_i(y), 1_m) = 
f_{m-1}(d_{q-1} d_j s_i(y), d_{q-1}(1_m))
\]

\NI{\bf Case} $p=j+1<q$
\[
\begin{aligned}
d_{j+1}(z_q) &= d_{j+1}\, f_m(d_q s_i(y), d_q s_j(1_m)) =
f_{m-1}(d_{j+1} d_q s_i(y), d_{j+1} d_q s_j(1_m)) \\
&= f_{m-1}(d_{j+1} d_q s_i(y), d_{q-1}(1_m))
\end{aligned}
\]
equals
\[
d_{q-1}(z_{j+1}) = d_{q-1}\, g_m(d_{j+1} s_i(y), 1_m) = 
f_{m-1}(d_{q-1} d_{j+1} s_i(y), d_{q-1}(1_m))
\]
This completes the verification that $\SETT{z_p:p \neq i,i+1} \in \BOX{i,i+1}{m+1}{D}$.

Next, the Compositions Lemma implies there is a uniquely defined $z \in D_{m+1}$ such that for all $p \neq i,i+1$, $d_p(z) = z_p$,
\MS

\NI {\bf Claim:} $d_i(z) = g_m(y,d_i s_j(1_m))$
\MS

\NI{\bf Verification:}

Both $d_i(z)$ and $g_m(y,d_i s_j(1_m))$ are $(m-1,i)$-compositions in $D$. To show they are equal it suffices to show that for all $p \neq i$, $d_p d_i(z) = d_p\,g_m(y,d_i s_j(1_m))$.

When $p<i$:
\[
\begin{array}{l}
d_p d_i(z) = d_{i-1}(z_p)=\\
\\
= \left\{
\begin{array}{l}
d_{i-1}\, f_m(d_p s_i(y), d_p s_j(1_m)) =  f_{m-1}(d_p(y), d_p d_i s_j(1_m)) \text{ if } p \neq j,j+1\\
\\
d_{i-1}\, g_m(d_p s_i(y),1_m) = f_{m-1}(d_p(y),d_{i-1}(1_m)) \text{ if } p=j,j+1
\end{array}
\right.
\end{array}
\]
equals
\[
\begin{array}{l}
d_p\,g_m(y,d_i s_j(1_m)) = g_{m-1}(d_p(y), d_p d_i s_j(1_m)) = \\
\\
= \left\{
\begin{array}{l}
f_{m-1}(d_p(y),d_p d_i s_j(1_m)) \text{ if } p \neq j,j+1\\ 
\\
g_{m-1}(d_p(y), d_{i-1} d_p s_j(1_m)) = f_{m-1}(d_p(y), d_{i-1}(1_m)) \text{ if } p=j,j+1
\end{array}
\right.
\end{array}
\]

When $p>i$:
\[
\begin{array}{l}
d_p d_i(z) = d_i(z_{p+1}) =\\
\\
= \left\{
\begin{array}{l}
d_i\, f_m(d_{p+1} s_i(y), d_{p+1} s_j(1_m)) = f_{m-1}(d_{p}(y), d_i d_{p+1} s_j(1_m)) \text{ if } p+1 \neq j,j+1\\
\\
d_i\, g_m(d_{p+1}s_i(y),1_m) = g_{m-1}(d_{p}(y),d_i(1_m)) \text{ if } p+1 = j,j+1
\end{array}
\right.
\end{array}
\]
equals
\[
\begin{array}{l}
d_p\,g_m(y,d_i s_j(1_m)) = g_{m-1}(d_p(y), d_p d_i s_j(1_m)) = g_{m-1}(d_p(y), d_i d_{p+1} s_j(1_m)) =\\
\\
= \left\{
\begin{array}{l}
f_{m-1}(d_p(y),d_i d_{p+1}s_j(1_m)) \text{ if } p+1 \neq j,j+1\\
\\
g_{m-1}(d_p(y), d_i(1_m)) \text{ if } p+1=j,j+1 
\end{array}
\right.
\end{array}
\]
This completes the proof that $d_i(z) = g_m(y,d_i s_j(1_m))$.
\MS

\NI {\bf Claim:} $d_{i+1}(z) = g_m(y,d_{i+1} s_j(1_m))$
\MS

\NI{\bf Verification:}

When $p<i$:
\[
\begin{array}{l}
d_p d_{i+1}(z) = d_i(z_p) =\\
\\
=\left\{
\begin{array}{l}
d_i\, f_m(d_p s_i(y), d_p s_j(1_m)) = f_{m-1}(d_p(y), d_i d_p s_j(1_m)) \text{ if } p \neq j,j+1\\
\\
d_i\, g_m(d_p s_i(y),1_m) = g_{m-1}(d_p(y),d_i(1_m)) \text{ if } p=j,j+1
\end{array}
\right.
\end{array}
\]
equals
\[
\begin{array}{l}
d_p\, g_m(y,d_{i+1}s_j(1_m)) = g_{m-1}(d_p(y),d_p d_{i+1}s_j(1_m)) =
g_{m-1}(d_p(y),d_i d_{p} s_j(1_m)) =\\
\\
= \left\{
\begin{array}{l}
f_{m-1}(d_p(y),d_i d_{p} s_j(1_m)) \text{ if } p \neq j,j+1\\
\\
g_{m-1}(d_p(y), d_i(1_m)) \text{ if } p=j,j+1
\end{array}
\right.
\end{array}
\]

When $p>i$:
\[
\begin{array}{l}
d_p d_{i+1}(z) = d_i(z_{p+1}) =\\ 
= \left\{
\begin{array}{l}
d_{i+1}\, f_m(d_{p+1} s_i(y), d_{p+1} s_j(1_m)) =  f_{m-1}(d_{p}(y), d_{i+1}d_{p+1} s_j(1_m)) \text{ if } p+1 \neq j,j+1\\
\\
d_{i+1}\, g_m(d_{p+1} s_i(y),1_m) = f_{m-1}(d_p(y), d_{i+1}(1_m)) \text{ if } p+1=j,j+1
\end{array}
\right.
\end{array}
\]
equals

\[
\begin{array}{l}
d_p\, g_m(y,d_{i+1} s_j(1_m)) = g_{m-1}(d_p(y), d_{i+1} d_{p+1} s_j(1_m)) =
\\
\\
= \left\{
\begin{array}{l}
f_{m-1}(d_p(y), d_{i+1} d_{p+1} s_j(1_m)) \text{ if } p+1 \neq j,j+1\\\
\\
f_{m-1}(d_p(y), d_{i+1}(1_m)) \text{ if } p+1 = j,j+1
\end{array}
\right. 
\end{array}
\]
This completes the proof that $d_{i+1}(z) = g_m(y,d_{i+1} s_j(1_m))$.
\MS

\NI To complete the proof of the lemma:
\[
d_i\, g_m(y,d_i s_j(1_m)) = d_i d_i(z) = d_i d_{i+1}(z) = d_i\, g_m(y,d_{i+1} s_j(1_m))
\]

\qed

\end{paragr}
\BS

\begin{paragr}
\label{g is simplicial}
{\bf Proof that $g$ is simplicial}
\MS

\NI{\bf Dimensions below $m-1$:}

$g$ is simplicial in dimensions below $m-1$ because $g_j=f_j$ when $j \leq m-2$.
\MS

\NI{\bf Dimension $m-1$: Verify $d_p g_{m-1} = g_{m-2}d_p$ and
$s_p g_{m-2} = g_{m-1} s_p$}
\MS

\NI For $d_p g_{m-1} = g_{m-2}d_p$:
\MS

Given $(x,\LAM) \in C_{m-1} \times \Delta[m]_{m-1}$, then $g_{m-1}(x,\LAM) = f_{m-1}(x,\LAM)$ unless $\LAM = \del_i = d_i(1_m)$. Therefore we need only verify
\[
d_p g_{m-1}(x,d_i(1_m)) \ISIT g_{m-2}(d_p(x), d_p d_i(1_m))
\]
We will use that $g_{m-1}(x,d_i(1_m)) \DFAS d_i \, g_m(s_i(x),1_m)$.
\MS

Case $p<i$:
\[
\begin{aligned}
d_p g_{m-1}(x,d_i(1_m)) &= d_p d_i \, g_m(s_i(x),1_m) = d_{i-1} d_p \, g_m(s_i(x),1_m)\\
&= d_{i-1}\, f_{m-1}(d_p s_i(x), d_p(1_m)) =
f_{m-2}(d_{i-1} d_p s_i(x), d_{i-1} d_p(1_m))\\
&= f_{m-2}(d_p(x), d_p d_{i}(1_m)) = g_{m-2}(d_p(x), d_p d_{i}(1_m))
\end{aligned}
\]

Case $p \geq i$:
\[
\begin{aligned}
d_p g_{m-1}(x,d_i(1_m)) &= d_p d_i \, g_m(s_i(x),1_m) = d_{i} d_{p+1} \, g_m(s_i(x),1_m)\\
&= d_i\, f_{m-1}(d_{p+1} s_i(x),d_{p+1}(1_m)) =
f_{m-2}(d_i d_{p+1} s_i(x),d_i d_{p+1}(1_m)) \\
&= f_{m-2}(d_{p}(x),d_{p} d_i (1_m)) = g_{m-2}(d_{p}(x),d_{p} d_i (1_m))
\end{aligned}
\]
Therefore $d_p g_m = g_{m-1} d_p$.
\BS

\NI For $s_p g_{m-2} = g_{m-1} s_p$:
\MS

Given $(t,\mu) \in C_{m-2} \times \Delta[m]_{m-2}$ then $\mu \in \Lambda^i[m]_{m-2}$ and $s_p(\mu) \in \Lambda^i[m]_{m-1}$. Therefore $s_p g_{m-2} = g_{m-1} s_p$ holds because $f$ is simplicial. 
\BS

\NI{\bf Dimension $m$:} Verify $d_p g_m = g_{m-1} d_p$
\MS

Let $(y,\gamma) \in C_m \times \Delta[m]_m$. The goal to prove 
\begin{equation}
\label{simp-at-m}
\ALL p \in [m], \quad d_p g_m(y,\gamma) = g_{m-1}(d_p(y),d_p(\gamma))
\end{equation}

If $\gamma \in \Lambda^i[m]_m$ then \eqref{simp-at-m} holds because $f$ is simplicial. We have already seen, above, that $d_p g_m(y,1_m) = g_{m-1}(d_p(y),d_p(1_m))$ in the definition of $g_m(y,1_m)$ (page \pageref{eqn-defining-g(y,1)}) and Corollary \refpage{eqn3 for g(y,1)}.
\MS

It remains to check \eqref{simp-at-m} when $\gamma = s_j d_i(1_m))$, $j \in [m-1]$.
\MS

Equation \eqref{simp-at-m} holds when $p \neq i$, by the definition of $g_m(y,\gamma)$ (page \pageref{def-g(y,sj di(1))}). Therefore, the only remaining step is to verify
\begin{equation}
d_i\, g_m(y,\gamma) = g_{m-1}(d_i(y), d_i(\gamma)), \quad \gamma = s_j d_i(1_m)
\end{equation}
This breaks into cases, as follows.
\MS

\NI {\bf Case} $0 \leq j<i$:

In this case $\gamma = s_j d_i(1_m) = d_{i+1} s_j(1_m)$. Using the technical lemma (indicated by ``$\overset{*}{=}$'') we have
\[
d_i\, g_m(y,\gamma) = d_i\, g_m(y,s_j d_i(1_m)) = d_i\,g_m(y,d_{i+1} s_j(1_m))
\overset{*}{=} d_i\,g_m(y,d_i s_j(1_m))
\]
Now $j<i \IMP d_i s_j(1_m) = s_j d_{i-1}(1_m) \in \Lambda^i[m]_m$. Therefore
\[
d_i\,g_m(y,d_i s_j(1_m)) = d_i\,g_m(y,s_j d_{i-1}(1_m)) = 
f_{m-1}(d_i(y),d_i s_j d_{i-1}(1_m))
\]
By the simplicial identities:
\[
d_i s_j d_{i-1} = d_i d_i s_j = d_i d_{i+1} s_j =d_i s_j d_{i}
\]
we obtain 
\[
d_i\, g_m(y,\gamma)= d_i\, g_m(y,s_j d_i(1_m)) = g_{m-1}(d_i(y), d_i s_j d_i(1_m)) = g_{m-1}(d_i(y), d_i (\gamma))
\]
as claimed.
\MS

\NI{\bf Case} $j=i$:

Here, $\gamma = s_i d_i(1_m) = d_i s_{i+1}(1_m)$. Applying the technical lemma to when $j=i+1$ we get
\[
\begin{aligned}
d_i\,g_m(y,s_i d_i(1_m)) &= d_i\,g_m(y,d_i s_{i+1}(1_m)) \overset{*}{=}
d_i\,g_m(y,d_{i+1} s_{i+1}(1_m))\\
&= d_i\,g_m(y,1_m) = g_{m-1}(d_i(y),d_i(1_m))
\end{aligned}
\]
and note that $d_i(1_m) = d_i s_i d_i(1_m) = d_i(\gamma)$.
\MS

\NI{\bf Case} $j>i$:

Here, $\gamma = s_j d_i(1_m) = d_i s_{j+1}(1_m)$. Applying the technical lemma to $j+1$ we use that $d_{i+1} s_{j+1}(1_m)= s_j d_{i+1}(1_m) \in \Lambda^i[m]_m$ and get
\[
\begin{aligned}
d_i\,g_m(y,\gamma) &= d_i\,g_m(y,d_i s_{j+1}(1_m)) \overset{*}{=} 
d_i\,g_m(y,d_{i+1} s_{j+1}(1_m))\\
&= f_{m-1}(d_i(y),d_i d_{i+1} s_{j+1}(1_m))
\end{aligned}
\]
By simplicial identities, we get
\[
d_i(\gamma) = d_i d_i s_{j+1}(1_m) = d_i d_{i+1} s_{j+1}(1_m)
\]
and therefore $d_i\,g_m(y,\gamma) = f_{m-1}(d_i(y),d_i d_{i+1} s_{j+1}(1_m)) = g_{m-1}(d_i(y),d_i (\gamma))$, as claimed.
\MS

This completes the verification of equation \eqref{simp-at-m}.
\BS

\NI {\bf Degeneracies:} Verify $s_q\, g_{m-1} = g_m s_q$
\MS

Let $(x,\LAM) \in C_{m-1} \times \Delta[m]_{m-1}$. The goal is to check
\begin{equation}
\label{degen-at-m}
\ALL q \in [m-1],\quad  s_q\, g_{m-1}(x,\LAM) = g_m(s_q(x), s_q (\LAM))
\end{equation}

Equation \eqref{degen-at-m} holds when $\LAM \neq \del_i = d_i(1_m))$ because then $\LAM \in \Lambda^i[m]_{m-1}$, $s_q\, g_{m-1}(x,\LAM) = s_q\, f_{m-1}(x,\LAM)$ and $f$ is simplicial. 

It remains to check equation \eqref{degen-at-m} when $\LAM = d_i(1_m)$.
\MS

\NI Now
\[
\begin{array}{l}
d_p s_q\, g_{m-1}(x,d_i(1_m)) = \\ \\=
\left\{
\begin{array}{l}
s_{q-1} d_p \, g_{m-1}(x,d_i(1_m)) = f_{m-1} (s_{q-1}d_p(x),s_{q-1}d_p d_i(1_m)) \text{ if } p<q \\
\\
g_{m-1}(x,d_i(1_m)) \text{ if } p=q,q+1\\
\\
s_q d_{p-1}\, g_{m-1}(x,d_i(1_m)) = f_{m-1}(s_q d_{p-1}(x),s_q d_{p-1}d_i(1_m)) \text{ if } p>q+1
\end{array}
\right.
\end{array}
\]
which, case by case, equals
\[
\begin{array}{l}
d_p \,g_m(s_q(x), d_q d_i(1_m)) = g_{m-1}(d_p s_q(x), d_p s_q d_i(1_m)) = \\
\\ 
= \left\{
\begin{array}{l}
f_{m-1}(s_{q-1} d_p(x), s_{q-1} d_p d_i(1_m)) \text{ if } p<q \\
\\
g_{m-1}(x,d_i(1_m)) \text{ if } p=q,q+1\\
\\
f_{m-1}(s_q d_{p-1}(x),s_q d_{p-1} d_i(1_m)) \text{ if } p>q+1
\end{array}
\right.
 
\end{array}
\]
This completes the verification of equation \eqref{degen-at-m} and the proof of Theorem \refpage{function complex theorem}.

\qed

\end{paragr}

\section{Comma-composers}
\label{comma-composers}

Given a small category $C$ and $x, y \in \OB{C}$,  then the hom-set $C(x,y)$ is implicated in several basic constructs: the hom functors $C(x,-)$ and $C(-,y)$, and the comma categories $C \downarrow y$ and $x \downarrow C$.

Statements concerning these ideas can be phrased in terms of the nerve of $C$ as statements about $0$- and $1$-simplices, and the associated nerves for $C \downarrow y$ and $x \downarrow C$.
\MS

Now if $C$ is, for the moment, an arbitrary simplicial set then similar statements can be made concerning simplices $x$ and $y$ of arbitrary dimension. If $x \in C_k$ and $y \in C_j$, we may consider any and all $z \in C_{j+k+1}$ which contain $x$ and $y$ as subfaces of $z$ {\em and} such that the vertices of $x$ and $y$, as sets of vertices of $z$, are disjoint. In a phrase, we will say ``$x$ and $y$ are {\em complementary subfaces of $z$}''. In this setting, the analog of the hom set $C(x,y)$ is the set, possibly empty, of all $z \in C_{j+k+1}$ in which $x$ and $y$ are complementary. The analog of the comma category $x \downarrow C$ is the set
\[
\bigcup_{j \geq 0} \left\{ z \in C_{j+k+1} : x \text{ is a subface of } z \right\}
\]
which we will denote by $C^x$.

We will verify the unsurprising fact that $C^x$ is a simplicial set analogous to the nerve of $x \downarrow C$ whose face and degeneracy maps are those of $C$ restricted to those not involving the subface $x$. (Theorem \refpage{C^x-is-simplicial-set}).
\MS

If $C$ is \NICOMP{n}{i} then what about $C^x$? We will show that it too is an \NICOMP{n}{i}.
\MS

This section develops all of these ideas in detail.

\subsection{Definitions and notation for complementary subfaces}

Let $C$ be an arbitrary simplicial set.

\begin{itemize}
\item 
\index{vertex index}
\index{$\VERT{t}(x)$}
Given any $k>0$, $x \in C_k$ and any $t \in [k]$ then vertex $t$ of $x$ is $d_0 \DDD{}\omit{d_t} \DDD{} d_k(x)$ and $t$ is the {\bf vertex index} of that vertex. We will use the notation $\VERT{t}(x) \DFAS \text{vertex } t \text{ of } x$.

\item
Given any $m>0$,  $z \in C_m$ and any proper non-empty subset \[
A = \SETT{u_0 \DDD{<} u_j} \subsetneq [m]
\]
let $x = d_{u_0} \DDD{} d_{u_j}(z)$.
\index{vertex index list}
Then
\[
\text{\bf vertex index list of } x \text{ \bf in } z \DFAS
\VL_z(x) \DFAS [m]-A
\]

\item
\index{sharp notation $\sharp$}
\index{flat notation $\flat$}
\index{$\sharp$}
\index{$\flat$}
\label{sharp-flat-notation}
{\bf Sharp ($\sharp$) and flat ($\flat$) notation:}

{\bf Sharp:} If $\DMAP{f}{k}{m}$ is a non-decreasing function then define the \SI\ function $\DMAP{\SHR{f}}{k}{m+k}$ by $\SHR{f}(t) \DFAS f(t)+t$.

That is, $f\SH$ is the pointwise sum of $f$ and $1_{[k]}$.

{\bf Flat:} If $\DMAP{g}{k}{m}$ is a \SI\ function, define the non-decreasing function $\DMAP{\FLT{g}}{k}{m-k}$ by $\FLT{g}(t) \DFAS g(t)-t$.

\item
\index{vertex index function}
Given any $m>0$, $z \in C_m$, and a proper subface $x \in C_k$ of $z$ with vertex index list $B = \SETT{w_0 \DDD{<} w_k} \SBS [m]$, then the {\bf vertex index function for $x$ in $z$} is the non-decreasing function $\mu_B : [k] \to [m-k]$ defined by $\mu_B(q) \DFAS w_q-q$ and characterized by the  condition that $\SHR{\mu}_B[k] = B$.

\item
{\bf Integer-interval notation:}

Given integers $p$ and $q$:
\[
[p,q] \DFAS \left\{
\begin{array}{ll}
\SETT{p ,p+1\DDD{,} q} & \text{if } p<q\\
\SETT{p} & \text{if } p=q\\
\emptyset & \text{if }p>q
\end{array}
\right.
\]
\BX
\end{itemize}

\begin{newdef}
{\bf (Complementary subfaces)}
\index{complementary subfaces}
{\rm
\SL

Let $C$ be a simplicial set, $m>0$ and $z \in C_m$. Let $x$ and $y$ be proper subfaces of $z$. Then $x$ and $y$ will be said to be {\bf complementary in $z$} if $\VL_z(x) = [m] -\VL_z(y)$.
} 
\BX
\end{newdef}

Note that given $m>0$ and $z \in C_m$ then any non-empty $A \subsetneq [m]$ determines exactly one pair of complementary subfaces of $z$.

To paraphrase the definition: if $x \in C_k$ and $y \in C_j$ are complementary in $z \in C_m$ and 
\[
\begin{aligned}
\VL_z(y) &= \SETT{u_0 \DDD{<} u_j}=A\\
\VL_z(x) &= \SETT{w_0 \DDD{<} w_k}=B
\end{aligned}
\]
then
\begin{enumerate}
\item 
$A$ and $B$ are disjoint and $A \cup B = [m] = [j+k+1]$.

\item
$d_{u_0} \cdots d_{u_j}(z)=d_A(z)=x$ and $d_{w_0} \cdots d_{w_k}(z)=d_B(z)=y$.

\item $\VERT{q}(x)$ is $\VERT{w_q}(z)$ and $\VERT{p}(y)$ is $\VERT{u_p}(z)$.

\end{enumerate}

\begin{newlem}
\label{subface-of-subface}
{\rm
\SL

Suppose $C$ is a simplicial set in which $u \in C_r$ is a subface of $x \in C_k$ and $x$ is a subface of $z \in C_m$ where
\[
\begin{aligned}
\VL_z(x)= \SETT{\SEQ{b}{q}} \SBS [m], & \quad [k] \XRA{\beta} [m], \quad \beta(q) \DFAS b_q\\
\VL_x(u) = \SETT{\SEQ{a}{r}} \SBS [k], & \quad [r] \XRA{\alpha}[k], \quad \alpha(t) \DFAS a_t\\
\VL_z(u) = \SETT{\SEQ{g}{r}} \SBS [k], & \quad [r] \XRA{\gamma} [m], \quad \gamma(t) \DFAS g_t
\end{aligned}
\]
Then $\gamma = \beta \circ \alpha$.
}
\end{newlem}
\Proof

Vertex $t$ of $u$ is vertex $a_t = \alpha(t)$ of $x$. Vertex $q$ of $x$ is vertex $\beta(q)$ of $z$. That is, vertex $t$ of $u$ is vertex $\beta(\alpha(t))$ of $z$.

\qed
\MS

\NI {\bf Remark:} Suppose $x \in C_k$ and $y \in C_j$ are complementary subfaces of $z \in C_m$ ($m=j+k+1$) where
\[
\begin{array}{l}
\VL_z(x) = \SETT{b_0 \DDD{<} b_k} \text{ where } [k] \XRA{\beta} [m], \quad \beta(q) \DFAS b_q\\
\VL_z(y) = \SETT{u_0 \DDD{<} u_j} \text{ where } [j] \XRA{\alpha} [m], \quad \alpha(p) \DFAS u_p
\end{array}
\]
Then in the simplicial set $\Delta[m]$, $\alpha \in \Delta[m]_j$ and $\beta \in \Delta[m]_k$ are complementary subfaces of $1_{[m]} \in \Delta[m]_m$.

\BX

\subsection{Simplicial encodings for complementary subfaces}

A pair of complementary subfaces of an $m$-simplex is described by a pair of vertex index lists which partitions $[m]$. In this section we examine that partition in terms of ``complementary vertex functions'', defined below.
\MS

The following notations will be convenient.

\begin{itemize}

\item 
Given any non-decreasing $\alpha : [m] \to [n]$ then we will denote $\alpha$ by its list of values $(a_0 \DDD{,} a_m)$, with $n$ understood from the context, where $a_t \DFAS \alpha(t)$. Then
\[
{\alpha}^\sharp = (a_0,a_1+1 \DDD{,} a_m+m)
\]

Given $\DMAP{\alpha}{m}{n}$, $p \in [m]$ and $r \in [n]$, then:
\begin{eqnarray*}
\alpha \del_p = d_p(\alpha) &=&\left\{
\begin{array}{ll}
(a_1 \DDD{,} a_m) & \text{if } p=0\\
(a_0 \DDD{,} a_{p-1},a_{p+1} \DDD{,} a_m) &\text{if } 0<p<m\\
(a_0 \DDD{,} a_{m-1}) & \text{if } p=m
\end{array}
\right. \\
\alpha \SIG_p = s_p(\alpha) &=& (a_0 \DDD{,} a_p,a_p \DDD{,} a_m)
\\
\del_r \alpha &=& \left\{
\begin{array}{ll}
(a_0+1 \DDD{,} a_m+1) & \text{if } r \leq a_0\\ 
(a_0 \DDD{,} a_{p-1},a_p+1 \DDD{,} a_m+1) & \text{if } a_{p-1}<r \leq a_p\\
(a_0 \DDD{,} a_m) & \text{if } r>a_m
\end{array}
\right. \\
\SIG_r \alpha &=& \left\{
\begin{array}{ll}
(a_0-1 \DDD{,} a_m-1) & \text{if } r<a_0\\
(a_0 \DDD{,} a_{p-1},a_p-1 \DDD{,} a_m-1) & \text{if } a_{p-1} \leq r<a_p\\
\LIST{a}{m} & \text{if } r \geq a_m 
\end{array}
\right.
\end{eqnarray*}

\BX
\end{itemize}

In the next definition and lemma we focus on describing a partition $[m] = A \cup B$ in terms of certain related non-decreasing functions.

\begin{newdef}
{\bf (Complementary vertex vunctions)}
\index{complementary vertex functions}
\label{comp-vert-functions}
{\rm
\SL

Given $j,k \geq 0$ and non-decreasing functions $\mu : [k] \to [j+1]$ and $\LAM : [j] \to [k+1]$, we say that $\mu$ and $\LAM$ are {\bf complementary vertex functions} if the sets
\[
A = \SETT{\LAM(p)+p: p \in [j]} \quad \text{and} \quad
B = \SETT{\mu(q)+q:q \in [k]}
\]
are disjoint and $A \cup B = [j+k+1]$. That is, $\SHR{\LAM}[j]=A$ and $\SHR{\mu}[k]=B$.

\BX
}
\end{newdef}
\NI Note that: 

(1) Any non-decreasing function $[j] \to [k+1]$ determines a unique non-decreasing function $[k] \to [j+1]$ which makes the pair complementary.

(2) The roles of $\LAM$ and $\mu$ are interchangeable. Therefore any true statement concerning $\mu$, $\LAM$, $k$ and $j$ remains true if the roles of $\mu$ and $k$ are interchanged with those of $\LAM$ and $j$.

(3) With $\LAM$ and $\mu$ as above (see page \pageref{del_A} for notation): 
\[
\begin{aligned}
\mu\SH &= \DMAP{\del_{\LAM\SH[j]}}{k}{j+k+1}\\
\LAM\SH &= \DMAP{\del_{\mu\SH[k]}}{j}{j+k+1}
\end{aligned}
\]
\MS

A pair of complementary vertex functions as in the definition above relate in several ways, as described in the following lemma. In particular, $\LAM$ has a certain minimal property with respect to $\mu$ (as does $\mu$, by symmetry,  with respect to $\LAM$).

\begin{newlem} 
\label{lambda-least}
{\rm
\SL

Suppose $j,k \geq 0$ and $\LAM:[j] \to [k+1]$ and $\mu : [k] \to [j+1]$ are complementary vertex functions. For each $p \in [j]$ and $q \in [k]$, denote $\SHR{\LAM}(p)$ by $u_p$, $\SHR{\mu}(q)$ by $w_q$, and let
\[
A \DFAS \SETT{u_p:p \in [j]}=\LAM\SH[j] \quad \text{and} \quad B \DFAS \SETT{w_q:q \in [k]} = \mu\SH[k]
\]
Then:
\begin{enumerate} 
\item 
 \begin{enumerate}
   \item 
   $\ALL p \in [j],\; \LAM(p) = \left| B \cap [0,u_p-1] \right|$
   
   \item
   $\ALL q \in [k],\; \mu(q) = \left| A \cap [0,w_q-1] \right|$
 \end{enumerate}
 
 \item 
 \label{al=lam}
  \begin{enumerate}
    \item 
    Define $\alpha:[j] \to [k+1]$ by 
    \[
    \alpha(p) = \left\{
    \begin{array}{ll}
    \text{least } q \text{ such that } u_p<w_q & \text{if } u_p<w_k\\
    k+1 & \text{if } u_p>w_k
 	\end{array}
    \right.
    \]
    Then $\alpha = \LAM$.
    
    \item
    Define $\beta:[k] \to [j+1]$ by 
    \[
    \beta(q) = \left\{
    \begin{array}{ll}
    \text{least } p \text{ such that } w_q<u_p & \text{if } w_q<u_j\\
    j+1 & \text{if } w_q>u_j
 	\end{array}
    \right.
    \]
    Then $\beta = \mu$.
  \end{enumerate}
 \item 
 \label{part3}
   \begin{enumerate}
      \item Given $p \in [j]$ and $q \in [k]$ then the following are equivalent:
      \begin{enumerate}
         \item $q$ is least such that $u_p<w_q$
         \item $q = \LAM(p)<k+1$
         \item $q$ is least such that $p<\mu(q)$
      \end{enumerate}
      
      \item Given $p \in [j]$ then the following are equivalent:
      \begin{enumerate}
         \item $u_p>w_k$
         \item $\LAM(p)= k+1$
         \item $p+1 > \mu(k)$
      \end{enumerate}

   \end{enumerate}

\item 
For all $p \in [j]$ and all $q \in [k]$, $p < \mu(q) \iff \LAM(p) \leq q$.

 \end{enumerate} 

}
\end{newlem}

\Proof

\begin{enumerate}
\item 
Given any $q \in [k]$ then, trivially, $w_q = \bigl| [0,w_q-1] \bigr|$ and 
\[
B \cap [0,w_q-1] = \left\{
\begin{array}{ll}
\SETT{\SEQ{w}{{q-1}}} & \text{if } q>0\\
\MT & \text{if } q=0
\end{array}
\right.
\]
Using that $A \cup B = [j+k+1]$ and $A \cap B = \MT$ then
\[
w_q = \Bigl|[0,w_q-1] \Bigr| = \Bigl|A \cap [0,w_q-1] \Bigr| + \Bigl| B \cap [0,w_q-1] \Bigr| = \Bigl|A \cap [0,w_q-1] \Bigr| + q
\]
Therefore $\mu(q) = w_q-q = \Bigl| A \cap [0,w_q-1] \Bigr|$, as claimed.

By interchanging $\mu$ with $\LAM$, the same argument proves that $\LAM(p) = u_p-p  = \Bigl| B \cap [0,u_p-1] \Bigr|$.

\item 
We will prove that $\alpha = \LAM$. The proof that $\beta = \mu$ follows by interchange.

Let $p \in [j]$. There are two cases for evaluating $\alpha(p)$: $u_p<w_k$ and $u_p>w_k$.

\NI {\bf Case $u_p<w_k$:} Let $q \in [k]$ be the least such that $u_p<w_q$. Then $\alpha(p)=q$ by definition, and
\[
B \cap [0,u_p-1] = \left\{
\begin{array}{ll}
\SETT{\SEQ{w}{{q-1}}} & \text{if } q>0 \\
\MT & \text{if } q=0
\end{array}
\right.
\]
Then, using reasoning similar to that in the previous item
\[
\begin{aligned}
u_p = \LAM(p)+p = \Bigl| [0,u_p-1] \Bigr| &= \Bigl| A \cap [0,u_p-1] \Bigr| + \Bigl| B \cap [0,u_p-1] \Bigr| \\
&= \Bigl| \SETT{ \SEQ{u}{{p-1}}} \Bigr| + q = p+q
\end{aligned}
\]
Therefore $\alpha(p) = q = u_p-p=\LAM(p)$ in this case.

\NI {\bf Case $u_p>w_k$:} By definition $\alpha(p)=k+1$ here, $B \cap [0,u_p-1] = B$ and $A \cap [0,u_p-1] = \SETT{\SEQ{u}{{p-1}}}$. Thus
\[
u_p = \LAM(p)+p = \Bigl| [0,u_p-1] \Bigr| = p + (k+1)
\]
and therefore $\LAM(p)=k+1=\alpha(p)$.

\item 
The result of part \ref{al=lam} is that for each $p \in [j]$\begin{itemize}
\item 
$\LAM(p)< k+1 \iff u_p<w_{\LAM(p)}$ with $\LAM(p)$ the least such.

\item
$\LAM(p)=k+1 \iff u_p>w_k$.
\end{itemize}

We will use this observation in the proofs of this part of the lemma.
\begin{enumerate}

\item 
Assume $p \in [j], q \in [k]$. (i.) $\iff$ (ii.) by part \ref{al=lam}.

(ii.) $\IMP$ (iii.): If $q = \LAM(p)<k+1$ then $u_p<w_q$. Rewriting this in terms of $\mu$ we have
$
\LAM(p)+p< \mu(\LAM(p))+\LAM(p)
$
which says $p<\mu(q)$ with $q$ the least such.

(iii.) $\IMP$ (i.): If $q$ is least such that $p<\mu(q)$ then, by adding $q$ we get $p+q<w_q$ with $q$ being least. There are two cases: $q>0$ and $q=0$.

If $q>0$ then $w_{q-1}<p+q<w_q$ which implies that $p+q = u_{p'}$ for some $p' \in [j]$. By part  \ref{al=lam}, $\LAM(p') = q$. Since $p+q = u_{p'} = \LAM(p')+p'$ then $p=p'$ and $w_{q-1}<u_p<w_q$. 

If $q=0$ then $p<\mu(0) = w_0$ implies $p=u_{p'}$ for some $p'$ where, by part \ref{al=lam}, $\LAM(p')=0$. Since $p=u_{p'} = \LAM(p')+p'$ then $p=p'$ and $u_p = u_{p'}<w_0$.

\item 
Assume $p \in [j]$. (i.) $\iff$ (ii.) by part \ref{al=lam}.

(ii.) $\IMP$ (iii.): If $\LAM(p)=k+1$ then $u_p>w_k$. In terms of $\mu$ this says $u_p=k+1+p>\mu(k)+k$ and therefore $p+1>\mu(k)$.

(iii.) $\IMP$ (i.): If $p+1>\mu(k)$ then $p+1+k>w_k$ and therefore $p+1+k = u_{p'}$ for some $p'$. By part \ref{al=lam} $\LAM(p')=k+1$. Since $p+1+k = u_{p'} = \LAM(p')+p'$ then $p=p'$ and therefore $u_p = u_{p'} > w_k$.

\end{enumerate}

\item 
Assume $p<\mu(q)$ where $q \in [k]$. Then $\LAM(p)<k+1$ because otherwise $\LAM(p)=k+1$ would imply $\mu(k) \leq p$ (part \ref{part3}) and $p<\mu(q) \leq \mu(k) \leq p$, which is impossible. By part \ref{part3}, $\LAM(p) \in [k]$ is least such that $p<\mu(\LAM(p))$. Therefore $p<\mu(q)$ implies $\LAM(p) \leq q$.

Conversely, if $\LAM(p) \leq q \in [k]$ then $\LAM(p)<k+1$ and (part \ref{part3} again) $\LAM(p)$ is the least element of $[k]$ such that $p<\mu(\LAM(p))$, hence $\mu(\LAM(p)) \leq \mu(q)$. Therefore $p< \mu(q)$.
\end{enumerate}
\qed

\begin{newlem} 
\label{cor-lambda-least}
{\rm
\SL

Suppose $j,k \geq 0$ and that $\LAM:[j] \to [k+1]$ and $\mu:[k] \to [j+1]$ are complementary vertex functions. Then for each $p \in [j]$
\[
\del_{\LAM\SH(p)} \mu^\sharp = \Bigl( \del_{p+1} \mu \Bigr)^\sharp \quad \text{and} \quad \SIG_{\LAM\SH(p)} \mu^\sharp = \Bigl( \SIG_p \mu \Bigr)^\sharp
\]
}
\end{newlem}

\Proof

Notation:  we'll denote $\mu(q)$ by $m_q$, $\mu^\sharp(q) \DFAS m_q+q$ by $w_q$ and $\LAM\SH(p)=\LAM(p)+p$ by $u_p$.
\MS

As a preliminary observation: 
\[
m_{q-1} \leq p<m_q \iff m_{q-1}-1<p<m_q \iff m_{q-1}<p+1 \leq m_q
\]
Now
\[
m_{q-1}-1<p<m_q \iff w_{q-1} < p+q < w_q
\]
and lemma \refpage{lambda-least} implies $p+q = u_{p'} = \LAM(p')+p'=p'+q$. That is $p+q=u_p$ and $w_{q-1} < u_p < w_q$.

Conversely: if $w_{q-1}<u_p<w_q$ then $q=\LAM(p)$ which says  $m_{q-1}-1<p<m_q$. 
\MS

In summary, the inequalities $m_{q-1}<p+1 \leq m_q$, $m_{q-1} \leq p<m_q$ and $w_{q-1}<u_p<w_q$ are equivalent.
\MS

Now to the assertions of the lemma.
\[
\del_{u_p} \mu^\sharp = \left\{
\begin{array}{ll}
(w_0+1 \DDD{,} w_k+1) & \text{if } u_p \leq w_0\\
(w_0 \DDD{,} w_{q-1},w_q+1 \DDD{,} w_k+1) & \text{if } w_{q-1}<u_p \leq w_q
\\
\LIST{w}{k} & \text{if } u_p>w_k 
\end{array}
\right.
\]
where $u_p \leq w_0 \IMP u_p<w_0$ and similarly, $u_p \leq w_q \IMP u_p<w_q$.
\[
\del_{p+1} \mu = \left\{
\begin{array}{ll}
(m_0+1 \DDD{,} m_k+1) & \text{if } p+1 \leq m_0\\
(m_0 \DDD{,} m_{q-1},m_q+1 \DDD{,} m_k+1) & \text{if } m_{q-1}<p+1 \leq m_q \\
\LIST{m}{k} & \text{if } p+1>m_k
\end{array}
\right.
\]
where, by lemma \refpage{lambda-least} and the preliminary remarks above, $p+1 \leq m_0 \iff u_p<w_0$, $m_{q-1}<p+1 \leq m_q \iff w_{q-1}<u_p<w_q$ and $p+1>m_k \iff u_p>w_k$.
\[
\Bigl(\del_{p+1} \mu \Bigr)^\sharp = \left\{
\begin{array}{ll}
(w_0+1 \DDD{,} w_k+1) & \text{if } p+1 \leq m_0\\ 
(w_0 \DDD{,} w_{q-1},w_q+1 \DDD{,} w_k+1) & \text{if } w_{q-1}<u_p < w_q \\
\LIST{w}{k} & \text{if } u_p>w_k
\end{array}
\right.
\]
This proves $\del_{u_p} \mu^\sharp = \Bigl(\del_{p+1} \mu \Bigr)^\sharp$.
\BS

By similar calculation:
\[
\SIG_{u_p} \mu^\sharp = \left\{
\begin{array}{ll}
(w_0-1 \DDD{,} w_k-1) & \text{if } u_p < w_0\\ 
(w_0 \DDD{,} w_{q-1},w_q-1 \DDD{,} w_k-1) & \text{if } w_{q-1} \leq u_p < w_q
\\
\LIST{w}{k} & \text{if } u_p \geq w_k
\end{array}
\right.
\]
where ``$w_{q-1} \leq u_p$'' and $``u_p \geq w_k$'' are necessarily strict inequalities.
\[
\SIG_p \mu = \left\{
\begin{array}{ll}
(m_0-1 \DDD{,} m_k-1) & \text{if } p<m_0\\
(m_0 \DDD{,} m_{q-1},m_q-1 \DDD{,} m_k-1) & \text{if } m_{q-1} \leq p< m_q \\
\LIST{m}{k} & \text{if } p \geq m_k 
\end{array}
\right.
\]
By lemma \refpage{lambda-least} and the remarks above, $p<m_0 \iff u_p<w_0$, $m_{q-1} \leq p< m_q \iff w_{q-1}<u_p<w_q$ and $p \geq m_k \iff u_p>w_k$.
Therefore:
\[
\Bigl(\SIG_p \mu \Bigr)^\sharp = \left\{
\begin{array}{ll}
(w_0-1 \DDD{,} w_k-1) & \text{if } u_p<w_0\\
(w_0 \DDD{,} w_{q-1},w_q-1 \DDD{,} w_k-1) & \text{if } w_{q-1} < u_p< w_q \\
\LIST{w}{k} & \text{if } u_p > w_k 
\end{array}
\right.
\]
This proves $\SIG_{u_p} \mu^\sharp = \Bigl(\SIG_p \mu \Bigr)^\sharp$.

\qed
\MS

\NI {\bf Comment:} The interchange version of this lemma is that for all $q \in [k]$, $\del_{w_q} \LAM^\sharp = \Bigl( \del_{q+1} \LAM \Bigr)^\sharp$ and $\SIG_{w_q} \LAM^\sharp = \Bigl( \SIG_q \LAM \Bigr)^\sharp$.

\BX
\BS

The gist of the next theorem is that when $x$ and $y$ are complementary subfaces of $z$, then:

(a) For each face of $z$, either $x$ or $y$ is a subface of that face, and there are formulas for the vertex functions derived from those for $x$ and $y$ in $z$.

(b) Each degeneracy of $z$ has complementary subfaces consisting either of $x$ and a degeneracy of $y$, or $y$ and a degeneracy of $x$. Again, there are formulas for the vertex functions in terms of those for $x$ and $y$.

\begin{newthm}
{\bf (Complementary Subface Theorem)}
\label{complementary-face-idents}
{\rm 
\SL

Suppose $x \in C_k$ and $y \in C_j$ are complementary subfaces of $z \in C_{j+k+1}$ with $\mu : [k] \to [j+1]$ and $\LAM:[j] \to [k+1]$ the vertex functions for $x$ and $y$ respectively. As above, denote $\LAM^\sharp(p)$ by $u_p$ and $\mu^\sharp(q)$ by $w_q$.

Then
\begin{enumerate}
\item 
For each $p \in [j]$:
\begin{enumerate}
\item 
$x$ and $d_p(y)$ are complementary in $d_{u_p}(z)$.

\item
The vertex function for $d_p(y)$ in $d_{u_p}(z)$ is $\LAM \del_p$.

\item
The vertex function for $x$ in $d_{u_p}(z)$ is $\SIG_p \mu$.

\end{enumerate}

\item 
For each $p \in [j]$:

\begin{enumerate}
\item 
$x$ and $s_p(y)$ are complementary in $s_{u_p}(z)$.

\item
The vertex function for $s_p(y)$ in $s_{u_p}(z)$ is $\LAM \SIG_p$.

\item
The vertex function for $x$ in $s_{u_p}(z)$ is $\del_{p+1} \mu$.

\end{enumerate}

\end{enumerate}

}
\end{newthm}

\Proof

\begin{enumerate}
\item 

\begin{enumerate}
\item 
$d_{u_p}(z)$ deletes vertex $u_p$ of $z$. Since that is vertex $p$ of $y$ then $x$ and $d_p(y)$ are complementary in $d_{u_p}(z)$.

\item
Since $\VL_z(y) = \SETT{u_0 \DDD{,} u_j}$ then, by the simplicial identities, 
\[
\VL_{d_{u_p}(z)}(d_p(y)) = \left\{
\begin{array}{ll}
\SETT{u_1-1 \DDD{,} u_j-1} & \text{if } p=0\\ 
\SETT{u_0 \DDD{,} u_{p-1}, u_{p+1}-1 \DDD{,} u_j-1} & \text{if } p \in [1,j-1]\\
\SETT{u_0 \DDD{,} u_{j-1}} & \text{if } p=j
\end{array}
\right.
\]
Now $\LAM = (u_0,u_1-1,u_2-2 \DDD{,} u_j-j)$. Therefore
\[
\LAM \del_p = 
\left\{
\begin{array}{ll}
(u_1-1,u_2-2 \DDD{,} u_j-j) & \text{if } p=0\\
&\\
(u_0,u_1-1 \DDD{,} u_{p-1}-(p-1), & \\
\qquad u_{p+1}-(p+1) \DDD{,} u_j-j)   & \\
& \text{if } p \in [1,j-1]\\
&\\
(u_0,u_1-1 \DDD{,} u_{j-1}-(j-1)) & \text{if } p=j
\end{array}
\right.
\]
and thus
\[
{(\LAM \del_p)}^\sharp = 
\left\{
\begin{array}{ll}
(u_1-1,u_2-1 \DDD{,} u_j-1) & \text{if } p=0\\
(u_0, u_1 \DDD{,} u_{p-1}, u_{p+1}-1 \DDD{,} u_j-1) & \text{if } p \in [1,j-1]\\
(u_0, u_1 \DDD{,} u_{j-1}) & \text{if } p=j
\end{array}
\right.
\]
This establishes the claim.

\item
By the simplicial identities:
\[
\VL_{d_{u_p}(z)}(x) = 
\left\{
\begin{array}{ll}
\SETT{w_0-1 \DDD{,} w_k-1} & \text{if } u_p<w_0 \\
\SETT{w_0 \DDD{,} w_{q-1},w_q-1 \DDD{,} w_k-1} & \text{if } w_{q-1}<u_p<w_q\\
\SETT{w_0 \DDD{,} w_k} & \text{if } u_p > w_k
\end{array}
\right.
\]
Now
\[
\SIG_{u_p} \mu^\sharp = 
\left\{
\begin{array}{ll}
(w_0-1 \DDD{,} w_k-1) & \text{if } u_p < w_0\\ 
(w_0 \DDD{,} w_{q-1},w_q-1 \DDD{,} w_k-1) & \text{if } w_{q-1} \leq u_p < w_q
\\
\LIST{w}{k} & \text{if } u_p \geq w_k
\end{array}
\right.
\]
By lemma \refpage{cor-lambda-least}, $\SIG_{u_p} \mu^\sharp = \Bigl( \SIG_p \mu \Bigr)^\sharp$ and therefore $\SIG_p \mu$ is the vertex function for $x$ in $d_{u_p}(z)$.

\end{enumerate}

\item 

\begin{enumerate}
\item 
$s_{u_p}(z)$ duplicates vertex $u_p$ of $z$ which is vertex $p$ of $y$. Thus $x$ and $s_p(y)$ are complementary in $s_{u_p}(z)$.

\item
Since $\VL_z(y) = \SETT{u_0 \DDD{,} u_j}$ the simplicial identities imply
\[
\VL{s_{u_p}(z)}(s_p(y)) = (u_0 \DDD{,} u_{p-1}, u_p, u_p+1 \DDD{,} u_j+1)
\]
Then
\[
\begin{aligned}
\LAM &= (u_0, u_1-1 \DDD{,} u_j-j)\\
\LAM \SIG_p &= (u_0,u_1-1 \DDD{,} u_p-p, u_p-p \DDD{,} u_j-j)\\
(\LAM \SIG_p)^\sharp &= (u_0,u_1 \DDD{,} u_p, u_p+1 \DDD{,} u_j+1)
\end{aligned}
\]
establishes the claim.

\item
By the simplicial identities 
\[
\VL_{s_{u_p}(z)}(x) =
\left\{
\begin{array}{ll}
\SETT{w_0+1 \DDD{,} w_{k}+1} & \text{if } u_p<w_0  \\
\SETT{w_0 \DDD{,} w_{q-1},w_q+1 \DDD{,} w_k+1} & \text{if } w_{q-1}<u_p<w_q\\
\SETT{w_0 \DDD{,} w_k} & \text{if } u_p>w_k
\end{array}
\right.
\]
Now
\[
\del_{u_p} \mu^\sharp = \left\{
\begin{array}{ll}
(w_0+1 \DDD{,} w_k+1) & \text{if } u_p < w_0\\
(w_0 \DDD{,} w_{q-1},w_q+1 \DDD{,} w_k+1) & \text{if } w_{q-1}<u_p < w_q
\\
\LIST{w}{k} & \text{if } u_p>w_k 
\end{array}
\right.
\]
By lemma \refpage{cor-lambda-least}, $\del_{u_p} \mu^\sharp = \Bigl(\del_{p+1} \mu \Bigr)^\sharp$. Therefore $\del_{p+1} \mu$ is the vertex function for $x$ in $s_{u_p}(z)$.

\end{enumerate}

\end{enumerate}
\qed

\subsection{Complementary complexes $C^x$}

In this subsection, $C$ denotes an arbitrary simplicial set.
\MS

\HD{Notations:}

\index{$C^{x,y}_{k,j}$}
\index{C$^{x,y}_{k,j}$}
\index{$C^x_j$}
\index{C$^x_j$}
\index{$C^x$}
\index{C$^x$}

Let $j,k \geq 0$ and suppose $x \in C_k$, $y \in C_j$. Then
\[
\begin{aligned}
C^{x,y}_{k,j} & \DFAS \SETT{z \in C_{j+k+1}: x,y \text{ are complementary in  } z}\\
C^x_j & \DFAS \SETT{z \in C_{j+k+1}: x \text{ is a subface of } z }\\
& = \bigcup_{y \in C_j} C^{x,y}_{k,j}\\
C^x & \DFAS \bigcup_{j \geq 0} C^x_j
\end{aligned}
\]
That is, $C^x_j$ consists of all $(j+k+1)$-simplices of $C$ with subface $x$ and complementary subface in $C_j$. And $C^x$ consists of those simplices of $C$ which have $x$ as a subface.

\BX

\begin{newdef}
{\bf ($\CPL(x,\mu;y,\LAM))$}
\label{complementary-subface-notations}
\index{$\CPL(x,\mu;y,\LAM)$}
\index{$\CPL(x,y)$}
\index{Cpl$(x,\mu;y,\LAM)$}
{\rm
\SL

Let $j,k \geq 0$ and suppose $x \in C_k$, $y \in C_j$. Given two complementary vertex functions $\mu:[k] \to [j+1]$ and $\LAM:[j] \to [k+1]$, we define
$\CPL(x,\mu;y,\LAM)$ to consist of those $z \in C^{x,y}_{k,j}$ such that $\mu$ is the vertex function for $x$ in $z$ and $\LAM$ is the vertex function for $y$ in $z$.

That is, $z \in \CPL(x,\mu;y,\LAM)$ iff $\VL_z(x) = \SHR{\mu}[k]$ and $\VL_z(y) = \SHR{\LAM}[j]$.
}

\BX
\end{newdef}
\MS

\NI Abbreviated notation: We sometimes abbreviate $\CPL(x,\mu;y,\LAM)$ by $\CPL(x,y)$.
\MS

Theorem \refpage{complementary-face-idents} states that if $z \in \CPL(x,\mu;y,\LAM)$ then, for each $p \in [j] = [\DIM(y)]$, and each $q \in [k] = [\dim(x)]$
\[
\begin{aligned}
d_{\LAM(p)+p}(z) & \in  \CPL(x,\SIG_p \mu;d_p(y), \LAM \del_p)\\
s_{\LAM(p)+p}(z) & \in  \CPL(x,\del_{p+1} \mu;s_p(y), \LAM \SIG_p)
\end{aligned}
\]
and
\[
\begin{aligned}
d_{\mu(q)+q}(z) & \in  \CPL(d_q(x),\mu \del_q;y, \SIG_q \LAM )\\
s_{\mu(q)+q}(z) & \in  \CPL(s_q(x),\mu \SIG_q ;y, \del_{q+1}\LAM)
\end{aligned}
\]

\begin{newdef}
{\rm
\SL

Given $k \geq 0$, $x \in C_k$, $j>0$, $z \in \CPL(x,\mu;y,\LAM) \SBS C^x_j$ and $p \in [j]$ then
define
\[
d^x_p(z) \DFAS d_{p+\LAM(p)}(z) = d_{\SHR{\LAM}(p)}(z)
\]
\MS

\NI For $j \geq 0$ and $p \in [j]$ then define
\[
s^x_p(z) \DFAS s_{p+\LAM(p)}(z) = s_{\SHR{\LAM}(p)}(z)
\]
\BX
}
\end{newdef}

It is a corollary of Theorem \refpage{complementary-face-idents} that if $z \in C^x_j$ then for all $p \in [j]$, $s^x_p(z) \in C^x_{j+1}$ and, if $j>0$ then $d^x_p(z) \in C^x_{j-1}$.
\MS

\begin{newthm} 
\label{C^x-is-simplicial-set}
{\rm
\SL

Fix $k \geq 0$ and $x \in C_k$. Then $C^x$ is a simplicial set with face operators $d^x_p$ and degeneracy operators $s^x_p$.
}
\end{newthm}

\Proof

We will verify each simplicial identity directly from the definition by transforming the alleged identity in $C^x$ into a corresponding one in $C$. For each of these, we start with an arbitrary $z \in \CPL(x,\mu;y,\LAM) \SBS C^x_j$. 

\MS

\NI {\bf $\bullet$ Check $d^x_p d^x_q \ISIT d^x_{q-1} d^x_p$} where $p<q$ in $[j]$ and $j \geq 2$:
\[
d^x_p d^x_q(z) = d^x_p(d_{q+\LAM(q)}(z)) = d_{p+(\LAM \del_q)(p)}(d_{q+\LAM(q))}(z)) = d_{p+\LAM(p)} d_{q+\LAM(q)}(z)
\]
and
\[
d^x_{q-1}d^x_p(z) = d^x_{q-1}(d_{p+\LAM(p)}(z)) = d_{q-1+(\LAM \del_p)(q-1)}d_{p+\LAM(p)}(z) = d_{q-1+(\LAM(q))}d_{p+\LAM(p)}(z)
\]
Therefore $d^x_p d^x_q(z)= d^x_{q-1} d^x_p(z)$ by the corresponding face identity in $C$.
\MS

\NI {\bf $\bullet$ Check $d^x_p s^x_q \ISIT s^x_{q-1} d^x_p$} where $p<q$ in $[j]$ and $j \geq 1$:
\[
d^x_p s^x_q(z) = d^x_p s_{q+\LAM(q)}(z) = d_{p+(\LAM \SIG_q)(p)}s_{q+\LAM(q)}(z) = d_{p+\LAM(p)}s_{q+\LAM(q)}(z) 
\]
and
\[
s^x_{q-1} d^x_p(z) = s^x_{q-1} d_{p+\LAM(p)}(z) = s_{q-1+(\LAM \del_p)(q-1)} d_{p+\LAM(p)}(z) = s_{q-1+\LAM(q)} d_{p+\LAM(p)}(z)
\]
The claim follows from the identity in $C$.
\MS

\NI {\bf $\bullet$ Check $d^x_p s^x_p \ISIT 1$:}
\[
d^x_p s^x_p(z) = d_{p+(\LAM \SIG_p)(p)} s_{p+\LAM(p)}(z) =
d_{p+\LAM(p)} s_{p+\LAM(p)}(z) =z
\]
\MS

\NI {\bf $\bullet$ Check $d^x_{p+1} s^x_p \ISIT 1$:}
\[
d^x_{p+1} s^x_p(z) = d_{p+1+(\LAM \SIG_p)(p+1)} s_{p+\LAM(p)}(z)
= d_{p+1+\LAM(p)} s_{p+\LAM(p)}(z) = z
\]
\MS

\NI {\bf $\bullet$ Check $d^x_p s^x_q \ISIT s^x_q d^x_{p-1}$ where $p>q+1$ is $[j]$:}
\[
d^x_p s^x_q(z) = d_{p+(\LAM \SIG_q)(p)} s_{q+\LAM(q)}(z) = 
d_{p+\LAM(p-1)} s_{q+\LAM(q)}(z) = s_{q+\LAM(q)} d_{p-1+\LAM(p-1)}(z)
\]
and
\[
s^x_q d^x_{p-1}(z) = s_{q+(\LAM \del_{p-1})(q)} d_{p-1+\LAM(p-1)}(z)
= s_{q+\LAM(q)} d_{p-1+\LAM(p-1)}(z)
\]
\MS

\NI {\bf $\bullet$ Check $s^x_p s^x_q \ISIT s^x_{q+1} s^x_p$} where $p \leq q$ in $[j]$:
\[
s^x_p s^x_q(z) = s_{p+(\LAM \SIG_q)(p)} s_{q+\LAM(q)}(z)
= s_{p+\LAM(p)} s_{q+\LAM(q)}(z) = s_{q+1+\LAM(q)}s_{p+\LAM(p)}(z)
\]
and
\[
s^x_{q+1} s^x_p(z) = s_{q+1+(\LAM \SIG_p)(q+1)} s_{p+\LAM(p)}(z)
= s_{q+1+\LAM(q)} s_{p+\LAM(p)}(z)
\]
\qed

\begin{newdef}
{\bf (Complementary complex)}
\label{complementary-complex}
\index{complementary complex}
{\rm
\SL

Given a simplicial set $C$, $k \geq 0$ and $x \in C_k$, we will call $C^x$ the {\bf complementary complex fixing $x$.}
}
\end{newdef}
\BX
\MS

\NI {\bf Comments:}
\MS

\NI 1. In the proof above we noted that $p<q$ implies $d^x_{p} d^x_q = d_{p+\LAM(p)} d_{q+\LAM(q)}$ because $\LAM \del_q(p) = \LAM(p)$ when $p<q$. This extends to any subface-permissible sequence $p_1 \DDD{<} p_r$:
\[
d^x_{p_1} \DDD{} d^x_{p_r} = d_{p_1+\LAM(p_1)} \DDD{} d_{p_r+\LAM(p_r)} = 
d_{\SHR{\LAM}(p_1)} \DDD{} d_{\SHR{\LAM}(p_r)}
\]
\MS

\NI A convenient notation: if $A=\SETT{p_1 \DDD{<} p_r}$ then $\LAM\SH(A) = \SETT{\LAM\SH(p_1) \DDD{<}\LAM\SH(p_r)}$ and so $d^x_A = d_{\LAM\SH(A)}$
\MS

\NI 2. It follows from Theorems \refpage{complementary-face-idents} and \refpage{C^x-is-simplicial-set} that 
\[
\begin{aligned}
z \in \CPL(x,\mu;y,\LAM) & \IMP \\
& \VERT{p}^x(z) \in \CPL(x,\mu';\VERT{p}(y),\VERT{p}(\LAM))
\end{aligned}
\]
where 
\[
\begin{aligned}
\VERT{p}(\LAM) &= \bigl( \LAM(p) \bigr):[0] \to [k+1] \\
\mu' &= \SIG_j \DDD{} \omit{\SIG_p} \DDD{} \SIG_0 \mu
\end{aligned}
\]
So $\VL_{\VERT{p}^x(z)}(\VERT{p}(y)) = \SETT{\LAM(p)}$ and $\VL_{\VERT{p}^x(z)}(x) =[k+1]-\SETT{\LAM(p)}$. Therefore
\[
\begin{aligned}
(\mu')^\sharp &= (0 \DDD{,} \LAM(p)-1,\LAM(p)+1 \DDD{,} k+1) \\
\mu' &= (0 \DDD{,} 0,1 \DDD{,} 1)
\end{aligned}
\]
that is,
\[
\mu'(q) = \left\{
\begin{array}{ll}
0 & \text{if } q \in [0,\LAM(p)-1]\\
1 & \text{if } q \in [\LAM(p),k+1]
\end{array}
\right.
\]

\begin{example}
{\rm
Let $x \in C_2$ and $y \in C_4$ be complementary subfaces of $z \in C^x_4 \SBS C_7$, where $\VL_z(x) =\SETT{1,4,5}$ and $\VL_z(y) =\SETT{0,2,3,6,7}$. Then $z \in \CPL(x,\mu;y,\LAM)$ where $\mu = (1,3,3):[2] \to [5]$ and $\LAM = (0,1,1,3,3):[4] \to [3]$. The partition of the vertices of $z$ into vertices of $x$ and vertices of $y$ in shown in the table:
\[
z = 
\begin{array}{|c|cccccccc|}
\hline
& 0 & 1 & 2 & 3 & 4 & 5 & 6 & 7 \\
\hline
x &&\bullet&&&\bullet&\bullet&&\\
\hline
y &\bullet&&\bullet&\bullet&&&\bullet&\bullet\\
\hline
\end{array}
\]To illustrate the previous comment,
$\VERT{1}^x(z)= d^x_0 d^x_2 d^x_3 d^x_4(z)= d_0 d_3 d_6 d_7(z) \in C^x_0 \SBS C_3$ and the corresponding table is
\[
\VERT{1}^x(z)=
\begin{array}{|c|cccc|}
\hline
& 0 & 1 & 2 & 3  \\
\hline 
x &\bullet&&\bullet&\bullet\\
\hline
\VERT{1}(y) &&\bullet&&\\
\hline
\end{array}
\]
Here, $\VERT{1}^x(z) \in \CPL(x,\mu';\VERT{1}(y),\LAM')$ where $\VL_{\VERT{1}^x(z)}(x) =\SETT{0,2,3}$, $\mu'=(0,1,1):[2] \to [1]$ and $\VL_{\VERT{1}^x(z)}(y) =\SETT{1}$ and $\LAM'=\VERT{1}(\LAM)=(1):[0] \to [3]$. 
}
\end{example}

\begin{newlem}
\label{E-lemma}
{\rm
\SL

Given a simplicial set $C$, $k \geq 0$ and $x \in C_k$, then there 
is a simplicial map $E^x : C^x \to C \times \Delta[k+1]$ defined at dimension $j$ by $E^x_j(z) \DFAS (y,\LAM)$ where $y$ is the subface of $z$ complementary to $x$ and $[j]\XRA{\LAM} [k+1]$ is the vertex function of $y$ in $z$.

That is, if $z \in \CPL(x,\mu;y,\LAM)$ then $E^x(z) = (y,\LAM)$.
}
\end{newlem}

\Proof

This follows immediately from the preceding theorem and the definitions.

\qed
\MS

\begin{newlem}
{\rm
\SL 

Suppose $C$ and $D$ are simplicial sets, and $F:C \to D$ a simplicial map. Given any $k \geq 0$ and any $x \in C_k$ then
\begin{enumerate}
\item 
Given any $z \in \CPL(x,\mu;y,\LAM)$ then $F_{j+k+1}(z) \in \CPL(F_k(x),\mu;F_j(y),\LAM)$.

\item
$F$ determines the simplicial map
\[
F^x : C^x \to D^{F_k(x)}, \quad \ALL j \geq 0,\;F^x_j(z) \DFAS F_{j+k+1}(z)
\]
\end{enumerate}
}
\end{newlem}

\Proof

\begin{enumerate}
\item 
Suppose $z \in \CPL(x,\mu;y,\LAM)$ where $y \in C_j$. Denote $F_{j+k+1}(z)$ by $z'$. Since $x$ is a subface of $z$ and $F$ is simplicial then $x' = F_k(x)$ is a subface of $z'$. Let $y'$ denote the complement of $x'$ in $z'$.

For each $q \in [k]$, $\VERT{q}(x) = \VERT{\mu(q)+q}(z)$. Therefore we have
\[
\VERT{q}(F_k(x)) = F_0(\VERT{q}(x)) = F_0(\VERT{\mu(q)+q}(z)) = \VERT{\mu(q)+q}(F_{j+k+1}(z))
\]
That is, for each $q \in [k]$, $\VERT{q}(x') = \VERT{\mu(q)+q}(z')$ i.e. $\mu$ is the vertex function for $x'$ in $z'$. It follows that $\LAM$ is the vertex function for $y'$ as a subface of $z'$. Finally,
\[
y' = d_{\mu(0)+0} \DDD{} d_{\mu(k)+k}(z') = 
F_j(d_{\mu(0)+0} \DDD{} d_{\mu(k)+k}(z)) = F_j(y)
\]
This proves the first claim.

\item

For each $j>0$ and each $p \in [j]$
\[
d^{F_k(x)}_p F^x_j = d_{\LAM(p)+p}F_{j+k+1} = F_{j+k} d_{\LAM(p)+p} = F^x_{j-1} d^x_p
\]
Similarly, for each $j \geq 0$ and each $p \in [j]$
\[
s^{F_k(x)}_p F^x_j = s_{\LAM(p)+p} F_{j+k+1} = F_{j+k+2} s_{\LAM(p)+p} = 
F^x_{j+1} s^x_p
\]
Therefore, $F^x$ is a simplicial map.
\end{enumerate}
\qed

\begin{newthm}
\label{cpl-commute}
{\rm
\SL

Let $C$ be a simplicial set. Suppose $j,k > 0$, $x \in C_k$ and $y \in C_j$. Then for all $p \in [j]$, all $q \in [k]$ and all $z \in C^{x,y}_{k,j}$:
\begin{eqnarray}
&& d^{d_p(y)}_q d^x_p(z) = d^{d_q(x)}_p d^y_q(z) \label{cplcommute-1}\\
&& s^{d_p(y)}_q d^x_p(z) = d^{s_q(x)}_p s^y_q(z)\label{cplcommute-2}
\end{eqnarray}
That is, the following diagrams commute.

\setlength{\unitlength}{1in}
\begin{picture}(4.5,1.3)
\put(0,1){\makebox(0,0){$C^{x,y}_{k,j}$}}
\put(1.5,1){\makebox(0,0){$C^{x,d_p(y)}_{k,j-1} $}}
\put(0,0){\makebox(0,0){$C^{d_q(x),y}_{k-1,j}$}}
\put(1.5,0){\makebox(0,0){$C^{d_q(x),d_p(y)}_{k-1,j-1}$}}
\put(.4,1){\vector(1,0){.55}}
\put(.6,1.1){\MB{d^x_p}}
\put(.5,0){\vector(1,0){.45}}
\put(.75,.15){\MB{d^{d_q(x)}_p}}
\put(0,.8){\vector(0,-1){.6}}
\put(-.2,.5){\MB{d^y_q }}
\put(1.5,.8){\vector(0,-1){.6}}
\put(1.75,.5){\MB{d^{d_p(y)}_q}}
\put(.75,.5){\makebox(0,0){\eqref{cplcommute-1}}}

\put(3,1){\makebox(0,0){$C^{x,y}_{k,j}$}}
\put(3,0){\makebox(0,0){$C^{s_q(x),y}_{k+1,j}$}}
\put(4.5,1){\makebox(0,0){$C^{x,d_p(y)}_{k,j-1}$}}
\put(4.5,0){\makebox(0,0){$C^{s_q(x),d_p(y)}_{k+1,j-1}$}}
\put(3.4,1){\vector(1,0){.55}}
\put(3.6,1.1){\MB{d^x_p}}
\put(3.5,0){\vector(1,0){.45}}
\put(3.75,.15){\MB{d^{s_q(x)}_p}}
\put(3,.8){\vector(0,-1){.6}}
\put(2.8,.5){\MB{s^y_q}}
\put(4.5,.8){\vector(0,-1){.6}}
\put(4.75,.5){\MB{s^{d_p(y)}_q}}
\put(3.75,.5){\makebox(0,0){\eqref{cplcommute-2}}}
\end{picture}
}
\end{newthm}
\hspace{1cm}

\Proof

For both claims, we start with an arbitrary $z \in \CPL(x,\mu;y, \LAM) \SBS C^{x,y}_{k,j}$. By Theorem \refpage{complementary-face-idents} we have
\renewcommand{\arraystretch}{1.5}
\[
\begin{array}{l}
d^x_p(z) = d_{\LAM(p)+p}(z) \in \CPL(x,\SIG_p \mu;d_p(y), \LAM \del_p)\\
d^y_q(z) = d_{\mu(q)+q}(z) \in \CPL(d_q(x),\mu \del_q ;y, \SIG_q \LAM)\\
d^{d_p(y)}_q d^x_p(z) = d_{(\SIG_p \mu)(q)+q}\; d_{\LAM(p)+p}(z)\\
d^{d_q(x)}_p d^y_q(z) = d_{(\SIG_q \LAM)(p)+p} \; d_{\mu(q)+q}(z)\\
\\
s^y_q(z) = s_{\mu(q)+q}(z) \in \CPL(s_q(x),\mu \SIG_q; y, \del_{q+1} \LAM)\\
s^{d_p(y)}_q d^x_p(z) = s_{(\SIG_p \mu)(q)+q}\; d_{\LAM(p)+p}(z)\\
d^{s_q(x)}_p s^y_q(z) = d_{(\del_{q+1} \LAM)(p)+p}\; s_{\mu(q)+q}(z)
\end{array}
\]
Next:
\[
(\SIG_p \mu)(q)+q = \left\{
\begin{array}{ll}
\mu(q)+q & \text{if } \mu(q) \leq p\\
\mu(q)-1+q & \text{if } \mu(q)>p
\end{array}
\right.
\]
\[
(\SIG_q \LAM)(p)+p = \left\{
\begin{array}{ll}
\LAM(p)+p & \text{if } \LAM(p) \leq q\\
\LAM(p)-1+p & \text{if } \LAM(p)>q
\end{array}
\right.
\]
\[
(\del_{q+1} \LAM)(p)+p = \left\{
\begin{array}{ll}
\LAM(p)+p & \text{if } \LAM(p)<q+1\\
\LAM(p)+1+p & \text{if } \LAM(p) \geq q+1
\end{array}
\right.
\]
Now we divide the argument into cases: $\mu(q) \leq p$ and $\mu(q)>p$.
\MS

\NI {\bf Case $\mu(q) \leq p$:}

By lemma \refpage{lambda-least}, $\mu(q) \leq p$ iff $\LAM(p)>q$. Also, $\mu(q)+q \leq p+q < \LAM(p)+p$, in this case.

For equation \eqref{cplcommute-1}:
\[
d^{d_p(y)}_q d^x_p(z) = d_{(\SIG_p \mu)(q)+q}\; d_{\LAM(p)+p}(z) = d_{\mu(q)+q}\; d_{\LAM(p)+p}(z)
\]
and
\[
d^{d_q(x)}_p d^y_q(z) = d_{(\SIG_q \LAM)(p)+p} \; d_{\mu(q)+q}(z) = 
d_{\LAM(p)-1+p} \; d_{\mu(q)+q}(z)
\]
Since $\mu(q)+q < \LAM(p)+p$, then $d^{d_p(y)}_q d^x_p(z) = d^{d_p(y)}_q d^x_p(z)$ is a face identity in $C$.
\BS

For equation \eqref{cplcommute-2} in {\em this} case:
\[
s^{d_p(y)}_q d^x_p(z) = s_{(\SIG_p \mu)(q)+q}\; d_{\LAM(p)+p}(z) = 
s_{\mu(q)+q}\; d_{\LAM(p)+p}(z)
\]
and
\[
d^{s_q(x)}_p s^y_q(z) = d_{(\del_{q+1} \LAM)(p)+p}\; s_{\mu(q)+q}(z) = 
d_{\LAM(p)+1+p}\; s_{\mu(q)+q}(z)
\]
Therefore equation \eqref{cplcommute-2} follows as a simplicial identity in $C$.
\MS

\NI {\bf Case $\mu(q)>p$:}

By lemma \refpage{lambda-least}, $\mu(q)>p$ iff $\LAM(p) \leq q$. And, in this case, $\mu(q)+q > p + \LAM(p)$.

For equation \eqref{cplcommute-1}:
\[
d^{d_p(y)}_q d^x_p(z) = d_{(\SIG_p \mu)(q)+q}\; d_{\LAM(p)+p}(z) = d_{\mu(q)-1+q}\; d_{\LAM(p)+p}(z)
\]
and
\[
d^{d_q(x)}_p d^y_q(z) = d_{(\SIG_q \LAM)(p)+p} \; d_{\mu(q)+q}(z) = 
d_{\LAM(p)+p} \; d_{\mu(q)+q}(z)
\]
So equation \eqref{cplcommute-1} is a simplicial identity in $C$.
\BS

For equation \eqref{cplcommute-2} in this case:
\[
s^{d_p(y)}_q d^x_p(z) = s_{(\SIG_p \mu)(q)+q}\; d_{\LAM(p)+p}(z) = 
s_{\mu(q)-1+q}\; d_{\LAM(p)+p}(z)
\]
and
\[
d^{s_q(x)}_p s^y_q(z) = d_{(\del_{q+1} \LAM)(p)+p}\; s_{\mu(q)+q}(z) = 
d_{\LAM(p)+p}\; s_{\mu(q)+q}(z)
\]
Therefore equation \eqref{cplcommute-2} also follows as a simplicial identity in $C$.

\qed

\begin{newcor}
\label{cpl-vertex-cor}
{\rm
\SL

Suppose $C$ is a simplicial set, $j,k > 0$, $x \in C_k$ and $y \in C_j$. Then for all $p \in [j]$, all $q \in [k]$ and all $z \in C^{x,y}_{k,j}$
\[
d^{\VERT{p}(y)}_q(\VERT{p}^x(z)) = \VERT{p}^{d_q(x)}(d^y_q(z))
\]
This relates the vertices of $z \in C^x_j$ (where $y$ is complementary to $x$ in $z$) with those of $d^y_q(z) \in C^{d_q(x)}_j$. 
That is: given $z \in C^{x,y}_{k,j}$, and, for each $p \in [j]$ we write $t_p = \VERT{p}^x(z)$ and $t'_p = \VERT{p}^{d_q(x)}(d^y_q(z))$ then $t'_p = d_q^{\VERT{p}(y)}(t_p)$.
}
\end{newcor}

\Proof

By concatenating diagram \refpage{cplcommute-1} we obtain the commutative diagram
\BS

\begin{center}
\begin{picture}(3,1.2)
\put(0,1){\makebox(0,0){$C^{x,y}_{k,j}$}}
\put(3,1){\makebox(0,0){$C^{x,\VERT{p}(y)}_{k,0}$}}
\put(0,0){\makebox(0,0){$C^{d_q(x),y}_{k-1,j}$}}
\put(3,0){\makebox(0,0){$C^{d_q(x),\VERT{p}(y)}_{k-1,0}$}}
\put(.4,1){\vector(1,0){1.9}}
\put(.5,0){\vector(1,0){1.8}}
\put(0,.8){\vector(0,-1){.6}}
\put(3,.8){\vector(0,-1){.6}}
\put(-.2,.5){\makebox(0,0){$d^y_q$}}
\put(3.4,.5){\makebox(0,0){$d^{\VERT{p}(y)}_q$}}
\put(1.4,1.2){\makebox(0,0){$d^x_0 \DDD{} \omit{d^x_p} \DDD{} d^x_j$}}
\put(1.4,-.2){\makebox(0,0){$d^{d_q(x)}_0 \DDD{} \omit{d^{d_q(x)}_p} \DDD{} d^{d_q(x)}_j$}}
\end{picture}
\end{center}
\MS

\NI showing that
\[
\begin{aligned}
d^{\VERT{p}(y)}_q\; d^x_0 \DDD{} \omit{d^x_p} \DDD{}d^x_j(z) &=
d^{\VERT{p}(y)}_q( \VERT{p}^x(z))\\
&= d^{d_q(x)}_0 \DDD{} \omit{d^{d_q(x)}_p} \DDD{}d^{d_q(x)}_j\; d^y_q(z)\\
&= \VERT{p}^{d_q(x)}(d^y_q(z))
\end{aligned}
\]
as claimed.

\vskip.3cm
\qed
\MS

Given $C$, $k \geq 0$, $x \in C_k$ and $z \in C^x_j$ then either $\VERT{j+k+1}(z) = \VERT{k}(x)$, or not. The distinction gives rise to a subcomplex of $C^x$.

\begin{newdef}
\label{Cxplus-def}
\index{$C^{x^+}$}
\index{C$^{x^+}$}
{\rm
\SL

Given a simplicial set $C$, $k \geq 0$ and $x \in C_k$ then for each $j \geq 0$ we define
\[
C^{x^+}_j \DFAS \SETT{z \in C^x_j : \VERT{j+k+1}(z) = \VERT{k}(x) }
\]
That is, if $z \in \CPL(x,\mu;y,\LAM) \SBS C^x_j$ then $z \in C^{x^+}_j$ iff $\mu(k)=j+1$ iff $\LAM(j) \leq k$.
}
\BX
\end{newdef}

\begin{newlem}
{\rm
\SL

Let $C$ be a simplicial set, $k \geq 0$ and $x \in C_k$. Then $C^{x+}$ is a subcomplex of $C^x$.
}
\end{newlem}

\Proof

Suppose $j>0$ and $z \in \CPL(x,\mu;y,\LAM) \SBS C^{x^+}_j$. Then for any $p \in [j]$ we have 
\[
d^x_p(z) \in \CPL(x,\SIG_p \mu;d_p(y), \LAM \del_p) \SBS C^x_{j-1}
\]
Since $(\SIG_p \mu)(k) = \SIG_p(j+1) = j = (j-1)+1$ then $d^x_p(z) \in C^{x+}_{j-1}$.
\MS

Next suppose $j \geq 0$. Then for any $p \in [j]$ we have 
\[
s^x_p(z) \in \CPL(x,\del_{p+1} \mu;s_p(y), \LAM \SIG_p) \SBS C^x_{j+1}
\]
Since $(\del_{p+1} \mu)(k) = \del_{p+1}(j+1) = j+2 = (j+1)+1$ then $s^x_p(z) \in C^{x^+}_{j+1}$.

\qed

\subsection{\NICOMP{n}{i} structure for $C^x$}

\begin{newthm}
\label{C^x-is-a-composer-thm}
{\rm
\SL

Suppose $C$ is an \NICOMP{n}{i}, $k \geq 0$ and $x \in C_k$. Then $C^x$ is also an \NICOMP{n}{i}.
}
\end{newthm}

\Proof

It suffices to show that $\phi^x_{n+1,i} : C^x_{n+1} \to \BOX{i}{n+1}{C^x}$ is an isomorphism. Since $C$ is an \NICOMP{m}{i} for all $m \geq n$ then the same approach would prove that $\phi^x_{m+1,i} : C^x_{m+1} \to \BOX{i}{m+1}{C^x}$ is an isomorphism for all $m>n$.
\MS

To begin, suppose
\[
( \OM{z}{n+1}{i} ) \in \BOX{i}{n+1}{C^x}
\]
To set notation, for each $p \in [n+1]-\SETT{i}$, $z_p \in \CPL(x,\mu_p;y_p,\LAM_p) \in C^x_n$ where $\mu_p:[k] \to [n+1]$, $\LAM_p:[n] \to [k+1]$.
\MS

\NI {\bf $\bullet$ Claim 1:} There exists $y \in C_{n+1}$ and $\LAM : [n+1] \to [k+1]$ such that for all $p \in [n+1]-\SETT{i}$, $\LAM \del_p = \LAM_p$ and $d_p(y) = y_p$.
\MS

\NI {\bf Reason:} By hypothesis, for all $p<q$ in $[n+1]-\SETT{i}$ we have $d^x_p(z_q) = d^x_{q-1}(z_p)$. Since
\[
\begin{aligned}
d^x_p(z_q) & \in  \CPL(x,\mu'_q;d_p(y_q), \LAM_q \del_p)\\
d^x_{q-1}(z_p) & \in  \CPL(x,\mu'_p;d_{q-1}(y_p), \LAM_p \del_{q-1})
\end{aligned}
\]
it follows that $d_p(y_q) = d_{q-1}(y_p)$ and $\LAM_q \del_p = \LAM_p \del_{q-1}$ for all such $p<q$. Then, for general combinatorial reasons (see page \pageref{delta-n-is-hypergroupoid}), there exists a unique $\LAM : [n+1] \to [k+1]$ such that $\LAM \del_p = \LAM_p$ for all $p \neq i$. Using $(n,i)$-composition in $C$, $y = \COMP{n}{i}(\OM{y}{n+1}{i})$ has the property that $d_p(y) = y_p$, all $p \neq i$.
\MS

\NI {\bf Notation:} In what follows, we will use the notations
\[
A = \SETT{\LAM^\sharp(p) : p \in [n+1]} = \SETT{u_0 \DDD{<} u_n}
\]
and
\[
B = \SETT{\mu^\sharp(q) : q \in [k]} = \SETT{w_0 \DDD{,} w_k}
\]
where $\mu$ is the unique vertex function $[k] \to [n+2]$ which is complementary to $\LAM$.
\MS

\NI {\bf $\bullet$ Claim 2:} Let $T \DFAS [n+k+2]-\SETT{\LAM^\sharp(p):p \in [n+1]-\{i\}}$. Then $\SETT{z_p:p \in [n+1]-\{i\}} \in \BOX{T}{n+k+2}{C}$.
\MS

\NI {\bf Reason:} For clarity, re-index $\SETT{z_p:p \in [n+1]-\{i\}}$ by denoting $z_p$ by $z'_{\LAM^\sharp(p)}$. The goal is to verify that whenever $\LAM^\sharp(p)<\LAM^\sharp(q)$ in $[n+k+2]$ then $d_{\LAM^\sharp(p)}(z'_{\LAM^\sharp(q)}) = d_{\LAM^\sharp(q)-1}(z'_{\LAM^\sharp(p)})$.

Now $\LAM^\sharp(p)<\LAM^\sharp(q) \iff p<q$ (since $\LAM^\sharp$ is strictly increasing) and so we have
\[
d^x_p(z_q) = d_{p+(\LAM \del_q)(p)}(z_q) = d_{p+\LAM(p)}(z_q) = d_{\LAM^\sharp(p)}(z'_{\LAM^\sharp(q)})
\]
and
\[
d^x_{q-1}(z_p) = d_{q-1+(\LAM \del_p)(q-1)}(z_p)= d_{q-1+\LAM(q)}(z_p) = 
d_{\LAM^\sharp(q)-1}(z'_{\LAM^\sharp(p)})
\]
The claimed equality follows from $d^x_p(z_q) = d^x_{q-1}(z_p)$.
\MS

\NI {\bf $\bullet$ Claim 3:} The $(n,i)$-structure of $C$ implies there exists a unique $z \in C_{n+k+2}$ such that for each $p \in [n+1]-\{i\}$, $d_{\LAM^\sharp(p)}(z) = z_p$. 
\MS

\NI {\bf Reason:} $|T| = k+2$. Write $T$ as $\SETT{t_1 \DDD{<} t_{k+2}}$. By definition, $\LAM^\sharp(i) = \LAM(i)+i \in T$. First we will calculate the index $q$ such that $t_q = \LAM(i)+i$.

Now $T = B \cup \SETT{\LAM(i)+i}$ and so
\[
B' \DFAS [0,\LAM(i)+i-1] \cap B = \text{ the set of predecessors of } \LAM(i)+i \text{ in } T
\]
which says that $q = 1+|B'|$. Next, if we set $A' \DFAS [0,\LAM(i)+i-1] \cap A$ then $[0,\LAM(i)+i-1] = A' \cup B'$ is a disjoint union. Therefore
\[
\LAM(i)+i = |[0,\LAM(i)+i-1]| = |A'|+|B'| = i+(q-1)
\]
and therefore $q = 1+\LAM(i)$. It then follows that $T$ is of type $q$ for $(n,i)$ because
\[
t_q = t_{1+\LAM(i)} = i+\LAM(i) = i+(1+\LAM(i))-1 = i+q-1
\]
By the Extended Compositions Theorem (Theorem \refpage{extended-comps-theorem}) there exists a unique $z \in C_{n+k+2}$ such that $d_{\LAM^\sharp(p)}(z) = z_p$ for all $p \in [n+1]-\{i\}$.
\MS

\NI {\bf $\bullet$ Claim 4:} $z \in \CPL(x,\mu;y,\LAM)$.
\MS

\NI {\bf Reason:} First, $x$ is a subface of $z$ because $x$ is a subface of $d_{\LAM^\sharp(p)}(z)$ for each $p \neq i$. Therefore $z \in \CPL(x,\mu';y',\LAM')$ for some complementary vertex functions $\mu' : [k] \to [n+2]$ and $\LAM':[n] \to [k+1]$. Now for each $p \in [n+1]-\{i\}$, $d^x_p(z) = z_p = d_{\LAM^\sharp(p)}(z)$. Since $d_{\LAM^\sharp(p)}(z) \in \CPL(x,\mu''_p;d_p(y'),\LAM' \del_p)$ and $z_p \in \CPL(x,\mu_p;d_p(y),\LAM \del_p)$, it follows that $y=y'$ because $d_p(y) = d_p(y')$ for all $p \neq i$, and $\LAM = \LAM'$ because $\LAM \del_p = \LAM' \del_p$ for all $p \neq i$. Therefore $\mu'=\mu$ because $\mu$ is the unique vertex function complementary to $\LAM$.
\MS

To conclude the proof, given any $z \in C^x_{n+1}$, claims 1-4 show that there is a $z' \in C^x$ such that $\phi^x_{n+1,i}(z') = \phi^x_{n+1,i}(z)$ and the uniqueness of $z'$ implies $z'=z$. That is, $\phi^x_{n+1,i}$ is an isomorphism.

\qed

\subsection{Low dimensional example of $C^x$}

In this subsection $C$ will be a \NICOMP{1}{1} ($n=1$; $C$ is a small category), $k=0$ and $x \in C_0$, fixed. According to Theorem \refpage{C^x-is-a-composer-thm}, $C^x$ will also be a \NICOMP{1}{1}.

For any $j \geq 0$ and $z \in C^x_j \SBS C_{j+1}$ then $z \in \CPL(x,\mu;y,\LAM)$ where $\mu : [0] \to [j+1]$, $\LAM : [j] \to [1]$ are the complementary vertex functions for $x$ and $y$ as subfaces of $z$. In this special situation, $\mu = (\mu(0))$ and $\mu^\sharp(0) = \mu(0)$. $B = \VL_z(x) = \SETT{\mu(0)}$ and $A = \VL_z(y) = [j+1]-\SETT{\mu(0)}$. There are exactly $j+2$ possible such complementary  vertex functions.
\MS

\NI Notation: For readability, given any $m$-simplex $w$ we will abbreviate $\VERT{t}(w)$ by $w^t$.
\MS

\NI When $j=0$:

$C^x_0 \SBS C_1$. Suppose $z \in \CPL(x,\mu;y,\LAM) \SBS C^x_0$.
The only two possibilities for $\mu$ and $\LAM$ are
\begin{enumerate}
\item 
$\mu = (0), \LAM = (1)$. Then $\VL_z(x) = \SETT{\mu^\sharp(0)} = \{0\}$ and $\VL_z(y) = \SETT{\LAM^\sharp(0)} \{ 1 \}$. Thus $z = x \to y$.

\item
$\mu = (1), \LAM = (0)$. Then $\VL_z(x) = \SETT{\mu^\sharp(0)} = \{1\}$ and $\VL_z(y) = \SETT{\LAM^\sharp(0)} \{ 0 \}$. Thus $z = y \to x$.

\end{enumerate}
\MS

\NI When $j=1$:

$C^x_1 \SBS C_2$ and $1$-simplices of $C^x$ are {\em maps} of $C^x$. Suppose $z \in \CPL(x,\mu;y,\LAM) \SBS C^x_1$. There are three possibilities for $\mu$ and $\LAM$, as follows.

\begin{enumerate}
\item 
$\mu=(0)$, $\LAM^\sharp=(1,2)$, $\LAM = (1,1)$. Then $z = x \to y^0 \to y^1$ with $d^x_0(z) = d_1(z) = x \to y^1$ and $d^x_1(z) = d_2(z) = x \to y^0$.

\item
$\mu=(1)$, $\LAM^\sharp=(0,2)$, $\LAM = (0,1)$. Then $z = y^0 \to x \to y^1$ with $d^x_0(z) = d_0(z) = x \to y^1$ and $d^x_1(z) = d_2(z) = y^0 \to x$.

\item
$\mu=(2)$, $\LAM^\sharp=(0,1)$, $\LAM = (0,0)$. Then $z = y^0 \to y^1 \to x$ with $d^x_0(z) = d_0(z) = y^1 \to x$ and $d^x_1(z) = d_2(z) = y^0 \to x$.

\end{enumerate}
\MS

\NI When $j=2$:

$C^x_2 \SBS C_3$ and $2$-simplices of $C^x$ are compositions in $C^x$. We will work through a specific illustration below. If $z \in \CPL(x,\mu;y,\LAM)$ then $y \in C_2$ and the complementary vertex functions are $\mu : [0] \to [3]$ and $\LAM : [2] \to [1]$. There are four possibilities, as follows.
\begin{enumerate}
\item 
$\mu = (0), \LAM^\sharp = (1,2,3), \LAM = (1,1,1)$. Then
\[
\begin{aligned}
z &= x \to y^0 \to y^1 \to y^2 \\
d^x_0(z) = d_1(z) &= x \to  y^1 \to y^2 \\
d^x_1(z) = d_2(z) &= x \to y^0 \to y^2 \\
d^x_2(z) = d_3(z) &= x \to y^0 \to y^1 
\end{aligned}
\]

\item 
$\mu = (1), \LAM^\sharp = (0,2,3), \LAM = (0,1,1)$. Then
\[
\begin{aligned}
z &= y^0 \to x \to y^1 \to y^2\\
d^x_0(z) = d_0(z) &= x \to y^1 \to y^2\\
d^x_1(z) = d_2(z) &= y^0 \to x \to y^2\\
d^x_2(z) = d_3(z) &= y^0 \to x \to y^1 
\end{aligned}
\]

\item
$\mu = (2), \LAM^\sharp = (0,1,3), \LAM = (0,0,1)$. Then
\[
\begin{aligned}
z &= y^0 \to y^1 \to x \to y^2\\
d^x_0(z) = d_0(z) &= y^1 \to x \to y^2\\
d^x_1(z) = d_1(z) &= y^0 \to x \to y^2\\
d^x_2(z) = d_3(z) &= y^0 \to y^1 \to x 
\end{aligned}
\]

\item
$\mu = (3), \LAM^\sharp = (0,1,2), \LAM = (0,0,0)$. Then
\[
\begin{aligned}
z &= y^0 \to y^1 \to y^2 \to x\\
d^x_0(z) = d_0(z) &= y^1 \to y^2 \to x\\
d^x_1(z) = d_1(z) &= y^0 \to y^2 \to x\\
d^x_2(z) = d_2(z) &= y^0 \to y^1 \to x 
\end{aligned}
\]

\end{enumerate}

Next, let's examine $(1,1)$-composition in $C^x$. Start with $(z_0,-,z_2) \in \BOX{1}{2}{C^x}$. That is, we're starting with $z_0 \in \CPL(x,\mu_0;y_0,\LAM_0) \SBS C^x_1$ and $z_2 \in \CPL(x,\mu_2;y_2,\LAM_2) \SBS C^x_1$ such that $d^x_0(z_2) = d^x_1(z_0)$. Since
\[
\begin{aligned}
d^x_0(z_2) & \in  \CPL(x, \mu'_2;d_0(y_2),\LAM_2 \del_0) \\
d^x_1(z_0) & \in  \CPL(x,\mu'_0;d_1(y_0),\LAM_0 \del_1)
\end{aligned}
\]
we know that $(y_0,-,y_2) \in \BOX{1}{2}{C}$ (i.e. $y_0$ and $y_2$ are composable) and that
\[
\LAM_2 \del_0 = (\LAM_2(1) ) = \LAM_0 \del_1 = ( \LAM_0(0))
\]
As noted in the proof of Theorem \refpage{C^x-is-a-composer-thm}, the equation $\LAM_2 \del_0 = \LAM_0 \del_1$ determines a unique $\LAM : [2] \to [1]$ such that $\LAM \del_0 = \LAM_0$ and $\LAM \del_2 = \LAM_2$, namely $\LAM = (\LAM_2(0),\LAM_2(1),\LAM_0(1))$. The unique $z \in C^x_2 \SBS C_3$ determine by composition is $z \in \CPL(x, \tilde{\mu}; y_2 y_0,\LAM)$ where $y_2 y_0$ denotes the composition of $y_2$ with $y_0$ in $C$ and $\tilde{\mu}$ is determined uniquely by $\LAM^\sharp = (\LAM_2(0),\LAM_2(1)+1,\LAM_0(1)+2)$.
\MS

\NI A specific numerical illustration of this composition:
\MS

Suppose $z_0 \in \CPL(x,\mu_0;y_0,(1,1))$ and $z_2 \in \CPL(x,\mu_2;y_2,(0,1))$. That is
\[
\begin{array}{l}
\LAM_0=(1,1), \quad \LAM^\sharp_0 = (1,2), \quad \mu_0 = (0)\\
\LAM_2=(0,1), \quad \LAM^\sharp_2 = (0,2), \quad \mu_2 = (1)
\end{array}
\]
where $y_2$ is composable with $y_0$ in $C$ corresponding to $y \in C_2$ with $y_1 = y_2 \circ y_0$

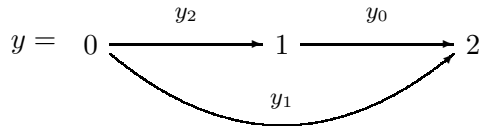
\begin{picture}(3,.5)
\put(1,0){\makebox(0,0){$0$}}
\put(2,0){\makebox(0,0){$1$}}
\put(3,0){\makebox(0,0){$2$}}
\qbezier(1.1,-.05)(2,-.8)(2.9,-.05)
\put(2.9,-.05){\vector(1,1){0}}
\put(1.1,0){\vector(1,0){.8}}
\put(2.1,0){\vector(1,0){.8}}
\put(1.5,.15){\MBS{y_2}}
\put(2.5,.15){\MBS{y_0}}
\put(2,-.3){\MBS{y_1}}
\put(.7,0){\makebox(0,0){$y=$}}
\end{picture}
\vskip.6in

\NI Then $\LAM = (0,1,1) : [2] \to [1]$, $\LAM^\sharp = (0,2,3) : [2] \to [3]$, $\tilde{\mu} = (1)$ and $z$ is $y^0 \to x \to y^1 \to y^2$ or, as a fancier diagram:

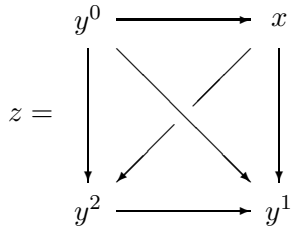
\begin{center}
\begin{picture}(1,1.5)
\put(0,1){\makebox(0,0){$y^0$}}
\put(1,1){\makebox(0,0){$x$}}
\put(0,0){\makebox(0,0){$y^2$}}
\put(1,0){\makebox(0,0){$y^1$}}
\put(.15,0){\vector(1,0){.7}}
\put(.15,1){\vector(1,0){.7}}
\put(0,.85){\vector(0,-1){.7}}
\put(1,.85){\vector(0,-1){.7}}
\put(.15,.85){\vector(1,-1){.7}}
\put(.85,.85){\line(-1,-1){.3}}
\put(.45,.45){\vector(-1,-1){.3}}
\put(-.3,.5){\makebox(0,0){$z=$}}
\end{picture}
\end{center}

\NI Here, $d^x_0(z) = d_0(z) = z_0$, $d^x_2(z) = d_3(z) = z_2$ and $d^x_1(z) = d_2(z)$ is the result of the composition in $C^x$.

\subsection{Special subcomplexes of $C^x$}

This section describes two subcomplexes of $C^x$, denoted $L^{x,t}$ and $R^{x,t}$, specializing to those $z \in C^x_j$ such that $\VL_z(x)$ contains a certain subinterval of $[0,j+k+1]$, where $x \in C_k$. We prove that if $C$ is an \NICOMP{n}{i} then so are $L^{x,t}$ and $R^{x,t}$. These \NICOMP{n}{i}s generalize the notion of comma-category, described in sections \refpage{comma-category-L} and \refpage{comma-category-R} below.

\subsubsection{$L^{x,t}$ and $R^{x,t}$}

\begin{newdef}
\label{Lxt and Rxt}
\index{$L^{x,t}$}
\index{$R^{x,t}$}
{\rm
\SL

Suppose $C$ is a simplicial set, $k \geq 0$, $x \in C_k$ and $t \in [k]$. 
\MS

For each $j \geq 0$:
\index{$L^{x,t}_j$}
\index{$L^{x,t}$}
\index{$R^{x,t}_j$}
\index{$R^{x,t}$}
\begin{eqnarray}
L^{x,t}_j & \DFAS & \SETT{ z \in C^x : [0,t] \SBS \VL_z(x)} \\
L^{x,t} & \DFAS & \bigcup_{j \geq 0} L^{x,t}_j\\
&& \\
R^{x,t}_j & \DFAS & \SETT{ z \in C^x : [j+k+1-t,j+k+1] \SBS \VL_z(x)}\\
R^{x,t} & \DFAS & \bigcup_{j \geq 0} R^{x,t}_j
\end{eqnarray}

\BX
}
\end{newdef}

\begin{newlem}
{\rm
\SL

Let $C$ be an arbitrary simplicial set. Suppose $x \in C_k$, $y \in C_j$ and $z \in \CPL(x,\mu;y,\LAM)$. Given any $t \in [k]$ the following statements are equivalent:
\begin{enumerate}
\item 
$[0,t] \SBS \VL_z(x)$

\item
For each $q \in [0,t]$, $\mu(q)=0$

\item
$\LAM(0)>t$.
\end{enumerate}
}
\end{newlem}

\Proof

$(1) \IMP (2)$: If $[0,t] \SBS \VL_z(x) = \mu\SH[k]$ then $\mu\SH(q) = \mu(q)+q = q$ for all $q \in [0,t]$. 

$(2) \IMP (3)$: Since $\VL_z(x)$ and $\VL_z(y)$ are disjoint, then $\LAM(0) =\LAM\SH(0)\notin [0,t]$ and therefore $\LAM(0)>t$. 

$(3) \IMP (1)$: If $\LAM(0)>t$ then $\LAM^\sharp(0)>t$ from which it follows that $[0,t] \cap \VL_z(y) = \MT$, hence $[0,t] \SBS \VL_z(x)$.

\qed
\MS

\begin{newlem}
{\rm
\SL

Let $C$ be an arbitrary simplicial set. Suppose $x \in C_k$, $y \in C_j$ and $z \in \CPL(x,\mu;y,\LAM)$. Let $t \in [k]$ and $m \DFAS j+k+1$. Then the following statements are equivalent:
\begin{enumerate}
\item 
$[m-t,\;m] \SBS \VL_z(x)$

\item
For each $q \in [k-t,k]$, $\mu(q)=j+1$.

\item
$\LAM(j) \leq k-t$
\end{enumerate}

}
\end{newlem}

\Proof

We will use that $\mu\SH[k]\; \cup \; \LAM\SH[j] = [m]$ is a disjoint union where $\DMAP{\mu\SH}{k}{m}$ and $\DMAP{\LAM\SH}{j}{m}$.
\MS

\NI $(1) \IMP (2)$: 
\MS

$[k] = [0,k-t-1] \cup [k-t,k]$ therefore $\mu\SH[k-t,k]=[m-t,m]$ since $\mu\SH$ is \SI. So for all $u \in [0,t]$, (setting $q = k-t+u$)
\[
\begin{aligned}
\mu\SH(k-t+u) &=m-t+u\\
\mu(k-t+u) &= m-t+u-(k-t+u) = m-k=j+1
\end{aligned}
\]

\NI $(2) \IMP (3)$: 
\MS

For each $u \in [0,t]$, $\mu(k-t+u) = j+1$ which implies $\mu\SH[k-t,k]=[m-t,m]$. Therefore since $\mu\SH[k] \cap \LAM\SH[j]=\MT$ then $[m-t,m] \cap \LAM\SH[j]=\MT$. It follows that
\[
\begin{aligned}
\LAM\SH(j) \leq m-t-1 & =j+k-t\\
\LAM(j) & \leq k-t
\end{aligned}
\]

\NI $(3) \IMP (1)$: 
\MS

$\LAM(j) \leq k-t$ implies $\LAM\SH(j) \leq k-t+j = m-t-1$ which implies $\LAM\SH[j] \SBS [0,m-t-1]$. Since $\LAM\SH[j] \; \cap \; \mu\SH[k]=\MT$ then 
$[m-t,m] \SBS \mu\SH[k]$ i.e. $[m-t,m] \SBS \VL_z(x)$.

\qed

\begin{newthm}
\label{G^{x,t} is subcomplex}
{\rm
\SL

Let $C$ be a simplicial set, $k \geq 0$, $x \in C_k$ and $t \in [k]$. Then $R^{x,t}$ and $L^{x,t}$ are subcomplexes of $C^x$.
}
\end{newthm}

\Proof

Let $z \in C^x_j$ where $z \in \CPL(x,\mu;y,\LAM)$ with vertex funxtion $\DMAP{\mu}{k}{j+1}$. Let $m \DFAS j+k+1$. For both claims we use that for all $p \in [j]$
\[
\begin{aligned}
d^x_p(z) & \in \CPL(x,\sigma_p \mu;d_p(y),\LAM \del_p) \\
s^x_p(z) & \in \CPL(x,\del_{p+1} \mu;s_p(y), \LAM \sigma_p)
\end{aligned}
\]

We will use lemma \refpage{cor-lambda-least} which states:
\[
(\sigma_p \mu)\SH = \sigma_{\LAM\SH(p)}\mu\SH \qquad \text{and} \qquad
(\del_{p+1} \mu)\SH = \del_{\LAM\SH(p)} \mu\SH
\]

\HD{That $R^{x,t}$ is a subcomplex of $C^x$:}
\MS

Suppose $z \in R^{x,t}_j$. By definition, $z \in \CPL(x,\mu;y,\LAM)$ and $\mu\SH[k] \SPS [m-t,m]$. We must show that for all $p \in [j]$ (and $j>0$), $d^x_p(z) \in R^{x,t}_{j-1}$ and, for all $j 
\geq 0$, $s^x_p(z) \in R^{x,t}_{j+1}$.
\MS

\HD{Face operators:}
\MS

We must show that $[m-1-t,m-1] \SBS (\sigma_p \mu)\SH[k]$. 
Since $\LAM\SH(p) \leq k-t<m-t$ then
\[
(\sigma_p \mu)\SH[k]=\sigma_{\LAM\SH(p)}\mu\SH[k] \SPS \sigma_{\LAM\SH(p)}[m-t,m] = [m-t-1,m-1]
\]
as claimed.
\MS

\HD{Degeneracy operators:}
\MS

To show $s^x_p(z) \in R^{x,t}_{j+1}$ we must show that $(\del_{p+1} \mu)\SH[k] \SPS [m-t+1,m+1]$.
Since $\LAM\SH(p)<m-t$ then 
\[
(\del_{p+1} \mu)\SH[k]=\del_{\LAM\SH(p)} \mu\SH[k] \SPS \del_{\LAM\SH(p)} \mu\SH[m-t,m]=[m-t+1,m+1]
\]
as claimed.
\BS

\HD{That $L^{x,t}$ is a subcomplex of $C^x$:}
\MS

By definition, $z \in \CPL(x,\mu;y,\LAM)$ such that $[0,t] \SBS \mu\SH[k]$.
\MS

\HD{Face operators:}
\MS

We must show that $(\sigma_p \mu)\SH[k] \SPS [0,t]$.
\[
(\sigma_p \mu)\SH[k] = \sigma_{\LAM\SH(p)}\mu\SH[k] \SPS \sigma_{\LAM\SH(p)}[0,t]
\]
But $\LAM\SH[j] \; \cap \; \mu\SH[k] = \MT$ implies $\LAM\SH(p)>t$ for all $p \in [j]$. Therefore $\sigma_{\LAM\SH(p)}[0,t] = [0,t]$.
\MS

\HD{Degeneracy operators:}
\MS

We must show that $(\del_{p+1} \mu)\SH[k] \SPS [0,t]$.
\[
(\del_{p+1} \mu)\SH[k] = \del_{\LAM\SH(p)} \mu\SH{k} \SPS \del_{\LAM\SH(p)}[0,t]
\]
But $\LAM\SH(p)>t$ (all $p$) implies $\del_{\LAM\SH(p)}[0,t] = [0,t]$.

\qed

\begin{newthm}
{\rm
\SL

Suppose $C$ is an \NICOMP{n}{i}, $k \geq 0$, $x \in C_k$ and $t \in [k]$. Then $L^{x,t}$ and $R^{x,t}$ are also \NICOMP{n}{i}s.
}
\end{newthm}

\Proof

It suffices to prove that $\BOX{i}{n+1}{G} \ISO G_{n+1}$ for $G = L^{x,t}$ and $G = R^{x,t}$. Since $C$ is also an \NICOMP{m}{i} for each $m>n$ the same argument would show that $\BOX{i}{m+1}{G} \ISO G_{m+1}$ for each $m>n$.
\MS

By Theorem \refpage{C^x-is-a-composer-thm}, $C^x$ is an \NICOMP{n}{i}. Now, for either $G=L^{x,t}$ or $G=R^{x,t}$, 
\[
(\OM{z}{n+1}{i} ) \in \BOX{i}{n+1}{G} \IMP (\OM{z}{n+1}{i} ) \in \BOX{i}{n+1}{C^x}
\]
Thus, there exists a unique $z \in C^x_{n+1}$, where $z \in \CPL(x,\mu;y,\LAM)$, such that $d^x_p(z) = z_p$ for each $p \in [n+1]-\{i\}$. Therefore it will suffice to show that $z \in G_{n+1}$ when $G=L^{x,t}$ and when $G=R^{x,t}$.
\BS

\NI {\bf Case $G=L^{x,t}$:}
\MS

By definition, and using the lemma for $L^{x,t}$ above:
\[
z \in L^{x,t}_{n+1} \iff [0,t] \SBS \VL_z(x) \iff \LAM(0)>t
\]
For each $p \neq i$, $d^x_p(z) = z_p \in \CPL(x,\mu';d_p(y),\LAM \del_p)$. By hypothesis, $z_p \in L^{x,t}_n$ and so, by the lemma, $(\LAM \del_p)(0)>t$. Now $n \geq 1 \IMP [n+1]-\SETT{0,i} \neq \MT$. So pick any $p \in [n+1]-\SETT{0,i}$. We get
\[
t<(\LAM \del_p)(0) = \LAM(0)
\]
and therefore $z \in L^{x,t}_{n+1}$.
\BS

\NI {\bf Case $G=R^{x,t}$:}
\MS

By the lemma above for $R^{x,t}$, $z \in R^{x,t}_{n+1} \iff \LAM(n+1)<k+1-t$. By hypothesis, for each $p \in [n+1] - \{i\}$, $z_p \in R^{x,t}_n$, and that implies $(\LAM \del_p)(n)<k+1-t$. Now $n \geq 1 \IMP [n+1]-\SETT{n+1,i} \neq \MT$. Pick any $p \in [n+1]-\SETT{n+1,i}$. Since $p \leq n$ then $k+1-t>(\LAM \del_p)(n) = \LAM(n+1)$, as required.

\qed

\subsubsection{$L^{x,0}$ when $x \in C_0$}
\label{comma-category-L}

For brevity, we will denote $L^{x,0}$ by $L^x$ in this subsection.
\MS

Suppose $z \in \CPL(x,\mu;y,\LAM) \SBS L^x_j \SBS C_{j+1}$, where $y \in C_j$, $\LAM : [j] \to [1]$ and $\mu : [0] \to [j+1]$.
Then $\mu(0)=0$ and we have $\VL_z(x) = \{0\}$ and $\VL_z(y) = \SETT{1,2 \DDD{,} j+1}$. That is, for all $p \in [j]$, $\LAM^\sharp(p)=p+1$ and $\LAM(p)=1$. For each $p \in [j]$, $d^x_p(z) = d_{p+1}(z)$.
\MS

If $j=0$ then $z$ is the 1-simplex of $C$: $x \to y$.
\MS

If $j=1$ then $z$ is a 2-simplex of $C$ looks like:

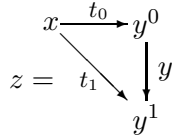
\begin{center}
\begin{picture}(.5,.5)
\put(0,.5){\makebox(0,0){$x$}}
\put(.5,.5){\makebox(0,0){$y^0$}}
\put(.5,0){\makebox(0,0){$y^1$}}
\put(.05,.5){\vector(1,0){.35}}
\put(.5,.4){\vector(0,-1){.3}}
\put(.05,.45){\vector(1,-1){.35}}
\put(.6,.25){\makebox(0,0){$y$}}
\put(-.1,.2){\makebox(0,0){$z = $}}
\put(.25,.57){\makebox(0,0){\footnotesize{$t_0$}}}
\put(.2,.2){\makebox(0,0){\footnotesize{$t_1$}}}
\end{picture}
\end{center}
where $d^x_0(z) = d_1(z)$ and $d^x_1(z) =d_2(z)$. If $C$ is a $(1,1)$-composer (i.e. a small category) then this 2-simplex is commutative.
\MS

If $j=2$ then $z$ is a 3-simplex of $C$:

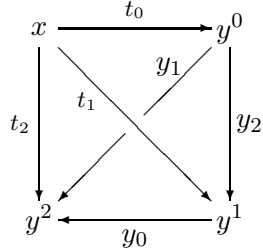
\begin{center}
\begin{picture}(1,1.1)
\put(0,1){\makebox(0,0){$x$}}
\put(1,1){\makebox(0,0){$y^0$}}
\put(1,0){\makebox(0,0){$y^1$}}
\put(0,0){\makebox(0,0){$y^2$}}
\put(0.1,1){\vector(1,0){.8}}
\put(.9,0){\vector(-1,0){.8}}
\put(0,.9){\vector(0,-1){.8}}
\put(1,.9){\vector(0,-1){.8}}
\put(.1,.9){\vector(1,-1){.8}}
\put(.9,.9){\line(-1,-1){.35}}
\put(.45,.45){\vector(-1,-1){.35}}
\put(.5,1.1){\makebox(0,0){\footnotesize{$t_0$}}}
\put(-.1,.5){\makebox(0,0){\footnotesize{$t_2$}}}
\put(.25,.62){\makebox(0,0){\footnotesize{$t_1$}}}
\put(1.1,.5){\makebox(0,0){$y_2$}}
\put(.5,-.1){\makebox(0,0){$y_0$}}
\put(.68,.8){\makebox(0,0){$y_1$}}
\end{picture}
\end{center}
where $y_p = d_p(y)$, $p \in [2]$. If $C$ is a $(1,1)$-composer then this is a commutative tetrahedron. Here, $d^x_p(z) = d_{p+1}(z)$. Note that the commutativity of $d^x_1(z)$, that is $y_1 \circ t_0 \ISIT t_2$ follows routinely from that of $d^x_0(z) = d_1(z)$ ($y_0 \circ t_1 = t_2$), $d^x_2(z) = d_3(z)$ ($y_2 \circ t_0 = t_1$) and that $y_0 \circ y_2 = y_1$.
\MS

Therefore in the case that $C$ is a $(1,1)$-composer, $L^x$ is the comma category $x \downarrow C$.

\subsubsection{$R^{x,0}$ when $x \in C_0$}
\label{comma-category-R}

This is similar to the discussion of $L^{x,0}$. For brevity we will denote $R^{x,0}$ by $R^x$.

Start with any $z \in \CPL(x,\mu;y,\LAM) \in R^{x,0}_j$, where $y \in C_j$, $\LAM : [j] \to [1]$ and $\mu : [0] \to [j+1]$. Here, $\VL_z(x) = \{ j+1 \}$ and $\VL_z(y) = \SETT{0 \DDD{,} j}$. That is, $\mu(0) = \mu^\sharp(0) = j+1$ and $\LAM^\sharp(p) = p$, $\LAM(p)=0$ for all $p \in [j]$. For each $p \in [j]$, $d^x_p(z) = d_p(z)$.
\BS

If $j=0$ then $z \in R^x_0$ is a 1-simplex $y \to x$ of $C$.
\MS

If $j=1$ then $z \in \CPL(x,\mu;y,\LAM) \SBS R^x_1$ is a 2-simplex of $C$ which looks like

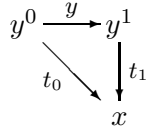
\begin{center}
\begin{picture}(.5,.5)
\put(0,.5){\makebox(0,0){$y^0$}}
\put(.5,.5){\makebox(0,0){$y^1$}}
\put(.5,0){\makebox(0,0){$x$}}
\put(.1,.5){\vector(1,0){.3}}
\put(.1,.4){\vector(1,-1){.3}}
\put(.5,.4){\vector(0,-1){.3}}
\put(.25,.57){\makebox(0,0){\footnotesize{$y$}}}
\put(.15,.2){\makebox(0,0){\footnotesize{$t_0$}}}
\put(.6,.25){\makebox(0,0){\footnotesize{$t_1$}}}
\end{picture}
\end{center}
with $d^x_0(z) = t_1$ and $d^x_1(z) = t_0$.
\BS

If $j=2$ then $z \in \CPL(x,\mu;y,\LAM) \SBS R^x_2$ is a 3-simplex of $C$ which looks like this:

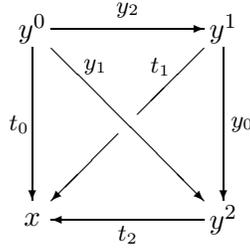
\begin{center}
\begin{picture}(1,1)
\put(0,1){\makebox(0,0){$y^0$}}
\put(1,1){\makebox(0,0){$y^1$}}
\put(1,0){\makebox(0,0){$y^2$}}
\put(0,0){\makebox(0,0){$x$}}
\put(.1,1){\vector(1,0){.8}}
\put(.9,0){\vector(-1,0){.8}}
\put(0,.9){\vector(0,-1){.8}}
\put(1,.9){\vector(0,-1){.8}}
\put(.1,.9){\vector(1,-1){.8}}
\put(.9,.9){\line(-1,-1){.35}}
\put(.45,.45){\vector(-1,-1){.35}}
\put(-.07,.5){\makebox(0,0){\footnotesize{$t_0$}}}
\put(.5,-.08){\makebox(0,0){\footnotesize{$t_2$}}}
\put(.67,.8){\makebox(0,0){\footnotesize{$t_1$}}}
\put(.5,1.1){\makebox(0,0){\footnotesize{$y_2$}}}
\put(1.1,.5){\makebox(0,0){\footnotesize{$y_0$}}}
\put(.33,.8){\makebox(0,0){\footnotesize{$y_1$}}}
\end{picture}
\end{center}
where $y_p = d_p(y)$.
\MS

Now if $C$ is a $(1,1)$-composer then the 1-simplices of $R^x$ are commutative triangles and the 2-simplices are commutative tetrahedra. Note that, with $z \in R^x_2$, the commutativity of $d_x^1(z)$ ($t2 \circ y1 \ISIT t_0$) follows routinely from the commutativity of the other faces of $z$. This is an instance of $\BOX{1}{2}{R^x} \ISO R^x_2$.
\MS

Therefore, in the case that $C$ is a $(1,1)$-composer, $R^x$ is the comma category $C \downarrow x$.

\subsection{$C^{x,F}$, a comma category generalization}

Suppose $C$ and $D$ are small categories i.e. $(1,1)$-composers and $F:D \to C$ is a functor. Given any $x \in \OB(C) = C_0$ then the comma category $x \downarrow F$ has as objects all pairs $(h,v) \in C_1 \times D_0$ such that $h : x \to F_0(v)$ is a map in $C$. A map in $x \downarrow F$ is $(h_0,v_0) \XRA{f} (h_1,v_1)$  where $f : v_0 \to v_1$ is in $D_1$ and $h_1 = F_1(f) h_0$. In the usual diagrams:
\setlength{\unitlength}{1in}
\begin{center}
\begin{picture}(4,.6)
\put(.3,.3){\makebox(0,0){Objects:}}
\put(1,.5){\MB{x}}
\put(1,0){\MB{F_0(v)}}
\put(1,.4){\vector(0,-1){.25}}
\put(1.1,.3){\MB{h}}
\put(2.5,.3){\makebox(0,0){Maps:}}
\put(3,0){\MB{F_0(v_0)}}
\put(3.5,.5){\MB{x}}
\put(4,0){\MB{F_0(v_1)}}
\put(3.4,.4){\vector(-1,-1){.25}}
\put(3.15,.34){\MBS{h_0}}
\put(3.6,.4){\vector(1,-1){.25}}
\put(3.85,.34){\MBS{h_1}}
\put(3.25,0){\vector(1,0){.5}}
\put(3.5,.1){\MBS{F_1(f)}}
\end{picture}
\end{center}
\MS

This can be generalized, as follows: replace the categories $C$ and $D$ with arbitrary simplicial sets, the functor $F$ with a simplicial map $F:D \to C$ and the object $x \in C_0$ with an element $x \in C_k$ for some given $k \geq 0$.

\begin{newdef}
{\bf ($C^{x,F}$)}
\index{$C^{x,F}$}
\index{C$^{x,F}$}
\label{C^x,F}
{\rm
\SL

Let $C$ and $D$ be simplicial sets, $F : D \to C$ a simplicial map, and $k \geq 0$. Fix $x \in C_k$. For each $j \geq 0$ define
\[
\begin{aligned}
C^{x,F}_j & \DFAS \SETT{
\; (z,v) \in C^x_j \times D_j : z \in \CPL(x;F_j(v))\; }\\
C^{x,F} & \DFAS \bigcup_{j \geq 0} C^x_j
\end{aligned}
\]
For each $j>0$ and each $p \in [j]$ 
\[
d^{x,F}_p : C^{x,F}_j \to C^{x,F}_{j-1}, \quad d^{x,F}_p(z,v) \DFAS (d^x_p(z),d_p(v))
\]
For each $j\geq 0$ and each $p \in [j]$ 
\[
s^{x,F}_p : C^{x,F}_j \to C^{x,F}_{j+1}, \quad s^{x,F}_p(z,v) \DFAS (s^x_p(z),s_p(v))
\]
}
\BX
\end{newdef}
It follows directly from the definitions that $C^{x,F}$ is a simplicial set.
\MS

Observe that given $m>1$ and $i \in [m+1]$ then
\[
\Bigl(
(z_0,v_0) \DDD{,} \underset{i}{-} \DDD{,} (z_{m+1},v_{m+1})
\Bigr) \in \BOX{i}{m+1}{C^{x,F}}
\]
implies
\[
(\OM{z}{m+1}{i}) \in \BOX{i}{m+1}{C^x} \; \text{and} \;
(\OM{v}{m+1}{i}) \in \BOX{i}{m+1}{D}
\]
This leads to:

\begin{newthm}
{\rm
\SL

Given $n \geq 1$ and $i \in [n+1]$, suppose $C$ and $D$ are \NICOMP{n}{i}s and $F : D \to C$ is a simplicial map. Then $C^{x,F}$ is also an \NICOMP{n}{i}.
}
\end{newthm}

\Proof

Suppose $m \geq n$ and 
\[
\Bigl(
(z_0,v_0) \DDD{,} \underset{i}{-} \DDD{,} (z_{m+1},v_{m+1})
\Bigr) \in \BOX{i}{m+1}{C^{x,F}}
\]
Then $(n,i)$-composition in $C^x$ and $D$ determine
\[
w = \COMP{m}{i}(\OM{z}{m+1}{i}) \in C^x_{m+1}
\]
and 
\[
y = \COMP{m}{i}(\OM{v}{m+1}{i}) \in D_{m+1}
\]
where, for all $p \neq i$, $d^x_p(w)= z_p$ and $d_p(y) = v_p$. Since $w$ and $y$ are determined uniquely by $(n,i)$-composition, to show that $\C^{x,F}_{m+1} \ISO \BOX{i}{m+1}{C^{x,F}}$, it suffices to show that $(w,y) \in C^{x,F}_{m+1}$; that is, $w \in \CPL(x;F_{m+1}(y))$.
\MS

Now suppose $\tilde{y} \in C_{m+1}$ is the complement of $x$ in $w$. For each $p \neq i$:
\[
\begin{aligned}
d^x_p(w) & \in \CPL(x;d_p(\tilde{y})) & \\
z_p & \in \CPL(x;F_m(v_p)) = \CPL(x;\, d_p(F_{m+1}(y)) & 
\end{aligned}
\]
It follows from $z_p = d^x_p(w)$ that $d_p(F_{m+1}(y)) = d_p(\tilde{y})$ for all $p \neq i$, and therefore that $\tilde{y}=F_{m+1}(y)$ and $(w,y) \in C^{x,F}_{m+1}$.

\qed

\subsection{Subcomplexes of $C^{x,F}$}

Given a simplicial map $F:D \to C$ we define two subcomplexes of $C^{x,F}$ derived from $L^{x,t}$ and $R^{x,t}$ above.

\begin{newdef}
\label{LxtF and RxtF}
\index{$L^{x,t,F}$}
\index{$R^{x,t,F}$}
{\rm
\SL

Let $F:D \to C$ be a simplicial map. Fix $k \geq 0$, $t \in [k]$ and $x \in C_k$. Define:
\[
\begin{aligned}
L^{x,t,F} & \DFAS \bigcup_{j \geq 0} \SETT{z \in L^{x,t}_j: \exists v \in D_j \;(z \in \CPL(x;F_j(v))\;} \\
R^{x,t,F} & \DFAS \bigcup_{j \geq 0} \SETT{z \in R^{x,t}_j: \exists v \in D_j \;(z \in \CPL(x;F_j(v))\;}
\end{aligned}
\]
}
\BX
\end{newdef}

The proof of Theorem \refpage{G^{x,t} is subcomplex} applies directly to the simplices of $L^{x,t,F}$ and $R^{x,t,F}$ and therefore $L^{x,t,F}$ and $R^{x,t,F}$ are subcomplexes of $C^{x,F}$. The proof that $C^{x,F}$ is an \NICOMP{n}{i} when $C$ and $D$ are \NICOMP{n}{i}s also applies. That is:

\begin{newcor}
{\rm
If $F:D \to C$ is a simplicial map of \NICOMP{n}{i}s, $k \geq 0$, $t \in [k]$ and $x \in C_k$ then $L^{x,t,F}$ and $R^{x,t,F}$ are also \NICOMP{n}{i}s.
}
\end{newcor}

\qed

\subsection{Formula relating complementary vertex functions}

This section will develop another description (Theorem \refpage{comp-vert-func-thm} below) of the relationship between complementary vertex functions.
\MS

Let
\[
\DMAP{f}{j}{k+1}, \qquad \DMAP{g}{k}{j+1}
\]
be complementary vertex functions.
Recall that any one of the functions $f, g, f\SH$ and $g\SH$ determines the other three.
\BS

\NI We will use the following notation and facts, consequences of $f$ and $g$ being non-decreasing:
\begin{enumerate}

\item 
For all $p \in [j+1]$, $g^{-1}(p) \SBS [k]$ is an interval, possibly empty.

For all $q \in [k+1]$, $f^{-1}(q) \SBS [j]$ is an interval, possibly empty.

\item
For all $p \in [j+1]$
\[
g\SH(g^{-1}(p)) = p + g^{-1}(p) = \SETT{p+q: g(q)=p}
\]
That is, $g\SH(g^{-1}(p))$ is the translation by $p$ of the interval $g^{-1}(p)$. 

Similarly, $f\SH(f^{-1}(q)) = q + f^{-1}(q)$.

\item
Some elementary observations:
\[
\begin{aligned}
p \neq p'\in [j+1] & \IMP g^{-1}(p) \cap g^{-1}(p') = \MT \\
 g\SH(g^{-1}(p)) \cap g\SH(g^{-1}(p')) &= \MT \quad\text{ (since $g\SH$ is 1-1)}\\
[k] & = \bigcup_{p \in [j+1]}g^{-1}(p) \quad \text{ (disjoint union)}\\
[j] & = \bigcup_{q \in [k+1]}f^{-1}(q) \quad \text{ (disjoint union)}
\end{aligned}
\]


\item
\label{A-B-notation}
Set $A \DFAS f\SH[j]$. Then
\[
A = f\SH \left( \bigcup_{q \in [k+1]}f^{-1}(q) \right) = \bigcup_{q \in [k+1]} f\SH(f^{-1}(q))
\]
is a disjoint union. Similarly
\[
B \DFAS g\SH[k] = \bigcup_{p \in [j+1]} g\SH(g^{-1}(p))
\]
is also a disjoint union.

\end{enumerate}

\begin{newdef}
\index{histogram} \index{h$^f$} \index{$h^f$}
{\rm
\SL

Given $f:[j] \to [k+1]$ non-decreasing, we define the {\bf histogram of $f$} to be $h^f: [k+1] \to [j+1]$ where
\[
h^f(q) \DFAS \left| f^{-1}(q) \right| = \left| \SETT{p \in [j]: f(p)=q} \right|
\]
(Note that $h^f$ is generally not monotonic).
}

\BX
\end{newdef}
\MS

Since $[j] = \bigcup_{q \in [k+1]} f^{-1}(q)$ then
\[
\Bigl| \, [j] \, \Bigr| = j+1 =  \sum_{q=0}^{k+1} h^f(q)
\]

\NI {\bf Notations:}
\index{$\int_0^q h^f$}
Given $q \in [k+1]$, we will denote $\sum_{t=0}^q h^f(t)$ by $\int_0^q h^f$, thinking of $h^f$ as a step function. Note that with this notation, $\int_0^0 h^f >0$ if $f^{-1}(0) \neq \MT$.
\MS

Given a pair of complementary vertex functions $\DMAP{f}{j}{k+1}$ and $\DMAP{g}{k}{j+1}$, define
\[
\hat{f}: [j] \to [k+1], \qquad \hat{f}(p) \DFAS \int_0^p h^g
= \sum_{t=0}^p h^g(t) = \sum_{t=0}^p |g^{-1}(t)|
\]
\[
\hat{g}:[k] \to [j+1], \qquad \hat{g}(q) \DFAS \int_0^q h^f
= \sum_{t=0}^q h^f(t) = \sum_{t=0}^q |f^{-1}(t)|
\]
\BX

Note that $\hat{f}$ is defined referring only to $g$, and $\hat{g}$ is defined referring only to $f$.
\MS

\begin{newthm}
\label{comp-vert-func-thm}
{\rm
\SL

Given any pair of complementary vertex functions $\DMAP{f}{j}{k+1}$ and $\DMAP{g}{k}{j+1}$, then $\hat{f} = f$ and $\hat{g}=g$.
}
\end{newthm}

\Proof

It suffices to show that $\hat{f} = f$.

First: since $g^{-1}(p)$ is an interval for each $p \in [j+1]$, $[k] = \bigcup_p g^{-1}(p)$ and $| g^{-1}(p)| = h^g(p)$ then:
\[
\begin{aligned}
g^{-1}(0) = & \Bigl[0,-1 + h^g(0) \Bigr] = \Bigl[0, -1+ \int_0^0 h^g \Bigr] = \Bigl[0,-1+\hat{f}(0) \Bigr]\\
g^{-1}(1) = & \Bigl[ \int_0^0 h^g,\, -1 + \int_0^1 h^g \Bigr] 
= \Bigl[\hat{f}(0), -1+\hat{f}(1) \Bigr]
\\
g^{-1}(2) = & \Bigl[ \int_0^1 h^g,\, -1 + \int_0^2 h^g \Bigr] =
\Bigl[ \hat{f}(1), -1 + \hat{f}(2) \Bigr] \\
\vdots & \\
g^{-1}(p) = & \Bigl[ \int_0^{p-1} h^g, \, -1 + \int_0^p h^g \Bigr] =
\Bigl[ \hat{f}(p-1), -1 + \hat{f}(p) \Bigr]
\end{aligned}
\]
and
\begin{equation} \label{B-interval-eqn}
g\SH(g^{-1}(p)) = \Bigl[ p+\hat{f}(p-1),\, p-1 + \hat{f}(p) \Bigr] \quad \text{for } p>0
\end{equation}
In the special case of $p=0$, $g\SH(g^{-1}(0)) = [0, -1+\hat{f}(0)]$. Note that $g\SH(g^{-1}(0)) \neq \MT$ when $\hat{f}(0)>0$.
\MS

\NI
Let $A \DFAS f\SH([j]), B \DFAS g\SH([k])$, as above.
Since $B = \bigcup_{p \in [j+1]} g\SH(g^{-1}(p))$ and $A$ is disjoint from $B$ then $A$ is the disjoint union of intervals between $g\SH(g^{-1}(p))$ and $g\SH(g^{-1}(p+1))$ for $0 \leq p \leq j$.
\BS

Next, we will examine those intervals which comprise $A$, by cases. There are three possibilities: (i) an interval of $A$ includes 0, which occurs when $g^{-1}(0) \neq \MT$ ($ \iff h^g(0)>0)$; (ii) an interval of $A$ is between two intervals of $B$; (iii) an interval of $A$ includes $j+k+1$.
\MS

\NI {\bf Case: $h^g(0)>0$}

Let $p'>0$ be the least such that $h^g(p')>0$. That is, $h^g(t)=0$ for any and all $t$ where $0<t<p'$. 

The interval of $A$ between $g\SH(g^{-1}(0)) = [0,-1+\hat{f}(0)]$ and
\[
g\SH(g^{-1}(p')) = \Bigl[p'+\hat{f}(p'-1), p'-1 + \hat{f}(p') \Bigr]
\]
is
\[
\Bigl[ \hat{f}(0), \, p'-1 + \hat{f}(p'-1) \Bigr]
\]
Since $h^g(t)=0$ when $0<t<p'$ then $\hat{f}(t) = \hat{f}(0)$ for all such $t$ and
\[
\Bigl[ \hat{f}(0), \, p'-1 + \hat{f}(p'-1)] = \Bigl[ \hat{f}(0), \, p'-1 + \hat{f}(0)] = \hat{f}(0) + \Bigl[ 0,p'-1
\Bigr]
\]
Thus in the case when $h^g(0)>0$, $A$ contains the interval consisting of all $\hat{f}(0) + t = \hat{f}(t) + t$ for each $t \in [0,p'-1]$.
\BS

\NI {\bf Case:} $0<p<p'$ such that $h^g(p)>0, h^g(p')>0$ and $h^g(t)=0$ for any and all $t$ such that $p<t<p'$.
\MS

In this case, $A$ contains the interval between $g\SH(g^{-1}(p))$ and $g\SH(g^{-1}(p'))$. Using equation \eqref{B-interval-eqn}, that interval is
\[
\Bigl[
p+\hat{f}(p), p'-1+\hat{f}(p'-1)
\Bigr] = 
\Bigl[
p+\hat{f}(p), p'-1+\hat{f}(p)
\Bigr] = \hat{f}(p) + [p,p'-1]
\]
using, again, that $h^g(t)=0$ for each $t \in [p+1,p'-1]$.
Since $\hat{f}(t) = \hat{f}(p)$ for each $t \in [p,p'-1]$ it follows that each element of $A$ in this interval has the form $\hat{f}(t) + t$, $t \in [p,p'-1]$.
\BS

\NI {\bf Case:} $h^g(p)>0$ and $h^g(t)=0$ for any and all $t \in [p+1,j+1]$.
\MS

In this case, equation \eqref{B-interval-eqn} implies that $A$ contains the interval
\[
\Bigl[
p+\hat{f}(p),j+k+1
\Bigr]= 
\hat{f}(p) + 
\Bigl[
p,j+k+1-\hat{f}(p)
\Bigr]
\]
Since $h^g(t)=0$ whenever $t \in [p+1,j+1]$, then for each $t>p$, $t + \hat{f}(p) = t + \hat{f}(t)$. Therefore, each element in this interval of $A$ has the form $t+\hat{f}(t)$.
\MS

\NI Collectively, these three cases show that
\[
A \DFAS \SETT{t+f(t):t \in [j]} = \SETT{t + \hat{f}(t) : t \in [j]}
\]
and therefore $\hat{f} = f$.
\MS

An identical analysis proves that $\hat{g} = g$.

\qed
\MS

\begin{example}
{\rm
Let $j=6, k=10$ and consider the complementary vertex functions $f:[6] \to [11]$ and $g:[10] \to [7]$ where
\[
\begin{aligned}
f & = (1,1,2,6,6,6,9) : [6] \to [11]\\
f\SH & = (1,2,4,9,10,11,15):[6] \to [17]\\
g\SH & = (0,3,5,6,7,8,12,13,14,16,17) : [10] \to [17]\\
g & = (0,2,3,3,3,3,6,6,6,7,7):[10] \to [7]
\end{aligned}
\]
Then
\[
\begin{array}{|c|c|c|c|c|c|c|c||c|}
\hline
t & 0 & 1 & 2 & 3 & 4 & 5 & 6 & 7 \\
\hline
g^{-1}(t) & [0,0] & \MT & [1,1] & [2,5] & \MT & \MT & [6,8] & [9,10]\\
\hline
h^g(t) = |g^{-1}(t)| & 1 & 0 & 1 & 4 & 0 & 0 & 3 & 2\\
\hline
\hat{f}(t) & 1 & 1 & 2 & 6 & 6 & 6 & 9 & -\\
\hline
\end{array}
\]
}
\end{example}

\section{Represented subcomplexes}

\subsection{Vertex-span subcomplexes of $\MC{S}$}

In this section, $C$ denotes an arbitrary simplicial set and $\MC{S}$ denotes the simplicial set of relations. 
\MS

We may associate with $C$ certain simplices of $\MC{S}$ by defining, for each $k \geq 0$, relations $R \SBS (C_0)^{k+1}$ and forming $y^R \in \MC{S}_k$, the minimal \SBSI\ $k$-simplex. The prequisite information is in section \refpage{min-simplex-of-rel}, definition \refpage{minimal-simplex}.

This construction yields a subcomplex $\MC{W}$ of $\MC{S}$ based on $C$ we will call the ``vertex-span'' subcomplex of $\MC{S}$ determined by $C$.
\MS

\begin{newdef}
{\bf (Sponsor)}
\label{sponsor-def}
\index{sponsor}
{\rm
\SL

Let $C$ be a simplicial set, $k \geq 0$. Then $z \in C_k$ will be said to {\bf sponsor} 
\[
\LIST{t}{k} \in \left( C_0 \right)^{k+1}
\]
if for each $p \in [k]$, $\VERT{p}(z) = t_p$. Note that $0$-simplices are self-sponsored.
\MS

\NI For each $k \geq 0$, define $\SPON{C}{k} \SBS \left( C_0 \right)^{k+1}$ by:
\[
\SPON{C}{k} \DFAS \SETT{\LIST{t}{k} \in \left( C_0 \right)^{k+1}: \LIST{t}{k} \text{ is sponsored.}}
\]

\BX
}
\end{newdef}

It follows immediately that if $z \in C_k$ sponsors $\LIST{t}{k}$ then for all $p \in [k]$, $d_p(z)$ sponsors $(\DLI{t}{k}{p})$ (when $k>0$) and, for all $k \geq 0$, $s_p(z)$ sponsors $(\DLIST{t}{k}{p})$.

Of course, a particular $\LIST{t}{k}$, may have no sponsors or more than one sponsor.
\MS

Given any non-empty $R \SBS \left( C_0 \right)^{k+1}$ such that $R \SBS \SPON{C}{k}$, let $y^R$ denote the minimal \SBSI\ $k$-simplex determined by $R$. That is:

\[
h_k \LIST{t}{k} \in y^R \iff \LIST{t}{k} \in R
\]
and for each subface-permissible $p_1 \DDD{<} p_m$,
\[
d_{p_1} \cdots d_{p_m}(y^R) = 
\SETT{e_{p_1} \cdots e_{p_m}h_k\LIST{t}{k} : \LIST{t}{k} \in R}
\]

From this, the simplicial set $C$ determines a subcomplex of $\MC{S}$, as follows.

\begin{newdef}
{\bf (Vertex-span subcomplex determined by $C$)}
\index{vertex-span subcomplex}
\label{vertex-span subcomplex}
{\rm
\SL

Let $C$ be an arbitrary simplicial set.
For each $k \geq 0$ define $\MC{W}_k \SBS \MC{S}_k$ as follows.
\MS

$\MC{W}_0 \DFAS $ the set of all unary relations $R \SBS \SPON{C}{0} = C_0$.
\MS

For $k>0$,
\[
\MC{W}_k \DFAS \SETT{y^R : R \SBS \SPON{C}{k}, \, R \neq \MT}
\]
}

\BX
\end{newdef}

\setlength{\unitlength}{1in}

\begin{example}
{\rm
To illustrate $\MC{W}_2 \SBS \MC{S}_2$: we have $(t_0,t_1,t_2) \in \SPON{C}{2}$ iff there exists $z \in C_2$ with

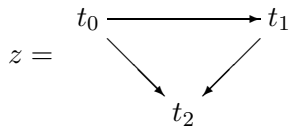
\begin{center}
\begin{picture}(1,.6)
\put(0,.5){\MB{t_0}}
\put(1,.5){\MB{t_1}}
\put(.5,0){\MB{t_2}}
\put(0.1,.4){\vector(1,-1){.3}}
\put(0.1,.5){\vector(1,0){.8}}
\put(0.9,.4){\vector(-1,-1){.3}}
\put(-.3,.3){\MB{z =}}
\end{picture}
\end{center}
Let $R \SBS \SPON{C}{2}$ be any non-empty subset. Then the element set of the 2-simplex $y^R \in \MC{W}_2$ is $\SETT{h_2(t_0,t_1,t_2) : (t_0,t_1,t_2) \in R}$ and the monic
\[
y^R  \to d_0(y^R) \times d_1(y^R) \times d_2(y^R) 
\]
is
\[
h_2(t_0,t_1,t_2) \mapsto \Bigl(h_1(t_1,t_2),\;h_1(t_0,t_2),\;h_1(t_0,t_1)\Bigr)
\]
Note that, for example,
\[
d_0(y^R) = \SETT{h_1(t_1,t_2) : \exists (t_0,t_1,t_2) \in R} = y^{d_0(R)}
\]
and the vertex-relation of $d_0(y^R)$ is
\[
\SETT{(t_1,t_2) \in C_0 \times C_0: \exists (t_0,t_1,t_2) \in R}
\]
with $d_0(y^R) \in \MC{W}_1$ since each $(t_1,t_2)$ in the vertex-relation of $d_0(y^R)$ is sponsored by $d_0(z)$ where $z$ is any sponsor of $(t_0,t_1,t_2) \in R$.
}
\end{example}
\MS

Following directly from these definitions:
\begin{newlem}
{\rm
\SL 

$\SETT{\MC{W}_k:k \geq 0}$ is a subcomplex of $\MC{S}$.}
\end{newlem}

\qed

\subsection{Subcomplex $\MC{R}(x,u)$ of representables determined by $C^x$}

If $x,y$ are $0$-simplices in a simplicial set $C$ then there are two distinguished unary relations on $C_1$ namely:
\[
\begin{aligned}
\SETT{z \in C_1 : d_0(z)=x,\; d_1(z)=y} \SBS C_1\\
\SETT{z \in C_1 : d_0(z)=y,\; d_1(z)=x} \SBS C_1\\
\end{aligned}
\]
If $C$ is the nerve of a small category then these two unary relations are, respectively, the hom sets $C(y,x) \SBS C_1$ and $C(x,y) \SBS C_1$. Given $x \in C_0$, these unary relations determine the hom functors $C(x,-):C \to \STS$ and $C(-,x):C^\text{op} \to \STS$.
\MS

This section generalizes these ideas, first to arbitrary simplicial sets $C$ where $x \in C_k$, $y \in C_j$ and where ``$z \in C(x,y)$'' is replaced by 
\[
z \in \CPL(x,\mu;y,\del_u \LAM)
\]
where $u \in [k+1]$ and
\[
[k] \XRA{\mu} [j+1]
\]
is the complementary vertex function of
\[
[j] \XRA{\LAM} [k] \XRA{\del_u} [k+1]
\]
The analog of hom-functor arises from a fixed $x \in C_k$ and $u \in [k+1]$ to yield a simplicial map $C \times \Delta[k] \to \MC{S}$. The choice of $u$ generalizes the notions of covariant and contravariant hom-functors.
\MS

\begin{example}
{\rm
To motivate this generalization, we first consider it when $k=0$ and for a fixed $x \in C_0$.

Fix $u \in [1]$. Then $x$ and $u$ determine a simplicial map $C \times \Delta[0] \to \MC{S}$ which we illustrate in dimensions 0 and 1 as follows. 
\MS

\NI {\bf Dimension 0: $C_0 \times \Delta[0]_0 \to \MC{S}_0$} 
\MS

Let $(y,\LAM) \in C_0 \times \Delta[0]_0$ and $u \in [1]$. Since $\Delta[0]_0 = \{ 1_{[0]} \}$, then the only choice of $\LAM$ is $1_{[0]}$. We get a 0-simplex of $\MC{S}$ using that $\del_u \LAM = \del_u$
and 
\[
\del_u \LAM = \del_u = \left\{
\begin{array}{cc}
 (1) & \text{if } u=0\\
 (0) & \text{if } u=1
\end{array}
\right.
\]
Then $x, y$ and $u$ determine the 0-simplex (unary relation) of $\MC{S}$
\[
\CPL(x,\mu;y,\del_u \LAM) \SBS C^x_0
\]
where, if $z \in \CPL(x,\mu;y,\del_u \LAM) \SBS C^x_0$ then $z \in C_1$ and
\[
\VL_z(y) = (\del_u \LAM)\SH[0] = \left\{
\begin{array}{cc}
\{1\} & \text{if } u=0\\
\{0\} & \text{if } u=1
\end{array}
\right.
\]
That is, 
\[
z \text{ is } \left\{
\begin{array}{ll}
 x \to y & \text{if } u=0\\
 y \to x & \text{if } u=1
\end{array}
\right.
\]
\MS

\NI {\bf Dimension 1: $C_1 \times \Delta[0]_1 \to \MC{S}_1$}
\MS

We associate to $(y,\LAM) \in C_1 \times \Delta[0]_1$ a 1-simplex of $\MC{S}$
\[
w \SBS d_0(w) \times d_1(w)
\]
where $d_0(w)$ and $d_1(w)$ are unary relations as described just above.

\NI and
\[
(z_0,z_1) \in w \overset{\text{def}}{\iff}
\exists z \in \CPL(x,\mu;y,\del_u \LAM)
\]
such that 
\[
\begin{aligned}
z_0 &= d^x_0(z) \in \CPL(x,\SIG_0 \mu;d_0(y),\del_u \LAM \del_0)\\
z_1 &= d^x_1(z) \in \CPL(x,\SIG_1 \mu;d_1(y),\del_u \LAM \del_1)
\end{aligned}
\]
$\del_u \LAM = [1] \XRA{\LAM} [0] \XRA{\del_u} [1]$ and therefore
\[
\del_u \LAM = \left\{
\begin{array}{ll}
(1,1) & \text{if } u=0\\
(0,0) & \text{if } u=1
\end{array}
\right.
\]
\[
(\del_u \LAM)\SH = \left\{
\begin{array}{ll}
(1,2) & \text{if } u=0\\
(0,1) & \text{if } u=1
\end{array}
\right.
\]
Also, $\del_u \LAM \del_p = [0] \XRA{\del_p} [1] \XRA{\LAM} [0] \XRA{\del_u} [1] =\del_u$.
\MS

\HD{Case $u=0$:}
\MS

Then $\del_0 \LAM = (1,1)$, $(\del_u \LAM)\SH = (1,2)$ and therefore $\VL_z(y) = \{1,2\}$.
\[
\begin{aligned}
d^x_0(z) &= d_{0+\del_0(0)}(z) = d_1(z) = x \XRA{d_1(z)} d_0(y) \\
d^x_1(z) &= d_{1+\del_0(1)}(z) = d_2(z) = x \XRA{d_2(z)} d_1(y)
\end{aligned}
\]
and $z$ is
\begin{center}
\begin{picture}(1,.6)
\put(0,.5){\MB{x}}
\put(1,.5){\MB{d_1(y)}}
\put(.5,0){\MB{d_0(y)}}
\put(.1,.4){\vector(1,-1){.25}}
\put(.9,.4){\vector(-1,-1){.25}}
\put(.1,.5){\vector(1,0){.7}}
\put(.85,.25){\MBS{y}}
\end{picture}
\end{center}
\MS

\HD{Case $u=1$:}
\MS

$\del_1 \LAM = (0,0)$, $(\del_1 \LAM)\SH = (0,1)$ and therefore $\VL_z(y) = \{0,1\}$.
\[
\begin{aligned}
d^x_0(z) &= d_{0+\del_1(0)}(z) = d_0(z) = d_0(y) \XRA{d_0(z)} x \\
d^x_1(z) &= d_{0+\del_1(1)}(z) = d_1(z) = d_1(y) \XRA{d_1(z)} x
\end{aligned}
\]
and $z$ is
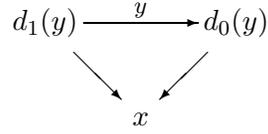
\begin{center}
\begin{picture}(1,.6)
\put(0,.5){\MB{d_1(y)}}
\put(1,.5){\MB{d_0(y)}}
\put(.5,0){\MB{x}}
\put(.15,.35){\vector(1,-1){.25}}
\put(.85,.35){\vector(-1,-1){.25}}
\put(.2,.5){\vector(1,0){.6}}
\put(.5,.57){\MBS{y}}
\end{picture}
\end{center}
\SL
}
\end{example}
\BS

\NI{\bf The general construction}
\MS

Suppose $C$ is an arbitrary simplicial set, with $k \geq 0$, $x \in C_k$  and $u \in [k+1]$ held fixed. The goal is to define a simplicial map $C \times \Delta[k] \to \MC{S}$ determined by $x$ and $u$. We will use the vertex-span subcomplex idea from the previous section applied to the simplicial set $C^x$.
\MS

\index{$\MC{W}^x$}
\index{W$^x$}
\NI {\bf Notation:} We will denote by $\MC{W}^x$ the full vertex-span subcomplex of $\MC{S}$ arising from $C^x$.

\BX
\MS

By definition, $\MC{W}^x_0$ consists of all unary relations $w \SBS C^x_0$. Each element $t \in w$ occurs as $t \in \CPL(x,\mu;y,\LAM) \SBS C^x_0$ for some $y \in C_0$ and complementary vertex functions $[k] \XRA{\mu} [1]$ and $[0] \XRA{\LAM} [k+1]$. $t \in C_{k+1}$.
\MS

For $j>0$, $\MC{W}^x_j$ consists of \SBSI\ relations $w \in \MC{S}_k$ whose fundamental signature is $\FSIG(C^x_0 \DDD{,} C^x_0)$. That is, given
\[
\LIST{t}{j} \in C^x_0 \DDD{\times} C^x_0
\]
then $h_j \LIST{t}{j} \in w$ iff there exists $z \in C^x_j$ which sponsors $\LIST{t}{j}$ i.e. for each $p \in [j]$, $\VERT{p}^x(z)=t_p$.
\BS

We will now define a particular $j$-simplex of $\MC{W}^x$, denoted $\HMU xuj{y,\LAM}$ corresponding to a choice of $(y,\LAM) \in C_j \times \Delta[k]_j$ and $u \in [k+1]$. We also define a certain subcomplex of $\MC{S}$, denoted $\MC{R}(x,u)$.

\begin{newdef}
{\bf ($\HMU xuj{y,\LAM}$)}
\label{rep-subcomplex}
\index{$\MC{R}(x,u)$}
\index{$H^{x,u}_j(y,\LAM)$}
\index{represented subcomplex}
\index{represented simplex}
{\rm
\SL

Given $(y,\LAM) \in C_j \times \Delta[k]_j$, we define a $j$-simplex in $\MC{W}^x$
\[
\HMU xuj{y,\LAM} \SBS \prod_{p=0}^j \HMU xu{j-1}{d_p(y),\LAM \del_p}
\]
and define 
\[
\MC{R}(x,u)_j \DFAS \bigcup_{(y,\LAM)} \HMU xuj{y,\LAM}
\qquad (y,\LAM) \in C_j \times \Delta[k]_j
\]
\MS

\HD{Dimension $j=0$:} 
\MS

Given any $(y,\LAM) \in C_0 \times \Delta[k]_0$, we define the 0-simplex
\[
\HMU xu0{y,\LAM} \SBS C^x_0
\]
by 
\[
\HMU xu0{y,\LAM}\DFAS \CPL(x,\mu;y,\del_{u} \LAM)
\]
where $[k] \XRA{\mu} [1]$ is the vertex function complementary to $[0] \XRA{\del_u \LAM}[k+1]$. Note that $z \in \HMU xu0{y,\LAM}$ implies $\VL_z(y) = \{ (\del_u \LAM)\SH(0) \}= \{ (\del_u \LAM)(0) \}$; that is, the 0=simplex $y$ is $\VERT{(\del_u\LAM)(0)}(z)$.
\MS

\NI We set $\MC{R}(x,u)_0 \DFAS \Bigl\{\Bigl(\HMU xu0{y,\LAM} \SBS C^x_0 \Bigr) : (y,\LAM) \in C_0 \times \Delta[k]_0\Bigr\}$.
\MS

\HD{Dimension $j>0$:}
\MS

Let $(y,\LAM) \in C_j \times \Delta[k]_j$. The fundamental entries of the proposed $\HMU xuj{y,\LAM}$ will be of the form
\[
t_p \in \HMU xu0{\VERT{p}(y), \VERT{p}(\LAM)} = \CPL(x,\mu_p;\VERT{p}(y),\del_u \VERT{p}(\LAM)), \quad p \in [j]
\]
where $\del_u \VERT{p}(\LAM) = [0] \XRA{\VERT{p}(\LAM)} [k] \XRA{\del_u} [k+1]$ and $\DMAP{\mu_p}{k}{1}$ is the complementary vertex function. 
\MS

Note that $(\VERT{p}(\del_u\LAM))(0) = \del_u(\LAM(p)) = (\del_u \VERT{p}(\LAM))(0)$. That is, $\VERT{p}(\del_u\LAM) = \del_u \, \VERT{p}(\LAM)$.
\MS

The entries of $\HMU xuj{y,\LAM}$ will be of the form $h_n\LIST{t}{j}$. The definition splits into cases: $(y,\LAM)$ is non-degenerate and $(y,\LAM)$ is degenerate.
\MS

Non-degenerate case: $h_n\LIST{t}{j} \in \HMU xuj{y,\LAM}$ iff there exists $z \in \CPL(x,\mu;y,\del_u \LAM) \SBS C^x_j$ where $\mu$ is complementary to $\del_u \LAM$, such that for all $p \in [j]$, 
\[
t_p = \VERT{p}^x(z) \in \CPL(x,\mu^{(p)}; \VERT{p}(y), \VERT{p}(\del_u \LAM))
\]
That is, $z$ sponsors $\LIST{t}{j}$.
\BS

Degenerate case: When $(y,\LAM)$ is degenerate
\[
(y,\LAM) = s_{m_1} \DDD{} s_{m_r}(y',\LAM')\qquad \text{for some }r \leq j
\]
and $(y',\LAM') \in C_{j-r} \times \Delta[k]_{j-r}$ is non-degenerate
then we define $h_j \LIST{t}{j} \in \HMU xuj{y,\LAM}$ iff there exists $h_{j-r} \LIST{t'}{{j-r}} \in \HMU xu{j-r}{y', \LAM'}$ such that
\[
h_j \LIST{t}{j} = c_{m_1} \DDD{} c_{m_r} h_{j-r} \LIST{t'}{{j-r}}
\]
Note: In terms of sponsors, this says that $\LIST{t'}{{j-r}}$ is sponsored by some $z' \in \CPL(x,\mu';y',\del_u \LAM')$. It then follows that $s^x_{m_1} \DDD{} s^x_{m_r}(z')$ sponsors $\LIST{t}{j}$.

\BX
}
\end{newdef}

\HD{Notes concerning the definition of $\HMU xuj{y,\LAM}$}
\begin{enumerate}
\item 
If $z$ sponsors an element of $\HMU xuj{y,\LAM}$ then for each $p \in [j]$, $d^x_p(z)$ sponsors an element of $\HMU xu{j-1}{d_p(y),\LAM \del_p}$ and $s^x_p(z)$ sponsors an element of $\HMU xu{j+1}{s_p(y),\LAM \SIG_p}$.

\item
Given $(y,\LAM) \in C_j \times \Delta[k]_j$ and $z \in \CPL(x,\mu;y,\del_u \LAM)$ then $(y,\LAM)$ is determined by $z$ because $y$ is the complementary subface of $x$ in $z$, $\del_u \LAM$ is the complementary vertex function of $\mu$ and $\LAM = \sigma_u \del_u \LAM$. Therefore $z$ sponsors the element $h_n\LIST{t}{j} \in \HMU xuj{y,\LAM}$ where for each $p \in [j]$, $t_p = \VERT{p}^x(z)$.

\item 
Let $j, k \geq 0$, $\DMAP{\LAM}{j}{k+1}$ and $u \in [k+1]$. We say ``$\LAM$ omits $u$'' if $u \notin \LAM[j]$. Then
\[
\begin{aligned}
\LAM \text{ omits } u & \iff \exists ! \DMAP{\gamma}{j}{k} \; (\LAM = \del_u \gamma) \\
\LAM \text{ omits } u & \IMP \ALL p \in [j]\; (\LAM \del_p \text{ omits } u) \\
\LAM \text{ omits } u & \IMP \ALL p \in [j]\; (\LAM \sigma_p \text{ omits } u)
\end{aligned}
\]
If $\LAM$ omits $u$ we denote the unique factorization through $\del_u$ by $\LAM = \del_u \LAM^u$. Note that $\LAM^u = \sigma_u \LAM$ in this case.
\MS

Given $x \in C_k$ and $u \in [k+1]$ there is a subcomplex $C^{x,u} \SBS C^x$ defined for each $j \geq 0$ by
\[
z \in C^{x,u}_j \iff z \in \CPL(x,\mu;y,\LAM)
\]
where $y$ is the complementary subface of $x$ in $z$ and the vertex function $\LAM$ of $y$ omits $u$. It follows that if $z \in C^{x,u}_j$ then $z$ sponsors an element of $\HMU 
xuj{y,\sigma_u \LAM}$.

\item 
Given
\[
\Bigl(
(y_0,\LAM_0) \DDD{,} (y_j,\LAM_j)
\Bigr) \in \Delta^\bullet(j)(C \times \Delta[k])
\]
so that $\LIST{y}{j} \in \Delta^\bullet(j)(C)$ and $\LIST{\LAM}{j} \in \Delta^\bullet(j)(\Delta[k])$, then there exists a unique $\LAM : [j] \to [k]$ such that for all $p \in [j]$, $\LAM_p = \LAM \del_p$. There may or may not exist $y \in C_j$ such that for all $p \in [j]$, $d_p(y) = y_p$.

If such a $y \in C_j$ does exist then, given any $u \in [k+1]$, there is a corresponding $j$-simplex of $\MC{R}(x,u)$, namely $\HMU xuj{y,\LAM} \SBS \prod_{p=0}^j \HMU xu{j-1}{y_p,\LAM_p}$.

There may be more than one such $y$, in general.

If $C$ is an \NICOMP{n}{i} $j \geq n+1$ and 
\[
\Bigl\{
(y_p,\LAM_p):p \in [j]-\{i\}
\Bigr\} \in \BOX{i}{j}{C \times \Delta[k]}
\]
then there is a unique $y$, determined by $(n,i)$-composition, and therefore just one $j$-simplex $\HMU{x}{u}{j}{y,\LAM}$ of $\MC{R}(x,u)_j$ whose faces are $\HMU xu{j-1}{y_p,\LAM_p}$, $p \in [j]-\{i\}$.

\end{enumerate}
\BX
\MS

\begin{example}
{\rm
The rather dense notation for elements of $\HMU xuj{y,\LAM}$ describes a relatively simple picture. Here is a low-dimensional example to show what that looks like.
\MS

Suppose $k=2$, $x \in C_2$, $j=3$, $y \in C_3$, $u=1 \in [3]=[k+1]$ and $(y,\LAM) \in C_3 \times \Delta[2]_3$
where $\LAM = (0,1,1,2) : [3] \to [2]$.

An element of $\HMU xu3{y,\LAM}$ is sponsored by some 
\[
z \in \CPL(x,\mu;y,\del_1 \LAM) \SBS C^x_3 \SBS C_6
\]
where $[2] \XRA{\mu} [4]$ is complementary to $[3] \XRA{\LAM} [2] \XRA{\del_1} [3]$. Direct calculations show:
\[
\begin{aligned}
\del_1 \LAM &= (0,2,2,3): [3] \XRA{\LAM} [2] \XRA{\del_1} [3]\\
(\del_1 \LAM)\SH &= (0,3,4,6):[3] \to [6]\\
\mu\SH &= (1,2,5) : [2] \to [6] \text{ since } \SETT{1,2,5} = [6]-\SETT{0,3,4,6}\\
\mu &= (1,1,3):[2] \to [4]
\end{aligned}
\]

So $x$ and $y$ are complementary in $z$ with $\VL_{z}(x) = \SETT{1,2,5}$ and $\VL_z(y) = \SETT{0,3,4,6}$. The picture showing just the ``outline'' (solid) edges of $x$ and $y$ as subfaces of $z$ is:
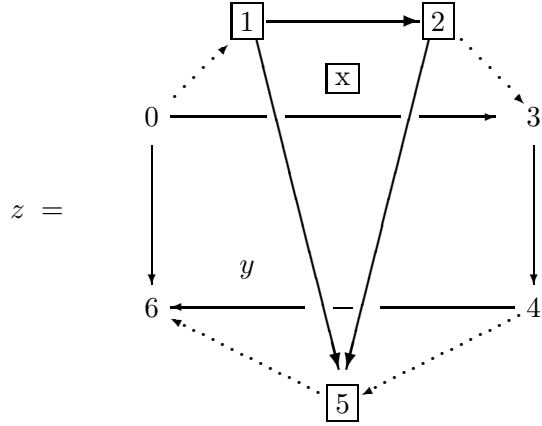
\begin{center}
\begin{picture}(2,2.2)
\put(-.6,1){\MB{z \; =}}
\put(0,1.5){\MB{0}}
\put(.5,2){\MB{\fbox{1}}}
\put(1.5,2){\MB{\fbox{2}}}
\put(2,1.5){\MB{3}}
\put(2,.5){\MB{4}}
\put(1,0){\MB{\fbox{5}}}
\put(0,.5){\MB{6}}
\thicklines
\put(.6,2){\vector(1,0){.8}} 
\put(1.45,1.9){\vector(-1,-4){.43}} 
\put(.55,1.9){\vector(1,-4){.43}} 
\thinlines
\put(.1,1.5){\line(1,0){.5}}
\put(1.4,1.5){\vector(1,0){.4}}
\put(2,1.35){\vector(0,-1){.73}}
\put(1.9,.5){\line(-1,0){.7}}
\put(.8,.5){\vector(-1,0){.7}}
\put(0,1.35){\vector(0,-1){.73}}
\put(.7,1.5){\line(1,0){.6}}
\put(1,.5){\MB{\line(1,0){.1}}}
\put(1,1.7){\MB{\fbox{x}}}
\put(.5,.7){\MB{y}}
\multiput(.1,1.6)(.05,.05){6}{.}
\put(.4,1.88){\vector(1,1){0}}
\multiput(1.6,1.9)(.05,-.05){7}{.}
\put(1.95,1.57){\vector(1,-1){0}}
\multiput(1.9,.45)(-.05,-.025){16}{.}
\put(1.1,.045){\vector(-2,-1){0}}
\multiput(.85,.05)(-.05,.025){16}{.}
\put(.1,.4425){\vector(-2,1){0}}
\end{picture}
\end{center}
\MS 

$z$ sponsors $h_3(t_0,t_1,t_2,t_3) \in \HMU x13{y,\LAM}$ where, for each $p \in [3]$, $t_p = \VERT{p}^x(z)$ is a fundamental entry of $h_3(t_0,t_1,t_2,t_3)$. 
\[
\begin{array}{l|c|c}
p & t_p & \VL_z(t_p)\\
\hline
0 & d^x_1 d^x_2 d^x_3(z) = d_3 d_4 d_6(z) & \SETT{0,1,2,5}\\
1 & d^x_0 d^x_2 d^x_3(z) = d_0 d_4 d_6(z) & \SETT{1,2,3,5}\\
2 & d^x_0 d^x_1 d^x_3(z) = d_0 d_3 d_6(z) & \SETT{1,2,4,5}\\
3 & d^x_0 d^x_1 d^x_2(z) = d_0 d_3 d_4(z) & \SETT{1,2,5,6}\\
\end{array}
\]

Components of $h_3(t_0,t_1,t_2,t_3)$ are sponsored by (sub)faces of $z$. For example, 
$
e_0(h_3(t_0,t_1,t_2,t_3)) = h_2(t_1,t_2,t_3) \in H^{x,1}_2(d_0(y),\LAM \del_0)
$
is sponsored by 
\[
d^x_0(z) = d_0(z) \in \CPL(x,\SIG_0 \mu;d_0(y), \del_1 \LAM \del_0)
\]
Here, 
\[
\begin{aligned}
\del_1 \LAM \del_0 & = (2,2,3): [2] \XRA{\del_0} [3] \XRA{\LAM} [2] \XRA{\del_1} [3]\\
(\del_1 \LAM \del_0)\SH &= (2,3,5):[2] \to [5]\\
(\SIG_0 \mu)\SH &= (0,1,4):[2] \to [5]\\
\SIG_0 \mu &= (0,0,2):[2] \to [3]
\end{aligned}
\]
Therefore
\[
\VL_{d^x_0(z)}(x) = \SETT{0,1,4}\quad \text{and} \quad \VL_{d^x_0(z)}(d_0(y)) = \SETT{2,3,5}
\]
The picture is:
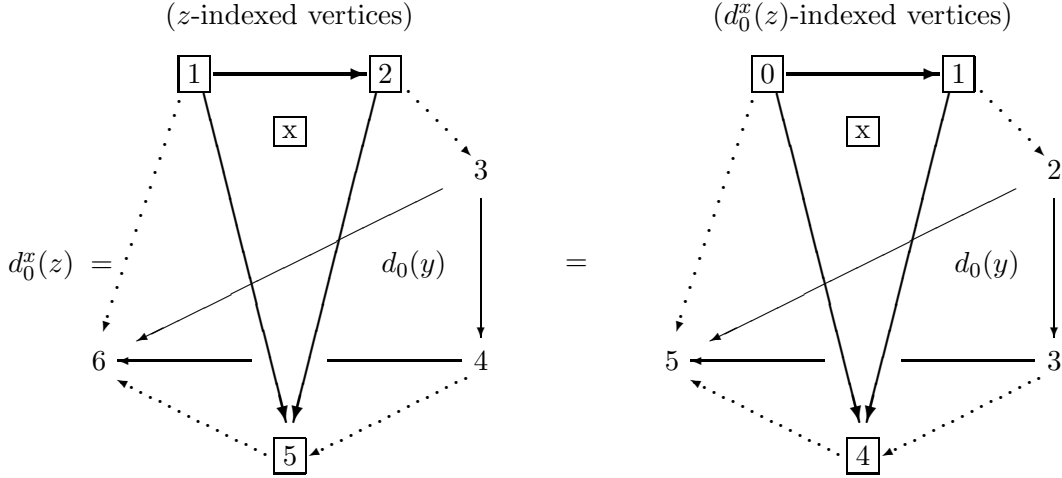
\begin{center}
\begin{picture}(6,2.4)
\put(1,2.3){\makebox(0,0){($z$-indexed vertices)}}
\put(-.2,1){\MB{d^x_0(z) \; =}}
\put(.5,2){\MB{\fbox{1}}}
\put(1.5,2){\MB{\fbox{2}}}
\put(2,1.5){\MB{3}}
\put(2,.5){\MB{4}}
\put(1,0){\MB{\fbox{5}}}
\put(0,.5){\MB{6}}
\thicklines
\put(.6,2){\vector(1,0){.8}} 
\put(1.45,1.9){\vector(-1,-4){.43}} 
\put(.55,1.9){\vector(1,-4){.43}} 
\thinlines
\put(2,1.35){\vector(0,-1){.73}}
\put(1.9,.5){\line(-1,0){.7}}
\put(.8,.5){\vector(-1,0){.7}}
\put(1.8,1.4){\vector(-2,-1){1.6}}
\put(1.65,1){\MB{d_0(y)}}
\put(1,1.7){\MB{\fbox{x}}}
\multiput(1.6,1.9)(.05,-.05){7}{.} 
\put(1.95,1.5725){\vector(1,-1){0}}
\multiput(1.9,.4)(-.05,-.025){16}{.} 
\put(1.1,0){\vector(-2,-1){0}}
\multiput(.85,0)(-.05,.025){16}{.}  
\put(.095,.4){\vector(-2,1){0}}
\multiput(.4,1.85)(-.025,-.075){16}{.} 
\put(.0225,.65){\vector(-1,-3){0}}
\put(2.5,1){\MB{=}}
\put(4,2.3){\makebox(0,0){($d^x_0(z)$-indexed vertices)}}
\put(3.5,2){\MB{\fbox{0}}} 
\put(4.5,2){\MB{\fbox{1}}} 
\put(4,0){\MB{\fbox{4}}} 
\thicklines
\put(3.6,2){\vector(1,0){.8}}
\put(4.45,1.9){\vector(-1,-4){.43}}
\put(3.55,1.9){\vector(1,-4){.43}}
\put(4,1.7){\MB{\fbox{x}}}
\thinlines
\put(5,1.5){\MB{2}} 
\put(5,.5){\MB{3}} 
\put(3,.5){\MB{5}} 
\put(4.8,1.4){\vector(-2,-1){1.6}}
\put(5,1.35){\vector(0,-1){.73}}
\put(4.9,.5){\line(-1,0){.7}}
\put(3.8,.5){\vector(-1,0){.7}}
\put(4.65,1){\MB{d_0(y)}}
\multiput(4.6,1.9)(.05,-.05){7}{.} 
\put(4.95,1.5725){\vector(1,-1){0}}
\multiput(4.9,.4)(-.05,-.025){16}{.} 
\put(4.1,0){\vector(-2,-1){0}}
\multiput(3.85,0)(-.05,.025){16}{.}  
\put(3.095,.4){\vector(-2,1){0}}
\multiput(3.4,1.85)(-.025,-.075){16}{.} 
\put(3.0225,.65){\vector(-1,-3){0}}

\end{picture}
\end{center}
\makebox(0,0){ }
}
\end{example}

\begin{newthm} 
\label{generalized-embedding}
{\rm
\SL

An arbitrary simplicial set $C$, $k \geq 0$, $u \in [k+1]$ and $x \in C_k$ determine a simplicial map $H^{x,u}:C \times \Delta[k] \to \mathcal{S}$.
}
\end{newthm}

\Proof

Given $(y,\LAM) \in C_j \times \Delta[k]_j$ and $j>0$ then for each $p \in [j]$, $\HMU xu{j-1}{d_p(y),\LAM \del_p} = d_p^\mathcal{S} \Bigl(\HM xj{y,\LAM} \Bigr)$ (where $d_p^\mathcal{S}$ denotes a face operator of $\mathcal{S}$). That is, $H^{x,u} : C \times \Delta[k] \to \mathcal{S}$ is a face map.

As for degeneracies, suppose $p \in [j]$ with $j \geq 0$. Then the signatures of $\HMU xu{j+1}{s_p(y), \LAM \SIG_p}$ and $s_p^\mathcal{S} \bigl( \HMU xuj{y,\LAM} \bigr)$ are equal. 

By the definition of degeneracies for $\mathcal{S}$, every element of $s_p^\mathcal{S} \bigl( \HM xj{y,\LAM} \bigr)$ is 
\[
c_p(h_j \LIST{t}{j})
\]
where
$h_j \LIST{t}{j} \in \HM xj{y,\LAM}$. By the degeneracy-case (b') and induction, $c_p h_j \LIST{t}{j} \in \HMU xu{j+1}{s_p(y), \LAM \SIG_p}$. That is, $s_p^\mathcal{S} \bigl( \HM xj{y,\LAM} \bigr) \SBS \HMU xu{j+1}{s_p(y), \LAM \SIG_p}$.
\MS

On the other hand, again by the degeneracy case, every element of $\HMU xu{j+1}{s_p(y), \LAM \SIG_p}$ has the form $c_p h_j \LIST{t}{j}$ where $h_j \LIST{t}{j} \in \HMU xuj{y,\LAM}$. Since $c_p h_j \LIST{t}{j} \in s_p^\mathcal{S} \bigl( \HM xj{y,\LAM} \bigr)$ then $\HMU xu{j+1}{s_p(y) \LAM \SIG_p} \SBS s_p^\mathcal{S} \bigl( \HM xj{y,\LAM} \bigr)$.
\MS

This verifies that $\HMU xu{j+1}{s_p(y), \LAM \SIG_p} = s_p^\mathcal{S} \bigl( \HM xj{y,\LAM} \bigr)$ as $(j+1)$-simplices.

\qed
\BS

\begin{newthm}
\label{represented-subcomplex}
{\rm
\SL

For each simplicial set $C$, $k \geq 0$ and $x \in C_k$, $\MC{R}(x,u)$ is a subcomplex of $\MC{S}$.
}
\end{newthm}

\Proof

Closure of $\MC{R}(x,u)$ under the simplicial operators of $\MC{S}$ follows directly from the previous theorem. 

\qed
\BS

\begin{newthm}
\label{Rxu-is-composer}
{\rm
\SL

Suppose $C$ is an \NICOMP{n}{i}, $k \geq 0$, $u \in [k+1]$ and $x \in C_k$. Then the representable subcomplex $\MC{R}(x,u) \SBS \mathcal{S}$ is also an \NICOMP{n}{i}.
}
\end{newthm}

\Proof

We will verify that $\BOX{i}{n+1}{\MC{R}(x,u)} \ISO \MC{R}(x,u)_{n+1}$ using the $(n,i)$-composition of $C$. Since the same argument using $(m,i)$-composition for all $m>n$ would show $\BOX{i}{m+1}{\MC{R}(x,u)} \ISO \MC{R}(x,u)_{m+1}$ it suffices to do the $m=n+1$ case.
\MS

So suppose
\[
\Bigl(
\HMU xun{y_0,\LAM_0} \DDD{,}  -  \DDD{,} \HMU xun{y_{n+1},\LAM_{n+1}}
\Bigr) \in \BOX{i}{n+1}{\MC{R}(x,u)}
\]
Since for $p<q$, $d_p(\HMU xun{y_q,\LAM_q}) = d_{q-1}(\HMU xun{y_p,\LAM_p})$, then $d_p(y_q)=d_{q-1}(y_p)$ and $\LAM_q \del_p = \LAM_p \del_{q-1}$. That is:
\[
\begin{aligned}
(\DLI{y}{n+1}{i} ) & \in \BOX{i}{n+1}{C}\\
(\DLI{\LAM}{n+1}{i} ) & \in \BOX{i}{n+1}{\Delta[k]}
\end{aligned}
\]
For general combinatorial reasons (see page \pageref{delta-n-is-hypergroupoid}) there exists a unique $\LAM \in\Delta[k]_{n+1}$ such that for all $p \in [n+1]-\SETT{i}$, $\LAM \del_p = \LAM_p$. Using $(n,i)$-composition in $C$, $y = \COMP{n}{i} (\DLI{y}{n+1}{i})$ is the unique $(n+1)$-simplex such that for all $p \in [n+1]-\SETT{i}$, $d_p(y) = y_p$ . Therefore 
\[
\Bigl(
\HMU xu{n+1}{y,\LAM} \SBS \prod_{p=0}^{n+1} \HMU xun{d_p(y), \LAM \del_p}
\Bigr) \in \MC{R}(x,u)_{n+1}
\]
is an $(n+1)$-simplex of $\MC{R}(x,u)$ such that
\[
\phi_{n+1,i}\Bigl( \HMU xu{n+1}{y,\LAM} \Bigr) = 
\Bigl(
\HMU xun{y_0,\LAM_0} \DDD{,}  -  \DDD{,} \HMU xun{y_{n+1},\LAM_{n+1}}
\Bigr)
\]
If 
\[
\Bigl(
\HMU xu{n+1}{y',\LAM'} \SBS \prod_{p=0}^{n+1} \HMU xun{d_p(y'), \LAM' \del_p}
\Bigr) \in \MC{R}(x,u)_{n+1}
\]
such that $\HMU xun{d_p(y'), \LAM' \del_p}= \HMU xun{d_p(y), \LAM \del_p}$ for each $p$ then for each $p$, $d_p(y')=d_p(y)$ and $\LAM \del_p = \LAM' \del_p$. This implies $y'=y$ and $\LAM'=\LAM$ by $(n,i)$-composition in $C$ and in $\Delta[k]$.
\MS

Hence $\phi_{n+1,i}: \MC{R}(x,u)_{n+1} \to \BOX{i}{n+1}{\MC{R}(x,u)}$ is an isomorphism.

\qed
\MS

We consider $H^{x,u}$ again below in section \refpage{D-times-delta-m}.

\section{Relating $C$ and $\mathcal{S}^C$}

In the subsections below we consider category structures associated with $\MC{S}$ (the simplicial set of relations) and $\MC{S}^C$, with $C$ an arbitrary simplicial set. Then we apply these ideas to the $k$-simplex $H^{x,u} : C \times \Delta[k] \to \MC{S}$ of \;$\MC{S}^C$ developed in the previous section. We will examine $x \mapsto H^{x,u}$ in terms of faces and degeneracies of $x \in C_k$.

\subsection{$\mathcal{S}$ as a simplicial object in {\it Cat}}

\begin{newdef}
{\bf (Map of relations)}
\label{map-of-relations}
\index{map of relations}
\index{domain function}
\index{signature function}
{\rm
\SL

Suppose $m,n \geq 0$ and $R \SBS \prod_{p=0}^m V_p$ and $R' \SBS \prod_{p=0}^n V'_p$ are relations such that $R, R' \neq \MT$. Then a {\bf map of relations}
\[
\Bigl( R \SBS \prod_{p=0}^m V_p \Bigr) \xrightarrow{(w,u)} \Bigl(R' \SBS \prod_{p=0}^n V'_p \Bigr)
\]
consists of functions $w : R \to R'$ and $u : \prod_{p=0}^m V_p \to \prod_{p=0}^n V'_p$ such that the following diagram commutes:

\begin{center}
\begin{picture}(2,.6)
\put(0,.5){\makebox(0,0){$R$}}
\put(0,0){\makebox(0,0){$R'$}}
\put(2,.5){\makebox(0,0){$\prod_{p=0}^m V_p$}}
\put(2,0){\makebox(0,0){$\prod_{p=0}^n V'_p$}}
\put(.2,.5){\vector(1,0){1.4}}
\put(1,.6){\makebox(0,0){\footnotesize{incl.}}}
\put(.2,0){\vector(1,0){1.4}}
\put(1,-.1){\makebox(0,0){\footnotesize{incl.}}}
\put(0,.375){\vector(0,-1){.25}}
\put(2,.375){\vector(0,-1){.25}}
\put(2.1,.25){\MBS{u}}
\put(-.1,.25){\MBS{w}}
\end{picture}
\end{center}
Since the horizontal maps are inclusions, this means simply that for all $t \in R$, $u(t) = w(t) \in R'$.

We will refer to $w$ as the {\bf ``domain function''} and to $u$ as the {\bf ``signature function''}.
} 

\BX
\end{newdef}
\MS

\NI  Clearly, the collection of relations forms a category whose maps are maps of relations in this sense. The next step is to specialize maps of relations to the $m$-simplices of $\mathcal{S}$ to define the notion of an ``$m$-simplex map'' taking into account the simplicial structure of $\mathcal{S}$.

\begin{newdef}
{\bf ($m$-simplex map)}
\label{m-simplex-map}
\index{m-simplex map}
{\rm
\SL

Given $m \geq 0$ and two $m$-simplices $y \SBS \prod_{p=0}^m d_p(y)$ and $y' \SBS \prod_{p=0}^m d_p(y')$ of $\MC{S}$ such that $y, y' \neq \MT$, an {\bf $m$-simplex map} 
\[
\SMAP{m}{p}{y}{y'}{w}{u}
\]
is a map of relations
\begin{center}
\begin{pictr}{2}{.6}{-.7}{.25}
\label{m-map-diag}
\put(0,.5){\makebox(0,0){$y$}}
\put(0,0){\makebox(0,0){$y'$}}
\put(2,.5){\makebox(0,0){$\prod_{p=0}^m d_p(y)$}}
\put(2,0){\makebox(0,0){$\prod_{p=0}^m d_p(y')$}}
\put(.2,.5){\vector(1,0){1.4}}
\put(1,.6){\makebox(0,0){\footnotesize{incl. = $i_y$}}}
\put(.2,0){\vector(1,0){1.4}}
\put(1,-.1){\makebox(0,0){\footnotesize{incl. = $i_{y'}$}}}
\put(0,.375){\vector(0,-1){.25}}
\put(-.1,.25){\MBS{w}}
\put(2,.375){\vector(0,-1){.25}}
\put(2.1,.25){\MBS{u}}
\end{pictr}
\end{center}
\MS

\NI specified inductively on $m$, as follows.
\MS

\NI {\bf Dimension $m=0$:}

A {\bf $0$-simplex map} $(y \SBS V) \xrightarrow{(w,u)} (y' \SBS V')$ is a commutative diagram:
\begin{center}
\begin{picture}(2,.6)
\put(0,.5){\makebox(0,0){$y$}}
\put(0,0){\makebox(0,0){$y'$}}
\put(2,.5){\makebox(0,0){$V$}}
\put(2,0){\makebox(0,0){$V'$}}
\put(.2,.5){\vector(1,0){1.6}}
\put(1,.6){\makebox(0,0){\MBS{i_y}}}
\put(.2,0){\vector(1,0){1.6}}
\put(1,-.1){\makebox(0,0){\MBS{i_{y'}}}}
\put(0,.375){\vector(0,-1){.25}}
\put(-.1,.25){\MBS{w}}
\put(2,.375){\vector(0,-1){.25}}
\put(2.1,.25){\MBS{u}}
\end{picture}
\end{center}
\MS

\NI {\bf Dimension $m>0$:}

An {\bf $m$-simplex map} 
is a map of relations (diagram (\ref{m-map-diag})\,) $(w,u)$ whose signature function $u : \prod_{p=0}^m d_p(y) \to \prod_{p=0}^m d_p(y')$ is defined by 
\begin{equation}
\label{def-of-u}
u \LIST{t}{m} \DFAS \bigl(
w_0(t_0) \DDD{,} w_m(t_m)
\bigr) \qquad t_p \in d_p(y)
\end{equation}
where, for each $k \in [m]$, $w_k : d_k(y) \to d_k(y')$ is the {\em domain function} of an $(m-1)$-simplex map 

\begin{center}
\begin{pictr}{1.5}{.6}{-.7}{.25}
\put(0,0){\MB{d_k(y')}}
\put(0,0.5){\MB{d_k(y)}}
\put(1.5,.5){\MB{\prod_{p=0}^{m-1} d_p d_k(y)}}
\put(1.5,0){\MB{\prod_{p=0}^{m-1} d_p d_k(y')}}
\put(.25,0){\vector(1,0){.7}}
\put(.25,.5){\vector(1,0){.7}}
\put(0,.35){\vector(0,-1){.2}}
\put(1.5,.35){\vector(0,-1){.2}}
\put(-.2,.25){\MBS{w_k}}
\put(1.7,.25){\MBS{u_k}}
\put(.55,.1){\MBS{i_{d_k(y')}}}
\put(.55,.6){\MBS{i_{d_k(y)}}}
\end{pictr}
\end{center}

We define $d^\MC{S}_k(w,u) \DFAS (w_k,u_k)$.
}

\BX
\end{newdef}

It follows from the definition that for all $q \in [m]$ the commutativity of the following diagram:
\begin{center}
\begin{picture}(3,1.1)
\put(0,1){\MB{y}}
\put(1,1){\MB{\prod_{p=0}^m d_p(y)}}
\put(3,1){\MB{d_q(y)}}
\put(0,0){\MB{y'}}
\put(1,0){\MB{\prod_{p=0}^m d_p(y')}}
\put(3,0){\MB{d_q(y')}}
\put(.2,1){\vector(1,0){.3}}
\put(1.4,1){\vector(1,0){1.3}}
\put(2,1.1){\MBS{\text{proj}_q}}
\put(.2,0){\vector(1,0){.3}}
\put(1.4,0){\vector(1,0){1.3}}
\put(2,.1){\MBS{\text{proj}_q}}
\put(0,.8){\vector(0,-1){.6}}
\put(-.1,.5){\MBS{w}}
\put(1,.8){\vector(0,-1){.6}}
\put(.9,.5){\MBS{u}}
\put(3,.8){\vector(0,-1){.6}}
\put(3.1,.5){\MBS{w_q}}
\end{picture}
\end{center}
\MS

\NI implies that given $t\in y$, then $w_q(e_q(t)) = e_q(w(t))$.
 Similarly, given any $\SETT{k_1 \DDD{,}k_j} \SBS [m]$ with $j<m$ then
\[
e_{k_1} \DDD{}e_{k_j}( w(t) ) = w_{e_{k_1} \DDD{}e_{k_j}}(e_{k_1} \DDD{}e_{k_j}(t))
\]
where $w_{e_{k_1} \DDD{}e_{k_j}} : d_{k_1} \DDD{}d_{k_j}(y) \to d_{k_1} \DDD{}d_{k_j}(y')$ is the domain function of an $(m-j)$-simplex map between the corresponding subfaces of $y$ and $y'$
\MS

The composition of $m$-simplex maps
\[
\Bigl(y \SBS \prod_{p=0}^m d_p(y) \Bigr) \XRA{(w,u)} \Bigl(y' \SBS \prod_{p=0}^m d_p(y') \Bigr) \XRA{(w',u')} \Bigl(y'' \SBS \prod_{p=0}^m d_p(y'') \Bigr)
\]
is 
\[
\SMAP{m}{p}{y}{y''}{w' w}{u'u}
\]
It follows directly that for each $k \in [m]$, $(u'u)_k = u'_k \, u_k$ and $(w'w)_k = w'_k\, w_k$. 
That is:
\begin{equation}
\label{d^S-functor}
d^\MC{S}_k(w'w,u'u) = d^\MC{S}_k(w',u')\, d^\MC{S}_k(w,u)
\end{equation}
\MS

Next, given any $m \geq 0$, we'll show that for any $m$-simplex map
\begin{equation} \label{given-smap}
\SMAP{m}{p}{y}{y'}{w}{u}
\end{equation}
and any $q \in [m]$ there are functions $s_q(w)$ and $s_q(u)$ (defined below) which comprise an $(m+1)$-simplex map:

\begin{center}
\begin{pictr}{2}{.6}{-.7}{.25}
\label{degen-m-map-diag}
\put(0,.5){\makebox(0,0){$s_q(y)$}}
\put(0,0){\makebox(0,0){$s_q(y')$}}
\put(2,.5){\makebox(0,0){$\prod_{p=0}^{m+1} d_p s_q(y)$}}
\put(2,0){\makebox(0,0){$\prod_{p=0}^{m+1} d_p s_q(y')$}}
\put(.25,.5){\vector(1,0){1.2}}
\put(.9,.6){\MBS{i_{s_q(y)}}}
\put(.25,0){\vector(1,0){1.2}}
\put(.9,.1){\MBS{i_{s_q(y')}}}
\put(0,.375){\vector(0,-1){.25}}
\put(-.2,.25){\MBS{s_q(w)}}
\put(2,.375){\vector(0,-1){.25}}
\put(2.2,.25){\MBS{s_q(u)}}
\end{pictr}
\end{center}
\smallskip

\NI The functions $s_q(w)$ and $s_q(u)$ are determined by induction and  the relationship between $w$ and $u$.
\MS

\NI {\bf Definition of $s_q(w)$:}

$s_q(w) \DFAS c^{y'}_q w (c^y_q)^{-1}$ \qquad i.e. \hspace{.5in}
\begin{picture}(1.5,.65)
\put(0,0){\MB{y}}
\put(0,.5){\MB{s_q(y)}}
\put(1,.5){\MB{s_q(y')}}
\put(1,0){\MB{y'}}
\put(.1,0){\vector(1,0){.8}}
\put(.2,.45){\vector(1,0){.6}}
\put(0,.1){\vector(0,1){.25}}
\put(1,.1){\vector(0,1){.25}}
\put(.5,.08){\MBS{w}}
\put(.5,.55){\MBS{s_q(w)}}
\put(-.15,.25){\MBS{c^y_q}}
\put(1.2,.25){\MBS{c^{y'}_q}}
\put(1.5,0){commutes.}
\end{picture} 
\MS

\NI Coordinate-wise:
\[
\begin{aligned}
s_q(w)(c_q(t)) &= c_q(w(t))\\
&=
\Bigl( c_{q-1}w_0(t_0) \DDD{,} c_{q-1}w_{q-1}(t_{q-1}),w(t),w(t),\\
& \qquad c_q w_{q+1}(t_{q+1}) \DDD{,} c_q w_m(t_m) \Bigr)
\end{aligned}
\]
so that
\[
e_p(s_q(w)(c_q(t)) = 
\left\{
\begin{array}{lc}
c_{q-1}w_p(t_p) & 0 \leq p <q\\
w(t) & p=q,q+1\\
c_q w_{p-1}(t_{p-1}) & q+1<p \leq m+1
\end{array}
\right.
\]
\BS

\NI {\bf Definition of $s_q(u): \prod_{p=0}^{m+1} d_p s_q(y) \to \prod_{p=0}^{m+1} d_p s_q(y)$:}
\MS

$s_q(u)$ is defined by requiring the commutativity for each $k \in [m+1]$ of

\begin{center}
\begin{picture}(1.5,.6)
\put(0,0){\MB{\prod d_p s_q(y')}}
\put(0,.5){\MB{\prod d_p s_q(y)}}
\put(1.5,.5){\MB{d_k s_q(y)}}
\put(1.5,0){\MB{d_k s_q(y')}}
\put (.4,0){\vector(1,0){.7}}
\put (.4,.5){\vector(1,0){.7}}
\put(0,.35){\vector(0,-1){.2}}
\put(1.5,.35){\vector(0,-1){.2}}
\put(-.2,.25){\MBS{s_q(u)}}
\put(1.9,.25){\MBS{\PR{k}(s_q(u))}}
\put(.75,.6){\MBS{\PR{k}}}
\put(.75,-.1){\MBS{\PR{k}}}
\end{picture}
\end{center}

\NI where

\begin{equation}
\label{s_q(u)-eqn}
\PR{k}(s_q(u)) \DFAS
\left\{
\begin{array}{lc}
s_{q-1}(w_k) &  0 \leq k < q\\
w & k=q, q+1\\
s_q(w_{k-1}) & q+1<k \leq m+1
\end{array}
\right.
\end{equation}
That is,
\[
s_q(u) = s_{q-1}(w_0) \DDD{\times} s_{q-1}(w_{q-1}) \times w \times w \times s_q(w_{q+1}) \DDD{\times} s_q(w_{m})
\]

The inductive nature of this arises from $w_p:d_p(y) \to d_p(y')$ being the domain map of an $(m-1)$-simplex map.
\MS

\NI {\bf Claim:} $s_q(u) i_{s_q(y)} = i_{s_q(y')} s_q(w)$; that is,

\begin{center}
\begin{pictr}{2}{.65}{-.7}{.25}
\label{degen-m-map-claim}
\put(0,0){\MB{s_q(y')}} 
\put(0,.5){\MB{s_q(y)}}
\put(0,0){\MB{s_q(y')}}
\put(1.3,0){\MB{\prod_{p=0}^{m+1}d_p s_q(y')}}
\put(1.3,.5){\MB{\prod_{p=0}^{m+1}d_p s_q(y)}}
\put(.2,0){\vector(1,0){.6}}
\put(.2,.5){\vector(1,0){.6}}
\put(0,.35){\vector(0,-1){.25}}
\put(1,.35){\vector(0,-1){.25}}
\put(.5,.1){\MBS{i_{s_q(y')}}}
\put(.5,.6){\MBS{i_{s_q(y)}}}
\put(-.2,.25){\MBS{s_q(w)}}
\put(1.2,.25){\MBS{s_q(u)}}
\put(2,.25){commutes.}
\end{pictr}
\end{center}

\NI Let $t = \LIST{t}{m} \in y$. Recall that
\[
c_q(t) = \Bigl(c_{q-1}(t_0) \DDD{,} c_{q-1}(t_{q-1}),t,t,c_q(t_{q+1}) \DDD{,} c_q(t_m) \Bigr)
\]

\NI We verify the commutativity of diagram (\ref{degen-m-map-claim}) coordinate-wise.
\MS

We have
\[
e_p \Bigl( \, s_q(u) \, i_{s_q(y)}\, c_q(t) \Bigr) = 
\left\{
\begin{array}{ll}
s_{q-1}(w_p)(c_{q-1}(t_p)) = c_{q-1}(w_p(t_p)) & 0 \leq p <q\\
i_{s_q(y)} w(t) & p=q,q+1\\
s_q(w_{p-1})(c_q(t_{p-1})) = c_q(w_{p-1}(t_{p-1})) & q+1<p \leq m+1
\end{array}
\right.
\]
and
\[
e_p \Bigl( i_{s_q(y')}\, s_q(w)(c_q(t)) \Bigr) = 
\left\{
\begin{array}{ll}
 c_{q-1}(w_p(t_p)) & 0 \leq p<q\\
 i_{s_q(y)} w(t) & p = q,q+1\\
 c_q(w_{p-1}(t_{p-1})) & q+1<p \leq m+1
\end{array}
\right.
\]
\BX
\MS

We define $s^\MC{S}_q(w,u) \DFAS (s_q(w),s_q(u))$.
%
\BS

\begin{newthm}
{\rm
\SL

For each $m \geq 0$, the $m$-simplices of $\MC{S}$ comprise the objects of a category whose maps are $m$-simplex maps. The maps $d^\MC{S}_k$ and $s^\MC{S}_k$ satisfy the simplicial identities. And, for each $k \in [m]$
\[
\begin{aligned}
d^\MC{S}_k : \MC{S}_m \to \MC{S}_{m-1} & \qquad (m>0)\\
s^\MC{S}_k : \MC{S}_m \to \MC{S}_{m+1} &
\end{aligned}
\] are functors.
}
\end{newthm}

\Proof

That $\MC{S}_m$ together with $m$-simplex maps comprises a category is immediate from the definitions, as is that the $d^\MC{S}_k$ and $s^\MC{S}_k$ satisfy the simplicial identities.
\MS

Consider a composition $(w'w,u'u)$ of $m$-simplex maps:

\begin{center}
\begin{pictr}{3}{.6}{-.7}{.25}
\put(0,0){\MB{y}} 
\put(1.5,0){\MB{y'}}
\put(3,0){\MB{y''}}
\put(0,.5){\MB{\prod_{p=0}^m d_p(y)}}
\put(1.5,.5){\MB{\prod_{p=0}^m d_p(y')}}
\put(3,.5){\MB{\prod_{p=0}^m d_p(y'')}}
\put(.15,0){\vector(1,0){1.2}}
\put(1.65,0){\vector(1,0){1.2}}
\put(.4,.5){\vector(1,0){.65}}
\put(1.9,.5){\vector(1,0){.65}}
\put(0,.15){\vector(0,1){.2}}
\put(1.5,.15){\vector(0,1){.2}}
\put(3,.15){\vector(0,1){.2}}
\put(.7,.1){\MBS{w}}
\put(2.2,.1){\MBS{w'}}
\put(.7,.6){\MBS{u}}
\put(2.2,.6){\MBS{u'}}
\end{pictr} 
\end{center}
\MS

As noted in equation \eqref{d^S-functor}, $d^\MC{S}_k$ respects $m$-simplex composition. We now check $s^\MC{S}_k$.
\MS

By definition, $s_k(w) c_k(t) = c_k(w(t))$. Then
\[
\begin{aligned}
 s_k(w') s_k(w) c_k(t) & = s_k(w') c_k(w(t)) = c_k(w'(w(t)))\\
 &= c_k((w'w)(t)) = s_k(w'w)(c_k(t))
\end{aligned}
\]
which shows $s_k(w'w) = s_k(w') s_k(w)$.

Similarly, by equation \eqref{s_q(u)-eqn} we have
\[
\PR{k} \Bigl( s_k(u'u) \Bigr) = 
\left\{
\begin{array}{ll}
s_{k-1}(w'_p w_p) = s_{k-1}(w'_p) s_{k-1}(w_p) & 0 \leq p<k\\
w'w & p=k, k+1\\
s_k(w'_{p-1} w_{p-1}) = s_k(w'_{p-1}) s_k(w_{p-1}) & k+1<p \leq m+1
\end{array}
\right.
\]
which shows that $s_k(u'u) = s_k(u') s_k(u)$.

Therefore, 
\[
s^\MC{S}_k(w'w,u'u)=(s_k(w'w),s_k(u'u)) = (s_k(w') s_k(w),s_k(u') s_k(u))
 = s^\MC{S}_k(w',u') s^\MC{S}_k(w,u)
\]
\qed

\subsection{$\mathcal{S}^C$ as a simplicial object in {\it Cat}}

In this section we will define a category structure on the set of $k$-simplices of $\MC{S}^C$ where $C$ is an arbitrary simplicial set.
\MS

Given a $k$-simplex $F:C \times \Delta[k] \to \MC{S}$ then for each $j > 0$ and each $(y,\LAM) \in C_j \times \Delta[k]_j$ we will denote the $j$-simplex $F_j(y,\LAM)$  of $\MC{S}$ by
\[
F_j(y,\LAM) \SBS \prod_{p=0}^j F_{j-1}(d_p(y),\LAM \del_p)
\]
and for $j=0$ write $F_0(y,\LAM)$ as
\[
F_0(y,\LAM) \SBS V^F_0(y,\LAM) \DFAS \text{ the signature of } F_0(y,\LAM)
\]

\begin{newdef}
{\bf ($k$-simplex tranform)}
\label{k-simplex-transform-in-S^C}
\index{k-simplex transform in $\mathcal{S}^C$}
{\rm
\SL

Suppose $C$ is a simplicial set and $\mathcal{S}$ is the simplicial set of relations. Given $k \geq 0$ and two $k$-simplices of $\MC{S}^C$ (simplicial maps)
\[
F, G : C \times \Delta[k] \to \MC{S}
\]
then a {\bf $k$-simplex transform of $\mathcal{S}^C$}, denoted $w : F \to G$, and defined inductively on dimension $j$ 
consists of:
\begin{itemize}
\item 
a family of $0$-simplex maps for all $(y,\LAM) \in C_0 \times \Delta[k]_0$

\begin{center}
\begin{pictr}{2}{.63}{-.7}{.25}
\put(0,0){\MB{F_0(y,\LAM)}} 
\put(2,0){\MB{G_0(y,\LAM)}}
\put(0,.5){\MB{V^F_0(y, \LAM)}}
\put(2,.5){\MB{V^G_0(y, \LAM)}}
\put(.35,0){\vector(1,0){1.3}}
\put(.35,.5){\vector(1,0){1.3}}
\put(0,.15){\vector(0,1){.2}}
\put(2,.15){\vector(0,1){.2}}
\put(.2,.25){\MBS{\text{incl.}}}
\put(2.2,.25){\MBS{\text{incl.}}}
\put(1,.1){\MBS{w(0,y,\LAM)}}
\put(1,.6){\MBS{u(0,y,\LAM)}}
\end{pictr}
 
\end{center}

\item
for each $j>0$, a family of $j$-simplex maps
\[
\Bigl\{
\bigl( w(j,y,\LAM),u(j,y,\LAM) \bigr):
(F_j(y,\LAM) \SBS \Sig(F_j(y,\LAM)) \to (G_j(y,\LAM) \SBS \Sig(G_j(y,\LAM))
\Bigr\}
\]
indexed by $(y,\LAM) \in C_j \times \Delta[k]_j$
where
\[
u(j,y,\LAM) = w\bigl(j-1,d_0(y),\LAM \del_0\bigr)\, \DDD{\times}\, w\bigl(j-1,d_j(y),\LAM \del_j\bigr)
\]

\begin{center}
\begin{pictr}{2.5}{.63}{-1.2}{.25}
\put(0,0){\MB{F_j(y,\LAM)}} 
\put(2.5,0){\MB{G_j(y,\LAM)}}
\put(0,.5){\MB{\prod_{p=0}^j F_{j-1}(d_p(y), \LAM \del_p)}}
\put(2.5,.5){\MB{\prod_{p=0}^j G_{j-1}(d_p(y), \LAM \del_p)}}
\put(.35,0){\vector(1,0){1.7}}
\put(.8,.5){\vector(1,0){.85}}
\put(0,.15){\vector(0,1){.2}}
\put(2.5,.15){\vector(0,1){.2}}
\put(1.25,.1){\MBS{w(j,y,\LAM)}}
\put(1.25,.6){\MBS{u(j,y,\LAM)}}
\put(.2,.25){\MBS{\text{incl.}}}
\put(2.7,.25){\MBS{\text{incl.}}}
\end{pictr}
\end{center}

\end{itemize}

\MS

The {\bf identity} $k$-simplex transform $F \to F$ consists of the identity $j$-maps for all $j \geq 0$.
}

\BX
\end{newdef}

That this defines a category structure on $(\MC{S}^C)_k$ follows from the discussion of maps of relations above.
\BS

Recall that the face and degeneracy operators for $\MC{S}^C$, $d^{\MC{S}^C}_q$ and $s^{\MC{S}^C}_q$ are defined as follows. We will simplify the notations $d^{\MC{S}^C}_q$ to $d_q$ and $s^{\MC{S}^C}_q$ to $s_q$.
\MS

Given any $k$-simplex $F: C \times \Delta[k] \to \mathcal{S}$ and any $q \in [k]$ then
\begin{itemize}
\item 
For $k>0$
\[
d_q(F) \text{ is: }\qquad C \times \Delta[k-1] \XRA{1 \times \delta^*_q} C \times \Delta[k] \XRA{F} \MC{S}
\]
where, for each $j \geq 0$, and each $(y,\LAM) \in C_j \times \Delta[k-1]_j$, \[
(1 \times \delta^*_q)_j(y,\LAM) \DFAS (y,\del_q \LAM)
\]
That is
\[
(d_q F)_{j}(y,\LAM) \DFAS F_j(y,\del_q \LAM)
\]
\MS

For $j=0$, $q \in [k]$ and $(y,\LAM) \in C_0 \times \Delta[k-1]_0$, $(d_qF)_0 (y,\LAM)$ is 
\[
F_0(y,\del_q \LAM) \SBS V^F_0(y,\del_q \LAM)
\]

For $j>0$, $q \in [k]$ and $(y,\LAM) \in C_j \times \Delta[k-1]_j$, $(d_qF)_j (y,\LAM)$ is
\[
(d_qF)_j (y,\LAM) \SBS \prod_{p=0}^j F_{j-1}(d_p(y), \del_q \LAM \del_p)
\]

\item
For $k \geq 0$
\[
s_q(F) \text{ is} \qquad C \times \Delta[k+1] \XRA{1 \times \sigma^*_q} C \times \Delta[k] \XRA{F} \MC{S}
\] 
where, for each $j \geq 0$, and each $(y,\LAM) \in C_j \times \Delta[k+1]_j$, \[
(1 \times \SIG^*_q)_j(y,\LAM) \DFAS (y,\SIG_q \LAM) 
\]
That is
\[
(s_q F)_{j}(y,\LAM) \DFAS F_j(y,\SIG_q \LAM)
\]
\MS

For $j=0$, $q \in [k]$ and $(y,\LAM) \in C_0 \times \Delta[k-1]_0$, $(s_qF)_0 (y,\LAM)$ is 
\[
F_0(y,\SIG_q \LAM) \SBS V^F_0(y,\SIG_q \LAM)
\]

For $j>0$, $q \in [k]$ and $(y,\LAM) \in C_j \times \Delta[k-1]_j$, $(s_qF)_j (y,\LAM)$ is
\[
(s_qF)_j (y,\LAM) \SBS \prod_{p=0}^j F_{j-1}(d_p(y), \SIG_q \LAM \del_p)
\]
\BX

\end{itemize}

\begin{newthm}
{\rm
\SL

Let $C$ be a simplicial set. Then for each $k \geq 0$, $\MC{S}^C_k$ is a category with respect to $k$-simplex transforms. Also, for each $q \in [k], k>0$, $d^{\MC{S}^C}_q: \MC{S}^C_k \to \MC{S}^C_{k-1}$ is a functor and for each $q \in [k], k\geq 0$, $s^{\MC{S}^C}_q: \MC{S}^C_k \to \MC{S}^C_{k+1}$ is a functor.
}
\end{newthm}

\Proof

That $\MC{S}^C_k$ is a category with respect to $k$-simplex transforms follows directly from the definitions and and compositions of $k$-simplex maps.
\MS

Suppose $F,G \in (\MC{S}^C)_k$, $k>0$ and $w : F \to G$ is a $k$-simplex transform of $\MC{S}^C$. For any $j \geq 0$ and any $(y,\LAM) \in C_j \times \Delta[k-1]_j$, we have $(d_q(F))_j(y,\LAM) \DFAS F_j(y,\del_q \LAM)$ and $(d_q(G))_j(y,\LAM) \DFAS G_j(y,\del_q \LAM)$. Then $(d_q(w))(j,y,\LAM) \DFAS w(j,y,\del_q \LAM)$ defines the claimed $(k-1)$-simplex transform $d_q(w) : d_q(F) \to d_q(G)$. It follows from the definition of composition of $k$-simplex maps that $d_q$ is a functor.
\MS

Similarly, for $k \geq 0$ and $q \in [k]$, we define $(s_q(w))(j,y,\LAM) \DFAS w(j,y,\SIG_q \LAM)$ to obtain a $(k+1)$-simplex transform $s_q(w) : s_q(F) \to s_q(G)$. That $s_q$ is a functor follows, again, from the definition of composition of $k$-simplex maps.

\qed

\subsection{$k$-simplex transforms involving $\MC{R}(x,u)$}

The next theorem refers to the subcomplex $\MC{R}(x,u)$ and its simplices $H^{x,u}$ defined above (page \pageref{rep-subcomplex}).

\begin{newthm}
\label{H-map-theorem}
{\rm
\SL

Suppose $C$ is a simplicial set, $k \geq 0$, $x \in C_k$, $q \in [k]$ and $u \in [k+1]$. Then:
\begin{enumerate}
\item 
If $u<q$ then there is a $k$-simplex transform $w: H^{x,u} \to s_{q-1}(H^{d_q(x),u})$.

\item
If $u>q+1$ then there is a $k$-simplex transform $w: H^{x,u} \to s_{q}(H^{d_q(x),u-1})$.
\end{enumerate}
}
\end{newthm}

\Proof

The proof of each statement follows from this: if $z \in \CPL(x,\mu;y, \del_u \LAM) \SBS C^x_j$ then for each $q \in [k]$, $d^y_q(z) \in \CPL(d_q(x), \mu \del_q;y, \SIG_q \del_u \LAM)$. By simplicial identities:
$\SIG_q \del_u \LAM = \del_u \SIG_{q-1} \LAM$ if $u<q$ and $\SIG_q \del_u \LAM = \del_{u-1} \SIG_{q} \LAM$ if $u>q+1$.
\MS

For each of the two claims, we define the domain functions of $w$ at dimension $j$ (denoted $w(j,y,\LAM)$, as before) starting at $j=0$, doing the $j>0$ definitions inductively, and then verifying the required commutativities.
\MS

\NI {\bf Case $u<q$:}

Given any $(y,\LAM) \in C_0 \times \Delta[k]_0$, the the elements of $H^{x,u}_0(y,\LAM)$ are $z \in \CPL(x,\mu;y,\del_u \LAM)$. For such $z$ and $q>u$ we have
\[
d^y_q(z) \in \CPL(d_q(x), \mu \del_q;y,\SIG_q \del_u \LAM) = 
\CPL(d_q(x), \mu \del_q;y,\del_u \SIG_{q-1} \LAM) 
\]
and so $d^y_q(z)$ sponsors an element of
\[
 H^{d_q(x),u}_0(y,\SIG_{q-1} \LAM) = s_{q-1}\bigl( H^{d_q(x),u} \bigr)_0(y,\LAM)
\]
We set $w(0,y,\LAM) \DFAS d^y_q$, and it is immediate that
\begin{center}
\begin{picture}(2,.6)
\put(0,.5){\MB{H^{x,u}_0(y,\LAM)}}
\put(0,0){\MB{H^{d_q(x),u}_0(y,\SIG_q \LAM)}}
\put(2,.5){\MB{C^x_0}}
\put(2,0){\MB{C^{d_q(x)}_0}}
\put(.6,.5){\vector(1,0){1.1}}
\put(1,.6){\makebox(0,0){\footnotesize{incl.}}}
\put(.6,0){\vector(1,0){1.1}}
\put(1,0.1){\makebox(0,0){\footnotesize{incl.}}}
\put(0,.4){\vector(0,-1){.25}}
\put(-.3,.25){\MBS{w(0,y,\LAM)}}
\put(2,.4){\vector(0,-1){.25}}
\put(2.15,.25){\MBS{d^y_q}}
\end{picture}
\end{center}
\MS

\NI commutes for all $(y,\LAM) \in C_0 \times \Delta[k]_0$.
\MS

Given $j>0$, $(y,\LAM) \in C_j \times \Delta[k]_j$ and any $h_j \LIST{t}{j} \in H^{x,u}_j(y,\LAM)$ sponsored by $z \in \CPL(x,\mu;y,\del_u \LAM)$, we define $w(j,y,\LAM)(h_j \LIST{t}{j})$ to be the element of $H^{d_q(x),u}_j(y, \SIG_{q-1} \LAM)$ sponsored by $d^y_q(z)$.
\MS

To verify the commutativity of
\begin{center}
\begin{picture}(3,.6)
\put(0,.5){\MB{\HMU xuj{y,\LAM}}}
\put(3,.5){\MB{\prod_{p=0}^j \HMU xu{j-1}{d_p(y), \LAM \del_p}}}
\put(0,0){\MB{\HMU {d_q(x)}uj{y,\SIG_{q-1}\LAM}}}
\put(3,0){\MB{\prod_{p=0}^j \HMU {d_q(x)}u{j-1}{d_p(y), \SIG_{q-1}\LAM \del_p}}}
\put(.6,.5){\vector(1,0){1.2}}
\put(1.4,.6){\makebox(0,0){\footnotesize{incl.}}}
\put(.7,0){\vector(1,0){1.2}}
\put(1.4,0.1){\makebox(0,0){\footnotesize{incl.}}}
\put(0,.4){\vector(0,-1){.25}}
\put(-.35,.25){\MBS{w(j,y,\LAM)}}
\put(2.8,.35){\vector(0,-1){.2}}
\put(3.7,.25){\MBS{\prod_{p=0}^j w(j-1,d_p(y),\LAM \del_p)}}
\end{picture}
\end{center}
\BS

\NI it suffices to show, for each $p \in [j]$ the commutativity of
\begin{center}
\begin{picture}(3,.6)
\put(0,.5){\MB{\HMU xuj{y,\LAM}}}
\put(3,.5){\MB{\HMU xu{j-1}{d_p(y), \LAM \del_p}}}
\put(0,0){\MB{\HMU {d_q(x)}uj{y,\SIG_{q-1}\LAM}}}
\put(3,0){\MB{\HMU {d_q(x)}u{j-1}{d_p(y), \SIG_{q-1}\LAM \del_p}}}
\put(.6,.5){\vector(1,0){1.5}}
\put(1.4,.6){\makebox(0,0){\footnotesize{proj$_p \circ$ incl.}}}
\put(.7,0){\vector(1,0){1.2}}
\put(1.4,0.1){\makebox(0,0){\footnotesize{proj$_p \circ$ incl.}}}
\put(0,.35){\vector(0,-1){.2}}
\put(-.35,.25){\MBS{w(j,y,\LAM)}}
\put(2.8,.35){\vector(0,-1){.2}}
\put(3.5,.25){\MBS{w(j-1,d_p(y),\LAM \del_p)}}
\end{picture}
\end{center}
\BS

\NI So suppose $h_j \LIST{t}{j} \in \HMU xuj{y, \LAM}$ is sponsored by $z \in \CPL(x,\mu;y,\del_u \LAM)$. Then:
\begin{itemize}
\item
$(\text{proj}_p \circ \text{ incl.})(h_j \LIST{t}{j})$ is sponsored by $d^x_p(z) \in \CPL(x,\SIG_p \mu;d_p(y),\del_u \LAM \del_p)$.

\item 
$\Bigl( w(j-1,d_p(y), \LAM \del_p) \circ \text{proj}_p \circ \text{ incl.} \Bigr)(h_j \LIST{t}{j})$ is sponsored by $d^{d_p(y)}_q d^y_p(z)$.

\item
$w(j,y,\LAM)(h_j \LIST{t}{j})$ is sponsored by $d^y_q(z)$.

\item
$\Bigl(
\text{proj}_p \circ \text{ incl.} \circ w(j,y,\LAM)
\Bigr)(h_j \LIST{t}{j})$ is sponsored by $d^{d_q(x)}_p d^y_q(z)$.
\end{itemize}

By Theorem \refpage{cpl-commute}, $d^{d_q(x)}_p d^y_q(z)=d^{d_p(y)}_q d^y_p(z)$. Therefore the diagrams above all commute, and this establishes the claimed $k$-simplex map $H^{x,u} \XRA{w} s_{q-1}(H^{d_q(x),u-1})$.
\BS

The claim in the case $u>q+1$ follows by an identical argument using the simplicial identity $\SIG_q \del_u = \del_{u-1} \SIG_q$.

\qed

\begin{newcor}
{\rm
\SL

Given a simplicial set $C$, $k \geq 0$, $q \in [k]$ and $u \in [k+1]$ then:
\begin{enumerate}
\item
If $u<q$ then

\begin{enumerate} 
\item 
If $k>0$ and $x \in C_{k-1}$ then there is a $k$ simplex transform $H^{s_q(x),u} \to s_{q-1}(H^{x,u})$.

\item
Given $x \in C_k$ there is a $k$-simplex transform $d_q(H^{x,u}) \to H^{d_q(x),u}$.

\item
If $q>0$, and $x \in C_k$ then there is a $k$-simplex transform $d_{q-1}(H^{x,u}) \to H^{d_q(x),u}$.
\end{enumerate} 

\item
If $u>q+1$ then

\begin{enumerate} 
\item 
If $k>0$ and $x \in C_{k-1}$ then there is a $k$-simplex transform $H^{s_q(x),u} \to s_{q}(H^{x,u-1})$.

\item
If $q<k$ and $x \in C_k$ then there is a $k$-simplex transform $d_{q+1}(H^{x,u}) \to H^{d_q(x),u-1}$. 

\item
If $x \in C_k$ then there is a $k$-simplex transform $d_{q}(H^{x,u}) \to H^{d_q(x),u-1}$. 

\end{enumerate} 

\end{enumerate}

}
\end{newcor}

\Proof

\begin{enumerate}
\item Assume $u<q$. Then $H^{s_q(x),u} \to s_{q-1}(H^{x,u})$ arises from applying the previous theorem to the $k$-simplex $s_q(x)$. 

The $k$-simplex transforms $d_q(H^{x,u}) \to H^{d_q(x),u}$ and $d_{q-1}(H^{x,u}) \to H^{d_q(x),u}$ arise from applying the operators $d_q$ and $d_{q-1}$ to the $k$-simplex map $H^{x,u} \to s_{q-1}(H^{d_q(x),u})$ from the previous theorem.

\item
In the case $u>q+1$ the same reasoning yields the claimed $k$-simplex transforms.
\end{enumerate}

\qed

\section{Raw material}
\label{raw-material}

The subsections below examine certain higher-dimensional versions of some concepts from category theory. The ideas here are included because they seem interesting, but they are provisional with regard to their appropriateness and relevance as generalizations. Even if relevant, they need further development.\footnote{Revisions of this paper may address these issues.}

\subsection{$D \times \Delta[m] \to C$ and $H^{x,u}$}
\label{D-times-delta-m}
In this subsection we will consider a higher-dimensional analog of the following dimension 1 notion: whenever $C$ and $D$ are categories, $F : D \to C$ is a functor and $x \in \OB (C)$, then there is a hom functor $H : D \to \STS$ dependant on $x$ and $F$ defined by
\[
\begin{aligned}
H(y) & \DFAS  C(x,F(y))\\
H(y \XRA{t}y') & \DFAS  C(x,F(y)) \to C(x,F(y')), \quad g \mapsto F(t) \circ g
\end{aligned}
\]
\BS

\NI We will generalize this, guided by the observation that a map $x\XRA{z} F(y)$ in $H(y)$ has $x$ and $F(y)$ as complementary subfaces. 
\MS

\NI For the generalization:
\MS

\NI $\bullet$ Replace the nerves of the small categories $C$ and $D$ with arbitrary simplicial sets $C$ and $D$.

\NI $\bullet$ Replace the functor $D \XRA{F} C$ (a simplicial map of the nerves, and a $0$-simplex of $C^D$) with a choice of $m \geq 0$ and $m$-simplex of $C^D$
\[
F : D \times \Delta[m] \to C
\]

\NI $\bullet$ Replace $x \in C_0$ with a choice of $k \geq 0$, $x \in C_k$ and a choice of $u \in [k+1]$.

\NI $\bullet$ Replace $1_{[0]} : [0] \to [0]$ with a fixed non-decreasing $\rho:[m] \to [k]$, $\rho \in \Delta[k]_m$.

\NI $\bullet$ Replace the functor $H:D \to \STS$ with the simplicial map
\[
D \times \Delta[m] \XRA{F*\rho} C \times \Delta[k] \XRA{H^{x,u}} \MC{S}
\]
to be defined below. (See page \pageref{rep-subcomplex} for the definition of $H^{x,u}$).
\BS

For the rest of this subsection, $C$ and $D$ denote arbitrary simplicial sets. We choose and fix $k, m \geq 0$, $x \in C_k$ and a simplicial map $F : D \times \Delta[m] \to C$.
\MS

\label{F-star-rho-def}
Given any $\rho : [m] \to [k]$ there is a simplicial map
\[
D \times \Delta[m] \XRA{F * \rho} C \times \Delta[k]
\]
defined at each $j \geq 0$ by
\[
\begin{aligned}
D_j \times \Delta[m]_j & \XRA{(F* \rho)_j} C_j \times \Delta[k]_j\\
(v,[j] \XRA{\theta} [m]) & \mapsto 
\Bigl(
F_j(v,\theta),\; [j] \XRA{\theta} [m] \XRA{\rho} [k]
\Bigr)
\end{aligned}
\]
Then $k$, $x$, $F$, any fixed $u \in [k+1]$ and $\rho$ determine a simplicial map
\begin{equation}
\label{G-rho-u}
G^{\rho,u} \DFAS H^{x,u} \circ (F * \rho) : D \times \Delta[m] \to C \times \Delta[k] \to \MC{S}
\end{equation}
For each $j \geq 0$, each element $h_j \LIST{t}{j}$ of $G^{\rho,u}_j(v,\theta) = \HMU xuj{F_j(v,\theta), \rho \theta}$ is sponsored by some $z \in \CPL(x,\mu; F_j(v,\theta), \del_{u} \rho\, \theta)$ where $\DMAP{\mu}{k}{j+1}$ is the complementary vertex function of $[j] \XRA{\theta} [m] \XRA{\rho} [k] \XRA{\del_u} [k+1]$.
\BS

The collection of sponsors of elements of $\HMU xuj{F_j(v,\theta), \rho \theta}$ for various $j$ and $(v,\theta)$ comprise a subcomplex of $C^{x,F}$ (page \pageref{C^x,F}) defined as follows.

\begin{newdef}
{\bf ($C^x(F,u,\rho)$)}
{\rm
\SL

Suppose $C$ and $D$ are a simplicial sets. Fix:
\[
\begin{aligned}
k \geq 0, & \quad x \in C_k\\
m \geq 0, & \quad \text{a simplicial map } F:D \times \Delta[m] \to C\\
u \in [k+1], & \quad \DMAP{\rho}{m}{k} \text{ non-decreasing}
\end{aligned}
\]
 Then for each $j \geq 0$, we define
\[
C^x(F,u,\rho)_j \DFAS \bigcup \CPL(x,\mu;F_j(v,\theta), \del_u \rho\, \theta)
\]
where the union is taken over all $(v,\theta) \in (D \times \Delta[m])_j$. 

That is, $z \in C^x_j$ is in $C^x(F,u,\rho)_j$ if the complementary subface in $z$ of $x$ is $F_j(v,\theta)$ for some $(v,\theta) \in (D \times \Delta[m])_j$ and the vertex function for $F_j(v,\theta)$ is $\del_u\, \rho \, \theta$.

We define
\[
C^x(F,u,\rho) \DFAS \bigcup_{j \geq 0} C^x(F,u,\rho)_j
\]
}

\BX
\end{newdef}

\begin{newlem}
{\rm
\SL

Given any $C,D,F, m,k,\rho,u$ and $x$ as in the definition above, then $C^x(F,u,\rho)$ is a subcomplex of $C^{x,F}$.
}
\end{newlem}

\Proof

Juxtaposing the definitions of $C^{x,F}_j$ (page \pageref{C^x,F}) and $C^x(F,u,\rho)$:
\[
\begin{aligned}
C^{x,F}_j & \DFAS \Bigl\{
\bigl(z,(v,\theta) \bigr) : z \in C^x_j, (v,\theta) \in (D \times \Delta[m])_j \text{ and } z \in \CPL(x;F_j(v,\theta)) \Bigr\}\\
C^x(F,u,\rho)_j & \DFAS \bigcup \CPL(x,\mu;F_j(v,\theta), \del_u \rho\, \theta)
\end{aligned}
\]
the map
\[
C^x(F,u,\rho)_j \to C^{x,F}_j, \quad z \mapsto (z,y)
\]
is monic where $y$ is the complementary subface of $x$ in $z$. That is, $y= \ F_j(v,\theta)$ for some $(v,\theta) \in (D \times \Delta[m])_j$.
\MS

Now if $z \in C^x(F,u,\rho)_j$ with $j>0$, with $z \in \CPL(x,\mu;F_j(v,\theta), \del_u \rho\, \theta)$ and $(v,\theta) \in (D \times \Delta[m])_j$, then for each $p \in [j]$, $j>0$:
\[
d^x_p(z) \in \CPL(x,\SIG_p \mu; F_{j-1}(d_p(v),\theta \del_p), \del_u \rho (\theta \del_p)) \SBS C^x(F,u,\rho)_{j-1}
\]
and, for any $j \geq 0$, 
\[
s^x_p(z) \in \CPL(x,\del_{p+1} \mu; F_{j+1}(s_p(v),\theta \SIG_p), \del_u \rho (\theta \SIG_p)) \SBS C^x(F,u,\rho)_{j+1}
\]
from which it follows that $C^x(F,u,\rho) \to C^{x,F}$ is a simplicial map.

\qed
\MS

\NI{\bf Notation:}

For $j \geq 0$, denote by $\MC{R}^{F,\rho}(x,u)_j$ the subset of $\MC{R}(x,u)_j$ consisting of $j$-simplices of $\MC{S}$ of the form $\HMU xuj{F_j(v,\theta,\rho \theta}$ for $(v,\theta) \in (D \times \Delta[k])_j$.

It follows from the definitions that $\MC{R}^{F,\rho}(x,u)$ is a subcomplex of $\MC{R}(x,u)$.

\BX
\BS

Now consider the special case where $C$ and $D$ are \NICOMP{n}{i}s and we consider $G^{\rho,u}$ defined in equation \eqref{G-rho-u} above where
\[
\begin{aligned}
m=k=n, & \quad x \in C_n\\
F:D \times \Delta[n] \to C & \quad \text{is a simplicial map}\\
\DMAP{\rho}{n}{n}& \quad \text{is non-decreasing}\\
G^{\rho,u} \text{ is} & \quad D \times \Delta[n] \XRA{F * \rho} C \times \Delta[n] \XRA{H^{x,u}} \MC{S}
\end{aligned}
\]
That is, for each $j \geq 0$ and each $(v,\theta) \in (D \times \Delta[n])_j$
\[
G^{\rho,u}(v,\theta) = H^{x,u}(F_j(v,\theta), \rho \theta)
\]
where each element of $G^{\rho,u}(v,\theta)$ is sponsored by some
\[
z \in \CPL(x,\mu;F_j(v,\theta), \del_u\, \rho \, \theta) \SBS C^{x,F}_j
\]

Now let $j \geq n+1$ and $(v,\theta) \in (D \times \Delta[n])_j$. We have for each $p \in [j]$
\[
d_p\Bigl( \HMU xuj{F_j(v,\theta),\del_u \rho\theta} \Bigr)
= \HMU xu{j-1}{F_{j-1}(d_p(v),\theta \del_p), \del_u \rho \theta \del_p}
\]
and Theorem \refpage{Rxu-is-composer} implies
\[
\Bigl(
\HMU xu{j-1}{F_{j-1}(d_0(v),\theta \del_0), \del_u \rho \theta \del_0} \DDD{,} \HMU xu{j-1}{F_{j-1}(d_j(v),\theta \del_j), \del_u \rho \theta \del_j}
\Bigr)
\]
is a $(j,i)$-composition in $\MC{S}$. Therefore $\MC{R}^{F,\rho}(x,u)$ is also an \NICOMP{n}{i}.

\subsection{Generalizing universality}
\label{singular-subsection}
Recall that if $F:D \to C$ is a functor of categories and $x \in \OB(C)$ then a ``universal map from $x$ to $F$'' consists of $y_0 \in \OB(D)$ and a map $u:x \to F_0(y_0)$ such that for all $y \in \OB(D)$ and all $v:x \to F_0(y)$ there exists a unique $f : y_0 \to y$ such that $v = F_1(f)u.$
\begin{center}
\begin{picture}(1,.5)
\put(0,0){\MB{F_0(y_0)}}
\put(1,0){\MB{F_0(y)}}
\put(.5,.5){\MB{x}}
\put(.4,.4){\vector(-1,-1){.25}}
\put(.6,.4){\vector(1,-1){.25}}
\put(.3,0){\vector(1,0){.4}}
\put(.25,.35){\FNS{\MB{u}}}
\put(.8,.35){\FNS{\MB{\ALL v}}}
\put(.5,.1){\FNS{\MB{\exists ! f}}}
\put(.5,-.1){\FNS{\MB{F_1(f)}}}
\end{picture}
\end{center}
\MS

Eventually, in section \refpage{universality in $C^{x,F}$} we will formulate a general version of this diagram.
\MS

The idea of universal map from $x$ to $F$ is an example of a general statement of the following form: 
\begin{quote}
A simplex $x$ in a simplicial set $C$ is ``universal'' for some type of configuration involving $x$ when, for each instance of such a configuration, there is a unique simplex involving $x$ in some specified way. 
\end{quote}
To make this vague statement precise, we have to say what is meant by a ``configuration'' and what ``involving $x$ in some specified way''means.

In the case of a universal map from $x$ to $F$, the simplicial set is the nerve of $x \downarrow F$
and the object (0-simplex) $x \XRA{u} F_0(y_0)$ has the property that for all ``configurations'' $x \XRA{v} F_0(y)$, there is a unique 1-simplex  $\tilde{f} \in (x \downarrow F)_1$ such that 
$\tilde{f} = u \XRA{F_1(f)} v$.

\BX
\MS

The following definition is one way to express the notion of uniqueness.

\begin{newdef}
{\bf ($\mu$-singular)}
\label{mu-singular}
\index{singular (with respect to ...)}
\index{$\mu$-singular}
\index{X-singular}
{\rm
\SL

Let $C$ be a simplicial set, $j, k \geq 0$, $x \in C_k$ and $\DMAP{\mu}{k}{j+1}$ be non-decreasing. 
We say {\bf $x$ is $\mu$-singular} 
if for all $y \in C_j$ there is exactly one $z \in C_{j+k+1}$ such that $\VL_z(x) = \mu\SH[k]$ and $d_{\mu\SH[k]}(z)=y$.

If $X \DFAS \mu\SH[k] \SBS [j+k+1]$ then we'll also say that {\bf $x$ is $X$-singular}.
}

\BX
\end{newdef}

\NI Notes: With $C$, $x$, $j,k$ and $\mu$ as above then $x$ being $\mu$-singular says that for all $y \in C_j$, $|\CPL(x,\mu;y,\LAM)|=1$ where $\DMAP{\LAM}{j}{k+1}$ is the complementary vertex function to $\mu$. In that case, the unique $z \in \CPL(x,\mu;y,\LAM)$ has $\VL_z(y) = \LAM\SH[j] = [j+k+1]-\mu\SH[k]$.
\MS

The definition is motivated by the following example.
\MS

\begin{example}
\label{initial-terminal-example}
{\rm
Given $C$, $j=k=0$, $x \in C_0$ and $\DMAP{\mu}{0}{1}$ then there are just two cases for $x$ being $\mu$-singular.
\MS

\NI If $\mu(0)=0$ then for all $y\in C_0$ there is a unique $z \in C_1$ such that $\VL_z(x)=\{0\}$. That is, for all $y\in C_0$ there is a unique 1-simplex $z=x \to y$. If $C$ is the nerve of a category, then $x$ is the initial object.
\MS

\NI If $\mu(0)=1$ then for all $y\in C_0$ there is a unique $z \in C_1$ such that $\VL_z(x)=\{1\}$. That is, for all $y\in C_0$ there is a unique 1-simplex $z=y \to x$. If $C$ is the nerve of a category, then $x$ is the terminal object.
}
\end{example}
\MS

Now suppose $C$ is an \NICOMP{n}{i} where $n \geq 2$. The question is: What is a (or the) useful version of a simplex of $C$ being ``singular'' in the context of the $(n,i)$-composor structure? The initial/terminal object example above suggests looking at cases where $j+k+1 = n$. That is, the unique $z \in \CPL(x,\mu;y,\LAM) \SBS C_n$ should be an $(n,i)$-factor. The two extreme cases are $k=n-1$ and $k=0$ which specialize to the motivating example above.
\MS

\NI {\bf Case:} $k=n-1,\, j=0$, $\DMAP{\mu}{n-1}{1}$ and $\DMAP{\LAM}{0}{n}$

Then 
\[
\begin{aligned}
\DMAP{\LAM\SH}{0}{n}\qquad & \text{Let }p \DFAS \LAM\SH(0).\\
\DMAP{\mu\SH}{n-1}{n} \qquad & \mu\SH=\del_p
\end{aligned}
\]

To say that the $(n,i)$-object $x \in C_{n-1}$ is $\mu$-singular means that for all $y \in C_0$ there is a unique $(n,i)$-factor $z \in C_n$ such that $d_p(z)=x$ and $\VL_z(y)=\{p\}$.
\MS

\NI {\bf Case:} $k=0,\, j=n-1$, $\DMAP{\mu}{0}{n}$, $\DMAP{\LAM}{n-1}{1}$

Then
\[
\begin{aligned}
\DMAP{\mu\SH}{0}{n}\qquad & \text{Let } p \DFAS \mu\SH(0).\\
\DMAP{\LAM\SH}{n-1}{n} \qquad & \LAM\SH = \del_p
\end{aligned}
\]
Then $x \in C_0$ is $\mu$-singular iff for all $(n,i)$-objects $y \in C_{n-1}$ there exists a unique $(n,i)$-factor $z \in C_n$ such that $\VL_z(x) = \{ p \}$ and $d_p(z)=y$.

\subsection{Trios}

\setlength{\unitlength}{1in}

As in the previous section we will use the definition of universal map from $x \in \OB(C)$ to the functor $F:D \to C$ as a motivating example.
\begin{center}
\begin{pictr}{1}{.5}{-1}{.2}
\label{univ-map-diag}
\put(0,0){\MB{F_0(y_0)}}
\put(1,0){\MB{F_0(y)}}
\put(.5,.5){\MB{x}}
\put(.4,.4){\vector(-1,-1){.25}}
\put(.6,.4){\vector(1,-1){.25}}
\put(.3,0){\vector(1,0){.4}}
\put(.25,.35){\FNS{\MB{u}}}
\put(.8,.35){\FNS{\MB{\ALL v}}}
\put(.5,.1){\FNS{\MB{\exists ! f}}}
\put(.5,-.1){\FNS{\MB{F_1(f)}}}
\end{pictr}
\end{center}
\BS

\NI The aim is to generalize this kind of diagram to the following situation: 
\begin{itemize}
\item 
$C$ and $D$ are simplicial sets and $F:D \to C$ is a simplicial map.

\item
$k,\, j_0,\, j \geq 0$ and $x \in C_k$, $y_0 \in D_{j_0}$, $y \in D_j$.

\item
The commutative triangle above, a 2-simplex in the nerve of the category $C$, is replaced by an $m$-simplex $w \in C_m$ where $m = j_0+j+k+2$ such that
the vertex index lists of $x, y_0$ and $y$ comprise a specified partition of $[m]$.

\end{itemize}

This configuration of $x, F_{j_0}(y_0), F_{j}(y)$ and their respective vertex index lists will be called a ``trio'' (the exact definition is below starting on page \pageref{combinatorial-trio}). The discussion of the generalization begins on page \pageref{universality in $C^{x,F}$}.

After a preparatory subsection concerning ``partial simplices'', we will develop the trio idea in two stages:
\begin{itemize}

\item 
Given $m\geq 2$, we will consider a partition $X \cup A \cup A'$ of $[m]$ and associated inclusion maps. (Section \refpage{combtrio-section})

\item
Specialize to $F$ being the identity simplicial map.
Given $j,j', k \geq 0$, $m = j+j'+k+2$, a simplicial set $C$, we will consider simplices $x \in C_k$, $y \in C_j$ and $y' \in C_{j'}$ and $z \in C_m$ (if any) such that $\VL_z(x)=X, \VL_z(y)=A$ and $\VL_z(y')=A'$ with $X \cup A \cup A' = [m]$ a disjoint union. (Definition \refpage{trio-def}).

\end{itemize}

\subsubsection{Partial simplices}

In this subsection, we generalize the open $A$-horn idea to define what we will call a ``partial simplex of dimension $m$''. (Definition \refpage{partial-m-simplex-def} below).
\BS

By ``subface of $z \in C_m$'', we mean any simplex of the form 
\[
d_{n_0} d_{n_1} \DDD{} d_{n_j}(z), \quad j \leq m-1
\]
The face identities imply that every such subface has the form
\[
d_{p_0} d_{p_1} \DDD{} d_{p_j}(z), \quad \text{where } p_0 \DDD{<} p_j
\]
\NI {\bf Notation:} \index{$d_A$} \index{$\del_A$}
In order to specify subfaces of a simplex $z \in C_m$ more concisely we will write $d_A(z)$ to mean $d_{p_0} d_{p_1} \DDD{} d_{p_j}(z)$ where $A = \SETT{p_0 \DDD{<} p_j}$.
\MS

\NI Dually, we will write $\del_A$ for $\del_{p_j} \DDD{} \del_{p_0}$.

\BX
\BS

We will consider certain families of simplices using the following observation as motivation. Given $m>0$ and $z \in C_m$, then a non-empty family $\mathcal{F}$ of proper non-empty subsets of $[m]$ determines a family $\SETT{d_A(z): A \in \mathcal{F}}$ of subfaces of $z$. To indicate {\em where} each subface $d_A(z)$ occurs in $z$ (i.e. that its vertex index set is $[m]-A$) we consider the set of pairs 
\[
z^\mathcal{F} = \Bigl\{(d_A(z),[m]-A): A \in \mathcal{F}\Bigr\}
\] 
Here, $d_A(z)$ is a subface of $z$ and $\VL_z(d_A(z)) = [m]-A$.
\MS

There are specific compatibility relationships among the $d_A(z)$ expressed by the relevant face identities. Given two subfaces $d_A(z)$ and $d_{A'}(z)$ that compatibility takes the following form: whenever you have $B \SBS [m - |A|]$ and $B' \SBS [m-|A'|]$ such that $d_B d_A = d_{B'} d_{A'}$ is a face identity then $d_B d_A(z) = d_{B'} d_{A'}(z)$. These face equations occur precisely when $d_A(z)$ and $d_{A'}(z)$ have intersecting vertex index sets. The case of $B$ or $B'$ being empty is allowed.

Briefly, a ``partial $m$-simplex'' is a family like $z^\mathcal{F}$ above but without mentioning or requiring a specific $z \in C_m$. The precise definition is below, based on a careful representation of the required face equations.
\MS

By way of preparation, we gather into one lemma some elementary facts about $\Delta$ and, given any $m > 0$, the simplicial set $\Delta[m]$.

\begin{newlem} 
{\bf (Face Identity lemma)}
\label{face-identity-lemma}
\index{face identity lemma}
{\rm
\SL

Let $m>0$.
\begin{enumerate}
\item 
Given any \SI\ function $\DMAP{g}{j}{m}$ then $g = \del_{[m]-g[j]}$ and so $g = d_{[m]-g[j]}(1_{[m]}) \in \Delta[m]_j$. 

\item
Given two \SI\ functions $[k] \XRA{f} [j] \XRA{g} [m]$ then $f=\del_{[j]-f[k]}$, $g = \del_{[m]-g[j]}$ and $gf = \del_{[m]-gf[k]}$. Therefore
\[
gf = d_{[j]-f[k]}\; d_{[m]-g[j]} (1_{[m]}) = d_{[m]-gf[k]}(1_{[m]})
\]
from which we get the face identity
\[
d_{[j]-f[k]}\; d_{[m]-g[j]} = d_{[m]-gf[k]}
\]

\item
Suppose $A_1$ and $A_2$ are non-empty subsets of $[m]$ such that
$A_1 \cup A_2 \subsetneq [m]$.
Then $A_1,A_2$ determine a commutative diagram 
\begin{center}
\setlength{\unitlength}{1in}
\begin{pictr}{2}{.5}{-.7}{.25}
\label{A1-A2-diag}
\put(1,0){\MB{[k]}}
\put(0,0){\MB{[j_1]}}
\put(2,0){\MB{[j_2]}}
\put(1,.5){\MB{[m]}}
\put(.2,.1){\vector(2,1){.6}} 
\put(.8,0){\vector(-1,0){.6}} 
\put(1.2,0){\vector(1,0){.6}} 
\put(1.8,.1){\vector(-2,1){.6}} 
\put(1,.15){\vector(0,1){.2}} 
\put(.5,.07){\MBS{\alpha_1}}
\put(1.5,.07){\MBS{\alpha_2}}
\put(.4,.3){\MBS{\beta_1}}
\put(1.6,.3){\MBS{\beta_2}}
\put(.925,.2){\MBS{\gamma}}
\end{pictr}
\end{center}
where
\[
j_1=m-|A_1|, \quad j_2=m-|A_2|\, \quad k=m-|A_1 \cup A_2|
\]
and the strictly increasing functions $\beta_1, \beta_2$ and $\gamma$ are defined by
\[
\begin{aligned}
\DMAP{\beta_1}{j_1}{m} & \qquad \text{such that }\beta_1[j_1] = [m]-A_1\\
\DMAP{\beta_2}{j_2}{m} & \qquad \text{such that }\beta_2[j_2] = [m]-A_2\\
\DMAP{\gamma}{k}{m} & \qquad \text{such that }\gamma[k] = [m]-(A_1 \cup A_2)
\end{aligned}
\]
and the strictly increasing functions $\alpha_1$ and $\alpha_2$ are defined by requiring
\[
\alpha_1[k] = \beta_1^{-1}\Bigl([m]-(A_1 \cup A_2\Bigr), \quad \alpha_2[k] = \beta_2^{-1}\Bigl([m]-(A_1 \cup A_2\Bigr)
\]
Then

\[
\begin{aligned}
d_{[j_1]-\alpha_1[k]} d_{[m]-\beta_1[j_1]} &=& d_{[m]-\gamma[k]} &=& 
d_{[j_2]-\alpha_2[k]} d_{[m]-\beta_2[j_2]} \\
d_{[j_1]-\alpha_1[k]}d_{A_1} &=& d_{A_1 \cup A_2} &=& d_{[j_2]-\alpha_2[k]} d_{A_2}
\end{aligned}
\]
\end{enumerate}
}
\end{newlem}
\Proof

Item (1) is by definition of ``$\del$'' and the simplicial structure of $\Delta[m]$. The remaining items follow directly from this.

\qed
\MS

\begin{newdef}
{\bf (partial $m$-simplex)}
\label{partial-m-simplex-def}
\index{partial $m$-simplex}
{\rm
\SL

Let $C$ be a simplicial set, $m>0$ and $\mathcal{F} \subsetneq 2^{[m]}$ a non-empty family of non-empty proper subsets of $[m]$. Then a family
\[
z^\mathcal{F} = \SETT{(y_A,[m]-A): A \in \mathcal{F}, y_A \in C_{m-|A|}}
\]
is a {\bf partial $m$-simplex} if
for all $A_1, A_2 \in \mathcal{F}$ such that $[m]-(A_1 \cup A_2) \neq \MT$, 
with the corresponding \SI\ functions (see the diagram below)
\[
\begin{aligned}
\DMAP{\gamma}{k}{m} \quad \text{defined by} \quad & \gamma[k] = [m]-(A_1 \cup A_2)\\
\DMAP{\beta_1}{j_1}{m} \quad \text{defined by} \quad & \beta_1[j_1] = [m]-A_1\\
\DMAP{\beta_2}{j_2}{m} \quad \text{defined by} \quad & \beta_2[j_2] = [m]-A_2\\
\DMAP{\alpha_1}{k}{j_1} \quad \text{defined by} \quad & \alpha_1[k] = \beta_1^{-1} \gamma[k]\\
\DMAP{\alpha_2}{k}{j_2} \quad \text{defined by} \quad & \alpha_2[k] = \beta_2^{-1} \gamma[k]\\
\end{aligned}
\]
comprising the commutative diagram:

\begin{center}
\begin{pictr}{2}{.6}{-.7}{.25}
\label{par-simp-diag-frag}
\put(0,0){\MB{[j_1]}}
\put(1,0){\MB{[k]}}
\put(2,0){\MB{[j_2]}}
\put(1,.5){\MB{[m]}}
\put(.2,.1){\vector(2,1){.65}}
\put(.5,.37){\MBS{\beta_1}}
\put(1.8,.1){\vector(-2,1){.65}}
\put(1.5,.37){\MBS{\beta_2}}
\put(1,.15){\vector(0,1){.2}}
\put(1.1,.25){\MBS{\gamma}}
\put(.75,0){\vector(-1,0){.5}}
\put(.6,.1){\MBS{\alpha_1}}
\put(1.25,0){\vector(1,0){.5}}
\put(1.4,.1){\MBS{\alpha_2}}
\end{pictr}
\end{center}
 
then $y_{A_1}$ and $y_{A_2}$ satisfy the face condition
\[
d_{[j_1]-\alpha_1[k]}(y_{A_1}) = d_{[j_2]-\alpha_2[k]}(y_{A_2})
\]
\BX
}
\end{newdef}

\begin{newdef}
{\bf (Sponsor)}
\label{sponsor-partial-simplex-def}
\index{sponsor of a partial simplex}
{\rm
\SL

A {\bf sponsor} of a partial $m$-simplex 
\[
z^\mathcal{F} = \SETT{(y_A,[m]-A): A \in \mathcal{F}}
\]
is any $z \in C_m$ such that for all $A \in \mathcal{F}$, $d_A(z) = y_A$.
}

\BX
\end{newdef}

\NI That is: if $(y_A,[m]-A) \in z^\MC{F}$ then $[m]-A$ is the vertex list of $y_A$ in any sponsor of $z^\MC{F}$. Referring back to diagram \eqref{par-simp-diag-frag} above, if $z$ sponsors $z^\MC{F}$ and $A_1, A_2 \in \MC{F}$ such that $A_1 \cup A_2 \subsetneq [m]$ then
\[
\begin{aligned}
d_{[m]-\gamma[k]}(z) & = d_{A_1 \cup A_2}(z) \\
& = d_{[j_1]-\alpha_1[k]}\, d_{[m]-\beta_1[j_1]}(z) \\
& = d_{[j_1]-\alpha_1[k]} d_{A_1}(z) = d_{[j_1]-\alpha_1[k]}(y_{A_1})
\end{aligned}
\]
So we have
\[
d_{[j_1]-\alpha_1[k]}(y_{A_1}) = d_{A_1 \cup A_2}(z) = d_{[j_2]-\alpha_2[k]}(y_{A_2})
\]
\BS

\begin{example}
{\rm
Let $m=6$, $A_1 = \SETT{0,1,4}, A_2 = \SETT{1,5}$. We will describe a partial 6-simplex consisting of:

$\bullet$\; a 3-simplex whose vertex index list is 
$\SETT{2,3,5,6} = [6]-A_1$

$\bullet$\; a 4-simplex whose vertex index list is $\SETT{0,2,3,4,6} = [6]-A_2$

Since these two vertex index lists intersect then there will be a face-compatibility requirement. In the notation of the definition above:

\[
\begin{aligned}
{[6]-A_1} = \SETT{2,3,5,6}, & \quad j_1 = 3\\
[6]-A_2 = \SETT{0,2,3,4,6}, & \quad j_2 = 4\\
[6]-(A_1 \cup A_2) = \SETT{2,3,6}, & \quad k=2
\end{aligned}
\]
Then 
\[
\begin{aligned}
\DMAP{\beta_1 = (2,3,5,6)}{3}{6}\\
\DMAP{\beta_2 = (0,2,3,4,6)}{4}{6}\\
\DMAP{\gamma = (2,3,6)}{2}{6}
\end{aligned}
\]
and in the following commutative diagram of \SI\ functions
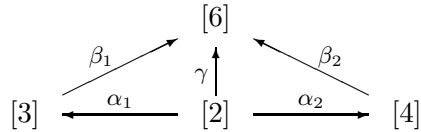
\begin{center}
\begin{picture}(2,.55)
\put(0,0){\MB{[3]}}
\put(1,0){\MB{[2]}}
\put(2,0){\MB{[4]}}
\put(1,.5){\MB{[6]}}
\put(.8,0){\vector(-1,0){.6}}
\put(1.2,0){\vector(1,0){.6}}
\put(1,.1){\vector(0,1){.25}}
\put(.2,.1){\vector(2,1){.6}}
\put(1.8,.1){\vector(-2,1){.6}}
\put(.5,.07){\MBS{\alpha_1}}
\put(1.5,.07){\MBS{\alpha_2}}
\put(.4,.3){\MBS{\beta_1}}
\put(1.6,.3){\MBS{\beta_2}}
\put(.925,.2){\MBS{\gamma}}
\end{picture}
\end{center}
\MS

\NI $\alpha_1 = (0,1,3)$ and $\alpha_2 = (1,2,4)$.

Now suppose $y_1 \in C_3$ and $y_2 \in C_4$. Then to say that
\[
\Bigl\{(y_1,[6]-A_1),(y_2,[6]-A_2)\Bigr\}
\]
is a partial 6-simplex means that the face equation $d_{[3]-\alpha_1[2]}(y_1) = d_{[4]-\alpha_2[2]}(y_2)$ holds. Since $[3]-\alpha_1[2] = \{2\}$ and $[4]-\alpha_2[2] = \{0,3\}$ then that equation is
\[
d_2(y_1) = d_0 d_3(y_2)
\]
representing the subface whose vertex index list is $\SETT{2,3,6}$.
Here are the vertex index lists in the partial $6$-simplex.
\[
\begin{array}{|r||c|c|c|c|c|c|c|c|}
\hline
 & 0 & 1 & 2 & 3 & 4 & 5 & 6 &\\
\hline
\VL(y_1) &  & & 0 & 1 & & 2 & 3 & \SETT{2,3,5,6} \\
\hline
\VL(y_2) & 0 & & 1 & 2 & 3 & & 4 & \SETT{0,2,3,4,6}  \\
\hline
\VL(d_2(y_1)) & & & 0 & 1 & & &2& \SETT{2,3,6}\\
\hline
\VL(d_0 d_3(y_2)) & & & 0 & 1 & & & 2 & \SETT{2,3,6} \\
\hline
\end{array}
\]
This can be visualized showing just the outline edges of $y_1$ (solid) and $y_2$ (dotted) by:
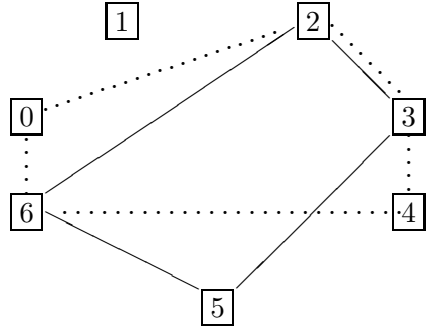
\begin{center}
\begin{picture}(3,1.7)
\put(.5,1){\MB{\fbox{0}}}
\put(1,1.5){\MB{\fbox{1}}}
\put(2,1.5){\MB{\fbox{2}}}
\put(2.5,1){\MB{\fbox{3}}}
\put(2.5,.5){\MB{\fbox{4}}}
\put(1.5,0){\MB{\fbox{5}}}
\put(.5,.5){\MB{\fbox{6}}}
\put(2.1,1.4){\line(1,-1){.3}} 
\put(2.4,.9){\line(-1,-1){.8}} 
\put(1.4,.1){\line(-2,1){.8}} 
\put(.6,.6){\line(3,2){1.3}} 
\multiput(.6,1.04)(.06,.02){21}{\makebox(0,0){.}} 
\multiput(2.1,1.48)(.05,-.05){8}{\makebox(0,0){.}} 
\multiput(2.5,.9)(0,-.07){5}{\makebox(0,0){.}} 
\multiput(2.45,.5)(-.07,0){26}{\makebox(0,0){.}} 
\multiput(.5,.6)(0,.07){5}{\makebox(0,0){.}} 
\end{picture}
\end{center}

A sponsor of $\SETT{(y_1,[6]-A_1),(y_2,[6]-A_2)}$ would be any $z \in C_6$ such that $d_{A_1}(z) = d_0 d_1 d_4(z) = y_1$ and $d_{A_2}(z) = d_1 d_5(z) = y_2$.
\MB{ }
}
\end{example}
\MS

\begin{example}
{\rm
Suppose $m>1$, $\MT \neq B \SBS[m]$ and $[m]-B$ has at least two elements. Consider an open $B$-horn
\[
\Bigl\{y_p : p \in [m]-B, y_p \in C_{m-1} \text{ and } \ALL p<q \text{ in } [m]-B, \; d_p(y_q) = d_{q-1}(y_p)\Bigr\}
\]
In order to show that an open $B$-horn is a partial $m$-simplex, it suffices to show that 
\[
\SETT{ (y_p,[m]-\{p\}),(y_q,[m]-\{q\}) }
\]
is a partial $m$-simplex for each $p<q$ in $[m]-B$.

In the notation used in the definition of partial $m$-simplex, $A_p \DFAS \{p\}$ and $A_q \DFAS \{q\}$. The corresponding commutative diagram of \SI\ functions is
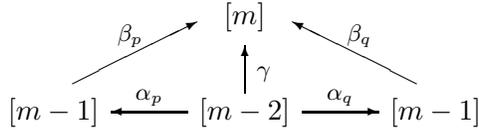
\begin{center}
\setlength{\unitlength}{1in}
\begin{picture}(2,.55)
\put(0,0){\MB{[m-1]}}
\put(1,0){\MB{[m-2]}}
\put(2,0){\MB{[m-1]}}
\put(1,.5){\MB{[m]}}
\put(.7,0){\vector(-1,0){.4}}
\put(1.3,0){\vector(1,0){.4}}
\put(1,.1){\vector(0,1){.25}}
\put(.1,.15){\vector(2,1){.65}}
\put(1.9,.15){\vector(-2,1){.65}}
\put(.4,.4){\FNS{\MB{\beta_p}}}
\put(1.6,.4){\FNS{\MB{\beta_q}}}
\put(.5,.075){\FNS{\MB{\alpha_p}}}
\put(1.5,.075){\FNS{\MB{\alpha_q}}}
\put(1.1,.2){\FNS{\MB{\gamma}}}
\end{picture}
\end{center}
where the maps are defined by:
\[
\begin{aligned}
\beta_p[m-1] = [m]-\{p\}, & \quad  \iff \DMAP{\beta_p = \del_p}{m-1}{m}\\
\beta_q[m-1] = [m]-\{q\}, & \quad \iff \DMAP{\beta_q = \del_q}{m-1}{m}\\
\gamma[m-2] = [m]-\{p,q\},& \quad \iff \DMAP{\gamma = \del_q \del_p = \del_p \del_{q-1}}{m-2}{m}
\end{aligned}
\]
which implies $\alpha_p = \del_{q-1}$ and $\alpha_q = \del_p$. 

The required face equation is $d_{[m-1]-\del_{q-1}[m-2]}(y_p) \ISIT d_{[m-1]-\del_p[m-2]}(y_q)$. This holds because
\[
d_{[m-1]-\del_{q-1}[m-2]} = d_{q-1}, \quad d_{[m-1]-\del_p[m-2]} = d_p
\]
and $d_{q-1}(y_p) = d_p(y_q)$.
\MB{ }
}
\end{example}
\MS

\begin{example}
{\rm
Suppose $x \in C_k$, $y \in C_j$ and $\DMAP{\mu}{k}{j+1}$ and $\DMAP{\LAM}{j}{k+1}$ are complementary vertex functions. Then
\[
\Bigl\{(x,\mu\SH[k]),(y,\LAM\SH[j])\Bigr\}
\]
is a partial $(j+k+1)$-simplex trivially since $\mu\SH[k] \cap \LAM\SH[j]= \MT$ and there is no face equation for $x$ and $y$. The set of sponsors of $\SETT{(x,\mu\SH[k]),(y,\LAM\SH[j])}$ is $\CPL(x,\mu;y,\LAM)$.
}
\end{example}

\subsubsection{Trios, combinatorial trios and sponsors}
\label{combtrio-section}

We are interested  in considering certain $m$-simplices with three specified subfaces whose vertex index lists partition $[m]$. We'll start by developing some notation and facts concerning such partitions.
\MS

\NI{\bf Notations for subsets of $[m]$:}
Any non-empty subset $B \SBS [m]$ with $k+1$ elements corresponds to a \SI\ function we will denote by $\HT{B} :[k] \to [m]$. That is
\[
B = \SETT{\HT{B}(q):a \in [k]} \text{ denoted } \HT{B}[k]
\]
\BX

\begin{newdef}
{\bf (Combinatorial-trios, Trios and Sponsors)}
\label{combinatorial-trio}
\index{combinatorial trio}
\index{comb-trio}
\index{trio (combinatorial)}
\label{trio-def} 
\index{trio} 
\index{trio sponsor} 
\index{sponsor of a trio}
{\rm
\SL
\MS

\NI {\bf 1. Combinatorial-trio:} 

Suppose $m \geq 2$, $X,A,A' \SBS [m]$ are non-empty subsets and 
\[
X \cup A \cup A' =[m]
\]
is a disjoint union where $|X|=k+1, |A|=j+1$ and $|A'|=j'+1$. Then we will refer to the partition $\SETT{X,A,A'}$ as a {\bf $(k,j,j')$-combinatorial trio of dimension $m$} (abbreviation: {\bf ``comb-trio''}). In this situation, $m = j+j'+k+2$.
\BS

\NI {\bf 2. Trio:} 

Given a $(k,j,j')$-comb-trio $\SETT{X,A,A'}$ of dimension $m$  and a simplicial set $C$ with $x \in C_k$, $y \in C_j$ and $y' \in C_{j'}$ then we will refer to the partial $m$-simplex 
\[
\Bigl\{(x,X), (y,A), (y',A')\Bigr\}
\]
as a {\bf $(k,j,j')$-trio (in $C$) of dimension $m$}, where $m=j+j'+k+2$.
Note that there are no face equations relating $x,y$ and $y'$ because $X,A,A'$ are pairwise disjoint.
\BS

\NI {\bf 3. Sponsor of a trio:}

 Given a trio $\SETT{(x,X), (y,A), (y',A')}$ of dimension $m$ in the simplicial set $C$, a {\bf sponsor} of it is any ${w} \in C_m$ such that:
\[
\begin{aligned}
d_{[m]-X}(w) = d_{A \cup A'}(w) = x \quad & \iff \quad \VL_w(x)=X\\
d_{[m]-A}(w) = d_{X \cup A'}(w) = y \quad & \iff \quad \VL_w(y)=A\\
d_{[m]-A'}(w) = d_{X \cup A}(w) = y' \quad & \iff \quad \VL_w(y')=A'
\end{aligned}
\]
}

\BX
\end{newdef}
\MS

\NI {\bf Complementary vertex function pairs associated to a comb-trio:}
\MS

Given a $(k,j,j')$-comb-trio $\SETT{X,A,A'}$ of dimension $m = k+j+j'+2$, there are several associated pairs of complementary vertex functions. 
\MS

First, $X \cap A = \MT$ and $|X \cup A|=j+k+1$. Therefore there are \SI\ functions $\mu\SH$ and $\LAM\SH$ making diagram \eqref{mu-lam-from-comb-trio} commute

\begin{center}
\begin{pictr}{2}{1.1}{-.7}{.5}
\label{mu-lam-from-comb-trio}
\put(0,0){\MB{[k]}}
\put(1,0){\MB{[j+k+1]}}
\put(2,0){\MB{[j]}}
\put(1,1){\MB{[m]}}
\put(.15,0){\vector(1,0){.45}}
\put(.4,.1){\MBS{\mu\SH}}
\put(1.85,0){\vector(-1,0){.45}}
\put(1.6,.1){\MBS{\LAM\SH}}
\put(1,.15){\vector(0,1){.7}}
\put(.1,.1){\vector(1,1){.75}}
\put(1.9,.1){\vector(-1,1){.75}}
\put(.4,.6){\MBS{\HT{X}}}
\put(1.6,.6){\MBS{\HT{A}}}
\put(.8,.4){\MBS{\HT{X\cup A}}}
\end{pictr}
\end{center}

\NI and $[k] \XRA{\mu} [j+1]$ and $[j] \XRA{\LAM} [k+1]$ comprise a complementary vertex function pair.
\MS

Applying the same observation to the disjoint pairs $X',A$ and $A,A'$ we get two other pairs of complementary vertex functions
\[
\DMAP{\mu'}{k}{j'+1}, \quad \DMAP{\LAM'}{j'}{k+1}
\]
and
\[
\DMAP{\gamma}{j}{j'+1}, \quad \DMAP{\gamma'}{j'}{j+1}
\]

All this is summarized by the following commutative diagram of strictly increasing functions.

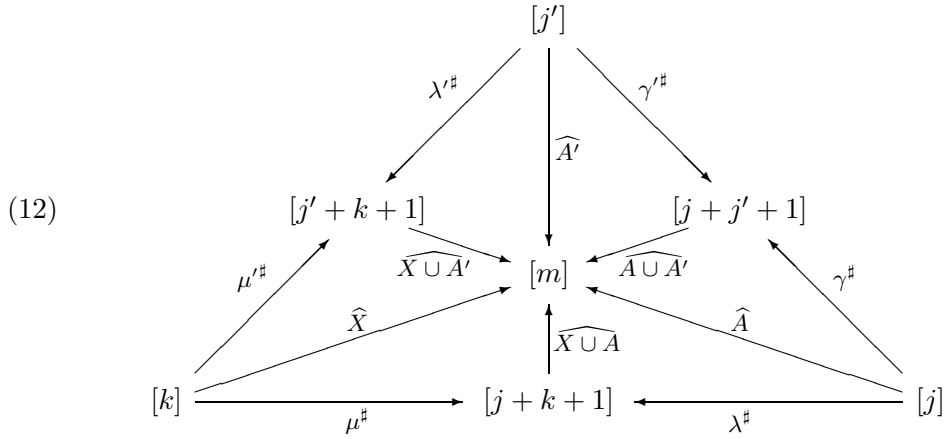
\begin{center}
\begin{pictr}{4}{2}{-.7}{1}
\label{comb-trio-diag-1}
\MBP{0}{0}{[k]}
\MBP{2}{0}{[j+k+1]}
\MBP{4}{0}{[j]}
\MBP{1}{1}{[j'+k+1]}
\MBP{2}{2}{[j']}
\MBP{3}{1}{[j+j'+1]}
\MBP{2}{.667}{[m]}
\put(.15,.15){\vector(1,1){.7}} 
\MBPS{.45}{.65}{{\mu'}\SH}
\put(.15,.05){\vector(3,1){1.65}} 
\MBPS{1}{.43}{\HT{X}}
\put(.15,0){\vector(1,0){1.4}}
\MBPS{1}{-.1}{\mu\SH}
\put(3.85,.15){\vector(-1,1){.7}}
\MBPS{3.55}{.65}{\gamma\SH} 
\put(3.85,.05){\vector(-3,1){1.65}} 
\MBPS{3}{.43}{\HT{A}}
\put(3.85,0){\vector(-1,0){1.4}}
\MBPS{3}{-.1}{\LAM\SH}
\put(1.85,1.85){\vector(-1,-1){.7}}
\MBPS{1.45}{1.65}{{\LAM'}\SH}
\put(2.15,1.85){\vector(1,-1){.7}}
\MBPS{2.55}{1.65}{{\gamma'}\SH}
\put(2,1.85){\vector(0,-1){1.03}}
\MBPS{2.1}{1.33}{\HT{A'}}
\put(2,.15){\vector(0,1){.36}}
\MBPS{2.2}{.33}{\HT{X \cup A}}
\put(1.27,.91){\vector(3,-1){.53}}
\MBPS{1.4}{.725}{\HT{X \cup A'}}
\put(2.73,.91){\vector(-3,-1){.53}}
\MBPS{2.55}{.725}{\HT{A \cup A'}}
\end{pictr}
\end{center}
\BS

\NI We will say that the complementary vertex function pairs $(\mu,\LAM), (\mu',\LAM')$ and $(\gamma, \gamma')$ ``belong to'' the comb-trio.
\MS

Note that this diagram is determined completely by, for example, just $\HT{X \cup A}$ and $\mu\SH$ because $A$ and $A'$ determine $\gamma\SH$ and $\gamma'{\SH}$, $X$ and $A'$ determine $\mu'{\SH}$ and $\LAM'{\SH}$ and:
\[
\begin{aligned}
\LAM\SH[j]&=[j+k+1]-\mu\SH[k]\\
X=\HT{X}[k] &=\HT{X \cup A}\circ \mu\SH[k]\\
A=\HT{A}[j]&=\HT{X \cup A}\circ \LAM\SH[j]\\
A' =\HT{A'}[j'] &=[m]-\HT{X \cup A}[j+k+1]\\
X \cup A'&=\HT{X \cup A}[j'+k+1]=[m]-A
\end{aligned}
\]
%
%
%
%
%
%
%
%
\BX
\MS

\begin{example}
{\rm
Let $k=2,j=3,j'=4$ and suppose $X = \SETT{1,5,10}, A= \SETT{0,2,6,9}$ and $A'=\SETT{3,4,7,8,11}$. Then $m=11$ and we have
\[
\begin{aligned}
\HT{X} &= \DMAP{(1,5,10)}{2}{11} \\
\HT{A} &= \DMAP{(0,2,6,9)}{3}{11}\\
\HT{A'} &= \DMAP{(3,4,7,8,11)}{4}{11}\\
\HT{X \cup A} &= (0,1,2,5,6,9,10):[6] \to [11]\\
\HT{X \cup A'} &= (1,3,4,5,7,8,10,11):[7] \to [11]\\
\HT{A \cup A'} &= (0,2,3,4,6,7,8,9,11):[8] \to [11]
\end{aligned}
\]
Since $\HT{X \cup A}^{-1}(A) = \SETT{0,2,4,5}$ then $\LAM\SH = \DMAP{(0,2,4,5)}{3}{6}$. Similar calculations show that 
\[
\begin{aligned}
\mu\SH = \DMAP{(1,3,6)}{2}{6} \\
{\mu'}\SH = \DMAP{(0,3,6)}{2}{7}\\
{\LAM'}\SH = \DMAP{(1,2,4,5,7)}{4}{7}\\
\gamma\SH = \DMAP{(0,1,4,7)}{3}{8}\\
{\gamma'}\SH = \DMAP{(2,3,5,6,8)}{4}{8}
\end{aligned}
\]
making for the an instance of diagram \eqref{comb-trio-diag-1}:

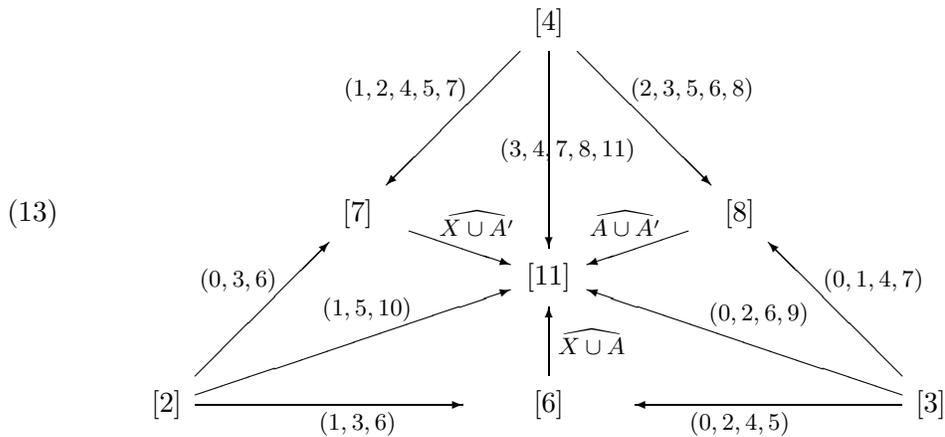
\begin{center}
\begin{pictr}{4}{2}{-.7}{1}
\MBP{0}{0}{[2]}
\MBP{2}{0}{[6]}
\MBP{4}{0}{[3]}
\MBP{1}{1}{[7]}
\MBP{2}{2}{[4]}
\MBP{3}{1}{[8]}
\MBP{2}{.667}{[11]}
\put(.15,.15){\vector(1,1){.7}} 
\MBPS{.375}{.65}{(0,3,6)}
\put(.15,.05){\vector(3,1){1.65}} 
\MBPS{1.05}{.5}{(1,5,10)}
\put(.15,0){\vector(1,0){1.4}}
\MBPS{1}{-.1}{(1,3,6)}
\put(3.85,.15){\vector(-1,1){.7}}
\MBPS{3.7}{.65}{(0,1,4,7)} 
\put(3.85,.05){\vector(-3,1){1.65}} 
\MBPS{3.1}{.47}{(0,2,6,9)}
\put(3.85,0){\vector(-1,0){1.4}}
\MBPS{3}{-.1}{(0,2,4,5)}
\put(1.85,1.85){\vector(-1,-1){.7}}
\MBPS{1.25}{1.65}{(1,2,4,5,7)}
\put(2.15,1.85){\vector(1,-1){.7}}
\MBPS{2.75}{1.65}{(2,3,5,6,8)}
\put(2,1.85){\vector(0,-1){1.03}}
\MBPS{2.1}{1.33}{(3,4,7,8,11)}
\put(2,.15){\vector(0,1){.36}}
\MBPS{2.225}{.33}{\HT{X \cup A}}
\put(1.27,.91){\vector(3,-1){.53}}
\MBPS{1.625}{.95}{\HT{X \cup A'}}
\put(2.73,.91){\vector(-3,-1){.53}}
\MBPS{2.4}{.95}{\HT{A \cup A'}}
\end{pictr}
\end{center}
\MB{ }
}
\end{example}
\BS

\NI {\bf Schematic diagram for a trio sponsor}
\MS

Suppose $C$ is a simplicial set, $x \in C_k, y\in C_j, y' \in C_{j'}$, $m = k+j+j'+2$, $\SETT{X,A,A'}$ is a $(k,j,j')$-comb-trio of dimension $m$ and $w \in C_m$ is a sponsor of the trio $\SETT{(x,X),(y,A),(y',A')}$. 
That is,
$d_{A \cup A'}(w) = x,\; d_{X \cup A'}(w) = y$ and $d_{X \cup A}(w) = y'$. We will represent $w$ schematically by
\begin{center}
\begin{pictr}{2}{.6}{-.7}{.25}
\label{schematic-diag}
\put(0,0){\MB{(y',A')}}
\put(2,0){\MB{(y,A)}}
\put(1,.5){\MB{(x,X)}}
\put(.2,.15){\line(2,1){.5}}
\put(1.8,.15){\line(-2,1){.5}}
\put(.3,0){\line(1,0){1.4}}
\put(.45,.45){\MBS{z'}}
\put(1.55,.45){\MBS{z}}
\put(1,.1){\MBS{z''}}
\end{pictr}
\end{center}
\MS

\NI where $z,z'$ and $z''$ are subfaces of $w$ defined by
\[
z = d_{A'}(w), \quad z' = d_A(w), \quad z'' = d_X(w)
\]
\BX

\begin{newlem}
\label{faces of sponsor lemma}
{\rm
\SL

Suppose $w$ sponsors the trio $\ATRIO{x}{X}{y}{A}{y'}{A'}$, where the comb-trio $\SETT{X,A,A'}$ has complementary vertex function pairs $(\mu,\LAM)$, $(\mu',\LAM')$ and $(\GAM,\GAM')$ as in diagram \eqref{comb-trio-diag-1} above and the schematic diagram \eqref{schematic-diag}.

Then
\[
z=d_{A'}(w) \in \CPL(x,\mu;y,\LAM), \quad z'=d_A(w) \in \CPL(x,\mu';y',\LAM') \text{ and } z'' =d_X(w) \in \CPL(y,\gamma;y',\gamma')
\]
}
\end{newlem}

\Proof

It suffices to show $d_{A'}(w) \in \CPL(x,\mu;y,\LAM)$. The other two claims follow by symmetry.
\MS

Apply the Face Identity lemma \refpage{face-identity-lemma} to the diagram:

\begin{center}
\begin{picture}(2,.6)
\put(0,0){\MB{[k]}}
\put(1,0){\MB{[j+k+1]}}
\put(2,0){\MB{[j]}}
\put(1,.5){\MB{[m]}}
\put(.2,0){\vector(1,0){.4}}
\put(1.8,0){\vector(-1,0){.4}}
\put(1,.15){\vector(0,1){.2}}
\put(.2,.1){\vector(2,1){.6}}
\put(1.8,.1){\vector(-2,1){.6}}
\put(1.1,.25){\MBS{\tau}}
\put(.4,.3){\MBS{\HT{X}}}
\put(1.6,.3){\MBS{\HT{A}}}
\put(.4,-.1){\MBS{\mu\SH}}
\put(1.6,-.1){\MBS{\LAM\SH}}
\end{picture}
\end{center}
\BS
where $\tau[j+k+1] = X \cup A$.
\MS

By definition, $x = d_{A \cup A'}(w) = d_{[m]-\HT{X}}(w)$. By the Face Identity lemma
\[
d_{[m]-\HT{X}} = d_{[j+k+1]-\mu\SH[k]}\; d_{[m]-\tau[j+k+1]}
\]
Since $[j+k+1]-\mu\SH[k] = \LAM\SH[j]$ and $[m]-\tau[j+k+1] = A'$ it follows that $x = d_{\LAM\SH[j]}\; d_{A'}(w)$.
\MS

Similarly, $y =d_{X \cup A'}(w) = d_{[m]-\HT{A}}(w)$, by definition. Applying the Face Identity lemma the same way we get
\[
d_{[m]-\HT{A}} = d_{[j+k+1]-\LAM\SH[j]}\; d_{[m]-\tau[j+k+1]}
\]
and since $d_{[j+k+1]-\LAM\SH[j]} = d_{\mu\SH[k]}$ we get
\[
y = d_{\mu\SH[k]} d_{A'}(w)
\]
That is, $w \in \CPL(x,\mu;y,\LAM)$, as claimed.

\qed

\subsubsection{Lemmas for strictly increasing functions}

Before continuing with comb-trios below (page \pageref{comb-trio constructions}), we'll record several elementary but pertinent lemmas concerning strictly increasing functions in $\Delta$.
We will use the following notations:
\MS

$\bullet$ Given a strictly increasing $f:[j] \to [n]$, define $\Delta f :[j] \to [n]$ by
\[
\Delta f(p) \DFAS \left\{
\begin{array}{ll}
f(0) & \text{if } p=0\\
f(p)-f(p-1) & \text{if } 0<p \leq j
\end{array}
\right.
\]
\index{$\Delta f$}
Trivially, for each $p \in [j]$, $f(p) = \sum_{t=0}^p \Delta f(t)$. 
Of course $\Delta f$ is {\em not}, in general, non-decreasing.
\MS

$\bullet$ Given strictly increasing functions $f:[j] \to [n]$ and $g : [j] \to [m]$ then write $\Delta f \leq \Delta g$ to mean that for all $p \in [j]$, $\Delta f(p) \leq \Delta g(p)$.

\BX
\MS

\begin{newlem} 
\label{extension-lemma}
\index{Extension Lemma}
\index{lemma!Extension Lemma}
{\bf (Extension Lemma)}
{\rm
\SL

Suppose $f:[j] \to [n]$ and $g:[j] \to [m]$ are strictly increasing and $n \leq m$. Then the following are equivalent:
\begin{enumerate}
\item 
There exists a strictly increasing function $[n] \XRA{h} [m]$ such that $hf=g$.

\item
$\Delta f \leq \Delta g$ and $n-f(j) \leq m-g(j)$.
\end{enumerate}
}
\end{newlem}

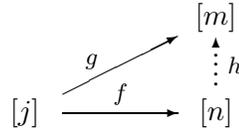
\begin{picture}(1,.4)(-3,0)
\put(0,0){\MB{[j]}}
\put(1,0){\MB{[n]}}
\put(1,.5){\MB{[m]}}
\put(.2,0){\vector(1,0){.6}}
\put(.2,.1){\vector(2,1){.6}}
\multiput(1,.15)(0,.05){4}{\MB{\cdot}}
\put(1,.4){\vector(0,1){0}}
\put(.5,.1){\FNS{\MB{f}}}
\put(.35,.27){\FNS{\MB{g}}}
\put(1.1,.25){\FNS{\MB{h}}}
\end{picture}

\Proof

$(1) \IMP (2)$ is immediate.

$(2) \IMP (1)$: Given $f$ and $g$ as assumed then for each $t \in [n]$ then either $t<f(0)$ or $t>f(j)$ or there exists a unique $p \in [0,j-1]$ such that $f(p) \leq t < f(p+1)$. Define a function $\DMAP{h}{n}{m}$ by
\[
h(t) \DFAS \left\{
\begin{array}{ll}
t & \text{if } t<f(0)\\
g(t)+t-f(p) & \text{if } f(p) \leq t< f(p+1)\\
g(j)+t-f(j) & \text{if } t>f(j)
\end{array}
\right.
\]
It follows from the assumptions on $f$ and $g$ that $h$ is well-defined, strictly increasing and that $hf=g$.

\qed
\MS

We note that in general, $h$ isn't unique.
\MS


\begin{newlem}
{\bf (Sum Lemma)}
\label{sum-lemma}
\index{Sum Lemma}
\index{lemma!Sum Lemma}
{\rm
\SL

Suppose $\DMAP{g}{n}{m}$ and $\DMAP{g'}{n'}{m}$ are strictly increasing such that $g[n] \cup g'[n'] = [m]$ and $|g[n] \cap g'[n']|=j+1$ for some $j \geq 0$. 

Then: 

(1) $m=n+n'-j$.   

(2) There exist unique \SI\ functions $\DMAP{f}{j}{n}$ and $\DMAP{f'}{j}{n'}$ such that $gf = g'f'$.    

(3) With $\DMAP{h}{j}{m}$ defined by $h = gf=g'f'$ then $h = (f+f')\FL$, equivalently $h\FL = f\FL + {f'}\FL$. 
}
\end{newlem}

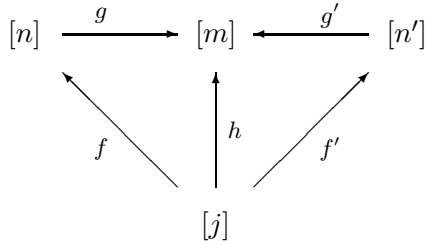
\begin{center}
\begin{picture}(2,1.2)
\MBP{0}{1}{[n]}
\MBP{1}{1}{[m]}
\MBP{2}{1}{[n']}
\MBP{1}{0}{[j]}
\put(.8,.2){\vector(-1,1){.6}}
\put(1.2,.2){\vector(1,1){.6}}
\put(1,.2){\vector(0,1){.6}}
\put(.2,1){\vector(1,0){.6}}
\put(1.8,1){\vector(-1,0){.6}}
\MBPS{1.1}{.5}{h}
\MBPS{.4}{.4}{f}
\MBPS{1.6}{.4}{f'}
\MBPS{.4}{1.1}{g}
\MBPS{1.6}{1.1}{g'}
\end{picture}
\end{center}

\Proof

\NI 1. For brevity, let $G = g[n]$, $G' = g'[n']$ and $H = G \cap G'$. Then $|G|=n+1$, $|G'|=n'+1$ and $|G \cap G'|=j+1$. Since $[m]$ is partitioned by $G-H$, $G'-H$ and $H$ we have
\[
\begin{aligned}
\Bigl|[m]\Bigr| &= |G-H| + |G'-H| + |H|\\
m+1 &= (n+1)-(j+1) + (n'+1)-(j+1) + (j+1)\\
\end{aligned}
\]
so that $m=n+n'-j$.
\MS

\NI 2. Define the \SI\ map $\DMAP{h}{j}{m}$ by $h[j] = H = g[n] \cap g'[n']$. Since $h[j] \SBS g[n]$ and $g$ is \SI\ then $g^{-1}(h[k]) \SBS [n]$ corresponds to exactly one \SI\ function $\DMAP{f}{j}{n}$ such that $gf=h$. The same reasoning shows that there is exactly one $\DMAP{f'}{j}{n'}$ such that $g'f'=h$.
\MS

\NI 3. Again we use $G \DFAS g[n]$, $G' \DFAS g'[n']$ and $H \DFAS h[j] = G \cap G'$.

For each $p \in [j]$ define
\[
I_p \DFAS [0,h(p)-1] = \SETT{0,1 \DDD{,} h(p)-1} \SBS [m]
\]
(We note that $I_p=\MT$ if $h(0)=0$.]

Then: $|I_p| = h(p)$ and $|I_p \cap H | = p$ because
\[
I_p \cap H = \SETT{ h(t):h(t)<h(p)}
\]
and $h$ is \SI.

Next, $I_p \cap G = \SETT{t \in [n]: g(t)<h(p)}$. Then
\[
t \in I_p \cap G \iff g(t)<h(p)=gf(p) \iff t<f(p)
\]
since $g$ is \SI. Therefore
\[
\Bigl| I_p \cap G \Bigr| = \Bigl| \SETT{t: 0 \leq t \leq f(p)-1}\Bigr| = f(p)
\]
By the same reasoning, $|I_p \cap G'| =f'(p)$.

Using that $G = (G-H) \cup H$ is a disjoint union then
\[
\begin{aligned}
I_p \cap G & = \bigl(I_p \cap(G-H)\bigr) \; \cup \; \bigl(I_p \cap H \bigr)\\
\Bigl| I_p \cap G \Bigr| & =\Bigl|I_p \cap(G-H)\Bigr| \; + \; \Bigl|I_p \cap H \Bigr|\\
f(p) &= \Bigl|I_p \cap(G-H)\Bigr| \; + p
\end{aligned}
\]
Therefore, $\Bigl|I_p \cap(G-H)\Bigr|=f(p)-p$. By the same reasoning, $\Bigl|I_p \cap(G'-H)\Bigr|=f'(p)-p$.

Finally, using that $G = (G-H) \cup (G'-H) \cup H$ is a disjoint union  then
\[
\begin{aligned}
h(p) = \Bigl| I_p \Bigr| &= \Bigl|I_p \cap(G-H) \Bigr| + \Bigl| I_p \cap (G'-H) \Bigr| + \Bigl| I_p \cap H \Bigr|\\
&= f(p)-p + f'(p)-p + p
\end{aligned}
\]
i.e. $h(p)=f(p)+f'(p)-p$ or, equivalently, 
\[
h\FL = f\FL + {f'}\FL \quad \text{and} \quad  h = (f+f')\FL
\]
\qed
\MS

\begin{example}
{\rm
Let $m=11$ and
\[
\begin{aligned}
G &= \SETT{1,2,4,5,7,9,10} \quad & \DMAP{g=(1,2,4,5,7,9,10)}{6}{11}\\
G' &= \SETT{0,2,3,6,7,8,9,10,11} \quad & \DMAP{g'=(0,2,3,6,7,8,9,10,11)}{8}{11}\\
H &= \SETT{2,7,9,10} \quad &   \DMAP{h=(2,7,9,10)}{3}{11}
\end{aligned}
\]
In the commutative diagram
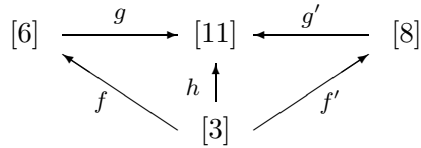
\begin{center}
\begin{picture}(2,.6)
\put(0,.5){\MB{[6]}} 
\put(1,.5){\MB{[11]}}
\put(2,.5){\MB{[8]}}
\put(1,0){\MB{[3]}}
\put(.2,.5){\vector(1,0){.6}}
\put(1.8,.5){\vector(-1,0){.6}}
\put(.8,0){\vector(-3,2){.6}}
\put(1.2,0){\vector(3,2){.6}}
\put(1,.15){\vector(0,1){.2}}
\put(.5,.6){\MBS{g}}
\put(1.5,.6){\MBS{g'}}
\put(.88,.23){\MBS{h}}
\put(.4,.15){\MBS{f}}
\put(1.6,.15){\MBS{f'}}
\end{picture}
\end{center}

\NI $f=\DMAP{(1,4,5,6)}{3}{6}$ and $f'=\DMAP{(1,4,6,7)}{3}{8}$.
\MS

The following chart is a summary:
\[
\begin{array}{c||c|c|c|c|c|c|c|c|c|c|c|c|}
 & 0 & 1 & 2 & 3 & 4 & 5 & 6 & 7 & 8 & 9 & 10 & 11\\
 \hline
\HT{G} &  & 0 & 1 &  &2 & 3 &  & 4 & & 5 & 6 & \\
\hline
\HT{G'} & 0 & &1 & 2 & & & 3 & 4 & 5 & 6 & 7 & 8\\
\hline
\HT{H} & & & 0 & & & & & 1 & & 2 & 3 & \\
\hline
\end{array}
\]
Then:
\[
\begin{array}{c||c|c|c|c|c|}
p & I_p & G \cap I_p & | G \cap I_p |=f(p) & G' \cap I_p & | G' \cap I_p |=f'(p) \\
\hline
0 & [0,1] & \{1\} & 1 & \{0\} & 1\\
\hline
1 & [0,6] & \{1,2,4,5\} & 4 & \{0,2,3,6\} & 4\\
\hline
2 & [0,8] & \{1,2,4,5,7\} & 5 & \{0,2,3,6,7,8\} & 6\\
\hline
3 & [0,9] & \{ 1,2,4,5,7,9\} & 6 & \{0,2,3,6,7,8,9 \} & 7\\
\hline
\end{array}
\]
and
\[
\begin{array}{c||c|c|c|c|}
p & h(p) & f(p) & f'(p) & f(p)+f'(p)-p\\
\hline
0 & 2 & 1 & 1 & 2 = 1+1-0\\
\hline
1 & 7 & 4 & 4 & 7=4+4-1\\
\hline
2 & 9 & 5 & 6 & 9=5+6-2\\
\hline
3 & 10 & 6 & 7 & 10=6+7-3 \\
\hline
\end{array}
\]
}
\end{example}
\BS

In the previous lemma, $f$ and $f'$ are determined by $g,\, g'$ and $h$ with $h(p)= f(p)+f'(p)-p$. There is a reverse construction, as follows. 

\begin{newlem} 
\label{sum-lemma-corollary}
{\rm
\SL

Suppose $j \geq 0$, $[j] \XRA{f} [n]$ and $[j] \XRA{f'} [n']$ are \SI, $m \DFAS n+n'-j$ and $\DMAP{h}{j}{m}$ is defined for each $p \in [j]$ by $h(p)=f(p)+f'(p)-p$. Then there exist \SI\ functions $[n] \XRA{g} [m]$ and $[n'] \XRA{g'} [m]$ such that $gf=h=g'f'$, $g[n] \cup g'[n'] = [m]$ and $g[n] \cap g'[n'] = h[j]$.

}
\end{newlem}

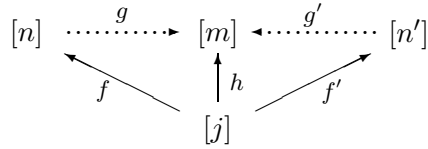
\begin{center}
\begin{picture}(2,.6)
\put(0,.5){\MB{[n]}}
\put(1,.5){\MB{[m]}}
\put(2,.5){\MB{[n']}}
\put(1,0){\MB{[j]}}
\put(1.2,.1){\vector(2,1){.6}}
\put(.4,.2){\FNS{\MB{f}}}
\put(.8,.1){\vector(-2,1){.6}}
\put(1.6,.2){\FNS{\MB{f'}}}
\multiput(.2,.5)(.05,0){12}{$.$}
\put(.8,.51){\vector(1,0){0}}
\multiput(1.8,.5)(-.05,0){13}{$.$}
\put(1.175,.51){\vector(-1,0){0}}
\put(1,.15){\vector(0,1){.25}}
\put(.5,.6){\FNS{\MB{g}}}
\put(1.5,.6){\FNS{\MB{g'}}}
\put(1.1,.25){\FNS{\MB{h}}}
\end{picture}
\end{center}

\Proof

First, we will apply the Extension Lemma (lemma \refpage{extension-lemma}) to show that at least one strictly increasing function $g : [n] \to [m]$ exists such that $gf=h$. 
To do so, we need to verify that $\Delta h \geq \Delta f$ and that $m-h(j) \geq n-f(j)$.
\MS

At $p=0$, $\Delta h(0) = f(0)+f'(0) \geq f(0) = \Delta f(0)$.

For $p>0$
\[
\begin{aligned}
h(p)-h(p-1) & = f(p)+f'(p)-p - f(p-1)-f'(p-1)+(p-1)\\
& = f(p)-f(p-1) \; + \; f'(p)-f'(p-1) -1\\
& \geq f(p)-f(p-1) = \Delta f(p)
\end{aligned}
\]
Finally, 
\[
\begin{aligned}
n+n'-j-h(j) & = n+n'-j-f(j)-f'(j)+j\\
& = n-f(j) \; + \; n'-f'(j)\\
& \geq n-f(j)
\end{aligned}
\]
Therefore, at least one strictly increasing $g:[n] \to [m]$ exists such that $gf=h$.

Working with any such $g$, let $G=g[n]$ and $H = h[j]$, as in the previous lemma. Define $G' = ([m]-G) \cup H$. Observe that
\[
\begin{aligned}
G \cup G' & = G \;\cup \; ([m]-G) \;\cup H = [m]\\
G \cap G' & = G \; \cap \; \Bigl([m]-G) \cup H\Bigr) = G \cap H = H
\end{aligned}
\]
Define $g':[n'] \to [m]$ by $g'[n'] = G'$.
Then, $g,g'$ and $h$ satisfy the hypotheses of the previous lemma. According to that lemma, there exists a unique strictly increasing function $f'' : [j] \to [n']$ such that $g' f'' = h$ and for all $p \in [j]$, $h(p) = f(p) + f''(p)-p$. It follows immediately that $f''=f'$.

\qed
\MS

\NI {\bf Note:} In the lemma, $[m] = (G-H) \cup (G'-H) \cup H$ is a disjoint union. That is, $\SETT{H,G-H,G'-H}$ forms a $(j,n-j,n'-j)$-comb-trio. We will apply this below.
\BX
\MS

\begin{newcor}
{\rm
\SL

Suppose $0 \leq j \leq \text{min}\{n,n'\}$ and $m = n+n'-j$. Let $\DMAP{h}{j}{m}$ and $\DMAP{f}{n}{j}$ be any \SI\ maps such that $\Delta h \geq \Delta f$ and $m-h(j) \geq n-f(j)$. Then: 

(1) There exists a \SI\ function $\DMAP{f'}{j}{n'}$ defined by $f'(t) \DFAS h(t)-f(t)+t$, $t \in [j]$. 

(2) There exist \SI\ functions $\DMAP{g}{n}{m}$ and $\DMAP{g'}{n'}{m}$ such that $gf=h=g'f'$, $g[n] \cup g'[n'] = [m]$ and $g[n] \cap g'[n'] = h[j]$.
}
\end{newcor}

\Proof

The solid arrows are given and the dotted arrows are claimed in:
\begin{center}
\begin{picture}(2,.6)
\put(0,.5){\MB{[n]}}
\put(1,.5){\MB{[m]}}
\put(2,.5){\MB{[n']}}
\put(1,0){\MB{[j]}}
\put(.8,.1){\vector(-2,1){.6}}
\multiput(1.2,.1)(.05,.025){12}{.}
\put(1.8,.4){\vector(2,1){0}}
\put(1,.1){\vector(0,1){.25}}
\put(.5,.15){\FNS{\MB{f}}}
\put(.93,.25){\FNS{\MB{h}}}
\put(1.5,.15){\FNS{\MB{f'}}}
\multiput(.2,.5)(.05,0){12}{.}
\put(.8,.505){\vector(1,0){0}}
\put(.5,.6){\FNS{\MB{g}}}
\multiput(1.8,.5)(-.05,0){12}{.}
\put(1.2,.505){\vector(-1,0){0}}
\put(1.5,.6){\FNS{\MB{g'}}}
\end{picture}
\end{center}
\MS

To see that $f'$ is \SI: for each $t \in [1,j]$ we have
\[
\begin{aligned}
f(t)-f(t-1) & = \Delta h(t) - \Delta f(t) + 1\\
& \geq 1
\end{aligned}
\]
$f'$ has the proper codomain because $f'(j) = h(j)-f(j)+j \leq m-n+j = n'$.
Since $f,f'$ and $h$ satisfy the hypotheses of Lemma \refpage{sum-lemma-corollary}, the claim (2) follows.
\qed
\MS

\begin{newlem} 
\label{double-extension-lemma}
{\bf (Joint Factorization Lemma)}
\index{Joint Factorization Lemma}
\index{lemma!Joint Factorization Lemma}
{\rm
\SL

Given $j,k \geq 0$, $n>j+k+1$ and strictly increasing maps $f,g,\alpha$ and $\beta$ in the diagram below such that $[j+k+1] = f[j] \cup g[k]$ is a disjoint union and $\alpha[j] \cap \beta[k] = \MT$, then a necessary and sufficient condition for the existence of a strictly increasing map $\gamma: [j+k+1] \to [n]$ such that $\gamma f = \alpha$ and $\gamma g = \beta$ is
\begin{equation} \label{lemma-condition}
\ALL p \in [j] \; \ALL q \in [k] \; \Bigl( f(p)<g(q)\; \iff \;\alpha(p)<\beta(q) \Bigr)
\end{equation}
\begin{center}
\begin{picture}(2,.5)
\put(0,0){\MB{[j]}}
\put(1,0){\MB{[j+k+1]}}
\put(2,0){\MB{[k]}}
\put(1,.5){\MB{[n]}}
\put(.2,0){\vector(1,0){.3}}
\put(.35,.1){\MBS{f}}
\put(1.8,0){\vector(-1,0){.3}}
\put(1.65,.1){\MBS{g}}
\put(.2,.15){\vector(2,1){.6}}
\put(.4,.35){\MBS{\alpha}}
\put(1.8,.15){\vector(-2,1){.6}}
\put(1.6,.35){\MBS{\beta}}
\multiput(1,.1)(0,.05){6}{.}
\put(1.02,.4){\vector(0,1){0}}
\put(.93,.23){\MBS{\gamma}}
\end{picture}
\end{center}
}
\end{newlem}

\Proof

The necessity of the condition \eqref{lemma-condition} is immediate.
\MS

Conversely: it follows from $[j+k+1]  = f[j] \cup g[k]$ being a disjoint union that there is a exactly one map $\gamma : [j+k+1] \to [n]$ such that $\gamma f = \alpha$ and $\gamma g = \beta$, namely
\[
\gamma (t) = 
\left\{
\begin{array}{ll}
\alpha(p) & \text{ if } t = f(p)\\
\beta(q) & \text{ if } t = g(q) 
\end{array}
\right.
\]
It remains to show that $\gamma$ is strictly increasing.

$\gamma$ is strictly increasing if and only if for each $t \in [0,j+k]$, $\gamma(t) < \gamma(t+1)$. Clearly, $\gamma(t)<\gamma(t+1)$ for those $t$ such that $t,t+1 \in f[j]$ and those $t$ such that $t,t+1 \in g[k]$.

Now if $t = f(p)$ and $t+1 = g(q)$ then $\gamma(t) = \alpha(p)$ and $\gamma(t+1) = \beta(q)$. Therefore, condition \eqref{lemma-condition} implies $\gamma(t) = \alpha(p)<\beta(q) = \gamma(t+1)$. 

Finally, if $t=g(q)<t+1 = f(p)$ then, by definition, $\gamma(t)= \beta(q)$ and $\gamma(t+1) = \alpha(p)$. By condition \eqref{lemma-condition}, and that $\alpha[j] \cap \beta[k] = \MT$ and $f[j] \cap g[k] = \MT $,
\[
f(p)>g(q) \iff \neg(f(p)<g(q)) \iff \neg(\alpha(p)<\beta(q)) \iff \beta(q)<\alpha(p)
\]
That is $\gamma(t)= \beta(q) < \alpha(p) = \gamma(t+1)$.

\qed

\begin{newthm}
\label{comb-trio-props}
{\rm
\SL

Suppose $\SETT{X,A,A'}$ is a $(k,j,j')$-comb-trio of dimension $m$. 
As before, suppose that the three pairs of complementary vertex functions associated with the comb-trio are, as in \DREF{comb-trio-diag-1}:
\[
\begin{aligned}
\DMAP{\mu}{k}{j+1}, \quad & \DMAP{\LAM}{j}{k+1} \quad \text{for } \HT{X \cup A}\\
\DMAP{\mu'}{k}{j'+1}, \quad & \DMAP{\LAM'}{j'}{k+1} \quad \text{for } \HT{X \cup A'} \\
\DMAP{\gamma}{j}{j'+1}, \quad & \DMAP{\gamma'}{j'}{j+1} \quad \text{for } \HT{A \cup A'}
\end{aligned}
\]
Then 
\[
\begin{aligned}
\ALL q \in [k], & \quad \HT{X}(q) = \mu\SH(q) + (\mu')\SH(q)-q\\
\ALL p \in [j], & \quad \HT{A}(p) = \LAM\SH(p) + \gamma\SH(p)-p\\
\ALL p' \in [j'], & \quad \HT{A'}(p') = (\LAM')\SH(p') + (\gamma')\SH(p')-p'\\
\end{aligned}
\]
And
\[
\begin{aligned}
\ALL p \in [j]\; \ALL q \in [k]\; &  \Bigl(\LAM\SH(p)<\mu\SH(q) \iff \HT{A}(p)<\HT{X}(q) \Bigr)\\
\ALL p' \in [j']\; \ALL q \in [k]\; &  \Bigl((\LAM')\SH(p')<(\mu')\SH(q) \iff \HT{A'}(p')<\HT{X}(q) \Bigr)\\
\ALL p \in [j]\; \ALL p' \in [j']\; & \Bigl(\gamma\SH(p)<(\gamma')\SH(p') \iff \HT{A}(p)<\HT{A'}(p') \Bigr)
\end{aligned}
\]
}
\end{newthm}

\Proof

Look at the portion of diagram \refpage{comb-trio-diag-1} involving $\HT{X}$, $\mu\SH$ and $\mu'{\SH}$:
\begin{center}
\begin{picture}(2,1.1)
\put(0,0){\MB{[j'+k+1]}}
\put(1,0){\MB{[k]}} 
\put(2,0){\MB{[j+k+1]}}
\put(1,1){\MB{[m]}}
\put(.85,0){\vector(-1,0){.4}}
\put(1.15,0){\vector(1,0){.4}}
\put(1,.15){\vector(0,1){.7}}
\put(0,.15){\vector(1,1){.8}}
\put(2,.15){\vector(-1,1){.8}}
\put(.65,-.1){\MBS{\mu'{\SH}}}
\put(1.35,-.1){\MBS{\mu{\SH}}}
\put(1.1,.5){\MBS{\HT{X}}}
\put(.2,.65){\MBS{\HT{X \cup A'}}}
\put(1.8,.65){\MBS{\HT{X \cup A}}}
\end{picture} 
\end{center}
\MS

\NI Since $X \cup A \SBS [m]$ and $X \cup A' \SBS [m]$ satisfy the hypotheses of the Sum Lemma above (page  \pageref{sum-lemma}), then for all $q \in[k]$, $\HT{X}(q)=\mu\SH(q)+\mu'{\SH}(q)-q$, as claimed.

The claims for $\HT{A}$ and $\HT{A'}$ follow the same way.
\MS

The remaining claims follow from the Joint Factorization Lemma \refpage{double-extension-lemma} above.

\qed

\subsubsection{Constructing comb-trios with given complementary function pairs}
\label{comb-trio constructions}

A comb-trio $\SETT{X,A,A'}$ of dimension $m$ determines three pairs of complementary vertex functions as in \DREF{comb-trio-diag-1}. In this section we consider the reverse: constructing a comb-trio of dimension $m$ given one or more pairs of complementary vertex functions.

\begin{newlem}
{\rm
\SL

Suppose $j,j',k \geq 0$, $m = j+j'+k+2$ and suppose that $\DMAP{\mu}{k}{j+1}$ and $\DMAP{\LAM}{j}{k+1}$ are complementary vertex functions.
\begin{enumerate}

\item
Let $X,A \SBS [m]$ such that $X \cap A=\MT$ and $A' \DFAS [m]-(X \cup A) \neq \MT$. Denote by $\DMAP{\HT{X}}{k}{m}$ and $\DMAP{\HT{A}}{j}{m}$ the \SI\ functions defined by $\HT{X}[k]=X$ and $\HT{A}[j] = A$ as in diagram \eqref{mu-lambda-X-A-lemma-diag} below.

Then the comb-trio $\SETT{X,A,A'}$ contains $(\mu,\LAM)$ iff
\[
\ALL q \in [k]\; \ALL p \in [j]\; \left(
\mu\SH(q)<\LAM\SH(p) \iff \HT{X}(q)<\HT{A}(p)
\right)
\]

\item
Suppose $\DMAP{\xi}{k}{m}$ is a strictly increasing function such that $\Delta\xi \geq \Delta \mu\SH$ and $m-\xi(k) \geq j+k+1-\mu\SH(k)$, as in diagram \eqref{mu-lambda-xi-lemma-diag} below. 

Then there exists a $(k,j,j')$-comb-trio $\SETT{X,A,A'}$ containing $(\mu,\LAM)$ such that $X =\xi[k]$.

\end{enumerate}
}
\end{newlem}

\Proof

\NI 1. By the Joint Factorization Lemma \refpage{double-extension-lemma} applied to the diagram
\begin{center}
\begin{pictr}{2}{.6}{-.7}{.25}
\label{mu-lambda-X-A-lemma-diag}
\put(0,0){\MB{[k]}}
\put(1,0){\MB{[j+k+1]}}
\put(2,0){\MB{[j]}}
\put(1,.5){\MB{[m]}}
\put(.2,0){\vector(1,0){.4}} 
\put(.4,-.1){\MBS{\mu\SH}}
\put(1.8,0){\vector(-1,0){.4}} 
\put(1.6,-.1){\MBS{\LAM\SH}}
\put(.2,.1){\vector(2,1){.6}} 
\put(.35,.28){\MBS{\HT{X}}}
\put(1.8,.1){\vector(-2,1){.6}} 
\put(1.7,.28){\MBS{\HT{A}}}
\multiput(1,.1)(0,.05){6}{.}
\put(1.02,.4){\vector(0,1){0}}
\put(1.1,.2){\MBS{\tau}}
\end{pictr}
\end{center}
\BS

\NI there exists a unique strictly increasing $\tau$ making the diagram commute. In fact, $\tau = \HT{X \cup A}$. As noted earlier, $\tau, \mu\SH$ and $\LAM\SH$ determine the comb-trio (as in diagram \refpage{comb-trio-diag-1}) $\SETT{X,A,A'}$
which contains the complementary vertex function pair $(\mu,\LAM)$.

\BS

\NI 2. The Extension Lemma \refpage{extension-lemma} implies there exists at least one $\tau:[j+k+1] \to [m]$ such that $\tau \mu\SH = \xi$. 
\begin{center}
\begin{pictr}{2}{.55}{-.7}{.25}
\label{mu-lambda-xi-lemma-diag}
\put(0,0){\MB{[k]}}
\put(1,0){\MB{[j+k+1]}}
\put(2,0){\MB{[j]}}
\put(1,.5){\MB{[m]}}
\put(.2,0){\vector(1,0){.4}}
\put(1.8,0){\vector(-1,0){.4}}
\put(1,.1){\vector(0,1){.3}}
\put(.2,.1){\vector(2,1){.6}}
\put(1.8,.1){\vector(-2,1){.6}}
\put(.2,.3){\MBS{\text{given }\xi}}
\put(1.12,.25){\MBS{\exists\, \tau}}
\put(1.6,.35){\MBS{\alpha}}
\put(.4,-.1){\MBS{\mu\SH}}
\put(1.6,-.1){\MBS{\LAM\SH}}
\end{pictr}
\end{center}
\BS

\NI $\tau$ determines a $(k,j,j')$-comb-trio $\SETT{X, A,A'}$ where $X \DFAS \xi[k]$, $A \DFAS \tau \LAM\SH[j]$ and $A' \DFAS [m]-(X \cup A)$.
Any such comb-trio includes a complementary vertex function pair 
\[
\Bigl([k] \XRA{\mu'} [j'+1], \; [j'] \XRA{\LAM'} [k+1] \Bigr)
\]
where, according to Theorem \refpage{comb-trio-props} above,
\[
\ALL q \in [k] \quad \xi(q) = \mu\SH(q) + {\mu'}\SH (q)-q
\]
which serves to define $\mu'$ and which does not mention $\tau$ explicitly. That ${\mu'}\SH$ is \SI\ follows from $\Delta \xi \geq \Delta \mu\SH$.

\qed

\begin{newthm}
\label{trio-from-two-pairs}
{\rm
\SL

Suppose $j,j',k \geq 0$ and $m = j+j'+k+2$. Let 
\[
\begin{aligned}
(\mu,\LAM), &  \quad \mu:[k] \to [j+1],\; \LAM:[j] \to [k+1]\\
(\mu',\LAM'), &  \quad \mu':[k] \to [j'+1],\; \LAM:[j'] \to [k+1]
\end{aligned}
\] 
be two pairs of complementary vertex functions, and define $\xi:[k] \to [m]$ by $\xi(q) \DFAS \mu\SH(q) + {\mu'}\SH(q)-q$, each $q \in [k]$.

Then there is a $(k,j,j')$-comb-trio to which the pairs $(\mu,\LAM)$ and $(\mu',\LAM')$ belong.
}
\end{newthm}

\Proof

Diagram (\ref{mu-lam-mu'-lam'-thm-diag}) below shows a partial diagram for the claimed comb-trio, where the solid arrows are given and the dotted arrows are to be deduced.

Apply Lemma \refpage{sum-lemma-corollary} with $\mu\SH, {\mu'}\SH$ and $\xi$ playing the roles of $f,f'$ and $h$ to obtain strictly increasing functions $\tau$ and $\tau'$ in the diagram below.

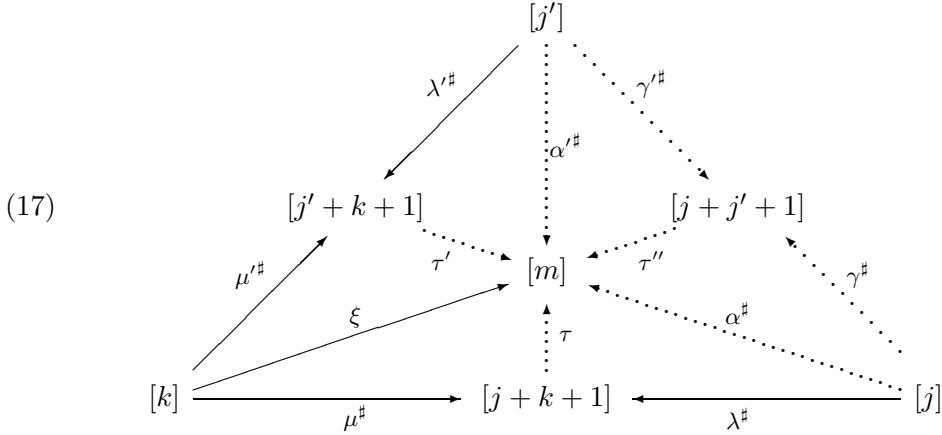
\begin{center}
\begin{pictr}{4}{2}{-.7}{1}
\label{mu-lam-mu'-lam'-thm-diag}
\MBP{0}{0}{[k]}
\MBP{2}{0}{[j+k+1]}
\MBP{4}{0}{[j]}
\MBP{1}{1}{[j'+k+1]}
\MBP{2}{2}{[j']}
\MBP{3}{1}{[j+j'+1]}
\MBP{2}{.667}{[m]}
\put(.15,.15){\vector(1,1){.7}} 
\MBPS{.45}{.65}{{\mu'}\SH}
\put(.15,.05){\vector(3,1){1.65}} 
\MBPS{1}{.43}{\xi}
\put(.15,0){\vector(1,0){1.4}}
\MBPS{1}{-.1}{\mu\SH}
\multiput(3.85,.25)(-.05,.05){12}{\MB{.}}
\put(3.25,.85){\vector(-1,1){0}}
\MBPS{3.63}{.65}{\gamma\SH} 
\multiput(3.85,.05)(-.051,.017){32}{\MB{.}}
\put(2.218,.594){\vector(-3,1){0}}
\MBPS{3}{.45}{{\alpha}\SH}
\put(3.85,0){\vector(-1,0){1.4}}
\MBPS{3}{-.1}{\LAM\SH}
\put(1.85,1.85){\vector(-1,-1){.7}}
\MBPS{1.45}{1.65}{{\LAM'}\SH}
\multiput(2.15,1.85)(.05,-.05){14}{\MB{.}}
\put(2.85,1.15){\vector(1,-1){0}}
\MBPS{2.55}{1.65}{{\gamma'}\SH}
\multiput(2,1.85)(0,-.05){21}{\MB{.}}
\put(2,.8){\vector(0,-1){0}}
\MBPS{2.1}{1.33}{{\alpha'}\SH}
\multiput(2,.15)(0,.05){7}{\MB{.}}
\put(2,.5){\vector(0,1){0}}
\MBPS{2.1}{.33}{\tau}
\multiput(1.36,.88)(.051,-.017){9}{\MB{.}}
\put(1.819,.727){\vector(3,-1){0}}
\MBPS{1.45}{.725}{\tau'}
\multiput(2.67,.89)(-.051,-.017){9}{\MB{.}}
\put(2.211,.737){\vector(-3,-1){0}}
\MBPS{2.55}{.725}{\tau''}
\end{pictr}
\end{center}
\BS

Let $T = \tau[j+k+1], T'=\tau'[j'+k+1]$ and $X = \xi[k]$. Then, according to the lemma, $T \cup T' = [m]$ and $T \cap T' = X$. Also
\[
[m] = (T-X) \cup (T'-X) \cup X \text{ is a disjoint union}
\]
Let $A = T-X$ and $A' = T'-X$.

Since $[j+k+1] = \mu\SH[k] \cup \LAM\SH[j]$  and $T=X \cup A$ are disjoint unions then
\[
\begin{aligned}
T & = \tau \mu\SH[k] \cup \tau \LAM\SH[j] \\
& = X \cup  \tau \LAM\SH[j] \\
& = X \cup A
\end{aligned}
\]
implies $\tau \LAM\SH[j] = A$. We define $\alpha\SH : [j] \to [m]$ by $\alpha\SH[j] = A$ so that $\alpha\SH = \tau \LAM\SH$.

By the same reasoning, $\tau'{\LAM'}\SH[j']= A'$ and if we define ${\alpha'}\SH : [j'] \to [m]$ by ${\alpha'}\SH[j'] = A'$ then $\tau' {\LAM'}\SH = {\alpha'}\SH$.

Define $\tau'' : [j+j'+1] \to [m]$ by $\tau''[j+j'+1] = A \cup A' = [m]-X$.

Finally, using that $A \cup A'$ is a disjoint union we may define $\gamma\SH:[j] \to [j+j'+1]$ and ${\gamma'}\SH:[j'] \to [j+j'+1]$ by
\[
\begin{aligned}
\gamma\SH(p)\DFAS t & \quad \text{ when } \tau''(t) = \alpha(p) \\
{\gamma'}\SH(p')\DFAS t' & \quad \text{ when } \tau''(t') = \alpha'(p') \\
\end{aligned}
\]
Then $\gamma$ and $\gamma'$ form a complementary vertex function pair because \\ $\gamma\SH[j] \cup {\gamma'}\SH[j'] = [j+j'+1]$ and $\gamma\SH[j] \cap {\gamma'}\SH[j'] = \MT$. 

This yields the complete specification of the claimed comb-trio.

\qed

The comb-trio is this theorem isn't unique because $\tau$ is not, in general, unique. Here's an example to illustrate that.
\MS

\begin{example}
{\rm
Take $k=j=j'=1$, $\mu\SH = (0,1):[1] \to [3]$ and $\xi = (0,2): [1] \to [5]$. Then $\LAM\SH = (2,3):[1] \to [3]$. 

One choice for $\tau : [3] \to [5]$ is $\tau = (0,2,3,4)$, and $\alpha = \tau \LAM\SH = (3,4)$. This determines the $(1,1,1)$-comb-trio $\SETT{X,A,A'}$ where $X=\{0,2\}, A = \{3,4\}, A'=\{1,5\}$.

Another choice for $\tau$ is $\tau = (0,2,3,5)$ and then $\tau \LAM\SH = (3,5)$. The corresponding comb-trio is $X=\{0,2\}, A = \{3,5\}, A'=\{1,4\}$.
}
\end{example}
\MS

Given three complementary vertex function pairs, there may be no comb-trio which contains all three of them, as the following example illustrates.
\MS

\begin{example}
{\rm
Let $k=j=j'=1$, $m=5$ and consider the following two vertex function pairs:
\[
\begin{aligned}
\DMAP{\mu=(0,0)}{1}{2} & \qquad \DMAP{\mu\SH = (0,1)}{1}{3} \\
\DMAP{\LAM = (2,2)}{1}{2} & \qquad \DMAP{\LAM\SH = (2,3)}{1}{3}\\
& \\
\DMAP{\mu' = (2,2)}{1}{2} & \qquad \DMAP{{\mu'}\SH = (2,3)}{1}{3}\\
\DMAP{\LAM' = (0,0)}{1}{2} & \qquad \DMAP{{\LAM'}\SH=(0,1)}{1}{3}\\
\end{aligned}
\]
We will see that there is exactly one comb-trio containing $(\mu,\LAM)$ and $(\mu',\LAM')$.
\MS

First, $\DMAP{\xi = (2,3)}{1}{5}$ is defined by the Sum Lemma as $\xi(q) = \mu\SH(q)+{\mu'}\SH(q)-q$. Then it is immediate that the only possible $\DMAP{\tau}{3}{5}$ such that $\tau \mu\SH = \xi$ is $\tau = (2,3,4,5)$. Then, following the construction in the proof of Theorem \refpage{sum-lemma-corollary}, $\DMAP{\tau'}{3}{5}$ is the strictly increasing function defined by
\[
\tau'[3] = \Bigl( [5]-\SETT{2,3,4,5} \Bigr) \cup \SETT{2,3} = \SETT{0,1,2,3}
\]
and this is the only possible function $\tau'$ such that $\tau' {\mu'}\SH = \xi$.
Then 
\[
\begin{aligned}
\DMAP{\alpha\SH= (2,3,4,5) \circ (2,3) = (4,5)}{1}{5}\\ 
\DMAP{{\alpha'}\SH = (0,1,2,3) \circ (0,1) = (0,1)}{1}{5}\\
\DMAP{\tau'' = (0,1,4,5)}{3}{5}
\end{aligned}
\]
and the only possible comb-trio is
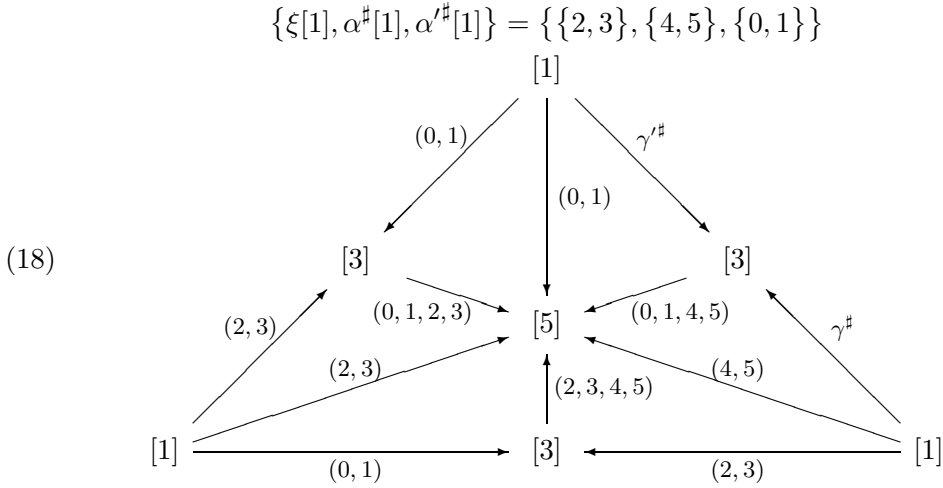
\[
\SETT{\xi[1], \alpha\SH[1], {\alpha'}\SH[1]}=
\SETT{\SETT{2,3}, \SETT{4,5}, \SETT{0,1}}
\]

\begin{center}
\begin{pictr}{4}{2}{-.7}{1}
\MBP{0}{0}{[1]}
\MBP{2}{0}{[3]}
\MBP{4}{0}{[1]}
\MBP{1}{1}{[3]}
\MBP{2}{2}{[1]}
\MBP{3}{1}{[3]}
\MBP{2}{.667}{[5]}
\put(.15,.15){\vector(1,1){.7}} 
\MBPS{.45}{.65}{(2,3)}
\put(.15,.05){\vector(3,1){1.65}} 
\MBPS{1}{.43}{(2,3)}
\put(.15,0){\vector(1,0){1.65}}
\MBPS{1}{-.1}{(0,1)}
\put(3.85,.15){\vector(-1,1){.7}}
\MBPS{3.55}{.65}{\gamma\SH} 
\put(3.85,.05){\vector(-3,1){1.65}} 
\MBPS{3}{.43}{(4,5)}
\put(3.85,0){\vector(-1,0){1.65}}
\MBPS{3}{-.1}{(2,3)}
\put(1.85,1.85){\vector(-1,-1){.7}}
\MBPS{1.45}{1.65}{(0,1)}
\put(2.15,1.85){\vector(1,-1){.7}}
\MBPS{2.55}{1.65}{{\gamma'}\SH}
\put(2,1.85){\vector(0,-1){1.03}}
\MBPS{2.2}{1.33}{(0,1)}
\put(2,.15){\vector(0,1){.36}}
\MBPS{2.3}{.33}{(2,3,4,5)}
\put(1.27,.91){\vector(3,-1){.53}}
\MBPS{1.35}{.725}{(0,1,2,3)}
\put(2.73,.91){\vector(-3,-1){.53}}
\MBPS{2.7}{.725}{(0,1,4,5)}

\end{pictr}
\end{center}
\BS

\NI From these calculations, $\DMAP{\gamma\SH = (2,3)}{1}{3}$ and $\DMAP{{\gamma'}\SH = (0,1)}{1}{3}$ are the only possible functions such that $\tau'' \gamma\SH = \alpha\SH$ and $\tau'' {\gamma'}\SH = {\alpha'}\SH$. 

Thus, given $(\mu,\LAM)$, $(\mu',\LAM')$ and any {\em other} choice of $(\gamma,\gamma')$, there is no comb-trio containing all three pairs. 
\MB{ }
}
\end{example}
\MS

\begin{newdef}
\label{compatible vertex function}
\index{compatible vertex function}
{\bf (Compatible vertex function)}
{\rm
\SL

Suppose $j,j',k \geq 0$. Given non-decreasing maps $\DMAP{\mu}{k}{j+1}$ and $\DMAP{\mu'}{k}{j'+1}$ we say {\bf the non-decreasing map $\DMAP{\gamma}{j}{j'+1}$ is compatible with $\mu$ and $\mu'$} if there is a comb-trio whose complementary vertex function pairs are $(\mu,\LAM)$, $(\mu',\LAM')$ and $(\gamma,\gamma')$, where $\LAM$, $\LAM'$ and $\gamma'$ are the complementary vertex functions of $\mu$, $\mu'$ and $\gamma$ respectively.
}

\BX
\end{newdef}
\MS

Given $\mu$ and $\mu'$ as in this definition, let $m \DFAS j+j'+k+2$ and $\DMAP{\xi}{k}{m}$ where $\xi(t) \DFAS \mu\SH(t) + {\mu'}\SH(t)-t$, and $X \DFAS \xi[k] \SBS [m]$. Then $\DMAP{\gamma}{j}{j'+1}$ is compatible with $\mu$ and $\mu'$ means that there is a partition $\SETT{X,A,A'}$ of $[m]$ whose complementary vertex functions include $\mu$, $\mu'$ and $\gamma$.

\subsubsection{Comb-trios and partial simplices}

\begin{newthm} 
\label{partial-simplex-from-two-pairs}
{\rm
\SL

Suppose $\SETT{X,A_1,A_2}$ is a $(k,j_1,j_2)$-comb-trio where $\tau_1= \HT{X \cup A_1}$, $\tau_2 = \HT{X \cup A_2}$ and $\tau_3=\HT{A_1 \cup A_2}$ as in the diagram:
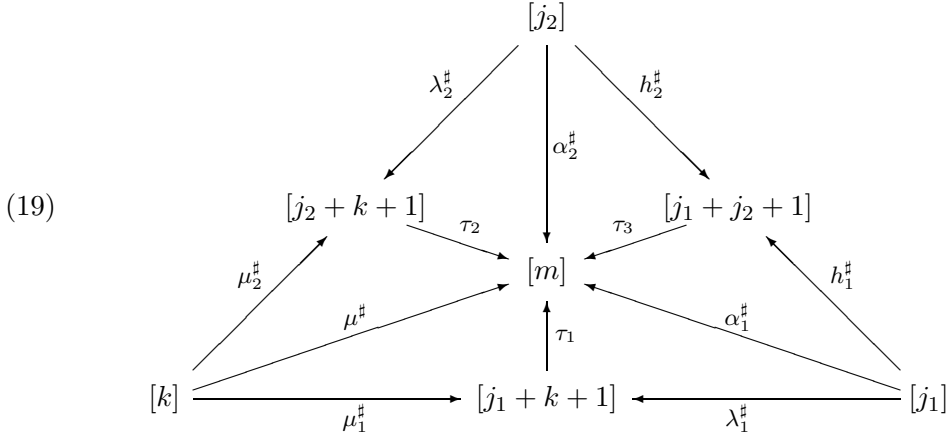
\begin{center}
\begin{pictr}{4}{2}{-.7}{1}
\MBP{0}{0}{[k]}
\MBP{2}{0}{[j_1+k+1]}
\MBP{4}{0}{[j_1]}
\MBP{1}{1}{[j_2+k+1]}
\MBP{2}{2}{[j_2]}
\MBP{3}{1}{[j_1+j_2+1]}
\MBP{2}{.667}{[m]}
\put(.15,.15){\vector(1,1){.7}} 
\MBPS{.45}{.65}{\mu_2\SH}
\put(.15,.05){\vector(3,1){1.65}} 
\MBPS{1}{.43}{\mu\SH}
\put(.15,0){\vector(1,0){1.4}}
\MBPS{1}{-.1}{\mu_1\SH}
\put(3.85,.15){\vector(-1,1){.7}}
\MBPS{3.55}{.65}{h_1\SH} 
\put(3.85,.05){\vector(-3,1){1.65}} 
\MBPS{3}{.43}{\alpha_1\SH}
\put(3.85,0){\vector(-1,0){1.4}}
\MBPS{3}{-.1}{\LAM_1\SH}
\put(1.85,1.85){\vector(-1,-1){.7}}
\MBPS{1.45}{1.65}{\LAM_2\SH}
\put(2.15,1.85){\vector(1,-1){.7}}
\MBPS{2.55}{1.65}{h_2\SH}
\put(2,1.85){\vector(0,-1){1.03}}
\MBPS{2.1}{1.33}{\alpha_2\SH}
\put(2,.15){\vector(0,1){.36}}
\MBPS{2.1}{.33}{\tau_1}
\put(1.27,.91){\vector(3,-1){.53}}
\MBPS{1.6}{.9}{\tau_2}
\put(2.73,.91){\vector(-3,-1){.53}}
\MBPS{2.4}{.9}{\tau_3}
\end{pictr}
\end{center}
\BS

Let $C$ be a simplicial set, $x \in C_k$, $y_1 \in C_{j_1}$ and $y_2 \in C_{j_2}$.

Then: If $z_1 \in \CPL(x,\mu_1;y_1,\LAM_1)$ and $z_2 \in \CPL(x,\mu_2;y_2,\LAM_2)$ then
\[
\SETT{(z_1,X \cup A_1),(z_2,X \cup A_2)}
\]
is a partial $m$-simplex. 
}
\end{newthm}

\Proof

Since $X \cup A_1 = [m]-A_2$ and $X \cup A_2 = [m]-A_1$ then the claim is that 
\[
\SETT{(z_1,[m]-A_2),(z_2,[m]-A_1)}
\]
is a partial $m$-simplex.
The relevant diagram corresponding to \DREF{par-simp-diag-frag} (page \pageref{par-simp-diag-frag}) is 
\begin{center}
\setlength{\unitlength}{1in}
\begin{picture}(2,.6)
\put(0,0){\MB{[j_2+k+1]}}
\put(1,0){\MB{[k]}}
\put(2,0){\MB{[j_1+k+1]}}
\put(1,.5){\MB{[m]}}
\put(.8,0){\vector(-1,0){.4}}
\put(1.2,0){\vector(1,0){.4}}
\put(1,.15){\vector(0,1){.25}}
\put(.1,.12){\vector(2,1){.65}}
\put(1.9,.12){\vector(-2,1){.65}}
\put(.4,.4){\FNS{\MB{\tau_2}}}
\put(1.65,.4){\FNS{\MB{\tau_1}}}
\put(.65,.12){\FNS{\MB{\mu_2\SH}}}
\put(1.3,.12){\FNS{\MB{\mu_1\SH}}}
\put(1.08,.25){\FNS{\MB{\HT{X}}}}
\end{picture}
\end{center}
\MS

\NI and we note that 
\[
\begin{aligned}
\tau_1[j_1+k+1] &= X \cup A_1= [m]-A_2\\
\tau_2[j_2+k+1] &= X \cup A_2 = [m]-A_1
\end{aligned}
\]
The Face Identity Lemma states
\[
d_{[j_1+k+1]-\mu_1\SH[k]} d_{[m]-\tau_1[j_1+k+1]} = 
d_{[j_2+k+1]-\mu_2\SH[k]} d_{[m]-\tau_2[j_2+k+1]}
\]
which in this case is
\[
d_{\LAM_1\SH[k]} d_{A_2} = d_{\LAM_2\SH[k]} d_{A_1}
\]
Therefore
the required face equation is 
\[
d_{[j_1+k+1]-\mu_1\SH[k]}(z_1)=d_{\LAM_1\SH[k]}(z_1) \ISIT d_{[j_2+k+1]-\mu_2\SH[k]}(z_2)=d_{\LAM_2\SH[k]}(z_2)
\] 
This holds because 
\[
d_{\LAM_1\SH[k]}(z_1) = x = d_{\LAM_2\SH[k]}(z_2)
\]
\qed
\BS

\begin{newdef} 
{\bf (Partial sponsor of a trio)}
\label{partial-sponsor} \index{partial sponsor}
{\rm
\SL

Let $C$ be a simplicial set, $j_1,j_2,k \geq 0$, $m=j_1+j_2+k+2$, and
\[
x \in C_k, \quad y_1 \in C_{j_1}, \quad y_2 \in C_{j_2}
\]
Suppose $z_1 \in \CPL(x,\mu_1;y_1,\LAM_1)$ and $z_2 \in \CPL(x,\mu_2;y_2,\LAM_2)$, and suppose $\SETT{X,A_1,A_2} $ is a $(k,j_1,j_2)$-comb trio having $(\mu_1,\LAM_1)$ and $(\mu_2,\LAM_2)$ as two of its complementary vertex function pairs. Then we'll call the partial $m$-simplex $\SETT{(z_1,X \cup A_1),(z_2,X \cup A_2)}$ a {\bf partial sponsor of the trio $\SETT{(x,X),(y_1,A_1),(y_2,A_2)}$ in $C$}. 
}

\BX
\end{newdef}
\MS

In effect, the partial sponsor is two-thirds of diagram \refpage{schematic-diag}:

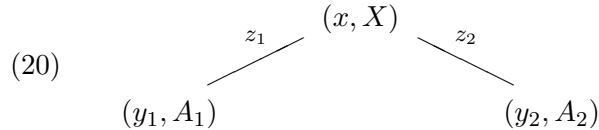
\begin{center}
\begin{pictr}{2}{.6}{-.7}{.25}
\put(0,0){\MB{(y_1,A_1)}}
\put(2,0){\MB{(y_2,A_2)}}
\put(1,.5){\MB{(x,X)}}
\put(.2,.15){\line(2,1){.5}}
\put(1.8,.15){\line(-2,1){.5}}
\put(.45,.4){\MBS{z_1}}
\put(1.55,.4){\MBS{z_2}}
\end{pictr}
\end{center}

\subsection{Trios, complementary vertex functions in $C^x$ and sponsors}
\label{trios in C^x}

In this section, $C$ is a simplicial set. Suppose:
\[
k \geq 0 ,\quad x \in C_k, \qquad  j_1 \geq 0, \quad y_1 \in C_{j_1}
\]
and $z_1 \in \CPL(x,\mu_1;y_1,\LAM_1)$ where
\[
\DMAP{\mu_1}{k}{j_1+1} \quad \text{and} \quad \DMAP{\LAM_1}{j_1}{k+1}
\]
are given complementary vertex functions.

Given any $m \geq j_1+k+2$ let $j_2 = m-(k+j_1+2)$ and let
\[
\DMAP{\mu_2}{k}{j_2+1}\quad \text{and} \quad \DMAP{\LAM_2}{j_2}{k+1}
\]
be any other pair of complementary vertex functions.
\MS

By Theorem \refpage{trio-from-two-pairs} there exists a comb-trio $\SETT{X,A_1,A_2}$ of dimension $m$ to which the pairs $(\mu_1,\LAM_1)$ and $(\mu_2,\LAM_2)$ belong. In effect, $z_1$ is one edge of a possible trio with schematic diagram (page \pageref{schematic-diag}):
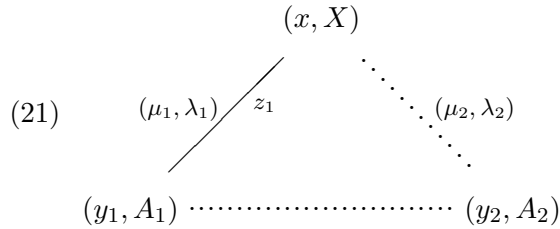
\begin{center}
\begin{pictr}{2}{1.1}{-.5}{.5}
\put(0,0){\MB{(y_1,A_1)}}
\put(1,1){\MB{(x,X)}}
\put(2,0){\MB{(y_2,A_2)}}
\put(.2,.2){\line(1,1){.6}}
\put(.25,.53){\MBS{(\mu_1,\LAM_1)}}
\put(.7,.55){\MBS{z_1}}
\multiput(.3,0)(.05,0){28}{.}
\multiput(1.2,.78)(.05,-.05){12}{.}
\put(1.8,.53){\MBS{(\mu_2,\LAM_2)}}
\end{pictr}
\end{center}
\MS

Given any $z_2 \in \CPL(x,\mu_2;y_2,\LAM_2)$ we ask whether there exists $w \in C_m$ which sponsors the partial trio. That is: does there exist $w \in C_m$ such that $d_{A_2}(w)=z_1$ and $d_{A_1}(w)=z_2$?
\MS

\NI Given this set-up, the assertion 
\begin{quote}
``For all $z_2  \in \CPL(x,\mu_2;y_2,\LAM_2)$ such a sponsor $w$ exists (or exists uniquely)''.
\end{quote}
amounts to a property of $z_1$.
\MS

To prepare for examining such a property, this section develops some facts concerning sponsors $w$ culminating in the Trio Sponsor Theorem \refpage{trio-sponsor-theorem}.
\BS

Let $\SETT{X,A_1,A_2}$ be a $(k,j_1,j_2)$-comb-trio whose complementary vertex function pairs are $(\mu_1,\LAM_1)$, $(\mu_2,\LAM_2)$ and $(h_{12},h_{21})$ as in the diagram of \SI\ functions: 

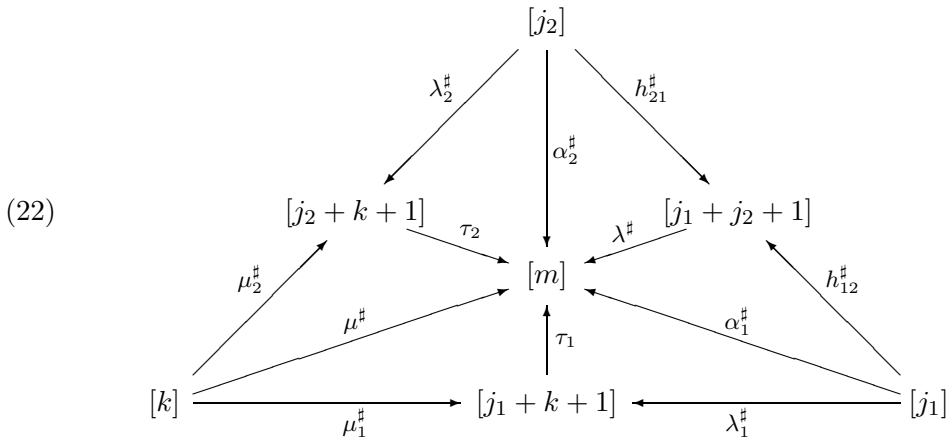
\begin{center}
\begin{pictr}{4}{2}{-.7}{1}
\label{comb-trio-diag-2}
\MBP{0}{0}{[k]}
\MBP{2}{0}{[j_1+k+1]}
\MBP{4}{0}{[j_1]}
\MBP{1}{1}{[j_2+k+1]}
\MBP{2}{2}{[j_2]}
\MBP{3}{1}{[j_1+j_2+1]}
\MBP{2}{.667}{[m]}
\put(.15,.15){\vector(1,1){.7}} 
\MBPS{.45}{.65}{\mu_2\SH}
\put(.15,.05){\vector(3,1){1.65}} 
\MBPS{1}{.43}{\mu\SH}
\put(.15,0){\vector(1,0){1.4}}
\MBPS{1}{-.1}{\mu_1\SH}
\put(3.85,.15){\vector(-1,1){.7}}
\MBPS{3.55}{.65}{h_{12}\SH} 
\put(3.85,.05){\vector(-3,1){1.65}} 
\MBPS{3}{.43}{\alpha_1\SH}
\put(3.85,0){\vector(-1,0){1.4}}
\MBPS{3}{-.1}{\LAM_1\SH}
\put(1.85,1.85){\vector(-1,-1){.7}}
\MBPS{1.45}{1.65}{\LAM_2\SH}
\put(2.15,1.85){\vector(1,-1){.7}}
\MBPS{2.55}{1.65}{h_{21}\SH}
\put(2,1.85){\vector(0,-1){1.03}}
\MBPS{2.1}{1.33}{\alpha_2\SH}
\put(2,.15){\vector(0,1){.36}}
\MBPS{2.1}{.33}{\tau_1}
\put(1.27,.91){\vector(3,-1){.53}}
\MBPS{1.6}{.9}{\tau_2}
\put(2.73,.91){\vector(-3,-1){.53}}
\MBPS{2.4}{.9}{\LAM\SH}
\end{pictr}
\end{center}
\BS

\NI where $m = k+j_1+j_2+2$ and
\[
\begin{aligned}
\mu\SH &= \HT{X}, \quad & \alpha_1\SH = \HT{A_1}, \quad \alpha_2\SH= \HT{A_2}\\
\tau_1 &= \HT{X \cup A_1}, \quad & T_1 \DFAS X \cup A_1= \tau_1[j_1+k+1] \\
\tau_2 &= \HT{X \cup A_2}, \quad & T_2 \DFAS X \cup A_2= \tau_2[j_2+k+1] \\
\LAM\SH &= \HT{A_1 \cup A_2}, \quad & A \DFAS A_1 \cup A_2 = \LAM\SH[j_1+j_2+1]
\end{aligned}
\]
When the Face Identity lemma is applied to the various commutative triangles in \DREF{comb-trio-diag-2} we get the following face identities:
\begin{eqnarray}
d_{A_1 \cup A_2} & = d_{\SHR{\LAM}_1[j_1]} d_{A_2} &= d_{\SHR{\LAM}_2[j_2]}\; d_{A_1}\\
d_{X \cup A_2} & = d_{\SHR{\LAM}_1[j_1]}\; d_{A_2} &= d_{\SHR{h}_2[j_2]}\; d_X \\
d_{X \cup A_1} & = d_{\SHR{h}_{12}[j_1]}\; d_X &= d_{\SHR{\mu}_2[k]}\; d_{A_1}
\end{eqnarray}
For example, $\HT{X}= \SHR{\mu} = \tau_1 \SHR{\mu}_1$. Therefore
\[
\begin{aligned}
d_{[m]-\SHR{\mu}[k]} &= d_{[k+j_1+1]-\SHR{\mu}_1[k]} d_{[m]-\tau_1[k+j_1+1]} \\
d_{\SHR{\LAM}[j_1+j_2+1]} &= d_{\SHR{\LAM}_1[j_1]}  d_{A_2}\\
d_{A_1 \cup A_2} &= d_{\SHR{\LAM}_1[j_1]}  d_{A_2}
\end{aligned}
\]
since $[j_1+k+1]-\mu_1{\SH}[k]=\SHR{\LAM}_1[j_1]$ and $[m]-\tau_1[j_1+k+1]=A_2$.
\MS

Now suppose $w \in C_m$ sponsors the trio $\ATRIO{x}{X}{y_1}{A_1}{y_2}{A_2}$. By definition this means:
\[
\begin{aligned}
\VL_w(x)=X & \iff & d_{A_1 \cup A_2}(w)=x \\
\VL_w(y_1)=A_1 & \iff & d_{X \cup A_2}(w)=y_1\\
\VL_w(y_1)=A_2 & \iff & d_{X \cup A_1}(w)=y_2
\end{aligned}
\]

\begin{center}
\begin{pictr}{2}{1.1}{-.5}{.5}
\put(0,0){\MB{(y_1,A_1)}}
\put(1,1){\MB{(x,X)}}
\put(2,0){\MB{(y_2,A_2)}}
\put(.2,.2){\line(1,1){.6}}
\put(.25,.53){\MBS{(\mu_1,\LAM_1)}}
\put(.7,.55){\MBS{z_1}}
\put(1.3,.55){\MBS{z_2}}
\put(.35,0){\line(1,0){1.3}}
\put(1.2,.8){\line(1,-1){.6}}
\put(1.8,.53){\MBS{(\mu_2,\LAM_2)}}
\end{pictr}
\end{center}
\MS

Observe that $w, d_{A_2}(w)$ and $d_{A_1}(w)$ all have $x$ as a subface. Therefore
\[
w \in C^x_{j_1+j_2+1}, \quad z_1 =d_{A_2}(w) \in C^x_{j_1}, \quad z_2=d_{A_1}(w) \in C^x_{j_2}
\]

Diagram \ref{sponsor-subface-diag} below portrays the face and subface relationships of the sponsor $w$. The symbol $\bullet \!\! \longrightarrow$ denotes (informally) a composition of face operators.

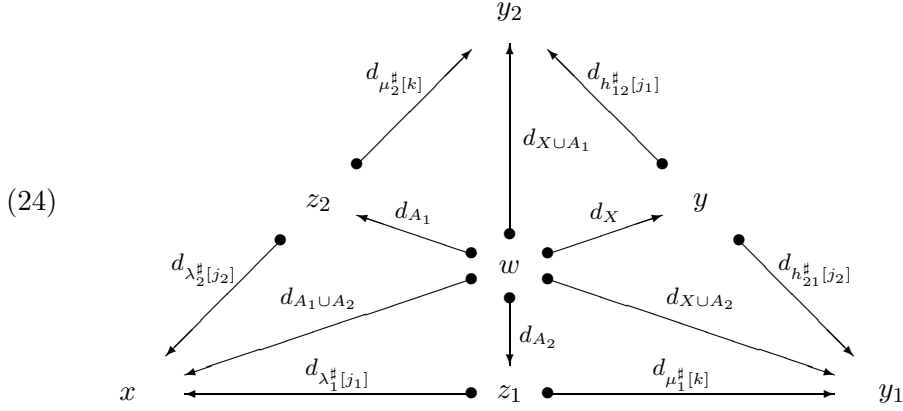
\begin{center}
\begin{pictr}{4}{2}{-.5}{1}
\label{sponsor-subface-diag}
\MBP{0}{0}{x}
\MBP{2}{0}{z_1}
\MBP{4}{0}{y_1}
\MBP{1}{1}{z_2}
\MBP{2}{2}{y_2}
\MBP{3}{1}{y}
\MBP{2}{.667}{w}
\put(1.8,.6){\vector(-3,-1){1.5}} 
\MBP{1.8}{.6}{\bullet}
\MBPS{1}{.5}{d_{A_1 \cup A_2}}

\put(2,.5){\vector(0,-1){.35}} 
\MBP{2}{.5}{\bullet}
\MBPS{2.15}{.3}{d_{A_2}}

\put(2.2,.6){\vector(3,-1){1.5}} 
\MBP{2.2}{.6}{\bullet}
\MBPS{3}{.5}{d_{X \cup A_2}}

\put(2.2,.734){\vector(3,1){.6}} 
\MBP{2.2}{.734}{\bullet}
\MBPS{2.5}{.95}{d_X}

\put(2,.834){\vector(0,1){1}} 
\MBP{2}{.834}{\bullet}
\MBPS{2.25}{1.33}{d_{X \cup A_1}}

\put(1.8,.734){\vector(-3,1){.6}} 
\MBP{1.8}{.734}{\bullet}
\MBPS{1.5}{.95}{d_{A_1}}

\put(1.8,0){\vector(-1,0){1.5}} 
\MBP{1.8}{0}{\bullet}
\MBPS{1.1}{.1}{d_{\LAM_1\SH[j_1]}}

\put(2.2,0){\vector(1,0){1.5}} 
\MBP{2.2}{0}{\bullet}
\MBPS{2.9}{.1}{d_{\mu_1\SH[k]}}

\put(.8,.8){\vector(-1,-1){.6}} 
\MBP{.8}{.8}{\bullet}
\MBPS{.4}{.65}{d_{\LAM_2\SH[j_2]}}

\put(1.2,1.2){\vector(1,1){.6}} 
\MBP{1.2}{1.2}{\bullet}
\MBPS{1.4}{1.65}{d_{\mu_2\SH[k]}}

\put(3.2,.8){\vector(1,-1){.6}} 
\MBP{3.2}{.8}{\bullet}
\MBPS{3.6}{.65}{d_{h_{21}\SH[j_2]}}

\put(2.8,1.2){\vector(-1,1){.6}} 
\MBP{2.8}{1.2}{\bullet}
\MBPS{2.6}{1.65}{d_{h_{12}\SH[j_1]}}

\end{pictr}
\end{center}

\vskip.5cm

For example, $y_1=d_{X \cup A_2}(w) = d_{\SHR{\mu}_1[k]}\; d_{A_2} (w) = d_{\SHR{\mu}_1[k]}(z_1)$.
\BS

\begin{newlem}
{\rm
\SL

Let $C$ be a simplicial set and $j_1,j_2, k \geq 0$. Suppose there are two complementary vertex functions pairs:
\[
\begin{aligned}
\DMAP{h_{12}}{j_1}{j_2+1}, \quad & \DMAP{h_{21}}{j_2}{j_1+1}\\
\DMAP{\mu}{k}{j_1+j_2+2}, \quad & \DMAP{\LAM}{j_1+j_2+1}{k+1}
\end{aligned}
\]
and $x \in C_k,\, y_1 \in C_{j_1},\, y_2 \in C_{j_2}$ and $y \in C_{j_1+j_2+1}$.

Given any $w \in \CPL(x,\mu;y,\LAM)$ such that $y \in \CPL(y_1,h_{12};y_2,h_{21})$ then $w$ and $y$ determine a $(k,j_1,j_2)$-comb-trio $\SETT{X,A_1,A_2}$ where
\[
X \DFAS \mu\SH[k], \quad A_1 \DFAS \LAM\SH h_{12}\SH[j_1], \quad A_2 \DFAS \LAM\SH h_{21}\SH[j_2]
\]
and $w$ is a sponsor of $\SETT{(x,X), (y_1,A_1), (y_2,A_2)}$.
}
\end{newlem}

\Proof

Let $m = j_1+j_2+k+2$. It follows immediately from the definitions that $X \cup A_1 \cup A_2=[m]$ is a disjoint union. The diagram of \SI\ functions for the  comb-trio is \DREF{comb-trio-diag-2} where
\[
\begin{aligned}
\alpha_1{\SH} \DFAS \LAM\SH \SHR{h}_1, \quad & \SHR{\alpha}_2 \DFAS \LAM\SH \SHR{h}_2\\
\tau_1[k+j_1+1]=X \cup A_1, \quad & \tau_2[k+j_2+1] = X \cup A_2\\
\SHR{\mu}_1[k] = \tau_1^{-1}(X), \quad & \SHR{\mu}_2[k] = \tau_2^{-1}(X)\\
\SHR{\LAM}_1[j_1] = \tau_1^{-1}(A_1), \quad & \SHR{\LAM}_2[j_2] = \tau_2^{-1}(A_2)
\end{aligned}
\]
To show that $w$ sponsors $\ATRIO{x}{X}{y_1}{A_1}{y_2}{A_2}$:
\MS

First, $w \in \CPL(x,\mu;y,\LAM)\IMP \VL_w(x) = \mu\SH[k]=X$.
\MS

Using that $y_1 = d_{h_{21}\SH[j_2]}(y)$, $h_{21}\SH[j_2] = [j_1+j_2+1]-h_1\SH[j_1]$, $y = d_{[m]-\mu\SH[k]}(w)$ and $\alpha_1\SH = \LAM\SH h_{12}\SH$, we get
\[
\begin{aligned}
y_1 & = d_{h_{21}\SH[j_2]}(y) \\
&= d_{[j_1+j_2+1]-h_{12}\SH[j_1]} \, d_{[m]-\mu\SH[k]}(w)\\
&= d_{[m]-A_1}(w)
\end{aligned}
\]
which shows that $\VL_w(y_1) = A_1$.
Similar reasoning shows $\VL_w(y_2)=A_2$.

\qed 
\BS

\begin{newthm}
{\bf (Trio Sponsor Theorem)}
\label{trio-sponsor-theorem}
\index{Trio Sponsor Theorem}
\index{Theorem!Trio Sponsor Theorem}
{\rm
\SL

Suppose $C$ is a simplicial set, $k,j_1,j_2 \geq 0$, $m=k+j_1+j_2+2$ and
\[
x\in C_k, \quad y_1 \in C_{j_1}, \quad y_2 \in C_{j_2}, \quad w \in C_m
\]
Then:
\begin{enumerate}
\item 
Let $\SETT{X,A_1,A_2}$ be a $(k,j_1,j_2)$-comb-trio with complementary vertex functions
\[
\begin{aligned}
\DMAP{\mu_1}{k}{j_1+1}, & \quad \DMAP{\LAM_1}{j_1}{k+1}\\
\DMAP{\mu_2}{k}{j_2+1}, & \quad \DMAP{\LAM_2}{j_2}{k+1}\\
\DMAP{h_{12}}{j_1}{j_2+1}, & \quad \DMAP{h_{21}}{j_2}{j_1+1}\\
\end{aligned}
\]
Then
\[
w \text{ sponsors } \ATRIO{x}{X}{y_1}{A_1}{y_2}{A_2} \IMP w \in \CPL^x(d_{A_2}(w),h_{12};d_{A_1}(w),h_{21})
\]

\item
Suppose $z_1 \in C^x_{j_1}, z_2 \in C^x_{j_2}$, and suppose
\[
\DMAP{h_{12}}{j_1}{j_2+1} \quad \text{and} \quad \DMAP{h_{21}}{j_2}{j_1+1}
\]
are complementary vertex functions. If $w \in \CPL^x(z_1,h_1;z_2,h_2)$ then there exists a uniquely determined comb-trio $\SETT{X,A_1,A_2}$, $y_1 \in C_{j_1}$ and $y_2 \in C_{j_2}$ such that $w$ sponsors $\ATRIO{x}{X}{y_1}{A_1}{y_2}{A_2}$.

\end{enumerate}

}
\end{newthm}

\Proof

The diagram for the trio $\SETT{X,A_1,A_2}$:
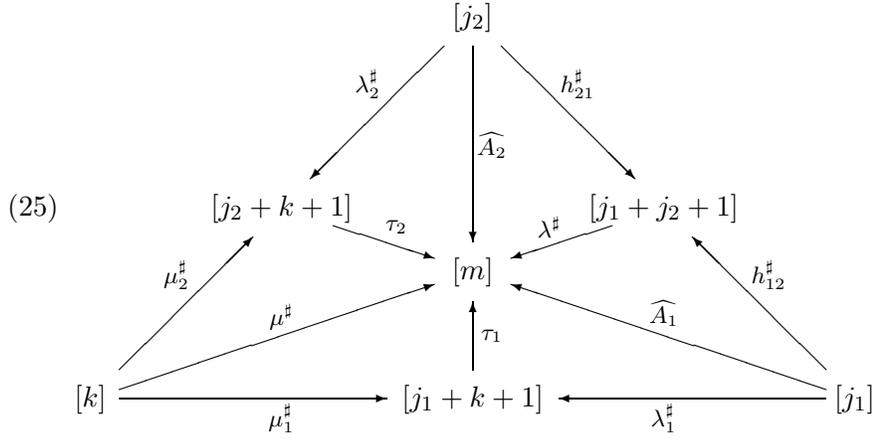
\begin{center}
\begin{pictr}{4}{2.1}{-.3}{1}
\label{comb-trio for trio sponsor}
\MBP{0}{0}{[k]}
\MBP{2}{0}{[j_1+k+1]}
\MBP{4}{0}{[j_1]}
\MBP{1}{1}{[j_2+k+1]}
\MBP{2}{2}{[j_2]}
\MBP{3}{1}{[j_1+j_2+1]}
\MBP{2}{.667}{[m]}
\put(.15,.15){\vector(1,1){.7}} 
\MBPS{.45}{.65}{\mu_2\SH}
\put(.15,.05){\vector(3,1){1.65}} 
\MBPS{1}{.43}{\mu\SH}
\put(.15,0){\vector(1,0){1.4}}
\MBPS{1}{-.1}{\mu_1\SH}
\put(3.85,.15){\vector(-1,1){.7}}
\MBPS{3.55}{.65}{h_{12}\SH} 
\put(3.85,.05){\vector(-3,1){1.65}} 
\MBPS{3}{.45}{\HT{A_1}}
\put(3.85,0){\vector(-1,0){1.4}}
\MBPS{3}{-.1}{\LAM_1\SH}
\put(1.85,1.85){\vector(-1,-1){.7}}
\MBPS{1.45}{1.65}{\LAM_2\SH}
\put(2.15,1.85){\vector(1,-1){.7}}
\MBPS{2.55}{1.65}{h_{21}\SH}
\put(2,1.85){\vector(0,-1){1.03}}
\MBPS{2.1}{1.33}{\HT{A_2}}
\put(2,.15){\vector(0,1){.36}}
\MBPS{2.1}{.33}{\tau_1}
\put(1.27,.91){\vector(3,-1){.53}}
\MBPS{1.6}{.9}{\tau_2}
\put(2.73,.91){\vector(-3,-1){.53}}
\MBPS{2.4}{.9}{\LAM\SH}
\end{pictr}
\end{center}
\BS

The schematic diagram for $w$:

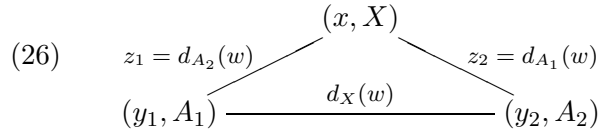
\begin{center}
\begin{pictr}{2}{.6}{-.7}{.3}
\put(0,0){\MB{(y_1,A_1)}}
\put(1,.5){\MB{(x,X)}}
\put(2,0){\MB{(y_2,A_2)}}
\put(.2,.1){\line(2,1){.6}}
\put(1.8,.1){\line(-2,1){.6}}
\put(.3,0){\line(1,0){1.4}}
\put(.1,.3){\MBS{z_1=d_{A_2}(w)}}
\put(1.9,.3){\MBS{z_2=d_{A_1}(w)}}
\put(1,.1){\MBS{d_X(w)}}
\end{pictr}
\end{center}
\BS

\NI Proof of claim 1.

First, $w \in \CPL(x,\mu;y,\LAM) \SBS C^x_{j_1+j_2+1}$ where 
\[
\mu\SH[k] = X, \quad \LAM\SH[j_1+j_2+1] = A_1 \cup A_2, \quad y = d_{\mu\SH[k]}(w)
\]
The claim that $w \in \CPL^x(d_{A_2}(w),h_{12};d_{A_1}(w),h_{21})$ says that $d^x_{h_{21}\SH[j_2]}(w) = z_1$ and, symmetrically, $d^x_{h_{12}\SH[j_1]}(w) = z_2$. By definition of face operators in $C^x$
\[
d^x_{h_{21}\SH[j_2]}(w) = d_{\LAM\SH h_{21}\SH[j_2]}(w) = d_{A_2}(w)
\]
as claimed. The same kind of calculation shows $d^x_{h_{21}\SH[j_2]}(w) =d_{A_1}(w)$.
\MS

\NI Proof of claim 2.

Assume $w \in \CPL^x(z_1,h_{12};z_2,h_{21})$ and let $y \in C_{j_1+j_2+1}$ denote the complementary subface of $x$ in $w$. As in part (1) above, let 
\[
\DMAP{\mu}{k}{j_1+j_2+2}, \quad \DMAP{\LAM}{j_1+j_2+1}{k+1}
\]
be the corresponding complementary vertex functions.  Then $\LAM\SH$ and $h_{12}\SH$ completely determine the comb-trio diagram (\ref{comb-trio for trio sponsor}) above. This, in turn, defines a trio $\SETT{X,A_1,A_2}$ where 
\[
X \DFAS \mu\SH[k],\quad A_1 \DFAS \LAM\SH h_{12}\SH[j_1],\quad A_2 \DFAS \LAM\SH h_{21}\SH[j_2]
\]
In particular, $d_{[m]-\mu\SH[k]}(w)=x$, which implies $\VL_w(x) = X$. 

Let $y_1 \DFAS d_{h_{21}\SH[j_2]}(y)$ and $y_2 \DFAS d_{h_{12}\SH[j_1]}(y)$; so that, by definition, $y \in \CPL(y_1,h_{12};y_2,h_{21})$.

Then, applying the Face Identity lemma
\[
\begin{aligned}
y_1 & = d_{h_{21}\SH[j_2]}(y) = d_{h_{21}\SH[j_2]}\; d_{\mu\SH[k]}(w)\\
&= d_{[j_1+j_2+1]-h_{12}\SH[j_1]}\; d_{[m]-\LAM\SH[j_1+j_2+1]}(w) =d_{[m]-\LAM\SH h_{12}\SH[j_1]}(w)\\
&= d_{[m]-A_1}(w)
\end{aligned}
\]
which implies $\VL_w(y_1)=A_1$. An identical calculation shows that $\VL_w(y_2)=A_2$. Therefore $w$ is the claimed sponsor.

\qed
\MS

\NI {\bf Note:} Given $w \in \CPL(x,\mu;y,\LAM)$, $y \in \CPL(y_1,h_1;y_2,h_2)$ and $z_1 \in \CPL(x,\mu_1;y_1,\LAM_1)$, then $y$ and $z_1$ resemble ``factors'' in the following schematic sense. Represent $w, y$ and $z_1$ schematically by 
\[
\SACP{x}{y}{w}, \quad \SACP{y_1}{y_2}{y} \quad \text{and} \quad \SACP{x}{y_1}{z_1}
\]
Then
\[
\begin{aligned}
\SACP{x}{y}{w} & \quad = \quad \SACP{x}{\left(\SACP{y_1}{y_2}{y}\right)}{w}\\
& \quad = \quad \SACP{\left(\SACP{x}{y_1}{z_1}\right)}{y_1}{w}
\end{aligned}
\]
\BX

\subsection{Universality in $C^{x,F}$}
\label{universality in $C^{x,F}$}

In this section, we will apply the concept of ``singular'' (in the sense of definition \refpage{mu-singular}) to the simplicial set $C^{x,F}$ (definition \refpage{C^x,F}) where $C$ and $D$ are simplicial sets, $F:D \to C$ is a simplicial map, $k \geq 0$ and $x \in C_k$. This will give a broad generalization of the ``universal map from $x$ to $F$'' concept when $C$ and $D$ are small categories and $F$ is a functor. Special cases will resemble the ``universal map from $x$ to $F$'' idea.
\MS

Recall that a $j$-simplex in $C^{x,F}$ is a pair $(z,v) \in C^x_j \times D_j$ such that $x$ and $F_j(v)$ are complementary subfaces of $z$. The face and degeneracy operators are 
\[
d^{x,F}_p(z,v) = (d^x_p(z), d_p(v))\quad (j>0) \quad \text{and} \quad s^{x,F}_p(z,v) = (s^x_p(z), s_p(v))
\]
\MS

\HD{Terminology:}

\begin{enumerate}
\item 
If $(z,v) \in C^{x,F}_j$ with $z \in \CPL(x,\mu;F_j(v),\LAM)$
we will say $\mu$ and $\LAM$ ``belong to'' or ``are associated with'' $(z,v)$.

\item
Suppose $(z_1,v_1) \in C^{x,F}_{j_1}$ has associated vertex function $\DMAP{\mu_1}{k}{j_1+1}$ and suppose $\DMAP{h_{12}}{j_1}{j_2+1}$ is a given non-decreasing function. Then we say 
\begin{quote}
\fbox{
{\bf $(z_2,v_2) \in C^{x,F}_{j_2}$ is compatible with $(z_1,v_1)$ and $h_{12}$}
}
\end{quote}
if the associated vertex function $\mu_2$ of $(z_2,v_2)$ is compatible with $\mu_1$ and $h_{12}$ in the sense of definition \refpage{compatible vertex function}.
\end{enumerate}
\BX
\BS

\begin{newdef}
{\bf (Singular in $C^{x,F}$)}
\label{singular in CxF}
{\rm
\SL

Let $F:D \to C$ be a simplicial map. Fix $k \geq 0$ and $x \in C_k$. Suppose $\DMAP{h_{12}}{j_1}{j_2+1}$ is non-decreasing and $\DMAP{h_{21}}{j_2}{j_1+1}$ is its complementary function. We say
\fbox{{\bf $(z_1,v_1) \in C^{x,F}_{j_1}$ is $h_{12}$-singular }}
if for all compatible $(z_2,v_2) \in C^{x,F}_{j_2}$ there exists a unique $(w,v) \in \CPL^{x,F}((z_1,v_1),h_{12};(z_2,v_2),h_{21})$.
}

\BX
\end{newdef}
\BS

With the Trio Sponsor Theorem \refpage{trio-sponsor-theorem} we can restate this in terms of the simplicial set $C$. For brevity, use the notations
\[
\begin{aligned}
y_1 \DFAS F_{j_1}(v_1), \quad & y_2 \DFAS F_{j_2}(v_2), \quad & y \DFAS F_{j_1+j_2+1}(v)\\
\end{aligned}
\]
Then $(w,v) \in C^{x,F}_{j_1+j_2+1}$, $(z_1,v_1) \in C^{x,F}_{j_1}$ and $(z_2,v_2) \in C^{x,F}_{j_2}$ translates to
\[
w \in \CPL(x,\mu;y,\LAM), \qquad z_1 \in \CPL(x,\mu_1;y_1,\LAM_1), \qquad z_2 \in \CPL(x,\mu_2;y_2,\LAM_2)
\]
for vertex functions $\mu$, $\mu_1$ and $\mu_2$ as in diagram \refpage{comb-trio-diag-2} 
where
\[
\mu\SH(t) \DFAS \mu_1\SH(t) + \mu_2\SH(t)-t
\]

Now $(w,v) \in \CPL^{x,F}\Bigl( (z_1,v_1),h_{12};(z_2,v_2),h_{21} \Bigr)$ implies
\[
\begin{aligned}
d^{x,F}_{h_{21}\SH[j_2]}(w,v) &= \left(d^{x}_{h_{21}\SH[j_2]}(w),d_{h_{21}\SH[j_2]}(v)\right) = (z_1,v_1)\\
d^{x,F}_{h_{12}\SH[j_1]}(w,v) &= \left(d^{x}_{h_{12}\SH[j_1]}(w),d_{h_{12}\SH[j_1]}(v)\right) = (z_2,v_2)
\end{aligned}
\]
and therefore $w \in \CPL^x(z_1,h_{12};z_2,h_{21})$,  $v_1 = d_{h_{21}\SH[j_2]}(v)$ and $v_2 = d_{h_{12}\SH[j_1]}(v)$.
Also
\[
d_{h_{21}\SH[j_2]}(y) = d_{h_{21}\SH[j_2]}(F_{j_1+j_2+1}(v)) =
F_{j_1}( d_{h_{21}\SH[j_2]}(v) ) = F_{j_1}(v_1)
\]
and, similarly, $d_{h_{12}\SH[j_1]}(y) = y_2$. Therefore
\[
y = F_{j_1+j_2+1}(v) \in \CPL(y_1,h_{12};y_2,h_{21})
\]
By part 2 of the Trio Sponsor Theorem above, $w \in C_{j_1+j_2+k+2}$ sponsors the trio $\ATRIO{x}{X}{F_{j_1}(v_1)}{A_1}{F_{j_2}(v_2)}{A_2}$ where
\[
X = \mu\SH[k], \quad A_1 = \LAM\SH h_{12}\SH[j_1], \quad A_2 = \LAM\SH h_{21}\SH[j_2]
\]
and the schematic picture is

\begin{center}
\begin{pictr}{2}{1.5}{-1}{.85}
\label{universal factor schematic general}
\put(0,.5){\MB{(F_{j_1}(v_1),A_1)}}
\put(1,1.5){\MB{(x,X)}}
\put(2,.5){\MB{(F_{j_2}(v_2),A_2)}}
\put(.5,.5){\line(1,0){1}}
\put(.2,.7){\line(1,1){.6}}
\put(1.8,.7){\line(-1,1){.6}}
\put(.3,1.1){\MBS{(z_1,v_1)}}
\put(2.3,1.1){\MBS{(z_2,v_2) \text{ (given compatible)}}}
\put(1,.35){\MBS{y = F_{j_1+j_2+1}(v)}}
\put(1,.9){\MBS{\exists ! (w,v)}}
\put(1,0){\MB{d_{A_2}(w)=z_1, \quad d_{A_1}(w)=z_2, \quad d_{X}(w)=y=F_{j_1+j_2+1}(v)}}
\end{pictr} 
\end{center}
\MS

\BX
\MS

\begin{example}
{\rm
It is instructive to look at the lowest dimension case: when $x \in C_0$. For $j \geq 0$, $(z,v) \in C^{x,F}_j$ means that $v \in D_j$ and $z \in \CPL(x,\mu;F_j(v),\LAM) \SBS C_{j+1}$ where $\DMAP{\mu}{0}{j+1}$ and $\DMAP{\LAM}{j}{1}$. In this special case, $\mu\SH = \mu$.
Simplices of $C^{x,F}$ in dimensions 0 and 1 are as follows.
\MS

\HD{Dimension $j=0$:}

Given $(z,v) \in C^{x,F}_0$ where $\DMAP{\mu}{0}{1}$ and $\DMAP{\LAM}{0}{1}$ are the vertex functions for the subfaces $x$ and $F_0(v)$ in $z$ then either $\mu\SH[0]=\{0\}$ or $\mu\SH[0]=\{1\}$.
\[
\begin{aligned}
\mu\SH(0)=0 & \IMP \LAM\SH(0)=1, & \; \text{and} \quad z = x \longrightarrow F_0(v)\\
\mu\SH(0)=1 & \IMP \LAM\SH(0)=0, & \; \text{and} \quad z = F_0(v) \longrightarrow x
\end{aligned}
\]

\HD{Dimension $j=1$;}

Here $(z,v) \in C^{x,F}_1$ means $v \in D_1$ and $z \in \CPL(x,\mu;F_1(v),\LAM) \SBS C_2$ where $\DMAP{\mu}{0}{2}$ and $\DMAP{\LAM}{1}{1}$. The three possibilities, with $y \DFAS F_1(v)$, are:
\BS

\setlength{\unitlength}{1cm}
\[
\mu\SH[0]=\{0\}, \; \LAM\SH[1]=\{1,2\}, \qquad
\makebox{
\begin{picture}(1.3,1)
\put(0,1){\MB{x}}
\put(.25,1){\vector(1,0){.5}}
\put(0,.75){\vector(0,-1){.5}}
\put(1,1){\MB{\bullet_1}}
\put(.75,.75){\vector(-1,-1){.5}}
\put(0,0){\MB{\bullet_2}}
\put(.7,.4){\MBS{y}}
\put(1.5,.5){\MB{=z}}
\end{picture}
}
\]
\MS

\[
\mu\SH[0]=\{1\}, \; \LAM\SH[1]=\{0,2\}, \qquad
\makebox{\begin{picture}(1.3,1)
\put(0,1){\MB{\bullet_0}}
\put(.25,1){\vector(1,0){.5}}
\put(0,.75){\vector(0,-1){.5}}
\put(1,1){\MB{x}}
\put(.75,.75){\vector(-1,-1){.5}}
\put(0,0){\MB{\bullet_2}}
\put(-.2,.5){\MBS{y}}
\put(1.5,.5){\MB{=z}}
\end{picture}}
\]
\MS

\[
\mu\SH[0]=\{2\}, \; \LAM\SH[1]=\{0,1\}, \qquad
\makebox{\begin{picture}(1.3,1)
\put(0,1){\MB{\bullet_0}}
\put(.25,1){\vector(1,0){.5}}
\put(0,.75){\vector(0,-1){.5}}
\put(1,1){\MB{\bullet_1}}
\put(.75,.75){\vector(-1,-1){.5}}
\put(0,0){\MB{x}}
\put(.5,1.2){\MBS{y}}
\put(1.5,.5){\MB{=z}}
\end{picture}}
\]
\MS

Next, we consider compatibility in the sense of definition \refpage{compatible vertex function}. 

The general setup here is: $k=0$, $j_1=j_2=j$ and $(z_1,v_1) \in C^{x,F}_j$ has vertex functions $\mu_1$ and $\LAM_1$. Let $\DMAP{h_{12}}{j}{j+1}$ and $\DMAP{h_{21}}{j}{j+1}$ be any pair of complementary vertex functions. The question is what $(z_2,v_2) \in C^{x,F}_j$ (with vertex functions $\mu_2$ and $\LAM_2$) are compatible with $(z_1,v_1)$ and the vertex function pair $h_{12}$ and $h_{21}$. This amounts to determining if there is a comb-trio whose vertex function pairs are $(\mu_1,\LAM_1)$, $(h_{12},h_{21})$ and $(\mu_2,\LAM_2)$. \MS

We'll examine this when $j=0$.
\MS

Let $(z_1,v_1) \in C^{x,F}_0$, $z_1 \in \CPL(x,\mu_1;F_0(v_1),\LAM_1)$. Suppose, as a specific example, that $\mu_1\SH[0] = \{0\}$, $\LAM_1\SH[0]=\{1\}$. We will examine {\em both} possibilities for $\DMAP{h_{12}\SH}{0}{1}$. For each possibility we will use the given choice of $\mu_1\SH$ and $h_{12}\SH$ to fill in the comb-trio diagram:

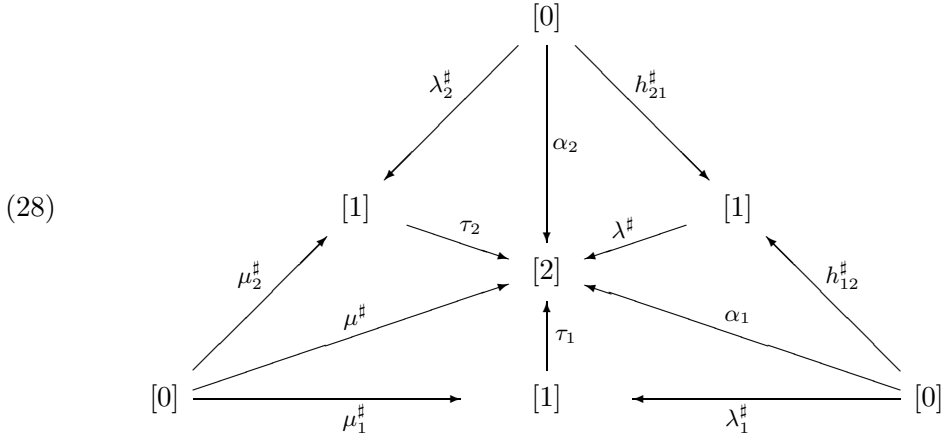
\begin{center}
\setlength{\unitlength}{1in}
\begin{pictr}{4}{2}{-.7}{1}
\label{low dim comb trio diag}
\MBP{0}{0}{[0]}
\MBP{2}{0}{[1]}
\MBP{4}{0}{[0]}
\MBP{1}{1}{[1]}
\MBP{2}{2}{[0]}
\MBP{3}{1}{[1]}
\MBP{2}{.667}{[2]}
\put(.15,.15){\vector(1,1){.7}} 
\MBPS{.45}{.65}{\mu_2\SH}
\put(.15,.05){\vector(3,1){1.65}} 
\MBPS{1}{.43}{\mu\SH}
\put(.15,0){\vector(1,0){1.4}}
\MBPS{1}{-.1}{\mu_1\SH}
\put(3.85,.15){\vector(-1,1){.7}}
\MBPS{3.55}{.65}{h_{12}\SH} 
\put(3.85,.05){\vector(-3,1){1.65}} 
\MBPS{3}{.43}{\alpha_1}
\put(3.85,0){\vector(-1,0){1.4}}
\MBPS{3}{-.1}{\LAM_1\SH}
\put(1.85,1.85){\vector(-1,-1){.7}}
\MBPS{1.45}{1.65}{\LAM_2\SH}
\put(2.15,1.85){\vector(1,-1){.7}}
\MBPS{2.55}{1.65}{h_{21}\SH}
\put(2,1.85){\vector(0,-1){1.03}}
\MBPS{2.1}{1.33}{\alpha_2}
\put(2,.15){\vector(0,1){.36}}
\MBPS{2.1}{.33}{\tau_1}
\put(1.27,.91){\vector(3,-1){.53}}
\MBPS{1.6}{.9}{\tau_2}
\put(2.73,.91){\vector(-3,-1){.53}}
\MBPS{2.4}{.9}{\LAM\SH}
\end{pictr}
\end{center}
\BS

\MS

\HD{Possibility \#1:  $h_{12}\SH[0]=\{0\}$, $h_{21}\SH[0]=\{1\}$}
\MS
 
\HD{Step 1:}
$\alpha_1(0) = \LAM_1\SH(0) + h_{12}\SH(0)=0=1$ which implies $A_1=\{1\}$ and $\tau_2 = (0,2)$.
\MS
 
\HD{Step 2:} 
$1 = \alpha_1(0) = \tau_1 \LAM_1\SH(0) = \tau_1(1)$. Since $\tau$ is \SI\ there is only one possibility: $\tau_1=(0,1)$. It follows that $\alpha_2(0)=2$, $A_2=\{2\}$, $X = \{0\}$ and $\mu\SH(0)=0$.
\MS

\HD{Step 3:}
$0 = \mu\SH(0) = \mu_1\SH(0)+\mu_2\SH(0)-0 = \mu_2\SH(0)$.
\MS

\HD{Summary:} 
Given $(z_1,v_1) \in C^{x,F}_0$, $z_1 \in \CPL(x,\mu_1;F_0(v_1),\LAM_1)$ and given that $h_{12}(0)=0$ then $(z_2,v_2) \in C^{x,F}_0$, $z_2 \in \CPL(x,\mu_2;F_0(v_2),\LAM_2)$ is compatible if $\DMAP{\mu_2}{0}{1}$ is $\mu_2(0)=0$.

To say that $(z_1,v_1)$ is $h_{12}$-singular is to say that for all compatible $(z_2,v_2) \in C^{x,F}_0$ there exists a unique $(w,v) \in C^{x,F}_1$ where 
\[
(w,v) \in \CPL^{x,F}((z_1,v_1),h_{12};(z_2,v_2),h_{21}), \qquad w \in \CPL^x(z_1,h_{12};z_2,h_{21})
\]
with $d_{A_2}(w) = d_2(w)=z_1$, $d_{A_1}(w) = d_1(w)=z_2$ and $d_X(w) = d_0(w) = F_1(v)$. In terms of simplices in $C$: 

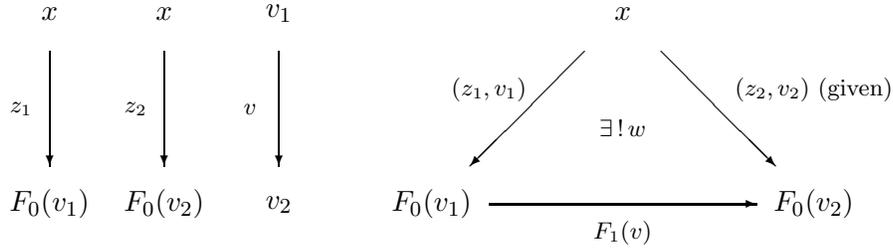
\setlength{\unitlength}{1in}
\begin{center}
\begin{picture}(4,1.1)
\put(0,1){\MB{x}}
\put(0,0){\MB{F_0(v_1)}} 
\put(0,.8){\vector(0,-1){.6}}
\put(-.15,.5){\MBS{z_1}}

\put(0.6,1){\MB{x}}
\put(0.6,0){\MB{F_0(v_2)}} 
\put(.6,.8){\vector(0,-1){.6}}
\put(.45,.5){\MBS{z_2}}

\put(1.2,1){\MB{v_1}}
\put(1.2,0){\MB{v_2}} 
\put(1.2,.8){\vector(0,-1){.6}}
\put(1.05,.5){\MBS{v}}

\put(2,0){\MB{F_0(v_1)}}
\put(4,0){\MB{F_0(v_2)}}
\put(3,1){\MB{x}}
\put(2.8,.8){\vector(-1,-1){.6}}
\put(3.2,.8){\vector(1,-1){.6}}
\put(2.3,0){\vector(1,0){1.4}}
\put(2.3,.6){\MBS{(z_1,v_1)}}
\put(4,.6){\MBS{(z_2,v_2) \text{ (given)}}}
\put(3,-.15){\MBS{F_1(v)}}
\put(3,.4){\MBS{\exists \, !\, w}}
\end{picture}
\end{center}
\MS

\NI which is the familiar case of $x \XRA{(z_1,v_1)} F_0(v_1)$ being the universal map from $x$ to $F$ when $C$ is a category and $F:D \to C$ a functor.
\BS

\HD{Possibility \#2: $h_{12}\SH[0]=\{1\}$, $h_{21}\SH[0]=\{0\}$}
\MS

\HD{Step 1:}
$\alpha_1(0) = \LAM_1\SH(0) + h_{12}\SH(0)-0 = 2$, which implies $A_1=\{2\}$ and $\tau_2=(0,1)$.
\MS

\HD{Step 2:}
$2 = \alpha_1(0) = \tau_2 \LAM_1\SH(0) = \tau_1(1)$. Thus there are {\em two} possible choices for $\tau_1$, namely $\tau_1=(0,2)$ (``choice A'') and $\tau_1=(1,2)$ (``choice B'').
\MS

\HD{Step 3A:}
$\tau_1=(0,2)$ which implies $A_2=\{1\}$ and $X = \{0\}$. Then
\[
0 = \mu\SH(0)=\mu_1\SH(0)+ \mu_2\SH(0)-0=\mu_2\SH(0)
\]

\HD{Summary ``A'':}
Given $(z_1,v_1) \in C^{x,F}_0$, $z_1 \in \CPL(x,\mu_1;F_0(v_1),\LAM_1)$ and given that $h_{12}(0)=1$ (and $\tau_1=(0,2)$) then $(z_2,v_2) \in C^{x,F}_0$, $z_2 \in \CPL(x,\mu_2;F_0(v_2),\LAM_2)$ is compatible if $\DMAP{\mu_2}{0}{1}$ is $\mu_2(0)=0$.

To say that $(z_1,v_1)$ is $h_{12}$-singular is to say that for all compatible $(z_2,v_2) \in C^{x,F}_0$ there exists a unique $(w,v) \in C^{x,F}_1$ where 
\[
(w,v) \in \CPL^{x,F}((z_1,v_1),h_{12};(z_2,v_2),h_{21}), \qquad w \in \CPL^x(z_1,h_{12};z_2,h_{21})
\]
with $d_{A_2}(w) = d_1(w)=z_1$, $d_{A_1}(w) = d_2(w)=z_2$ and $d_X(w) = d_0(w) = F_1(v)$. The simplex diagram in $C$ is:

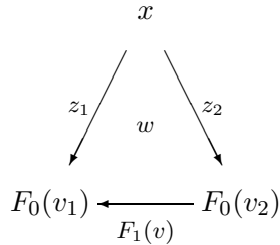
\begin{center}
\begin{picture}(1,1)
\put(0,0){\MB{F_0(v_1)}}
\put(1,0){\MB{F_0(v_2)}}
\put(.5,1){\MB{x}}
\put(.4,.8){\vector(-1,-2){.3}}
\put(.6,.8){\vector(1,-2){.3}}
\put(.75,0){\vector(-1,0){.5}}
\put(.15,.5){\MBS{z_1}}
\put(.85,.5){\MBS{z_2}}
\put(.5,-.15){\MBS{F_1(v)}}
\put(.5,.4){\MBS{w}}
\end{picture}
\end{center}
\MS

\HD{Step 3B:}
$\tau_1=(1,2)$ which implies $A_2=\{0\}$, $X = \{1\}$ and $A_1=\{2\}$. Then
\[
1 = \mu\SH(0)=\mu_1\SH(0)+ \mu_2\SH(0)-0=\mu_2\SH(0)
\]
and $\tau_3=(0,2)$
\MS

\HD{Summary ``B'':}
Given $(z_1,v_1) \in C^{x,F}_0$, $z_1 \in \CPL(x,\mu_1;F_0(v_1),\LAM_1)$ and given that $h_{12}(0)=1$ (and $\tau_1=(0,2)$) then $(z_2,v_2) \in C^{x,F}_0$, $z_2 \in \CPL(x,\mu_2;F_0(v_2),\LAM_2)$ is compatible if $\DMAP{\mu_2}{0}{1}$ is $\mu_2(0)=1$.

To say that $(z_1,v_1)$ is $h_{12}$-singular is to say that for all compatible $(z_2,v_2) \in C^{x,F}_0$ there exists a unique $(w,v) \in C^{x,F}_1$ where 
\[
(w,v) \in \CPL^{x,F}((z_1,v_1),h_{12};(z_2,v_2),h_{21}), \qquad w \in \CPL^x(z_1,h_{12};z_2,h_{21})
\]
with $d_{A_0}(w) = d_0(w)=z_1$, $d_{A_2}(w) = d_1(w)=z_2$ and $d_X(w) = d_1(w) = F_1(v)$. The simplex diagram in $C$ is:

\begin{center}
\begin{picture}(1,1)
\put(0,0){\MB{F_0(v_1)}}
\put(1,0){\MB{x}}
\put(.5,1){\MB{F_0(v_2)}}
\put(.4,.8){\vector(-1,-2){.3}}
\put(.6,.8){\vector(1,-2){.3}}
\put(.75,0){\vector(-1,0){.5}}
\put(.05,.5){\MBS{F_1(v)}}
\put(.85,.5){\MBS{z_2}}
\put(.5,-.15){\MBS{z_1}}
\put(.5,.4){\MBS{w}}
\end{picture}
\end{center}
\MS

Note that in this case, the location of $x$ as a subface of $z_1$ is different from that in $z_2$. This is a consequence of $h_{12}$.
}
\end{example}

\subsubsection{Universal $(n,i)$-factor}

The definition of ``singular'' in $C^{x,F}$ for a given simplicial map $F:D \to C$ involves a general choice of:

$\BLT$\; numerical parameters $k, j_1, j_2 \geq 0$

$\BLT$\; a simplex $x \in C_k$

$\BLT$\; a non-decreasing function $\DMAP{h_{12}}{j_1}{j_2+1}$
\MS

In this section we will specialize these to choices which seems apt in the context of $(n,i)$-composition.
\MS

Referring back to definition \refpage{singular in CxF} of ``singular in $C^{x,F}$\,'', fix notation as follows: a comb-trio $\SETT{X,A_1,A_2}$ of dimension $m \DFAS j_1+j_2+k+2$ with complementary vertex functions
\[
\begin{aligned}
\DMAP{\mu_1}{k}{j_1+1}, \qquad & \DMAP{\LAM_1}{j_1}{k+1}\\
\DMAP{\mu_2}{k}{j_2+1}, \qquad & \DMAP{\LAM_2}{j_2}{k+1}\\
\DMAP{h_{12}}{j_1}{j_2+1}, \qquad & \DMAP{h_{21}}{j_2}{j_1+1}
\end{aligned}
\]
and the trio $\ATRIO{x}{X}{F_{j_1}(v_1)}{A_1}{F_{j_2}(v_2)}{A_2}$ where $x \in C_k$, $v_1 \in D_{j_1}$ and $v_2 \in D_{j_2}$. The definition of ``singular'' posits the existence of a unique sponsor $(w,v) \in C^{x,F}_{j_1+j_2+1}$ where $d_{A_2}(w,v)=(z_1,v_1)$, $d_{A_1}(w,v) = (z_2,v_2)$ and $d_X(w,v) = F_{j_1+j_2+1}(v)$. The schematic diagram for such a sponsor would be:

\begin{center}
\begin{pictr}{2}{1}{-1}{.5}
\label{trio for sing discussion}
\put(0,0){\MB{(F_{j_1}(v_1),A_1)}}
\put(2,0){\MB{(F_{j_2}(v_2),A_2)}}
\put(1,1){\MB{(x,X)}}
\put(.2,.2){\line(1,1){.6}}
\put(1.8,.2){\line(-1,1){.6}}
\put(.45,0){\line(1,0){1}}
\put(.4,.7){\MBS{(z_1,v_1)}}
\put(1.6,.7){\MBS{(z_2,v_2)}}
\put(1,-.1){\MBS{F_{j_1+j_2+1}(v)}}
\put(1,.4){\MBS{(w,v)}}
\end{pictr}
\end{center}
\BS

Abbreviate $F_{j_1}(v_1)$ by $y_1 \in C_{j_1}$, $F_{j_2}(v_2)$ by $y_2 \in C_{j_2}$ and $F_{j_1+j_2+1}(v)$ by $y \in C_{j_1+j_2+1}$.
\MS

\NI Recall that
\[
z_1 \in \CPL(x,\mu_1;y_1,\LAM_1), \quad z_2 \in \CPL(x,\mu_2;y_2,\LAM_2), \quad
y \in \CPL(y_1,h_{12};y_2,h_{21})
\]
It follows from the discussion in section \refpage{trios in C^x} that $w \in \CPL(x,\mu;y,\LAM)$ where $\DMAP{\mu\SH}{k}{m}$ is defined by 
\[
\mu\SH(t) \DFAS \mu_1\SH(t)+\mu_2\SH(t)-t
\]
and that
\[
\begin{aligned}
w  & \in \CPL^x(z_1,h_{12};z_2,h_{21})\\
(w,v) & \in \CPL^{x,F}((z_1,v_1),h_{12};(z_2,v_2),h_{21})
\end{aligned}
\]
\MS

In order to motivate the proposed specialization below, first consider the next example which connects the general description of ``singular'' with the familar definition of ``universal map from $x$ to $F$''.
\MS

\begin{example}
{\rm
Let $C$ and $D$ be small categories i.e. \NICOMP{1}{1}s, $k=0$, $x \in C_0$, $j_1=j_2=0$ and let $\mu_1,\mu_2:[0] \to [1]$ be $\mu_1(0)=\mu_2(0)\DFAS 0$. Then $(z_1,v_1) \in C^{x,F}_0$ is the 1-simplex $x \XRA{z_1} F_0(v_1)$ in $C$. If $\DMAP{h_{12}}{0}{1}$ is $h_{12}(0)\DFAS 0$ then $x \XRA{z_1} F_0(v_1)$ is ``the universal map from $x$ to $F$'' iff for all $(z_2,v_2) \in C^{x,F}_0$, (that is, $x \XRA{z_2} F_0(v_2)$) there is a unique 
\[
(w,v) \in C^{x,F}_1, \qquad (w,v) = \CPL^{x,F}(z_1,v_1),h_{12};(z_2,v_2),h_{21})
\]
whose schematic trio diagram is 

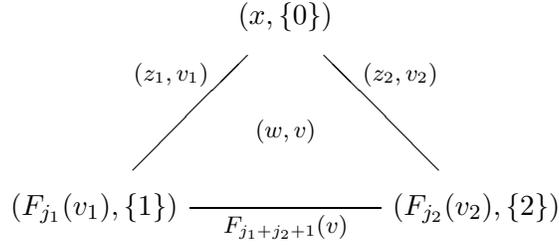
\begin{center}
\begin{picture}(2,1)
\put(0,0){\MB{(F_{j_1}(v_1),\{1\})}}
\put(2,0){\MB{(F_{j_2}(v_2),\{2\})}}
\put(1,1){\MB{(x,\{0\})}}
\put(.2,.2){\line(1,1){.6}}
\put(1.8,.2){\line(-1,1){.6}}
\put(.5,0){\line(1,0){1}}
\put(.4,.7){\MBS{(z_1,v_1)}}
\put(1.6,.7){\MBS{(z_2,v_2)}}
\put(1,-.1){\MBS{F_{j_1+j_2+1}(v)}}
\put(1,.4){\MBS{(w,v)}}
\end{picture}
\end{center}
\MS

\NI and whose diagram in $C$ is 

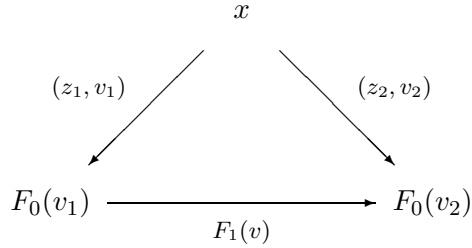
\begin{center}
\begin{picture}(2,1)
\put(1,1){\MB{x}}
\put(0,0){\MB{F_0(v_1)}}
\put(2,0){\MB{F_0(v_2)}}
\put(.8,.8){\vector(-1,-1){.6}}
\put(1.2,.8){\vector(1,-1){.6}}
\put(.3,0){\vector(1,0){1.4}}
\put(.2,.6){\MBS{(z_1,v_1)}}
\put(1.8,.6){\MBS{(z_2,v_2)}}
\put(1,-.15){\MBS{F_1(v)}}
\end{picture}
\end{center}
\SL
}
\end{example}
\MS

This example suggests a version of ``universal'' appropriate for the case when $C$ and $D$ are \NICOMP{n}{i}s. We require both $(z_1,v_1) \in C^{x,F}_{j_1} \SBS C_{k+j_1+1}$ and $(z_2,v_2) \in C^{x,F}_{j_2} \SBS C_{k+j_2+1}$ to be $(n,i)$-factors; i.e $j_1+k+1=j_2+k+1=n$. We also require that $\mu_1=\mu_2:[k] \to [j+1]$ where $z_1 \in \CPL(x,\mu_1;F_j(v_1),\LAM_1)$ and $z_2 \in \CPL(x,\mu_2;F_j(v_2),\LAM_2)$. 
\MS

\NI That is $j_1=j_2=j$ and $j+k+1=n$. The corresponding comb-trio $\SETT{X,A_1,A_2}$ of dimension $m \DFAS 2j+k+2 = n+j+1$ has complementary vertex function diagram:

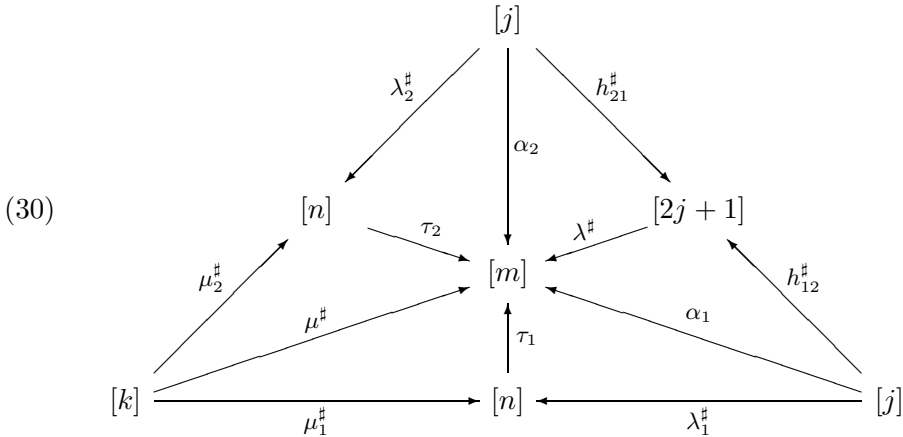
\begin{center}
\begin{pictr}{4}{2}{-.5}{1}
\label{k,j,j comb-trio}
\MBP{0}{0}{[k]}
\MBP{2}{0}{[n]}
\MBP{4}{0}{[j]}
\MBP{1}{1}{[n]}
\MBP{2}{2}{[j]}
\MBP{3}{1}{[2j+1]}
\MBP{2}{.667}{[m]}
\put(.15,.15){\vector(1,1){.7}} 
\MBPS{.45}{.65}{\mu_2\SH}
\put(.15,.05){\vector(3,1){1.65}} 
\MBPS{1}{.43}{\mu\SH}
\put(.15,0){\vector(1,0){1.7}}
\MBPS{1}{-.1}{\mu_1\SH}
\put(3.85,.15){\vector(-1,1){.7}}
\MBPS{3.55}{.65}{h_{12}\SH} 
\put(3.85,.05){\vector(-3,1){1.65}} 
\MBPS{3}{.45}{\alpha_1}
\put(3.85,0){\vector(-1,0){1.7}}
\MBPS{3}{-.1}{\LAM_1\SH}
\put(1.85,1.85){\vector(-1,-1){.7}}
\MBPS{1.45}{1.65}{\LAM_2\SH}
\put(2.15,1.85){\vector(1,-1){.7}}
\MBPS{2.55}{1.65}{h_{21}\SH}
\put(2,1.85){\vector(0,-1){1.03}}
\MBPS{2.1}{1.33}{\alpha_2}
\put(2,.15){\vector(0,1){.36}}
\MBPS{2.1}{.33}{\tau_1}
\put(1.27,.91){\vector(3,-1){.53}}
\MBPS{1.6}{.9}{\tau_2}
\put(2.73,.91){\vector(-3,-1){.53}}
\MBPS{2.4}{.9}{\LAM\SH} 
\end{pictr}
\end{center}
\MS

\NI where $X = \mu\SH[k]$, $A_1 = \alpha_1[j]$ and $A_2 = \alpha_2[j]$.
\MS

Starting with a given $(z_1,v_1) \in C^{x,F}_j$ asserted to have a singularity property, one might regard the compatibility requirement in one of the following ways:
\begin{enumerate}
\item 
Specify a non-decreasing $\DMAP{h_{12}}{j}{j+1}$ and assert that $(z_1,v_1)$ is $h_{12}$-singular with respect to all compatible $(z_2,v_2)$.

\item
Specify a family of $j$-simplices $(z_2,v_2)$ of $C^{x,F}$ (i.e. certain $n$-simplices of $C$) with respect to which the given $\DMAP{h_{12}}{j}{j+1}$ is compatible i.e. there exists for each $(z_2,v_2)$ in the family, a comb-trio containing the vertex functions $\DMAP{\mu_1}{k}{j+1}$ (arising from $z_1$), $\DMAP{\mu_2}{k}{j+1}$ (arising from $z_2$) and $h_{12}$.

In this case, the singular property would be that $(z_1,v_1)$ is $h_{12}$-singular with respect to the specified family of $j$-simplices in $C^{x,F}$.
\end{enumerate}
\MS

\NI Notes:

A $(k,j,j)$-comb trio as in diagram \eqref{k,j,j comb-trio} above has 
\[
m=k+2j+2=n+j+1
\]
Therefore $m \equiv k \mod 2$ and $n-k=j+1$.
Given $\DMAP{\mu_1}{k}{k+j+1}$, $\DMAP{\mu_2}{k}{k+j+1}$ with $\mu_1=\mu_2$, then
$\DMAP{\mu\SH}{k}{n+1}$, defined by
\[
\mu\SH(t) \DFAS \mu_1\SH(t) + \mu_2\SH(t)-t = 2 \mu_1\SH(t)-t
\]
has the property that its values are alternately even then odd. (This is not true in general of $\LAM\SH$ where $\DMAP{\LAM}{j}{k+1}$ is the complementary vertex function of $\mu$).

\BX
\MS

\begin{example}
{\rm
Let $C$ and $D$ be \NICOMP{n}{i}s and $F:D \to C$ a simplicial map. Fix $x \in C_{n-1}$ and suppose $(z_1,v_1), (z_2,v_2) \in C^{x,F}_0$ where $z_1 \in \CPL(x,\mu_1;F_0(v_1),\LAM_1)$ and $z_2 \in \CPL(x,\mu_2;F_0(v_1),\LAM_2)$. We will consider a trio arising from the case when $\mu_1\SH = \mu_2\SH \DFAS \DMAP{\del_i}{n-1}{n}$. In that case, 
\[
\LAM_1\SH = \LAM_2\SH = \DMAP{(i)}{0}{n}
\]

Diagram \ref{k,j,j comb-trio} becomes:
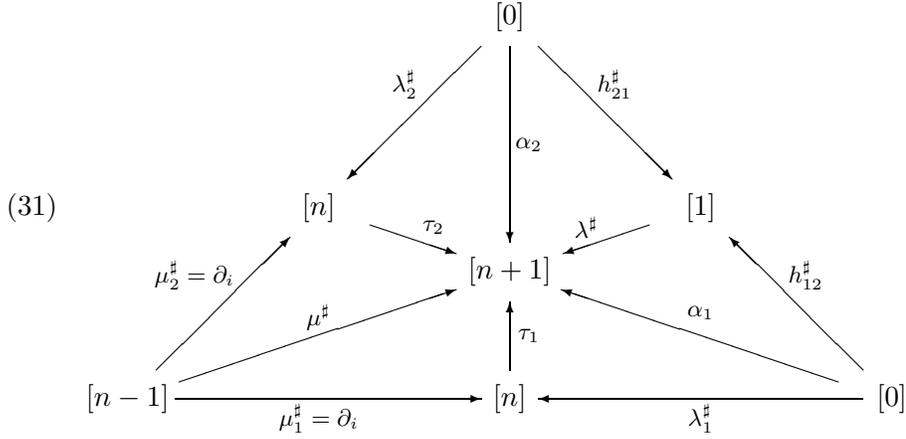
\begin{center}
\begin{pictr}{4}{2}{-.5}{1}
\label{n-1,0,0 comb-trio}
\MBP{0}{0}{[n-1]}
\MBP{2}{0}{[n]}
\MBP{4}{0}{[0]}
\MBP{1}{1}{[n]}
\MBP{2}{2}{[0]}
\MBP{3}{1}{[1]}
\MBP{2}{.667}{[n+1]}
\put(.15,.15){\vector(1,1){.7}} 
\MBPS{.35}{.65}{\mu_2\SH=\del_i}
\put(.27,.09){\vector(3,1){1.45}} 
\MBPS{1}{.43}{\mu\SH}
\put(.25,0){\vector(1,0){1.6}}
\MBPS{1}{-.1}{\mu_1\SH=\del_i}
\put(3.85,.15){\vector(-1,1){.7}}
\MBPS{3.55}{.65}{h_{12}\SH} 
\put(3.72,.09){\vector(-3,1){1.45}} 
\MBPS{3}{.45}{\alpha_1}
\put(3.85,0){\vector(-1,0){1.7}}
\MBPS{3}{-.1}{\LAM_1\SH}
\put(1.85,1.85){\vector(-1,-1){.7}}
\MBPS{1.45}{1.65}{\LAM_2\SH}
\put(2.15,1.85){\vector(1,-1){.7}}
\MBPS{2.55}{1.65}{h_{21}\SH}
\put(2,1.85){\vector(0,-1){1.03}}
\MBPS{2.1}{1.33}{\alpha_2}
\put(2,.15){\vector(0,1){.36}}
\MBPS{2.1}{.33}{\tau_1}
\put(1.27,.91){\vector(3,-1){.45}}
\MBPS{1.6}{.9}{\tau_2}
\put(2.73,.91){\vector(-3,-1){.45}}
\MBPS{2.4}{.9}{\LAM\SH} 
\end{pictr}
\end{center}
\MS

It follows from the general formula for $\mu\SH$ that
\[
\begin{aligned}
\mu\SH & = (0 \DDD{,} i-1,i+2 \DDD{,} n+1):[n-1] \to [n+1]\\
\LAM\SH &= (i,i+1):[1] \to [n+1]
\end{aligned}
\]
There are two possible choices for $\DMAP{\tau_1}{n}{n+1}$. They are
\[
\begin{aligned}
\tau_1 & = (0 \DDD{,} i-1,i,i+2\DDD{,}n+1)=\del_{i+1} \quad \text{(the minimal choice)}\\
\tau_1 & = (0 \DDD{,} i-1,i+1,i+2\DDD{,}n+1)=\del_{i}
\end{aligned}
\]
Taking $\tau_1$ to be the minimal choice, then
\[
\begin{aligned}
\alpha_1 & = \tau_1 \LAM_1\SH = \DMAP{(i)}{0}{n+1}\\
\alpha_2 &= \DMAP{(i+1)}{0}{n+1}
\end{aligned}
\]
It follows from $\alpha_1 = \LAM\SH h_{12}\SH$ that $h_{12}\SH = \DMAP{(0)}{0}{1}$ and $h_{21}\SH = \DMAP{(1)}{0}{1}$

The trio is $\ATRIO{x}{X}{F_0(v_1)}{A_1}{F_0(v_2)}{A_2}$ where
\[
\begin{aligned}
X &= \mu\SH[k] = \SETT{0\DDD{,}i-1,i+2 \DDD{,} n+1}\\
A_1 &= \alpha_1[0] = \{i \}\\
A_2 &= \alpha_2[0] = \{i+1\}
\end{aligned}
\]
If $(w,v) \in C^{x,F}_1$ is a sponsor of this trio then the schematic diagram for $(w,v)$ is
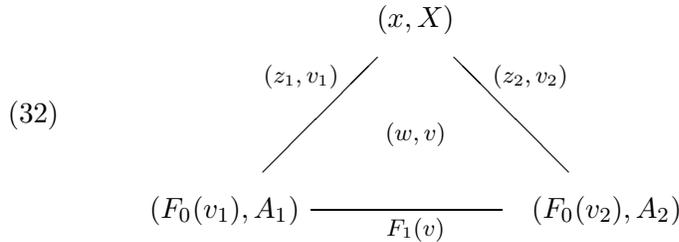
\begin{center}
\begin{pictr}{2}{1}{-1}{.5}
\label{trio for n-1,0,0}
\put(0,0){\MB{(F_{0}(v_1),A_1)}}
\put(2,0){\MB{(F_{0}(v_2),A_2)}}
\put(1,1){\MB{(x,X)}}
\put(.2,.2){\line(1,1){.6}}
\put(1.8,.2){\line(-1,1){.6}}
\put(.45,0){\line(1,0){1}}
\put(.4,.7){\MBS{(z_1,v_1)}}
\put(1.6,.7){\MBS{(z_2,v_2)}}
\put(1,-.1){\MBS{F_{1}(v)}}
\put(1,.4){\MBS{(w,v)}}
\end{pictr}
\end{center}
\MS

\NI where $(z_1,v_1) = d_{A_2}(w,v)$, $(z_2,v_2) = d_{A_1}(w,v)$ and $F_{j_1+j_2+1}(v) = d_X(w,v)$.
\MS

Note that $w \in C^x_1 \SBS C_{n+1}$ sponsors the trio $\ATRIO{x}{X}{F_0(v_1)}{A_1}{F_0(v_2)}{A_2}$ construed as a trio for the simplicial set $C^x$.
\MS

Finally, to say that $(z_1,v_1) \in C^{x,F}_0$ is $h_{12}$-universal means that given any compatible $(z_2,v_2)$ there is a unique $(w,v) \in C^{x,F}_1$ which sponsors $\ATRIO{x}{X}{F_0(v_1)}{A_1}{F_0(v_2)}{A_2}$ such that $(z_1,v_1) = d_{A_2}(w,v)$, $(z_2,v_2) = d_{A_1}(w,v)$ and $F_{j_1+j_2+1}(v) = d_X(w,v)$.

Observe that $w$, as a composition in $C$, has $d_i(w) = d_{A_1}(w) = z_2$ and $d_{i+1}(w) = d_{A_2}(w) = z_1$.
}
\end{example}

\subsection{Addendum: Context for trios}

Complementary subfaces come in pairs. Trios specify three subfaces of a potential $m$-simplex. A similar description applies to any number of subfaces (of a potential $m$-simplex) with disjoint vertex lists, as follows.
\BS

\begin{newdef}
{\bf ($(m,n)$-partition)}
\index{(m,n)-partition}
\index{(m,n)-partial simplex}
\index{sponsor of an (m,n)-partial simplex}
{\rm
\SL

Let $n \geq 1$, and $m \geq n$. 
\begin{enumerate}
\item 
A partition of $[m]$ into $n+1$ non-empty sets
\[
A_0 \DDD{\cup} A_n = [m]
\]
together with the indicated ordering of the sets will be called an {\bf \PART{m}{n}}. Each subset $A_t$ corresponds to a \SI\ map $[k_t] \to [m]$ where $k_t \DFAS |A_t|-1$. It follows that $m = n + \sum_{t=0}^n k_t$.

\item
Let $C$ be a simplicial set. Given an \PART{m}{n} $\SETT{A_0 \DDD{,} A_n}$, then an {\bf \PSMP{m}{n} of $C$} is a set of pairs
\[
\SETT{(y_0,A_0) \DDD{,} (y_n,A_n)}
\]
where, for each $t \in [n]$, $y_t \in C_{j_t}$.

Note: $\SETT{(y_0,A_0) \DDD{,} (y_n,A_n)}$ is a partial $m$-simplex; there are no required subface compatibility requirements among the $y_t$ because the $A_t$ are pairwise disjoint.

\item 
A {\bf sponsor} of an \PSMP{m}{n}\ $\SETT{(y_0,A_0) \DDD{,} (y_n,A_n)}$ is any simplex $w \in C_m$ such that for each $t \in [n]$, $\VL_w(y_t)=A_t$, equivalently $y_t = d_{[m]-A_t}(w)$.

\end{enumerate}
}
\BX
\end{newdef}

\begin{example}
An \PART{m}{1}\ $A_0 \cup A_1 =[m]$ with $j_0=|A_0|-1$ and $j_1 = |A_1|-1$
corresponds to the complementary vertex function pair $\DMAP{f_0}{j_0}{j_1+1}$ and $\DMAP{f_1}{j_1}{j_0+1}$ where $f_0\SH[j_0] = A_0$ and $f_1\SH[j_1]=A_1$. 
\MS

An \PART{m}{2} is a comb-trio of dimension $m$.
\end{example}

\begin{newdef}
\label{subpartition}
\index{subpartition}
{\bf (Subpartition)}
{\rm
\SL

Suppose $\{A_0 \DDD{,} A_n\}$ is an \PART{m}{n} and $B \subsetneq [n]$ where 
\[
B = \SETT{p_0 \DDD{<} p_k} \neq \MT
\]
Let
\[
A_B \DFAS \bigcup_{i \in B} A_i = \bigcup_{t=0}^k A_{p_t}, \quad \text{and } m_B \DFAS |A_B|-1
\]
and define the strictly increasing function $\DMAP{g_B}{m_B}{m}$ by the condition $g_B([m_B]) = A_B$. For each $t \in [k]$ let $A'_{p_t} \SBS [m_B]$ be the pullback in the diagram:

\begin{center}
\begin{picture}(1,1)
\put(0,0){\MB{A'_{p_t}}} 
\put(0,1){\MB{[m_B]}}
\put(1,0){\MB{A_{p_t}}}
\put(1,1){\MB{[m]}}
\put(.2,0){\vector(1,0){.6}}
\put(.25,1){\vector(1,0){.5}}
\put(0,.2){\vector(0,1){.6}}
\put(1,.2){\vector(0,1){.6}}
\put(.5,1.15){\MBS{g_B}}
\put(1.2,.5){\MBS{\text{incl.}}}
\put(-.2,.5){\MBS{\text{incl.}}}
\end{picture}
\end{center}
\MS

\NI That is, $r \in A'_{p_t} \iff g_B(r) \in A_{p_t}$.
\MS

The \PART{m_B}{k} $\SETT{A'_{p_t}: t \in [k]}$ will be called the {\bf subpartition of $\{A_0 \DDD{,} A_n\}$ determined by $B$}.
}

\BX
\end{newdef}
\MS

\begin{example}
{\rm
Consider the \PART{16}{3} $\{A_0,A_1,A_2,A_3\}$ ($m=16$ and $n=3$) where
\[
\begin{array}{l|l}
i  & A_i\\
\hline
0  & \{1,2,8,10\} \\
1  & \{4,12,14,16\} \\
2  & \{0,5,7,13,15\} \\
3  & \{ 3, 6, 9, 11\}
\end{array}
\]

For the subpartition defined by $B=\{0,2\} \SBS [3]$, we have
\[
\begin{aligned}
k & =1 \\
A_B &= A_0 \cup A_2 =\SETT{0,1,2,5,7,8,10,13,15}\\
m_B &= 8\\
g_B &= (0,1,2,5,7,8,10,13,15):[8] \to [16]\\
\end{aligned}
\]
The relevant pullbacks are ($i=0$, $i=2$):

\begin{center}
\begin{picture}(2,1.3)
\put(0,0){\MB{A'_i}} 
\put(0,1){\MB{[8]}}
\put(1,0){\MB{A_i}}
\put(1,1){\MB{[16]}}
\put(.2,0){\vector(1,0){.6}}
\put(.25,1){\vector(1,0){.5}}
\put(0,.2){\vector(0,1){.6}}
\put(1,.2){\vector(0,1){.6}}
\put(.5,1.1){\MBS{g_B}}
\put(1.2,.5){\MBS{\text{incl.}}}
\put(-.2,.5){\MBS{\text{incl.}}}
\end{picture}
\end{center}
\MS

\NI Then: $A_0 = \{ 1,2,8,10\} \IMP A'_0=\{ 1,2,5,6\}$
since $g_B(\{ 1,2,5,6\}) = \{ 1,2,8,10\}$.
$A_2=\{ 0,5,7,13,15\} \IMP A'_2=\{ 0,3,4,7,8\}$
since $g_B(\{ 0,3,4,7,8\}) = \{ 0,5,7,13,15\}$.
Therefore the subpartition is $\SETT{A'_0,A'_2}$, an \PART{8}{1} of $[8]$.
}
\end{example}
\MS

\begin{newdef}
{\bf (Complementary vertex functions of a partition)}
{\rm
\SL

Suppose $\SETT{A_0 \DDD{,} A_n}$ is a \PART{m}{n} where $n \geq 1$. For each $i \in [n]$ let $k_i \DFAS |A_i|-1$ and let $\DMAP{\hat{A_i}}{k_i}{m}$ be the strictly increasing map defined by the condition that $\hat{A_i}([k_i])=A_i$.

For any $i \neq j$ in $[n]$ let $m_{ij} \DFAS |A_i \cup A_j|-1 = k_i+k_j+1$ and let $\DMAP{g_{ij}}{m_{ij}}{m}$ be the strictly increasing map defined by the condition $g_{ij}[m_{ij}]=A_i \cup A_j$. Consider the $(m_{ij},1)$-subpartition $\SETT{A'_i,A'_j}$ of $[m_{ij}]$. The subsets $A'_i \SBS [m_{ij}]$ and $A'_i \SBS [m_{ij}]$ are represented by strictly increasing maps $h_{ij}$ and $h_{ji}$ in what we will call ``the $i,j$-diagram'':

\begin{center}
\begin{pictr}{2}{1.2}{-.3}{.5}
\label{subpart diagram}
\put(0,0){\MB{[k_i]}}
\put(2,0){\MB{[k_j]}}
\put(1,1){\MB{[m]}}
\put(1,0){\MB{[k_i+k_j+1]}}
\put(.2,0){\vector(1,0){.3}}
\put(1.8,0){\vector(-1,0){.3}}
\put(1,.2){\vector(0,1){.6}}
\put(.2,.2){\vector(1,1){.6}}
\put(1.8,.2){\vector(-1,1){.6}}
\put(.4,.55){\MBS{\hat{A_i}}}
\put(1.6,.55){\MBS{\hat{A_j}}}
\put(1.1,.5){\MBS{g_{ij}}}
\put(.3,-.1){\MBS{h_{ij}}}
\put(1.7,-.1){\MBS{h_{ji}}}
\end{pictr}
\end{center}
\MS

\NI and these correspond to a pair of complementary vertex functions 
\[
\DMAP{\FLT{h}_{ij}}{k_i}{k_j+1} \quad \text{and}\quad \DMAP{\FLT{h}_{ji}}{k_j}{k_i+1}
\]
We will refer to the set of maps $\SETT{h_{ij}:i,j \in [n],\; i \neq j}$ as {\bf belonging to the partition}.
}

\BX
\end{newdef}
\MS

\begin{newthm}
\label{extended-sum-theorem}
\index{Theorem!Extended Sum Theorem}
{\rm
{\bf (Extended Sum Theorem)}

Suppose $\SETT{A_0 \DDD{,} A_n}$ is an \PART{m}{n}. Then for each $i \in [n]$
\[
A_i\FL = \sum_{t \neq i} h_{it}\FL
\]
}
\end{newthm}

\Proof

It suffices, by symmetry, to prove the claim when $i=0$. We use induction on $n$. 
\MS

In the (trivial) case $n=1$ we have $m = k_0+k_1+1$, the ``$(0,1)$-diagram''
is
\begin{center}
\begin{pictr}{2}{1.2}{-.5}{.5}
\MBP{0}{0}{[k_0]}
\MBP{1}{0}{[k_0+k_1+1]}
\MBP{2}{0}{[k_1]}
\MBP{1}{1}{[m]}
\put(.2,0){\vector(1,0){.3}}
\MBPS{.35}{-.1}{h_{01}}
\put(1.8,0){\vector(-1,0){.3}}
\MBPS{1.65}{-.1}{h_{10}}
\put(1,.2){\vector(0,1){.6}}
\MBPS{1.15}{.4}{1_{[m]}}
\put(.2,.2){\vector(1,1){.6}}
\MBPS{.4}{.6}{A_0}
\put(1.8,.2){\vector(-1,1){.6}}
\MBPS{1.6}{.6}{A_1}
\end{pictr}
\end{center}
and the claim is immediate.
\MS

\NI In the case $n=2$ the \PART{m}{2} is a comb-trio and the Sum Lemma (Lemma \refpage{sum-lemma}) applies.
\MS

Now suppose $n \geq 3$. From the given \PART{m}{n} we form the \PART{m}{n-1} $\SETT{A'_0 \DDD{,} A'_{n-1}} = \SETT{A_0 \DDD{,} A_{n-2},\; A_{n-1} \cup A_n}$.
We set $k'_i = |A'_i|-1$ so that $k'_i=k_i$ for $i=0 \DDD{,} n-2$ and $k'_{n-1} = |A_{n-1} \cup A_n|-1 = k_{n-1}+k_n+1$.
\MS

For each $i =1 \DDD{,} n-2$, the corresponding ``$(0,i)$-diagram'' is 
\begin{center}
\begin{pictr}{2}{1}{-.5}{.5}
\MBP{0}{0}{[k'_0]}
\MBP{1}{0}{[k'_0+k'_i+1]}
\MBP{2}{0}{[k'_i]}
\MBP{1}{1}{[m]}
\put(.2,0){\vector(1,0){.3}}
\MBPS{.35}{-.1}{h'_{0i}}
\put(1.8,0){\vector(-1,0){.3}}
\MBPS{1.65}{-.1}{h'_{i0}}
\put(1,.2){\vector(0,1){.6}}
\MBPS{1.25}{.4}{A'_0 \cup A'_i}
\put(.2,.2){\vector(1,1){.6}}
\MBPS{.4}{.6}{A'_0}
\put(1.8,.2){\vector(-1,1){.6}}
\MBPS{1.6}{.6}{A'_i}
\end{pictr}
\end{center}
\BS

\NI where $A'_i=A_i$, $k'_i=k_i$, $h'_{0\,i}=h_{0\,i}$ and $h'_{i\,0}=h_{i\,0}$.
\MS

By the induction hypothesis, 
\[
(A'_0)\FL = (h'_{01})\FL \DDD{+} (h'_{0\, n-2})\FL + (h'_{0\, n-1})\FL
\]
Since $A'_i=A_i$ for $0 \leq i \leq n-2$ we get
\[
A_0\FL = h_{01}\FL \DDD{+} h_{0\,n-2}\FL + (h'_{0\, n-1})\FL
\]

Now the ``$(0,n-1)$-diagram'' is
\begin{center}
\begin{pictr}{4}{1}{-.5}{.5}
\MBP{0}{0}{[k'_0]}
\MBP{1}{0}{[k'_0+k'_{n-1}+1]}
\MBP{2}{0}{[k'_{n-1}]}
\MBP{1}{1}{[m]}
\put(.2,0){\vector(1,0){.3}}
\MBPS{.35}{-.15}{h'_{0\,{n-1}}}
\put(1.8,0){\vector(-1,0){.3}}
\MBPS{1.65}{-.15}{h'_{{n-1}\,0}}
\put(1,.2){\vector(0,1){.6}}

\put(.2,.2){\vector(1,1){.6}}
\MBPS{.4}{.6}{A_0}
\put(1.8,.2){\vector(-1,1){.6}}
\MBPS{1.75}{.6}{A_{n-1} \cup A_n} 
\MBP{2.5}{.5}{=}
\MBP{3}{0}{[k_0]}
\MBP{4}{0}{[m']}
\MBP{5}{0}{[k'_{n-1}]}
\MBP{4}{1}{[m]}
\put(3.2,0){\vector(1,0){.6}}
\MBPS{3.4}{-.15}{h'_{0\, n-1}}
\put(4,.2){\vector(0,1){.6}}
\MBPS{4.6}{-.15}{h'_{n-1\,0}}
\put(4.8,0){\vector(-1,0){.6}}
\put(3.2,.2){\vector(1,1){.6}}
\MBPS{3.4}{.6}{A_0}
\put(4.8,.2){\vector(-1,1){.6}}
\MBPS{4.75}{.6}{A_{n-1} \cup A_n}
\end{pictr}
\end{center}
\vskip1cm

\NI where $m' = k'_0+k'_{n-1}+1 =k_0 + k_{n-1}+k_n+2$.

Consider, then, the comb-trio $\SETT{A_0, A_{n-1},A_n}$ with diagram (as in \DREF{comb-trio-diag-1} on page \pageref{comb-trio-diag-1}):
\BS

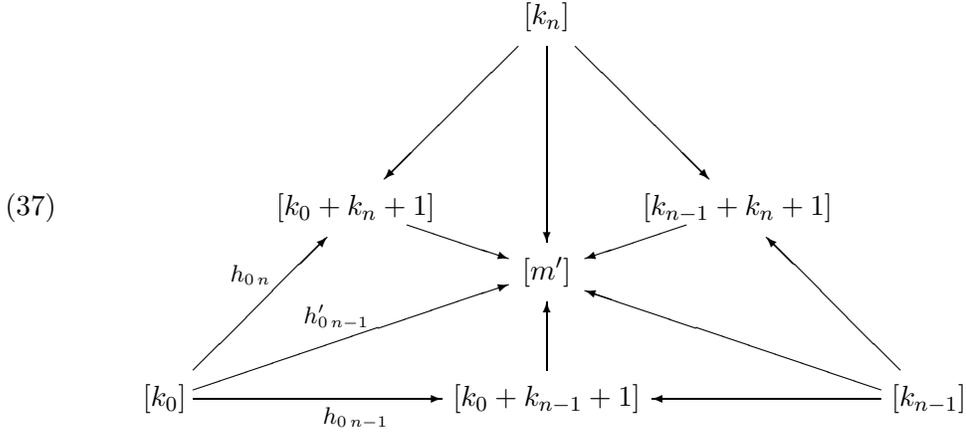
\begin{center}
\begin{pictr}{4}{2}{-.7}{1}
\MBP{0}{0}{[k_0]}
\MBP{2}{0}{[k_0+k_{n-1}+1]}
\MBP{4}{0}{[k_{n-1}]}
\MBP{1}{1}{[k_0+k_n+1]}
\MBP{2}{2}{[k_n]}
\MBP{3}{1}{[k_{n-1}+k_n+1]}
\MBP{2}{.667}{[m']}

\put(.15,.15){\vector(1,1){.7}} 
\MBPS{.45}{.65}{h_{0\,n}}
\put(.15,.05){\vector(3,1){1.65}} 
\MBPS{.9}{.43}{h'_{0\, n-1}}
\put(.15,0){\vector(1,0){1.3}}
\MBPS{1}{-.1}{h_{0\,n-1}}
 
\put(3.85,.15){\vector(-1,1){.7}}
\put(3.75,.05){\vector(-3,1){1.55}} 
\put(3.75,0){\vector(-1,0){1.2}}

\put(1.85,1.85){\vector(-1,-1){.7}}
\put(2.15,1.85){\vector(1,-1){.7}}
\put(2,1.85){\vector(0,-1){1.03}}

\put(2,.15){\vector(0,1){.36}}
\put(1.27,.91){\vector(3,-1){.53}}
\put(2.73,.91){\vector(-3,-1){.53}}
\end{pictr}
\end{center}
\BS

\NI It follows from the Sum Lemma for comb-trios that $(h'_{0\,n-1})\FL = h_{0\,n-1}\FL + h_{0\,n}\FL$ and the claim of the theorem follows from that.

\qed
\BS

\begin{example}
{\rm
Consider the following \PART{19}{4} and associated maps $h_{0i}$ and $h_{i0}$ belonging to it:
\[
\begin{array}{c|c|l|l|l}
i & k_i & A_i & \DMAP{h_{0i}}{2}{3+k_i} & \DMAP{h_{i0}}{k_i}{3+k_i} \\
\hline
0 & 2 & \{1,8,9\} & - & - \\
1 & 1 & \{4,12\} & \DMAP{(0,2,3)}{2}{4} & \DMAP{(1,4)}{1}{4}\\
2 & 3 & \{ 6,11,14,19 \} & \DMAP{(0,2,3)}{2}{6} & \DMAP{(1,4,5,6)}{3}{6}\\
3 & 5 & \{ 0,5,7,15,16,18 \} & \DMAP{(1,4,5)}{2}{8} & \DMAP{(0,2,3,6,7,8)}{5}{8}\\
4 & 4 & \{ 2,3,10,13,17 \} & \DMAP{(0,3,4)}{2}{7} & \DMAP{(1,2,5,6,7)}{4}{7}
\end{array}
\]
Then
\[
\begin{array}{c|l}
i & \DMAP{h_{0i}\FL}{2}{k_i+1} \\
\hline
1 & \DMAP{(0,1,1)}{2}{2}\\
2 & \DMAP{(0,1,1)}{2}{4}\\
3 & \DMAP{(1,3,3)}{2}{6}\\
4 & \DMAP{(0,2,2)}{2}{5}\\
\hline
\text{sum}= & \DMAP{(1,7,7)}{2}{17}
\end{array}
\]
and $A_0\FL = \DMAP{(1,8,9)\FL}{2}{17} = \DMAP{(1,7,7)}{2}{17}$.
}
\end{example}
\BS

The next theorem extends Theorem \refpage{trio-from-two-pairs} by showing that just some of the $h_{i,j}$ suffice to construct an \PART{m}{n} to which the given $h_{i,j}$ belong.

\begin{newthm}
{\rm
\SL

Suppose $n \geq 2$, $k_0 \DDD{,} k_n$ are non-negative, $j \in [n]$ fixed and $n$ \SI\ functions 
\[
\DMAP{h_{j,i}}{k_j}{k_i+k_j+1}, \quad i \in [n]-\{j\}
\]
are given. Then there exists an \PART{m}{n} to which the given $h_{j,i}$ belong, where $m = n + \sum_{p=0}^n k_p$.
}
\end{newthm}

\Proof

It suffices to take $j=0$ with $\DMAP{h_{0,i}}{k_0}{k_0+k_i+1}$, $1 \leq i \leq n$, the given \SI\ functions. The proof is by induction on $n$ with the base case $n=2$ for comb-trios established in Theorem \refpage{trio-from-two-pairs}.

Assume the claim holds for $n$ and consider any set of \SI\ functions
\[
\DMAP{h_{0,i}}{k_0}{k_0+k_i+1}, \quad 1 \leq i \leq n+1
\]
By induction, the set of functions $h_{0,i}, 1 \leq i \leq n$ belongs to an \PART{m'}{n} $\SETT{A'_0 \DDD{,} A'_n}$ where $m' = n + \sum_{p=0}^n k_p$. Denote the corresponding \SI\ subset functions by $\DMAP{\hat{A'}_i}{k_i}{m'}$. From the previous theorem we know that
\[
{\FLT{\hat{A'}_0}} = {h\FL_{0,1}}\DDD{+} {h\FL_{0,n}}
\]

Let $m \DFAS (n+1)+\sum_{p=0}^{n+1} k_p$ so that $m = m' + k_{n+1}+1$. We define the \SI\ function $\DMAP{\hat{A_0}}{k_0}{m}$ to be
\[
\FLT{\hat{A_0}} \DFAS h\FL_{0,1} \DDD{+} h\FL_{0,n+1} = {A'}\FL_0 + h\FL_{0,n+1}
\]
which is forced by the claim of the theorem. That is, for each $p \in [k_0]$,
\begin{equation}
\label{AA'-eqn}
\hat{A_0}(p) = \FLT{\hat{A'_0}}(p)+ h_{0,n+1}(p)-p
\end{equation}

Consider the diagram:
\begin{center}
\begin{pictr}{1}{1.1}{-.7}{.5}
\MBP{0}{0}{[k_0]} 
\MBP{1}{0}{[m']}
\MBP{1}{1}{[m]}
\put(.2,0){\vector(1,0){.6}}
\put(.2,.2){\vector(1,1){.6}}
\multiput(1,.2)(0,.05){12}{\MB{.}}
\put(1,.8){\vector(0,1){0}}
\MBPS{.4}{.6}{\hat{A_0}}
\MBPS{.5}{.1}{\hat{A'_0}}
\MBPS{1.15}{.5}{U_0}
\end{pictr}
\end{center}
We will apply the Extension Lemma \refpage{extension-lemma} to show that a \SI\ function $U_0$ exists such that $U_0 \hat{A'_0} = \hat{A_0}$.
\MS

\NI To check:
\MS

(1) $\hat{A'_0}(0) \leq \hat{A_0}(0)$

(2) For $0 <p \leq k_0$, $\hat{A'_0}(p)-\hat{A'_0}(p-1) \leq \hat{A_0}(p)-\hat{A_0}(p-1)$

(3) $m'-\hat{A'_0}(k_0) \leq m-\hat{A_0}(k_0)$.
\MS

\NI Condition (1) is immediate from equation \eqref{AA'-eqn}.
\MS

\NI For condition (2):
\[
\begin{aligned}
\hat{A_0}(p)-\hat{A_0}(p-1) &= \hat{A'_0}(p)-\hat{A'_0}(p-1)+h_{0,n+1}(p)-h_{0,n+1}(p-1)-1\\
& \geq \hat{A'_0}(p)-\hat{A'_0}(p-1)
\end{aligned}
\]

\NI For condition (3): 

Since $\DMAP{h_{0,n+1}}{k_0}{k_0+k_{n+1}+1}$ then $h_{0,n+1}(k_0) \leq k_0+k_{n+1}+1$. Therefore, using that $m=m' + k_{n+1}+1$ and equation \eqref{AA'-eqn}, we get:
\[
\begin{aligned}
m-\hat{A_0}(k_0) &= m'+k_{n+1}+1-\hat{A'_0}(k_0)-h_{0,n+1}(k_0)+k_0\\
&=m'-\hat{A'_0}(k_0)+k_0+k_{n+1}+1-h_{0,n+1}(k_0)\\
& \geq m'-\hat{A'_0}(k_0)
\end{aligned}
\]
Therefore the claimed $\DMAP{U_0}{m'}{m}$ exists.

To complete the proof, we define $\hat{A_i} \DFAS U_0 \hat{A'_i}$ for $1 \leq i \leq n$ and $A_{n+1} \DFAS [m]-(A_0 \DDD{\cup} A_n)$. This yields an \PART{m}{n+1} of $[m]$ containing the given $h_{0,i}$. 

\qed
\MS

\begin{example}
{\rm
$n=3$ in this example. Suppose the functions $h_{0,i}$ and $h_{i,0}$ are:
\[
\begin{array}{l|l|c|c|c}
i & k_i & \DMAP{h_{0,i}}{k_0}{k_0+k_i+1} & \DMAP{h\FL_{0,i}}{k_0}{k_i+1} & \DMAP{h_{i,0}}{k_i}{k_0+1}\\
\hline
0 & 2 & - & - & - \\
1 & 1 & \DMAP{(0,1,4)}{2}{4} & \DMAP{(0,0,2)}{2}{2} & \DMAP{(2,3)}{1}{4}\\
2 & 1 & \DMAP{(0,3,4)}{2}{4} & \DMAP{(0,2,2)}{2}{2} & \DMAP{(1,2)}{1}{4} \\
3 & 2 & \DMAP{(1,3,5)}{2}{5} & \DMAP{(1,2,3)}{2}{3} & \DMAP{(0,2,4)}{2}{5} 
\end{array}
\]
In this example, $m = 9$ and we will construct a \PART{9}{3} containing $h_{0,1}, h_{0,2}$ and $h_{0,3}$. The proof of the theorem above gives a procedure to do so.

First we construct a comb-trio containing $h_{0,1}$ and $h_{0,2}$, as follows.

Work with $k_0, k_1$ and $k_2$. Then $m' \DFAS 2+k_0+k_1+k_2 = 6$. Since $h\FL_{0,1} + h\FL_{0,2} = (0,2,4)$ we define $\hat{A'_0} \DFAS (0,2,4)\SH = \DMAP{(0,3,6)}{2}{6}$.
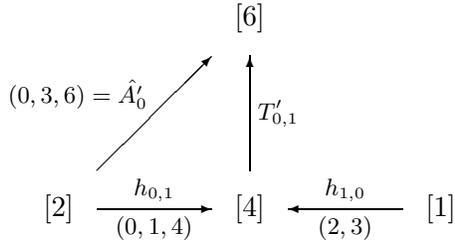
\begin{center}
\begin{picture}(2,1.2)
\MBP{0}{0}{[2]}
\MBP{1}{0}{[4]}
\MBP{2}{0}{[1]}
\MBP{1}{1}{[6]}
\put(.2,0){\vector(1,0){.6}}
\put(1.8,0){\vector(-1,0){.6}}
\put(1,.2){\vector(0,1){.6}}
\put(.2,.2){\vector(1,1){.6}}
\MBPS{.5}{-.1}{(0,1,4)}
\MBPS{.5}{.1}{h_{0,1}}
\MBPS{1.5}{-.1}{(2,3)}
\MBPS{1.5}{.1}{h_{1,0}}
\MBPS{.1}{.6}{(0,3,6)=\hat{A'_0}}
\MBPS{1.15}{.5}{T'_{0,1}}
\end{picture}
\end{center}
\MS

\NI where $T'_{0,1}$ is $\DMAP{(0,3,4,5,6)}{4}{6}$. Then $\hat{A'_2}$ must be $[6]-\{0,3,4,5,6\} = \{1,2\}$. That is, $\hat{A'_2}$ is $\DMAP{(1,2)}{1}{6}$. Then $\hat{A'_1}$ is $\DMAP{(4,5)}{1}{6}$. So the comb-trio defined is the \PART{6}{2} $A'_0 \cup A'_1 \cup A'_2 = [6]$.

Now the \PART{9}{3} to be defined must have
\[
\FLT{\hat{A_0}} = h\FL_{0,1} + h\FL_{0,2} + h\FL_{0,3} = (1,4,7)
\]
and therefore $\hat{A_0} = (1,4,7)\SH = \DMAP{(1,5,9)}{2}{9}$.
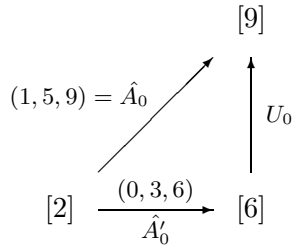
\begin{center}
\begin{picture}(1,1.1)
\MBP{0}{0}{[2]}
\MBP{1}{0}{[6]}
\MBP{1}{1}{[9]}
\put(.2,0){\vector(1,0){.6}}
\put(.2,.2){\vector(1,1){.6}}
\put(1,.2){\vector(0,1){.6}}
\MBPS{.5}{-.1}{\hat{A'_0}}
\MBPS{.5}{.1}{(0,3,6)}
\MBPS{.1}{.6}{(1,5,9)=\hat{A_0}}
\MBPS{1.15}{.5}{U_0}
\end{picture}
\end{center}
\BS

\NI where one (but not the only) choice of $U_0$ is $\DMAP{(1,2,3,5,6,7,9}{6}{9}$. Then we define $\hat{A_1} \DFAS U_0 \hat{A'_1} = \DMAP{(6,7)}{1}{9}$, $\hat{A_2} \DFAS U_0 \hat{A'_2} = \DMAP{(2,3)}{1}{9}$ which forces $\DMAP{\hat{A_3}=(0,4,8)}{2}{9} $. It is easy to check directly that this \PART{9}{3} contains the functions $h_{0,1}, h_{0,2}$ and $h_{0,3}$.
\makebox(0,0){ }
}
\end{example}


\section{Appendix: Simplicial algebra background}
\label{simpAlg}

This section contains background for the simplicial algebra used in this
paper.
Some standard references are \cite{maclane}, \cite{curtis-e-b} and \cite{may-j-p}.
\MS

\NI {\bf Contents:}

\begin{tabular}{lcr}
The simplicial category & ..... & \pageref{the-simplicial-category} \\
Simplicial object & ..... & \pageref{simplicial-object} \\
Generating maps of $\Delta$, faces and subfaces & ..... &  \pageref{generating-maps-of-delta} \\
$n$-simplex geometry & ..... &  \pageref{n-simplex-geometry} \\
Simplicial identities & ..... &  \pageref{simplicial-identities} \\
Simplicial images and degenerate simplices & ..... &
\pageref{simplicial-images-and-degenerate simplices} \\
Generators and simplicial maps & ..... &  \pageref{generators-of-a-simplicial-set} \\
Simplicial maps in terms of non-degenerate simplices
 & ..... &  \pageref{simp-maps-in-terms-of-non-degens} \\
Truncated simplicial sets
 & ..... &  \pageref{truncated-simplicial-sets} \\
Simplicial kernel
 & ..... &  \pageref{simplicial-kernel} \\
Open $i$-horn
 & ..... &  \pageref{open-i-horn} \\
The standard $m$-simplex $\Delta[m]$ and standard $i$-horn $\Lambda^i[m]$
 & ..... &  \pageref{standard-m-simplex} \\
Simplices in $\Delta[n]$
 & ..... &  \pageref{delta-m-simplices} \\
Function complexes
 & ..... &  \pageref{function-complexes}
\end{tabular}
\BX

\begin{enumerate}

\item {\bf The simplicial category} 
\index{simplicial category}
\label{the-simplicial-category}

The rules of simplicial algebra are encoded in the {\bf simplicial
category} $\Delta$ whose objects are the finite linearly ordered sets
$[n] = \lbrace 0,1,2, \cdots ,n \rbrace$, $n \in \mathbb{N}$, and whose maps are the
non-decreasing functions $f : [n] \to [m]$.

\item {\bf Simplicial object} 
\index{simplicial object}
\label{simplicial-object}

Given any category $\mathcal{C}$, a {\bf simplicial object} in
$\mathcal{C}$ is a functor $X: \Delta^{\text{op}} \to \mathcal{C}$. It may
visualized as a diagram in $\mathcal{C}$ consisting of objects $X(n)$
(or $X_n$, typically) and, for each $f : [n] \to [m]$ in $\Delta$, a map
$X(f) : X(m) \to X(n)$. In particular, a simplicial set is a functor $X
: \Delta^{\text{op}} \to \text{Sets}$. An element of $X_n$ in this case is
called an $n$-simplex of $X$.

A {\bf simplicial map} is a natural transformation of such functors.
\index{simplicial map}

\item {\bf Generating maps of $\Delta$, faces and subfaces}
\label{generating-maps-of-delta}

The maps in $\Delta$ are generated by the following family of maps.

\begin{enumerate}

\item
For each $n \ge 1$ and each $i = 0,\cdots,n$:
\[
\del_i : [n-1] \to [n]
\text{ is the map whose image omits } i
\]
That is
\[
\del_i(t) = \left\{
\begin{array}{ll}
t & \text{ if } 0 \leq t \leq i-1\\
t+1 & \text{ if } i \leq t \leq n-1
\end{array}
\right.
\]

\item
For each $n \ge 0$ and each $i = 0,\cdots,n$:
\[
\sigma_i : [n+1] \to [n]
\text{ is the map whose image repeats } i
\]
That is
\[
\sigma_i(t) = \left\{
\begin{array}{ll}
t & \text{ if } 0 \leq t \leq i\\
t-1 & \text { if } i+1 \leq t \leq n+1 
\end{array}
\right.
\]
\end{enumerate}

The $\del_p$ and $\SIG_q$ satisfy the following identities:
\begin{eqnarray*}
\del_q \del_p &=& \del_p \del_{q-1} \quad \text{if } p<q\\
\SIG_q \del_p &=& \left\{
\begin{array}{ll}
\del_p \SIG_{q-1} & \text{if } p<q\\
1 & \text{if } p=q,q+1\\
\del_{p-1} \SIG_q & \text{if } p>q+1 
\end{array}
\right.\\
\SIG_q \SIG_p &=& \SIG_p \SIG_{q+1} \quad \text{if } p \leq q
\end{eqnarray*}

Every map $f$ of $\Delta$ factors uniquely in the form 
$\del_{i_1} \cdots \del_{i_q}\, \sigma_{j_1} \cdots \sigma_{j_p}$ where $p \geq 0$, $q \geq 0$ and $i_1 \DDD{>} i_q$ if $q \geq 2$ and $j_1 \DDD{<} j_p$ if $p \geq 2$.
\MS

If $X: \Delta^{\text{op}} \to \mathcal{C}$ is any simplicial object,
$X(\del_i)$ and $X(\sigma_i)$ are denoted by $d_i : X_{n} \to X_{n-1}$
and $s_i : X_n \to X_{n+1}$ respectively. The maps $d_i$ are called face
maps and the maps $s_i$ are called degeneracy maps. Precise notation
would require ``\; $d^{n}_i$\;'' but that never seems to be really necessary.
\MS

\label{face-of-x}
\label{subface-of-x}
\label{subface-permissible-for-n}
\label{del_A}
\label{d_A}
Given $n>0$ and an $n$-simplex $x$, then $d_0(x) \DDD{,} d_n(x)$ are the {\bf faces} of $x$. The expression $d_{p_0} \DDD{} d_{p_k}(x)$ makes sense if $k<n$, $p_k \leq n$, $p_{k-1} \leq n-1$ etc. That is, for each $0 \leq j \leq k$, $p_j \leq n-k+j$. We'll refer to $\SEQ{p}{k}$ as {\bf subface-permissible} (``for $n$'', if necessary) and call $d_{p_0} \DDD{} d_{p_k}(x)\in X_{n-k-1}$ a {\bf subface} of $x$. 

If $A=\SETT{p_0 \DDD{<} p_k}$ is subface-permissible for $n$ then abbreviate
\[
d_A \DFAS d_{p_0} \DDD{}d_{p_k}:X_n \to X_{n-k-1}
\]
and 
\[
\del_A \DFAS \del_{p_k} \DDD{} \del_{p_0}:[n-k-1] \to [n]
\]
\index{face of $x$}
\index{subface of $x$}
\index{subface-permissible for $n$}
\index{$\del_A$}
\index{$d_A$}

\item {\bf $n$-simplex geometry}
\label{n-simplex-geometry}

If $X$ is a simplicial set, then the combinatorics of $\Delta$ may be
portrayed geometrically by representing $x \in X_n$ as an $n$-dimensional oriented polyhedron with $n+1$ vertices $v_0, \cdots,v_n$. For each $i =
0,\cdots, n$, $d_i(x)$ is the $i$'th face of $x$ (the face ``opposite''
$i$) and is the $(n-1)$-simplex whose vertices are $v_0,\cdots,
v_{i-1},v_{i+1},\cdots,v_n$ (with suitable adjustments in this list in
the cases $i=0$ and $i=n$). We will sometimes say that $d_i(x)$ is
``spanned'' by those vertices.

\item {\bf ``omit $i$'' notation}\index{omit $i$ notation}

There will be many occasions to write a list of things ``$t$'' indexed
by $0,1,\cdots , i-1,i+1,\cdots,n$. The standard notation for such a
list is
\[
t_0, \cdots, \omit{t_i}, \cdots, t_n
\]
{\em without} any presumption that ``$t_i$'' exists. We will also sometimes write:
\[
\OM{t}{n}{i}
\]

\item {\bf Simplicial identities}
\label{simplicial-identities}
\index{simplicial identities}

The identities satisfied by the $\del_p$ and $\SIG_q$ in $\Delta$ carry over to the so-called simplicial identities for
the $d_p$ and $s_q$ in any simplicial object. For the record, they are:
\begin{eqnarray*}
d_p d_q &=& d_{q-1}d_p \quad \text{if } p < q\\
d_p s_q &=&
\left\{
\begin{array}{ll}
s_{q-1}d_p & \text{if }  p<q \\
1 & \text{if } p=q \text{ or } q+1\\
s_q d_{p-1} & \text{if } p>q+1
\end{array}
\right. \\
s_p s_q &=& s_{q+1} s_p \quad \text{if } p \leq q
\end{eqnarray*}

Because of the role played by the maps $\del_p$ and $\sigma_q$ in $\Delta$,
a simplicial object $X : \Delta^{\text{op}} \to \mathcal{C}$ may be visualized as a set of objects $X_n$, $n \ge 0$, and maps $X_{n}
\xleftarrow{d_p} X_{n+1} \xleftarrow{s_q} X_n$ for each $n \ge 0$,
for each $p \in [n+1]$ and each $q \in [n]$.

\item {\bf Simplicial images and degenerate simplices}
\label{simplicial-images-and-degenerate simplices}
\index{simplicial image} \index{degenerate simplex}

If $x$ is any simplex (of an appropriate dimension) of a simplicial set, then a simplex of
the form $s_{i_1} \cdots s_{i_m}d_{j_1}\cdots d_{j_n}(x)$ (allowing also
the cases where $m=0$ or $n=0$) will be called a {\bf simplicial image}
of $x$.

A simplex $x$ is {\bf degenerate} if for some $i$ and $j$, $x =
s_id_jx$. Otherwise, $x$ is said to be {\bf non-degenerate}. Here are some basic facts about degenerate simplices.
\begin{enumerate}
\item 
Suppose $u$ and $v$ are non-degenerate simplices such that 
\[
s_{j_0} \cdots s_{j_t}(u) = s_{k_0} \cdots s_{k_r}(v)
\]
Then $t=r$ and $u=v$.

The proof is as follows: The given information implies that
\[
 u = d_{j_t} \cdots d_{j_0}\;  s_{k_0} \cdots s_{k_r}(v)
\]
That $u$ is non-degenerate implies that $t \geq r$. By symmetry, $r \geq t$. Therefore $t=r$. Now apply the face-degeneracy simplicial identities the right hand side. If $d_{j_t} \cdots d_{j_0}\;  s_{k_0} \cdots s_{k_r} \neq 1$ then $u$ would be a degenerate image of $v$ which would contradict the assumption that $u$ is non-degenerate.

\item
Suppose $x,y \in C_m$ such that $s_p(x) = s_q(y)$. If either $x$ or $y$ is non-degenerate then $x=y$ and $p=q$. If both $x$ and $y$ are degenerate then there is a non-degenerate simplex $u$ and values
$j_0 , \cdots , j_t$ and $k_0 , \cdots , k_t$ such that
\[
x =  s_{j_0} \cdots s_{j_t}(u) \quad \text{and} \quad y= s_{k_0} \cdots s_{k_t}(u)
\]

The proof is as follows:
First, $s_p(x) = s_q(y)$ implies $x = d_p s_p(x) = d_p s_q(y)$. If $p<q$ or $p>q+1$ then the simplicial identities would imply that $x$ is a degenerate image of a face of $y$. Therefore, if $x$ is non-degenerate then $x=y$ follows. Also, $p=q$ because $p=q+1$ would imply $x = d_{q+2}s_{q+1}(x) = d_{q+2}s_q(y) = s_q d_{q+1}(y)$, contradicting that $x$ is non-degenerate.

Otherwise, if $x$ and $y$ are both degenerate then there exist non-degenerate simplices $u$ and $v$ such that 
\[
x = s_{k_0} \cdots s_{k_t}(u) \text{ and }
y = s_{j_0} \cdots s_{j_r}(v)
\]
and therefore
\[
s_p \; s_{k_0} \cdots s_{k_t}(u) =
s_q \; s_{j_0} \cdots s_{j_r}(v)
\]
By the first part of the lemma, this implies $r=t$ and $u=v$.

\end{enumerate}

\item{\bf Generators and simplicial maps}
\label{generators-of-a-simplicial-set}
\index{generators of a simplicial set}

Suppose $C$ is a simplicial set and $G \SBS C$ is any non-empty set of simplices of $C$. Then $G$ determines a subcomplex of $C$, denoted (here) by $C(G)$ defined to be the set of all simplicial images of elements of $G$. $G$ will be termed {\bf a set of generators of $C(G)$}.

If $f : C \to C'$ is any simplicial map and if $G \SBS \C$ is non-empty, let $G' \DFAS \SETT{f(x): x \in G}$. (Notational abbreviation: $f(x)$ is $f_m(x)$ if $x \in G \cap C_m$). Then $f$ restricts to a (surjective) simplicial map $f(G) : C(G) \to C'(G')$. 
\MS

More generally, if $C$ is a simplicial set and $G \SBS C$ is a non-empty set of simplices such that every simplex of $C$ is a simplicial image of some simplex in $G$, then we say $G$ {\bf generates} $C$. Any simplicial map $f: C \to D$ is determined by its values on any set of generators of $C$.

\item {\bf Simplicial maps in terms of non-degenerate simplices}
\label{simp-maps-in-terms-of-non-degens}

Suppose $C$ and $D$ are simplicial sets and there is a function
\[
f : \SETT{y \in C : y \text{ is non-degenerate}} \to D
\]
such that $\DIM_D(f(y)) = \DIM_C(y)$ and for all $p \in [\DIM_C(y)]$, $d_p(f(y)) = f(d_p(y))$. Then $f$ extends to a unique simplicial map $F : C \to D$.
\MS

\Proof

Suppose $y \in C_m$ is degenerate with $y = s_{u_1} \cdots s_{u_k}(y')$ and $y' \in C_{m-k}$ non-degenerate. Then we define $F(y) \DFAS s_{u_1} \cdots s_{u_k}(f(y'))$. It remains to verify that $F$ is well-defined and that $F$ is actually a simplicial map.
\MS

Clearly the definition of $F$ given above is the only possible one. That $F$ is well defined follows from the uniqueness of the representation of $y$ as $s_{u_1} \cdots s_{u_k}(y')$.
\MS

It only remains to check that $F(d_p(y)) = d_p(F(y))$ and $F(s_p(y)) = s_p(F(y))$. So let $p \in [m]$. There are two possibilities according to simplicial identities.
\begin{eqnarray*}
\text{Either } \quad d_p s_{u_1} \cdots s_{u_k} &=& s_{v_1} \cdots s_{v_k} d_q\\
\text{or } \quad d_p s_{u_1} \cdots s_{u_k} &=& s_{v_1} \cdots s_{v_{k-1}}
\end{eqnarray*}
In the first case we have
\[
\begin{aligned}
d_p(F(y)) &= F(s_{v_1} \cdots s_{v_k} d_q (y'))\\
&= s_{v_1} \cdots s_{v_k} f(d_q(y'))\\
&= s_{v_1} \cdots s_{v_k} d_q(f(y'))\\
&= s_{v_1} \cdots s_{v_k} d_q(F(y'))
\end{aligned}
\]
On the other hand
\[
\begin{aligned}
F(d_p(y)) &= F(d_p s_{u_1} \cdots s_{u_k}(y'))\\
&= F(s_{v_1} \cdots s_{v_k} d_q(y'))\\
&= s_{v_1} \cdots s_{v_k} f(d_q(y'))
\end{aligned}
\]
So $d_p(F(y)) = F(d_p(y))$, as required.
\MS

In the second case, where $d_p s_{u_1} \cdots s_{u_k} = s_{v_1} \cdots s_{v_{k-1}}$,
\[
F(d_p(y)) = s_{v_1} \cdots s_{v_{k-1}} f(y')
\]
and
\[
d_p(F(y)) = d_p s_{u_1} \cdots s_{u_k} f(y')= s_{v_1} \cdots s_{v_{k-1}} f(y')
\]
Again, $d_p(F(y)) = F(d_p(y))$, as required.
\MS

Finally
\[
F(s_p(y)) = F(s_p s_{u_1} \cdots s_{u_k} (y')) = s_p s_{u_1} \cdots s_{u_k}(f(y')) = s_p(F(y))
\]
as required.

\item {\bf Truncated simplicial sets}
\label{truncated-simplicial-sets}
\index{truncated simplicial set}

Given a simplicial object $X$ and any $n \ge 0$, the truncation of $X$
to dimension $n$, denoted $\text{Tr}_n(X)$, arises from restricting $X$
to the full subcategory of $\Delta^{\text{op}}$ whose objects are $[0], [1], \cdots,
[n]$. One may visualize this as a finite diagram of objects $X_0,
\cdots, X_n$ and finitely many maps generated by the $d_i$ and $s_j$ among those objects.

$n$-truncated simplicial sets form a category in the obvious way, and
$\text{Tr}_n$ is a functor.

\item {\bf Simplicial kernel}
\label{simplicial-kernel}
\index{simplicial kernel} 
\index{$\Delta^\bullet(n)(X)$} 
\index{coskeleton}
\index{$\text{cosk}^n(\cdot)$}
\index{$\text{Cosk}^n(\cdot)$}

Suppose $X$ is a simplicial set truncated to dimension $n$. The
set
\[
\SET{(x_0,\cdots , x_{n+1}) \in X_n^{n+2}}{\forall i,j \in [n+1]( i < j \Rightarrow d_ix_j = d_{j-1}x_i}
\]
is called the $(n+1)$'st {\bf simplicial kernel}
and we will denote this set $\Delta^\bullet(n+1)(X)$. An element of $\Delta^\bullet(n+1)(X)$ may be visualized as a list of $n$-simplices whose $(n-1)$-dimensional subfaces have compatible faces so that the list comprises the set of faces of what looks like an $(n+1)$-simplex. For example, one may visualize $(x_0,x_1,x_2) \in \Delta^\bullet(2)(X)$ as three 1-simplices (directed line segments) whose endpoints match so that they form the edges of an oriented triangle, the boundary of what looks like a 2-simplex of $X$.

The projection maps $p_i : \Delta^\bullet(n+1)(X) \to X_n$ satisfy 
$d_ip_j = d_{j-1}p_i$ whenever $0 \le i<j \le n$. Also, given any $x \in
X_n$, and any $j = 0 ,\cdots n$, we define
\[
q_j(x) = \Bigl( s_{j-1}d_0(x),\cdots, s_{j-1}d_{j-1}(x),x,x,s_jd_{j+1}(x),\cdots,s_jd_m(x) \Bigr) 
\]
It is easily verified that $q_j(x) \in \Delta^\bullet(n+1)(X)$.
Furthermore, the $q_j$ act like degeneracy maps with respect to the face and degeneracy maps of the truncated simplicial set. Thus, $X$ truncated to dimension $n$, together with $\Delta^\bullet(n+1)(X)$ and the $p_j$ and $q_j$ form a simplicial set truncated to dimension $n+1$. One may repeat this construction in all higher dimensions to obtain a simplicial set denoted $\text{cosk}^n(X)$, the ``$n$-coskeleton'' of $X$. Note therefore that for any $m>n$, $\text{cosk}^m( \text{cosk}^n(X)) = \text{cosk}^n(X)$.

The process is functorial and given any simplicial set $X$, there is a unique simplicial map $\phi : X \to \text{cosk}^n(\text{Tr}_n(X))$ which is the identity map in dimensions $\le n$. It is defined inductively in dimension $m>n$ by setting $\phi_{n+1}(x) \DFAS
(d_0(x), \cdots, d_{n+1}(x))$ and using the iterative nature of the definition of simplicial kernel to define $\phi_m$ for $m>n+1$. 
\MS

A concise way to say that a simplicial object $X$ consists of simplicial
kernels from dimension $n+1$ on up is to say $X =
\text{cosk}^n(\text{Tr}_n(X))$.

Notation: $\text{cosk}^n(\text{Tr}_n(X)) = \text{Cosk}^n(X)$.

Let $\text{SimpSets}$ denote the category of simplicial sets and simplicial maps, and let $\text{SimpSets}_{\leq n}$ denote the category of $n$-truncated simplicial sets and $n$-truncated simplicial maps. Then 
\[
\text{Tr}_n : \text{SimpSets} \to \text{SimpSets}_{\leq n}
\]
is a functor whose right adjoint is
\[
\text{cosk}^n : \text{SimpSets}_{\leq n} \to \text{SimpSets}
\]
The simplicial map $X \to \text{cosk}^n(\text{Tr}_n(X))$ is the unit of that adjunction evaluated at $X$.

Now suppose $Y$ is any simplicial set with $Y_m = \Delta^\bullet(m)(Y)$ for all $m >n$; that is $Y =\text{cosk}^n(\text{Tr}_n(Y)$. If $X$ is any simplicial set, then any truncated simplicial map $\text{Tr}_n(X) \to \text{Tr}_n(Y)$ extends to a unique simplicial map $X \to \text{cosk}^n(\text{Tr}_n(Y))$, namely
\[
X \to \text{cosk}^n(\text{Tr}_n(X)) \to \text{cosk}^n(\text{Tr}_n(Y)) = Y
\]

\item {\bf Open $i$-horn} 
\label{open-i-horn}
\index{open $i$-horn} 
\index{$i$-horn}
\index{$\Lambda^i(n+1)(X)$}

Given a simplicial set $X$, $n \ge 0$ and $i \in  [n]$ the set
$\Lambda^i(n+1)(X)$ consists of all 
\[
(\OM{x}{n+1}{i}) \in \prod_{j=0, j \neq i}^{n+1} X_j
\]
such that whenever $j,k \in
[n+1]-\{i\}$ with $j<k$ then $d_j(x_k) = d_{k-1}(x_j)$.

An element of $\Lambda^i(n+1)(X)$ is called an {\bf ``open $i$-horn of
dimension $n+1$''} and is like an element of $\Delta^\bullet(n+1)(X)$
whose $i$'th coordinate is missing. Note the sequence of functions:\MS

$\begin{array}{lllll}
X_{n+1} & \longrightarrow & \Delta^\bullet(n+1)(X) & \longrightarrow & \Lambda^i(n+1)(X) \\
x & \mapsto & (d_0(x), \cdots , d_{n+1}(x)) & \mapsto & 
(d_0(x), \cdots , \omit{d_i(x)}, \cdots ,d_{n+1}(x))
\end{array}$

We will denote $X_{n+1} \to \Lambda^i(n+1)(X)$ by $\phi_{n+1,i}$ (omitting $X$ from the notation).

\item
{\bf The standard $m$-simplex $\Delta[m]$ and standard $i$-horn $\Lambda^i[m]$}
\label{standard-m-simplex}
\index{standard $m$ simplex}
\label{standard-i-horn}
\index{standard $i$-horn}

Let $m \geq 0$ and consider the {\bf standard $m$-simplex}, the simplicial set $\Delta[m]$, defined as follows:
\[
\begin{aligned}
(\Delta[m])_k & \DFAS \SETT{[k] \to [m]: \text{ non-decreasing}}\\
\Delta[m] & \DFAS \bigcup_{k \geq 0} (\Delta[m])_k
\end{aligned}
\]
The face and degeneracy operators for $\Delta[m]$ are, given $\DMAP{f}{k}{m}$ in $\Delta[m]_k$:
\[
\begin{aligned}
d_p(f) & \DFAS \DMAP{f \del_p}{k-1}{m}, \quad k>0 \text{ and } p \in [k]\\
s_p(f) & \DFAS \DMAP{f \SIG_p}{k+1}{m}, \quad p \in [k]
\end{aligned}
\]
In particular, the $p$'th vertex $\VERT{p}(f)$ of $f$ is 
\[
\VERT{p}(f)=\DMAP{f \del_k \DDD{} \omit{\del_p} \DDD{} \del_0}{0}{m},\; 0 \mapsto f(p)
\]
and given $[k] \XRA{f} [m] \XRA{g} [n]$ then
\[
\DMAP{\VERT{p}(gf)}{0}{n} = \DMAP{(gf(p))}{0}{n} = g\,  \VERT{p}(f)
\]

Note that $\Delta[m]$ is generated by $1_{[m]} \in \Delta[m]_m$.

Suppose $m \geq 2$ and $i \in [m]$. Consider
\[
G \DFAS \SETT{\DMAP{\del_p}{m-1}{m}: p \neq i} \SBS \Delta[m]_{m-1}
\]
The subcomplex of $\Delta[m]$ generated by $G$ is denoted $\Lambda^i[m]$. It has the following useful properties based on the fact\footnote{with apologies for the notation} that 
\[
(\OM{\del}{m}{i}) \in \BOX{i}{m}{\Lambda^i[m]} \SBS \BOX{i}{m}{\Delta[m]}
\]
Given any simplicial map $g:\Lambda^i[m] \to C$, then
\[
\Bigl(
g_{m-1}(\del_0) \DDD{,} \underset{i}{-} \DDD{,} g_{m-1}(\del_m)
\Bigr) \in \BOX{i}{m}{C}
\]
Reversing: any $(\OM{y}{m}{i}) \in \BOX{i}{m}{C}$ determines uniquely the simplicial map $g: \Lambda^i[m] \to C$ defined (on generators) for each $p \neq i$ by $g(\del_p) \DFAS y_p$. 

This shows there is a 1-1 correspondence between $\SETT{g:\Lambda^i[m] \to C}$ and the $i$-horns $\BOX{i}{m}{C}$.
\MS

Next, let $f:\Delta[m] \to C$ be any simplicial map. This corresponds, by Yoneda's lemma to $y = f_m(1_{[m]}) \in C_m$. Consider the composite simplicial map
\[
\Lambda^i[m] \XRA{\text{incl}} \Delta[m] \XRA{f} C
\]
For each generator $\del_p \in \Lambda^i[m]_{m-1}$, we have $f_{m-1}(\del_p) = d_p f_m(1_{[m]}) = d_p(y)$. That is, $f \circ \text{ incl}$ corresponds to $\phi_{m,i}(y)$.
\MS

Finally suppose $(\OM{y}{m}{i}) \in \BOX{i}{m}{C}$ corresponds to the simplicial map $g : \Lambda^i[m] \to C$. To say there exists $y \in C_m$  such that $\phi_{m,i}(y) = (\OM{y}{m}{i})$ is to say that $g$ factors through $f$:
\begin{center}
\begin{picture}(1,.6)
\put(0,0){\MB{\Lambda^i[m]}}
\put(0,.5){\MB{\Delta[m]}}
\put(1,.5){\MB{C}}
\put(0,.15){\vector(0,1){.2}}
\put(-.2,.25){\MBS{\text{incl}}}
\put(.2,.5){\vector(1,0){.6}}
\put(.5,.6){\MBS{f}}
\put(.2,.1){\vector(2,1){.6}}
\put(.5,.15){\MBS{g}}
\end{picture}
\end{center}
where $f : \Delta[m] \to C$ is defined by $f(1_{[m]})=y$.

\item {{\bf Simplices in $\Delta[n]$}}
\label{delta-m-simplices}

Suppose $\lambda \in \Delta[m]_k$, i.e. $\lambda: [k] \to [m]$ is non-decreasing. We may denote $\lambda$ (with ambiguity about $m$) by listing its values:
\[
\lambda = \Bigl( \lambda(0), \cdots , \lambda(k) \Bigr)
\] 
For any $p \in [k]$ we have
\[
(d_p(\LAM))(t) = (\lambda \del_p)(t) = \left\{
\begin{array}{ll}
\lambda(t) & \text{ if } t \in [0,p-1]\\
\lambda(t+1) & \text{ if } t \in [p,k-1]
\end{array}
\right.
\]
So
\[
d_p(\LAM)=\lambda \del_p = \Bigl(  \lambda(0), \cdots , \lambda (p-1),\; \lambda(p+1), \cdots , \lambda(k) \Bigr)
\]
It follows that the $i$'th vertex of $\lambda$ is $( \lambda(i) ) \in \Delta[m]_0$.
\BS

If $\mu \in \Delta[m]_{k+1}$ and $\mu \del_p = \lambda$ then
\[
\mu = \Bigl( \lambda(0) \DDD{,} \lambda(p-1),\; \mu(p),\; \lambda(p), \lambda(p+1) \DDD{,} \lambda(k) \Bigr)
\]
There are $\lambda(p)- \lambda(p-1)+1$ possible $\mu$ such that $\mu \del_p = \lambda$.
\BS

{\bf Now suppose $\lambda_p, \lambda_q \in \Delta[m]_k$, ($k>0$) such that $p<q$ in $[k]$ and $\lambda_q \del_p = \lambda_p \del_{q-1}$. Then there is {\em exactly one} $\mu \in \Delta[m]_{k+1}$ such that $\mu \del_p = \lambda_p$ and $\mu \del_q = \lambda_q$.}
The reason is as follows.
\MS

By definition, we have
\[
(\lambda_q \del_p)(t)= \left\{
\begin{array}{ll}
\lambda_q(t) & \text{ if } t \in [0,p-1]\\
\lambda_q(t+1) & \text{ if } t \in [p,k-1] 
\end{array}
\right.
\]
and
\[
(\lambda_p \del_{q-1})(t) = \left\{
\begin{array}{ll}
\lambda_p(t) & \text{ if } t \in [0,q-2]\\
\lambda_p(t+1) & \text{ if } t \in [q-1, k-1]
\end{array}
\right.
\]
Thus 
\[
\lambda_p(t) = \lambda_q(t) \text{ for } t \in [0,p-1] \cap [0,q-2] = [0,p-1]
\]
\[
\lambda_p(t+1) = \lambda_q(t+1) \text{ for } t \in [p,k-1] \cap [q-1,k-1] = [q-1,k-1]
\]
Now $[0,k-1] = [0,p-1] \cup [p,q-2] \cup [q-1,k-1]$ is a disjoint union where $[p,q-2]$ might be empty (in case $p=q-1$). If $q > 1+p$ then, for $t \in [p,q-2]$, $\lambda_q(t+1) = \lambda_p(t)$. That is
\[
\begin{aligned}
\lambda_q(p+1) &= \lambda_p(p)\\
\lambda_q(p+2) &= \lambda_p(p+1)\\
& \vdots  \\
\lambda_q(q-1) &= \lambda_p(q-2)
\end{aligned}
\]
\MS

As we saw above, any $\mu : [k+1] \to [m]$ with $\mu \del_p = \lambda_p$ must be
\[
\mu = \Bigl( \lambda_p(0) \DDD{,} \lambda_p(p-1), \mu(p), \lambda_p(p) \DDD{,} \lambda_p(k) \Bigr)
\]
If $\mu$ also has the property that $\mu \del_q = \lambda_q$ then
\[
\lambda_q(t) = (\mu \del_q)(t) = \left\{
\begin{array}{ll}
\mu(t) & \text{ if } t \in [0,q-1]\\
\mu(t+1) & \text{ if } t \in [q,k]
\end{array}
\right.
\]
Since $p \leq q-1$ then $\lambda_q(p) = (\mu \del_q)(p) = \mu(p)$. Therefore, if $\mu$ exists then
\[
\mu = \Bigl(\lambda_p(0) \DDD{,} \lambda_p(p-1),\; \lambda_q(p), \;\lambda_p(p) \DDD{,} \lambda_p(k) \Bigr)
\]
is the {\em only} possibility for $\mu$. We just need to check that 
\[
\lambda_p(p-1) \leq \lambda_q(p) \leq \lambda_p(p)
\]
As observed above $\lambda_p(p-1) = \lambda_q(p-1)$ therefore $\lambda_p(p-1) = \lambda_q(q-1) \leq \lambda_q(p)$. Also (by the calculation above) $\lambda_q(p) \leq \lambda_q(p+1) = \lambda_p(p)$. This verifies that the sequence above given for $\mu$ is non-decreasing. And, by construction, $\mu \del_p = \lambda_p$ and $\mu \del_q = \lambda_q$.
\MS

A related fact:
\label{delta-n-is-hypergroupoid}

\NI {\bf Lemma:}
Suppose $m \geq 0$ and $n \geq 2$. Then
\[
\phi_{n+1} : \Delta[m]_{n+1} \to \Delta^\bullet(n+1)(\Delta[m])
\]
is an isomorphism.

\Proof

Suppose $\beta = \LIST{b}{n+1} \in \Delta[m]_{n+1}$. Then
\[
\phi_{n+1}(\beta) = 
\left(
\beta \del_0 \DDD{,} \beta \del_{n+1} 
\right)
\]
where $\beta \del_p = (\OM{b}{n+1}{p})$. It follows immediately from this that $\phi_{n+1}$ is monic.

As a visual aid, we may represent $\phi_{n+1}(\beta)$ as an array whose rows are the values of $\beta \del_q$, $q \in [n+1]$. That array looks like
\[
\begin{array}{c|cccccc}
 & 0 & 1 & 2 & \cdots & n-1 & n\\
\hline
\beta \del_0 =& b_1 & b_2 & b_3 & \cdots & b_n & b_{n+1}\\
\beta \del_1 =& b_0 & b_2 & b_3 & \cdots & b_n & b_{n+1}\\
\beta \del_2 =& b_0 & b_1 & b_3 & \cdots & b_n & b_{n+1}\\
\vdots & \vdots & \vdots & \vdots & \cdots & \vdots & \vdots \\
\beta \del_{n+1} =& b_0 & b_1 & b_2 & \cdots & b_{n-1} & b_n
\end{array}
\]

To show that $\phi_{n+1}$ is surjective, suppose
\[
\LIST{\alpha}{n+1} \in \Delta^\bullet(n+1)(\Delta[m])
\]
where, for each $p \in [n+1]$, $\alpha_p=(a_{p\,0} \DDD{,} a_{p\,n})$. For each $p \in [1,n+1]$, set $b_p = a_{0\,p-1}$, and set $b_0 = a_{1\,0}$. We claim that $\phi_{n+1}\LIST{b}{n+1} = \LIST{\alpha}{n+1}$.

The $\alpha_p$ satisfy the face identities 
\[
\alpha_q \del_p = d_p(\alpha_q) = d_{q-1}(\alpha_p) = \alpha_p \del_{q-1}
\]
for all $p<q$ in $[n+1]$. In particular, the face identities $\alpha_{p+1} \del_p = \alpha_p \del_p$ for $p \in [0,n]$ imply that we may construct the array for $\LIST{\alpha}{n+1}$ sequentially starting with row 0 to give
\[
\begin{array}{r|cccccc}
 & 0 & 1 & 2 & 3 & \cdots & n \\
\hline
\alpha_0 =& b_1& b_2 & b_3 & b_4 & \cdots & b_{n+1} \\
\alpha_1 =& b_0& b_2 & b_3 & b_4 & \cdots & b_{n+1} \\
\alpha_2 =& b_0& b_1 & b_3 & b_4 & \cdots & b_{n+1} \\
\vdots & \vdots &\vdots &\vdots &\vdots & \cdots \vdots \\
\alpha_{n+1} =& b_0 & b_1 & b_2 & b_3 & \cdots & b_n
\end{array}
\]
It follows directly that $\LIST{\alpha}{n+1} = \phi_{n+1}\LIST{b}{n+1}$.

\qed

\item {\bf Function complexes}
\label{function-complexes}
\index{function complex}

If $C$ and $D$ are simplicial sets then the function complex $C^D$ is the simplicial set whose $k$-simplices (by Yoneda's lemma) are the simplicial maps
\[
D \times \Delta[k] \XRA{F} C
\]
If $k>0$ and $q \in [k]$ then
\[
d_q(F) \DFAS D \times \Delta[k-1] \XRA{1 \times \delta^*_q} D \times \Delta[k] \XRA{F} C
\]
That is, for each $j \geq 0$ and each $(x,\LAM) \in D_j \times \Delta[k-1]_j$ then $(d_q(F))_j(x,\LAM) = F_j(x,\del_q \LAM)$.
\MS

If $k \geq 0$ and $q \in [k]$ then
\[
s_q(F) \DFAS D \times \Delta[k+1] \XRA{1 \times \sigma^*_q} D \times \Delta[k] \XRA{F} C
\]
That is, for each $j \geq 0$ and each $(x,\LAM) \in D_j \times \Delta[k+1]_j$ then $(s_q(F))_j(x,\LAM) = F_j(x,\sigma_q \LAM)$.

\end{enumerate}
\newpage

\label{indexstart}
\addtocontents{toc}{\NI \NI \bf Index \hfill \pageref{indexstart}}

\printindex  

\addtocontents{toc}{\newline \NI \bf References \hfill \pageref{refs}}
\bibliographystyle{plain} \label{refs}
\bibliography{nicomp2021}

\end{document}